\documentclass[
    a4paper,
    english,
    oneside]{book}
    
\usepackage[margin=4cm]{geometry}

\usepackage{expl3}

\usepackage{amssymb, amsfonts, amsthm, mathtools,amsmath}

\usepackage[
    hidelinks, 
    unicode,
    colorlinks=true,
    linkcolor=blue,
    urlcolor=magenta,
    citecolor=blue]{hyperref}

\newcommand{\quot}[2]{{\begin{center}\textit{#1}\end{center}

\begin{flushright}#2\end{flushright}}}

\AtBeginDocument{%
\let\Im\relax \DeclareMathOperator\Im{\mathrm{im}}%
}

\usepackage{makeidx}
\makeindex
\newcommand{\idx}[1]{\textbf{#1}\index{#1}}

\newcommand{\ev}{\operatorname{ev}}
\newcommand{\triv}{{\operatorname{triv}}}
\newcommand{\GMod}{G\textrm{-Mod}}
\newcommand{\cGMod}{G\textrm{-Mod}_{\mcM}}
\newcommand{\bull}{\bullet}
\newcommand{\RMod}{R\textrm{-Mod}}
\newcommand{\RRMod}{\R\textrm{-Mod}}
\newcommand{\ModR}{\Mod_R}
\newcommand{\CondMod}{\mathbf{CondMod}}
\newcommand{\CondModR}{\CondMod_R}
\newcommand{\Mon}{\textrm{Mon}}
\newcommand{\CondRing}{\mathbf{CondRing}}

\newcommand{\DRMod}{D(R\textrm{-Mod})}
\newcommand{\cMod}{\RMod_{\mcM}}
\newcommand{\cDMod}{D(\RMod)_{\mcM}}
\newcommand{\cDRMod}{D(\RMod)_{\mcM}}
\newcommand{\liq}{\mathrm{liq}}

\newcommand{\cab}{\mathbf{CondAb}}
\newcommand{\Cond}{\mathbf{Cond}}
\newcommand{\Lan}{\operatorname{Lan}}
\newcommand{\Ran}{\operatorname{Ran}}
\newcommand{\Coim}{\operatorname{coim}}
\newcommand{\ihom}{\underline{\hom}}
\newcommand{\Ext}{\operatorname{Ext}}
\newcommand{\iExt}{\underline{\operatorname{Ext}}}
\newcommand{\Ch}{\mathrm{Ch}}
\newcommand{\tot}{\mathrm{tot}}
\newcommand{\sSet}{\mathbf{sSet}}
\newcommand{\Cat}{\mathrm{Cat}}
\newcommand{\Ani}{\mathrm{Ani}}
\newcommand{\Fun}{\mathrm{Fun}}
\newcommand{\Fin}{\mathrm{Fin}}
\newcommand{\LCA}{\mathrm{LCA}}
\newcommand{\qc}{\mathrm{qc}}
\newcommand{\EM}{\mathrm{EM}}

\usepackage{todonotes}
\usepackage{tikz-cd} 
\usepackage{quiver}

\RequirePackage{enumerate}
\RequirePackage[autostyle]{csquotes} 
\RequirePackage{xcolor}
\RequirePackage{mathrsfs}
\RequirePackage{amsfonts}
\RequirePackage{euscript}

\newcommand{\Rhom}{\mathrm{RHom}}
\newcommand{\RHom}{\mathrm{RHom}}
\newcommand{\iRHom}{\underline{\mathrm{RHom}}}

\renewcommand{\Im}{\mathrm{Im}}
\renewcommand{\hom}{\mathrm{hom}_{\mathcal{C}}}
\newcommand{\coker}{\mathrm{coker}\,}
\newcommand{\sub}{\subseteq}

\newcommand{\PSh}{\mathrm{PSh}}

\newcommand{\M}{\mathcal{M}}

\newcommand{\fin}{\mathbf{Fin}}
\newcommand{\profin}{\mathbf{Pro}(\mathbf{Fin})}
\newcommand{\prof}{\profin}
\newcommand{\prolam}{{\profin}_{\lambda}}
\newcommand{\prolak}{{\profin}_{\kappa}}
\newcommand{\f}{\mathbb{F}_2}
\newcommand{\bool}{\mathbf{Bool}}
\newcommand{\stone}{\mathbf{Stone}}
\newcommand{\CHaus}{\mathbf{CHaus}}

\newcommand{\extr}{\mathbf{Extr}}

\newcommand{\condAb}{\mathrm{CondAb}}

\newcommand{\bdis}{\beta\mathbf{Disc}}
\newcommand{\Top}{\mathbf{Top}}
\newcommand{\cond}{\mathbf{Cond}}
\newcommand{\condk}{{\mathbf{Cond}_{\kappa}}}
\newcommand{\set}{\mathbf{Set}}
\newcommand{\qcqs}{\mathbf{qcqs}}
\newcommand{\qs}{\mathbf{qs}}
\newcommand{\CAb}{\mathbf{CondAb}}
\newcommand{\Sh}{\mathrm{Sh}}
\newcommand{\Ab}{\mathbf{Ab}}
\newcommand{\eps}{\epsilon}
\newcommand{\Ind}{\operatorname{Ind}}
\newcommand{\sInd}{\operatorname{sInd}}
\renewcommand{\hom}{\mathrm{Hom}}

\newcommand{\Id}{\mathrm{Id}}

\newcommand\cA{\mathcal{A}}
\newcommand\cB{\mathcal{B}}
\newcommand\cC{\mathcal{C}}
\newcommand\cD{\mathcal{D}}
\newcommand\cE{\mathcal{E}}
\newcommand\cF{\mathcal{F}}

\newcommand\cI{\mathcal{I}}
\newcommand\cJ{\mathcal{J}}

\newcommand\cM{\mathcal{M}}
\newcommand\cN{\mathcal{N}}

\newcommand\cP{\mathcal{P}}

\newcommand\cS{\mathcal{S}}
\newcommand\cT{\mathcal{T}}
\newcommand\cU{\mathcal{U}}
\newcommand\cV{\mathcal{V}}
\newcommand\cW{\mathcal{W}}

\newcommand\mcA{{\cA}}
\newcommand\mcB{\cB}
\newcommand\mcC{{\cC}}
\newcommand\mcD{{\cD}}
\newcommand\mcE{{\cE}}
\newcommand\mcF{{\cF}}

\newcommand\mcI{\cI}
\newcommand\mcJ{\cJ}

\newcommand\mcM{\cM}
\newcommand\mcN{\cN}

\newcommand\mcP{\cP}

\newcommand\mcS{\cS}
\newcommand\mcT{\cT}
\newcommand\mcU{\cU}
\newcommand\mcV{\cV}
\newcommand\mcW{\cW}

\AtBeginDocument{
  \ifdefined\up
    
    \renewcommand{\up}{\mathrm{p}}
  \else 
    \newcommand{\up}{\mathrm{p}}
  \fi
}

\newcommand{\N}{{\mathbb{N}}} 
\newcommand{\Z}{\mathbb{Z}}
\newcommand{\Q}{\mathbb{Q}}
\newcommand{\R}{\mathbb{R}}
\newcommand{\C}{\mathbb{C}}

\newcommand{\T}{\mathbb{T}}

\newcommand{\fix}{\operatorname{fix}}

\newcommand{\id}{\operatorname{id}}
\newcommand{\im}{\operatorname{im}}

\newcommand{\orb}{\operatorname{orb}}

\newcommand{\res}{\operatorname{res}}

\renewcommand{\Im}{\operatorname{Im}}

\newcommand{\Aut}{\operatorname{Aut}}

\newcommand{\End}{\operatorname{End}}
\newcommand{\GL}{\operatorname{GL}} 

\newcommand{\Hom}{\operatorname{Hom}}

\newcommand{\cof}{\operatorname{cof}}

\newcommand{\Grp}{\mathbf{Grp}}

\newcommand{\Mod}{\mathbf{Mod}} 
\newcommand{\Ring}{\mathbf{Ring}}
\newcommand{\Set}{\mathbf{Set}}

\newcommand{\Vect}{\mathbf{Vect}} 

\newcommand{\op}{\mathrm{op}}

\let\oldforall\forall
\let\oldexists\exists
\renewcommand{\forall}{\,\oldforall\,}
\renewcommand{\exists}{\,\oldexists\,}

\renewcommand{\geq}{\geqslant}

\renewcommand{\epsilon}{\varepsilon}

\renewcommand{\phi}{\varphi}

\renewcommand{\theta}{\vartheta}

\renewcommand{\subset}{\subseteq}

\renewcommand{\tilde}[1]{\widetilde{#1}}
\renewcommand{\hat}[1]{\widehat{#1}}

\let\longimplies\implies
\renewcommand{\implies}{\longimplies}

\newcounter{proofstep}
\newcommand{\step}[1][]{
  \refstepcounter{proofstep}
  \ifthenelse{\equal{#1}{}}{%
    \par\textbf{\mystepname\,\theproofstep:}
  }{
    \par\textbf{\mystepname\,\theproofstep\,{\normalfont(#1)}:}\unskip
  }
}

\DeclareFontFamily{U}{matha}{\hyphenchar\font45}
\DeclareFontShape{U}{matha}{m}{n}{
      <5> <6> <7> <8> <9> <10> gen * matha
      <10.95> matha10 <12> <14.4> <17.28> <20.74> <24.88> matha12
      }{}
\DeclareSymbolFont{matha}{U}{matha}{m}{n}
\DeclareFontSubstitution{U}{matha}{m}{n}

\DeclareFontFamily{U}{mathx}{\hyphenchar\font45}
\DeclareFontShape{U}{mathx}{m}{n}{
      <5> <6> <7> <8> <9> <10>
      <10.95> <12> <14.4> <17.28> <20.74> <24.88>
      mathx10
      }{}
\DeclareSymbolFont{mathx}{U}{mathx}{m}{n}
\DeclareFontSubstitution{U}{mathx}{m}{n}

\DeclareMathDelimiter{\vvvert}{0}{matha}{"7E}{mathx}{"17}

\allowdisplaybreaks[2] 

\usepackage{amsthm}

\theoremstyle{plain}

  \newtheorem{theorem}{Theorem}[section]
\newtheorem{lemma}[theorem]{Lemma}
\newtheorem{corollary}[theorem]{Corollary}
\newtheorem{proposition}[theorem]{Proposition}
\newtheorem{conjecture}[theorem]{Conjecture}

\theoremstyle{definition}

\newtheorem{remark}[theorem]{Remark}
\newtheorem{warning}[theorem]{Warning}

\theoremstyle{definition}

\newtheorem{definition}[theorem]{Definition}
\newtheorem{question}[theorem]{Question}

\newtheorem{example}[theorem]{Example}

\theoremstyle{plain}

\newtheorem*{theorem*}{Theorem}
\newtheorem*{lemma*}{Lemma}
\newtheorem*{corollary*}{Corollary}
\newtheorem*{proposition*}{Proposition}
\newtheorem*{conjecture*}{Conjecture}

\newtheorem*{satz*}{Satz}
\newtheorem*{folgerung*}{Folgerung}

\theoremstyle{remark}

\newtheorem*{reminder*}{Reminder}
\newtheorem*{remark*}{Remark}
\newtheorem*{notation*}{Notation}
\newtheorem*{motivation*}{Motivation}
\newtheorem*{assumption*}{Assumption}

\theoremstyle{definition}

\newtheorem*{definition*}{Definition}
\newtheorem*{question*}{Question}
\newtheorem*{problem*}{Problem}
\newtheorem*{exercise*}{Exercise}
\newtheorem*{example*}{Example}
\newtheorem*{construction*}{Construction}

\newcommand{\lean}[1]{}
\newcommand{\discussion}[1]{}
\newcommand{\leanok}{}

\newcommand{\chapone}{}
\newcommand{\chaptwo}{}
\newcommand{\chapthree}{}
\newcommand{\chapfour}{}

\ExplSyntaxOn
\NewDocumentCommand{\uses}{m}
 {\clist_map_inline:nn{#1}{\vphantom{\ref{##1}}}%
  \ignorespaces}
\NewDocumentCommand{\proves}{m}
 {\clist_map_inline:nn{#1}{\vphantom{\ref{##1}}}%
  \ignorespaces}
\ExplSyntaxOff

\usepackage{adjustbox}
\usepackage[nottoc]{tocbibind}

\begin{document}

\pagenumbering{arabic}

\begingroup
\cleardoublepage

\begin{titlepage}
\parbox{.75\linewidth}{%
Universität Tübingen \\ 
Fachbereich Mathematik\\
Arbeitsgemeinschaft Funktionalanalysis}
\par\vspace{5cm}
\begin{center}\Large\bfseries
	\begin{tabular}{c}
	\huge{Aspects of Condensed Mathematics} \\
        \Large{From Abstract Nonsense to Ergodic Theory} \\
	\end{tabular}%
	\par\normalsize\normalfont\vspace*{2cm}%
	\begin{tabular}{c}
		\large{Noa Bihlmaier, Nick Ruoff, Philipp Schmale} \\
		 \\[20pt]
		March 14th 2025
	\end{tabular}\par\vspace*{\fill}
\end{center}
\end{titlepage}
\setcounter{page}{2}

\thispagestyle{empty}    
\cleardoublepage
\thispagestyle{empty}

\pagestyle{empty}

\tableofcontents
\cleardoublepage
\begin{center}
\emph{\enquote{What can you prove with exterior algebra that you cannot prove without it?}}
\end{center}

{\raggedright
Whenever you hear this question raised about some new piece of mathematics, be assured that you are likely to be in the presence of something important. In my time, I have heard it repeated for random variables, Laurent Schwartz' theory of distributions, ideles and Grothendieck's schemes, to mention only a few. A proper retort might be:
}

\begin{center}
\emph{\enquote{You are right. There is nothing in yesterday's mathematics that you can prove with exterior algebra that could not also be proved without it. Exterior algebra is not meant to prove old facts, it is meant to disclose a new world. Disclosing new worlds is as worthwile a mathematical enterprise as proving old conjectures.}}

\begin{flushright}
Gian-Carlo Rota in \enquote{Indiscrete Thoughts} \cite{Rota_1997}.
\end{flushright}

\end{center}

\cleardoublepage

\thispagestyle{empty}

\subsubsection*{The story behind this text}

As members of the functional analysis working group \enquote{Arbeitsgemeinschaft Funktionalanalysis} (AGFA) in Tübingen,
we were originally interested in questions arising in areas of functional analysis such as locally convex \cite{schaefer1971Topological},
Banach lattice \cite{schaefer1974Banachlattices} and spectral theory,
as well as in structured ergodic theory \cite{Eisner2016}.

However, studying these topics, it was often impossible to employ arguments that work well algebraically
due to issues with the nature of highly non-separated topological spaces;
see section~\ref{subsec:problems_in_top} for a more exhaustive list of these problems with topology.

Thus, we were very interested when we stumbled upon a blog post%
\footnote{\url{https://xenaproject.wordpress.com/2020/12/05/liquid-tensor-experiment/}}
of Peter Scholze where he -- together with Dustin Clausen -- proposed an entirely new framework, so-called \emph{condensed mathematics},
promising to resolve those issues.
Because we were already accustomed to extremally disconnected spaces,
seeing that they also play a central role in their theory
was additional motivation to take a closer look at the lecture notes \cite{scholze2019condensed, scholze2019Analytic, scholze2022complex} by Clausen and Scholze.
Fortunately, they are extremely well-written (including many explanations and intuitive discussions) and gave us the feeling of understanding the rough ideas.
Soon, we were completely convinced that the approach of Clausen and Scholze will
yield a beautiful and potent theory and that, by understanding this theory,
we should be able to tackle even old classical questions and resolve many of our problems.

Having agreed on the necessity to understand the theory, we stumbled upon a practical problem:
not knowing any category theory.
Indeed, we were already scared by the first precise definition (definition~1.2 in \cite{scholze2019condensed}) of a condensed set:
\begin{center}
A \emph{condensed set} is a sheaf of sets on $\ast_{\textrm{pro\'{e}t}}$, the pro-\'{e}tale site of a point.
\end{center}

This led to the first task:
understand and collect all the category theory necessary for understanding condensed mathematics.
Because in our working group no one (and, we feel, generally in our community few people) really enjoy category theory,
we decided to learn the needed topics on our own, which resulted in a first write-up (unfortunately in German).

After that, we were able to follow some elementary aspects and somewhat understand how to work with condensed sets.
This motivated us to approach our ultimate goal of employing this theory in the pursuit of functional analysis and ergodic theory.
Inevitably, we stumbled upon the problem of not knowing any homological algebra, and being afraid of even thinking about derived $\infty$-categories.
Having taught ourselves some of the basics of homological algebra
we finally moved towards condensing algebraic theories, talking about completeness, Banach spaces, and elementary ergodic theory.

We were thrilled to experience a huge interest in condensed mathematics from many mathematicians inside and outside of our community.
For example, we were given the opportunity to hold several talks in Wuppertal, a short lecture series in Kiel, talked to many people at several conferences, etc.
However, it always seemed to us that although many people are interested,
they often did not follow through with studying the theory for lack of some text collecting the necessary background
and that, at the same time, presents enough of the theory to convince someone coming from, e.g.,
functional analysis or ergodic theory of the applicability of these methods.

Supplying such a document was our motivation for writing this text,
which essentially also forms the joint master's thesis of the first two authors.
We hope it may be helpful to some people.

\clearpage

\textbf{Note that we are no experts, but rather master students who taught themselves the presented theory.
  There will most likely be some mathematical errors and misunderstandings on our side.
  No expert of the field has ever proofread the text.
We are extremely grateful for all sorts of feedback, remarks or corrections concerning this first version,
in particular with regard to the many questions and conjectures we have included.}
\vfill

\textsc{Noa Bihlmaier}, \href{mailto:nobi@fa.uni-tuebingen.de}{nobi@fa.uni-tuebingen.de}

\textsc{Nick Ruoff}, \href{mailto:niru@fa.uni-tuebingen.de}{niru@fa.uni-tuebingen.de}

\textsc{Philipp Schmale}, \href{mailto:philipp.p.schmale@gmail.com}{philipp.p.schmale@gmail.com}

\clearpage

\subsubsection*{Acknowledgements}
This text would not have been possible without the constant support and help of many people.

The first person we want to mention is our advisor, Rainer Nagel.
Inside and outside mathematics his constant support throughout the past years cannot be overestimated.
We are deeply grateful not only for the joint breakfast every Tuesday and Thursday, but also for the support in going to conferences, (almost) finishing papers, trying to understand \enquote{useless nonsense no one understands} (and apply it to \enquote{real problems}), and being interested not only in mathematics but also in the rest of the world, just to mention a few things.
Without him our mathematical and personal development would have been completely different.

The next person we are deeply indebted to is our second mathematical leader, Asgar Jamneshan.
Not only are we pursuing many joint projects at the intersection between category theory, topos theory and ergodic theory, his interest in the presented theory and the possibility of applications thereof is in fact one of the key motivations of this text.
He has enabled the first author to do a research visit with him in Istanbul,
go to \enquote{high-end} conferences in Chicago, Mondov\`{i} and many more, and always encouraged us to continue with our \enquote{abstract nonsense}.
His mathematical positivity and optimism helped us to not give up on seemingly too complicated topics.
Lastly, he taught us to focus on the important things -- without him this text would not have been finished for many more months.
We are very much looking forward to many joint projects in the next years.  

Furthermore, we thank Roland Derndinger for always being there whenever we needed someone to talk to -- regardless of whether the topic is condensed mathematics or personal.
His interest and positive feedback was extremely important to us.

We are grateful to Frank Loose for his support and good advice throughout our studies and many inspiring lectures.

We thank Markus Haase for the opportunity to give a miniworkshop about condensed mathematics in ergodic theory at the University of Kiel, and for many fruitful discussions about topos theory.

Moreover, we thank Ulrich Groh for his continuous guidance and many good hints to the classical literature on functional analysis and topology.

We thank our friends and fellow AGFA members Christian Alber, Robert Boehringer, Leon Duensing, Theodor M\"{u}ller, Sophia Schmidt and Michael Zimmermann
for many productive discussions and feedback on earlier versions of this text.

Lastly, we thank Christian Krause for technical support in adapting the lean-blueprint-tool from Patrick Massot, on which our {\LaTeX}-layout is based.

\thispagestyle{empty}
\cleardoublepage
\endgroup


\addtocontents{toc}{\protect\thispagestyle{empty}}

\chapter{Introduction to Category Theory}
\quot{Category theory is generalized abstract nonsense.}{Steenrod \cite{MacLane1997}}

In this chapter we give an introduction to those parts of category theory that are used in the main parts on condensed mathematics.
Most of the first sections can be found in \cite{mac2013categories},
in the Stacks project \cite{SPA2023}\footnote{\url{https://stacks.math.columbia.edu/}}
or on the nLab\footnote{\url{https://ncatlab.org/nlab/show/mathematics}}.
Other classical references include \cite{Borceux1994, Borceux2008, Borceux2008a, Mitchell1965,Freyd1964,kashiwara2005categories, Lane1992SheavesIG, Gelfand2003, Riehl}.

Further suggestions of reading include \cite{Vakil2023, Gabber2019, Lurie2009,Lurie2017, Johnstone}.
 More explanations and deeper insights can be found in \cite{Vakil, Vakila, thomas2000derived, Urbanik2019, baez1997introduction, Murfet2006, Johnstone, Bergner2009,Dowker1950,forster1999lectures, Godemont1958, griffiths2014principles, hatcher2002algebraic, Johnstone1982, metzler2003topological, Niemiro, sanchez2020homological, Vakila, Vistoli2007, Gelfand2003, shulman2016homotopy, shulman2017homotopy, riehl2023type, Riehl_Verity_2022, Hartshorne1997, Grothendieck1957}.

Before getting started, we want to clarify which axiomatic framework we employ.
This remark can be skipped with no harm, and we suggest to do so.

\begin{remark}[About axioms]
	We mostly behave as if we were working in \emph{von Neumann-Bernays-Gödel} set theory (NBG) and working with proper classes.
	However, silently we believe that this is purely an abuse of language and we work in ZFC with the additional assumption of two universes
	(the existence of uncountable strongly inaccessible cardinal numbers), the \enquote{small} and the \enquote{large} universe, and the \enquote{outer} universe.

	The first \enquote{small} universe is the universe of sets, and the second \enquote{large} universe is the universe of classes.
	Hence, whenever we say \emph{class of sets} we mean the set of small sets, and whenever we talk about a \emph{class of classes} or
	a \emph{metaclass of classes} we mean the outer set of all large sets in the middle universe.
	This allows us to be relaxed whenever we abuse language and use classes as if they were sets.

	We remark that we find it very tempting to forget about ZFC-based axiomatics and work in an $\infty$-categorical
	setting developed by Cisinski et al. in \cite{Cisinski2025} or in \cite{riehl2023type} (being a more homotopy-type-theoretical approach).
	But since we think that this is still work in progress, we stick to the classical axiomatics just described.
\end{remark}

\section{Preliminary definitions}
\quot{Perhaps the purpose of categorical algebra is to show that that which is trivial is trivially trivial.}
{Peter J.\ Freyd, allegedly}
\subsection{Categories, functors and natural transformations}

The main idea of category theory is that in many fields in mathematics one considers objects, defined by some properties,
together with structure preserving maps, so-called morphisms between them.

\begin{definition}[Category]\label{def:category}\chapone
A (locally small) \idx{category} $\mcC$ is a class of \idx{objects} $\mathrm{ob}\,\mcC$,
usually denoted $A,B,C,\dots$, together with
\begin{enumerate}[(i)]
\item an assignment $(A,B)\mapsto \hom(A,B)$ of sets of \textbf{morphisms}\index{morphism}
        (also called \textbf{arrows}\index{arrow})
        for every two objects $A,B$,
        whose elements usually are denoted by $f,g,h,\dots$,
\item \textbf{identity morphisms}\index{identity morphism} $1_A\in \hom(A,A)$
        for every object $A$, and
\item a \idx{composition map} $\circ\colon \hom(B,C)\times\hom(A,B)\to \hom(A,C)$
        for every three objects $A,B,C$ such that
        \begin{enumerate}[(1)]
            \item the identity behaves like an identity,
                which means that for all $f\in \hom(A,B)$ and $g\in\hom(B,C)$ the diagram
\begin{center}
\begin{tikzcd}
	A \\
	B & B \\
	& C
	\arrow["f"', from=1-1, to=2-1]
	\arrow["f", from=1-1, to=2-2]
	\arrow["{1_B}", from=2-1, to=2-2]
	\arrow["g"', from=2-1, to=3-2]
	\arrow["g", from=2-2, to=3-2]
\end{tikzcd}
\end{center}
                is commutative, i.e., $1_B\circ f=f$ and $g\circ 1_B=g$, and
            \item composition is associative, i.e. for all objects $A,B,C,D$
                and $f\in\hom(A,B),\; g\in \hom(B,C),\; h\in \hom(C,D)$, the diagram
\begin{center}
\begin{tikzcd}
	A \\
	B & C \\
	& D
	\arrow["f"', from=1-1, to=2-1]
	\arrow["{h\circ g}"', from=2-1, to=3-2]
	\arrow["{g\circ f}", from=1-1, to=2-2]
	\arrow["h", from=2-2, to=3-2]
	\arrow["g", from=2-1, to=2-2]
\end{tikzcd}
\end{center}
is commutative, meaning that $(h\circ g)\circ f=h\circ (g\circ f)$.
\end{enumerate}
\end{enumerate}
\end{definition}

\begin{remark}
We write $f\in \mcC$ or $C\in \mcC$ to denote any morphism or object in the category $\mcC$.
For $f\in\hom(A,B)$ we also write $f\colon A\to B$ and call $A$ the \idx{domain} of $f$
and $B$ its \idx{codomain}.

Sometimes, the identity morphism $1_A$ on $A$ is also denoted by $\Id_A$.
It is uniquely determined, using $1_{A}=1_{A}1_{A}'=1_A'$.

If confusion could arise, we sometimes also write $\hom_{\mcC}(A,B)$ for $\hom(A,B)$ if $A$ and $B$ are objects of $\mcC$.
\end{remark}

Some authors define a general category using \textit{classes} of morphisms $\hom(A,B)$ instead of \emph{sets}.
In this terminology, our definition of a category is called \idx{locally small}.
We usually omit the prefix \enquote{locally small} and just call this a \emph{category}.
For the other case, we define the following.

\begin{definition}\label{def:small_cat}\uses{def:category}
A \idx{large category} is a category with proper classes of morphisms.
A \idx{small category} is a (locally small) category with a \emph{set} of objects.
\end{definition}

Many familiar definitions of objects in mathematics define categories,
and hence there is a huge number of examples.

\begin{example}[Examples of categories]\label{ex:examples_of_categories}
The following is a list of basic examples of categories.
\begin{itemize}
\item The category $\Set$ whose objects are the sets and whose morphisms are the functions between them is the most fundamental category.
\item Classical algebraic structures like (abelian) groups, (commutative) rings, modules, algebras, over a ring or a field, etc., with their corresponding homomorphisms.
\item Topological structures such as topological, uniform, metrisable, compact, Hausdorff, etc., spaces together with continuous maps as morphisms.
\item Order theoretical structures as posets (partially ordered sets), directed sets, lattices, complete lattices, Boolean algebras with monotone or order continuous morphisms. We mention them separately even if some of these cases belong to the first item.
\item Any combination of the \enquote{mother structures} above,
e.g. locally convex vector spaces with continuous linear maps as morphisms,
or topological groups, or $C^*$-algebras, etc.
\end{itemize}

Some kind of mathematical objects directly give rise to a category, such as the following.
\begin{itemize}
    \item A \idx{semigroupoid} can be seen as a small category and conversely.
    \item A \idx{monoid} (a semigroup with unit) can be viewed as a small category with exactly one object and conversely.
    \item A \idx{groupoid} can be interpreted as a small category,
        where every morphism is invertible (see below for invertibility) and conversely.
    \item A \idx{group} $G$ can be viewed as a small category $\mathbf{B}G$
        with exactly one object, where every morphism is invertible.
    \item A \idx{partially ordered set} $\mathbf{P}$ identifies with a category
        where $\hom(p,q)$ has at most one element.
\end{itemize}
For a larger list of some categories used throughout this text,
see the table before the index.
\end{example}

Having defined the objects of study as \enquote{categories}, next we define the relevant notion of \enquote{morphism} between categories, the so-called functors.
This analogy of viewing categories as \enquote{objects} and functors as \enquote{morphisms} can be made explicit by defining the \enquote{meta-category} or \enquote{2-category} of categories,
which is a common point of view in higher category theory,
but for the moment we restrict ourselves to the classical view and regard this as a pure analogy.

\begin{definition}[Functor]\label{def:functor}\uses{def:category}
    A \idx{functor} $F\colon \mcC\to \mcD$ between two categories $\mcC$ and $\mcD$
    is a map that assigns to each object $A\in \mcC$
    an object $F(A)\in \mcD$ and to each morphism $f\in \hom_{\mcC}(A,B)$
    a morphism $F(f)\in \hom_{\mcD}(F(A),F(B))$ such that
\begin{enumerate}[(i)]
	\item the identity morphism is preserved,
	i.e. $F(1_A)=1_{F(A)}$ for all objects $A\in \mcC$, and
	\item composition is preserved,
	i.e. $F(g\circ f)=F(g)\circ F(f)$ for all morphisms $f\in \hom_{\mcC}(A,B)$
	and $g\in \hom_{\mcC}(B,C)$.
\end{enumerate}
\end{definition}

Again, many familiar constructions that translate certain structures into other structures are functors.

\begin{example}\label{ex:functors}\uses{def:functor}
	Some functors are given by the following examples.
	\begin{enumerate}[(i)]

		\item \textbf{Forgetful functors}\index{forgetful functor}, often denoted by \enquote{?} or left implicit, are functors that, as the name indicates,
		forget some structures,
		e.g., one can forget the multiplication of an algebra to obtain modules,
		and further forget the scalar multiplication to obtain a group,
		and further forget all algebraic structure to obtain a set.
		These functors map the morphisms to themselves,
		e.g., any algebra homomorphism is a map of modules.
		Similarly, every compact Hausdorff space is a completely regular space,
		and every completely regular space is a set.

		\item \textbf{Free functors}\index{free functor} like the free (abelian) group
		functor, free module functor, the Stone-\v{C}ech compactification functor, etc.,
		are functors usually starting in the category of sets and assigning to any set
		a free object in a certain category,
		as well as to any morphism its (unique) extension.
		We will make this more precise later in terms of adjunctions,
		see Section~\ref{sec:adjoint-functors}.
	\end{enumerate}
\end{example}

\begin{remark}
Often functors will be defined just in terms of the underlying map of objects, and the map on morphisms is left implicit.
\end{remark}

After having defined the relevant structures (categories) and morphisms (functors)
and having given some examples,
a typical next step is to define basic constructions.
These constructions (except for the opposite category) form a special type of (co)limits in the 2-category of categories, which we will introduce later (see \ref{def:strict_two_category}).

\begin{definition}[Elementary constructions of categories]\label{def:elementary_constructions_of_categories}\uses{def:category}
Let $\mcC$, $\mcD$ and $\mcC_i$ be categories for $i$ in some set $I$.
\begin{enumerate}
	\item The most simple construction is the \idx{opposite category}, or \idx{dual category}, $\mcC^{\mathrm{op}}$ of a category $\mcC$.
	It has the same objects but the direction of the morphisms is reversed, i.e.,
    \[
        \hom_{\mcC^{\mathrm{op}}}(A,B)=\hom_{\mcC}(B,A).
    \]
	Therefore we also have to reverse the composition of $\mcC$,
	so that for $f\in \hom_{\mcC^{\mathrm{op}}}(A,B)$ and
	$g\in \hom_{\mcC^{\mathrm{op}}}(B,C)$
	we define the composition
	\[
	    g\circ_{\mcC^{\mathrm{op}}}f\coloneqq g\circ_{\mcC} f.
	\]
	Clearly, the opposite of the opposite category is the original category,
	i.e. $(\mcC^{\mathrm{op}})^{\mathrm{op}}=\mcC$.
	Passing to the opposite category yields the so-called
	\idx{dual} concept of any concept, as we will explore later.
    
	\item The \idx{product category} $\prod_{i\in I}\mcC_i$ has as objects the tuples $(A_i)_{i\in I}$ with $A_i\in \mcC_i$, and morphisms $(A_i)_{i\in I}\to (B_i)_{i\in I}$ are given by tuples $(f_i)_{i\in I}$ with $f_i\in \hom_{\mcC_i}(A_i,B_i)$, equipped with coordinatewise composition.
	
	\item A \idx{subcategory} $\mcD$ of $\mcC$ is a category such that $\mathrm{Ob}(\mcD)\subseteq \mathrm{Ob}(\mcC)$, and for any $A,B\in\mcD$ we have $\hom_{\mcD}(A,B)\subseteq \hom_{\mcC}(A,B)$.
	Furthermore we demand that the identity morphisms agree, and that the composition on $\mcD$ is the restriction of the composition of $\mcC$.
	 
	\item Assume that for every two objects $A,B$ one is given an equivalence relation $\sim$ on $\hom_{\mcC}(A,B)$, and that the collection of these equivalence relations is stable under composition,
	i.e., for $f\sim f'$ we have $hfg\sim h f'g$ for all morphisms $g,h$, whenever the composition is defined.
	Then one can build the \idx{quotient category} $\mcC/\sim$ of $\mcC$ by the relation $\sim$ by taking the same objects,
	but defining $\hom_{\mcC/\sim}(A,B)$ to be $\hom_{\mcC}(A,B)/\sim$.

\end{enumerate}
\end{definition}

The construction of opposite categories allows to define another very central class of examples of functors, the \textbf{$\hom$-functors}\index{Hom-functors}.
\begin{example}\label{ex:hom_functor}
	\begin{enumerate}[(i)]

\item Let $\mcC$ be a category and $C\in\mcC$.
	    Then the \idx{covariant $\hom$-functor} is the functor $y_C = \hom_{\mcC}(C,-)$
	    from $\mcC$ to $\Set$ which maps any object $A$ of $\mcC$ to $\hom_{\mcC}(C,A)$
	    and a morphism $f\colon A\to B$ to its \idx{pushforward}
	    \[
	        f_*\colon \hom_{\mcC}(C,A)\to\hom_{\mcC}(C,B),\quad h\mapsto f\circ h.
	    \]
	    Using the axioms of a category,
	    it is readily checked that this assignment satisfies (i) and (ii) in the Definition~\ref{def:functor} of functors.

	    \item Dually, one defines the \idx{contravariant $\hom$-functor}
	    $y^C = \hom_{\mcC}(-,C)$ as the functor from $\mcC^\mathrm{op}$ to $\Set$
	    by assigning to an object $A$ of $\mcC$ the set $\hom_{\mcC}(A,C)$
	    and to any morphism $f\colon A\to B$ (direction as in $\mcC$) its \idx{pullback}
	    \[
	        f^*\colon \hom_{\mcC}(B,C)\to\hom_{\mcC}(A,C),\quad h\mapsto h\circ f.
	    \]
		More generally, a \idx{contravariant functor} from $\mcC$ to $\mcD$
	is a functor $\mcC^{\mathrm{op}}\to\mcD$.
		\item Combining these, the \idx{ $\hom$-functor} $\hom_{\mcC}(-,-)$ is the functor
		\[
		    \hom_{\mcC}(-,-)\colon \mcC^{\mathrm{op}}\times\mcC\to\Set,\quad (A,B)\mapsto \hom_{\mcC}(A,B).
		\]
		\end{enumerate}
\end{example}

Often a priori different constructions (functors) can be transformed into each other via morphisms on the objects in the target category.
This is the idea of natural transformations between functors, which intuitively form morphisms or homotopies between functors.
Again, this could be made precise using the language of 2-categories.

\begin{definition}[Natural transformation]\label{def:natural_transformation}\uses{def:category, def:functor}
    Let $F,\, G$ be functors between two categories $\mcC$ and $\mcD$.
    A natural transformation $\eta=(\eta_A)_{A\in \mcC}$ from $F$ to $G$ consists of
    morphisms $\eta_A\in\hom(FA,GA)$ for all $A\in \mcC$ such that
    for any morphism $f\in \hom_{\mcC}(A,B)$ the diagram
    \begin{center}
    \begin{tikzcd}
	FA & GA \\
	FB & GB
	\arrow["{\eta_A}", from=1-1, to=1-2]
	\arrow["{F(f)}"', from=1-1, to=2-1]
	\arrow["{G(f)}", from=1-2, to=2-2]
	\arrow["{\eta_B}"', from=2-1, to=2-2]
    \end{tikzcd}
    \end{center}
    is commutative, i.e., $\eta_B\circ F(f)=G(f)\circ \eta_A$.
	Often a natural transformation $\eta$ from $F$ to $G$ will be denoted
	$\eta\colon F\Rightarrow G$ or $\eta\colon F\to G$.
\end{definition}

\begin{example}
Let $\beta$ be the functor from $\Set$ to $\Top$,
which assigns to every set $S$ its Stone-\v{C}ech compactification $\beta S$
and $D$ the functor from $\Set$ to $\Top$
which puts on every set the discrete topology.
Then the embedding $S\to\beta S$ which maps $s\in S$ to the ultrafilter fixed on $s$,
is a natural transformation from the functor $D$ to the functor $\beta$.
\end{example}

Since functors and natural transformations morally are some sort of maps,
it makes sense to ask for composition of functors and natural transformations.

\begin{remark}\label{rem:composition_of_natural_transformations_and_functors}\uses{def:functor,def:natural_transformation}
	The composition (so-called \emph{vertical composition})
	of natural transformations is defined componentwise
	and yields a natural transformation,
	\begin{center}
	\begin{tikzcd}
	FA & GA & HA \\
	FB & GB & HB
	\arrow["{{\eta_A}}", from=1-1, to=1-2]
	\arrow["{{F(f)}}"', from=1-1, to=2-1]
	\arrow["{\eps_A}", from=1-2, to=1-3]
	\arrow["{{G(f)}}", from=1-2, to=2-2]
	\arrow["{H(f)}", from=1-3, to=2-3]
	\arrow["{{\eta_B}}"', from=2-1, to=2-2]
	\arrow["{\eps_B}"', from=2-2, to=2-3]
    \end{tikzcd},
	\end{center}
	and the composition of functors is defined by composition of the underlying maps
	and yields a functor.
	Furthermore, we can compose (so-called \emph{whiskering}) a natural transformation
	$\eta \colon F\to G$ between functors $F,G\colon \mcC\to\mcD$ with any functors
	$H\colon \mcE\to\mcC$ and $J\colon \mcD\to\mcE$ via
	\[
	    \eta.H\colon FH\to GH,\quad (\eta.H)_E=\eta_{HE}
    \]
	and
	\[
    	J.\eta\colon JF\to JH,\quad (J.\eta)_A=J(\eta_A).
	\]
	More generally, we can \emph{horizontally compose} natural transformations
	$\eta\colon F\to G$ for $F,G\colon\mcC\to\mcD$ and $\nu\colon H\to J$ for
	$H,J\colon\mcD\to\mcE$
	to a natural transformation $\nu\ast\eta\colon  H\circ F\to J\circ G$ via
	\[
    	(\nu\ast\eta)_A \coloneq \nu_{GA}\circ H(\eta_A) = J(\eta_A)\circ\nu_{FA},
	\]
	as in the diagram
	\begin{center}
	\begin{tikzcd}
	HFA & JFA \\
	HGA & JGA
	\arrow["{\nu_{FA}}", from=1-1, to=1-2]
	\arrow["{H(\eta_A)}"', from=1-1, to=2-1]
	\arrow["{J(\eta_A)}", from=1-2, to=2-2]
	\arrow["{\nu_{GA}}"', from=2-1, to=2-2]
\end{tikzcd}
	\end{center}
	Both horizontal and vertical composition are associative
	and have a unit (which is the same for both compositions), i.e.,
	the identity natural transformation $1_F\colon F\to F$.
	Furthermore, they are compatible.
	If we have three categories, six functors and four natural transformations,
	\begin{center}
	\begin{tikzcd}
	{\mcC} &&& {\mcD} &&& {\mcE}
	\arrow[""{name=0, anchor=center, inner sep=0}, curve={height=-30pt}, from=1-1, to=1-4]
	\arrow[""{name=1, anchor=center, inner sep=0}, curve={height=30pt}, from=1-1, to=1-4]
	\arrow[""{name=2, anchor=center, inner sep=0}, from=1-1, to=1-4]
	\arrow[""{name=3, anchor=center, inner sep=0}, curve={height=30pt}, from=1-4, to=1-7]
	\arrow[""{name=4, anchor=center, inner sep=0}, curve={height=-30pt}, from=1-4, to=1-7]
	\arrow[""{name=5, anchor=center, inner sep=0}, from=1-4, to=1-7]
	\arrow["\eps", shorten <=4pt, shorten >=4pt, Rightarrow, from=0, to=2]
	\arrow["\mu", shorten <=4pt, shorten >=4pt, Rightarrow, from=2, to=1]
	\arrow["\nu", shorten <=4pt, shorten >=4pt, Rightarrow, from=4, to=5]
	\arrow["\eta", shorten <=4pt, shorten >=4pt, Rightarrow, from=5, to=3]
    \end{tikzcd},
    \end{center}
	then
	\[
	    (\eta\circ\nu)\ast (\mu\circ \eps)=(\eta\ast\mu)\circ (\nu\ast\eps).
    \]
    For more see chapter II.5 in \cite{mac2013categories}.
\end{remark}

The first part of the remark implies that the collection of functors
between two categories forms a category, the so-called \idx{functor category}.

\begin{definition}[Functor categories]\label{def:functor_category}\uses{def:category, def:functor, def:natural_transformation}
    The class of functors between two categories $\mcC$ and $\mcD$ is denoted by
    $[\mcC,\mcD]$, $\hom(\mcC,\mcD), \mcD^\mcC$ or $\mathrm{Fun}(\mcC,\mcD)$.
	It becomes a large category with natural transformations as morphisms,
	where the identity morphisms are given by the trivial natural transformations.

	The special case $\mathrm{Fun}(\mcC^{\mathrm{op}},\Set)$ is called the category of \idx{presheaves} on $\mcC$ and denoted $\PSh(\mcC)$ or $\hat{\mcC}$.
\end{definition}

\begin{remark}
    If the category $\mcC$ is small,
    then the functor category $[\mcC,\mcD]$ is locally small.    
    Very often, although, functor categories are large categories.

    We have defined a contravariant functor from $\mcC$ to $\mcD$ as a functor
    $\mcC^\mathrm{op}\to\mcD$.
    This is equivalent to saying that it is a functor $\mcC\to\mcD^\mathrm{op}$
    However, the functor categories $[\mcC^{\op},\mcD]$ and $[\mcC,\mcD^{\op}]$
    are dual, i.e.
    \[
        [\mcC^{\op},\mcD]=[\mcC,\mcD^{\op}]^{\op}.
    \]
\end{remark}

\begin{example}[Yoneda embedding]\label{ex:yoneda}\uses{ex:functors, def:functor, def:functor_category}
\begin{enumerate}[(i)]
\item    For $\mcD = \Set$ every object $C\in\mcC$ determines a functor $y_C=\hom(C,-)$
    in $\mathrm{Fun}(\mcC,\Set)$, see example~\ref{ex:hom_functor}.
    This assignment, called the \idx{covariant Yoneda embedding}, determines a functor $\mcC^{\mathrm{op}}\to \mathrm{Fun}(\mcC,\mathrm{Set})$
    by mapping a morphism $f\colon A\to B$ (direction denoted as it were in $\mcC$)
    to its \emph{pullback} $f^*\colon y_B\to y_A$ given on $C$ by
    \[
        f^*_C\colon\hom_{\mcC}(B,C)\to\hom_{\mcC}(A,C),\quad h\mapsto h\circ f.
    \]
    One easily checks that $f^*$ is a natural transformation from $y_B$ to $y_A$,
    meaning that for every morphism $g\colon C\to D$ the diagram
    \begin{center}
    \begin{tikzcd}
	{y_B(C)} & {y_A(C)} \\
	{y_B(D)} & {y_A(D)}
	\arrow["{f^*_C}", from=1-1, to=1-2]
	\arrow["{g_*}"', from=1-1, to=2-1]
	\arrow["{g_*}", from=1-2, to=2-2]
	\arrow["{f^*_D}"', from=2-1, to=2-2]
    \end{tikzcd}
    \end{center}
    is commutative.

\item    Dually, we have the \idx{contravariant Yoneda embedding}
    $\mcC\to \mathrm{Fun}(\mcC^{\mathrm{op}},\mathrm{Set})=\PSh(\mcC)$
    as the functor sending an object $C$ to
    $y^C\coloneqq \hom_{\mcC}(-,C)\in \mathrm{Fun}(\mcC^{\mathrm{op}},\mathrm{Set})$,
    see~\ref{ex:hom_functor}, and a morphism $f\colon A\to B$ to its \emph{pushforward}
    $f_*\colon y^A\to y^B$, given on $C$ by
    \[
       (f_*)_C\colon\hom_{\mcC}(C,A)\to\hom_{\mcC}(C,B),\quad h\mapsto f\circ h.
    \]
    By the same argument as before, this determines a natural transformation.
\end{enumerate}    
Note that the covariant Yoneda embedding is a contravariant functor mapping to covariant functors, and conversely the contravariant Yoneda embedding is a covariant functor mapping to contravariant functors.\end{example}
\quot{The Yoda embedding, contravariant it is.}{The Rising Sea \cite{Vakil2023}, 1.3.12}

\begin{definition}\label{def:representable-object}\uses{ex:yoneda}
A presheaf that is isomorphic (in $\PSh(\mcC)$) to some $y^C$
is called \idx{representable}.
\end{definition}

We conclude that the basic structures of category theory are categories,
with functors as maps between categories and natural transformations as maps between functors,
which together forms a \emph{(strict) 2-category} (ignoring size-issues,
which can be solved).

\begin{definition}[strict 2-category]\label{def:strict_two_category}\uses{def:category}
	A \textbf{(strict) 2-category}\index{strict 2-category}\index{2-category} consists of a category $\mcC$
	and for any two morphisms $f,g\in \mcC$ a set of 2-morphisms $\hom_2(f,g)$,
	equipped with two associative compositions (vertical and horizontal),
    such that there exists a simultaneous identity 2-morphism $1_f$ and such that
    \begin{center}
    \begin{tikzcd}
	X &&& Y &&& Z
	\arrow[""{name=0, anchor=center, inner sep=0}, curve={height=-30pt}, from=1-1, to=1-4]
	\arrow[""{name=1, anchor=center, inner sep=0}, curve={height=30pt}, from=1-1, to=1-4]
	\arrow[""{name=2, anchor=center, inner sep=0}, from=1-1, to=1-4]
	\arrow[""{name=3, anchor=center, inner sep=0}, curve={height=30pt}, from=1-4, to=1-7]
	\arrow[""{name=4, anchor=center, inner sep=0}, curve={height=-30pt}, from=1-4, to=1-7]
	\arrow[""{name=5, anchor=center, inner sep=0}, from=1-4, to=1-7]
	\arrow["\eps", shorten <=4pt, shorten >=4pt, Rightarrow, from=0, to=2]
	\arrow["\mu", shorten <=4pt, shorten >=4pt, Rightarrow, from=2, to=1]
	\arrow["\nu", shorten <=4pt, shorten >=4pt, Rightarrow, from=4, to=5]
	\arrow["\eta", shorten <=4pt, shorten >=4pt, Rightarrow, from=5, to=3]
    \end{tikzcd}
    \end{center}
	commutes, i.e., $(\eta\circ\nu)\ast (\mu\circ \eps)=(\eta\ast\mu)\circ (\nu\ast\eps)$.
\end{definition}

In particular for all $A,B$, $(\hom(A,B),\hom_2(-,-))$ again is a category, and we could have defined a 2-category using this enrichment.
The collection of all small (due to size issues) categories forms a strict 2-category
by the observations made in
Remark~\ref{rem:composition_of_natural_transformations_and_functors}.
For more on this see \cite{baez1997introduction,baez2023}
and chapter XII in \cite{mac2013categories}.

\subsection{Morphisms}

For now lets turn back to ordinary categories.
First we will define special types of morphisms.

\begin{definition}[Special morphisms]\label{def:special_morphisms}\uses{def:category}
    A morphism $f\colon A\to B$ in a category $\mcC$ is called
    \begin{enumerate}[(i)]
	    \item a \idx{monomorphism} or \idx{monic}, if for any two morphisms $g,h\colon C\to A$ the equality $f\circ g=f\circ h$ implies $g=h$,
    	\item a \idx{split monomorphism} or a \idx{section}, if there exists a morphism $g\colon B\to A$ such that $g\circ f=1_A$ (which is called a \idx{retraction} of $f$),
    	\item an \idx{epimorphism} or \idx{epic}, if for any two morphisms $g,h\colon B\to D$ the equality $g\circ f=h\circ f$ implies $g=h$,
	    \item a \idx{split epimorphism} or a \idx{retraction}, if there exists a morphism $g\colon B\to A$ such that $f\circ g=1_B$ (which is called a \idx{section} of $f$),
    	\item a \idx{bimorphism}, if it is both a monomorphism and an epimorphism,
    	\item an \idx{isomorphism}, if there exists a morphism $g\colon B\to A$ such that $f\circ g=1_B$ and $g\circ f=1_A$.
	    If there exists an isomorphism $A\to B$ we call $A$ and $B$ \idx{isomorphic} and write $A\simeq B$.
\end{enumerate}
\end{definition}

\begin{remark}
    In the case of (ii) and (iv), one also says that $A$ is a \idx{retract} of $B$.
    
    In the case of (vi), $g$ is referred to as the \idx{inverse} of $f$.
    As usual, an inverse is, if it exists, uniquely determined.
\end{remark}

\begin{remark}
	\begin{enumerate}[(i)]
		\item
    Clearly, the properties of being split monic, split epic or being an isomorphism are preserved under arbitrary functors,
	while the classes of epimorphisms, of monomorphisms and of bimorphisms do not need to be preserved.
    For example, a continuous map with dense image is epimorphic in the category of Hausdorff spaces with continuous maps,
    but applying the forgetful functor to $\Set$ yields a non-surjective map
    which is henceforth not epimorphic.

		\item
	One can easily check that every split monomorphism is monic,
	every split epimorphism is epic,
	and that every isomorphism is a bimorphism.
    However the converse does not hold,
    as e.g. in the category of topological spaces a continuous bijection
    does not need to be a homeomorphism, i.e., have a continuous inverse.
\item The properties of (split) monics and (split) epimorphisms are dual concepts, i.e., a morphism $f$ is (split) monic in $\mcC$ precisely if the reversed arrow is (split) epimorphic in $\mcC^\op$.
\item
    In the category of sets these subtle differences disappear
    and one recovers the usual notion of injective, surjective and bijective maps.
	There, the notions of monomorphism and split monomorphism coincide,
	and agree with the injective maps.
	The class of epimorphisms agrees with the split epimorphisms and the surjective maps,
	which is equivalent to the axiom of choice.
	The isomorphisms agree with the bimorphisms and the bijective maps.
\end{enumerate}
\end{remark}

Categories, where bimorphisms are invertible, have a special name.

\begin{definition}\label{def:balanced-category}\uses{def:special_morphisms}
A category $\mcC$ is called \idx{balanced} if every bimorphism is an isomorphism.
\end{definition}

In non-balanced categories, we can only expect an inverse of a bimorphism
if there already exists a splitting in one direction.

\begin{lemma}[Monic and split epi is iso]\label{lem:split_epi_mono_is_iso}\uses{def:special_morphisms}
	For any morphism $f\colon A\to B$ in a category $\mcC$ the following conditions are equivalent.
\begin{enumerate}[(a)]
	\item The morphism $f$ is an isomorphism.
	\item The morphism $f$ is a split monomorphism and an epimorphism.
	\item The morphism $f$ is a split epimorphism and a monomorphism.
\end{enumerate}
\end{lemma}

\begin{proof}
  Clearly an isomorphism is both a split monomorphism and a split epimorphism.

	Let $f\colon A\to B$ be a morphism that is both a split monomorphism and an epimorphism.
	Then there exists a morphism $g\colon B\to A$ such that $g\circ f=1_A$ and for any morphisms $h,\,k\colon B\to D$ the equality $h\circ f=k\circ f$ implies $h=k$.
	Now the equation
\[(f\circ g)\circ f=f\circ (g\circ f)=f=1_B\circ f\]
implies $f\circ g=1_B$, and thus $g$ is the inverse to $f$.

Analogously, if $f$ is a split epimorphism and a monomorphism, then there exists a morphism $g\colon B\to A$ such that $f\circ g=1_B$
and for any morphisms $h,\,k\colon C\to A$ the equality $f\circ h=f\circ k$ implies $h=k$.
Again,
	\[f\circ (g\circ f)=(f\circ g)\circ f=f= f\circ 1_A\]
	leads to $g\circ f=1_A$ and thus $f$ is an isomorphism with inverse $g$.
\end{proof}

Moreover, certain \emph{nice} retracts preserve isomorphism,
as the following lemma shows, which will be useful later.

\begin{lemma}\label{lem:retract-diagram-iso}
In any category $\mcC$, if in the diagram
\begin{center}
\begin{tikzcd}
	A && B \\
	\\
	C && D
	\arrow["f", from=1-1, to=1-3]
	\arrow["s", shift left, from=1-1, to=3-1]
	\arrow["t", shift left, from=1-3, to=3-3]
	\arrow["p", shift left, from=3-1, to=1-1]
	\arrow["g", from=3-1, to=3-3]
	\arrow["q", shift left, from=3-3, to=1-3]
\end{tikzcd}
\end{center}
both squares commute, i.e. $t\circ f = g\circ s$ and $f\circ p = q\circ g$,
and if furthermore $p\circ s = 1_A$, $q\circ t = 1_B$,
i.e. $A$ is a retract of $C$ and $B$ is a retract of $D$,
then if $g$ is an isomorphism, so is $f$.
\end{lemma}

\begin{proof}
Let $g$ be an isomorphism with inverse $g^{-1}$.
We show that $p\circ g^{-1}\circ t$ is an inverse of $f$.
First,
\[
    p\circ g^{-1}\circ t\circ f = p\circ g^{-1}\circ g\circ s = p\circ s = 1_A
\]
and second,
\[
    f\circ p\circ g^{-1}\circ t = q\circ g\circ g^{-1}\circ t = q\circ t = 1_B.
\]
\end{proof}

\begin{definition}\label{def:subobject}\uses{def:special_morphisms}
    A \idx{subobject} of an object $A$ in a category $\mcC$ is a monomorphism
    $B\hookrightarrow A$ for an object $B$.
\end{definition}

\begin{remark}
    The class of subobjects of a fixed object $A$ yields a category by defining
    morphisms from $B\hookrightarrow A$ to $C\hookrightarrow A$ to be morphisms
    $B\to C$ such that the diagram
    \begin{center}
    \begin{tikzcd}
	& A \\
	B && C
	\arrow[hook, from=2-1, to=1-2]
	\arrow[from=2-1, to=2-3]
	\arrow[hook', from=2-3, to=1-2]
    \end{tikzcd}
    \end{center}
    is commutative.
    In particular, $B\to C$ is a monomorphism.
\end{remark}

\begin{warning}
	Sometimes we identify a whole isomorphism class of subobjects
    to be a subobject, which is a common abuse of language,
    see chapter V.7 in \cite{mac2013categories}.
\end{warning}
This has the advantage that although often there is a class of monics, there is just a set of equivalence classes thereof.

\begin{definition}[well powered]
We say that a category $\mcC$ is \idx{well-powered},
if every object $A$ in $\mcC$ has only a (\emph{small}) set
of isomorphism classes of subobjects.
\end{definition}

The category of subobjects of a given object is a special case of the more general concept of comma categories which we will now introduce.

\begin{definition}[Comma category]\label{def:comma_category}\uses{def:category}
	Let $\mcC$ and $\mcD$ be categories and $F\colon \mcC\to \mcD$ a functor,
	and $D\in \mcD$.
	The \idx{comma category} $(F\downarrow D)$  has as objects tuples $(C, f)$,
	where $C\in \mcC$ and $f\colon F(C)\to D$,
	and morphisms $(C,f)\to (C',f')$ are given by $g\colon C\to C'$ such that the diagram
    \begin{center}\begin{tikzcd}
	{F(C)} && {F(C')} \\
	& D
	\arrow["{F(g)}", from=1-1, to=1-3]
	\arrow["f"', from=1-1, to=2-2]
	\arrow["{f'}", from=1-3, to=2-2]
\end{tikzcd}\end{center}
	is commutative.

	Dually, the comma category $(D\downarrow F)$ has as objects tuples $(C, f)$,
	where $C\in \mcC$ and $f\colon D\to F(C)$,
	and morphisms $(C,f)\to (C',f')$ are given by $g\colon C\to C'$ such that the diagram
\begin{center}\begin{tikzcd}
	& D \\
	{F(C)} && {F(C')}
	\arrow["f"', from=1-2, to=2-1]
	\arrow["{f'}", from=1-2, to=2-3]
	\arrow["{F(g)}"', from=2-1, to=2-3]
\end{tikzcd}\end{center}
is commutative.
\end{definition}

\begin{remark}
  Note that two morphisms $f,g$ in the comma category are distinct as soon as they are distinct morphisms in $\mcC$, although the morphisms $F(f)$ and $F(g)$ might coincide as morphisms in $\mcD$.

Often, especially when $F$ is the identity functor on $\mcC$ and $C$ an object of $\mcC$,
we write $F/C = (F\downarrow C)$ and $C/F = (C\downarrow F)$.
(In the case $F=\id_{\mcC}$ we write $\mcC/C$ resp.\ $C/\mcC$)
\end{remark}
	
\begin{example}
The category of subobjects of an object $C\in \mcC$ can be recovered using the identity functor on the subcategory of $\mcC$ having only the monomorphisms of $\mcC$ as morphisms.
\end{example}

A definition of injectivity and surjectivity for functors is slightly more involved,
since here we are dealing with maps between objects and maps between sets of morphisms.

\begin{definition}[Full, faithful, essentially surjective functors]\label{def:full_faithful_essentially_surjective_functors}\uses{def:functor, def:special_morphisms}
	A functor $F\colon \mcC\to \mcD$ is called
	\begin{enumerate}[(i)]
		\item \idx{full}, if for any two objects $A,B\in \mcC$ the map $\hom_{\mcC}(A,B)\to \hom_{\mcD}(F(A),F(B))$ is surjective,
		\item \idx{faithful}, if this map is injective,
		\item \idx{fully faithful} if it is both full and faithful,
		\item \idx{essentially surjective}, if for any object $D\in \mcD$ there exists an object $C\in \mcC$ such that $F(C)\cong D$.
		\item an \idx{embedding} if it is faithful and injective on objects.
	\end{enumerate}
\end{definition}

\begin{warning}
	Note that while the naive image of a functor
	(meaning the set-theoretic image of the object maps as well as the $\hom$ set maps)
	does not have to be a category,
	the image of a full functor always yields a category.
	This motivates the definition of the \idx{essential image} to be the smallest subcategory of $\mcD$ containing the naive image and being closed under isomorphisms.
\end{warning}

\begin{proof} To see that the image of a functor is not necessarily a category consider e.g. the category $\mcC$ induced by
\begin{center}\begin{tikzcd}
	\bullet & \bullet & \bullet & \bullet \\
	{}
	\arrow[from=1-1, to=1-2]
	\arrow[from=1-3, to=1-4]
\end{tikzcd}\end{center}
as well as the category $\mcD$ via
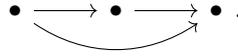
\begin{center}\begin{tikzcd}
	\bullet & \bullet & \bullet
	\arrow[from=1-1, to=1-2]
	\arrow[from=1-2, to=1-3]
	\arrow[curve={height=18pt}, from=1-1, to=1-3]
\end{tikzcd}.\end{center}
Consider the functor $F\colon\mcC\to\mcD$ induced by the dashed arrows as follows
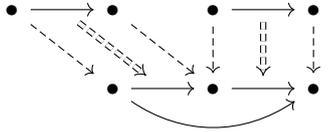
\begin{center}\begin{tikzcd}
	\bullet & \bullet & \bullet & \bullet \\
	& \bullet & \bullet & \bullet
	\arrow[""{name=0, anchor=center, inner sep=0}, from=2-2, to=2-3]
	\arrow[""{name=1, anchor=center, inner sep=0}, from=2-3, to=2-4]
	\arrow[curve={height=18pt}, from=2-2, to=2-4]
	\arrow[""{name=2, anchor=center, inner sep=0}, from=1-1, to=1-2]
	\arrow[""{name=3, anchor=center, inner sep=0}, from=1-3, to=1-4]
	\arrow[dashed, from=1-1, to=2-2]
	\arrow[dashed, from=1-2, to=2-3]
	\arrow[dashed, from=1-3, to=2-3]
	\arrow[dashed, from=1-4, to=2-4]
	\arrow[shorten <=7pt, shorten >=7pt, Rightarrow, dashed, from=2, to=0]
	\arrow[shorten <=4pt, shorten >=4pt, Rightarrow, dashed, from=3, to=1]
\end{tikzcd}.\end{center}
Clearly, the naive image of ${F}$ is not a category, because the composition of the two arrows in the image is no longer in the image.

If on the other hand the functor is full, for any 
\[F(f)\in F(\Hom(A,B))=\Hom(F(A),F(B))\]
 and \[F(g)\in F(\Hom(C,D))=\Hom(F(C),F(D))\] with $F(C)=F(B)$ there exists a preimage $j\in \Hom(B,D)$ with $F(j)=F(g)$ and hence
for $F(g)\circ F(f)\in \Hom(F(A),F(D))=F(\Hom(A,D))$ there exists a preimage $h=j \circ f\in \Hom(A,D)$ such that $F(g)\circ F(f)=F(h)\in F(\Hom(A,D))$, which implies that the restriction of the composition is well-defined in the image, hence the image forms a subcategory.
\end{proof}

Next, as an example of (full) embeddings of categories, we will show the standard way of embedding a category into a functor category on this category.
This is one of the most fundamental concepts in category theory and has applications and appearances in many other fields of mathematics.
The reader is encouraged to not think too hard about this lemma until one has seen it being used, as we have the feeling that one understands it best while applying it (it took us quite a while until this result felt natural).
The proof is surprisingly straight-forward, and can be found in any standard text on category theory, e.g. in \cite[2.2.4, 2.2.8]{Riehl}, in chapter III.2 in \cite{mac2013categories}, in chapter 1.4 in \cite{kashiwara2005categories} and in chapter 1.3 in \cite{Borceux2008}.

\begin{lemma}[\idx{Yoneda lemma}]\label{lem:yoneda}\uses{def:functor_category, def:functor, def:natural_transformation, def:special_morphisms, def:category, ex:functors}
\begin{enumerate}[(i)]
	\item Consider the {covariant Yoneda embedding} (see example~\ref{ex:yoneda} (i))
	\[
    	\mcC^{\mathrm{op}}\to \mathrm{Fun}(\mcC,\mathrm{Set}),\quad C\mapsto y_C.
	\]
	For any functor $F\colon \mcC\to \mathrm{Set}$ there is an isomorphism
	 \[
	     F(A)\simeq \hom(y_A,F),
	\]
	natural in $A$ and in $F$.
	In particular, taking $F = y_B$ yields $\hom(y_A,y_B)\simeq\hom_{\mcC}(B,A)$
	so that the covariant Yoneda embedding is fully faithful.
	
	\item Similarly, consider the {contravariant Yoneda embedding} (see example~\ref{ex:yoneda} (ii))
	\[
    	\mcC\to\PSh(\mcC),\quad C\mapsto y^C
	\]
	For any presheaf $F\colon \mcC^{\mathrm{op}}\to\mathrm{Set}$
	there is an isomorphism
	\[
    	F(A)\simeq \hom(y^A,F)
	\]
    natural in $A$ and in $F$.
	In particular, taking $F=y^B$ yields $\hom(y^A,y^B)\simeq\hom_{\mcC}(B, A)$
	so that the contravariant Yoneda embedding is fully faithful.
\end{enumerate}
\end{lemma}

\begin{proof}
In both cases, the isomorphism is given by
\[
    (\eta\colon y_A\Rightarrow F)\mapsto \eta_A(1_A)\in F(A)
\]
(resp. $y^A$ in the contravariant case).
If one is already familiar with the Yoneda lemma,
the proof is straightforward.
If one is not yet familiar with the Yoneda lemma,
to get a better understanding,
we recommend doing all of the needed verifications yourself.
\end{proof}

While one can install the naive term of isomorphism on functors (corresponding to the 1-metacategory of categories),
the presence of natural transformations as morphisms between functors allows us to define a slightly more refined version of isomorphism,
so-called equivalence (which corresponds to the 2-categorical point of view).

\begin{definition}[Equivalent and isomorphic categories]\label{def:equivalent_isomorphic_categories}\uses{def:category, def:functor, def:natural_transformation}
\begin{enumerate}[(i)]
	\item A functor $F\colon \mcC\to\mcD$ is called \idx{isomorphism}, if there exists a functor $G\colon \mcD\to\mcC$ with $F\circ G=1_{\mcD}$ and $G\circ F=1_{\mcC}$.
	In this case, we call the two categories \idx{isomorphic}.
\item More importantly, two functors $F, G\in \mathrm{Fun}(\mcC,\mcD)$ are called \textbf{isomorphic} (short $F\simeq G$), if they are isomorphic in the functor category,
i.e. there exists $\eta\colon F\to G$ and $\nu\colon G\to F$ with $\nu\circ \eta=1_F$ and $\eta\circ \nu=1_G$.

	Two categories are called \idx{equivalent} if there exist functors $R\colon \mcC\to \mcD$ and $L\colon \mcD\to \mcC$ such that $R\circ L\simeq 1_{\mcD}$ and $L\circ R\simeq 1_{\mcC}$.
	In this case, one calls $L$ and $R$ \idx{essentially inverse} to each other.
\end{enumerate}
\end{definition}

\begin{example}
	\begin{enumerate}[(i)]
		\item The category $\{1\to 2\}$ is isomorphic to the category $\{2\to 3\}$, although not equal.

		\item The category
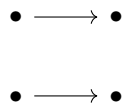
\begin{center}\begin{tikzcd}
	\bullet & \bullet \\
	\bullet & \bullet
	\arrow[from=1-1, to=1-2]
	\arrow[from=2-1, to=2-2]
\end{tikzcd}\end{center}
		is not equivalent to the category $\{\bullet\to \bullet\}$, but the category
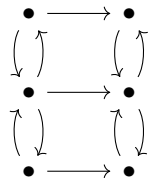
\begin{center}\begin{tikzcd}
	\bullet & \bullet \\
	\bullet & \bullet \\
	\bullet & \bullet
	\arrow[from=1-1, to=1-2]
	\arrow[from=2-1, to=2-2]
	\arrow[curve={height=6pt}, from=1-1, to=2-1]
	\arrow[curve={height=6pt}, from=2-1, to=1-1]
	\arrow[curve={height=6pt}, from=1-2, to=2-2]
	\arrow[curve={height=6pt}, from=2-2, to=1-2]
	\arrow[curve={height=-6pt}, from=2-1, to=3-1]
	\arrow[curve={height=-6pt}, from=3-1, to=2-1]
	\arrow[from=3-1, to=3-2]
	\arrow[curve={height=-6pt}, from=3-2, to=2-2]
	\arrow[curve={height=-6pt}, from=2-2, to=3-2]
\end{tikzcd}\end{center}
is equivalent to  $\{\bullet\to \bullet\}$, although not isomorphic.

		\item The category of finite dimensional $\C$-vector spaces is equivalent to the category $\{\C^n\,:\,n\in \N_0\}$ with linear maps.
		\item The category of Hilbert spaces is equivalent to the category of $\ell^2(I)$-spaces for varying index sets $I$ (assuming a strong enough axiom of choice).
		\item Skip if you're an algebraist: The category of compact Hausdorff spaces is \textbf{dual}\index{dual category} to the category of commutative unital $C^*$-algebras with unital $*$-algebra homomorphisms, meaning that $\CHaus^{\mathrm{op}}$ is equivalent to the category $C^*$-alg.
		 This is the classical Gelfand-Naimark theorem.
		The category $\CHaus$ is also equivalent to the dual to the category of unital AM-spaces with lattice homomorphisms (this is Gelfand-Kakutani, see e.g. \cite{schaefer1974Banachlattices}), and to the category of $C(K)$-spaces with Koopman-operators (pullback-operators).
		Restricting this to hyperstonean spaces, one obtains analogous duality to unital commutative $W^*$-algebras, $L^\infty$-spaces or order complete unital AM-spaces.
		\item Stone-duality will be another important example of a duality of categories,
		see \ref{subsec:stone_duality} below.
	\end{enumerate}
\end{example}

The following characterisation yields a good insight into the nature of equivalences of categories, i.e. that they are some sort of isomorphisms after taking isomorphism classes of objects in the categories.

\begin{proposition}[Characterisation of equivalences]\label{prop:characterisation_of_equivalences}\uses{def:equivalent_isomorphic_categories, def:special_morphisms, def:small_cat}
	A functor $F\colon \mcC\to \mcD$ where $\mcD$ is a small category is an equivalence if and only if it is fully faithful and essentially surjective.
\end{proposition}

\begin{proof}
See \cite[IV.4.1]{mac2013categories} or \cite[1.5.9]{Riehl}.
\end{proof}

\begin{remark}
	The construction of the essential inverse heavily relies on the axiom of choice.
	If one assumes a stronger version of it,
	e.g. by assuming the axiom of global choice,
	one can obtain versions of this result that do not need smallness of the categories.
\end{remark}

The following result is very useful for replacing an a priori very big category with a suitable smaller one.

\begin{corollary}[Equivalence of a category to a skeleton]\label{cor:equivalence_to_skeleton}\uses{prop:characterisation_of_equivalences}
	Any small category $\mcC$ is equivalent to a \idx{skeleton} of itself, i.e. a full subcategory $\mcC'\subset \mcC$ such that any element in $\mcC$ is isomorphic to precisely one object in $\mcC'$.
\end{corollary}

Often smallness conditions can be dropped, as soon as the skeleton is explicit enough.

\begin{warning}
	We will often switch to equivalent categories without mentioning and leave these steps implicit; we hope that this simplifies more than it creates confusion.
	In particular by switching to a skeleton we sometimes omit isomorphisms and directly write equalities instead of isomorphy.

	A better version, however, would be to always mean isomorphism,
	whenever we say equality.
	This is an instance of \emph{univalence}, see \cite{Cisinski2025}.
\end{warning}


\section[Limits]{Limits and Colimits}
\quot{A systematic treatment of all possible properties of limits was contained in a manuscript by Chevalley on category theory; the manuscript was unfortunately lost by some shipping company}{Mac Lane, p232 in \cite{mac2013categories}}
\subsection{Defining (co)limits}

One of the most important concepts in category theory are limits and colimits,
describing some sort of \emph{universal} objects with respect to some diagram of other objects.

\begin{warning}[skip if you're an algebraist]
	Note that, although there is the classical notion of limit in topology,
	we strongly suggest to see the analogy of topological with categorical limits
	in purely the names.
	See e.g. \url{https://math.stackexchange.com/q/62800} for a setup
    that yields topological limits as categorical ones.	
\end{warning}

\begin{definition}[Diagram]\label{def:diagram}\uses{def:category, def:functor}
A \idx{diagram} in a category $\mcC$ is a functor $D\colon \mcI\to \mcC$
from a small category $\mcI$ to $\mcC$.
Morphisms of diagrams are given by natural transformations between the functors.

Moreover, for any cardinal number $\kappa$, a diagram $\mcI\to \mcC$ is said to be \idx{$\kappa$-small}, if the cardinality of the set of all morphisms in $\mcI$ is less than $\kappa$.
\end{definition}

At first sight, this definition might seem too general for the concrete concept of a diagram.
One should think of the category $\mcI$ as the \idx{form} or the \idx{index category} of the diagram,
noting that one can always represent small categories as (directed) graphs,
and the functor is simply filling in the form by putting a concrete object of $\mcC$ at every node and a concrete morphism to every arrow.

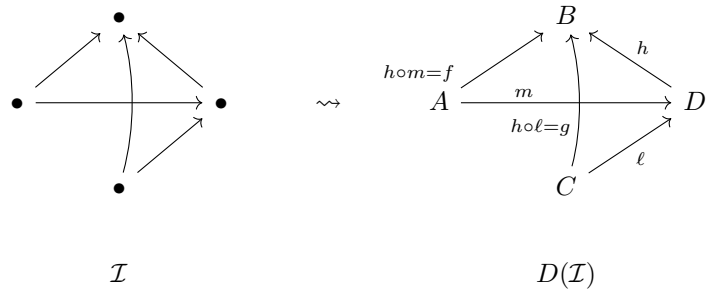
\begin{center}
\begin{tikzcd}
	& \bullet &&&& B \\
	\bullet && \bullet & \leadsto & A && D \\
	& \bullet &&&& C \\
	& \mcI &&&& {D(\mcI)}
	\arrow[from=2-1, to=1-2]
	\arrow[from=2-1, to=2-3]
	\arrow[from=2-3, to=1-2]
	\arrow["{{h\circ m=f}}"{pos=0}, from=2-5, to=1-6]
	\arrow["m"{pos=0.3}, from=2-5, to=2-7]
	\arrow["h"', from=2-7, to=1-6]
	\arrow[curve={height=6pt}, from=3-2, to=1-2]
	\arrow[from=3-2, to=2-3]
	\arrow["{{h\circ\ell=g}}"{pos=0.2}, curve={height=6pt}, from=3-6, to=1-6]
	\arrow["\ell"', from=3-6, to=2-7]
\end{tikzcd}
\end{center}

A particularly easy diagram one can build is given by always filling in the same object of the category, and just identity morphisms on the morphisms.
This is done by using a constant functor $\mcI\to\mcC$, and called a \idx{constant diagram}.
Often one collapses constant diagrams to just a single node (as the commutativity of any diagram does not change by this operation).

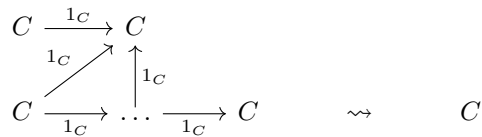
\begin{center}
\begin{tikzcd}
	C & C \\
	C & \dots & C & \leadsto & C
	\arrow["{1_C}", from=1-1, to=1-2]
	\arrow["{1_C}", from=2-1, to=1-2]
	\arrow["{1_C}"', from=2-1, to=2-2]
	\arrow["{1_C}"', from=2-2, to=1-2]
	\arrow["{1_C}"', from=2-2, to=2-3]
\end{tikzcd}
\end{center}
Usually, we denote a constant diagram $\mcI\to\mcC$ just by the object $C$.

\begin{definition}[Cone]\label{def:cone}\uses{def:category, def:diagram}
Let $D\colon \mcI\to \mcC$ be a diagram of nonempty form $\mcI$ in a category $\mcC$.
A \idx{cone} over $D$ is a morphism (natural transformation)
from a constant diagram $C\in\mcC$ to $D$.
Explicitly, after collapsing the constant diagram as described before,
this consists of the object $C$ and a family of morphisms $(\pi_i\colon C\to D(i))_{i\in \mcI}$
such that for all morphisms $f\colon i\to j$ in $\mcI$ the diagram
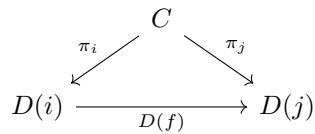
\begin{center}
\begin{tikzcd}
	& C \\
	{D(i)} && {D(j)}
	\arrow["{{\pi_i}}"', from=1-2, to=2-1]
	\arrow["{{\pi_j}}", from=1-2, to=2-3]
	\arrow["{{D(f)}}"', from=2-1, to=2-3]
\end{tikzcd}
\end{center}
is commutative.
Dually, a \idx{cocone} over $D$ is a morphism from $D$ to a constant diagram $C\in\mcC$,
i.e., the object $C$ and a family of morphisms $(\tau_i\colon D(i)\to C)_{i\in \mcI}$
such that for all morphisms $f\colon i\to j$ in $\mcI$ the diagram
\begin{center}
\begin{tikzcd}
	& C \\
	{D(i)} && {D(j)}
	\arrow["{{\tau_i}}", from=2-1, to=1-2]
	\arrow["{{D(f)}}"', from=2-1, to=2-3]
	\arrow["{{\tau_j}}"', from=2-3, to=1-2]
\end{tikzcd}
\end{center}
is commutative.
\end{definition}

Having defined cones and diagrams, we can define the notions of limits and colimits as \emph{universal} cones, resp. cocones.

\begin{definition}[Limit and colimit]\label{def:limit}\uses{def:category, def:diagram, def:cone}
Let $D\colon \mcI\to \mcC$ be a diagram in a category $\mcC$.
A \idx{limit} of $D$ is a cone $(L, (\pi_i)_{i\in \mcI})$ over $D$
such that for every cone $(C, (p_i)_{i\in \mcI})$ over $D$
there exists a unique morphism $u\colon C\to L$ in $\mcC$ such that for all $i\in \mcI$ the diagram
\begin{center}
\begin{tikzcd}
	C && L \\
	& {D(i)}
	\arrow["u", dashed, from=1-1, to=1-3]
	\arrow["{{p_i}}"', from=1-1, to=2-2]
	\arrow["{{\pi_i}}", from=1-3, to=2-2]
\end{tikzcd}
\end{center}
is commutative.
We denote the limit by
\[
    L=\varprojlim D=\varprojlim\limits_{i\in \mcI}D(i),
\]
and call the components $\pi_i\colon  L\to D(i)$ the \idx{projection onto the coordinate} $i$
(although this morphism is not necessarily epic or a projection).

A \idx{colimit} is a cocone $(L,(\tau_i)_{i\in \mcI})$
such that for all cocones $(C,(e_i)_{i\in \mcI})$ there exists a unique
$u\colon L\to C$ such that for all $i\in \mcI$

\begin{center}
\begin{tikzcd}
	C && L \\
	& {D(i)}
	\arrow["u"', dashed, from=1-3, to=1-1]
	\arrow["{{e_i}}", from=2-2, to=1-1]
	\arrow["{{\tau_i}}"', from=2-2, to=1-3]
\end{tikzcd}
\end{center}
is commutative.
We denote colimits by
\[
    L=\varinjlim D=\varinjlim\limits_{i\in \mcI}D(i)
\]
and call the components $\tau_i\colon D(i)\to L$ \idx{embedding into the coordinate} $i$
(although this morphism is not necessarily monic or an embedding).

A category $\mcC$ is called \idx{complete} if all diagrams in $\mcC$ have a limit, and \idx{cocomplete} if all diagrams in $\mcC$ have a colimit.
If it is both, we call the category \idx{bicomplete}.
\end{definition}

The following diagrams summarize the definition of limits and cones
\begin{center}\begin{tikzcd}
	{D:} & {D(i)} & {D(j)} & \dots & \dots \\
	\\
	&& {\varprojlim D} && C
	\arrow[from=1-2, to=1-3]
	\arrow[from=1-3, to=1-4]
	\arrow["{\pi_i}", from=3-3, to=1-2]
	\arrow["{\pi_j}"'{pos=0.4}, from=3-3, to=1-3]
	\arrow[from=3-3, to=1-4]
	\arrow[from=3-3, to=1-5]
	\arrow[from=3-5, to=1-2]
	\arrow[from=3-5, to=1-3]
	\arrow[from=3-5, to=1-4]
	\arrow[from=3-5, to=1-5]
	\arrow["{\exists !}", dashed, from=3-5, to=3-3]
\end{tikzcd}\end{center}

as well as colimits and cocones
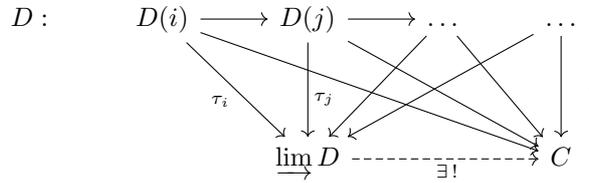
\begin{center}\begin{tikzcd}
	{D:} & {D(i)} & {D(j)} & \dots & \dots \\
	\\
	&& {\varinjlim D} && C
	\arrow[from=1-2, to=1-3]
	\arrow["{\tau_i}"', from=1-2, to=3-3]
	\arrow[from=1-2, to=3-5]
	\arrow[from=1-3, to=1-4]
	\arrow["{\tau_j}"{pos=0.6}, from=1-3, to=3-3]
	\arrow[from=1-3, to=3-5]
	\arrow[from=1-4, to=3-3]
	\arrow[from=1-4, to=3-5]
	\arrow[from=1-5, to=3-3]
	\arrow[from=1-5, to=3-5]
	\arrow["{\exists !}"', dashed, from=3-3, to=3-5]
\end{tikzcd}.\end{center}

\begin{remark}
A limit is sometimes also referred to as \idx{inverse limit} or \idx{projective limit},
and a colimit is also called \idx{direct limit} or \idx{inductive limit}.

One could also rephrase the definition of colimits in terms of limits in the opposite category $\mcC^{\mathrm{op}}$.
\end{remark}

\begin{definition}[Initial and terminal objects]\label{def:init_term}\uses{def:limit}
For the empty diagram we define cones and cocones to be given by the above explicit description rather than by functors, i.e. they are given by single objects.
The limit over the empty diagram is called \idx{terminal} object of the category.
This is an object $\ast$, such that for every other object $A\in\mcC$ there exists a unique arrow $A\to \ast$.

Dually, the colimit over the empty diagram is called \idx{initial object},
i.e. an object $\emptyset$ with a unique arrow $\emptyset\to A$ to every $A\in\mcC$.
\end{definition}

We remark some very important general facts about limits and colimits,
before diving into concrete examples.

\begin{remark}[Large (Co)limits]\label{rem:large_limits}\uses{def:limit, def:limit}
If the category $I$ is not small, one can still define diagrams and limits over $I$, which we will call \idx{large (co)limits}.
However, most of the relevant limits and colimits will be small, and in particular it is not reasonable to expect all large limits to exist.
    This is due to the fact that such a category automatically would be \idx{thin}, i.e. hom-sets have at most one element,
which follows by contradiction by assuming that there are two distinct morphisms $X\to Y$ and concluding that the cardinality
of $\hom(X,\prod_I Y)$ is unbounded for a proper class $I$.
    The same reasoning shows that a small category that is complete is thin, by remarking that otherwise the cardinality of $\hom(X,\prod_{\hom(\mcC)}Y)$ would be larger than the set of all morphisms $\hom(\mcC)$ of $\mcC$.
\end{remark}

\begin{remark}[(Co)limits are unique up to isomorphism]\label{rem:unique_up_to_isomorphism}\uses{def:limit, def:limit}
Limits are, like every other universal constructions,
uniquely determined up to unique isomorphisms.
For any two limits $(L_1,(\pi_i^1)) \, (L_2,(\pi_i^2))$ of a diagram $D$
there exists a unique morphism $f\colon L_1\to L_2$
such that for all $i$ the  diagram
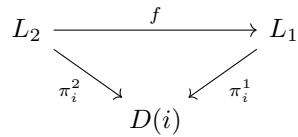
\begin{center}
\begin{tikzcd}
	{L_2} && {L_1} \\
	& {D(i)}
	\arrow["f", from=1-1, to=1-3]
	\arrow["{{\pi_i^2}}"', from=1-1, to=2-2]
	\arrow["{{\pi_i^1}}", from=1-3, to=2-2]
\end{tikzcd}
\end{center}
is commutative, and this morphism is an isomorphism.
The (easy) argument for this can be found in every standard book
and essentially reduces to using the existence part in the universal property
of both limits for the existence of two morphisms $f\colon L_1\to L_2$ and $g\colon L_2\to L_1$,
and the uniqueness to see that $fg=1$ and $gf=1$.
The same discussion applies to colimits.
\end{remark}

This justifies talking about \emph{the} limit of a diagram.
Often, one can compute limits and colimits of diagrams by just considering certain subdiagrams.

\begin{definition}[final and initial functor]\label{def:final}\uses{def:category}
A functor $E\colon \mcI\to\mcJ$ is called \idx{final} (the old term for this is \idx{cofinal}), if for any $j\in\mcJ$ the comma category $(j\downarrow E)$ is connected, which in this case means that
\begin{enumerate}[(i)]
	\item it is nonempty, i.e. there exists an $i\in \mcI$ with a morphism $j\to E(i)$,
	    and
	\item every two morphisms $f\colon j\to E(i)$ and $g\colon j\to E(i')$ can be connected by a \emph{zig-zag},
	meaning that there exist finitely many $i_0 = i,i_1,i_2,\dots, i_n,i_{n+1}=i'$ with $f_{i_k}\colon j\to E(i_k)$ such that for any consecutive pair $i_k,i_{k+1}$ there exists
	 an arrow $h_k\colon i_{k}\to i_{k+1}$
	 or an arrow $h_k\colon i_{k+1}\to i_k$ such that
\begin{center}
\begin{tikzcd}
	& j \\
	{E(i_k)} && {E(i_{k+1})}
	\arrow["{{f_k}}"', from=1-2, to=2-1]
	\arrow["{{f_{k+1}}}", from=1-2, to=2-3]
	\arrow["{{E(h_k)}}"', from=2-1, to=2-3]
\end{tikzcd},
    resp.
\begin{tikzcd}
	& j \\
	{E(i_k)} && {E(i_{k+1})}
	\arrow["{{f_k}}"', from=1-2, to=2-1]
	\arrow["{{f_{k+1}}}", from=1-2, to=2-3]
	\arrow["{{E(h_k)}}", from=2-3, to=2-1]
\end{tikzcd}
\end{center}
    is commutative.
\end{enumerate}

A functor $E\colon \mcI\to\mcJ$ is called \idx{initial},
if for any $j\in\mcJ$ the comma category $(E\downarrow j)$ is connected, meaning that
\begin{enumerate}[(i)]
	\item there exists a $i\in \mcI$ with a morphism $ E(i)\to j$, and
	\item every two morphisms $f\colon  E(i)\to j$ and $g\colon  E(i')\to j$ can be connected by a zig-zag, meaning that there exist finitely many $i_0=i,i_1,i_2,\dots, i_n,i_{n+1}=i'$ with $f_i\colon j\to E(i_n)$ such that
	for any consecutive pair $i_k,i_{k+1}$ there exists
    an arrow $h_k\colon i_{k}\to i_{k+1}$
    or an arrow  $h_k\colon i_{k+1}\to i_k$ such that
\begin{center}
\begin{tikzcd}
	{E(i_k)} && {E(i_{k+1})} \\
	& j
	\arrow["{{{E(h_k)}}}", from=1-1, to=1-3]
	\arrow["{{{f_k}}}"', from=1-1, to=2-2]
	\arrow["{{{f_{k+1}}}}", from=1-3, to=2-2]
\end{tikzcd},
    resp.
\begin{tikzcd}
	{E(i_k)} && {E(i_{k+1})} \\
	& j
	\arrow["{{{f_k}}}"', from=1-1, to=2-2]
	\arrow["{{{E(h_k)}}}"', from=1-3, to=1-1]
	\arrow["{{{f_{k+1}}}}", from=1-3, to=2-2]
\end{tikzcd}
\end{center}
    is commutative.
\end{enumerate}
A \idx{subdiagram} (meaning the restriction of some Diagram $\mcJ\to \mcC$ onto some subcategory $\mcI\subset\mcJ$) is called \textbf{final}\index{final subdiagram}/\textbf{initial}\index{initial subdiagram},
if the inclusion functor is final/initial.
\end{definition}

\begin{example}[Final subdiagrams of linearly ordered sets]\label{ex:final_subsets_lin_ordered}
Consider any \idx{totally ordered set} $\mcJ$ (also called linearly ordered),
regarded as a category with morphism $a\to b$ iff $a\le b$ as in example~\ref{ex:examples_of_categories}.
Since for every $j$ and any subcategory $\mcI$ condition (ii) for $(j\downarrow\mcI)$
is fulfilled,
the final full subcategories agree with the usual cofinal subsets,
being subsets $\mcI\subseteq \mcJ$ such that for any $j\in \mcJ$
there exists $i\in \mcI\subseteq\mcJ$ with $j\le i$.
\end{example}

\begin{lemma}[Colimits along final functors]\label{lem:limits_along_final_subdiagrams}\uses{def:limit, def:final}
A functor $E\colon \mcI\to\mcJ$ is final precisely if for every category $\mcC$
and any diagram $D\colon \mcJ\to \mcC$ such that $DE$ has a colimit,
the diagram $D$ has a colimit and the natural map
\[
    \varinjlim_i DE(i)\to \varinjlim_j D(j)
\]
is an isomorphism.

Dually, the functor is initial if and only if for every category $\mcC$
and any diagram $D\colon \mcJ\to \mcC$ such that $DE$ admits a limit,
the diagram $D$ has a limit and the natural map
\[
    \varprojlim_j D(j)\to \varprojlim_i DE(i)
\]
is an isomorphism.
\end{lemma}

\begin{proof}
See \cite[Theorem IX.3.1]{mac2013categories} (and Exercise 5 for the converse),
or \cite[Proposition 2.5.2]{kashiwara2005categories}.
See also Lemma 4.17.2 in \cite[04E7]{SPA2023} or final functor in nLab.
\footnote{\url{https://ncatlab.org/nlab/show/final+functor}}
\end{proof}

\begin{example}\label{ex:colimit-evaluation-at-point}
Let $j\in\mcJ$. Then the discrete subdiagram $\{j\}\subset\mcJ$ is final
if and only if $j$ is a terminal object in $\mcJ$.
In view of the previous lemma~\ref{lem:limits_along_final_subdiagrams},
this means that one can compute colimits over a diagram with terminal object
by just evaluating at a terminal object.

The dual statement holds for limits,
those are given by evaluating the diagram at the initial object, if it exists.
\end{example}

If one restricts the target category instead of the index category, the following basic fact often helps.

\begin{lemma}[Limits in full subcategories]\label{lem:lim_in_full_sub}\uses{def:limit, def:full_faithful_essentially_surjective_functors}
Consider any diagram $D\colon \mcI\to \mcD$ in any full subcategory $\mcD\sub\mcC$.
If the (co)limit, computed in $\mcC$, exists and again is contained in $\mcD$ (more precisely, in the essential image of the inclusion), then it automatically is a (co)limit in $\mcD$.
\end{lemma}

\begin{proof}
This is immediately verified.
\end{proof}

\subsection{Examples of (co)limits}
\quot{What do you call someone who reads a paper on category theory?\\
Answer: A coauthor}{The Rising Sea \cite{Vakil2023}, 1.4.11}

Now we consider different important examples of limits and colimits.

\begin{definition}[Special limits and colimits]\label{def:special_limits_and_colimits}\uses{def:limit, def:limit}
In this definition we give a summary of the most important limits and colimits we encounter, the ones in brackets will not be important to us, we include them for the sake of completeness.
\begin{table}[h]
\begin{center}
\begin{adjustbox}{max width=\textwidth}
\begin{tabular}{|ll||ll|}
\hline
Concept & & Dual concept& \\ \hline
Form $\mcI$ & Limit & Form $\mcI$ & Colimit \\ \hline\hline
 Discrete & \idx{Product} $\prod$ & Discrete & \idx{Coproduct} \\&&& or \idx{disjoint union}  $\coprod$ \\ \hline
\begin{tikzcd}
	\bullet & \bullet
	\arrow[shift left, from=1-1, to=1-2]
	\arrow[shift right, from=1-1, to=1-2]
\end{tikzcd} & \idx{Equalizer} & \begin{tikzcd}
	\bullet & \bullet
	\arrow[shift left, from=1-1, to=1-2]
	\arrow[shift right, from=1-1, to=1-2]
\end{tikzcd}  & \idx{Coequalizer} \\ \hline

( \begin{tikzcd}
	\bullet & \bullet
	\arrow["{\text{(common retraction)}}"', shift right=2, from=1-1, to=1-2]
	\arrow[shift left=2, from=1-1, to=1-2]
	\arrow[from=1-2, to=1-1]
\end{tikzcd} )& (\idx{coreflexive equalizer})
 &
\begin{tikzcd}
	\bullet & \bullet
	\arrow["{\text{(common section)}}"', shift right=2, from=1-1, to=1-2]
	\arrow[shift left=2, from=1-1, to=1-2]
	\arrow[from=1-2, to=1-1]
\end{tikzcd}&\idx{Reflexive coequalizer}\\ \hline
\begin{tikzcd}
	& \bullet \\
	\bullet & \bullet
	\arrow[from=1-2, to=2-2]
	\arrow[from=2-1, to=2-2]
\end{tikzcd}& \idx{Pullback} or \idx{fibre product} &
\begin{tikzcd}
	\bullet & \bullet \\
	\bullet
	\arrow[from=1-1, to=1-2]
	\arrow[from=1-1, to=2-1]
\end{tikzcd} & \idx{Pushout} \\ \hline
\idx{Cofiltered diagram}& \idx{Cofiltered limit} & \idx{Filtered diagram}	&  \idx{Filtered colimit}\\ \hline
(\idx{Cosifted diagram})&(\idx{Cosifted limit})&\idx{Sifted diagram}& \idx{Sifted colimit}\\ \hline
\end{tabular}
\end{adjustbox}
\end{center}
\end{table}
\end{definition}

In more detail, let $\mcC$ be a category, and $D$ a diagram of the form $\mcI$.
\begin{enumerate}[(i)]
\item If $\mcI$ is discrete, we identify $\mcI$ with its underlying set of objects.
	Then the diagram $D\colon \mcI\to \mcC$ is simply an $\mcI$-indexed family of objects $(D_i)_{i\in \mcI}$ in $\mcC$.
A cone now simply consists of an object $C\in \mcC$ together with an arbitrary family of morphisms $(f_i\colon C\to D_i)$ with no further required compatibilities.
A \idx{product} of $(D_i)_{i\in\mcI}$ now is an object $\prod_{i\in \mcI} D_i$ together with \emph{coordinate projections}
\[
    \pi_j\colon \prod_{i\in \mcI} D_i\to D_j
\]
for all $j\in \mcI$
(which in general do not have to be projections or epimorphisms)
such that for all objects $C\in \mcC$ and all families
$(f_i\colon C\to D_i)_{i\in \mcI}$ of morphisms there exists a unique \idx{gluing}
$\left(\prod f_i\right)\colon C\to \prod D_i$, making the diagram
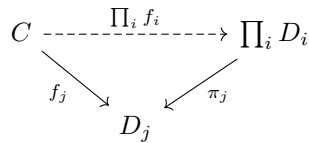
\begin{center}
\begin{tikzcd}
	C && {\prod_i D_i} \\
	& {D_j}
	\arrow["{{\prod_i f_i}}", dashed, from=1-1, to=1-3]
	\arrow["{{f_j}}"', from=1-1, to=2-2]
	\arrow["{{\pi_j}}", from=1-3, to=2-2]
\end{tikzcd}
\end{center}
commutative for all $j\in \mcI$.
The gluing intuitively can be thought of as writing $f_i$ in the $i$-th coordinate.

A \idx{coproduct} or \idx{disjoint union} consists of an object $\coprod_{i\in\mcI} D_i$ together with a family of \emph{embeddings}
\[
    e_i\colon D_i\to \coprod_{i\in\mcI} D_i
\]
(again, they do not have to be monics),
such that for any family $(f_i\colon D_i\to C)$ there exists a unique $\sqcup_i f_i \colon \coprod D_i\to C$ such that for all $j\in\mcI$ the diagram
	
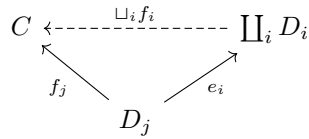
\begin{center}
\begin{tikzcd}
	C && {\coprod_i D_i} \\
	& {D_j}
	\arrow["{{\sqcup_i f_i}}"', dashed, from=1-3, to=1-1]
	\arrow["{{f_j}}", from=2-2, to=1-1]
	\arrow["{e_i}"', from=2-2, to=1-3]
\end{tikzcd}
\end{center}
is commutative.
The morphism can be thought of as the morphism that is given on each summand
$D_i$ by $f_i$.

\item Let us assume that the diagram consists of two parallel arrows $f,g\colon A\to B$.

A cone over such a diagram a priori is an object $C$ together with arrows $h_A\colon C\to A$ and $h_B\colon C\to B$ such that
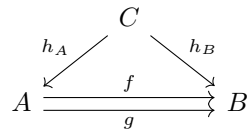
\begin{center}\begin{tikzcd}
	& C \\
	A && B
	\arrow["{h_A}"', from=1-2, to=2-1]
	\arrow["{h_B}", from=1-2, to=2-3]
	\arrow["f", shift left, from=2-1, to=2-3]
	\arrow["g"', shift right, from=2-1, to=2-3]
\end{tikzcd}\end{center}
commutes.
But the commutativity reduces to
	 \[h_B=fh_A=gh_A\]
	and hence $h_B$ is uniquely determined by $h_A$ and the property that $h_A$ \idx{equalizes} $f$ and $g$, meaning $fh_A=gh_A$.
	Thus one abbreviates the cone diagram to
\begin{center}\begin{tikzcd}
	C & A & B
	\arrow["{h_A}", from=1-1, to=1-2]
	\arrow["f", shift left, from=1-2, to=1-3]
	\arrow["g"', shift right, from=1-2, to=1-3]
\end{tikzcd}.\end{center}
	 Having understood cones over these diagrams, we call their limit the \idx{equalizer} of the diagram and often denote it by $\mathrm{eq}(f,g)$.
	Abbreviating cones as explained above, the equalizer is a morphism $e\colon E\to A$ with $fe=ge$ such that all other $h\colon F\to A$ with $fh=gh$ uniquely factors through $e$,
	meaning that there exists exactly one $j\colon F\to E$ with $h=ej$.
\begin{center}
\begin{tikzcd}
	C \\
	E & A & B
	\arrow["j"', dashed, from=1-1, to=2-1]
	\arrow["h", from=1-1, to=2-2]
	\arrow["e"', from=2-1, to=2-2]
	\arrow["f", shift left, from=2-2, to=2-3]
	\arrow["g"', shift right, from=2-2, to=2-3]
\end{tikzcd}
\end{center}

	Intuitively, the equalizer describes the largest subobject of $A$ such that on this subobject the two morphisms $f$ and $g$ agree.
	In particular, one can show that due to the uniqueness part of the universal property, an equalizer is always monic (note that this is an example where the \emph{coordinate projection} is often not epic).
	It is an isomorphism precisely if the two arrows coincide.

Dually, cocones along such diagrams are given by morphisms $h\colon B\to C$
\begin{center}\begin{tikzcd}
	A & B & C
	\arrow["f", shift left, from=1-1, to=1-2]
	\arrow["g"', shift right, from=1-1, to=1-2]
	\arrow["h", from=1-2, to=1-3]
\end{tikzcd}\end{center}
	such that $h$ \idx{coequalizes} $f$ and $g$, i.e. $hg=hf$.
	The colimit of this diagram is called \idx{coequalizer},
	and is given by a morphism $c\colon B\to C$ with $cf=cg$ such that every other $h\colon B\to D$ with $hf=hg$ uniquely \idx{factors through} $c$, meaning that there exists a unique $j\colon C\to D$ with $jc=h$.

\begin{center}\begin{tikzcd}
	&& D \\
	A & B & C
	\arrow["g"', shift right, from=2-1, to=2-2]
	\arrow["f", shift left, from=2-1, to=2-2]
	\arrow["h", from=2-2, to=1-3]
	\arrow["c"', from=2-2, to=2-3]
	\arrow["j"', dashed, from=2-3, to=1-3]
\end{tikzcd}\end{center}

	This approximately corresponds to taking the quotient by the smallest equivalence relation that pointwise identifies the images of $f$ and $g$.
	Again, due to the uniqueness part of the universal property, any coequalizer is always an epimorphism.
	It is an isomorphism precisely if the two arrows coincide.

\item A \idx{reflexive coequalizer} is a slightly stronger notion of a coequalizer,
	where we add the property that both maps $f,g\colon A\to B$ have a common section $s\colon B\to A$ (meaning $fs=1_B$ and $gs=1_B$). Whether one takes the colimit of the diagram with or without the section is irrelevant, only the existence of such a section is important.
	
	 The intuition of reflexive coequalizers is closer to being a generalisation of taking factors after an equivalence relation.
One may also call this coequalizer a \idx{quotient}.

	In fact, for any classical equivalence relation on sets, $E\sub B\times B$ one can consider the two coordinate projections $p_1,p_2\colon E\to B$,
	and the common section is given by the diagonal map, which is precisely the reflexivity property of the equivalence.

\item Now we investigate diagrams with three objects and two morphisms.

Consider two morphisms with the same codomain
\begin{center}\begin{tikzcd}
	& B \\
	C & D
	\arrow["g", from=1-2, to=2-2]
	\arrow["f"', from=2-1, to=2-2]
\end{tikzcd}.\end{center}
	A cone to this diagram consists of an object $A$ together with morphisms $A\to B$, $A\to C$ and $A\to D$, as depicted in the diagram
	\begin{center}
\begin{tikzcd}
	A & B \\
	C & D
	\arrow["{h_B}", from=1-1, to=1-2]
	\arrow["{h_C}"', from=1-1, to=2-1]
	\arrow["{h_D}", from=1-1, to=2-2]
	\arrow["f", from=1-2, to=2-2]
	\arrow["g"', from=2-1, to=2-2]
\end{tikzcd}.
		\end{center}
Again, the morphism $A\to D$ can be recovered via $A\to C\to D$ or $A\to B\to D$, and hence contains no further information and may be left implicit.
	Thus a cone can be described via the morphisms $h_B\colon A\to B$ and $h_C\colon A\to C$ such that $gh_C=fh_B$, i.e. the commutative diagram
\begin{center}\begin{tikzcd}
	A & B \\
	C & D
	\arrow["{h_B}", from=1-1, to=1-2]
	\arrow["{h_C}"', from=1-1, to=2-1]
	\arrow["f", from=1-2, to=2-2]
	\arrow["g"', from=2-1, to=2-2]
\end{tikzcd}.\end{center}
	The limit of this diagram if called the \idx{pullback}, \idx{fibre product}, \idx{cartesian square} or \idx{intersection} of $f$ and $g$, and is denoted by
\begin{center}\begin{tikzcd}
	{C\times_DB} & B \\
	C & D
	\arrow["{p_g}", from=1-1, to=1-2]
	\arrow["{p_f}"', from=1-1, to=2-1]
	\arrow["f", from=1-2, to=2-2]
	\arrow["g"', from=2-1, to=2-2]
\end{tikzcd},\end{center}
	where $p_f$ is the pullback of $f$ along $g$ (and $p_g$ the pullback of $g$ along $f$), which is sometimes also denoted by $f^*$ and $g^*$.
\footnote{Note that for any morphism $f\colon A\to B$ the precomposition $f^*\colon \hom(B,-)\to \hom(A,-)$ is also called the pullback of $f$. These notions can be related but we suggest to think of the term and notion pullback to have two different meanings, depending on the context.}
The universal property is given by
\begin{center}\begin{tikzcd}
	A \\
	& {C\times_DB} & B \\
	& C & D
	\arrow["{\exists ! \,k}"{pos=0.8}, dashed, from=1-1, to=2-2]
	\arrow["{h_B}", curve={height=-12pt}, from=1-1, to=2-3]
	\arrow["{h_C}"', curve={height=12pt}, from=1-1, to=3-2]
	\arrow["{p_g}", from=2-2, to=2-3]
	\arrow["{p_f}"', from=2-2, to=3-2]
	\arrow["f", from=2-3, to=3-3]
	\arrow["g"', from=3-2, to=3-3]
\end{tikzcd},\end{center}
	meaning that for every $h_B\colon A\to B, h_C\colon A\to C$ with $fh_B=gh_C$ there exists a unique $k\colon A\to B\times_D C$ with $p_f k=h_B$ and $p_g k=h_C$.
	A fibre product roughly collects all tuples of points $(b,c)\in B\times C$ for which $f(b)=g(c)$, yielding some sort of intersection relative to $D$.
	A special case of this is the classical \idx{intersection} of two subsets $B,C$ of $D$, yielding $B\times_D C=B\cap D$.
	Another example is taking \idx{preimages}, by considering a map $f\colon B\to D$
	and a subset $C\subset D$, leading to $B\times_{D} C=f^{-1}(C)$.

Dually, consider two morphisms with common domain,
\begin{center}
\begin{tikzcd}
	A & B \\
	C
	\arrow["f", from=1-1, to=1-2]
	\arrow["g"', from=1-1, to=2-1]
\end{tikzcd}.\end{center}
			  A cocone is given by an object $D$ and two morphisms $h_B\colon B\to D$, $h_C\colon C\to D$ such that the diagram

\begin{center}\begin{tikzcd}
	A & B \\
	C & D
	\arrow["f", from=1-1, to=1-2]
	\arrow["g"', from=1-1, to=2-1]
	\arrow["{h_B}", from=1-2, to=2-2]
	\arrow["{h_C}"', from=2-1, to=2-2]
\end{tikzcd}\end{center}
is commutative.
			   The colimit over this diagram is called the \idx{pushout} or \idx{cocartesian square} of $f$ and $g$, denoted by $C\sqcup_A B$, equipped with the universal property
\begin{center}\begin{tikzcd}
	A & B \\
	C & {C\sqcup_AB} \\
	&& D
	\arrow["f", from=1-1, to=1-2]
	\arrow["g"', from=1-1, to=2-1]
	\arrow[from=1-2, to=2-2]
	\arrow["{h_B}", curve={height=-12pt}, from=1-2, to=3-3]
	\arrow[from=2-1, to=2-2]
	\arrow["{h_C}"', curve={height=12pt}, from=2-1, to=3-3]
	\arrow["{\exists!}"'{pos=0.4}, dashed, from=2-2, to=3-3]
\end{tikzcd}\end{center}
			  A pushout morally is gluing $B$ and $C$ together along $A$,
  		   i.e. being the coproduct of $B$ and $C$ but with values of $f$ being identified with the corresponding values of $g$.
			  This generalizes the \idx{union} of subobjects, when taking $A$ as their intersection (i.e. the union is defined to be the pushout of the pullback of two subobjects). \footnote{In general, the union should rather be defined as the coproduct in the comma category, but these notions often coincide.}

\item A \idx{filtered category} $\mcJ$ is a nonempty category fulfilling one of the following equivalent conditions.
	\begin{enumerate}[(a)]
		\item 
		\begin{enumerate}[(i)]
		    \item For any two $i,j\in\mcJ$ there exists $k\in\mcJ$
		    with maps $i\to k$ and $j\to k$, and
		    \item for every two parallel arrows $f,g\colon i\to j$ for some $i,j$
		    there exists a $k\in\mcJ$ and a morphism $h\colon j\to k$
		    equalizing both, i.e. $hf=hg$.
        \end{enumerate}
		\item For evey \idx{finite category} $\mcF$
		(meaning $\mcF$ consists of finitely many objects and finitely many morphisms)
        the diagonal functor $\mcJ\to\mcJ^\mcF$ is final.
	\end{enumerate}
	A \idx{filtered colimit} is a colimit along a diagram with filtered index category (a \idx{filtered diagram}).
	The most important example of filtered categories are directed sets.
	In fact, there are not much more filtered categories,
	since every filtered index category admits a thin directed subcategory
	(which, up to size issues, can be identified with a directed set)
	such that the inclusion is final, i.e the colimits don't change,
	see \cite{Andreka1982} or lemma 4.21.5 in \cite[0032]{SPA2023}.
	Thus usually we assume implicitly that filtered colimits are taken over directed sets.
	Clearly, every category with finite colimits is filtered.

	For filtered categories the final full subcategories are particularly easy to describe.
	Similarly to example~\ref{ex:final_subsets_lin_ordered} one can easily check that
	a full subcategory $\mcI\sub\mcJ$ of a filtered category $\mcJ$ is final precisely if 
	it contains to every $j\in\mcJ$ an $i\in\mcI$ with arrow $j\to i$,
	because the second condition is automatic.

        Here is a nice property of filtered colimits whose proof is an excellent exercise.
        If the arrows in a filtered diagram are all epic,
        then the \enquote{embeddings} into the filtered colimit are also epic.

	More generally, for any cardinal $\kappa$,
	an index category $\mcJ$ is	\idx{$\kappa$-filtered},
	if for any $\kappa$-small subdiagram $D\subseteq \mcJ$,
	there exists a \emph{cocone} in $\mcJ$,
	meaning an object $j\in \mcJ$ with arrows $h_d\colon d\to j$ for all $d\in D$
	such that for all arrows $d\to d'$ in $\mcD$, the diagram
\begin{center}
\begin{tikzcd}
	& j \\
	d && {d'}
	\arrow["{{h_d}}", from=2-1, to=1-2]
	\arrow[from=2-1, to=2-3]
	\arrow["{{h_{d'}}}"', from=2-3, to=1-2]
\end{tikzcd}
\end{center}
	is commutative.
	This is equivalent to saying that the diagonal functor $\mcJ\to\mcJ^{D}$ is final for all such $D$.
	Clearly, the usual filtered categories are just the $\aleph_0$-filtered categories.
	
	\item A nonempty index category $\mcJ$ is called \idx{sifted}, if the diagonal functor $\mcJ\to\mcJ\times\mcJ$ is final.
	A colimit along a sifted category is called \idx{sifted colimit}.
	In particular, all reflexive coequalizers and all filtered colimits are sifted.
	Moreover, a category with finite coproducts is sifted.
	These constitute the most important examples,
	due to Lemma \ref{lem:comm_with_sifted}.
	For more on sifted colimits see \cite{Adamek2010}, \cite{Adamek2001} and \cite{Adamek2010a}, where their close relation to algebraic theories is highlighted.
\end{enumerate}

An example of quite large colimits are unions of categories.
\begin{definition}[Unions of categories]\label{def:cat-unions}
Consider any totally ordered class $I$, and a family of categories $\mcC_i$ indexed by $i\in I$.
Furthermore for any $i\le j$ in $I$ consider a fully faithful embedding $\tau_i^j\colon \mcC_i\to \mcC_j$ (for $i=j$ being the identity functor), commuting in the sense that $\tau_j^k\tau_i^j=\tau_i^k$ for $i\le j\le k$.
Hence we can view this collection as a large diagram in the (large) category of categories, $\Cat$.

\begin{center}\begin{tikzcd}
	\dots & {\mcC_i} & {\mcC_j} & {\mcC_k} & \dots
	\arrow[hook, from=1-1, to=1-2]
	\arrow[hook, from=1-2, to=1-3]
	\arrow[hook, from=1-3, to=1-4]
	\arrow[hook, from=1-4, to=1-5]
\end{tikzcd}\end{center}

The colimit $\varinjlim_i \mcC_i$ of this diagram, which we call the \idx{union of the categories $C_i$}, is given by the category whose
	objects are given by all objects of $\mcC_i$ for all $i\in I$,
	and morphisms between an object $A\in \mcC_i$ to $B\in \mcC_j$ are given by
	\[\hom(A,B)=\begin{cases}\hom_{\mcC_i}(A, \tau_j^i B)& j\le i\\\hom_{\mcC_j}(\tau_i^j A, B)&i<j\end{cases}\]

	One can easily check that this forms a colimit in $\Cat$, i.e. for any family of functors $F_i\colon \mcC_i\to \mcD$ that are compatible with the embeddings $\tau_i^j$,
	there exists a unique functor $F\colon \varinjlim \mcC_i \to \mcD$ such that $F$ on the restriction to $\mcC_i$ agrees with $F_i$
\end{definition}
	
\begin{warning}
It would be better to define this in an appropriate 2-categorical sense
where everything commutes up to natural isomorphy.
However, when working with embeddings we can restrict to diagrams that "really" commute.

Moreover, note the set-theoretic problem in this definition, which can be solved by using the Grothendieck universe approach, since then a class (a set in the middle universe) of classes (again in the middle universe) has a well defined class as union (using the regularity of the outer strongly inaccessible cardinal).
\end{warning}

\subsection{Computing (co)limits}

As of now, limits and colimits might seem like abstract constructions
which are hard to describe.
This is not the case, since in most categories, they have an explicit description.
As an example, we compute all limits and colimits in the category of sets, see e.g. 
Chapter 2.1 and Proposition 2.4.1 of \cite{kashiwara2005categories}
and Chapter II.3,4 and V.2 of \cite{mac2013categories}.

\begin{lemma}[$\Set$ is bicomplete]\label{lem:set_bicomplete}\uses{def:limit, prop:criterion_for_completeness}
For any diagram $D\colon \mcI\to \mathrm{Set}$,
the limit $\varprojlim D$ exists and is given by the subset of the set-theoretic product
$\prod_{i\in\mcI}D(i)$ consisting of all tuples $(d_i)_{i\in\mcI}$
such that for all morphisms $f\colon i\to j$ one has $d_j=D(f)(d_i)$. 
The corresponding universal cone is given by restriction of the coordinate projections. 
In particular the set-theoretic (cartesian) product equipped with coordinate projections agrees with the categorical product.

The colimit $\varinjlim D$ exist and is given by taking the quotient of the disjoint union
$\dot\bigcup_{i\in\mcI} D(i)$ by the smallest equivalence relation identifying
$d_i\sim d_j$ for $d_i\in D(i), \, d_j\in D(j)$ if $D(f)(d_i)= d_j$ for some $f\colon i\to j$.
The corresponding cocone morphisms are given by the set-inclusions followed by the factor map of the equivalence relation.

Again, the set-theoretic disjoint union equipped with coordinate embeddings agrees with the categorical coproduct, and will from now on mostly be denoted by $\bigsqcup$.

\end{lemma}

\begin{proof}
This is a nice exercise and the verification is elementary.
\end{proof}

Therefore, in $\Set$, any limit is a subset of a product,
and any colimit is a quotient of a coproduct.
There is a more general principle of describing (co)limits
via (co)products and reflexive (co)equalizers.
The following lemma describes the general situation
and shows how one can compute (co)limits out of simpler ones.

\begin{lemma}[Limits as combination of special limits]\label{lem:lim_via_special_lim}\uses{def:special_limits_and_colimits}
\begin{enumerate}[(i)]
	\item Any (finite) colimit can be written as the reflexive coequalizer of two (finite) coproducts,
	assuming that all relevant colimits exist.
	\item Every colimit can be written as the filtered colimit of the colimits over all finite subdiagrams,
	assuming the relevant colimits exist.
    \item Every colimit can be written as the sifted colimit of the coproducts over all finite discrete subdiagrams,
    assuming the relevant colimits exist.
\end{enumerate}
The dual statements also hold, in particular any (finite) limit can be written as the coreflexive equalizer of two (finite) products, assuming that all relevant limits exist.
\end{lemma}

\begin{proof}
	For (i) and its dual statement see e.g. chapter V.2 in \cite{mac2013categories}, noting that there is an obvious section.
	The second assertion is proved (for coproducts) in Theorem 1 in Chapter IX.1 in \cite{mac2013categories}.
	The proof directly generalizes to arbitrary colimits. 	The last claim is analogous to the third.
\end{proof}

The lemma implies the following.

\begin{corollary}[Criterion for (co)completeness]\label{prop:criterion_for_completeness}\uses{def:limit, lem:lim_via_special_lim}
A category is (finitely) complete precisely if it admits all coreflexive equalizers
of parallel pairs of morphisms and all (finite) products.
Dually, a category is (finitely) cocomplete precisely if it admits all reflexive
 coequalizers of parallel pairs of morphisms and all (finite) coproducts.
\end{corollary}

Another good class of examples of (co)complete categories can be found by looking at functor categories.

\begin{lemma}\label{lem:limits_in_functor_categories}\uses{def:functor_category, def:limit}
Assuming that all relevant (co)limits in the codomain category exist,
(co)limits in functor categories are computed pointwise.
That is, the limit of any diagram $F\colon \mcI\to\mathrm{Fun}(\mcC,\mcD)$ exists
and it is given by the functor $\varprojlim_{i\in \mcI} F_i$
defined on objects $A\in\mcC$ by
\[
    (\varprojlim_{i\in \mcI} F_i)(A) = \varprojlim_{i\in\mcI} (F_i A)
\]
and on morphisms $f\colon A\to B$ by the induced universal arrow given in the diagram
\begin{center}
\begin{tikzcd}
	{(\varprojlim F)(A)} &&&&& {(\varprojlim F)(B)} \\
	& {F_i(A)} &&& {F_i(B)} \\
	& {F_j(A)} &&& {F_j(B)}
	\arrow["{(\varprojlim F)(f)}", dashed, from=1-1, to=1-6]
	\arrow[from=1-1, to=2-2]
	\arrow[from=1-1, to=3-2]
	\arrow[from=1-6, to=2-5]
	\arrow[from=1-6, to=3-5]
	\arrow["{F_i(f)}", from=2-2, to=2-5]
	\arrow["{(F(i\to j))_A}", from=2-2, to=3-2]
	\arrow["{(F(i\to j))_B}"', from=2-5, to=3-5]
	\arrow["{F_j(f)}"', from=3-2, to=3-5]
\end{tikzcd}.
\end{center}
The corresponding limiting cone morphisms $\varprojlim_i F_i \to F_i$ are the natural transformations induced by the limiting cone morphisms
\[
    \varprojlim_i F_i(A)\to F_i(A).
\]
Dually, one can compute colimits of functors pointwise,
defining
\[
    (\varinjlim F)(A)=\varinjlim_i F_i(A).
\]
\end{lemma}

\begin{proof}
See Chapter V.3 in \cite{mac2013categories}.
\end{proof}

\begin{corollary}[Completeness of functor categories]\label{cor:completeness_of_functor_categories}\uses{def:functor_category, lem:limits_in_functor_categories , def:limit}
If $\mcD$ is a (co)complete category and $\mcC$ any category,
then the functor category $\mathrm{Fun}(\mcC,\mcD)$ is (co)complete.
\end{corollary}

An important special case of this corollary concerns the case when $\mcD = \Set$.
Then the presheaf category $\PSh(\mcC)$ of all contravariant functors from $\mcC$ to $\Set$ is bicomplete by the previous corollary.
So in particular, every diagram in $\PSh(\mcC)$ consisting of objects arising from the
Yoneda embedding $\mcC\to\PSh(\mcC)$, i.e., consisting of representable presheaves,
has a colimit.
What presheaves can arise as those colimits?
With a slight restriction on the domain category $\mcC$,
the following lemma gives an answer to this question

\begin{lemma}\label{lem:colimits-representable}
Let $\mcC$ be a small category and $F\colon\mcC^\mathrm{op}\to\Set$ a presheaf.
Then $F$ is a colimit of representable presheaves, i.e.,
of images of the Yoneda embedding $C\mapsto y^C=\hom(-,C)$.
\end{lemma}

\begin{proof}
This is Theorem III.7.1 in~\cite{mac2013categories}.
\end{proof}

\begin{remark}
	We describe the diagram and cocone yielding $F$ as its colimit.
	In doing so,
	the Yoneda lemma has multiple occurrences, especially in the incarnation $F(C)=\hom_{\PSh}(y^C,F)$.
	
	As objects of our diagram,
	we take the presheaves $y^C$ for $C\in\mcC$ but we include each $y^C$ as many times as there
	are morphisms $y^C\to F$ (which are just elements of $F(C)$)
	so that we can index the $y^C$'s by these morphisms.
	As morphisms we take those $h\colon y^C_{f}\to y^D_{g}$ ($f\in F(C)$, $g\in F(D)$)
	for which the obvious diagram also involving $F$ commutes.
	(That is, $f=g\circ h$.)
	
	In other words,
	the full subcategory of the comma category $(\PSh(\mcC)\downarrow F)$ with objects the representable
	presheaves is the index category and the functor describing the diagram is the restriction of
	\[
		(\PSh(\mcC)\downarrow F)\to\PSh(\mcC),\quad (G,f)\mapsto G
	\]
	to this subcategory.
	We also denote this by
	\[F=\varinjlim_{y^C\to F}y^C.\]
	The description of the cone is trivial.
	For $C\in\mcC$ and $f\colon y^C\to F$ we simply use $f$ as the morphism $y^C_f \to F$.
\end{remark}

In other words, the Yoneda embedding (freely) generates the presheaf category under colimits.
We will investigate this notion in the context of cocompletions
more closely in section~\ref{ssec:completions}.

\subsection{Permanence of (co)limits}

Having dealt with different notions of projectivity, we will dive deeper into different properties and definitions of limits and colimits.
We suggest the reader to skip the rest of the section at first reading
and continue directly with completions in section~\ref{ssec:completions}.

\begin{lemma}[Monics and equalizers]\label{lem:eq_comp_monic}\uses{def:special_limits_and_colimits}
For any two morphisms $f,g\colon A\to B$ their equalizer $e\colon E\to A$ is also the equalizer of $mf$ and $mg$ for any monic $m\colon B\to C$,
\begin{center}
\begin{tikzcd}
	E & A & B & C
	\arrow["e", from=1-1, to=1-2]
	\arrow["g"', shift right, from=1-2, to=1-3]
	\arrow["f", shift left, from=1-2, to=1-3]
	\arrow["m", hook, from=1-3, to=1-4]
\end{tikzcd}.
\end{center}
Dually, coequalizers of $f$ and $g$ agree with coequalizers of $fp$ and $gp$
for any epic $p$.
\end{lemma}

\begin{proof}
This follows directly from noting that for any $h\colon X\to A$
the condition $mfh = mgh$ implies $fh = gh$.
\end{proof}

\begin{lemma}[Pullback of monic is monic]\label{lem:pull_of_mon_is_mon}\uses{def:special_limits_and_colimits}
The pullback $m^{*}\colon C\times_D B\to C$ of a monomorphism $m\colon B\to D$ along any morphism $ C\to D$ is monic,
\begin{center}
\begin{tikzcd}
	{C\times_DB} & B \\
	C & D
	\arrow["{h^*}", from=1-1, to=1-2]
	\arrow["{m^*}"', hook', from=1-1, to=2-1]
	\arrow["m", hook', from=1-2, to=2-2]
	\arrow["h"', from=2-1, to=2-2]
\end{tikzcd}.
\end{center}
Dually, every pushout $e_*$ of an epimorphism $e$ is epic, i.e.
\begin{center}
\begin{tikzcd}
	A & B \\
	C & {B\times_AC}
	\arrow["e", two heads, from=1-1, to=1-2]
	\arrow["k"', from=1-1, to=2-1]
	\arrow["{k_*}", from=1-2, to=2-2]
	\arrow["{e_*}"', two heads, from=2-1, to=2-2]
\end{tikzcd}.
\end{center}
\end{lemma}

\begin{proof}
For $f,g\colon A\to C\times_D B$ with $m^*f = m^*g$
we have that $mh^*f = hm^*f = hm^*g = mh^*g$ and since $m$ is monic, $h^*f = h^*g$.
The universal property of the pullback then implies $f = g$.
\end{proof}

\begin{lemma}\label{lem:colim-preserve-epis}
If $D,E\colon\mcI\to\mcC$ are diagrams such that their colimits exist
and $\eta\colon D\Rightarrow E$ is an epimorphism
(i.e. every component $\eta_i$ is epic),
then the induced morphism $f\colon \varinjlim D\to\varinjlim E$ is epimorphic.

Dually, if $\eta$ is a monomorphism,
then the induced morphism $f\colon \varprojlim D\to\varprojlim E$ is monic,
if the limits exist.
\end{lemma}

\begin{proof}
Consider the case of colimits, where $\eta$ is an epimorphism.
Let $h,k\colon\varinjlim E\to C$ such that $hf = kf$.
In particular, $h\tau^E_i\eta_i = hf\tau^D_i = kf\tau^D_i = h\tau^E_i\eta_i$
and thus $h\tau^E_i = k\tau^E_i$.
Uniqueness in the colimit property then implies $h = k$.
\end{proof}

Generally, colimits do not preserve monomorphisms,
neither do limits preserve epimorphisms.
For example, consider the diagram
\begin{center}
\begin{tikzcd}
	\emptyset & {\{0,1\}} && {\{0,1\}} \\
	{\{0,1\}} & {\{0,1\}} && {*}
	\arrow["{\mathrm{eq}}", from=1-1, to=1-2]
	\arrow[from=1-1, to=2-1]
	\arrow["0", shift left, from=1-2, to=1-4]
	\arrow["1"', shift right, from=1-2, to=1-4]
	\arrow["{\mathrm{id}}"', two heads, from=1-2, to=2-2]
	\arrow[two heads, from=1-4, to=2-4]
	\arrow["{\mathrm{eq}}"', from=2-1, to=2-2]
	\arrow[shift left, from=2-2, to=2-4]
	\arrow[shift right, from=2-2, to=2-4]
\end{tikzcd}
\end{center}
in $\Set$ where the two rows are equalizer diagrams
and the two vertical arrows on the right are the morphisms between the diagrams
which are epimorphic, however the induced morphism $\emptyset\to\{0,1\}$ is not epic.

Later we will get to know a large class of categories where certain types of (co)limits,
however, have this property.
The proofs of all of the preceding lemmas also follow from the fact
that monomorphisms can be described in terms of limits,
and epimorphisms in terms of colimits,
see remark~\ref{rem:epi_as_colimit},
and that (co)limits commute with (co)limits,
see lemma~\ref{lem:limits_commute_with_limits}.

\begin{lemma}[Pasting law for pullbacks]\label{lem:pasting_lemma}\uses{def:special_limits_and_colimits}
Consider any commutative diagram of the form
\begin{center}
\begin{tikzcd}
	\bullet & \bullet & \bullet \\
	\bullet & \bullet & \bullet
	\arrow[from=1-1, to=1-2]
	\arrow[from=1-1, to=2-1]
	\arrow[from=1-2, to=1-3]
	\arrow[from=1-2, to=2-2]
	\arrow[from=1-3, to=2-3]
	\arrow[from=2-1, to=2-2]
	\arrow[from=2-2, to=2-3]
\end{tikzcd}.
\end{center}
\begin{enumerate}[(i)]
	\item If the right square is a pullback, the outer rectangle is a pullback if and only if the left square is a pullback.
	\item Dually, if the left square is a pushout, the outer rectangle is a pushout if and only if the left square is a pushout.
\end{enumerate}
\end{lemma}

\begin{proof}
See e.g. pasting law for pullbacks in nLab
\footnote{\url{https://ncatlab.org/nlab/show/pasting+law+for+pullbacks}}.
\end{proof}

The following is a further commutation property some special limits and colimits might have.
Let $\mcI$ be an index category and $D$ a diagram over $\mcI$ with colimit $Y$.
Furthermore, let $f\colon X\to Y$ and assume that all required pullbacks exist.
Then there exists a diagram $D\times_Y X$ over $\mcI$
given by the universal property of the pullback:
\begin{center}
\begin{tikzcd}
	{D(i)\times_Y X} \\
	& {D(j)\times_Y X} & X \\
	& {D(j)} & Y \\
	& {D(i)}
	\arrow[dashed, from=1-1, to=2-2]
	\arrow[shift left=2, from=1-1, to=2-3]
	\arrow[from=1-1, to=4-2]
	\arrow[from=2-2, to=2-3]
	\arrow[from=2-2, to=3-2]
	\arrow[from=2-3, to=3-3]
	\arrow[from=3-2, to=3-3]
	\arrow[from=4-2, to=3-2]
	\arrow[from=4-2, to=3-3]
\end{tikzcd}.
\end{center}
Thus we have constructed a cocone $D\times_Y X \Rightarrow X$.

\begin{definition}[Stable under base change]\label{def:col_stable_under_base_change}\uses{def:special_limits_and_colimits}
    We say that colimits of the form $\mcI$ are \idx{stable under base change}
    if the natural transformation $D\times_Y X\Rightarrow X$ is a colimiting cocone
    of $D\times_Y X$ for every diagram $D$ over $\mcI$ with colimit $Y$
    and every morphism $f\colon X\to Y$.
\end{definition}

\begin{remark}
Some authors instead say that colimits of this form are \emph{universal},
e.g.~\cite{Borceux2008}.
\end{remark}

\begin{lemma}
Assume that colimits of form $\mcI$ are stable under base change.
If $f\colon X\to Y$ is a morphism, $D$ a diagram over $\mcI$ with colimit $C$
and $D\Rightarrow Y$ a cocone over $Y$,
the following diagram is a pullback square:
\begin{center}
\begin{tikzcd}
	{\varinjlim_{\mcI}(D(i)\times_Y X)} & X \\
	{\varinjlim_{\mcI} D(i)} & Y.
	\arrow[from=1-1, to=1-2]
	\arrow[from=1-1, to=2-1]
	\arrow[from=1-2, to=2-2]
	\arrow[from=2-1, to=2-2]
\end{tikzcd}
\end{center}
We also say the colimit \emph{distributes} over the pullback.
\end{lemma}
\begin{proof}
For any $i\in\mcI$ we can form the double pullback diagram
\begin{center}
\begin{tikzcd}
	{D(i)\times_C(C\times_Y X)} & {C\times_Y X} & X \\
	{D(i)} & C & Y
	\arrow[from=1-1, to=1-2]
	\arrow[from=1-1, to=2-1]
	\arrow[from=1-2, to=1-3]
	\arrow[from=1-2, to=2-2]
	\arrow[from=1-3, to=2-3]
	\arrow[from=2-1, to=2-2]
	\arrow[from=2-2, to=2-3]
\end{tikzcd}.
\end{center}
By the pasting lemma, we have
$D(i)\times_C(C\times_Y X) \cong D(i)\times_Y X$.
Now since colimits over $\mcI$ are stable under base change,
we obtain
\[
          (\varinjlim_{\mcI} D(i))\times_Y X
    \cong C\times_Y X
    \cong \varinjlim_{\mcI} (D(i)\times_C(C\times_Y X))
    \cong \varinjlim_{\mcI} (D(i)\times_Y X).
\]
\end{proof}

\begin{example}
The most important cases are those of coproducts and of coequalizers.
We say that \emph{coproducts are stable},
if colimits over every discrete diagram are stable (under base change).
We say that \emph{coequalizers are stable},
if colimits over the diagram
\begin{center}
\begin{tikzcd}
	\bullet && \bullet
	\arrow[shift left, from=1-1, to=1-3]
	\arrow[shift right, from=1-1, to=1-3]
\end{tikzcd}
\end{center}
are stable (under base change).
\end{example}

\begin{example}\label{ex:CHaus_stable_base_change}
In the category $\CHaus$, finite coproducts are stable under base change.
\end{example}

\subsection{Continuous functors}

Having a definition of limit, one has a canonical definition of continuity and preservation of limits.

Consider a functor $F\colon \mcC\to\mcD$ and a diagram $D\colon\mcI\to\mcC$.
Then the composition $FD\colon\mcI\to\mcD$ is a diagram in $\mcD$.
If $(C,(f_i)_{i\in\mcI})$ is a cone for $D$,
then after applying $F$, we obtain a cone $(F(C),(F(f_i))_{i\in\mcI})$ for $FD$.
This is essentially whiskering the natural transformation $f\colon C\Rightarrow D$
with the functor $F$,
which therefore yields a natural transformation $F.f\colon F(C)\Rightarrow FD$.

\begin{definition}[Continuous functor]\label{def:(co)continuous_functor}\uses{def:functor, def:limit, def:limit}
Let $F\colon\mcC\to\mcD$ and $D\colon\mcI\to\mcC$ as before.
\begin{enumerate}[(i)]
	\item We say that $F$ \idx{preserves} (or \idx{commutes with})
	    the limit $(L,(\pi_i)_{i\in\mcI})$ of $D$,
	    if the cone $(F(L),(F(\pi_i))_{i\in\mcI})$ is a limit of $FD$,
	    and we write
        \[
            \varprojlim_{i\in\mcI} F(D(i))=F(\varprojlim_i D(i)).
        \]
	\item Dually, $F$ preserves the colimit $(L,\tau_i)$,
	    if the cocone $(F(L),F(\tau_i))$ is a colimit of $FD$, i.e.
        \[
            \varinjlim_{i\in\mcI} F(D(i))=F(\varinjlim_i D(i)).
        \]
	\item We say that $F$ \textbf{reflects}\index{reflect} the (co)limit of $D$,
	    if for any (co)cone $(C,f_i)$ to $D$,
	    if $(F(C),F(f_i))$ is a (co)limit of $FD$, then $(C,f_i)$ is a (co)limit of $D$.
	\item The functor $F$ \idx{creates} the (co)limit of $D$,
	    if it reflects it and the existence of a (co)limit of $FD$
	    implies the existence of a (co)limit of the form $(F(C),F(f_i))$.
\end{enumerate}
A functor $F\colon \mcC\to\mcD$ is called \textbf{(co)continuous}\index{continuous functor}\index{cocontinuous functor},
if it preserves all (co)limits that exist in $\mcC$.
\end{definition}

The intuition of preserving colimits is clear, reflecting means that as soon as we have a cone in $\mcC$ that is mapped (by $F$) to a limiting cone in $\mcD$, then the cone in $\mcC$ already was a limiting cone.
Creating colimits means that as son as we have a colimiting cone to some diagram in $\mcD$ coming from $\mcC$, there exists a colimiting cone in $\mcC$ that is the preimage of the colimiting cone in $\mcD$.

\begin{remark}
Note that some authors, like \cite[Chapter V.1]{mac2013categories},
define creation of limits with an additional uniqueness property.
We follow \cite[Def 3.3.1]{Riehl}.
\end{remark}


One very important example of continuous functors is given by $\hom$-functors.

\begin{lemma}[Continuity of Hom-functor]\label{lem:hom_functor_preserves_limits}\uses{def:(co)continuous_functor, def:limit, lem:set_bicomplete}
The covariant $\hom$-functor $\hom(A,-)\colon\mcC\to\Set$ is continuous,
i.e. for every diagram $D\colon\mcI\to\mcC$ one obtains
\[
    \hom(A,\varprojlim_{i\in\mcI} D(i))=\varprojlim_{i\in\mcI} \hom(A,D(i)).
\]
Dually, the contravariant $\hom$-functor $\hom(-,A)\colon\mcC^\mathrm{op}\to\Set$
is continuous, i.e. for every diagram $D\colon\mcI\to\mcC$ one obtains
\[
    \hom(\varinjlim_{i\in\mcI} D(i),A)=\varprojlim_{i\in\mcI} \hom(D(i),A).
\]
\end{lemma}

\begin{remark}
Note that for the contravariant $\hom$-functor the colimit $\varinjlim_i D(i)$
is taken in the category $\mcC$,
so that it becomes the limit $\varprojlim_{i} D(i)$ in the opposite category $\mcC^\mathrm{op}$
and thus $\hom(-,A)$ is continuous (takes limits to limits) as a functor from $\mcC^\mathrm{op}$ to $\Set$.
\end{remark}

\begin{proof}
This is precisely the definition of (co)limits.
\end{proof}

Together with the description of (co)limits in functor categories,
lemma~\ref{lem:limits_in_functor_categories}, this yields the following corollary.

\begin{corollary}[Continuity of the Yoneda embedding]\label{cor:cocontinuity_of_yoneda}\uses{lem:hom_functor_preserves_limits, def:functor_category, lem:limits_in_functor_categories , def:(co)continuous_functor, lem:yoneda}
The contravariant Yoneda embedding
$\mcC^\mathrm{op}\to\mathrm{Fun}(\mcC, \mathrm{Set})$
and the covariant Yoneda embedding
$\mcC\to\PSh(\mcC)$ are continuous.
\end{corollary}

\begin{example}\label{ex:comma-cat-creates-limits}
Let $F\colon\mcC\to\mcD$ be a functor and $D\in\mcD$.
Then we have a forgetful functor $P\colon (F\downarrow D)\to \mcC$
mapping an object $(C,f)$ to $C$ and a morphism
$g\colon (C,f)\to(C',f')$ to $g\colon C\to C'$.
If $F$ preserves colimits, then $P$ creates colimits.
Dually, the forgetful functor $Q\colon (D\downarrow F)\to \mcC$
creates limits if $F$ preserves limits,
see Chapter V.6 in \cite{mac2013categories}.

In particular, if $F$ is cocontinuous and $\mcC$ is cocomplete,
then $(F\downarrow D)$ is cocomplete as well.
The dual holds for $(D\downarrow F)$.
\end{example}

\begin{example}
  Forgetful functors from \enquote{algebraic} categories
  such as $\Ab\to\Set$, $\ModR\to\Set$, ... or even $\Ab(\Top)\to\Top$
  (the former being the category of topological abelian groups)
  always create limits.
  We will prove precise (and shockingly general) versions of this statement in \ref{sec:univ-alg}.
\end{example}

Another example of continuous functors is given as follows.

\begin{example}\label{ex:pullback-functor}
Suppose $\mcC$ is a category with pullbacks and
where colimits of the form $\mcI$ are stable under base change.
Then for every morphism $f\colon X\to Y$ in $\mcC$,
the induced \idx{pullback functor} $f^\ast\colon \mcC/Y\to\mcC/X$
between the comma categories preserves colimits of the form $\mcI$.
\end{example}

Next, we introduce an important criterion for continuity of functors,
which follows from proposition~\ref{prop:criterion_for_completeness}.

\begin{corollary}[Criterion for (co)continuity)]\label{prop:criterion_for_continuity}\uses{def:(co)continuous_functor,def:limit, lem:lim_via_special_lim}
A functor with (finitely) complete domain is (finitely)
(co)continuous precisely if it preserves all (reflexive) (co)equalizers and all (finite) (co)products.
\end{corollary}

The following result explains the intuition of sifted colimits being essentially just filtered colimits and reflexive coequalizers.

\begin{lemma}[Commuting with sifted colimits]\label{lem:comm_with_sifted}\uses{def:special_limits_and_colimits}
A functor with finitely cocomplete domain preserves sifted colimits
precisely if it preserves reflexive coequalizers and filtered colimits.
\end{lemma}

\begin{proof}
See \cite[Theorem 2.1]{Adamek2010}.
\end{proof}

The defect of the $\hom$-functor $\hom(A,-)$ to preserve colimits
gives rise to many of the most important special classes of objects.
We define various classes of objects in a category $\mcC$
for which the covariant $\hom$-functor commutes with more colimits.

\begin{definition}[Special objects]\label{def:special_objects}\uses{def:category, lem:hom_functor_preserves_limits, lem:comm_with_sifted}
Let $\mcC$ be a category.
\begin{enumerate}[(i)]
	\item A \idx{tiny object} is an object $T\in\mcC$
	such that $\hom(T,-)$ preserves colimits.
	\item A ($\kappa$-)\idx{compact object} is an object $K\in\mcC$ such that
	$\hom(K,-)$ preserves ($\kappa$-)filtered colimits.
	\item A \idx{small object} is an object that is $\kappa$-compact for some regular cardinal $\kappa$.
	\item An \idx{effective projective object} is an object $P\in\mcC$
	such that $\hom(P,-)$ commutes with reflexive coequalizers.
	\item A \idx{compact projective object} is an object $K\in\mcC$
	such that $\hom(K,-)$ preserves sifted colimits.
	\item A \idx{connected object} is an object $C\in\mcC$
	such that $\hom(C,-)$ preserves coproducts.
\end{enumerate}
\end{definition}

\begin{remark}\label{rem:special_objects_other_names}
In some cases, especially in algebraic contexts, e.g. in \cite{Adamek2010a},
a compact object is also called \emph{finitely presentable}
and a compact projective object is called \emph{perfectly presentable}
or \emph{strongly finitely presentable}.

By lemma~\ref{lem:comm_with_sifted},
an object is compact projective if and only if it is compact and effective projective.
\end{remark}

The tiny objects play a central role in the presheaf category $\PSh(\mcC)$ of a category $\mcC$.

\begin{lemma}[Representable is tiny]\label{lem:representable_tiny}\uses{lem:yoneda,def:special_objects,lem:limits_in_functor_categories,lem:retract-diagram-iso}
    Let $\mcC$ be a small category.
	A presheaf $F\in\PSh(\mcC)$ is tiny
	precisely if it is a retract of a representable presheaf.
	\end{lemma}
\begin{proof}
    Since the canonical map
    \[
        \varinjlim\hom(F,G_i) \to \hom(F,\varinjlim G_i)
    \]
    is natural in $F$, the class of tiny objects is closed under retracts
	by lemma~\ref{lem:retract-diagram-iso}.
	Thus it suffices to show that any presheaf of the form $y^X=\hom(-,X)$
	for $X\in\mcC$ is tiny in $\PSh(\mcC)$.
	Hence consider any diagram $D\colon\mcI\to \PSh(\mcC)$ with colimit $\varinjlim D$,
	which by lemma \ref{lem:limits_in_functor_categories} is computed pointwise.
	Now, using the Yoneda lemma twice,
	\[
	\hom(y^X,\varinjlim_i D_i)
	= (\varinjlim_i D_i)(X)
	= \varinjlim_i (D_i(X))
	= \varinjlim_i \hom(y^X, D_i).
	\]
	On the other hand, let $F$ be tiny.
	Every presheaf $F$ is the colimit of representable presheaves $y^X$,
	see lemma~\ref{lem:colimits-representable}.
	Therefore we obtain
	\[
	    \hom(F,F) = \hom(F,\varinjlim y^{X_i}) = \varinjlim\hom(F,y^{X_i})
    \]
    where the colimit runs over every $f\colon y^X\to F$ for every $X\in\mcC$.
    In particular, there exists $g\colon F\to y^X$ such that
    $1_F = fg$.
    So $F$ is a retract of $y^X$.
\end{proof}

\begin{remark}
  The smallness condition on $\mcC$ cannot be dropped.
\end{remark}

We are going to prove more stability properties of the introduced objects later
in lemma~\ref{lem:stability-compact-projective}.
For now, let us give a characterisation of compact objects.

\begin{lemma}\label{lem:comp_hom_endl}\uses{lem:lim_via_special_lim}
    Let $\mcC$ be a finitely cocomplete category.
    An object $K$ of $\mcC$ is compact precisely if for every diagram $D\colon\mcJ\to \mcC$ which has a colimit in $\mcC$,
    every map $f\colon K\to \varinjlim_{\mcJ} D$ factors through the colimit
    of a finite subdiagram $\mcF\subset\mcJ$, i.e.
    \[
        f\colon K\to \varinjlim_{\mcF}D\to \varinjlim_{\mcI} D.
    \]
\end{lemma}

\begin{proof}
    By lemma~\ref{lem:lim_via_special_lim} (iii) we have that
    \[
        \varinjlim_{\mcJ} D \simeq \varinjlim_{\mcF\subset\mcJ}\varinjlim_{\mcF} D
    \]
    where $\mcF\subset\mcJ$ is a finite full subcategory
    and hence the outer colimit is filtered since (even directed on the nose).
    If now $K$ is compact, we obtain
    \[
        \hom(K,\varinjlim_{\mcJ} D)
        = \hom(K,\varinjlim_{\mcF\subset\mcJ}\varinjlim_{\mcF} D)
        = \varinjlim_{\mcF\subset\mcJ}\hom(K,\varinjlim_{\mcF} D)
    \]
    and therefore $f$ factors as required.

    For the other direction let $\mcJ$ be a filtered category.
    By assumption one has
    \[
        \hom(K,\varinjlim_{\mcJ}D)
        = \varinjlim_{\mcF\subset\mcJ}\hom(K,\varinjlim_{\mcF} D).
    \]
    Now the diagonal functor $\mcJ\to{\mcJ}^{\mcF}$ is final and hence
    there exists $j\in\mcJ$ such that $\varinjlim_{\mcF} D\to \varinjlim_{\mcJ} D$
    factors as $\varinjlim_{\mcF} D\to D(j)\to \varinjlim_{\mcJ} D$ which implies
    \[
        \hom(K,\varinjlim_{\mcJ}D) = \varinjlim_{\mcJ}\hom(K,D).
	\]
\end{proof}

\begin{remark}\label{rem:kappa_commutes_into_subdiag}
The same assertion holds for compact projective objects with discrete finite subdiagrams for $\mcF$,
and for $\kappa$-compact objects with $\kappa$-small subdiagrams.
\end{remark}

\begin{remark}\label{rem:comp-factorization}
    The arrow \(\varinjlim_{\mcF}D\to \varinjlim_{\mcJ}D\)
    in this factorization need not be monic nor epic.
    For example, if \(\varinjlim_{\mcJ} D = \coprod_{\mcJ} D(j)\) and \(K = D(j_0)\)
    then we factorize
    \[
        K \to D(j_0) \to \coprod_{\mcJ} D(j)
    \]
    but the latter morphism is not epimorphic.

    Similarly when taking a diagram over \(\N^\mathrm{op}\) then its colimit is the value of the diagram at \(1\),
    so if we set this value to be \(*\) the second arrow in the factorization is in general not monic.
\end{remark}

\subsection{Interchanging limits and colimits}

Another important example of continuous functors are limit functors themselves,
meaning that the order of iterated (co)limits is irrelevant.

Let $\mcI$ and $\mcJ$ be index categories and $D\colon\mcI\times\mcJ\to\mcC$ a diagram
in a category $\mcC$.
Thus, we can interpret the diagram $D$ as a diagram over $\mcJ$ in $[\mcI,\mcC]$.
If all the involved limits exist, one obtains a diagram
\[
    \varprojlim_{j\in\mcJ}D(-,j)\colon\mcI\to\mcC,\quad i\mapsto\varprojlim_{j\in\mcJ} D(i,j),
\]
which essentially comes down to computing the limit
\[
    \varprojlim_{j\in\mcJ} D(i,j)
\]
in $\mcC$ for every $i\in\mcI$ by lemma~\ref{lem:limits_in_functor_categories}.
Now one can again compute the limit
\[
    \varprojlim_{i\in\mcI}\varprojlim_{j\in\mcJ}D(i,j)
\]
in $\mcC$ of the limit functor $\varprojlim_{j\in\mcJ}D(-,j)$.
Analogously, one defines $\varprojlim_{j\in\mcJ}\varprojlim_{i\in\mcI}D(i,j)$.
Now, using the universal properties of limits,
one obtains morphisms $\mu_i$ and $\nu_j$ for every $i\in\mcI$ and $j\in\mcJ$
which in turn yield the canonical morphisms $\mu$ and $\nu$,
as depicted in the diagram.
\begin{center}
\begin{tikzcd}
	{D(i,j)} & {\varprojlim\limits_{j\in\mcJ}D(i,j)} & {\varprojlim\limits_{i\in\mcI}\varprojlim\limits_{j\in\mcJ}D(i,j)} \\
	{\varprojlim\limits_{i\in\mcI}D(i,j)} \\
	{\varprojlim\limits_{j\in\mcJ}\varprojlim\limits_{i\in\mcI}D(i,j)}
	\arrow["{\pi_{i,j}}"', from=1-2, to=1-1]
	\arrow["{\pi_i}"', from=1-3, to=1-2]
	\arrow["{\nu_j}"{pos=0.4}, dashed, from=1-3, to=2-1]
	\arrow["\nu", shift left=4, curve={height=-12pt}, dashed, from=1-3, to=3-1]
	\arrow["{\tau_{i,j}}", from=2-1, to=1-1]
	\arrow["{\mu_i}"', curve={height=6pt}, dashed, from=3-1, to=1-2]
	\arrow["\mu", shift right, curve={height=12pt}, dashed, from=3-1, to=1-3]
	\arrow["{\tau_j}", from=3-1, to=2-1]
\end{tikzcd}
\end{center}
Furthermore, $\mu$ and $\nu$ are natural in $D$,
meaning that for every natural transformation $\eta\colon D\Rightarrow E$,
the diagram

\begin{adjustbox}{max width=\textwidth}
\begin{tikzcd}
	& {E(i,j)} && {\varprojlim\limits_{j\in\mcJ}E(i,j)} && {\varprojlim\limits_{i\in\mcI}\varprojlim\limits_{j\in\mcJ}E(i,j)} \\
	{D(i,j)} && {\varprojlim\limits_{j\in\mcJ}D(i,j)} && {\varprojlim\limits_{i\in\mcI}\varprojlim\limits_{j\in\mcJ}D(i,j)} \\
	& {\varprojlim\limits_{i\in\mcI}E(i,j)} && {\varprojlim\limits_{j\in\mcJ}\varprojlim\limits_{i\in\mcI}E(i,j)} \\
	{\varprojlim\limits_{i\in\mcI}D(i,j)} && {\varprojlim\limits_{j\in\mcJ}\varprojlim\limits_{i\in\mcI}D(i,j)}
	\arrow["{p_{i,j}}"', from=1-4, to=1-2]
	\arrow["{p_i}"', from=1-6, to=1-4]
	\arrow["{\eta_{i,j}}", from=2-1, to=1-2]
	\arrow["{\eta_{i}}", dashed, from=2-3, to=1-4]
	\arrow["{\pi_{i,j}}"'{pos=0.7}, from=2-3, to=2-1]
	\arrow["\vartheta", dashed, from=2-5, to=1-6]
	\arrow["{\pi_i}"'{pos=0.2}, from=2-5, to=2-3]
	\arrow["{t_{i,j}}"'{pos=0.7}, from=3-2, to=1-2]
	\arrow["{m_i}"{pos=0.3}, from=3-4, to=1-4]
	\arrow["m"', curve={height=30pt}, from=3-4, to=1-6]
	\arrow["{t_j}"'{pos=0.7}, from=3-4, to=3-2]
	\arrow["{\tau_{i,j}}", from=4-1, to=2-1]
	\arrow["{\eta_j}", dashed, from=4-1, to=3-2]
	\arrow["{\mu_i}"'{pos=0.3}, from=4-3, to=2-3]
	\arrow["\mu"', curve={height=30pt}, from=4-3, to=2-5]
	\arrow["\theta", dashed, from=4-3, to=3-4]
	\arrow["{\tau_j}"', from=4-3, to=4-1]
\end{tikzcd}
\end{adjustbox}
is commutative (and the same for $\nu$).
Moreover, the morphisms $\mu$ and $\nu$ are inverse to each other,
as the following lemma explains.

\begin{lemma}[Limits commute with limits]\label{lem:limits_commute_with_limits}\uses{def:limit, def:functor}
Limits commute with limits, i.e.,
in the situation above,
as soon as one of the objects exists,
\[
    \varprojlim_{j\in\mcJ}\varprojlim_{i\in\mcI}D(i,j)
    \simeq\varprojlim_{i\in\mcI}\varprojlim_{j\in\mcJ}D(i,j)
    \simeq \varprojlim_{(i,j)\in\mcI\times\mcJ}D(i,j),
\]
where all isomorphisms are given by the canonical morphisms.
Similarly, colimits commute with colimits, i.e.,
\[
    \varinjlim_{j\in\mcJ}\varinjlim_{i\in\mcI}D(i,j)
    \simeq\varinjlim_{i\in\mcI}\varinjlim_{j\in\mcJ}D(i,j)
    \simeq \varinjlim_{(i,j)\in\mcI\times\mcJ}D(i,j),
\]
if one of the objects exists.
\end{lemma}

\begin{proof}
See Chapters IX.2 and IX.8 in \cite{mac2013categories}.
\end{proof}

Now one may ask, what happens if we mix limits and colimits.
If we have, as before, a diagram $D\colon\mcI\times\mcJ\to\mcC$,
one can also form the colimit
and then one can form the colimit
\[
    \varinjlim_{i\in\mcI}\varprojlim_{j\in\mcJ} D(i,j).
\]
of the diagram $\varprojlim_{j\in\mcJ}D(-,j)$ over $\mcI$.
On the other hand, we can first compute the colimits over $\mcJ$
and then take the limit over $\mcJ$,
\[
    \varprojlim_{j\in\mcJ}\varinjlim_{i\in\mcI} D(i,j).
\]
Then one obtains by the universal property of limits,
for every $i\in\mcI$ a morphism
\[
    \mu_i\colon\varprojlim_{j\in\mcJ} D(i,j)\to\varprojlim_{j\in\mcJ}\varinjlim_{i\in\mcI}D(i,j)
\]
and this in turn yields a canonical morphism
\[
    \varphi\colon\varinjlim_{i\in\mcI}\varprojlim_{j\in\mcJ} D(i,j)\to\varprojlim_{j\in\mcJ}\varinjlim_{i\in\mcI} D(i,j).
\]
The situation is summarised in the following diagram.
\begin{center}
\begin{tikzcd}
	{D(i,j)} & {\varprojlim\limits_{j\in\mcJ}D(i,j)} & {\varinjlim\limits_{i\in\mcI}\varprojlim\limits_{j\in\mcJ}D(i,j)} \\
	{\varinjlim\limits_{i\in\mcI}D(i,j)} & {\varprojlim\limits_{j\in\mcJ}\varinjlim\limits_{i\in\mcI}D(i,j)}
	\arrow["{\tau_{i,j}}"', from=1-1, to=2-1]
	\arrow["{\pi_{i,j}}"', from=1-2, to=1-1]
	\arrow["{\tau_j}", from=1-2, to=1-3]
	\arrow["{\mu_i}", dashed, from=1-2, to=2-2]
	\arrow["\varphi", shift left=3, shorten >=5pt, dashed, from=1-3, to=2-2]
	\arrow["{\pi_j}", from=2-2, to=2-1]
\end{tikzcd}
\end{center}
As before, the canonical morphism is natural in $D$.
However, in general we cannot expect a morphism in the other direction,
making the diagram commutative.

\begin{warning}
It is not true that in general \idx{colimits commute with limits},
meaning that the canonical morphism
\[
    \varphi\colon\varinjlim_{i\in\mcI}\varprojlim_{j\in\mcJ} D(i,j)\to\varprojlim_{j\in\mcJ}\varinjlim_{i\in\mcI} D(i,j).
\]
is an isomorphism, see example~\ref{ex:filtered-colimit-limit-not-commute} below.
There are, however, situations when we can construct an inverse.
\end{warning}

\begin{proposition}[Characterisation of filtered and sifted categories]\label{prop:filtered_sifted_commutes_set}\uses{def:special_limits_and_colimits}
\begin{enumerate}[(i)]
    \item An index category $\mcI$ is filtered precisely if
        taking colimits of diagrams over $\mcI$ commutes with finite limits in $\Set$,
		i.e. for any finite category $\mcJ$ and diagram $D\colon \mcI\times\mcJ\to\Set$
		the canonical morphism
		\[
    		\varinjlim_{i\in\mcI} \varprojlim_{j\in\mcJ} D(i,j)\simeq \varprojlim_{j\in\mcJ}\varinjlim_{i\in\mcI} D(i,j).
		\]
		is an isomorphism.
		More generally, $\kappa$-filtered colimits are precisely those which commute with $\kappa$-small limits in $\Set$.

    \item An index category $\mcI$ is sifted precisely if
        taking colimits of diagrams with this form commute with finite products in $\Set$,
        i.e. for any finite set $\mcJ$ and diagram $D\colon \mcI\times\mcJ\to\Set$
        the canonical morphism
        \[
            \varinjlim_{i\in\mcI} \prod_{j\in\mcJ} D(i,j)\simeq \prod_{j\in\mcJ}\varinjlim_{i\in\mcI} D(i,j).
        \]
        is an isomorphism.
\end{enumerate}
\end{proposition}

\begin{proof}
We only proof (i) because the proof of (ii) is similar.
Let $\mcI$ be a filtered index category, $\mcJ$ finite and $D\colon\mcI\times\mcJ\to\Set$
a diagram.
We have to show that $\varphi$ is an isomorphism.
By the explicit construction of colimits in $\Set$, one has for every $j\in\mcJ$
\[
    \varinjlim_{i\in\mcI} D(i,j) = \left(\bigsqcup_{i\in\mcI} D(i,j)\right) / \sim
\]
where two elements $x_i\in D(i,j)$, $y_{i'}\in D(i',j)$ are equivalent precisely if
there exist $f\colon i\to i_0$, $g\colon i'\to i_0$ such that
\[
    D(f,1_j)x_i = D(g,1_j)x_{i'}.
\]
Since $\mcI$ is filtered, this indeed is an equivalence relation
(transitivity is using the filtered condition).
On the other hand, it contains the forced identifications in the colimit,
described in lemma~\ref{lem:set_bicomplete},
and clearly is the smallest equivalence relation containing those.
Thus, an element of
\[
    \varprojlim_{j\in\mcJ}\varinjlim_{i\in\mcI} D(i,j)
\]
is a tuple $([x_j,i_j])_{j\in\mcJ}$ of equivalence classes with $x_j\in D(i_j,j)$,
such that for every morphism $h\colon j\to j'$ one has
$D(i_j,j') \ni D(1_{i_j},h)(x_j,i_j) \sim (x_{j'},i_{j'}) \in D(i_{j'},j')$.
So there exist $f^j_h\colon i_j\to i_h$, $f^{j'}_h\colon i_{j'}\to i_h$ such that
\[
    D(f^j_h,h)(x_j,i_j) = D(f^{j'}_h,1_{j'})(x_{j'},i_{j'}).
\]
Now the $i_j$ for $j\in\mcJ$ together with the $i_h$ for $h\in\mathrm{mor}\mcJ$
and morphisms $f^j_h$ constitute a finite subdiagram $\mcF$ of $\mcI$.
Since the latter is filtered, there exists a cone $\eta\colon\mcF\Rightarrow i_0$ over this finite subdiagram.
Let $z_j = D(\eta_{i_j},1_j)(x_j,i_j) \in D(i_0,j)$.
Then we have for $h\colon j\to j'$, that
\begin{align*}
    D(1_{i_0},h)z_j
    &= D(\eta_{i_j},h)(x_j,i_j) \\
    &= D(\eta_{i_h},1_{j'})D(f^j_h,h)(x_j,i_j) \\
    &= D(\eta_{i_h},1_{j'})D(f^{j'}_h,1_{j'})(x_{j'},i_{j'}) \\
    &= D(\eta_{i_{j'}},1_{j'})(x_{j'},i_{j'}) \\
    &= z_{j'}.
\end{align*}
In particular, $(z_j)_{j\in\mcJ}\in\varprojlim_{j\in\mcJ}D(i_0,j)$.
Taking the equivalence class, one obtains an element of
$\varinjlim_{i\in\mcI}\varprojlim_{j\in\mcJ} D(i,j)$.
It remains to check that this element is independent of the choices made
and therefore determines a map
\[
    \varprojlim_{j\in\mcJ}\varinjlim_{i\in\mcI} D(i,j)\to \varinjlim_{i\in\mcI}\varprojlim_{j\in\mcJ} D(i,j),
\]
which is the desired inverse of $\varphi$.

If, on the other hand, colimits over $\mcI$ commute with finite limits,
then for a finite subdiagram $F\colon\mcF\to\mcI$ define a diagram
\[
    D\colon\mcI\times\mcF^\mathrm{op}\to\Set,\quad (i,t)\mapsto \hom_{\mcI}(F(t),i)
\]
and the values of $D$ on morphisms are given by pushforward and pullback, respectively.
Now, for any $t\in\mcF^\mathrm{op}$, we have
\[
    \varinjlim_{i\in\mcI} D(i,t) = \left(\bigsqcup_{i\in\mcI} \hom_{\mcI}(F(t),i) \right) / \sim
\]
and for $(f,i)\in \bigsqcup_{i\in\mcI} \hom_{\mcI}(F(t),i)$
one has $D(1_t,f)(1_{F(t)},F(t)) = (f,i)$,
so that any two elements are identified, yielding
\[
    \varinjlim_{i\in\mcI} D(i,t) = \ast
\]
for every $t\in\mcF^\mathrm{op}$.
Together with the assumption that the limit commutes with the colimit,
this implies that
\[
    \varinjlim_{i\in\mcI}\varprojlim_{t\in\mcF^\mathrm{op}} D(i,t)
    = \varprojlim_{t\in\mcF^\mathrm{op}}\varinjlim_{i\in\mcI} D(i,t) = \ast.
\]
In particular, there exists an $i_0\in\mcI$ such that
\[
    \varprojlim_{t\in\mcF^\mathrm{op}} D(i_0,t)\neq\emptyset.
\]
Hence there exists a tuple $(f_t)_{t\in\mcF^\mathrm{op}}$ of morphisms
$f_t\in\hom_{\mcI}(F(t),i_0)$ such that, for every $h\colon s\to t$ one has
$f_t = f_s\circ F(h)$.
We have constructed a cocone $f\colon F\Rightarrow i_0$.
\end{proof}

\begin{example}\label{ex:filtered-colimit-limit-not-commute}
This is a special property of $\Set$ and not true in every category.
Consider the category of all closed subsets of $\alpha\N=\N\cup \{\infty\}$ as partially ordered set, see example~\ref{ex:examples_of_categories}.
Then $[n]=\{1,\dots, n\}$ yields an increasingly filtered diagram $([n])_{n\in\N}$ with colimit $\alpha\N=\overline{\bigcup_n [n]}$.
Now, since the product is the intersection in our category,
this yields
\[
    \varinjlim([n]\times\{\infty\})=\varinjlim \emptyset =\emptyset
\]
but
\[
    (\varinjlim[n])\times \{\infty\}=\alpha\N\times\{\infty\}=\{\infty\}.
\]
Hence filtered colimits do not commute with products in this category.
\end{example}
This commutativity allows the following easy consequence.
\begin{lemma}\label{lem:kappa_comp_from_comp}
Any $\kappa$-small colimit of compact objects is an $\kappa$-compact object, and dually any $\kappa$-small limit of cocompact objects is $\kappa$-cocompact.
\end{lemma}
\begin{proof}
	Consider any $\kappa$-small colimit $\varinjlim X_i$ for compact objects $X_i$.
	Now, consider any $\kappa$-filtered colimit $\varinjlim T_j$.
	Then
	\[\hom(\varinjlim X_i, \varinjlim T_j)=\varprojlim_i \hom(X_i, \varinjlim T_j)=\varprojlim_i \varinjlim_j\hom(X_i,T_j),\]
	but now this commutes in set and we obtain
	\[\varinjlim_j\varprojlim_i\hom(X_i,T_j)=\varinjlim_j\hom(\varinjlim X_i, T_j).\]
	\end{proof}

\begin{warning}
It is not true that if the diagonal functor $\mcI\to{\mcI}^{\mcJ}$ is final,
then the canonical morphism
\[
    \varprojlim_{\mcJ}\varinjlim_{\mcI} D \to \varinjlim_{\mcI}\varprojlim_{\mcJ} D
\]
is an isomorphism for any diagram $D\colon\mcJ\times\mcI\to\Set$.
For this, consider the form $\mcI$ given by
\begin{center}
\begin{tikzcd}
	C && D \\
	& B \\
	& A
	\arrow["h", from=2-2, to=1-1]
	\arrow["k"', from=2-2, to=1-3]
	\arrow["f", curve={height=-6pt}, from=3-2, to=2-2]
	\arrow["g"', curve={height=6pt}, from=3-2, to=2-2]
\end{tikzcd}.
\end{center}
where $f$ and $g$ are distinct, but $hf = hg$ and $kf = kg$.
Then the diagonal functor $\mcI\to\mcI^{\mcJ}$,
where $\mcJ$ is the form
\begin{center}
\begin{tikzcd}
	j_0 & j_1
	\arrow[shift left, from=1-1, to=1-2]
	\arrow[shift right, from=1-1, to=1-2]
\end{tikzcd}
\end{center}
is final.
However, for the diagram
\[
    D\colon\mcI\times\mcF\to\Set,\quad (X,j)\mapsto
    \begin{cases}
         \hom_{\mcI}(A,X),& j = j_1 \\
         \hom_{\mcI}(B,X),& j = j_0
    \end{cases},
\]
the limit over the colimits is the point $\ast$,
but
\[
    \varinjlim_{X\in\mcI}\varprojlim_{E\in\mcF} D(X,E) = \{ h,k \}.
\]

\end{warning}

Let us end with the announced stability properties of the special objects of the previous chapter.

\begin{lemma}\label{lem:stability-compact-projective}
Let $\mcC$ be a category.
Then the following hold.
\begin{enumerate}[(i)]
    \item The class of compact objects is closed under finite colimits and retracts.
    \item The class of effective projective objects is closed under finite coproducts and retracts.
    \item The class of compact projective objects is closed under finite coproducts and retracts.
\end{enumerate}
\end{lemma}

\begin{proof}
The proof for the closure under retracts of all three classes works
the same as the one for tiny objects.
The other assertions follow from
propositions~\ref{prop:filtered_sifted_commutes_set} and~\ref{lem:hom_functor_preserves_limits}.
\end{proof}

\subsection{Distinguished morphisms and objects}

Using limits and colimits, we can refine our understanding of morphisms.

\begin{definition}[Kernel pair]\label{def:kernel_pair}\uses{def:category, def:special_limits_and_colimits}
    The \textbf{kernel}\index{kernel pair} and \idx{cokernel pair} of a morphism $f\colon A\to B$ are defined (if they exist)
    as the projection maps $k_{1}$, $k_{2}$ from to pullback of $f$ with itself (to $A$)
    and maps $r_{1}$, $r_{2}$ from the pushout of $f$ with itself,
\begin{center}\begin{tikzcd}
	{A\times_BA} & A \\
	A & B \\
	B & {B\sqcup_A B}
	\arrow["{k_2}", from=1-1, to=1-2]
	\arrow["{k_1}"', from=1-1, to=2-1]
	\arrow["f", from=1-2, to=2-2]
	\arrow["f", from=2-1, to=2-2]
	\arrow["f"', from=2-1, to=3-1]
	\arrow["{r_1}"', from=3-1, to=3-2]
	\arrow["{r_2}", from=2-2, to=3-2]
\end{tikzcd}.\end{center}
	The \idx{quotient} of $B$ by $A$ corresponding to $f\colon A\to B$ is defined as the coequalizer of the kernel pair, and often denoted by $B/A$.
\end{definition}

\begin{remark}\label{rem:epi_as_colimit}
A morphism $f\colon A\to B$ is epic precisely if $1_B,1_B$ is a cokernel pair of $f$, i.e.

\begin{center}\begin{tikzcd}
	A & B \\
	B & B
	\arrow["f", from=1-1, to=1-2]
	\arrow["f"', from=1-1, to=2-1]
	\arrow["{1_B}", from=1-2, to=2-2]
	\arrow["{1_B}"', from=2-1, to=2-2]
\end{tikzcd}\end{center}
is a pushout diagram.
Hence epimorphisms can be described with colimits and in particular any colimit preserving functor preserves epis.
Specializing to the colimit functor this also yields that any epic natural transformation between diagrams induces an epimorphism between the colimits.

Dually, a morphism is monic if and only if $1_A,1_A$ is a kernel pair of $f$, i.e.

\begin{center}\begin{tikzcd}
	A & A \\
	A & B
	\arrow["{1_A}", from=1-1, to=1-2]
	\arrow["{1_A}"', from=1-1, to=2-1]
	\arrow["f", from=1-2, to=2-2]
	\arrow["f"', from=2-1, to=2-2]
\end{tikzcd}\end{center}
is a pullback diagram, also implying that all limit preserving functors preserve monics.
This in particular implies that any monomorphism between diagrams induces a monic of limits,
see lemma~\ref{lem:colim-preserve-epis}.
\end{remark}

\begin{definition}[Special mono-/epimorphisms]\label{def:special_monoepi}\uses{def:kernel_pair,def:category, def:special_morphisms}
An epimorphism $f\colon X\to Y$ is called
\begin{enumerate}[(i)]
\item \textbf{regular}\index{regular epimorphism},
if it is the coequalizer of a pair of morphisms
\begin{center}\begin{tikzcd}
	K & X & Y
	\arrow[shift left, from=1-1, to=1-2]
	\arrow[shift right, from=1-1, to=1-2]
	\arrow["f", from=1-2, to=1-3]
\end{tikzcd},\end{center}
\item \textbf{effective}\index{effective epimorphism},
if it is the coequalizer of its kernel pair
\begin{center}\begin{tikzcd}
	{X\times_Y X} & X & Y
	\arrow[shift left, from=1-1, to=1-2]
	\arrow[shift right, from=1-1, to=1-2]
	\arrow["f", from=1-2, to=1-3]
\end{tikzcd}.\end{center}

\end{enumerate}
Dually, a monomorphism $f\colon X\to Y$ is called
\begin{enumerate}[(i)]
\item \textbf{regular}\index{regular monomorphism}, if it is the equalizer of a pair of morphisms
\begin{center}\begin{tikzcd}
	X & Y & K
	\arrow["f", from=1-1, to=1-2]
	\arrow[shift right, from=1-2, to=1-3]
	\arrow[shift left, from=1-2, to=1-3]
\end{tikzcd},\end{center}
\item \textbf{effective}\index{effective monomorphism}, if it is the equalizer of its cokernel pair
\begin{center}\begin{tikzcd}
	X & Y & {Y\sqcup_X Y}
	\arrow["f", from=1-1, to=1-2]
	\arrow[shift right, from=1-2, to=1-3]
	\arrow[shift left, from=1-2, to=1-3]
\end{tikzcd}.\end{center}

\end{enumerate}
\end{definition}

\begin{example}
  In $\Top$,
  the monomorphisms are precisely the injective maps.
  Regular monomorphisms are those monomorphisms $A\hookrightarrow B$
  for which $A$ carries the subspace topology.
\end{example}

\begin{remark}\label{rem:special_monoepi_reflexive}
Every effective epimorphism is in particular a reflexive coequalizer,
since there exists a common section by
\begin{center}\begin{tikzcd}
	X \\
	& {X\times_YX} & X \\
	& X & Y
	\arrow[dashed, from=1-1, to=2-2]
	\arrow["1", from=1-1, to=2-3]
	\arrow["1"', from=1-1, to=3-2]
	\arrow[from=2-2, to=2-3]
	\arrow[from=2-2, to=3-2]
	\arrow["f"', from=2-3, to=3-3]
	\arrow["f", from=3-2, to=3-3]
\end{tikzcd}\end{center}
and this section doesn't affect the colimit.
In particular, it is the colimit of a sifted diagram.
So one could also call an effective epimorphism a reflexive epimorphism.
We avoid this terminology and stick to effective.

The same holds for effective monomorphism.
\end{remark}

\begin{example}\label{example:set-epi-reflective}
    In \(\Set\) every epimorphism is effective.
    For this consider any surjective map \(f\colon X\to Y\).
    With the axiom of choice, every set is projective so we obtain a section
    \(s\colon Y\to X\) such that \(fs = 1_Y\).
    Now let \(g\colon X\to Z\) be any map (with \(Z\) any set) such that \(gp_1 = gp_2\).
    Define \(h = g\circ s\).
    Since \(f\circ s\circ f = f\) we obtain that for any \(x\in X\) the tuple
    \(t = (x,s(f(x)))\) is an element of \(X\times_Y X\) and therefore
    \[
        g(x) = g(p_1(t)) = g(p_2(t)) = g(s(f(x))) = h(f(x))
    \]
    which means that \(g = h\circ f\).
    This implies the commutativity of the diagram
    \begin{center}
        \begin{tikzcd}
	       {X\times_Y X} && X && Y \\
	       \\
	       &&&& Z
	       \arrow["{p_1}", shift left, from=1-1, to=1-3]
	       \arrow["{p_2}"', shift right, from=1-1, to=1-3]
	       \arrow["f", two heads, from=1-3, to=1-5]
	       \arrow["g"', from=1-3, to=3-5]
	         \arrow["h", from=1-5, to=3-5]
        \end{tikzcd}.
    \end{center}
    Uniqueness of \(h\) follows by surjectivity of \(f\).
    Thus we have shown that \(f\) is the coequalizer of its kernel pair.
    Similarly, every monomorphism is effective in $\Set$.
\end{example}

\begin{lemma}[Characterisation of special mono-/epimorphisms]\label{lem:characterisation_of_special_monoepi}\uses{def:special_monoepi}
In a category with pullbacks,
the notions of effective and regular epimorphism agree.

Dually, in a category with pushouts,
the notions of effective and regular monomorphism agree.
\end{lemma}

\begin{proof}
Let $f\colon X\to Y$ be regular epimorphic,
so it is the coequalizer of a diagram
\begin{center}\begin{tikzcd}
	K & X.
	\arrow[shift left, from=1-1, to=1-2]
	\arrow[shift right, from=1-1, to=1-2]
\end{tikzcd}\end{center}
If the pullback $X\times_Y X$ exists, clearly there is a unique morphism $k\colon K\to X\times_Y X$ translating the pairs of arrows.
But now any $X\to Z$ equalizing both morphisms $p_1,p_2\colon X\times_Y X\rightrightarrows X$
clearly equalizes $p_1k$ and $p_2k$, and hence both arrows $K\rightrightarrows X$.
Thus there exists a factorization through $f$, and thereby $f$ is the coequalizer of its kernel pair.
If the epi is effective, then it clearly is regular.
\end{proof}

\begin{lemma}\label{lem:square_with_diag}
Consider a commutative square
\begin{center}\begin{tikzcd}
	X & Y \\
	Z & W.
	\arrow[two heads, from=1-1, to=1-2]
	\arrow[from=1-1, to=2-1]
	\arrow[from=1-2, to=2-2]
	\arrow[hook, from=2-1, to=2-2]
\end{tikzcd}
\end{center}
If the monomorphism or the epimorphism is regular,
then there exists a diagonal arrow $Y\to Z$ making the diagram
\begin{center}
\begin{tikzcd}
	X & Y \\
	Z & W
	\arrow[two heads, from=1-1, to=1-2]
	\arrow[from=1-1, to=2-1]
	\arrow[dashed, from=1-2, to=2-1]
	\arrow[from=1-2, to=2-2]
	\arrow[hook, from=2-1, to=2-2]
\end{tikzcd}
\end{center}
commutative.
\end{lemma}

\begin{proof}
Assume that $Z\to W$ is regular, i.e. the equalizer of two arrows $W\to H$.
\begin{center}\begin{tikzcd}
	X & Y \\
	Z & W & H
	\arrow[two heads, from=1-1, to=1-2]
	\arrow[from=1-1, to=2-1]
	\arrow[from=1-2, to=2-2]
	\arrow[hook, from=2-1, to=2-2]
	\arrow[shift left, from=2-2, to=2-3]
	\arrow[shift right, from=2-2, to=2-3]
\end{tikzcd}\end{center}
Now clearly, the arrow $X\to W$ equalizes the arrows $W\to H$ as it factors through $Z$.
Since $X\to Y$ is epic, this implies that $Y\to W$ equalizes both arrows $W\to H$, and now the universal property of the equalizer implies the existence of the arrow $Y\to Z$ such that the lower triangle commutes.
That the upper triangle commutes now follows since composed with the monic $Z\to W$ both arrows coincide.

A dual argumentation yields the result if $X\to Y$ is regular.
\end{proof}

\begin{definition}\label{def:projectivity}\uses{def:special_objects, def:special_monoepi}
An object $P$ is called \idx{projective} if for every regular epimorphism $X\to Y$
and every morphism $P\to Y$ there exists a lift $P\to X$,
\begin{center}
\begin{tikzcd}
	& W \\
	& X \\
	P & Y.
	\arrow[shift right, from=1-2, to=2-2]
	\arrow[shift left, from=1-2, to=2-2]
	\arrow[two heads, from=2-2, to=3-2]
	\arrow[dashed, from=3-1, to=2-2]
	\arrow[from=3-1, to=3-2]
\end{tikzcd}
\end{center}
\end{definition}

\begin{warning}
Phrased differently, the $\hom$-functor $\hom(P,-)$ preserves regular epimorphisms.
Often a projective object is referred to as an object that has the lifting property for \emph{all} epimorphisms.
However, this is not the correct definition as we show in the next remark~\ref{rem:injectivity}.
In practice, we only work in categories where every epimorphism is regular
and then the notions coincide.
\end{warning}

\begin{example}
In the category of compact Hausdorff spaces and its full subcategory of Stone spaces
the projective objects are precisely the extremally disconnected spaces.
See \cite{Gleason1958} and lemma~\Ref{lem:CHaus_enough_proj}.
\end{example}

\begin{warning}
Recall that an object $P$ is called \emph{effective projective}
if $\hom(P,-)$ preserves reflexive coequalizers, see definition~\ref{def:special_objects}.
Every effective projective object is projective
since reflexive coequalizers yield regular epimorphisms.
However, the converse does not hold.
For example, in the category $\Set$, every object is projective,
but the effective projectives are precisely the finite sets, see remark 3.4 in~\cite{Pedicchio2000}.

However, when one adds the assumption of $P$ being compact
(as in definition~\ref{def:special_objects})
then in many categories (such as so-called \emph{exact} categories)
projectivity implies effective projectivity.
This means that there the compact projective objects are precisely the objects
that are compact and projective. See lemma~\ref{lem:comp-proj} for more.
\end{warning}

\begin{lemma}\label{lem:comp-proj}
Let $\mcC$ be an \idx{exact category},
i.e. $\mcC$ is finitely complete,
has coequalizers of kernel pairs,
regular epimorphisms are stable under pullback
and every equivalence relation is effective (see definition~\ref{def:equivalence_relation}).
Then the following are equivalent for an object $K$ of $\mcC$.
\begin{enumerate}[(a)]
    \item The object $K$ is compact and projective.
    \item The object $K$ is compact projective.
\end{enumerate}
\end{lemma}

\begin{proof}
This follows form theorem 18.1 in \cite{Adamek2010a},
see also corollary 18.3.
\end{proof}

\begin{example}\label{ex:CHaus_exact}
	The category $\CHaus$ of compact Hausdorff spaces is exact.
\end{example}

\begin{remark}[Injective object]\label{rem:injectivity}
Dual to the definition of projective objects one defines an \idx{injective object}
to be an object $I$ such that for every regular monomorphism $X\to Y$
every morphism $X\to I$ extends to a morphism $Y\to I$ as in the diagram
\begin{center}
\begin{tikzcd}
	& W \\
	& Y \\
	I & X
	\arrow[shift left, from=2-2, to=1-2]
	\arrow[shift right, from=2-2, to=1-2]
	\arrow[dashed, from=2-2, to=3-1]
	\arrow[hook, from=3-2, to=2-2]
	\arrow[from=3-2, to=3-1]
\end{tikzcd}
\end{center}
Equivalently, $\hom(-,I)$ preserves regular epimorphisms (as a functor from $\mcC^\mathrm{op}$ to $\Set$).
\end{remark}

\begin{example}
The Hahn-Banach-Schaefer theorem is the statement that
the Banach space $C(S)$ for extremally disconnected compact Hausdorff $S$
is an injective object in the category of Banach spaces.
This would not be true anymore if we would consider injectivity with respect to
\emph{all} monomorphisms, which are precisely the injective maps
in the category of Banach spaces.
Then the embedding $\ell^1\to c_0$ is monic, but not every continuous functional
on $\ell^1$ extends to one on $c_0$.
Restricting to regular monomorphisms yields only closed embeddings
where then the classical result applies.

In more generality, these spaces are injective with respect to closed embeddings in complete locally convex vector spaces
(we are confident that it is even injective, however, we were unable to find a precise reference for this, see \ref{conj:M1_vs_lctvs} for more).

The Tietze extension theorem asserts that the unit interval $[0,1]$
is injective in the category of normal spaces.
Here, again it is necessary to work with regular monomorphisms.
\end{example}

\begin{lemma}[Stability of projective objects]\label{lem:stability_of_proj}\uses{def:projectivity, def:enough_projectives}
The class of projective objects is closed under retracts and coproducts.
\end{lemma}

\begin{proof}
See \cite[Proposition 4.6.3 and 4.6.4]{Borceux2008},
but note the slightly different definition of projective there.
Nevertheless, the proof is the same.
\end{proof}

In many categories, projective objects are precisely those where every \enquote{surjection onto them} admits a section.

\begin{lemma}[Projectivity via sections]\label{lem:projectivity_via_sections}\uses{def:special_objects}
    Let $\mcC$ be a category with pullbacks
    and with the property that pullbacks of regular epimorphisms are again regular epimorphisms.
    Then the following are equivalent.
    \begin{enumerate}[(a)]
        \item The object $P$ is projective.
        \item Every regular epimorphism $X\twoheadrightarrow P$ admits a section.
    \end{enumerate}
\end{lemma}

\begin{proof}
    If $P$ is projective and $f\colon X\to P$ a regular epimorphism,
    then lifting the identity $P\to P$ along $f$ yields the desired section $P\to X$.
    
	If (b) holds, let \(f\colon X\to Y\) be regular epic and \(g\colon P\to Y\).
	Then in the pullback diagram
\begin{center}
	\begin{tikzcd}
		X\times_Y P && P \\
		\\
		X && Y
		\arrow["p", from=1-1, to=1-3]
		\arrow["q", from=1-1, to=3-1]
		\arrow["g", from=1-3, to=3-3]
		\arrow["f", two heads, from=3-1, to=3-3]
	\end{tikzcd}
\end{center}
the morphism $p$ is regular epic.
Therefore there exists a section $s\colon P\to X\times_Y P$ of $p$.
Then $qs$ is the desired lift of $g$ along $f$ because $fqs = gps = g$.
Hence $P$ is projective.
\end{proof}

\begin{example}
This is fulfilled in the category $\CHaus$ of compact Hausdorff spaces or in the category of Stone spaces $\stone$.
\end{example}

\begin{definition}[Enough projectives]\label{def:enough_projectives}\uses{def:special_objects,def:cogenerating_set}
A category is said to have \idx{enough projectives},
if for any object $C$ there exists a projective object $P$
and a regular epimorphism $P\to C$.
\end{definition}

Often mathematical objects can be constructed out of basic building blocks.
We make this notion precise for categories.

\begin{definition}[Generators and separators]\label{def:cogenerating_set}\uses{def:limit, def:category}
Let $\mcC$ be a category and $\mcA$ a subclass of objects of $\mcC$.
	\begin{enumerate}[(i)]
		\item The class $\mcA$ is called \idx{separating}
		if for any two arrows $f,g\colon B\to C$ in $\mcC$ with $f\ne g$
		there exists an object $A\in \mcA$ and a morphism
		$h\colon A\to B$ with $fh\ne gh$.
		\item Dually, $\mcA$ is called \idx{coseparating}
		if for any two arrows $f,g\colon B\to C$ in $\mcC$ with $f\ne g$
		there exists $A\in \mcA$ and a morphism $h\colon C\to A$ with $hf\ne hg$.
		\item We call $\mcA$ \idx{densely generating} (or \idx{dense})
		if every object $C\in\mcC$ is a (small) colimit of a diagram using only objects in $\mcA$.
		\item Dually, $\mcA$ is \idx{densely cogenerating} (or \idx{codense})
		if every object $C\in\mcC$ is a (small) limit of a diagram using only objects in $\mcA$.
	\end{enumerate}
\end{definition}

\begin{remark}\label{rem:dense_is_separating}
If coproducts in $\mcC$ exist, then if $\mcA$ is dense, it is separating.
For this let $f,g\colon B\to C$ with $f\neq g$.
By the general description of colimits by coequalizers and coproducts,
the object $B$ admits a surjection
\[
    h\colon\coprod_{i\in I} A_i \to B
\]
with $A_i\in\mcA$.
Therefore $fh \neq gh$.
In particular there exists $A_i$ such that
$fhe_i \neq ghe_i$ for $e_i\colon A_i\to\coprod_{i\in I} A_i$ the embedding.
The dual holds for codense classes.

Note that many authors, e.g.~\cite{mac2013categories},
call separating sets generating (and coseparating cogenerating, respectively).
We follow \cite[A1.2.4, B2.4.1]{Johnstone}.
\end{remark}

\begin{lemma}[Charactersation of separating sets]\label{lem:cosep_coprod}\uses{def:cogenerating_set,def:limit,def:limit}
    Let $\mcC$ be a category with (small) coproducts
    and $\mcA$ a \emph{set} of objects of $\mcC$.
    Then the following are equivalent.
    \begin{enumerate}[(a)]
        \item The set $\mcA$ is separating.
        \item For every object $X$ there exists a family $(A_i)_{i\in I}$
        of objects $A_i\in\mcA$ and an epimorphism
            \[
                \coprod_{i\in I} A_i \to X.
            \]
        \item For every object $X$ the canonical morphism
            \[
                \coprod_{\substack{{f\colon A\to X} \\ {A\in\mcA}}} A\to X
            \]
            is epic.
    \end{enumerate}
\end{lemma}

\begin{proof}
This is evident.
\end{proof}

In many cases, the separating family consists of projective objects.
Then one already has a complete characterisation of all projective objects of the category.

\begin{lemma}\label{lem:enough-projectives-separating}
Let $\mcC$ be a category with coproducts, where every epimorphism is regular,
and let $\mcA$ be a (small) separating family of projectives.
Then $\mcC$ has enough projectives
and every projective object is a retract of coproducts of objects in $\mcA$.
\end{lemma}

\begin{proof}
See \cite[Proposition 4.6.5]{Borceux2008}.
\end{proof}

In remark~\ref{rem:dense_is_separating}
we saw that every dense set is separating in categories with coproducts.
What about the converse:
are separating sets already dense?
Say that $\mcC$ is a category in which every epimorphism is regular
and $\mcA$ be a separating set of objects.
For $X\in\mcC$ we get a surjection
\[
    \coprod_{i\in I} A_i \to X
\]
by lemma \ref{lem:cosep_coprod} since this epimorphism is regular and post-composition with epimorphisms does not change coequalizers, see \ref{lem:eq_comp_monic},
we have that the diagram
\begin{center}
\begin{tikzcd}
	{\coprod\limits_{j\in J}A_j} & Y & {\coprod\limits_{i\in I}A_i} & X
	\arrow[two heads, from=1-1, to=1-2]
	\arrow["f", shift left, from=1-2, to=1-3]
	\arrow["g"', shift right, from=1-2, to=1-3]
	\arrow[two heads, from=1-3, to=1-4]
\end{tikzcd}
\end{center}
is a coequalizer diagram.
However, one does not a priori now if the morphisms between the coproducts
arise from morphisms $A_j\to A_i$ between the individual objects.
But if we add the assumption that coproducts are stable,
we obtain the following.

\begin{lemma}\label{lem:cosep_is_cogen}
Let $\mcC$ be a category in which every epimorphism is regular,
pullbacks and (small) coproducts exist and coproducts are stable under pullback.
Then the following are equivalent for a \emph{set} $\mcA$ of objects of $\mcC$.
\begin{enumerate}[(a)]
    \item The family $\mcA$ is separating.
    \item The family $\mcA$ is densely generating.
\end{enumerate}
\end{lemma}

\begin{proof}
See  \cite[Proposition 4.5.6]{Borceux2008}
\end{proof}

\begin{lemma}\label{lem:colimit-injective-transition-maps}
Let $\mcC$ be a category such that the class $\mcC^\mathrm{c}$ of compact objects is separating
and let $F\colon \mcI\to\mcC$ be a filtered diagram with colimit $C$.
Then the following are equivalent.
\begin{enumerate}[(a)]
    \item Every embedding $\tau_i\colon F(i)\to C$ is a monomorphism.
    \item For every $i\to j$ in $\mcI$, the transition map $F(i)\to F(j)$ is monic.
\end{enumerate}
\end{lemma}

\begin{proof}
If every embedding $\tau_i\colon F(i)\to C$ is a monomorphism,
then clearly the transition map $f_i^j\colon F(i)\to F(j)$ is monic
since $\tau_j\circ f_i^j = \tau_i$.

Now assume that every transition map $f_i^j$ is a monomorphism.
We have to show that $\tau_i$ is monic.
Since the category $\mcC$ admits a separating class of compact objects,
it suffices to show that for $g,h\colon K\to F(i)$ with $K$ compact,
$\tau_ig = \tau_ih$ implies $g = h$.
Because $K$ is compact, we have
\[
    \hom(K,C) = \varinjlim_{i\in\mcI}\hom(K,F(i)).
\]
Now $g,h\in\hom(K,F(i))$ are equivalent on the right hand side,
since they are the same on the left hand side.
By the explicit description of colimits in $\Set$,
one obtains that there exists a $j\in\mcI$ and morphisms $a,b\colon i\to j$
such that $F(a)g = F(b)h$.
By the filteredness condition, there exists $c\colon j\to k$ such that
$ca = cb$ and thus $F(ca)g = F(cb)h$.
Since $F(ca) = F(cb)$ is monic, we obtain $g = h$.
\end{proof}

\subsection[(Co)completions]{Completions and cocompletions}\label{ssec:completions}
Some categories do not yet have all the (co)limits that we would like to consider.
(Co)completions are a formal way to freely adjoin some class of (co)limits to a category, see~\cite{kashiwara2005categories} for many theorems about (mainly $\mathrm{Ind}$-) completions.
We refer to 5.3.6 in \cite{Lurie2009} for an $\infty$-categorical version.
Our main source, however, is~\cite{Velebil2001}.

\begin{definition}[(Co)completion]\label{def:completion}\uses{def:category, def:functor, def:limit}
Consider any class of forms $\mcW$ that we want to adjoin and another class $\mcM$ of forms that we want to preserve.
	The $\mcM$-continuous $\mcW$-completion of any category $\mcC$ is a category $\mcD$ that possesses all limits of forms in $\mcW$,
    together with a fully faithful embedding\footnote{
      It is,
      in fact,
      irrelevant whether or not one simply asks for a functor or a fully faithful embedding.
      If the completion exists,
      this functor will always be fully faithful.
    }
    $\mcC\to\mcD$ preserving all limits of form in $\mcM$,
    such that for all other $\mcM$-continuous functors $\mcC\to \mcE$ into $\mcW$-complete categories,
    there exists a $\mcW$-continuous functor, unique up to natural isomorphy, $\mcD\to\mcE$ such that
\begin{center}
\begin{tikzcd}
	& \mcE \\
	\mcC & \mcD
	\arrow[""{name=0, anchor=center, inner sep=0}, from=2-1, to=1-2]
	\arrow[from=2-1, to=2-2]
	\arrow[from=2-2, to=1-2]
\end{tikzcd}\end{center}
commutes.
Dually,
$\mcM$-cocontinuous $\mcW$-cocompletions are defined.

	Special cases of this definition include \idx{free} $\mcW$-(co)completions, meaning that $\mcM=\emptyset$:
	\begin{itemize}
			\item The \idx{free (co)completion} of $\mcC$, using $\mcW$ being all small categories and $\mcM=\emptyset$.
			\item The \idx{free $\mathrm{Ind}$-completion} $\mathrm{Ind}(\mcC)$ of $\mcC$ is the cocompletion of $\mcC$ with $\mcW$ being all (small) filtered categories and $\mcM=\emptyset$.
			\item The \idx{free $\mathrm{Pro}$-completion} $\mathrm{Pro}(\mcC)$ of $\mcC$ is the completion with $\mcW$ being all (small) cofiltered categories and $\mcM=\emptyset$.
			\item The \idx{free $\mathrm{sInd}$-completion (sifted completion)} $\mathrm{sInd}(\mcC)$ of $\mcC$  is the cocompletion with $\mcW$ being all (small) sifted categories and $\mcM=\emptyset$.
	\end{itemize}
\end{definition}

For $\mcM\subseteq\mcW$,
this can be thought of as a special case of a 2-adjunction when considering the forgetful functor from $\mcW$-complete categories to $\mcM$-complete categories
with $\mcW$- resp.\ $\mcM$-continuous functors.
(Of course,
this only makes sense after specialization to the case that $\mcC$ is $\mcM$-complete.)

As example of the appearance of cocompletions we consider the process of 1-animation (see e.g. Chapter 4 in \cite{Adamek2010a} for more),
which is a preliminary version of animation, and is explained by the preceeding discussion.

\begin{definition}[1-animation]\label{def:one_ani}
Let $\mcC$ be small cocomplete, and denote by $\mcC^{\mathrm{cp}}$ the full subcategory of $\mcC$ spanned by the compact projective objects.

If $\mcC$ is generated by $\mcC^{\mathrm{cp}}$, then
\[\sInd(\mcC^{\mathrm{cp}})\simeq \mcC.\]
If $\mcC^{\mathrm{cp}}$ is small, then
\[\sInd(\mcC^{\mathrm{cp}})\]
can be described as the category of functors ${\mcC^{\mathrm{cp}}}^{\op}\to \Set$ mapping finite coproducts in $\mcC^{\mathrm{cp}}$ to finite products.
\end{definition}

One can always recover completions of cocompletions (and vice versa).
\begin{lemma}[Completions via cocompletions]\label{lem:comp_via_cocomp}\uses{def:completion}
Consider any category $\mcC$ as well as classes $\mcM$ and $\mcW$.
Then the $\mcM$-continuous $\mcW$-completion of $\mcC$ can be calculated
by taking the opposite category of the $\mcM^{\mathrm{op}}$-cocontinuous $\mcW^{\mathrm{op}}$-cocompletion of $\mcC^\mathrm{op}$.
(By $\mcM^{\op}$,
we mean $\{\mcI^{\op}\mid \mcI\in\mcM\}$ and similarly for $\mcW$.)

In particular $\mathrm{Pro}(\mcC)=\mathrm{Ind}(\mcC^{\mathrm{op}})^\mathrm{op}$.
\end{lemma}
\begin{proof}
	Note that a diagram of form $\mcI$ in $\mcC$ with limit $X$ is equivalently a diagram $I^{\op}\to \mcC^{\op}$ with colimit $X$.
	Hence a $\mcW$-complete category $\mcD$ is equivalently a $\mcW^{\op}$-cocomplete category $\mcD^{\op}$,
    a $\mcM$-continuous functor $\mcC\to \mcD$ is equivalently a $\mcM^{\op}$-cocontinuous functor $\mcC^{\op}\to \mcD^{\op}$
	and thus, denoting by $\mcE$ the $\mcM$-continuous $\mcW$-completion of $\mcC$, $\mcE^{\op}$ yields a $\mcM^{\op}$-cocontinuous $\mcW^{\op}$-cocompletion of $\mcC^\op$.
\end{proof}

\begin{definition}[Closure under $\mcW$-limits]
	Let $\mcC$ be any category admitting all $\mcW$-colimits, and $\mcA\sub \mcC$ a class ob objects.
	The $\mcW$-closure of $\mcA$ in $\mcC$ is defined as the smallest full subcategory $\mcB$ of $\mcC$ containing all colimits
    of forms in $\mcW$ and containing $\mcA$ (and closed under isomorphisms).
\end{definition}

\begin{lemma}[Constructing free completions]\label{lem:completions_via_Yoneda}\uses{def:completion,lem:yoneda, lem:limits_in_functor_categories, lem:hom_functor_preserves_limits}
 For any category $\mcC$ define the category of \idx{$\mcW$-accessible Presheaves} of $\mcC$ are defined to be the closure with respect to forms in $\mcW$ of the subcategory of representable presheaves,

and dually define the accessible functors in $\mathrm{Fun}(\mcC, \Set)$ as defined as the closure with respect to limits of the forms in $\mcW$ of the image of the covariant Yoneda embedding.

If $\mcW$ is the class of all small categories, we simply call these functors \idx{accessible}.
Note that, due to the yoneda lemma, the category of accessible presheaves in $\PSh(\mcC)$ is locally small whenever $\mcC$ is.
If $\mcC$ is small, all presheaves are accessible.
Often, accessible presheaves form a superior replacement to the large category of presheaves.

Define the category $\mcD$ to consist of all $\mcW$-accessible functors $\mathrm{Fun}(\mcC,\Set)$.
Then the Yoneda embedding $\mcC\to\mcD$ incorporates $\mcD$ as a free $\mcW$-completion of $\mcC$.
\end{lemma}

Recall, that the Yoneda embedding $\mcC\to\PSh(\mcC)$ is continuous but does in general not commute with colimits, as we allow too many presheaves (and hence allow to compute colimits pointwise).
However, individually $y(X)=\hom(-,X)$ maps all colimits to limits, and hence forms a presheaf of the nicest possible form.

This leads to the following more general theorem.
\begin{theorem}[General cocompletions of small categories]\label{thm:constructing_cocompletions}\uses{def:completion}
If $\mcC$ is any small category, $\mcM$ and $\mcW$ are classes of forms, then define $\mcD$ to be the full subcategory of all $\mcM$-cocontinuous presheaves (functors $\mcC^{\mathrm{op}}\to\Set$ that preserve all limits in $\mcC^{\mathrm{op}}$ of diagrams indexed by $\mcI^{\op}$ for forms $\mcI\in\mcM$).
As all representable presheaves are of this form, we can define $\mcE$ as the closure of representable presheaves in $\mcD$ under colimits of diagrams of representable presheaves with forms in $\mcW$.
Now, the restriction of the Yoneda embedding exhibits $\mcE$ as the $\mcM$-continuous $\mcW$-completion of $\mcC$.
\end{theorem}
\begin{proof}
	See Thm 14 of \cite{Velebil2001}.
	See also chapter 5.7 in \cite{Kelly1982}.
\end{proof}
If $\mcC$ is not small, then one can pass to the presheaf category with values in a larger universe, do the same construction there, and obtains a
(almost always) large category with the same property.
See \cite{Velebil2001} and \cite{Kelly1982} for more.

\begin{remark}
  Often,
  we are in the convenient situation
  that to obtain the $\mcW$-closure we only need to \enquote{add $\mcW$-limits} once,
  i.e.,
  $\mcW$-limits of $\mcW$-limits of objects in our subcategory can be written simply as $\mcW$-limits
  of objects in this subcategory.
  Specifically,
  in the cases relevant to the construction,
  this is the case for
  \begin{itemize}
    \item $\mcW$ the class of all small categories,
    \item $\mcW$ the class of small filtered categories (see \cite[chapter 6]{kashiwara2005categories}),
    \item $\mcW$ the class of small sifted categories (see \cite[section 2]{Adamek2001}).
  \end{itemize}
  Of course,
  the dual cases work the same way. See also the discussion before section 5.8 in \cite{Kelly1982}.

	Then another way to define a free $\mcW$-completion is to define a category $\mcE$ whose objects are diagrams $D\colon \mcI\to\mcC$ in $\mcC$ of forms $\mcI$ in $\mcW$,
denoted by $"\varprojlim_i" D(i)$ and whose sets of morphisms are defined as
	\[\hom("\varprojlim_i" D(i), "\varprojlim_j" E(j))=\varprojlim_j\varinjlim_i \hom(D(i),E(j)).\]

	The corresponding way to define a free $\mcW$-cocompletion is to define a category $\mcE$ whose objects are diagrams $D\colon \mcI\to\mcC$ in $\mcC$ of forms $\mcI$ in $\mcW$,
denoted by $"\varprojlim_i" D(i)$ and whose sets of morphisms are defined as
	\[\hom("\varprojlim_i" D(i), "\varprojlim_j" E(j))=\varprojlim_i\varinjlim_j \hom(D(i),E(j))\]
\end{remark}

\begin{remark}
	Note that the condition of representing iterated colimits as single colimits is a special feature of our concrete setup.
	As example, consider the full subcategory $\{\ast\}\sub \CHaus$.
	Then all colimits that can be achieved with $\{\ast\}$ are coproducts; and all coproducts $\coprod_{i\in I}\ast$ are described as $\beta I$.
	However, every compact Hausdorff space is a colimit of spaces of the form $\beta I$, thus every compact Hausdorff space $K$ may be described as
	\[K=\varinjlim_{j}\varinjlim_{i\in I_j}\ast=\varinjlim_j \beta I_j.\]
	But clearly, not every compact Hausdorff space is given by a $\beta I$; showing that here, iterated colimits
	of a point may not always be written as single colimits of a point.
\end{remark}

Sometimes it will be nice to consider certain cocompletions of larger categories, that are given as (large) union of small categories.
\begin{lemma}[Cocompletions of unions]\label{lem:cocomp_union}\uses{def:completion}
Assume that $\mcC$ is an increasing union of full subcategories $\mcC_i$, indexed by some totally ordered class $\mcI$, such that
\begin{itemize}
	\item all inclusions $\mcC_i\to \mcC_j$ are $\mcM$-continuous,
	\item the $\mcC_i$ are small, and
	\item $\mathrm{cof}(I)>|J|$ for all $J\in \mcM$ or $J\in \mcW$ (see \ref{def:cofinality} for a definition of cofinality).
\end{itemize}
Then the $\mcM$-continuous $\mcW$-cocompletion of $\mcC$ exists, and can either be constructed by the construction for small categories, or by the union over all completions of the $\mcC_i$ along the extensions of the inclusions.
\end{lemma}
\begin{proof}
	We can construct the cocompletion of $\mcC$ as the union of the cocompletions of the $\mcC_i$.
	Denote by $\mcC_i\hookrightarrow \overline{\mcC_i}$ the $\mcM$-continuous $\mcW$-cocompletion of $\mcC_i$.

Then, by the universal property of the cocompletion, we obtain a unique $\mcW$-cocontinuous functor $\overline{\mcC_i}\to \overline{\mcC_j}$ for $i\le j$, with

\begin{center}\begin{tikzcd}
	{\mcC_i} && {\overline{\mcC_i}} \\
	{\mcC_j} && {\overline{\mcC_j}} \\
	\vdots && \vdots \\
	\\
	{\varinjlim \mcC_i=\mcC} && {\varinjlim\overline{\mcC_i}}.
	\arrow[hook, from=1-1, to=1-3]
	\arrow[hook', from=1-1, to=2-1]
	\arrow[dashed, from=1-3, to=2-3]
	\arrow[hook, from=2-1, to=2-3]
	\arrow[hook', from=2-1, to=3-1]
	\arrow[dashed, from=2-3, to=3-3]
	\arrow[hook', from=3-1, to=5-1]
	\arrow[from=3-3, to=5-3]
\end{tikzcd}\end{center}

	This in turn glues to a unique functor $\varinjlim \mcC_i\to \varinjlim \overline{\mcC_i}$.
	We claim that this is indeed an $\mcM$-continuous $\mcW$-cocompletion.
    For the $\mcM$-continuity note that for any $\mcJ\in \mcM$ and any diagram $\mcJ\to \varinjlim \mcC_i$ with colimit $L$, the diagram and its limit is contained in some level $\mcC_i$,
    as the cofinality of $I$ is bigger than $|\mcJ|$.
	Hence the diagram and colimit $L$ is contained in some $\mcC_i$, and the functor $\varinjlim \mcC_i\to \varinjlim \overline{\mcC_i}$ maps the
	limit of the diagram to the limit of the diagram in $\overline{\mcC_i}$.
	As the map $\mcC_i\to \overline{\mcC_j}$ is always $\mcM$-continuous, the same is true for all $\overline{\mcC_j}$.
    Hence the colimit of the diagram in $\overline{\mcC_i}$ is the colimit of the diagram in $\varinjlim \overline{\mcC_i}$.

	A similar reasoning shows, that $\varinjlim \overline{\mcC_i}$ is $\mcW$-complete.

    Now, consider any $\mcM$-continuous functor $\varinjlim \mcC_i\to \mcD$ for any $\mcW$-complete $\mcD$.
    Then, individually there exist unique $\mcW$-continuous functors $\overline{\mcC_i}\to \mcD$ that make the diagram

\begin{center}\begin{tikzcd}
	{\mcC_i} && {\overline{\mcC_i}} \\
	{\mcC_j} && {\overline{\mcC_j}} \\
	\vdots && \vdots \\
	& \mcD \\
	{\varinjlim \mcC_i=\mcC} && {\varinjlim\overline{\mcC_i}}
	\arrow[hook, from=1-1, to=1-3]
	\arrow[hook', from=1-1, to=2-1]
	\arrow[hook, from=1-3, to=2-3]
	\arrow[dashed, from=1-3, to=4-2]
	\arrow[hook, from=2-1, to=2-3]
	\arrow[hook', from=2-1, to=3-1]
	\arrow[hook, from=2-3, to=3-3]
	\arrow[dashed, from=2-3, to=4-2]
	\arrow[hook', from=3-1, to=5-1]
	\arrow[hook, from=3-3, to=5-3]
	\arrow[from=5-1, to=4-2]
	\arrow[hook, from=5-1, to=5-3]
	\arrow[dashed, from=5-3, to=4-2]
\end{tikzcd}\end{center}
	commute.
	These glue to a unique $\mcW$-continuous functor $\varinjlim \overline{\mcC_i}\to \mcD$.

	This shows, that the cocompletion exists and is locally small. \cite{Velebil2001} implies that in this case it is described by the usual way of presheaves on $\varinjlim_i \mcC_i$.
	\end{proof}

\begin{definition} For any class $\mcW$ define the \idx{orthogonal class} $\mcW^\perp$ of forms as all those small categories $I$,
	such that colimits of forms in $\mcW$ commute with limits of form $I$ in $\Set$.

	\end{definition}
See e.g. 10.4 in \cite{Hoermann2019} for more.
\begin{example}
	\begin{tabular}{c|c}
		$\mcW$ & $\mcW^\perp$\\
		\hline
		$\{\ast\}$ & all small cats\\
		sifted & finite discrete sets\\
		filtered & finite categories\\
		all small cats & $\{\ast\}$
		\end{tabular}
	\end{example}

\begin{lemma}
	Assume that $\mcW$ is stable under concatenation as explained above.
	The $\emptyset$-continuous $\mcW$-cocompletion of a small category agrees with the $\mcW^\perp$-cocontinuous $\mcW$-cocompletion.
	\end{lemma}
\begin{proof}
	Define $\mcM$ to be $\mcW^{\perp}$.
	We have to show that the embedding is $\mcM$-continuous.
	Hence take any form $I\in \mcW$ and $J\in \mcM$.
	Now, take any colimit $T=\varinjlim_{i\in I}y(X_i)$ of $X_i\in \mcC$, as well as any $S=\varinjlim_j S_j$ in $\mcC$.
	We want to see that 
	\[\hom(\varinjlim S_j, \varinjlim X_i)=\varinjlim_i \hom(\varinjlim S_j,X_i).\]
	The left hand side by Yoneda is given by $T(S)$.	
	As the colimit is taken in presheaves, it is taken pointwise in sets,
	\[T(S)=\varinjlim y(X_i)(S)\]
By Yoneda,
	\[y(X_i)(\varinjlim_{J} S_j)=\hom(\varinjlim S_j, T_i)=\varprojlim T_i(S_j).\]
	We conclude
	\[T(S)=\varinjlim_i \varprojlim_j T_i(S_j)=\varprojlim_j \varinjlim_i T_i (S_j)=\varprojlim_j T(S_j)\]
	And thus $T$ maps $\mcM$-colimits in $\mcC$ to limits in $\Set$.

	\end{proof}
\begin{question}
	We believe that the stability of $\mcW$ under concatenation can be dropped, as long as the category of presheaves is obtained by transfinitely taking colimits of forms in $\mcW$ of objects in the last step, starting with representables.
	If we assume enough universes or $\mcW$ to be small this condition should hence become redundant.

	We are not sure whether the $\mcM$-cocontinuous $\mcW$-cocompletion can be dualised similarly.
\end{question}

This in particular implies the following classical theorems.
\begin{corollary}[Preserved colimits into completions]\label{lem:completion_preserves_limits}\uses{lem:completions_via_Yoneda}
The embedding $\mcC\to \mathrm{Ind}(\mcC)$ commutes with finite colimits, and $\mcC\to\mathrm{Pro}(\mcC)$ with finite limits.
The embedding $\mcC\to\mathrm{sInd}(\mcC)$ commutes with finite coproducts.
\end{corollary}

If now $\mcW^\perp$ and $\mcW$ together create all limits, the $\mcW$-completion in fact is cocomplete.
\begin{lemma}\label{lem:sind_via_cont}
	Let $\mcW$ again be stable under concatenation, and furthermore assume $\mcC$ to posess all colimits of forms in $\mcW^{\perp}$.
	If every colimit can be decomposed into colimits of form $\mcW$ and of form $\mcW^\perp$, then the free $\mcW$-cocompletion agrees with the $\mcW^\perp$-continuous cocompletion.
\end{lemma}
\begin{proof}
	Define $\mcM=\mcW^\perp$, and $\mcD$ as the $\mcM$-cocontinuous $\mcW$-cocompletion constructed as in the theorem.
We only need to show cocompleteness of $\mcD$, as in this case $\mcD$ would agree with the construction of the $\mcM$-continuous cocompletion.
	Hence consider any colimit $\varinjlim_h D_h$ in $\mcD$ of representables $D_h$.
	Decompose this colimit
	\[\varinjlim_h D_h=\varinjlim_{i\in \mcI}\varinjlim_{j\in \mcJ_i} D_j\]
	where $\mcI\in \mcW$ and $\mcJ_i\in \mcM$.

	Then the $\varinjlim_{j\in \mcJ_i} D_j$ exist, as the embedding is $\mcM$-cocontinuous and all $\mcM$-colimits in $\mcC$ exist,
	and now by $\mcW$-cocompleteness the outer colimit also exists, yielding the existence of $\varinjlim_h D_h$.

\end{proof}

\begin{remark}
	The condition of $\mcW^{\perp}$ and $\mcW$ together generate all limits/colimits is true for the classes of all small diagrams, of filtered diagrams and of sifted diagrams.
	We do not think that this will be true in general, but we are highly uncertain about this.
\end{remark}
This yields at least the last parts of the following well known results, see \cite{kashiwara2005categories} 6.1.17 and 6.1.18.
\begin{lemma}[Pro-completion is complete and admits finite coproducts]\label{lem:pro-completion_is_complete}\uses{def:completion}
If $\mcC$ admits (finite) limits, then $\mathrm{Ind}(\mcC)$ admits (finite) limits and $\mcC\to\mathrm{Ind}(\mcC)$ is continuous.
If $\mcC$ admits (finite) colimits, then $\mathrm{Pro}(\mcC)$ admits (finite) colimits and $\mcC\to\mathrm{Pro}(\mcC)$ is cocontinuous.

Furthermore, if $\mcC$ admits finite colimits, then $\mathrm{Ind}(\mcC)$ admits all small colimits and dually if $\mcC$ is finitely complete, then $\mathrm{Pro}(\mcC)$ is complete.
\end{lemma}

Sometimes, it will be useful to generalise the notion of $\mcM$-continuous cocompletions to $\mcM$ not being a class of forms, but rather a class of diagrams in $\mcC$.
We call this the \idx{non-uniform} cocompletion.
\begin{definition}\label{def:non-uniform-coco}
	Take a class of diagrams $\mcD$ in $\mcC$, and assume that all colimits of elements in $\mcD$ exist (this might be dropped).
	The non-uniform $\mcD$-continuous cocompletion of $\mcC$ is defined to be a cocomplete category $\mcE$ together with an embedding $\mcC\to\mcE$ preserving all the colimits of diagrams in $\mcD$,
	such that every other such functor $\mcC\to \mcE'$ factors uniquely cocontinuously through $\mcE$.
\end{definition}
\begin{remark}
In particular, an analogous definition of $\mcD$-continuous $\mcW$-cocompletions can be given
and the existence and description of the cocompletions also holds for the non-uniform $\mcD$-cocontinuous $\mcW$-cocompletion, if $\mcD\sub \mcW$ and everything is small.
This is essentially 6.23 in \cite{Kelly1982}.

Furthermore, a suitable replacement for the commutativity of cocompletions and unions also works for a non-uniform cocompletion, since we did not use the fact that the colimits in $\mcM$ contain all colimits of the same forms.
\end{remark}

\begin{question}
	How much of the above results hold for non-uniform cocompletions?
\end{question}

\begin{remark}
Developing further the notion of \enquote{non-uniform} cocompletions should lead to the general theory of \idx{sketches}, see e.g \cite{Kelly1982} or \cite{Johnstone} for more.
However, since we only need a very special case of this, we will not discuss this further.
\end{remark}


\section{Adjoint functors}\label{sec:adjoint-functors}

\quot{We did not then regard [category theory] as a field for further research efforts, but just as a language and an orientation --
  a limitation which we followed for a dozen years or so, till the advent of adjoint functors}{Saunders Mac Lane, \enquote{Concepts and categories in perspective}}
The next important concept in category theory is that of adjoint functors.
\begin{definition}[Three equivalent definitions of adjunctions]\label{def:adjunction}\uses{def:functor, def:natural_transformation}
    Consider two categories $\mcC$ and $\mcD$.
    An \idx{adjunction} between $\mcC$ and $\mcD$ consists of
    \begin{itemize}
        \item a functor $R\colon\mcD\to\mcC$, called the \idx{right adjoint},
        \item a functor $L\colon\mcC\to\mcD$, called the \idx{left adjoint},
        \item a natural isomorphism,
        \[\hom(L-,-)\simeq \hom(-,R-)\]
    \item a natural transformation $\eta\colon 1_{\mcC}\to RL$, called the \idx{unit},
    \item a natural transformation $\eps\colon LR\to 1_{\mcD}$, called the \idx{counit},
\end{itemize}
    so that the triangle identities

\begin{center}\begin{tikzcd}
	& LRL &&& RLR \\
	L & L && R & R
	\arrow["{\eps.L}", from=1-2, to=2-2]
	\arrow["{R.\eps}", from=1-5, to=2-5]
	\arrow["{L.\eta}", from=2-1, to=1-2]
	\arrow["{1_L}"', from=2-1, to=2-2]
	\arrow["{\eta.R}", from=2-4, to=1-5]
	\arrow["{1_R}"', from=2-4, to=2-5]
\end{tikzcd}\end{center}
hold.

Up to natural isomorphy, an adjunction can be reconstructed from only part of these data:
\begin{itemize}
    \item Given the right adjoint $R$, the left adjoint $L$ and the $\hom$-isomorphism $\phi\colon\hom(L-,-)\to \hom(-,R-)$,
    one defines the unit $(\eta_C)_{C\in\mcC}$ via $\eta_C=\phi(1_{LC})$ and the counit $(\eps_D)_{D\in\mcD}$ via $\eps_D=\phi^{-1}(1_{RD})$.
    \item Given $R,L$, the unit $(\eta_C)$ and counit $(\eps_D)$, fulfilling the triangle identities,
     the composition

\begin{center}\begin{tikzcd}
	{\hom(L-,-)} & {\hom(RL-,R-)} & {\hom(-,R-)}
	\arrow["R", from=1-1, to=1-2]
	\arrow["{\eta^*}", from=1-2, to=1-3]
\end{tikzcd}\end{center}
    yields a natural isomorphism with inverse induced by

\begin{center}\begin{tikzcd}
	{\hom(-,R-)} & {\hom(L-,RL-)} & {\hom(L-,-)}
	\arrow["L", from=1-1, to=1-2]
	\arrow["{\eps_*}", from=1-2, to=1-3]
\end{tikzcd}\end{center}
\item Assume that one has given the functor $R\colon\mcD\to\mcC$ such that for every object $C\in\mcC$ there exists an object $LC\in \mcD$  together with
    a morphism $\eta_C\colon C\to R(LC)$ such that for all morphisms $f\colon C\to RD$ for any $D\in\mcD$ there exists a unique morphism $\phi(f)\colon LC\to D$ with
\begin{center}\begin{tikzcd}
	C & {R(LC)} & LC \\
	& RD & D
	\arrow["{\eta_C}", from=1-1, to=1-2]
	\arrow["f"', from=1-1, to=2-2]
	\arrow["{R(\phi(f))}", from=1-2, to=2-2]
	\arrow["{\exists! \,\phi(f)}", dashed, from=1-3, to=2-3]
\end{tikzcd}.\end{center}
    Then the obvious choice of our notation defines the object function of the functor $L$, as well as the unit $\eta$ and the isomorphism $\phi$.
    The hom-functions of $L$ are now determined by

\begin{center}\begin{tikzcd}
	C & {R(LC)} & LC \\
	C' & RLC' & LC'.
	\arrow["{\eta_C}", from=1-1, to=1-2]
	\arrow["f"', from=1-1, to=2-1]
	\arrow["{\eta_{C'} f}"', shift left=2, from=1-1, to=2-2]
	\arrow["{R(\phi(f))}", from=1-2, to=2-2]
	\arrow["{\exists! \,\phi(\eta_{C'} f)=L(f)}", dashed, from=1-3, to=2-3]
	\arrow["{\eta_{C'}}"', from=2-1, to=2-2]
\end{tikzcd}\end{center}
\end{itemize}
For proofs, explanations and further characterisations,
see the article on adjoint functors in nLab.
\footnote{\url{https://ncatlab.org/nlab/show/adjoint+functor}}
Adjoint functors are denoted by $L\dashv R$, or
\begin{center}\begin{tikzcd}
	\mcC & \mcD
	\arrow[""{name=0, anchor=center, inner sep=0}, "L", curve={height=-12pt}, from=1-1, to=1-2]
	\arrow[""{name=1, anchor=center, inner sep=0}, "R", curve={height=-12pt}, from=1-2, to=1-1]
	\arrow["\dashv"{anchor=center, rotate=-90}, draw=none, from=0, to=1]
\end{tikzcd}.\end{center}
\end{definition}

\begin{remark}
  Note that this definition can be made in any 2-category, see adjoint in nLab.
  \footnote{\url{https://ncatlab.org/nlab/show/adjunction}}
\end{remark}

\begin{warning}[Algebraists, skip!]
    Again, the classical notion of adjoint operators in functional analysis can be seen as a special case of adjoint functors (see e.g. \cite{Baez1996, Palmquist1974}),
	but we strongly suggest to regard both notions as completely different and just notationally similar.
\end{warning}
\begin{lemma}[Uniqueness of adjoints]\label{lem:uniqueness_of_adjoints}\uses{def:adjunction}
	Any adjunction $(R,L,\phi,\eta,\eps)$ is, up to isomorphy, uniquely determined by the functor $R$ or, analogously, by the functor $L$.
\end{lemma}
\begin{proof}
  Clear, by the latter characterisation above.
\end{proof}
This explains, why often large parts of adjunctions are left implicit,
e.g. it is a reasonable question to ask for a fixed functor $R$ whether it is a right adjoint,
i.e. admits any left adjoint, which by this uniqueness result can justifiably be called \enquote{the} left adjoint of $R$.

One important example of adjunctions is the case of forgetful functors possessing a left adjoint.
\begin{definition}[Reflective subcategory]\label{def:reflective_subcategory}\uses{def:adjunction}
  A \idx{reflective subcategory} of a category $\mcC$ is a full subcategory $\mcD\subseteq \mcC$ such that the inclusion (forgetful functor) admits a left adjoint $L$.
\end{definition}

\begin{remark}
  A left adjoint to a forgetful functor to the category of sets (or sometimes also to another fixed base category modelling the category of sets one works with)
  is often referred to as a \idx{free functor}.
\footnote{Usually endings such as \enquote{-ification} or \enquote{-ization} are good hints to recognize left adjoints to forgetful functors.}
\end{remark}
\begin{example}
	\begin{enumerate}[(i)]
		\item All free functors from example \ref{ex:functors} are left adjoints to forgetful functors, i.e.:
		\begin{itemize}
			\item The \idx{free abelian group functor} $\Z[-]$ is the left adjoint to $\Ab\to\Set$.
			\item The \idx{free group functor} $F(-)$ is the left adjoint to $\Grp\to\Set$.
			\item The left adjoint to the forgetful functor from $\R$-vector spaces $\Mod_{\R}\to\Set$ is given by assigning any set $S$ the vector space $\R^S$.
		  \item The forgetful functor $\Top\to \Set$ admits the discrete functor $(-)_{\mathrm{disc}}$
				(installing the discrete topology on any set) as left adjoint,
				and the indiscrete topology functor $(-)_{\mathrm{triv}}$ as right adjoint.
		\end{itemize}
	  \item The Cauchy completion of a metric space is the left adjoint to the forgetful functor from complete metric spaces with uniformly continuous maps as morphisms
			to all metric spaces equipped with uniformly continuous maps.
		\item The Stone-\v{C}ech compactification is the left adjoint to the forgetful functor from compact Hausdorff spaces to all topological spaces $\CHaus\to\Top$.
		\item The abelianisation of a group is the left adjoint to the forgetful functor $\Grp\to\Ab$.
		\item The Hausdorffification of a topological space is the left adjoint to the forgetful functor from Hausdorff spaces to all topological spaces.
	  \item Although this is not a precise statement,
		the definition of completions is motivated by defining a (2-categorical) left adjoint
		to the forgetful functor from the 2-metacategory of $\mcW$-complete categories with $\mcW$-continuous functors to $\mcM$-complete categories with $\mcM$-continuous morphisms.
	  \item A wealth of adjoints (especially simultaneous left and right adjoints) can be found in representation theory.
	\end{enumerate}
	The reader is strongly encouraged to check that the corresponding \enquote{extension results} in fact are nothing but the universal property of adjoints.
\end{example}
\quot{Adjoint functors arise everywhere.}
{Mac Lane, \cite{mac2013categories}}

By composing the natural $\hom$-set isomorphisms one obtains the following result.
\begin{lemma}[Composition of Adjunctions]\label{lem:comp_of_adjoints}\uses{def:adjunction}
Consider any two adjunctions

\begin{center}\begin{tikzcd}
	\mcC & \mcD & \mcE
	\arrow[""{name=0, anchor=center, inner sep=0}, "{L_1}", curve={height=-12pt}, from=1-1, to=1-2]
	\arrow[""{name=1, anchor=center, inner sep=0}, "{R_1}", curve={height=-12pt}, from=1-2, to=1-1]
	\arrow[""{name=2, anchor=center, inner sep=0}, "{L_2}", curve={height=-12pt}, from=1-2, to=1-3]
	\arrow[""{name=3, anchor=center, inner sep=0}, "{R_2}", curve={height=-12pt}, from=1-3, to=1-2]
	\arrow["\dashv"{anchor=center, rotate=-90}, draw=none, from=0, to=1]
	\arrow["\dashv"{anchor=center, rotate=-90}, draw=none, from=2, to=3]
\end{tikzcd}\end{center}

Then, by composing the natural $\hom$-isomorphisms, one obtains a $\hom$-isomorphism
\[\hom(L_2L_1-,-)\simeq \hom(L_1-,R_2-)\simeq \hom(-,R_1R_2-)\]
And thus one obtains that the composed functors again form an adjunction.
\begin{center}\begin{tikzcd}
	\mcC & \mcE
	\arrow[""{name=0, anchor=center, inner sep=0}, "{L_2L_1}", curve={height=-12pt}, from=1-1, to=1-2]
	\arrow[""{name=1, anchor=center, inner sep=0}, "{R_1R_2}", curve={height=-12pt}, from=1-2, to=1-1]
	\arrow["\dashv"{anchor=center, rotate=-90}, draw=none, from=0, to=1]
\end{tikzcd}\end{center}

Clearly $\eta^{L_2L_1}_{C}=\phi^{12} (1_{L_2L_1 C})=\phi^{1}(\phi^2(1_{L_2L_1C}))$.
Alternatively,
the counit and unit are given by

\begin{center}\begin{tikzcd}
	C & {R_1L_1 C} & {R_1R_2L_2L_1 C} && {L_2L_1R_1R_2D} & {L_2R_2D} & D
	\arrow["{\eta^1_C}", from=1-1, to=1-2]
	\arrow["{R_1\eta^2_{L_1C}}", from=1-2, to=1-3]
	\arrow["{L_2\eps^1_{R_2D}}", from=1-5, to=1-6]
	\arrow["{\eps^2_D}", from=1-6, to=1-7]
\end{tikzcd}\end{center}
\end{lemma}

The following preservation property is 4.4.6 in \cite{Riehl}.

\begin{lemma}[Adjunctions in functor categories]\label{lem:adj_in_functor_cats}\uses{lem:limits_in_functor_categories}
For any adjunction

\begin{center}\begin{tikzcd}
	\mcC & \mcD
	\arrow[""{name=0, anchor=center, inner sep=0}, "L", curve={height=-12pt}, from=1-1, to=1-2]
	\arrow[""{name=1, anchor=center, inner sep=0}, "R", curve={height=-12pt}, from=1-2, to=1-1]
	\arrow["\dashv"{anchor=center, rotate=-90}, draw=none, from=0, to=1]
\end{tikzcd}\end{center}
	any small category $\mcJ$ and any category $\mcE$, the induced \idx{pullback functors} and the \idx{pushforward functors} induce adjunctions

\begin{center}\begin{tikzcd}
	{\mcC^{\mcJ}} & {\mcD^{\mcJ}} && {\mcE^{\mcC}} & {\mcE^{\mcD}}
	\arrow[""{name=0, anchor=center, inner sep=0}, "{L_*}", curve={height=-12pt}, from=1-1, to=1-2]
	\arrow[""{name=1, anchor=center, inner sep=0}, "{R_*}", curve={height=-12pt}, from=1-2, to=1-1]
	\arrow[""{name=2, anchor=center, inner sep=0}, "{R^*}", curve={height=-12pt}, from=1-4, to=1-5]
	\arrow[""{name=3, anchor=center, inner sep=0}, "{L^*}", curve={height=-12pt}, from=1-5, to=1-4]
	\arrow["\dashv"{anchor=center, rotate=-90}, draw=none, from=0, to=1]
	\arrow["\dashv"{anchor=center, rotate=-90}, draw=none, from=2, to=3]
\end{tikzcd}\end{center}
\end{lemma}

In fact, as (co)limits can be seen as special cases of adjunctions (see e.g. 6.3.10 in \cite{Riehl} and in the complete setting, see e.g. 4.5.1 in \cite{Riehl}),
the above result \ref{lem:comp_of_adjoints} in combination with the present lemma can be seen as some generalisation of (co)limits commuting with (co)limits \ref{lem:limits_commute_with_limits}.
In the same vein, that (co)limits in functor categories can be computed pointwise is a simple special case of this lemma
(together with currying of functor categories).

That (co)limits commute with (co)limits is more directly a special case of the following.
(The formation of limits is right adjoint and therefore commutes with limits.)
\begin{lemma}[Adjoints commute with (co)limits]\label{lem:adjoints_commute_with_limits}\uses{def:adjunction, def:limit, def:limit}
Any left adjoint functor commutes with all colimits, and any right adjoint functor commutes with all limits.
\end{lemma}
This elementary result has far-reaching consequences for the computation and existence of (co)limits,
as soon as one can install adjunctions.
\begin{remark} A nice trick for proving adjunctions is, that if one is able to prove cocontinuity of the potential left adjoint,
 afterward it suffices to check the adjunction (in the left hand side) on a generating class of the category.
	\end{remark}

Next we will show the general reason for universal morphisms into free constructions (units) being monomorphisms;
e.g., the embedding of a set into its Stone-\v{C}ech compactification or of a set into the free abelian group over it is monic.
\begin{lemma}[(Co)unit isomorphism via fully faithfulness]\label{lem:unit_via_fully_faithfullness}\uses{def:full_faithful_essentially_surjective_functors,lem:split_epi_mono_is_iso,def:adjunction, def:special_morphisms}
Consider any adjunction $L\dashv R$.
	\begin{enumerate}
		\item The right adjoint $R$ is faithful precisely if all components of the counit $\eps_D\colon LRD\to D$ are epic.
		\item The right adjoint $R$ is full precisely if all components of the counit are split monic.
		\item Consequently (using \ref{lem:split_epi_mono_is_iso}), $R$ is fully faithful precisely if all components of the counit are isomorphisms.
		\item The left adjoint $L$ is faithful precisely if all components of the unit $\eta_C\colon C\to RLC$ are monic.
		\item The left adjoint $L$ is full precisely if all components of the unit $\eta_C\colon C\to RLC$ are split epic
		\item Consequently (using \ref{lem:split_epi_mono_is_iso}), $L$ is fully faithful precisely if all $\eta_C\colon C\to RLC$ are isomorphisms.
	\end{enumerate}
	Combining this we obtain that adjunctions are a generalisation of equivalences of categories:
	Equivalences are precisely those adjunctions such that the left and the right adjoint are fully faithful.
\end{lemma}
\begin{example}
  The above theorem directly implies the cocompleteness of any reflective subcategory of a cocomplete category
  because the (fully faithful) inclusion functor creates colimits.
  A special case of this result is that colimits in $\CHaus$ exist and may be computed as Stone-\v{C}ech compactifications of the colimits in $\Top$.
  Another example is that the completion of metric spaces commutes with coproducts.
\end{example}

Let us further remark that adjunctions can be glued along unions of categories.
\begin{lemma}\label{lem:tower-adjs}
	Assume that we have a commutative tower of adjunctions

\begin{center}\begin{tikzcd}
	\vdots & \vdots \\
	{\mcC_i} & {\mcD_i} \\
	{\mcC_{i+1}} & {\mcD_{i+1}} \\
	\vdots & \vdots \\
	{\varinjlim C_i} & {\varinjlim{D_i}}
	\arrow[hook', from=1-1, to=2-1]
	\arrow[hook', from=1-2, to=2-2]
	\arrow[""{name=0, anchor=center, inner sep=0}, "{L_i}", curve={height=-6pt}, from=2-1, to=2-2]
	\arrow[hook', from=2-1, to=3-1]
	\arrow[""{name=1, anchor=center, inner sep=0}, "{R_i}", curve={height=-6pt}, from=2-2, to=2-1]
	\arrow[hook', from=2-2, to=3-2]
	\arrow[""{name=2, anchor=center, inner sep=0}, "{L_{i+1}}", curve={height=-6pt}, from=3-1, to=3-2]
	\arrow[hook', from=3-1, to=4-1]
	\arrow[""{name=3, anchor=center, inner sep=0}, "{R_{i+1}}", curve={height=-6pt}, from=3-2, to=3-1]
	\arrow[hook', from=3-2, to=4-2]
	\arrow[hook', from=4-1, to=5-1]
	\arrow[hook', from=4-2, to=5-2]
	\arrow["L", curve={height=-6pt}, from=5-1, to=5-2]
	\arrow["R", curve={height=-6pt}, from=5-2, to=5-1]
	\arrow["\dashv"{anchor=center, rotate=-90}, draw=none, from=0, to=1]
	\arrow["\dashv"{anchor=center, rotate=-90}, draw=none, from=2, to=3]
\end{tikzcd}\end{center}
	where the downward arrows are fully faithful inclusions.
	Then $L=\varinjlim L_i$ and $R=\varinjlim R_i$ are adjoints.
	\end{lemma}
\begin{proof}
	As the embeddings are fully faithful, this follows immediately.
	Explicitly, take any $A\in \varinjlim \mcC_i$ and $B\in \varinjlim \mcD_i$.
	Then there exists $i\in I$ such that $A\in\mcC_i$ and $j\in I$ such that $B\in \mcC_j$, and after embedding the smaller index fully faithful into the bigger one, we may assume that $i=j$.
	But now
	\[\hom(LA,B)=\hom(L_i A,B)=\hom(A, R_i B)=\hom(A, B).\]
\end{proof}

\begin{lemma}[Criterion for adjoints to preserve projective objects]\label{lem:adjoints_preserve_projective_objects}\uses{def:special_morphisms,def:adjunction, def:special_objects}
If the right adjoint $R$ of an adjunction
\begin{center}\begin{tikzcd}
	{ \mcC} & \mcD
	\arrow[""{name=0, anchor=center, inner sep=0}, "L", curve={height=-12pt}, from=1-1, to=1-2]
	\arrow[""{name=1, anchor=center, inner sep=0}, "R", curve={height=-12pt}, from=1-2, to=1-1]
	\arrow["\dashv"{anchor=center, rotate=-90}, draw=none, from=0, to=1]
\end{tikzcd}\end{center}
maps regular epimorphisms to regular epimorphisms, then the left adjoint maps projective objects to projective objects.
\end{lemma}
\begin{proof}
	This is completely straight-forward.
	Consider any projective object $P\in\mcC$ and any regular epimorphism $p\colon A\twoheadrightarrow B$ in $\mcD$, as well as any morphism $f\colon LP\to B$

\begin{center}\begin{tikzcd}
	& A \\
	LP & B.
	\arrow["p", two heads, from=1-2, to=2-2]
	\arrow["f", from=2-1, to=2-2]
\end{tikzcd}\end{center}

	Applying the $\hom$-isomorphism to $f$ yields a morphism $\tilde f\colon P\to RB$, and applying $R$ to $p$ yields a regular epimorphism $Rp\colon RA\twoheadrightarrow RB$ by assumption.
	By projectivity of $P$ we obtain a lift $\tilde g\colon P\to RA$,

\begin{center}\begin{tikzcd}
	& RA \\
	P & RB.
	\arrow["Rp", two heads, from=1-2, to=2-2]
	\arrow["{\tilde g}", dashed, from=2-1, to=1-2]
	\arrow["{\tilde f}", from=2-1, to=2-2]
\end{tikzcd}.\end{center}
	Applying the $\hom$-adjunction once again, we obtain a morphism $g\colon LP\to A$.
	The naturality of the isomorphism $\phi\colon \hom(P,R-)\to \hom(LP,-)$ in the second coordinate precisely yields that

\begin{center}\begin{tikzcd}
	{\tilde g\in \hom(P, RA)} & {g\in \hom(LP, A)} &&& A \\
	{\tilde f\in \hom(P, RB)} & {f=pg\in \hom(LP,B)} & {\text{i.e.}} & LP & B
	\arrow["{\phi_A}", maps to, from=1-1, to=1-2]
	\arrow["{(Rp)_*}", maps to, from=1-1, to=2-1]
	\arrow["{p_*}", maps to, from=1-2, to=2-2]
	\arrow["p", two heads, from=1-5, to=2-5]
	\arrow["{\phi_B}"', from=2-1, to=2-2]
	\arrow["g", from=2-4, to=1-5]
	\arrow["f"', from=2-4, to=2-5]
\end{tikzcd}\end{center}
commutes.
\end{proof}

\begin{remark}
  Clearly,
  instead of regular epimorphisms one could just as well consider all epimorphisms and get a completely analogous result.
\end{remark}

\begin{example}
  One toy example of this theorem is that the Stone-\v{C}ech compactification of a discrete set is a projective compact Hausdorff space
  (as the left adjoints to $\CHaus\to\Top\to\Set$ compose).
  Another example is that free objects form projective objects, e.g., free abelian groups or free $R$-modules are projective.
\end{example}

Having seen some basic properties of adjunctions, we turn the focus towards calculating them.
\begin{lemma}[Computation of Adjoints]\label{lem:computation_of_adjoints}\uses{def:limit}
Consider any functor $R\colon \mcD\to\mcC$.
Then $R$ admits a left adjoint if and only if the following hold.
	\begin{enumerate}
		\item $R$ is continuous.
		\item The following limit exists:
		\[L(C)=\varprojlim_{C\to R(D)\in C/R} D.\]
			The diagram is given by the projection $C/R\to\mcC$,
			i.e., its form is the comma category $C/R$.

	\end{enumerate}
In this case clearly the left adjoint is defined objectwise via the formula above, and the morphism $\eta_C\colon C\to RL(C)$
	is the universal arrow of the limit $RL(C)=\varprojlim_{C\to R(D)}R(D)$ induced by the cone of all $C\to R(D)$.
\end{lemma}
\begin{proof}
	See, e.g., Prop. 3.6 in adjoint functor in nLab.
\footnote{\url{https://ncatlab.org/nlab/show/adjoint+functor}}
	This is a special case of the formula for pointwise Kan extensions.
\end{proof}

Using this, we can easily deduce some further criteria for a functor to admit a left adjoint.
The usual way to do so is to use the following criterion, see 4.6.1 in \cite{Riehl}.
\begin{lemma}[Adjoints via initial objects in comma categories]\label{lem:adjoints_comma_initial}\uses{lem:computation_of_adjoints}
	A functor $R\colon\mcD\to\mcC$ has a left adjoint if and only if all comma categories $C/R$ have an initial object.
\end{lemma}
However, we will avoid the comma category and rather give a direct proof of the following theorem. 
\begin{theorem}[Freyd's Adjoint Functor Theorem (\idx{FAFT})]\label{thm:FAFT}\uses{def:adjunction, lem:limits_along_final_subdiagrams, lem:limits_along_final_subdiagrams, lem:computation_of_adjoints}
For a complete category $\mcD$, a functor $R\colon \mcD\to\mcC$ admits a left adjoint if and only if it is continuous and it satisfies the following \idx{solution set condition}:
For any $C\in\mcC$ there exists a set(!!!) of morphisms $\{f_i \colon C\to R(D_i)\}_{i\in I}$ such that for any morphism
$f\colon C\to RD$ for at least one of the $f_i$'s there exists a $D_i \to D$ with

\begin{center}\begin{tikzcd}
	C \\
	{R(D_i)} & {R(D)} \\
	{D_i} & D
	\arrow["{f_i}"', from=1-1, to=2-1]
	\arrow["f", from=1-1, to=2-2]
	\arrow[from=2-1, to=2-2]
	\arrow[dashed, from=3-1, to=3-2]
\end{tikzcd}\end{center}

This theorem is also called the \idx{General Adjoint Functor Theorem}.
\end{theorem}
\begin{proof}
  We show that the solution set spans a small initial subcategory in the comma category $C/R$,
  i.e., that the inclusion functor is initial.
	The solution set condition is a simple reformulation of the first defining property of initial functors.
	The second property of initial functors follows using the following diagram
	(it is a good thing to know that comma categories of complete categories along continuous functors are complete, but not necessary)
\begin{center}\begin{tikzcd}
	{R(D_k)} && C \\
	& {R(D_i\times_DD_j)} && {R(D_i)} \\
	& {R(D_j)} && {R(D)} \\
	{D_k} & {D_i\times_DD_j} && {D_i} \\
	& {D_j} && D
	\arrow[dashed, from=1-1, to=2-2]
	\arrow["{f_k}"', from=1-3, to=1-1]
	\arrow[dashed, from=1-3, to=2-2]
	\arrow["{f_i}", from=1-3, to=2-4]
	\arrow["{f_j}"{pos=0.7}, from=1-3, to=3-2]
	\arrow["f"'{pos=0.7}, from=1-3, to=3-4]
	\arrow[from=2-2, to=2-4]
	\arrow[from=2-2, to=3-2]
	\arrow[from=2-4, to=3-4]
	\arrow[from=3-2, to=3-4]
	\arrow[from=4-2, to=4-4]
	\arrow[from=4-2, to=5-2]
	\arrow[from=4-4, to=5-4]
	\arrow[from=5-2, to=5-4]
	\arrow[from=4-1, to=4-2]
\end{tikzcd}\end{center}
thus by \ref{lem:limits_along_final_subdiagrams} leading to the existence of the required limit of \ref{lem:computation_of_adjoints}.
(The arrow $C\to R(D_{i}\times_{D}D_{j})$ is given by the universal property of $R(D_{i})\times_{R(D)}R(D_{j})$;
after all, $R$ is continuous by assumption.
The arrow $D_{k}\to D_{i}\times_{D}D_{j}$ comes from the solution set condition.)

		The \enquote{only} if direction follows by observing that the limit is an initial object in the index category and hence a solution set.
	\end{proof}
In some cases (e.g. when using $\mcD=\CHaus$) the solution set criterion is automatic.
\begin{corollary}[Special adjoint functor theorem (\idx{SAFT})]\label{cor:SAFT}\uses{thm:FAFT, def:cogenerating_set}
Let $\mcD$ be a complete category admitting all (not necessarily small!) intersections of subobjects of any fixed object $D\in\mcD$, and assume that $\mcD$ admits a small coseparating set of objects.
	Then a functor $\mcD\to\mcC$ is continuous if and only if it admits a left adjoint.
\end{corollary}
The condition of all intersections to exist may be dropped as soon as subobjects of any given object form a set (e.g. when using $\mcD=\Set, \Ab, \CHaus$).
In that case, the existence of all intersections is immediate from the completeness of $\mcD$.

\begin{definition}
  We say a category $\mcD$ is \idx{well-powered} if for any $D\in\mcD$
  tha class of subobjects of $D$ modulo isomorphism (of subobjects) is a set.
  In other words,
  there is always small skeleton of the full subcategory of the comma category $\mcD/D$
  spanned by monomorphisms $D'\hookrightarrow D$.
\end{definition}

\begin{corollary}[SAFT, version 2]\label{cor:SAFT2}\uses{cor:SAFT}
	If $\mcD$ is a complete, well-powered category with a small coseparating set, then a functor $\mcD\to\mcC$ admits a left adjoint if and only if it is continuous.
	\end{corollary}
As one example we will show a new proof of a classical Gelfand duality type result.
\begin{example}[Gelfand as adjunction]\label{ex:gelfand_adjunction}
	The Baire functor from $\CHaus$ to the category of $\sigma$-Boolean algebras, mapping every space to its Baire $\sigma$-algebra and every continuous function to the pullback morphism, is continuous.
	Thus by the SAFT, it admits a left adjoint $K\colon \sigma\mathbf{-Bool}^{\op}\to \CHaus$, mapping any $\sigma$-algebra $\Sigma$ to a compact Hausdorff space $K_\Sigma$, which we call the \idx{Stone representation space} of $\Sigma$.
	The $\hom$-isomorphism now reads as the natural isomorphism
	\[ C(K_\Sigma, T)\simeq \hom(\Sigma, \mathrm{Baire}(T)).\]
	Using $\Sigma$ as the measure algebra of any measure space and $T=[-n,n]\subseteq\R$, and taking the colimit $n\to\infty$,
	one indeed recovers a (more general version of) the classical Gelfand-type result
	\[C(K_\Sigma)\simeq L^\infty(\Sigma).\]
	Playing around with this example indeed yields a good generalisation of many duality results, by purely abstract reasons.
\end{example}
Other examples of the application of adjoint functor theorems include all examples above, e.g. the existence of Stone-\v{C}ech compactifications, of a Hausdorffification, free abelian groups etc.


\section[Kan extensions]{Kan extensions}

\begin{center}
\quot{All concepts are Kan extensions}{Mac Lane, Chapter X.7 of \cite{mac2013categories}.}
\end{center}
Having seen limits and adjunctions we are finally prepared to unify these concepts into one exciting idea: \emph{Kan extensions}.
One could have started category theory with Kan extensions,
and once having understood those,
all concepts we have presented so far fit beautifully into the world of Kan extensions.
However, the definition of Kan extensions seemed quite unnatural to us
when we saw it for the first time,
and it only began to make sense after having understood limits and adjunctions.
So we decided to stick to this more pedagogical approach
and hope that by now the reader feels well prepared for Kan extensions.
See chapter X of \cite{mac2013categories}, chapter 1 of \cite{Riehl2014}
or chapter 6 of \cite{Riehl} for most of the following statements.

The idea of Kan extensions is to extend a functor $\mcC\to \mcE$
to another category $\mcD$ as best as possible along some information $\mcC\to \mcD$
telling us how the category $\mcC$ can be found inside $\mcD$.

\begin{definition}[Kan extension]\label{def:kan_extensions}\uses{def:functor, def:natural_transformation}
Let $F\colon \mcC\to \mcE$ be a functor between categories $\mcC$ and $\mcE$ that we want to \idx{extend} along some other functor $G\colon \mcC\to\mcD$.

A \idx{left Kan extension} of $F$ along $G$ is a pair $(\Lan_G F, \eta)$ consisting of a functor $\Lan_{G} F\colon \mcD\to \mcE$ and a natural transformation $\eta\colon F\Rightarrow (\Lan_G F)\circ G$ as in the diagram
\begin{center}
\begin{tikzcd}
	\mcE \\
	\mcC & \mcD
	\arrow[""{name=0, anchor=center, inner sep=0}, "F", from=2-1, to=1-1]
	\arrow["G"', from=2-1, to=2-2]
	\arrow["{{{\Lan_G F}}}"', shift right, from=2-2, to=1-1]
	\arrow["\eta"{pos=0.4}, shorten <=2pt, Rightarrow, from=0, to=2-2]
\end{tikzcd}
\end{center}
such that for any other pair $(L, \nu)$ of a functor $L\colon \mcD\to \mcE$ and a natural transformation $\nu\colon F\Rightarrow L\circ G$,
there exists a unique natural transformation $\alpha\colon \Lan_{G} F\Rightarrow L$ making the following diagram of natural transformations commute

\begin{center}
\begin{tikzcd}
	\mcE \\
	\\
	\mcC && \mcD
	\arrow[""{name=0, anchor=center, inner sep=0}, "F", from=3-1, to=1-1]
	\arrow["G"', from=3-1, to=3-3]
	\arrow[""{name=1, anchor=center, inner sep=0}, "L"', curve={height=24pt}, from=3-3, to=1-1]
	\arrow[""{name=2, anchor=center, inner sep=0}, "{{\Lan_G F}}"{description, pos=0.8}, curve={height=-6pt}, from=3-3, to=1-1]
	\arrow["\nu"{description, pos=0.7}, shift right=5, curve={height=-6pt}, shorten <=5pt, Rightarrow, from=0, to=1]
	\arrow["\eta"{description}, curve={height=6pt}, shorten <=12pt, shorten >=6pt, Rightarrow, from=0, to=3-3]
	\arrow["\alpha", shift left=5, curve={height=-6pt}, shorten <=5pt, shorten >=5pt, Rightarrow, from=1, to=2]
\end{tikzcd}
\end{center}
This means that the diagram
\begin{center}
\begin{tikzcd}
	{\Lan_GF} & F & {(\Lan_GF)\circ G} \\
	L && {L G}
	\arrow["{\exists ! \alpha}"', dashed, Rightarrow,  from=1-1, to=2-1]
	\arrow["\eta",Rightarrow, from=1-2, to=1-3]
	\arrow["\nu"',Rightarrow, from=1-2, to=2-3]
	\arrow["{\alpha.G}",Rightarrow, from=1-3, to=2-3]
\end{tikzcd}\end{center}
is commutative.
This can also be rephrased as saying that the map $\alpha \mapsto (\alpha.G)\circ\eta$ is a natural bijection
\[
\mathrm{Nat}(\Lan_G F, -) \to \mathrm{Nat}(F, - \circ G).
\]

A \idx{right Kan extension} of $F$ along $G$ is a pair $(\Ran_G F, \eps)$ consisting of a functor $\Ran_G F\colon \mcD\to \mcE$ and a natural transformation
$\eps\colon F\Rightarrow (\Ran_G F)\circ G$
\begin{center}\begin{tikzcd}
	& \mcE \\
	\mcC & \mcD
	\arrow[""{name=0, anchor=center, inner sep=0}, "F", from=2-1, to=1-2]
	\arrow["G"', from=2-1, to=2-2]
	\arrow["{\Ran_G F}"'{pos=0.4}, from=2-2, to=1-2]
	\arrow["\eps"'{pos=0.3}, shorten >=2pt, Rightarrow, from=2-2, to=0]
\end{tikzcd}\end{center}
such that for every other pair $(R\colon\mcD\to\mcE, \mu\colon F\Rightarrow RG)$
there exists a unique $\beta\colon R\Rightarrow \Ran_G F$ such that the diagram
\begin{center}\begin{tikzcd}
	{\Ran_G F} & F & {(\Ran_GF)\circ G} \\
	R && RG
	\arrow["\mu", Rightarrow, from=2-3, to=1-2]
	\arrow["\eps"', Rightarrow, from=1-3, to=1-2]
	\arrow["{\exists ! \beta}", dashed,Rightarrow, from=2-1, to=1-1]
	\arrow["{\beta.G}"',Rightarrow, from=2-3, to=1-3]
\end{tikzcd}\end{center}
is commutative.
In other words, the map $\beta\mapsto \eps\circ (\beta.G)$ is a natural isomorphism
\[
    \mathrm{Nat}(-, \Ran_G F)\simeq \mathrm{Nat}(-\circ G, F).
\]
\end{definition}

\begin{remark}
Similar to adjunctions, this concept can be defined in any 2-category, and Kan extensions are unique up to unique isomorphism, justifying the term \emph{the} Kan extension.
\end{remark}

\begin{remark}[Kan as adjoint]\label{rem:Kan-adjoint}
The characterisation of Kan extensions via isomorphisms of $\hom$-sets
between functors shows that for a functor $G\colon \mcC\to\mcD$ and a category $\mcE$,
if for all functors $F\colon\mcC\to\mcE$ the left Kan extension of $F$ along $G$ exist,
then the functor 
\[
    \Lan_G\colon [\mcC,\mcE]\to [\mcD, \mcE],\quad K\mapsto\Lan_G K
\]
is the left adjoint of the pullback functor
\[
    G^*\colon[\mcD, \mcE]\to [\mcC,\mcE],\quad H\mapsto H\circ G.
\]
The dual statement holds for the right Kan extension which is then the right adjoint
of the pullback $G^*$.
See also 6.1.5 in \cite{Riehl}.
\end{remark}

First, we show that Kan extensions are transitive.

\begin{lemma}[Kan extensions along compositions]\label{lem:kan_along_comp}\uses{def:kan_extensions}
Let $G\colon \mcC\to\mcD$, $H\colon \mcD\to \mcF$ and $F\colon \mcC\to \mcE$ be functors.
Assume that the left Kan extension $(\Lan_G F,\eta)$ exists.
If furthermore $(\Lan_{H}(\Lan_G F),\mu)$ exists,
it is the left Kan extension $(\Lan_{HG} F, (\mu.G)\circ \eta)$ of $F$ along $HG$.
\begin{center}
\begin{tikzcd}
	\mcE \\
	\\
	\mcC && \mcD && \mcF
	\arrow[""{name=0, anchor=center, inner sep=0}, "F", from=3-1, to=1-1]
	\arrow["G"', from=3-1, to=3-3]
	\arrow[""{name=1, anchor=center, inner sep=0}, "{{{\Lan_G F}}}"{description}, from=3-3, to=1-1]
	\arrow["H"', from=3-3, to=3-5]
	\arrow["{{{\Lan_H(\Lan_G F)}}}"', from=3-5, to=1-1]
	\arrow["\eta"{description}, shorten <=6pt, shorten >=6pt, Rightarrow, from=0, to=3-3]
	\arrow["\mu"{description}, shorten <=18pt, shorten >=9pt, Rightarrow, from=1, to=3-5]
\end{tikzcd}
\end{center}
\end{lemma}

\begin{proof}
Since $(\Lan_G F, \eta)$ is the left Kan extension of $F$ along $G$,
the map
\[
    \hom(\Lan_G F, L)\to\hom(F,L\circ G),\quad \alpha\mapsto (\alpha.G)\circ\eta
\]
is an isomorphism (bijection) natural in $L\colon\mcD\to\mcE$.
Similarly, $(\Lan_H(\Lan_G F),\mu)$ being the left Kan extension of $\Lan_G F$ along $H$
is the statement that the map
\[
    \hom(\Lan_H(\Lan_G F),K)\to\hom(\Lan_G F,K\circ H),\quad \beta\mapsto (\beta.H)\circ\mu
\]
is a bijection natural in $K\colon\mcF\to\mcE$.
Composing these two maps yields the bijection
\[
    \hom(\Lan_H(\Lan_G F),K)\to\hom(F,K\circ (HG)),\quad
    \beta\mapsto ((\beta.H)\circ\mu).G\circ\eta
        = (\beta.HG)\circ(\mu.G\circ\eta)
\]
which is natural in $K$ since composition of natural transformations
(in the functor categories) yields natural transformations.
In particular, $(\Lan_H \Lan_G F, (\mu.G)\circ \eta)$ is the left Kan extension of $F$ along $HG$.
\end{proof}

A similar statement holds for right Kan extensions by reversing the direction of the natural transformations.
This result can be combined with another direction of preservation of Kan extensions,
the postcomposition with a fixed functor.

\begin{definition}[Preservation of Kan extensions]\label{def:preserves_kan}\uses{def:kan_extensions}
Let $(\Ran_G F,\eps)$ be the right Kan extension
of a functor $F\colon\mcC\to\mcD$ along $G\colon\mcC\to\mcD$,
and $H\colon\mcE\to \mcF$ be another functor.
We say that $H$ \idx{preserves the right Kan extension} $\Ran_G F$
if $(H\Ran_G F,H.\eps)$ is the right Kan extension of $FH$ along $G$:
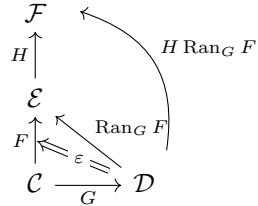
\begin{center}
  \begin{tikzcd}
	\mcF \\
	\mcE \\
	\mcC & \mcD
	\arrow["H", from=2-1, to=1-1]
	\arrow[""{name=0, anchor=center, inner sep=0}, "F", from=3-1, to=2-1]
	\arrow["G"', from=3-1, to=3-2]
	\arrow["{H\Ran_G F}"', shift right=4, curve={height=18pt}, from=3-2, to=1-1]
	\arrow["{\Ran_GF}"', from=3-2, to=2-1]
	\arrow["\eps"{description}, shorten <=5pt, Rightarrow, from=3-2, to=0]
\end{tikzcd}
\end{center}
An analogous definition applies to left Kan extensions.
\end{definition}

There are special cases of this preservation which have names.

\begin{definition}[Pointwise and absolute Kan extensions]\label{def:pointwise_kan_extensions}\uses{def:kan_extensions, def:preserves_kan}\chapone
A right Kan extension $\Ran_G F$ of $F\colon \mcC\to\mcE$ along $G\colon\mcC\to\mcD$ is called
\begin{enumerate}
	\item \textbf{pointwise}\index{pointwise Kan extension}, if it is preserved by all representable functors $y_E=\hom_{\mcE}(E, -)\colon\mcE\to\Set$,
	\item \textbf{absolute}\index{absolute Kan extension}, if it is preserved by every functor.
\end{enumerate}
\end{definition}

\begin{remark}
The version for left Kan extensions that are pointwise
naturally uses the contravariant $\hom$-functor $y^E = \hom(-,E)$,
leading to the necessity of passing to the opposite direction of some 2-morphisms.

A left Kan extension $\Lan_ G F\colon \mcD\to\mcE$ of $F\colon\mcC\to\mcD$ along $G\colon\mcD\to\mcE$
is precisely the right Kan extension of $F\colon C^{\op}\to\mcE^{\op}$ along $G\colon \mcC^{\op}\to \mcD^{\op}$.
Thus, a left Kan extension is called \textbf{pointwise},
if for every $E\in\mcE$ the functor $y^E=\hom_{\mcE}(-,E)\colon\mcE^{\mathrm{op}}\to \Set$ preserves the right Kan extension $\Lan_G F\colon \mcD^{\op}\to\mcE^{\op}$.
A left Kan extension is called \textbf{absolute}, if it is preserved by every functor.
\end{remark}

\begin{example}\label{ex:adjoints_preserve_kan}
Left adjoint functors $L$ preserve left Kan extensions,
and right adjoint functors preserve right Kan extensions,
see~Theorem X.5.1 in \cite{mac2013categories}.

In particular, if $\mcE$ has small copowers
(meaning small coproducts over copies of the same object,
$\coprod_{i\in I} E$ with $E\in\mcE$),
then right Kan extensions are pointwise.
Dually, if $\mcE$ admits all small powers ($\prod_{i\in I} E$ with $E\in\mcE$),
then left Kan extensions are pointwise.

We urge the reader to think through that $\Set\to\mcE$, $I\mapsto\coprod_{I}E$ is left adjoint to $\hom_{\mcE}(E,-)$.
\end{example}

Unfortunately, Kan extensions in general need not be preserved by many functors,
although in practice they are often preserved by representable functors.
Thus, many authors claim (e.g. \cite{Riehl2014}) that the better concept of Kan extensions
are pointwise Kan extensions, due to the
\enquote{failure to find a single instance where a [non-pointwise] Kan extension plays any mathematical role whatsoever} (Kelly in \cite[§4]{Kelly1982}).
For us,
this is convincing due to the following lemma and theorem~\ref{thm:computation_of_kan_extensions} which we will frequently use when working with Kan extensions.

Since \emph{All Concepts Are Kan Extensions}, clearly adjoints and limits form special cases of Kan extensions, see e.g. 6.5.1, 6.5.2 in \cite{Riehl} or chapter X.7 in \cite{mac2013categories} for the following result.

\begin{lemma}[Adjunctions and limits are absolute Kan extensions]\label{lem:adjunctions_and_limits_are_kan_extensions}\uses{def:kan_extensions, def:adjunction, def:limit}
   \begin{enumerate}[(i)]
\item Consider an adjunction $L\dashv R$ with unit $\eta$ and counit $\eps$.

Then the \emph{left} adjoint $(L, \eps)$ is the absolute \emph{right} Kan extension of the identity $1_{\mcD}$ along $R$,
\begin{center}\begin{tikzcd}
	& \mcD \\
	\mcD & \mcC
	\arrow[""{name=0, anchor=center, inner sep=0}, "{1_D}", from=2-1, to=1-2]
	\arrow["R"', from=2-1, to=2-2]
	\arrow["L"', from=2-2, to=1-2]
	\arrow["\eps"', shorten <=2pt, Rightarrow, from=2-2, to=0]
\end{tikzcd},\end{center}
and the \emph{right} adjoint $(R,\eta)$ is the absolute \emph{left} Kan extension of $1_{\mcC}$ along $F$,
\begin{center}\begin{tikzcd}
	& \mcC \\
	\mcC & \mcD
	\arrow[""{name=0, anchor=center, inner sep=0}, "{1_C}", from=2-1, to=1-2]
	\arrow["L"', from=2-1, to=2-2]
	\arrow["R"', from=2-2, to=1-2]
	\arrow["\eta", shorten <=2pt, Rightarrow, from=0, to=2-2]
\end{tikzcd}.
\end{center}
Conversely, if $(R,\eta)$ is the left Kan extension of $1_C$ along $L$ and $L$ preserves this Kan extension, then $L\dashv R$ (and similarly for $L$).

\item The colimit of any diagram $D\colon \mcI\to\mcC$ exists precisely if the left Kan extension of $D$ along the functor $\mcI\to \ast$ to the terminal category exists.
In this case, the value of the left Kan extension on the point is the colimit.
Dually, a limit of $D$ is the same thing as a right Kan extension of $D$ along $\mcC\to\ast$.
\end{enumerate}
\end{lemma}

Note that, since Kan extensions are unique, this in particular implies the known uniqueness results of limits and adjoints.
Furthermore, in view of the lemma, example~\ref{ex:adjoints_preserve_kan} can be seen as a generalisation of the fact that left adjoints commute with colimits and right adjoint commute with limits.
The lemma also shows that the natural transformations in Kan extensions need not be isomorphisms.
We will now investigate when this is the case.

\begin{lemma}[Kan extension extends functor]\label{lem:Kan_extension_extends}\uses{def:kan_extensions}
Pointwise Kan extensions along a fully faithful functor $G\colon\mcC\to\mcD$ have natural isomorphisms as unit and counit, i.e.,
\[
    (\Ran_G F)\circ G\simeq F\simeq (\Lan_G F)\circ G
\]
for functors $F\colon\mcC\to\mcE$ and $G\colon\mcC\to\mcD$.
If the functor $G$ is a full inclusion, i.e., $\mcC$ is a full subcategory of $\mcD$,
then the isomorphisms are equalities.
\end{lemma}

\begin{proof}
See Corollary 3 and 4 in chapter X.3. of \cite{mac2013categories}.
\end{proof}

This means that extending along a fully faithful functor really yields an \emph{extension}.

Having seen some elementary properties of (pointwise) Kan extensions,
we look for ways to compute them and criteria for their existence.

\begin{theorem}[Computation of Kan extensions]\label{thm:computation_of_kan_extensions}\uses{def:kan_extensions} \chapone
Let $F\colon\mcC\to\mcE$ and $G\colon\mcC\to\mcD$ be functors.
The right Kan extension $\Ran_G F$ is pointwise if any only if it is given by
\[
    \Ran_G F(X)=\varprojlim_{(f,C)\in (X\downarrow G)} F(C)
\]
for every $X\in\mcD$.
Dually, the left Kan extension $\Lan_G F$ is pointwise if and only if it is computed by
\[
    \Lan_G F(X)=\varinjlim_{(f,C)\in (G\downarrow X)} F(C)
\]
for all $X\in\mcD$.
If the relevant (co)limits exist in $\mcE$,
the given formulas define pointwise left resp.\ right Kan extensions.

In these cases, one can recover the counit $\eps$ (resp. the unit $\eta$)
and the extended functor on morphisms analogously to \ref{lem:computation_of_adjoints}.
\end{theorem}

\begin{proof}
See 6.3.7 in \cite{Riehl}.
\end{proof}

This in particular explains the terminology \emph{pointwise}.
Note that the theorem implies the formula for computation of adjoints,
compare the formula there in \ref{lem:computation_of_adjoints}.

Next, similarly to the adjoint functor theorem,
we give a criterion for the existence of the occurring (co)limits.

\begin{theorem}[Existence criterion]\label{thm:crit_for_kan}\uses{thm:computation_of_kan_extensions}
If the category $\mcC$ is small and $\mcE$ is (co)complete,
then the pointwise right (resp. left) Kan extension of any functor $F\colon \mcC\to\mcE$
along any functor $G\colon \mcC\to\mcD$ exists.
\end{theorem}

\begin{remark}
  The smallness condition on $\mcC$ can be replaced by a solution set condition
  ensuring that the relevant (co)limits can be computed as small (co)limits.
  This is a straightforward generalisation of \ref{thm:FAFT}.
\end{remark}

In view of remark~\ref{rem:Kan-adjoint},
this also assures a right resp. left adjoint to the pullback functor $G^*$.
The theorem also implies the commutativity of limits with limits
and the composition result of adjunctions.

We have already observed that often, Kan extensions are not absolute.
However, as we will see in Chapter 3, in homotopical categories,
they are absolute.
In this case, the following lemma is of great importance.

\begin{lemma}[Absolute Kan preserves adjoints]\label{lem:Kan_preserves_adjoints}\uses{def:pointwise_kan_extensions}
Let $L\dashv R$ be an adjunction and $G\colon \mcC\to\tilde \mcC$, $H\colon \mcD\to\tilde\mcD$ be functors,
and assume that the \emph{absolute} right Kan extension $\tilde L=\Ran_{G}HL$ of $HL$ along $G$
and the \emph{absolute} left Kan extension $\tilde R=\Lan_{H}GR$ of $GR$ along $H$ exist.
Then they are again adjoint, $\tilde L\dashv \tilde R$.
\begin{center}
\begin{tikzcd}
	\mcD & {\tilde\mcD} \\
	\mcC & {\tilde\mcC}
	\arrow["H", from=1-1, to=1-2]
	\arrow[""{name=0, anchor=center, inner sep=0}, "R", curve={height=-6pt}, from=1-1, to=2-1]
	\arrow[""{name=1, anchor=center, inner sep=0}, "{{\tilde R}}", curve={height=-6pt}, from=1-2, to=2-2]
	\arrow[""{name=2, anchor=center, inner sep=0}, "L", curve={height=-6pt}, from=2-1, to=1-1]
	\arrow["G"', from=2-1, to=2-2]
	\arrow[""{name=3, anchor=center, inner sep=0}, "{{\tilde L}}", curve={height=-6pt}, from=2-2, to=1-2]
	\arrow["\dashv"{anchor=center}, draw=none, from=2, to=0]
	\arrow["\dashv"{anchor=center}, draw=none, from=3, to=1]
\end{tikzcd}
\end{center}
\end{lemma}

\begin{proof}
Denote by $\eta\colon 1\to RL$ the unit and by $\varepsilon\colon LR\to 1$
the counit of the adjunction $L\dashv R$.
Let $\mu\colon\tilde{L}G\Rightarrow HL$ be the universal natural transformation of $\tilde{L}$
and $\nu\colon GR\Rightarrow\tilde{R}H$ be the one for $\tilde{R}$.
Since the Kan extensions are absolute, $(\tilde{R}\tilde{L},\tilde{R}.\mu)$
is the right Kan extension of $\tilde{R}HL\colon\mcC\to\tilde{\mcC}$ along $G$.
Hence for $\xi = \nu.L\circ G.\eta\colon 1_{\tilde{\mcC}}G\Rightarrow GRL\Rightarrow\tilde{R}HL$
there exists a unique $\alpha\colon 1_{\tilde{\mcC}}\Rightarrow\tilde{R}\tilde{L}$
such that $\tilde{R}.\mu\circ\alpha.G = \xi$.
Similarly, $(\tilde{L}\tilde{R},\tilde{L}.\nu)$ is the left Kan extension of $\tilde{L}GR\colon\mcD\to\tilde{\mcD}$ along $H$.
For $\zeta = H.\varepsilon\circ\mu.R\colon \tilde{L}GR\Rightarrow HLR\Rightarrow 1_{\tilde{\mcD}}H$ there exists a unique $\beta\colon \tilde{L}\tilde{R}\Rightarrow 1_{\tilde{\mcD}}$ such that $\beta.H\circ\tilde{L}.\nu = \zeta$.

We now show that $\alpha$ and $\beta$ are the unit and counit making $\tilde{L}$ and $\tilde{R}$ adjoint to each other.
By the triangle identities it suffices to show that $(\beta.\tilde{L})\circ(\tilde{L}.\alpha) = 1_{\tilde{L}}$ and $(\tilde{R}.\beta)\circ(\alpha.\tilde{R}) = 1_{\tilde{R}}$.
This reduces to showing that the triangles
\begin{center}
\begin{tikzcd}
	& {\tilde{L}G} &&& {\tilde{R}H} \\
	{\tilde{L}G} && HL & GR && {\tilde{R}H}
	\arrow["\mu", Rightarrow, from=1-2, to=2-3]
	\arrow["{((\tilde{R}.\beta)\circ(\alpha.\tilde{R})).H}", Rightarrow, from=1-5, to=2-6]
	\arrow["{((\beta.\tilde{L})\circ(\tilde{L}.\alpha)).G}", Rightarrow, from=2-1, to=1-2]
	\arrow["\mu"', Rightarrow, from=2-1, to=2-3]
	\arrow["\nu", Rightarrow, from=2-4, to=1-5]
	\arrow["\nu"', Rightarrow, from=2-4, to=2-6]
\end{tikzcd}
\end{center}
are commutative, because then the triangle identities follow from the uniqueness part
in the definition of Kan extension.
Now these are straightforward verifications where one needs to be careful
with the 2-categorical structure on natural transformations.
We have by the triangle identity for the adjoint pair $L\dashv R$, that
\begin{align*}
        \mu
    =   (H.1_L)\circ\mu
    =   (H.(\varepsilon.L\circ L.\eta))\circ\mu
    =   H.\varepsilon.L\circ HL.\eta\circ\mu
\end{align*}
Now since $\mu$ is a natural transformation and vertical composition intertwines with horizontal composition, we obtain
\begin{align*}
        H.\varepsilon.L\circ HL.\eta\circ\mu
    =   H.\varepsilon.L\circ\mu.RL\circ\tilde{L}G.\eta
    =   (H.\varepsilon\circ\mu.R).L\circ\tilde{L}G.\eta
    =   \zeta.L\circ\tilde{L}G.\eta
\end{align*}
Now by the characterisation of $\zeta$, one can deduce
\begin{align*}
        \zeta.L\circ\tilde{L}G.\eta
    =   (\beta.H\circ\tilde{L}.\nu).L\circ\tilde{L}G.\eta
    =   \beta.HL\circ\tilde{L}.\nu.L\circ\tilde{L}G.\eta
    =   \beta.HL\circ\tilde{L}.(\nu.L\circ G.\eta)
    =   \beta.HL\circ\tilde{L}.\xi.
\end{align*}
Because $\xi$ factors over $\tilde{R}.\mu$, we obtain
\begin{align*}
        \beta.HL\circ\tilde{L}.\xi
    =   \beta.HL\circ\tilde{L}.(\tilde{R}.\mu\circ\alpha.G)
    =   \beta.HL\circ\tilde{L}\tilde{R}.\mu\circ\tilde{L}.\alpha.G
\end{align*}
Since $\beta$ is a natural transformation, we interchange again
\begin{align*}
        \beta.HL\circ\tilde{L}\tilde{R}.\mu\circ\tilde{L}.\alpha.G
    =   \mu\circ\beta.\tilde{L}G\circ\tilde{L}.\alpha.G
\end{align*}
which is exactly what we wanted.
The second verification works almost the same,
but here one starts with
\begin{align*}
    \nu = \nu\circ G.1_R = \nu\circ G.(R.\varepsilon\circ\eta.R)
\end{align*}
using the second triangle identity for the adjunction $L\dashv R$.
Then again using the factorisation for $\xi$ and then for $\zeta$
and interchanging vertical with horizontal composition,
one arrives at the desired result.
\end{proof}

To end the section about Kan extensions,
let us remark how they behave under unions of categories.

\begin{lemma}\label{lem:kan-extension-fully-faithful}
Consider a union $\varinjlim_{i\in I}\mcC_i$ of categories $\mcC_i$
along fully faithful inclusions $\tau_i^j \colon C_i\to C_j$,
where $I$ is a totally ordered class.
Let $F\colon \mcC_k\to\mcD$ be a functor for some $k\in I$,
and assume that for all $j\ge k$ the left Kan extension $(\Lan_{\tau_{k}^j}F,\eta_j)$ exists.
Then the induced functor $\tilde{F}\colon\varinjlim_{i\in I}\mcC_i\to\mcD$
is the left Kan extension of $\Lan_{\tau_{k}^j}F$ along the embedding $\tau_j\colon\mcC_j\to\varinjlim _{i\in I}\mcC_i$ for every $j\geq k$.
\begin{center}
\begin{tikzcd}
	{\mcC_k} \\
	{\mcC_j} && D \\
	\vdots \\
	{\varinjlim\mcC_i}
	\arrow["{\tau_k^j}"', hook', from=1-1, to=2-1]
	\arrow["F", from=1-1, to=2-3]
	\arrow["{\Lan_{\tau_k^j}F}"', from=2-1, to=2-3]
	\arrow[hook', from=2-1, to=3-1]
	\arrow[hook', from=3-1, to=4-1]
	\arrow["{\tilde{F}}"', dashed, from=4-1, to=2-3]
\end{tikzcd}
\end{center}
\end{lemma}

\begin{proof}
First of all, by lemma~\ref{lem:kan_along_comp},
for $i\geq j\geq k$, we have that $\Lan_{\tau_j^i}(\Lan_{\tau_k^j}F) = \Lan_{\tau_k^i}F$,
so without loss of generality we may assume $j = k$.
Take any functor $G\colon\varinjlim\mcC_i\to D$ with natural transformation
$\nu\colon F\Rightarrow G\circ \tau_k$.
By restriction, for $j\geq k$ we have a natural transformation
$\nu_j\colon F\Rightarrow(G\circ\tau_j)\circ\tau_k^j$.
Here, $\Lan_{\tau_k^j}F$ exists and thus there is a unique natural transformation
\[
    \alpha_j\colon\tilde{F}\circ\tau_j = \Lan_{\tau_k^j}F\Rightarrow G\circ\tau_j
\]
such that $\alpha_j.\tau_k^j\circ\eta_j = \nu_j$.
The $\alpha_j$ together yield a natural transformation $\alpha\colon\tilde{F}\Rightarrow G$, which is the unique natural transformation such that $\alpha.\tau_k\circ\eta = \nu$.
\end{proof}

\begin{remark}
The same argument can be applied to the right Kan extension of $F$.
\end{remark}

For further uses and implications of Kan extensions,
e.g. a more general Yoneda lemma \cite[6.5.4]{Riehl}, we refer to the literature.

Kan extensions form a special case of universal arrows of cocompletions.
\begin{question}
	For any small category $\mcC$ with $\mcM$-continuous cocompletion $\tau\colon \mcC\to\overline\mcC$, and any $\mcM$-continuous functor $G\colon \mcC\to \mcD$ with $\mcD$ being
	cocomplete, right/left Kan extension $\Ran_\tau G$ yields the unique cocontinuous functor

\begin{center}\begin{tikzcd}
	\mcC & {\overline{\mcC}} \\
	\mcD
	\arrow["\tau", from=1-1, to=1-2]
	\arrow["G"', from=1-1, to=2-1]
	\arrow["{\Ran_{\tau} G}", from=1-2, to=2-1]
\end{tikzcd}\end{center}

\end{question}

\section{Sheaves and topoi}\label{sec:sheaves-topoi}

In this chapter we refine the understanding of functor categories and find especially good presheaves, so-called \emph{sheaves}.
In fact,
this potential to be a sheaf is the reason why they are called presheaves.
The appearance of sheaf categories, so-called \emph{Grothendieck Topoi},
unifies multiple ideas at once.
One can view them as non-free cocompletions, leading to excellent categorical properties, which in turn completely classify categories of sheaves.
On the other hand, one can view a sheaf directly as a presheaf
preserving certain limits,
that encode the geometric structure of the category in terms of \emph{coverings}.

The classical approach to sheaf categories,
which we follow first, is to start with an idea of which collections of morphisms
with common codomain one wants to be coverings.
This information is gathered in a so-called \emph{Grothendieck topology} on the category.
Now sheaves are those functors that respect coverings,
i.e., one can recover the information of the sheaf on a space
by knowing its values on a covering of the space.
This idea arose from considering the category $\mathbf{Ouv}(X)$
of open sets of a topological space $X$,
where a family $U_i\subset U$ covers $U$ iff $U = \bigcup_{i\in I} U_i$.
There, the presheaf $U\mapsto C(U,\R)$ which assigns to every open set $U\subset X$
the set of continuous real valued functions on $U$ is a sheaf.
Indeed, one can recover a continuous function $f\colon U\to\R$
from the information on an open cover of $U$.
This sheaf is called the \emph{sheaf of continuous functions}.

Now the classical notion of coverings on a topological space lead to the
definition of a Grothendieck topology by relaxing the fact that
all maps $U_i\to U$ need to be open embeddings,
and just considering general morphisms $X\to Y$.

There exists a good amount of literature, we refer to
\cite{Johnstone, Lane1992SheavesIG, Borceux1994, kashiwara2005categories, Vakil2023}.
In particular,
\cite[chap.\ 2]{Vakil2023} contains a lovely introduction to the main ideas;
another is at \cite{Urbanik2019}.

\subsection{Sites and sheaves}\label{ssec:sites-sheaves}

First we discuss the notion of \emph{covering},
which is similar to the usual way of defining a topology:
To every object (subset) $X$ one assigns a collection $\mcT(X)$ of families of morphisms 
(subsets) with codomain $X$, which we call \emph{coverings} of $X$,
satisfying certain axioms.

\begin{definition}[Coverage]\label{def:coverage}\uses{def:category}
A \idx{coverage} on a category $\mcC$ consists of a map $\mcT$ assigning
to each $X\in\mcC$ a class of families
\footnote{The families itself need not be small.
This is often dealt with by working in universes.}
$(f_i\colon X_i\to X)_{i\in\mcI}$ of morphisms with codomain $X$,
which are called \textbf{covering families}\index{covering} of $X$, such that
\begin{enumerate}[]
    \item for every covering $(f_i\colon X_i\to X)$ and every morphism $g\colon Y\to X$
        there exists a covering $(h_j\colon Y_j\to Y)$ of $Y$
        such that all $gh_j$ factor through a morphism $f_i$:
        \begin{center}
        \begin{tikzcd}
        	{Y_j} & Y \\
        	{X_i} & X
        	\arrow["{h_j}", from=1-1, to=1-2]
        	\arrow[dashed, from=1-1, to=2-1]
        	\arrow["g", from=1-2, to=2-2]
        	\arrow["{f_i}"', from=2-1, to=2-2]
        \end{tikzcd}
        \end{center}
\end{enumerate}
\end{definition}

The intuition behind this condition is that if the $f_i$ cover $X$,
then the preimages of the images of all $f_i$ under $g$ cover the domain of $g$.
Since the required pullbacks need not exist,
we allow to refine the covering by the $h_j$.
However, often one can check the criterion by trying one of the following two possible coverings.

\begin{lemma}[Criterion for coverage with pullbacks]\label{lem:coverage_with_pullback}\uses{def:coverage, def:special_limits_and_colimits}
Consider an assignment $\mcT$ mapping $X$ to collections of families $(f_i\colon X_i\to X)$.
Then one of the following conditions is sufficient for $\mcT$ being a coverage.
\begin{enumerate}[(i)]
		\item If $\mcC$ has pullbacks, and for every covering $(f_i\colon X_i\to X)$
		and every morphism $g\colon Y\to X$, the pullback family
		$(f_i^*\colon X_i\times_X Y\to Y)$ is a covering of $Y$.
		\item For every covering $S=(f_i\colon X_i\to X)$
		and every morphism $g\colon Y\to X$ the \idx{pullback sieve} $g^{*}(S)$,
		defined as
        \[
            g^*(S) = \{h_j\colon Y_j\to Y: 
                \,\text{there exists}\, f_i\in\mcC,\, k\colon Y_j\to X_i
                \,\text{with}\, f_{i}k=gh_j \},
        \]
		is a covering.
\end{enumerate}
If furthermore $\mcT$ is \idx{upwards closed},
i.e., any family $(f_i\colon X_i\to X)_{i\in I}$ that contains a covering
(an element of $\mcT(X)$), is itself a covering of $X$,
then condition (ii) is equivalent to $\mcT$ being a coverage.
\end{lemma}

\begin{definition}[Site]\label{def:site}\uses{def:coverage}
A \idx{site} is a category $\mcC$ equipped with a coverage.
A \idx{small site} is a site on a small category $\mcC$.
\end{definition}

Having defined what coverings are,
we can define those functors that respect this structure.

\begin{definition}[Preliminary definition of sheaves]\label{def:sheaf}\uses{def:site}
Let $(\mcC,\mcT)$ be a site and $F\colon \mcC^{\op}\to \Set$ a presheaf.
\begin{enumerate}[(i)]
	\item For a family of morphisms $(f_i\colon X_i\to X)_{i\in\mcI}$ we define a \idx{compatible family} to be a family $(s_i)_{i\in\mcI}$ with $s_i\in F(X_i)$ such that for all $h\colon V\to X_i$ and $g\colon V\to X_j$ with $hf_i=gf_j$,
\begin{center}\begin{tikzcd}
	V & {X_j} \\
	{X_i} & X
	\arrow["g", from=1-1, to=1-2]
	\arrow["h"', from=1-1, to=2-1]
	\arrow["{f_j}", from=1-2, to=2-2]
	\arrow["{f_i}"', from=2-1, to=2-2]
\end{tikzcd},\end{center}
we have $F(h)(s_i)=F(g)(s_j)\in F(V)$.
Under the existence of pullbacks in $\mcC$ this is equivalent to $s_i$ and $s_j$ agreeing on the intersection, i.e. for all $i,j\in\mcI$ with pullback diagram
\begin{center}\begin{tikzcd}
	{X_{i}\times_X X_j} & {X_j} \\
	{X_i} & X
	\arrow["{\tau_i^j}", from=1-1, to=1-2]
	\arrow["{\tau_j^i}"', from=1-1, to=2-1]
	\arrow["{f_j}", from=1-2, to=2-2]
	\arrow["{f_i}"', from=2-1, to=2-2]
\end{tikzcd}\end{center}
one has $F(\tau_i^j)s_j=F(\tau_j^i)s_{i}$.
Two elements $s, t\in F(X)$ \idx{locally agree} on $X$,
if there exists a covering $S$ of $X$ such that $F(f)s=F(f)t$ for all $f\in S$.

\item We say that $F$ satisfies the \idx{sheaf condition} with respect to the family
$(f_i\colon X_i\to X)_{i\in\mcI}$ if for all compatible families $(s_i)_{i\in\mcI}$ there exists a unique $s\in F(X)$ with $F(f_i)(s)=s_i$ for all $s_i$.
Now with the explicit description of limits in $\Set$,
see lemma~\ref{lem:set_bicomplete},
this can be rephrased to $(F(X),(F(f_i))_{i})$ being a limit cone in the (large) diagram obtained by applying $F$ to the (large) diagram
\begin{center}
\begin{tikzcd}
	V & {X_i} \\
	W & {X_j} && X \\
	\dots & {X_k}
	\arrow[from=1-1, to=1-2]
	\arrow[from=1-1, to=2-2]
	\arrow["{{f_i}}", from=1-2, to=2-4]
	\arrow[shift left, from=2-1, to=2-2]
	\arrow[shift right, from=2-1, to=2-2]
	\arrow["{{f_j}}"', from=2-2, to=2-4]
	\arrow[from=3-1, to=2-2]
	\arrow[from=3-1, to=3-2]
	\arrow["{{f_k}}"', from=3-2, to=2-4]
\end{tikzcd}
\end{center}
consisting of all the $X_i$ and all $V$ with morphisms to $X_i,X_j$ building a commutative square.
Note that we include the case with $X_i=X_j$, and in this case the commutative square
        in fact are just two parallel arrows.
        Again,
        if $\mcC$ admits pullbacks,
        we can require the diagram to only consist of the $X_{i}$ and $X_{i}\times_{X}X_{j}$ (in the \enquote{second layer}).

  \item If $(\mcC,\mcT)$ is a small\footnote{
        There are ways of avoiding this smallness condition and retaining a workable equalizer diagram.
        For example,
        one could take $F\colon \mcC^{\op}\to\Set_{\mcU}$ and then take the products in $\Set_{\mcV}$
        for a larger Grothendieck universe $\mcV$,
        working with product classes and only using theexplicit description of equalizers,
        assuming that all relevant families $(f_{i}\mid i\in \mcI)$ are indexed by sets (this is rarely the case)
        or assuming the cone above to always possess an initial sub\textit{set} (this is often the case)
        and requiring the equalizer condition for this restricted diagram.
        }
        site (and hence $\mcI$ is small as well),
    the presheaf $F$ satisfies the sheaf condition for $(f_i\colon X_i\to X)_{i\in\mcI}$
    precisely if the diagram
\begin{center}\begin{tikzcd}
	{F(X)} & {\prod\limits_{f_i}F(X_i)} & {\prod\limits_{(h_V^i, g_V^j)} F(V)}
	\arrow["f", from=1-1, to=1-2]
	\arrow["h", shift left, from=1-2, to=1-3]
	\arrow["g"', shift right, from=1-2, to=1-3]
\end{tikzcd}\end{center}
is an equalizer diagram,
where the first product is taken over all morphisms $f_i$.
The morphism $f$ is the one induced by all $F(f_i)$
from the universal property of the product.
The second product is taken over to all tuples
$(h_V^i\colon V\to X_i,\, g_V^j\colon V\to X_j)$ which satisfy $f_ih_V^i=f_jg_V^j$.
The morphism $h$ is induced by the coordinates
\[
    F(h_V^i)\circ\pi_i\colon \prod_{f_i} F(X_i)\to  F(X_i)\to F(V)
\]
and $g$ is defined similarly by gluing $F(g_V^j)\circ\pi_j$ together.
If $\mcC$ admits fibre products
\begin{center}
\begin{tikzcd}
	{X_i\times_X X_j} & {X_j} \\
	{X_i} & X
	\arrow["{\tau_i^j}", from=1-1, to=1-2]
	\arrow["{\tau_j^i}"', from=1-1, to=2-1]
	\arrow["{f_j}", from=1-2, to=2-2]
	\arrow["{f_i}"', from=2-1, to=2-2]
\end{tikzcd},
\end{center}
this in turn is equivalent to the diagram
\begin{center}
\begin{tikzcd}
	{F(X)} & {\prod\limits_{f_i}F(X_i)} & {\prod\limits_{(f_i, f_j)} F(X_i\times_X X_j)}
	\arrow["f", from=1-1, to=1-2]
	\arrow["{h}", shift left, from=1-2, to=1-3]
	\arrow["{g}"', shift right, from=1-2, to=1-3]
\end{tikzcd}
\end{center}
being an equalizer diagram.
Here, $h$ is the gluing of all morphisms
\[
    F(\tau_j^i)\circ\pi_i\colon \prod_{f_i} F(X_i)\to F(X_i\times_X X_j)
\]
and $g$ the gluing of all $F(\tau_i^j)\circ \pi_j$.

	\item Finally, a \idx{sheaf} on $\mcC$ is a presheaf $F\colon \mcC^{\op}\to \Set$
	that satisfies the sheaf condition with respect to all coverings of every object.
\end{enumerate}
The full subcategory of $\PSh(\mcC)$ consisting of all sheaves on $\mcC$ is called the \idx{category of sheaves} on $\mcC$ and denoted by $\Sh(\mcC)$.
(Of course, a more orthodox notation would be $\Sh_{\mcT}(\mcC)$, but $\mcT$ is usually left implicit.)
\end{definition}

\begin{remark}
One should read this definition as follows:
the $X_i$ are some subspaces of $X$ with inclusions $f_i$,
the $X_i\times_X X_j$ are the intersections of $X_i$ and $X_j$,
together with an embedding $\tau_j^i$ of the intersection $X_i\cap X_j$ into the first subspace $X_i$,
and an embedding $\tau_i^j$ of the intersection into the second space $X_j$.
The functor $F$ can be imagined as assigning to each subspace a set of functions on this space,
and turns every inclusion into restriction of functions along these inclusions.
Hence the function $f$ in the equalizer diagram in (iii) of the definition
assigns to a globally defined function the tuple of restrictions $(f\lvert_{X_i})$.
The functions $h$ and $g$ of this equalizer diagram take any tuple of functions living on $X_i$ and restrict them further to twofold intersections $X_i\cap X_j$.
Now the equalizer condition reads as follows:
any tuple of functions defined on all subspaces $X_i$ individually,
that agree on overlaps $X_i\cap X_j$,
glue in a unique way to a globally defined function.
\end{remark}

We now extend the definition of a sheaf to define the category of $\mcD$-valued sheaves,
$\Sh(\mcC,\mcD)$, for any category $\mcD$ in the obvious way, using the large limit.

\begin{definition}[$\mcD$-valued sheaf]\label{def:D-sheaf}\uses{def:sheaf}
Let $(\mcC,\mcT)$ be a site and $\mcD$ a category.
\begin{enumerate}[(i)]
\item Let $S = (f_i\colon X_i\to X)_{i\in\mcI}$ a family of morphisms with codomain $X$.
Then we build a diagram $C_S$ in $\mcC$ which consists of every $X_i$ for $i\in\mcI$
and if $V$ is an object in $\mcC$ such that there exist $i,j\in\mcI$
and two morphisms $g\colon V\to X_i$ and $h\colon V\to X_j$,
then $V$ is also an object of the diagram, and so are $g$ and $h$.
A small portion of diagram $D_{S}:=F\circ C_S$ can be pictured as follows.
\begin{center}
\begin{tikzcd}
	{F(W)} && {F(X_i)} \\
	\\
	{F(V)} && {F(X_j)}
	\arrow["{F(k)}"', shift right, from=1-3, to=1-1]
	\arrow["{F(l)}", shift left, from=1-3, to=1-1]
	\arrow["{{F(g)}}", from=1-3, to=3-1]
	\arrow["{{F(h)}}", from=3-3, to=3-1]
\end{tikzcd}.
\end{center}
\item A functor $F\colon\mcC^{\op}\to \mcD$ satisfies the sheaf condition for $S$,
if the cone $(F(X),(F(f_i))_{i\in\mcI})$ is a limiting cone for $D_S$ in $\mcD$.

\item A functor $F\colon\mcC^{\op}\to \mcD$ is called a $\mcD$-valued \idx{sheaf},
if $F$ satisfies the sheaf condition for every $X\in\mcC$ and all coverings $S$ of $X$.
\end{enumerate}
We denote by $\Sh(\mcC,\mcD)$ the full subcategory of $[\mcC^\mathrm{op},\mcD]$
whose objects are $\mcD$-valued sheaves.
\end{definition}

In general, the category of sheaves on a site is a large category,
except for when $\mcC$ is a small site.

\begin{remark}
At first sight it seems that this definition might be too general
or too abstract to be computable.
In practice, one often has a (forgetful) functor $\mcD\to\Set$
that creates limits.
In this case, the sheaf condition can be computed in $\Set$
after forgetting the additional structure of $\mcD$.
This will be the case for (concrete) algebraic categories
as for example the category of abelian groups.
And, of course, for $\mcD = \Set$,
this definition yields the preliminary definition~\ref{def:sheaf}
and we write $\Sh(\mcC,\Set) = \Sh(\mcC)$.
In another direction of specialisation,
for $\mcT=\emptyset$,
we get $[\mcC^{\op},\mcD]=\Sh_{\mcT}(\mcC,\mcD)$ as a special case.

For more general situations,
especially if intersections of more than two \enquote{subspaces} $X_i\cap X_j$
are of interest,
one should switch to hypercoverings and hypersheaves to obtain a slightly stronger theory.
To us,
it appears that this is more elegantly handled in the $\infty$-categorical setting.
In any case,
the $1$-categorical $1$-sheaves suffices for now.
\end{remark}

\begin{definition}[Grothendieck topos]\label{def:grothendieck_topos}\uses{def:sheaf}
A \idx{Grothendieck topos} is a category that is equivalent to a category $\Sh(\mcC)$ for a \emph{small} site $(\mcC,\mcT)$.
\end{definition}

\begin{example}[Examples of sheaves]
Consider the category of open sets $\mathbf{Ouv}(X)$ of a topological space $X$.
For an object $U$ we say that a family of subsets $(U_i\subset U)_{i\in I}$ covers $U$
iff $\bigcup_{i\in I} U_i = U$.
This is the so-called $\sup$-topology on $\mathbf{Ouv}(X)$.
Then for any topological space $Y$ the functor
\[
    C(-,Y)\colon U\mapsto C(U,Y)=\hom_{\mathrm{Top}}(U, Y)
\]
(and mapping inclusions to restrictions) is a sheaf.
If the space $X$ has the structure of a smooth manifold,
another example is the sheaf of holomorphic functions $U\mapsto \mathrm{Hol}(U, \C)$,
or the sheaf of (real) differentiable functions $U\mapsto C^1(X,\C)$.
If $X$ has the structure of an algebraic variety,
then the sheaf of regular functions (local quotients of polynomials) is a sheaf.
\end{example}

\subsection{Stability of the sheaf condition}\label{ssec:stab-sheaf-cond}

Next, we will show that the sheaf condition is preserved under certain modifications of the coverage,
allowing us to impose further properties on the coverage without changing the category of sheaves.

\begin{definition}[Sieve]\label{def:sieve}\uses{def:category}
A \idx{sieve} on an object $X$ is a right ideal in the class of all homomorphisms with target $X$,
i.e., it is a collection $S=(f_i\colon X_i\to X)$ such that for all $h\colon V\to X_i$ one has $f_ih\in S$.
We say the sieve is \idx{stable under precomposition}.

Equivalently, a sieve on $X$ is a \idx{subfunctor} of the $\hom$-functor $\hom(-,X)$
(meaning a subobject in the functor category).

The \idx{generated sieve} of a family $(f_i\colon X_i\to X)$ is the smallest sieve containing this family.
It consists of all compositions $f_i h$ for every $h\colon V\to X_i$.
\end{definition}

If $S$ is a sieve on $X$, then the sheaf condition in definition~\ref{def:D-sheaf}~(ii)
has another description.
For this, denote by $(S\downarrow X)$ the full subcategory of $\mcC/X$
whose objects are the elements of $S$
and whose morphisms are factorisations of elements of $S$ by another element of $S$
(i.e., usually do not lie in $S$ themselves as their codomain may not be $X$).
Then we have a canonical forgetful functor $P\colon (S\downarrow X)\to \mcC$
which yields a (large) diagram in $\mcC$.

\begin{lemma}\label{lem:sheaf-cond-colimit}
Let $F\colon\mcC^\mathrm{op}\to\mcD$ be a $\mcD$-valued presheaf
and let $S$ be a sieve on $X$.
Then $F$ satisfies the sheaf condition for $S$ if and only if
$(F(X),(F(f))_{f\in S})$ is a limiting cone of the diagram $FP$.
\end{lemma}

\begin{proof}
If $F$ satisfies the sheaf condition for $S$,
let $(D,(\pi_f)_{f\in S})$ be a cone for $FP$.
Now, if we have morphisms $g\colon V\to X_i$, $h\colon V\to X_j$
for $f_i\colon X_i\to X$ and $f_j\colon X_j\to X$ in $S$ such that $f_ig = f_jh$,
Since $f_ig = f_jh$ is an element of $S$,
and therefore $g\colon (V,f_ig)\to (X_i,f_i)$ is a morphism in the diagram $P$,
we have that $\pi_{f_ig} = F(g)\pi_{f_i}$ and similarly $\pi_{f_jh} = F(h)\pi_{f_j}$
But now $f_ig = f_jh$ implies $F(g)\pi_{f_i} = F(h)\pi_{f_j}$.
Thus we have extended the cone for $FP$ to a cone for $D_S$ and since there $F(X)$
is the limit, we obtain that it is also the limit of $FP$.

For the converse, let $(D,(\pi_f)_{f\in S})$ be a cone for $D_S$.
Now, if $f\colon Y\to X$ and $g\colon Z\to X$ are in $S$ and $h\colon (Y,f)\to (Z,g)$
is a morphism in $(S\downarrow X)$, then $f1_Y = gh$ and so $\pi_f = F(h)\pi_g$.
Thus, the cone for $D_S$ is also a cone for $FP$ and hence $F(X)$ is also a limit of $FP$.
\end{proof}

\begin{definition}[Sifted coverage]\label{def:sifted-coverage}\uses{def:sieve}
A coverage $\mcT$ on a category $\mcC$ is called \idx{sifted},
\footnote{
  As with set-theoretic filters and filtered categories,
  there are obvious similarities between sifted coverages and sifted categories.
}
if all its coverings are sieves.

If $\mcT$ is a collection $(\mcT(X))_{X\in\mcC}$ of families of morphisms with codomain $X$,
then we define $\overline{\mcT}(X)$ to consist of all the generated sieves of elements of $\mcT$
and call the collection $(\overline{\mcT}(X))_{X}$ the \idx{sifting} of $\mcT$.
Clearly,
$\mcT$ is a coverage precisely if $\overline{\mcT}$ is one.

Recall that a coverage is called \emph{upwards closed}, if supersets of coverings are coverings.
Define the \idx{upwards closure} of any coverage by defining a family to cover as soon as it contains a covering.
\end{definition}

\begin{lemma}[Stability of the sheaf property]\label{lem:stab_sh}\uses{def:sheaf, def:sieve}
Let $F$ be a presheaf on a site $(\mcC,\mcT)$ with values in $\mcD$.
Then the following hold.
\begin{enumerate}[(i)]
    \item For any family of morphisms $S = (f_i\colon X_i\to X)$,
	the presheaf $F$ fulfills the sheaf condition with respect to $S$
	precisely if it satisfies the sheaf condition for the sieve $\overline{S}$
	generated by $S$.	
	\item If $f\colon Y\to X$ is an isomorphism,
	then $F$ satisfies the sheaf condition for $(f)$.
\end{enumerate}
For the next two points, assume that $F$ is a sheaf for $\mcT$
and let $(f_i\colon X_i\to X)_{i\in I}$ and $(f_{ij}\colon X_{ij}\to X_i)_{j\in J_i}$
be covering families for $\mcT$.
\begin{enumerate}[(i)]
\setcounter{enumi}{2}
	\item The sheaf $F$ satisfies the sheaf condition for $(f_if_{ij}\colon X_{ij}\to X)_{i\in I,j\in J_i}$.
	\item If $(g_j\colon Y_j\to X)$ is a family such that each $f_i$ factors through a $g_j$,
	then $F$ satisfies the sheaf condition for the family of the $g_j$.
\end{enumerate}
\end{lemma}

\begin{proof}
This is C2.1.3 and C2.1.5-C2.1.7 in \cite{Johnstone}
but because we have a slightly more general definition of sheaf,
we give the proofs.

Proof of (i).
If $F$ satisfies the sheaf condition for $S$,
then it clearly satisfies the sheaf condition for $\overline{S}$:
for if $(A,(\pi_f)_{f\in\overline{S}})$ is a cone for $D_{\overline{S}}$,
then it also is a cone for $D_S$ and hence there exists a unique morphism
$m\colon A\to F(X)$ with $\pi_{f_i} = F(f_i)m$.
Then for $f\in\overline{S}$ we have $f = f_ig$ for some morphism $g$
and therefore $\pi_f = F(g)\pi_{f_i} = F(g)F(f_i)m = F(f)m$.

If, on the other hand, $F$ satisfies the sheaf condition for $\overline{S}$,
then we can extend any cone $(A,(\pi_{f_i})_{f_i\in S})$ over $D_S$ to a cone
$(A,(\pi_f)_{f\in S})$ over $D_{\overline{S}}$ by setting $\pi_f = F(g)\pi_{f_i}$
for $f = f_ig$.
This is well-defined by definition of $D_S$.
Therefore, $F$ satisfies the sheaf condition for $S$.

The second assertion is clear as $F(f)$ is an isomorphism.

To prove (iii), let $(A,(\pi_{ij})_{i,j})$ be a cone for $D_S$
where $S = (f_if_{ij}\colon X_{ij}\to X)_{i,j}$.
For any $i$, the cone $(A,(\pi_{ij})_{j})$ is a cone for $D_{R_i}$
where $R_i = (f_if_{ij}\colon X_{ij}\to X)_j$.
For morphisms $g\colon V\to X_{ij}$, $h\colon V\to X_{ij'}$
with $f_{ij}g = f_{ij'}h$ we have $f_if_{ij}g = f_if_{ij'}h$
and therefore $F(g)\pi_{ij} = F(h)\pi_{ij'}$.
The presheaf $F$ satisfies the sheaf condition for $R_i$
and so there exists a unique morphism $m_i\colon A\to F(X_i)$
such that $F(f_{ij})m_i = \pi_{ij}$ for every $j$.
Again, for $h\colon W\to X_i$, $\ell\colon W\to X_{i'}$ with $f_ik = f_{i'}\ell$,
there exists a covering $(g_k\colon Y_j\to W)_{j\in J_i}$
such that every composition $g_kh$ and $g_k\ell$ factors over $f_{ij}$ and $f_{i'j}$ respectively (use the property of coverages three times), as in
\begin{center}
\begin{tikzcd}
	{Y_k} && {X_{i'j}} \\
	& W & {X_{i'}} \\
	{X_{ij}} & {X_i} & X
	\arrow["{f'_{i'j}}", dashed, from=1-1, to=1-3]
	\arrow["{g_k}", from=1-1, to=2-2]
	\arrow["{f'_{ij}}"', dashed, from=1-1, to=3-1]
	\arrow["{f_{i'j}}", from=1-3, to=2-3]
	\arrow["\ell", from=2-2, to=2-3]
	\arrow["h"', from=2-2, to=3-2]
	\arrow["{f_{i'}}", from=2-3, to=3-3]
	\arrow["{f_{ij}}"', from=3-1, to=3-2]
	\arrow["{f_i}"', from=3-2, to=3-3]
\end{tikzcd}.
\end{center}
Now from the commutativity in the diagram $f_if_{ij}f'_{ij} = f_{i'}f_{i'j}f'_{i'j}$,
we obtain by applying the functor $F$, that
\[
    F(g_k)F(h)m_i = F(f'_{ij})F(f_{ij})m_i = F(f'_{ij})\pi_{ij} =
    F(f'_{i'j})\pi_{i'j} = F(f'_{i'j})F(f_{i'j})m_{i'} = F(g_k)F(\ell)m_{i'}.
\]
By the sheaf condition, $F(W)$ is the limit over the diagram generated by the $F(Y_k)$
and so we get $F(h)m_i = F(\ell)m_{i'}$.
Thus $(D,(m_i)_{i\in I})$ is a cone for the diagram $D_R$ where $R = (f_i\colon X_i\to X)$.
Hence there exists a unique $m\colon D\to F(X)$ such that $F(f_i)m = m_i$.
This implies $F(f_if_{ij}) = F(f_{ij})m_i = \pi_{ij}$ and $m$ is unique with this property.
So $F$ satisfies the sheaf condition for $S$.

For (iv), let $(D,(\pi_j))$ be a cone over the diagram of the covering $(g_j\colon Y_j\to X)$.
This extends to a cone over the diagram of the covering $(f_i\colon X_i\to X)$.
Set $\pi_{i} = \pi_jF(f_{ij})$ for $f_i = g_jf_{ij}$ where
$f_{ij}\colon X_i\to Y_j$.
This is well defined by definition of the diagram for the covering of the $g_j$.
Now for if $h\colon V\to X_i$ and $k\colon V\to X_{i'}$ such that $f_ih = f_{i'k}$,
then there exist $j,j'$ and morphisms $f_{ij}\colon X_i\to Y_j$, $f_{i'j'}\colon X_{i'}\to Y_{j'}$ such that $f_i = g_jf_{ij}$ and $f_{i'} = g_{j'}f_{i'j'}$.
In particular, $g_jf_{ij}h = g_{j'}f_{i'j'}k$.
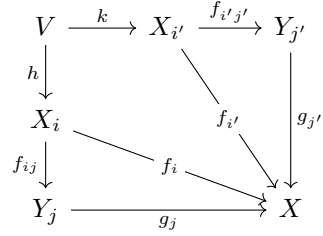
\begin{center}
\begin{tikzcd}
	V & {X_{i'}} & {Y_{j'}} \\
	{X_i} \\
	{Y_j} && X
	\arrow["k", from=1-1, to=1-2]
	\arrow["h"', from=1-1, to=2-1]
	\arrow["{f_{i'j'}}", from=1-2, to=1-3]
	\arrow["{f_{i'}}"{description}, from=1-2, to=3-3]
	\arrow["{g_{j'}}", from=1-3, to=3-3]
	\arrow["{f_{ij}}"', from=2-1, to=3-1]
	\arrow["{f_i}"{description}, from=2-1, to=3-3]
	\arrow["{g_j}"', from=3-1, to=3-3]
\end{tikzcd}
\end{center}
This means
\[
    F(h)\pi_i = F(h)F(f_{ij})\pi_j = F(k)F(f_{i'j'})\pi_{j'} = F(k)\pi_{i'}.
\]
By the sheaf condition there exists a unique $m\colon D\to F(X)$ such that
$F(f_i)m = \pi_i$.
Now the $f_i$ cover $X$, so there exists a cover $(h_k\colon Z_k\to Y_j)$
such that each $g_{j_0}h_k$ factors through $f_i$ (for a fixed $j_0$)
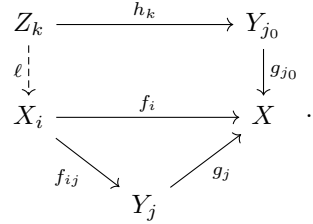
\begin{center}
\begin{tikzcd}
	{Z_k} && {Y_{j_0}} \\
	{X_i} && X \\
	& {Y_j}
	\arrow["{{h_k}}", from=1-1, to=1-3]
	\arrow["\ell"', dashed, from=1-1, to=2-1]
	\arrow["{{g_{j_0}}}", from=1-3, to=2-3]
	\arrow["{{f_i}}", from=2-1, to=2-3]
	\arrow["{f_{ij}}"', from=2-1, to=3-2]
	\arrow["{g_j}"', from=3-2, to=2-3]
\end{tikzcd}.
\end{center}
Then
\[
    F(h_k)F(g_{j_0})m = F(\ell)F(f_i)m = F(\ell)\pi_i = F(\ell)F(f_{ij})\pi_j = F(k)\pi_{j_0} 
\]
where the last equality follows from the fact that $g_jf_{ij}\ell = f_i\ell = g_{j_0}k$.
Since the $h_k$ cover $Y_j$ it follows that $F(g_{j_0})m = \pi_{j_0}$
which was to be shown.
\end{proof}

\begin{remark}
In view of part (i) of the lemma, one could also formulate the assertions (ii)-(iv)
in terms of sieves and sifted coverings and immediately obtain an equivalence
between these two points of view.
For more details see C2.1 in \cite{Johnstone}.
\end{remark}

In particular, one can assume for any coverage $\mcT$ that isomorphisms are coverings,
that it is sifted,
and that if for a sieve $S$ on $X$, there exists a covering $R$ of $X$
such that for every $f\in R$ the pullback sieve $f^*(S)$ is a covering,
then $S$ itself is a covering of $X$.
All this does not change the category of sheaves.

Collecting these additional properties of coverings,
we obtain the following definition.

\begin{definition}[Grothendieck Topology]\label{def:grothendieck_topology}\uses{def:sieve}
A \idx{Grothendieck topology} $\mcT$ on a category $\mcC$ assigns to every object $X$ a collection of sieves $\mcT(X)$ on $X$, again called coverings, such that
\begin{enumerate}[(i)]
	\item for every object the maximal sieve,
	consisting of all morphisms with codomain $X$, is a covering,
	\item if $R$ covers $X$ and $f\colon Y\to X$, then $f^{*}(R)$ covers $Y$, and
	\item if $R$ covers $X$ and $S$ is a sieve on $X$ such that for all $f\colon Y\to X\in R$ the sieve $f^*(S)$ covers $Y$, then $S$ covers $X$.
\end{enumerate}
In particular, Grothendieck topologies are upwards closed (with respect to sieves).
Furthermore, the \idx{generated Grothendieck topology} of a coverage $\mcT$ is defined as the
intersection of all Grothendieck topologies containing $\mcT$.
\end{definition}

\begin{corollary}[Sheaves on generated Grothendieck topologies]\label{cor:sh_on_gen_groth_top}\uses{def:grothendieck_topology, lem:stab_sh}
The passage to the generated Grothendieck topology does not affect the sheaves,
i.e., a presheaf is a sheaf for a coverage precisely if it is a sheaf for the generated Grothendieck topology.
\end{corollary}

Often one implicitly identifies the coverage with its generated Grothendieck topology.
We now look at some special Grothendieck topologies on a category.

\begin{example}[Subcanonical]\label{ex:subcanonical}\uses{def:site}
A coverage $\mcT$ on a category $\mcC$ is called \idx{subcanonical}
if all representable presheaves are sheaves.
In this case, the Yoneda embedding allows us to identify $\mcC$ with a full subcategory of $\Sh(\mcC)$.
The largest sifted subcanonical coverage is called the \idx{canonical} coverage.
\end{example}

There is also a description of subcanonical coverages only in terms of the category $\mcC$.

\begin{lemma}\label{lem:subcanonical-colimit}
A coverage $\mcT$ on a category $\mcC$ is subcanonical precisely if
for any covering $S = (f_i\colon X_{i}\to X)$ the diagram $C_S$
(equivalently, the diagram $P$) has $X$ as colimit,
i.e., the covering families form colimiting cocones.
\end{lemma}

\begin{proof}
In view of lemma~\ref{lem:sheaf-cond-colimit}
(taking $\mcD = \mcC^\mathrm{op}$ and $F$ the identity functor)
the claims for the diagram $C_S$ and for the diagram $P$ are equivalent.
The other equivalence follows from the definition of limits via the $\hom$-functor,
see~\ref{lem:hom_functor_preserves_limits}.
\end{proof}

Thus, a site is subcanonical precisely if the condition for $F\colon \mcC^{\op}\to \mcD$
to be a sheaf reduces to $F$ mapping certain colimits in $\mcC$ to limits in $\mcD$.

Note the similarity to $\mcM$-continuous cocompletions: there, one uses presheaves
preserving certain limits.
In fact, one often is able to identify categories of sheaves with $\mcM$-continuous cocompletions.
However, a priori the class of limits that need to be preserved is not purely specified by their shape,
other than in the case of $\mcM$-continuous cocompletions.
In the other direction, however, there might be a connection.

\begin{conjecture}\label{con:non-uniform-cocompletion}
One can define cocompletions more generally by not specifying the index categories
but rather the diagrams themselves (see~\ref{def:non-uniform-coco}).
We believe that this would yield stronger results that allow to describe subcanonical 
Grothendieck topoi as \emph{non-uniform} cocompletions.
\end{conjecture}

\subsection{Properties of categories of sheaves}\label{ssec:properties-sheaves}

We have seen that there are certain operations modifying the coverage on a site
without changing the sheaf category.
In this section we establish an important theorem that
allows us to substantially restrict the underlying category of a site without changing the sheaf category.
For this we need a notion of subsite.
We assume that all coverages $\mcT$ are Grothendieck topologies.

\begin{definition}[Dense subsites]\label{def:dense_subsite}\uses{def:grothendieck_topology}
Let $(\mcC,\mcT)$ be a site and $\mcD$ a full subcategory of $\mcC$.
The \idx{induced Grothendieck topology} $\mcT_{\mcD}$ on $\mcD$ is given
by taking coverings of $X\in\mcD$ to be $S\cap\hom(\mcD)$ for $S\in\mcT(X)$.

We call $(\mcD,\mcT_{\mcD})$ a \idx{dense subsite}
if every object in $\mcC$ has a covering sieve generated by morphisms with domain in $\mcD$.
\end{definition}

It is immediately verified that if $\mcT$ is a Grothendieck topology,
then so is $\mcT_{\mcD}$.
Before stating the main theorem,
we make some observations, see C2.2.2 in \cite{Johnstone}.

\begin{lemma}\label{lem:dense-subsite-properties}
Let $(\mcC,\mcT)$ be a site and $(\mcD,\mcT_{\mcD})$ a dense subsite.
Then the following hold.
\begin{enumerate}[(i)]
    \item A sieve $S$ on $X\in\mcD$ is in $\mcT_{\mcD}(X)$ if and only if
            the generated sieve $\overline{S}$ in $\mcC$ is in $\mcT(X)$.
    \item Let $F$ be an $\mcE$-valued sheaf on $(\mcC,\mcT)$.
            Then the restriction of $F$ to $\mcD$ is a sheaf on $(\mcD,\mcT_{\mcD})$.
\end{enumerate}
\end{lemma}

\begin{proof}
Proof of (i).
If $S$ is in $\mcT_{\mcD}(X)$,
then there exists a covering sieve $R$ on $X$ in $\mcC$,
such that $S = R\cap\mathrm{mor}(\mcD)$.
Now for $f\colon Y\to X$ in $R$, let $T$ be the covering sieve on $Y$ generated by
morphisms with domain in $\mcD$.
Then the sieve $f^*(\overline{S})$ contains $T$ and therefore is a covering sieve.
Hence, $\overline{S}$ is a covering sieve.
If $\overline{S}$ is covering, then $S = \overline{S}\cap\mathrm{mor}(\mcD)$ is
covering in $\mcD$.

Proof of (ii).
Let $S$ be a covering sieve in $\mcD$ of $X$.
We need to show that the cone $(F(X),(F(f))_{f\in S})$ is a limit cone
of the diagram $D_S$ in $\mcD$ as defined in definition~\ref{def:D-sheaf}.
For this it suffices, since $F$ is a sheaf on $\mcC$,
that the we can extend any cone $(E,(\pi_f)_{f\in S})$ over $D_S$
to a cone over $D'_S$, where $D'_S$ is now generated in $\mcC$, rather than in $\mcD$
(since by lemma~\ref{lem:stab_sh} it suffices to check the sheaf condition on a generator).
Let $f\colon Z\to X$ and $f'\colon Z'\to X$ be two morphisms in $S$
and $g\colon V\to Z$, $h\colon V\to Z'$ morphisms in $\mcC$ with $fg = f'h$.
We need to show that $F(g)\pi_f= F(h)\pi_{f'}$.
For this, take a covering $(f_i\colon V_i\to V)$ in $\mcC$ where $V_i\in\mcD$.
Then, since the compositions $gf_i$ and $hf_i$ are morphisms in $\mcD$,
and $(E,(\pi_f)_{f\in S})$ is a cone over $D_S$ in $\mcD$,
we have $F(f_i)F(g)\pi_f = F(f_i)F(h)\pi_{f'}$.
But now since $F$ is a sheaf on $\mcC$, it follows that $F(g)\pi_f = F(h)\pi_{f'}$
which was to be shown.
\end{proof}

With this we can now proof the promised comparison lemma.
See C.2.2.3 in \cite{Johnstone}, Corollary 3 in Chapter 4 of the Appendix in \cite{Lane1992SheavesIG} and \cite{Vistoli2007}
for different versions of the following lemma, which can be generalised to non-full subcategories.

\begin{theorem}[Comparison lemma]\label{thm:comparison_lemma}\uses{def:grothendieck_topos,cor:sh_on_gen_groth_top,  def:kan_extensions,prop:characterisation_of_equivalences}
Let $(\mcC,\mcT)$ be a site and $(\mcD,\mcT_{\mcD})$ be a \emph{small} dense subsite.
Let $\mcE$ be a complete category.
Then restriction to $\mcD$ induces an equivalence
\begin{center}
\begin{tikzcd}
	{\Sh(\mcC,\mcE)} && {\Sh(\mcD,\mcE)}
	\arrow["\simeq", from=1-1, to=1-3]
\end{tikzcd}
\end{center}
between categories.
The inverse (and right adjoint) is given by right Kan extension along the inclusion $\mcD\to\mcC$.
\end{theorem}

\begin{proof}
First of all, restriction is well-defined by lemma~\ref{lem:dense-subsite-properties}.
Furthermore, if $F\colon\mcD^\mathrm{op}\to\mcE$ is a sheaf,
then the pointwise right Kan extension of $F$ along the inclusion
$\mcD^\mathrm{op}\to\mcC^\mathrm{op}$ exists by theorem~\ref{thm:crit_for_kan}
because $\mcE$ is complete and $\mcD$ is small.

Next, we show that this right Kan extension $\Ran F$ is a sheaf on $\mcC$.
For this, let $S$ be a covering sieve of $X$ in $\mcC$
and $(E,(\pi_f)_{f\in S})$ a cone for $FP$
where $P\colon (S\downarrow X)\to\mcC$ is the corresponding diagram.
Let $(D,g)\in\mcD/X$.
Because $g^*(S)$ covers $D$ in $\mcC$,
the intersection $R = g^*(S)\cap\mathrm{mor}(\mcD)$ covers $D$ in $\mcD$.
For $(h,D')\in (R\downarrow D)$, $(k,D'')\in (R\downarrow D)$ and
$\ell\colon (k,D'')\to (h,D')$, we have that $gh\ell = gk$.
Since $gh,gk\in S$, we obtain that $\pi_{gk} = F(\ell)\pi_{gh}$.
Hence $(E,(\pi_{gh})_{(h,D')\in (R\downarrow D)})$ forms a cone over $FQ$
with $Q\colon (R\downarrow D)\to\mcD$.
Because $F$ is a sheaf, there exists a unique arrow $\tau_{(D,g)}\colon E\to F(D)$
such that $F(h)\tau_{(D,g)} = \pi_{gh}$.
If now $f\colon (D,g)\to (D_0,g_0)$ is a morphism in $\mcD/X$,
then $g_0fk = gk\in S$, so $fk\in g_0^*(S)$
which means that
\[
    F(k)F(f)\tau_{(D_0,g_0)}
    = F(fk)\tau_{(D_0,g_0)}
    = \pi_{g'fk}
    = \pi_{gk}
    = F(k)\tau_{(D,g)}.
\]
The uniqueness of $\tau_{(D,g)}$ implies that $\tau_{(D,g)} = F(f)\tau_{(D_0,g_0)}$.
In particular, $(E,(\tau_{(D,g)})_{(D,g)\in \mcD/X})$ is a cone.
Since
\[
    (\Ran F)(X) = \varprojlim_{(D,g)\in\mcD/X} F(D),
\]
there exists a unique morphism $E\to \Ran F(X)$ making the whole diagram commutative.
Therefore, $\Ran F$ is a sheaf on $\mcC$.

Having remark~\ref{rem:Kan-adjoint} in mind, it follows that
the right Kan extension is the right adjoint of the restriction.
It remains to show that the two functors are essentially inverse to each other.
By lemma~\ref{lem:kan-extension-fully-faithful},
the restriction of the right Kan extension of a sheaf on $\mcD$ is the sheaf itself.
For the other direction, let $F$ be a sheaf on $\mcC$ and denote by $G$ the right
Kan extension of its restriction.
For $X\in\mcC$ we have
\[
    G(X) = \varprojlim_{(D,f)\in\mcD/X} F(D)
\]
by the formula for pointwise Kan extension.
Since $\mcD/X$ contains a covering of $X$ by the denseness of $\mcD$,
the functor $F$ satisfies the sheaf condition for $\mcD/X$ by lemma~\ref{lem:stab_sh}~(iv).
But in view of lemma~\ref{lem:sheaf-cond-colimit} this is exactly the assertion
that $G(X) = F(X)$.
\end{proof}

\begin{remark}
Note that the smallness condition on $\mcD$ is only needed to assure the existence
of the right Kan extension of a sheaf.
In any setting where they exist,
one can omit the smallness condition.
\end{remark}

Now we turn our attention to important properties of categories of sheaves.
For this,
let us first explain how to turn $\mcT(X)$ ($\mcT$ a Grothendieck topology on $\mcC$, $X\in\mcC$)
into a category.
As it is naturally ordered by inclusion,
we can turn it into a category by taking reverse inclusions as morphisms.
This results in a filtered category $\mcT(X)$.

\begin{theorem}[Sheafification]\label{thm:sheafification}\uses{def:sheaf}
  For any site $\mcC$,
  a left adjoint to the inclusion $\Sh(\mcC)\to\PSh(\mcC)$ is called \idx{sheafification}, often denoted $\Sh$.
  For any Grothendieck topos $\Sh(\mcC)$ (i.e., small $\mcC$) there exists a sheafification.
  This can be computed in the following ways:
  \begin{enumerate}[(i)]
    \item In the case that $\mcC$ has pullbacks, one can use
          \[\Sh(F)(X)=\varinjlim_{S\in\mcT(X)}\left\{(s_i)\in\prod_{f_{i}\in S}X_{i}:  \text{ locally agree on } X_i\times_X X_j \text{ for all }f_{i},f_{j}\in S\right\}.\]
    \item Classically, one defines the \idx{Grothendieck $F^+$ construction}
          \[F^+(X)=\varinjlim_{S\in\mcT(X)}\{(s_i)\text{ compatible family}\}.\]
          Then $\Sh(F)=F^{++}$ (note that we have to do this process twice; in the first step one only obtains a separated presheaf).
  \end{enumerate}
\end{theorem}

See \cite{Johnstone}, \cite{Lane1992SheavesIG} and the blog post \emph{How to sheafify in one go}
\footnote{\url{https://web.stanford.edu/~dkim04/blog/sheafification/}}.

\begin{question}
  How can one obtain the formula of (i) from the general formula for left adjoints?
  Can this theorem be generalized to $\mcE$- (i.e., not necessarily $\Set$-)valued sheaves?
\end{question}

\begin{corollary}[Limits and colimits in $\Sh(\mcC)$]\label{lem:limits_pointwise}\uses{def:sheaf,lem:lim_in_full_sub, lem:limits_commute_with_limits,  lem:adjoints_commute_with_limits, thm:sheafification, lem:limits_in_functor_categories}
  The category of sheaves on any site $\mcC$ is complete and limits may be calculated pointwise.
  After all,
  limits (in $\PSh(\mcC)$) of sheaves are sheaves and $\Sh(\mcC)$ is a full subcategory of $\PSh(\mcC)$.

  On the other hand,
  $\Sh(\mcC)$ is cocomplete and colimits can be computed as sheafification of the pointwise colimit (as left adjoints such as $\Sh$ preserve colimits).
\end{corollary}

The following is a direct generalization of the definition of accessible presheaves
as presheaves are sheaves for a subcanonical site (namely the empty Grothendieck topology).
\begin{definition}
  Let $(\mcC,\mcT)$ be a subcanonical site.
  Then a sheaf $X\in\Sh(\mcC)$ is an \idx{accessible sheaf} if it is a small colimit of representables.
  The category of accessible sheaves is locally small.
\end{definition}
Frequently, the category of accessible sheaves is a superior replacement to the large category of sheaves.
If $\mcC$ is small, every sheaf is accessible.

\subsection{Elementary and Grothendieck topoi}\label{ssec:el-topoi}

This subsection is aimed at building a better understanding.
It is not necessary for most of what follows.

Categories of sheaves, and especially Grothendieck topoi,
are famous for having excellent categorical properties, which make them,
categorically speaking, \emph{almost as good as sets}.
The main idea of \emph{topos theory} is to take this seriously,
and replace the category of sets by any topos,
see \cite{Johnstone, Lane1992SheavesIG, Borceux1994, Lurie2009}
for some of the main references.
Since the categorical properties of topoi and sets are similar,
one can develop an \emph{internal logic} of a topos,
where one then obtains transfer principles:
an (easy enough) theorem, that can be proven for sets,
is also true internally in the topos.
In \cite{Jamneshan2013, jamneshan2019sheaves, carl2019transfer} topos theory helps in giving a new approach to measure theory.
And in \cite{leroy2013theorielamesuredans}, one is able to circumvent paradoxes arising in measure theory by extending the theory to a topos.
See \cite{Bell, Cunningham1973} for more logical applications.
That the statement be \emph{easy enough} is critical for the transfer.
The better the topos, the more axioms hold internally
and hence the stronger the internal logic,
yielding more classical results.
Our main subject, the category of \emph{condensed sets},
is almost a topos of the best possible sort --
size issues are the only thing holding it back from actually being a topos.
Therefore many arguments from topos theory can be applied there almost verbatim.
The most important (and, in a sense, the only) thing condensed sets do not have that topoi do have is the following.

\begin{definition}[Subobject classifier and exponentiable object]\label{def:subobj_class}\uses{def:subobject}
\begin{enumerate}[(i)]
	\item Consider a category $\mcC$ with pullbacks.
An object $\Omega$ together with a monomorphism $\ast\hookrightarrow\Omega$
is called a \idx{subobject classifier},
if for any monomorphism $f\colon A\hookrightarrow B$
there exists exactly one morphism $\chi_A\colon B\to\Omega$,
the \idx{characteristic function} of \(A\),
such that the diagram
\begin{center}
	\begin{tikzcd}
		A && \ast \\
		\\
		B && \Omega
		\arrow["!", from=1-1, to=1-3]
		\arrow["f"', hook, from=1-1, to=3-1]
		\arrow["{\chi_A}"', from=3-1, to=3-3]
		\arrow[hook, from=1-3, to=3-3]
	\end{tikzcd}
\end{center}
is a pullback diagram.

\item Now assume that $\mcC$ admits finite products.
Then an object $A$ of $\mcC$ is called \idx{exponentiable},
if the endofunctor $B\mapsto B\times A$ admits a right adjoint, denoted $(-)^A$.
Spelling out the universal property, one has a natural isomorphism
\[
    \hom(B\times A, C)\simeq \hom(B, C^A).
\]
The counit of the adjunction $\eps_B\colon B^A\times A\to B$ is called \idx{evaluation}.
A category with finite products is called \idx{cartesian closed}
if all its objects are exponentiable.
A finitely product preserving functor $F$ between cartesian closed categories
is called \textbf{cartesian closed}\index{cartesian closed functor}
if the natural morphism $F(B^A)\to F(B)^{F(A)}$ is an isomorphism.

\item Combining those two points, an \idx{elementary topos}
is a finitely complete cartesian closed category with subobject classifier.
\end{enumerate}
In this case, the domain of the subobject classifier is the terminal object of $\mcC$.
\end{definition}

Next we collect some good properties every elementary topos has.
See \cite{Johnstone} and \cite{Lane1992SheavesIG} for proofs of these results.
First we define some of the properties we observe in the category of sets and that we would like to have.

\begin{definition}[Power object]\label{def:power_objects}\uses{def:special_limits_and_colimits}
For an object $X$, a \idx{power object} is a pair $(PX, \in_X\hookrightarrow X\times PX)$
consisting of an object $PX$ and a subobject of $X\times PX$, called $\in_X$,
such that for every object $B$ and every subobject $R\hookrightarrow B\times X$
there exists a unique map $g\colon B\to PX$ inducing a pullback square
\begin{center}\begin{tikzcd}
	R & {\in_X} \\
	{B\times X} & {PX\times X}
	\arrow[from=1-1, to=1-2]
	\arrow[from=1-1, hook, to=2-1]
	\arrow[from=1-2, to=2-2]
	\arrow["{g\times 1_X}"', from=2-1, to=2-2]
\end{tikzcd}\end{center}
\end{definition}

\begin{definition}[Disjointness of coproducts]\label{def:coprod_disj}\uses{def:special_limits_and_colimits}
We say that a coproduct \(\coprod X_i\) is \textbf{disjoint}\index{disjoint coproduct},
    if for any pair \(X_i,X_j\) of components of the coproduct,
    the pullback
    \begin{center}
        \begin{tikzcd}
            	P & {X_j} \\
            	{X_i} & {\coprod X_i}
            	\arrow[from=1-1, to=1-2]
            	\arrow[from=1-1, to=2-1]
            	\arrow["{e_j}", from=1-2, to=2-2]
            	\arrow["{e_i}"', from=2-1, to=2-2]
        \end{tikzcd}
    \end{center}
is the initial object.
We call two subobjects \textbf{disjoint}\index{disjoint subobjects}
if their intersection (the pullback) is the initial object.
\end{definition}

\begin{definition}[(Co)image]\label{def:image-coimage}
Let $\mcC$ be a category.
\begin{enumerate}[(i)]
    \item Define the \idx{image} of a morphism $f\colon A\to B$ as the subobject $m\colon M\hookrightarrow B$
        such that $f$ factors over $m$ and for any $h\colon N\hookrightarrow B$ such that $f$ factors over $h$,
        $m$ factors over $h$ as well.
    \item The \idx{coimage} is defined dually: it is an epimorphism $e\colon A\twoheadrightarrow M$ such that $f$ factors over $e$ and for any $f=iq$ with $q$ epic there is another epic $w$ such that $e=wq$.
\end{enumerate}
\end{definition}

\begin{remark}\label{rem:image-epi}
By definition, the image and the coimage of a morphism are,
if they exist, unique up to isomorphism.
Every monomorphism is its own image and every epimorphism is its own coimage.

If the category $\mcC$ has equalizers,
then the morphism $e$ in the image factorisation $f=me$ is an epimorphism:
for if $g,h\colon M\to C$ are so that $ge = he$,
then taking the equalizer $E$ of $g$ and $h$ yields a subobject of $B$
over which $f$ factors.
But $E$ is also a subobject of $M$, so $E\simeq M$.
In particular, $g=h$.

Dually, if $\mcC$ has coequalizers, then the morphism $m$ in the coimage factorisation $f=me$ is a monomorphism.
\end{remark}

\begin{definition}\label{def:stable-image}
Let $\mcC$ be a category such that every morphism has an image and pullbacks exist.
We say that \idx{images are stable under pullback} if for every morphism $f\colon A\to B$
with image $m\colon M\to B$ and every morphism $g\colon C\to B$,
the morphism $m^*\colon M\times_B C\to C$ is the image of $f^*\colon A\times_B C\to C$.
\begin{center}
\begin{tikzcd}
	{A\times_B C} & {M\times_B C} & C \\
	A & M & B
	\arrow[from=1-1, to=1-2]
	\arrow["{{{f^*}}}", curve={height=-12pt}, from=1-1, to=1-3]
	\arrow["{{{g^*}}}"', from=1-1, to=2-1]
	\arrow["{{{m^*}}}"', hook, from=1-2, to=1-3]
	\arrow[from=1-2, to=2-2]
	\arrow["g", from=1-3, to=2-3]
	\arrow[from=2-1, to=2-2]
	\arrow["f"', curve={height=12pt}, from=2-1, to=2-3]
	\arrow["m", hook, from=2-2, to=2-3]
\end{tikzcd}
\end{center}
\end{definition}

The following is an extension of the notion of kernel pair established in definition \ref{def:kernel_pair}.

\begin{definition}[Equivalence relation]\label{def:equivalence_relation}\uses{def:kernel_pair}
Let $\mcC$ be a finitely complete category.
An \idx{equivalence relation} is a pair of parallel morphisms
$\partial_0,\partial_1\colon R\to A$ such that it is
\begin{enumerate}
		\item a \idx{relation}, i.e., the induced arrow $(\partial_0,\partial_1)\colon R\to A\times A$ is monic,
		\item \idx{reflexive}, i.e., there exists $s\colon A\to R$ with $\partial_0s=\partial_1s=1_A$.
		\item \idx{symmetric}, i.e., there exists $r\colon R\to R$ with $\partial_0r=\partial_1$ and $\partial_1r=\partial_0$
		\item \idx{transitive}, i.e., there exists $t\colon R\times_A R\to R$
		for the pullback
\begin{center}\begin{tikzcd}
	{R\times_AR} & R \\
	R & A
	\arrow["p", from=1-1, to=1-2]
	\arrow["q"', from=1-1, to=2-1]
	\arrow["a", from=1-2, to=2-2]
	\arrow["b"', from=2-1, to=2-2]
\end{tikzcd}\end{center}
fulfilling $\partial_0 t=\partial_0q$ and $\partial_1t=\partial_1p$.
\end{enumerate}
The kernel pair of a morphism $f$ is an equivalence relation
and we call an equivalence relation arising this way \idx{effective}.
We say that a category has \idx{effective equivalence relations} if every equivalence relation is an effective.

Given an equivalence relation on $A$, the \idx{quotient by the equivalence relation}
$A/R$ is the coequalizer of the parallel arrows $\partial_0,\partial_1$.
\end{definition}

Recall that a category is called \emph{balanced}
if every monomorphism that is also an epimorphism is already an isomorphism,
see definition~\ref{def:balanced-category}.
With these definitions at hand, we can collect all of the special properties of elementary topoi.

\begin{theorem}[Properties of elementary topoi]\label{thm:elem_props_topos}\uses{def:equivalence_relation, def:special_monoepi, def:image-coimage}
Let $\mcC$ be an elementary topos.
Concerning objects, the following properties hold.
\begin{enumerate}[(i)]
    \item There exist {power objects}, given by $PX=\Omega^X$
     and $\in_X\hookrightarrow PX\times X$ is induced by the evaluation $PX\times X\to\Omega$.
    \item Conversely, a finitely complete cartesian closed category with power objects is an elementary topos.
    \item All power objects $PX$ are injective.
        In particular, there exist enough injectives.
 	\item It is finitely cocomplete.
    \item Coproducts are disjoint.
    \item The initial object is \idx{strict},
    meaning that any morphism with codomain being the initial object is an isomorphism.
\end{enumerate}
Concerning morphisms, we furthermore have the following properties.
\begin{enumerate}[(i)]
\setcounter{enumi}{6}
    \item The category is balanced.
    \item Every epimorphism is effective.
	    The pullback of an epimorphism is again epimorphic
        and the pullback square then also is a pushout square.
	\item Dually, every monomorphism is effective.
	    Moreover, the pushout of a monomorphism is a monomorphism
	    and the pushout diagram is a pullback diagram as well.
    \item If \(A_i\hookrightarrow B\) are monic and pairwise disjoint,
    then \(\coprod A_i\hookrightarrow B\) is monic as well.
    \item If \(A_i\hookrightarrow B_i\) are monic, then \(\coprod A_i\hookrightarrow\coprod B_i\) is monic as well.
    \item If $f\colon A\to B$ and $g\colon C\to D$ are epimorphisms,
        then so is $f\times g\colon A\times C\to B\times D$.
    \item Every equivalence relation if effective.
    \item Every morphism $f$ has an image which is given by the equalizer of the cokernel pair of $f$.
    Moreover, images are stable under pullback.
    Dually, every morphism $f$ has a coimage which is given by the coequalizer of the kernel pair of $f$.
    In particular, every morphism $f\colon A\to B$ can be monic-epic factorized as
        \[
            A\twoheadrightarrow M\hookrightarrow B,
        \]
    and this monic-epic factorization is unique up to isomorphism.
\end{enumerate}
\end{theorem}

\begin{proof}
All of the claims (except the pushout claim in (viii)) can be found in chapter IV of \cite{Lane1992SheavesIG},
in chapter 5 of \cite{Borceux1994} and throughout \cite{Johnstone}.

To illustrate how to work with the subobject classifier, we spell out proofs of two of the assertions.

(vii)
Let $f\colon A\to B$ be a bimorphism.
Consider the pullback diagram of the characteristic function
\begin{center}\begin{tikzcd}
	A & \ast \\
	B & \Omega
	\arrow[from=1-1, to=1-2]
	\arrow[hook', two heads, from=1-1, to=2-1]
	\arrow[hook', from=1-2, to=2-2]
	\arrow["\chi_A"', from=2-1, to=2-2]
\end{tikzcd}.\end{center}
Now note that, for $h\colon B\to\ast\to \Omega$ we have that $hf = \chi_Af$
since the diagram
\begin{center}\begin{tikzcd}
	A & \ast \\
	B & \Omega
	\arrow[from=1-1, to=1-2]
	\arrow[hook', two heads, from=1-1, to=2-1]
	\arrow[hook', from=1-2, to=2-2]
	\arrow[from=2-1, to=1-2]
	\arrow[from=2-1, to=2-2]
\end{tikzcd}\end{center}
is commutative.
But now the outer rectangle is a pullback diagram,
which means that $f$ is the equalizer of $h$ and $\chi_A$.
Because $f$ is epic we conclude $h=\chi_A$.
But the equalizer of two identical arrows is an isomorphism, so $f$ is an isomorphism.

(viii)
For the pushout property,
use that coequalizers are stable under pullback.
This follows from topoi being cartesian closed.
Now write both epics in the pullback as coequalizers of arrows fitting into a pullback diagram.
Then the existence of an arrow out of the bottom right vertex follows from the universal property of coequalizers
and uniqueness comes for free since there are enough epimorphisms in the diagram.

(xiv)
Consider the equalizer $m\colon M\hookrightarrow B$ of the cokernel pair $(r_1,r_2)$ of $f$.
Since $r_1\circ f=r_2\circ f$ we know that $f$ factors through $m$, $f=me$ for some $e\colon A\to M$.
For any other factorisation $f=hg$ with monic $h$, the monic is equalizer of its cokernel pair $s,t\colon B\to C$ by (ix).
\begin{center}\begin{tikzcd}
	A & M & B & {B\sqcup_AB} \\
	A & N & B & C
	\arrow["e"', from=1-1, to=1-2]
	\arrow["f", curve={height=-12pt}, from=1-1, to=1-3]
	\arrow["m"', hook, from=1-2, to=1-3]
	\arrow["{r_1}", shift left, from=1-3, to=1-4]
	\arrow["{r_2}"', shift right, from=1-3, to=1-4]
	\arrow["1", tail reversed, from=2-1, to=1-1]
	\arrow["g", from=2-1, to=2-2]
	\arrow["f"', curve={height=12pt}, from=2-1, to=2-3]
	\arrow["h", hook, from=2-2, to=2-3]
	\arrow["1", tail reversed, from=2-3, to=1-3]
	\arrow["s", shift left, from=2-3, to=2-4]
	\arrow["t"', shift right, from=2-3, to=2-4]
\end{tikzcd}\end{center}
As $s,t$ equalize $h$ and hence $f$, by the pushout property of $B\sqcup_A B$ there exists $u$ with
\begin{center}\begin{tikzcd}
	A & M & B & {B\sqcup_AB} \\
	A & N & B & C
	\arrow["e"', from=1-1, to=1-2]
	\arrow["f", curve={height=-12pt}, from=1-1, to=1-3]
	\arrow["m"', hook, from=1-2, to=1-3]
	\arrow["{r_1}", shift left, from=1-3, to=1-4]
	\arrow["{r_2}"', shift right, from=1-3, to=1-4]
	\arrow["u", from=1-4, to=2-4]
	\arrow["1", tail reversed, from=2-1, to=1-1]
	\arrow["g", from=2-1, to=2-2]
	\arrow["f"', curve={height=12pt}, from=2-1, to=2-3]
	\arrow["h", hook, from=2-2, to=2-3]
	\arrow["1", tail reversed, from=2-3, to=1-3]
	\arrow["s", shift left, from=2-3, to=2-4]
	\arrow["t"', shift right, from=2-3, to=2-4]
\end{tikzcd}\end{center}
This yields $sm=tm$ and thereby a $k\colon M\to N$ with $hk=m$
\begin{center}\begin{tikzcd}
	A & M & B & {B\sqcup_AB} \\
	A & N & B & C
	\arrow["e"', from=1-1, to=1-2]
	\arrow["f", curve={height=-12pt}, from=1-1, to=1-3]
	\arrow["m"', hook, from=1-2, to=1-3]
	\arrow["k"', from=1-2, to=2-2]
	\arrow["{r_1}", shift left, from=1-3, to=1-4]
	\arrow["{r_2}"', shift right, from=1-3, to=1-4]
	\arrow["u", from=1-4, to=2-4]
	\arrow["1", tail reversed, from=2-1, to=1-1]
	\arrow["g", from=2-1, to=2-2]
	\arrow["f"', curve={height=12pt}, from=2-1, to=2-3]
	\arrow["h", hook, from=2-2, to=2-3]
	\arrow["1", tail reversed, from=2-3, to=1-3]
	\arrow["s", shift left, from=2-3, to=2-4]
	\arrow["t"', shift right, from=2-3, to=2-4]
\end{tikzcd}.\end{center}
Thus the image has the universal property.
The morphism $e$ is epic by remark~\ref{rem:image-epi}.

For the stability under pullbacks, consider the following diagram
\begin{center}
\begin{tikzcd}
	{A\times_B C} & {M\times_B C} & C & {C\sqcup_{A\times_B C}C} \\
	A & M & B & {B\sqcup_AB}
	\arrow[from=1-1, to=1-2]
	\arrow["{{{f^*}}}", curve={height=-12pt}, from=1-1, to=1-3]
	\arrow["{{{g^*}}}"', from=1-1, to=2-1]
	\arrow["{{{m^*}}}"', hook, from=1-2, to=1-3]
	\arrow[from=1-2, to=2-2]
	\arrow["{{s_1}}", shift left, from=1-3, to=1-4]
	\arrow["{{s_2}}"', shift right, from=1-3, to=1-4]
	\arrow["g", from=1-3, to=2-3]
	\arrow["h", dashed, from=1-4, to=2-4]
	\arrow[from=2-1, to=2-2]
	\arrow["f"', curve={height=12pt}, from=2-1, to=2-3]
	\arrow["m", hook, from=2-2, to=2-3]
	\arrow["{{r_1}}", shift left, from=2-3, to=2-4]
	\arrow["{{r_2}}"', shift right, from=2-3, to=2-4]
\end{tikzcd}
\end{center}
where $(r_1,r_2)$ is the cokernel pair of $f$, $(s_1,s_2)$ is the cokernel pair of $f^*$
and $h$ is the morphism from the universal property of the pushout
\begin{center}
\begin{tikzcd}
	{A\times_B C} & C \\
	C & {C\sqcup_{A\times_B C}C} & B \\
	& B & {B\sqcup_A B}
	\arrow["{{f^*}}", from=1-1, to=1-2]
	\arrow["{{f^*}}"', from=1-1, to=2-1]
	\arrow["{{s_2}}"', from=1-2, to=2-2]
	\arrow["g", from=1-2, to=2-3]
	\arrow["{{s_1}}", from=2-1, to=2-2]
	\arrow["g"', from=2-1, to=3-2]
	\arrow["h"{description}, dashed, from=2-2, to=3-3]
	\arrow["{{r_2}}", from=2-3, to=3-3]
	\arrow["{{r_1}}"', from=3-2, to=3-3]
\end{tikzcd}
\end{center}
since $r_1gf^* = r_1f^*g = r_2f^*g = r_2gf^*$.
For any morphism $k\colon D\to C$ with $s_1k = s_2k$
there exists a unique $\ell\colon D\to M$ such that $gk = m\ell$
because $r_1gk = hs_1k = hs_2k = r_2gk$.
This implies the existence of a unique $D\to M\times_B C$ which is the desired arrow,
such that $m^*\colon M\times_B C\to C$ actually is the image of $f^*$.

The proof for the coimage follows the same argument essentially,
but we use that every epimorphism is the coequalizer of its kernel pair by (viii).

Now let $f = me = m'e'$ be two factorizations of $f$ over $M$ and $M'$
where $m$, $m'$ are monic and $e$, $e'$ are epic.
We can assume that $m$ is the image of $f$ as constructed before.
So there exists a morphism $s\colon M\to M'$ making the diagram
\begin{center}
\begin{tikzcd}
	& M \\
	A && B \\
	& {M'}
	\arrow["m", hook, from=1-2, to=2-3]
	\arrow["s"', dashed, from=1-2, to=3-2]
	\arrow["e", two heads, from=2-1, to=1-2]
	\arrow["{e'}"', two heads, from=2-1, to=3-2]
	\arrow["{m'}"', hook, from=3-2, to=2-3]
\end{tikzcd}
\end{center}
commutative.
Uniqueness of $s$ follows from $m'$ being a monomorphism.
Now $s$ is both a monomorphism as well as an epimorphism
(since $e'$ is epic and $m$ is monic).
Now (vii) implies that $s$ is an isomorphism.
In particular, the codomain of the coimage is isomorphic to the domain of the image.

\end{proof}

Grothendieck topoi, as defined in \ref{def:grothendieck_topos},
are a special kind of topos.
Remark that it is not trivial that Grothendieck topoi are indeed elementary topoi,
see also chapter III.7 in \cite{Lane1992SheavesIG}.
But in fact, Grothendieck topoi have way stronger properties.
Surprisingly, it is possible to extract those properties singling out Grothendieck topoi from the elementary topoi.
These are the so-called \idx{Giraud Axioms}.

Recall that colimits are called \emph{stable under base change}, definition~\ref{def:col_stable_under_base_change},
if for every diagram $D\colon\mcI\to\mcC$ and morphism $f\colon X\to \varinjlim_i D$
the induced pullback diagram $D(i)\times_{\varinjlim_i D} X$ yields a colimiting cocone for $X$.

A class $\mcA$ of objects of $\mcC$ is called \emph{separating}
if, for any two distinct morphisms $f,g\colon B\to C$, there exists an object $A$ in $\mcA$
and a morphism $h\colon A\to B$ such that $fh \neq gh$.
It is called \emph{densely generating}
if every object of $\mcC$ is a colimit of objects in $\mcA$, see also definition~\ref{def:cogenerating_set}.

\begin{theorem}[Giraud axioms]\label{thm:gireaud_characterisation}\uses{def:grothendieck_topos}
A category $\mcC$ is a Grothendieck topos if and only if
\begin{enumerate}[(i)]
    \item it is finitely complete,
    \item it has (small) coproducts,
            and they are disjoint and stable under base change,
    \item coequalizers are stable under base change,
    \item it has effective equivalence relations quotients of them,
    \item every epimorphism is effective, and
    \item it admits a \emph{small} separating set.
\end{enumerate}
\end{theorem}

We are not going to prove this theorem.
Note that, however, it is easy to see that all the conditions are necessary for $\mcC$ to be a Grothendieck topos:
they all are true in $\Set$,
and so they hold in any presheaf category $\PSh(\mcE)$ for a small category $\mcE$,
since there (co)limits are computed pointwise
and the representable presheaves separate $\PSh(\mcE)$.
As sheafification is cocontinuous and preserves finite limits,
all of the properties also hold in the sheaf category $\Sh(\mcE)$.
For a more detailed argument see the proof of proposition~\ref{prop:cond-morphisms}.
The proof of the theorem can be found in
the Appendix of \cite{Lane1992SheavesIG},
Theorem C2.2.8 in \cite{Johnstone},
Theorem 3.6.1 in \cite{Borceux1994},
or 6.1.0.1 in \cite{Lurie2009}.

\begin{remark}
Another formulation of Giraud's axioms goes as follows:
a category $\mcC$ is a Grothendieck topos precisely if
\begin{enumerate}[(i)]
    \item it is (small) cocomplete,
    \item colimits are stable under base change,
    \item coproducts are disjoint,
    \item it has effective equivalence relations, and
	\item it has a \emph{small} densely generating set $S$ of objects
	    and there exists a regular cardinal $\lambda$,
	    such that all objects of $S$ are $\lambda$-compact,
	    and they generate $\mcC$ under $\lambda$-filtered colimits.
\end{enumerate}
If a cocomplete category satisfies the last point,
it is also called \idx{locally presentable}.
\end{remark}

For more about topos theory,
especially for the topic of geometric morphisms, we refer to the literature.

\clearpage{\thispagestyle{empty}\cleardoublepage}

\chapter{Condensed Topology}
\quot{Point set topology is a disease from which the human race will soon recover.}
{Henri Poincar{\'e}, quoted in Chapter 3 of \cite{Eisner2016}.}

In this chapter we introduce condensed sets,
the basic notion on which condensed mathematics is built, as developed by Clausen and Scholze in the series of lecture notes \cite{scholze2019Analytic, scholze2019condensed, scholze2022complex, Scholze2023} as well as two series of lectures \cite{Scholze2024}%
\footnote{\url{https://www.youtube.com/playlist?list=PLAMniZX5MiiLXPrD4mpZ-O9oiwhev-5Uq}}
and \cite{ScholzeNov2020}%
\footnote{\url{https://www.youtube.com/playlist?list=PLx5f8IelFRgGmu6gmL-Kf_Rl_6Mm7juZO}}.
Very similar, though not the same concepts
have also been introduced by Barwick and Haine in \cite{barwick2019pyknotic}

Further work on condensed mathematics can be found, e.g., in \cite{asgeirsson2021foundations, Commelin, mair2021animated, Zhu2023, Aparicio2021, Morgan, Rasekh2023, yamazaki2022condensed, tong2021topologization, Link2023, Commelin2023, gregoric2024stone, bhatt2014proetale, Lurie2019, Garcia2022, LeStum2024} and some rather advanced applications in \cite{fargues2024geometrization, camargo2024analytic, mann2022padic}.

In the first section we lay out our rough intuitive understanding of condensed sets (it is, indeed, devoid of formal content).
In the next section we will intensively study profinite sets,
before we move towards $\kappa$-condensed sets
with $\kappa$ being some size restriction on (the \enquote{complexity} of) the condensed set.
Afterward, we will get rid of the $\kappa$ and thereby obtain the definition of condensed sets.
Then, we will explain basic properties of the category of condensed sets, discuss the embedding of topological spaces,
and develop some basic topological notions.
Many of the results obtained there will be proven and studied again in more general versions later in chapter~4.
However, we have decided to provide some elementary proofs,
which hopefully helps in gaining intuition.

\section{A very rough sketch}
\quot{Categorification is the left adjoint to forgetting category theory.}
{Sophia}

In this chapter we try to convince the functional analytically minded person that condensed mathematics is worth spending time on.
We mainly focus on two questions.
\begin{itemize}
  \item What are problems with topology that we would love to eliminate/resolve?
  \item What are the main ideas of condensed mathematics relevant for this purpose, that is, how are these problems resolved in condensed mathematics?
\end{itemize}
Of course, this chapter is \emph{highly subjective} and we just present what has convinced us.
In particular, the term \enquote{main ideas} should not be taken too seriously, as we merely understand a fraction of the theory that Clausen and Scholze have developed.

\subsection{Problems with topology}\label{subsec:problems_in_top}

\quot{It is never worth a first-class man's time to express a majority opinion.
By definition, there are plenty of others to do that.}{G.H. Hardy, quoted by C.P. Snow}

We have the feeling that many of the problems with topology arise from the fact that the notion of topological space is not in essence intrinsically motivated.
In other words, in our minds,
the cardinal sin of topology is faulty generalisation.
One observes the behaviour of \enquote{open sets} in nice metric spaces and axiomatizes this,
not necessarily considering at first
how well this is suited to other (say, non-metrisable) \enquote{spaces} that one would like to consider.%
\footnote{
  Of course,
  all of this is, historically speaking, unfair, as many of these interesting spaces would never have been discovered without the concept of topological space.
  But as mathematicians,
  we want the best, most reasonable theory, not one riddled with quirks that can only be justified by explaining their history.
}
Specifically, $\R$ with its euclidean topology is seen as the most natural object in the world (nothing but a line, after all) and thus the phenomena of $\R$ dictate the path of generalisation.
The first step takes us to $\R^n$ with the real-world intuition of balls and provides a feeling for Banach spaces;
then one proceeds further to metric spaces and topological spaces, where every step makes total sense as long as one has the image of open sets as (generated by open) balls in mind.
This \enquote{top-down} approach is certainly useful for applied mathematics; but it causes problems whenever we want to combine
topological methods with algebraic structures, which are defined in a bottom-up way:
in algebra, special emphasis is placed not only on the definitions generalising day-to-day examples as special cases but often more on the \enquote{naturality} of the definitions themselves
(seeing them from a \enquote{Platonistic} perspective),
leading to such pristine definitions as those of groups or modules and giving excellent general theories.

One aspect of this problem is that topology is heavily set-theory based (i.e., uses quintessentially notions such as intersection or union, even the powerset in the very definition),
whereas for algebra, only points (classically) or homomorphisms (categorically) are needed.
In particular, in the non-Hausdorff case (or, even more drastically, in the presence of non-closed points, i.e., non-T1 spaces),
topologically one cannot suitably distinguish points.
If one cannot distinguish points with continuous functions, or, in other words, there are not enough continuous functions, this implies that the spaces behave badly with category theory.
As a concrete example,
in the category of topological spaces,
the two morphisms $\ast\to\{0,1\}_{\mathrm{indisc.}}$ are actually indistinguishable by the opens of $\{0,1\}_{\mathrm{indisc.}}$!
This implies that such spaces behave bad with algebraic structures, where single points are of great importance and which are more category-theory based.
For example, many algebraic theorems for abelian groups do not hold in topological abelian groups.
For a specific example of \enquote{mono $+$ epi $=$ iso} failing,
take $\R_{\textrm{disc.}}\to\R$.

Put as a slogan: topology mixes badly with algebra.
But also without additional algebraic structure, there are few non-trivial statements about topological spaces that do not assume
some separation (or compactness) condition, leading to a whole zoo of different separation axioms.%
\footnote{By skimming through this section,
  one might object that for condensed sets,
  there are also different separation conditions.
  However,
  there are two important differences.
  First of all,
  we claim that the notions we discuss are more conceptually clear and differentiable (and in particular, just fewer).
  More importantly,
  these notions only come up seldomly when we actually work with condensed sets!
  It is not at all the case that in chapter~4,
  practically every theorem or proposition starts with a bag of assumptions on the condensed sets involved in the statements.
}
This creates the feeling that there are some topological spaces without any good properties,
i.e., spaces that should not exist -- both in the sense that they form pathological counterexamples to perfectly reasonable theorems
and in the sense that they do not serve any purpose in the theory.
So on the one hand, there seem to be too many bad spaces in the non-Hausdorff world -- some of them should not exist.

On the other hand, the necessity of working in Hausdorff spaces causes another problem:
One can only factor by closed subspaces.%
\footnote{Indeed,
  in the category of Hausdorff spaces,
  the pushout of $\ast\leftarrow A\hookrightarrow B$ is (automatically) $B/\overline{A}$.}
Hence, whenever one has, e.g., a subspace in a Banach space one wants to quotient away, one has to take its closure first,
but this often destroys good algebraic descriptions of the subspace.
This causes troubles in particular in classical problems in spectral theory, where this \enquote{taking closures} regularly annihilates the relevant information.
Another simple example is that of a topological group action $G\times X\to X$.
Whenever one wants to consider the orbit space, it is common to take the quotient by the closure of the (algebraic, point-set) orbit
    \[
            \mathrm{orb}(x) = \overline{\{g.x:g\in G\}}.
      \]
But if, e.g., the action has dense orbits, this means that the orbit space becomes trivial, although the action might be highly non-transitive (consider an irrational rotation on the torus).

In order to obtain a good theory that can handle these cases, we need to be able to take quotients by non-closed subspaces, leading to a new variant of non-Hausdorff behaviour.
Spelling this out in some examples, we have the embedding of Banach spaces $\ell^1\hookrightarrow \ell^2$.
This is injective, i.e., a monomorphism, but clearly no isomorphism of Banach spaces.
Hence we would like to have a cokernel $\ell^2/\ell^1$, which has to be highly non-trivial, as the topologies on $\ell^1$ and $\ell^2$ are quite different.
More drastically, consider the monomorphism $\R_{d}=\R_{\mathrm{disc}}\hookrightarrow \R$ (the latter being equipped with the usual euclidean topology).
We need a cokernel, which has to be a space like $\R/\R_{d}$, i.e., having as underlying set just a point.
But this space has to be highly non-trivial!
This also explains why the classical topological quotient space cannot yield the correct space -- there is just one topology on a point.
We regard this as a defect of using set theory: if the complexity of the space is bounded by the power set of the underlying set,
we cannot build the complicated factor spaces we need.

This leads to the feeling that highly non-separated spaces are \enquote{incorrectly} defined and that one should prefer a category-theoretic definition over a set-theoretic one.

Summarising the above, some problems we want to resolve in topology are the following.
\begin{itemize}
  \item Some topological spaces should exist, but do not. (Think: $\R/\R_{d}$.)
  \item Some topological spaces should not exist, but do so. (Think: indiscrete spaces.)
  \item Algebraic theory does not work well combined with topological structure as soon as one does not work in compact Hausdorff spaces.
  \item Favourable categorical aspects of algebra get lost when working with topological spaces. (Think: $A\to B$ iso iff monic and epic for abelian groups.)
\end{itemize}

\subsection{How do condensed sets work?}
\quot{The mathematician's patterns, like the painter's or the poet's,
  must be beautiful; the ideas like the colours or the words, must fit together in a harmonious way.
  Beauty is the first test: there is no permanent place in the world for ugly mathematics.}
{G.H. Hardy in A Mathematicians Apology (1941)}

In order to approach the problems mentioned above, condensed mathematics was conceived as a wonderful combination of well-known ideas in new shape.

The name \enquote{condensed set} describes perfectly the \enquote{bottom-up} route we will take,
starting with the simplest possible spaces and step-by-step building \enquote{bigger and more complex} spaces,
always aiming to take the simplest, most canonical steps.
We start with the \enquote{discrete} world of finite sets.
Finite sets come naturally equipped with just one reasonable topology (recall that non-T1 spaces in our opinion are not reasonable) -- the discrete topology.
This is the first thing we can all agree on,
giving us a shared base of \enquote{molecules} that so far do not interact or \enquote{get close to one another}.
Now we ask ourselves what we can naturally build out of just those finite discrete sets.

Which spaces have their topology fully determined by finite sets?
The answer is rather unsurprising:
The universal property of products and, more generally, closed subspaces of products
force some spaces to carry a specific topology.
In fact, take an arbitrarily big product of finite spaces $S_i$,
\[\prod S_i.\]
Then the universal property of product spaces already fully determines the topological properties that $\prod S_i$ must have,
by just using that we know what the correct topology on finite sets is!
(For example,
a map from any other space $X\to\prod S_{i}$ is continuous if and only if all composites $X\to\prod S_{i}\to S_{i}$ are continuous.
And conversely, such continuous maps $X\to S_{i}$ uniquely determine one map $X\to\prod S_{i}$.)

And knowing that we do not have any choice in the topology of $\prod S_i$, we can take any closed subspace of this, and
rewrite the universal property -- again this shows that there is no choice in the topology.
(For those who want to object that using the notion \enquote{closed subspace} here is somewhat self-referential,
note that closed subspaces of products are precisely those that can be delineated by equations of their coordinates,
i.e., this class of subsets indeed allows for an intrinsic description.)

This means that having decided that the correct topology on finite sets has to be the discrete topology,
we already know what the topology on sets of the form $\prod S_i$ and closed subspaces thereof must look like.

Let us have a look at the topological spaces we have created in this way.
Clearly, by Tychonov, arbitrary products of finite discrete sets are compact Hausdorff spaces,
and since all the coordinates are discrete, they form totally disconnected spaces (or to be precise, $0$-dimensional spaces).
This also passes to closed subspaces, hence we conclude that all the topological spaces we obtain by limits of finite sets
are compact Hausdorff totally disconnected spaces (for short: Stone spaces).
It is now one part of the classical Stone duality (theorem~\ref{thm:stone}) that conversely every
totally disconnected compact Hausdorff space is a limit of finite sets, hence the topological spaces we have created are precisely those.

But how does this procedure work formally, without ever referring to topological spaces and just using the category $\fin$ of finite sets?

Categorically, products of finite sets and closed subspaces of these are just limits of finite sets, in fact all limits of finite sets.
But they are not finite anymore, hence not contained in the category of finite sets.
We rather want these limits to be adjoined freely to the category of finite sets.
(The \enquote{freely} is important here.
We do not want to just take the limit in $\Set$ which could only ever give discrete objects.)
This suggests the idea of passing from $\fin$ to the free completion of $\fin$.
But this is not quite what we want -- we still want that finite limits of finite sets are preserved,
as these clearly should be given by taking the limit in finite sets and equipping this with the discrete topology.
Therefore, we should rather look at the finitely continuous completion of $\fin$.
This is precisely the pro-completion (see section~\ref{ssec:completions} (co)completions).
Hence, we pass from finite sets purely via abstract nonsense (see again section~\ref{ssec:completions} for a purely categorical construction of this)
to the free pro-completion of finite sets, $\profin$.

It is now a classical version of Stone duality that these constructions agree.
(We will prove this in the next section.)

Returning to our process of condensation, we started with finite discrete sets $\fin$, and by pure abstract nonsense obtained the category $\profin$ of
totally disconnected compact Hausdorff spaces.
So there was already some magic involved -- the spaces are no longer discrete, but rather \enquote{gaseous} and slightly non-discrete.
The topology is still totally disconnected and hence quite close to a discrete topology (we have many clopens),
hence these spaces like $\beta\N$ or $\alpha\N$ of course look quite fuzzy and gaseous,
but clearly this is the closest to being discrete that a non-discrete space can be
(this is not entirely true, even closer are the extremally disconnected compact Hausdorff spaces $\extr$, which to us will be at least as important as the profinite sets).

Having unlocked this first achievement of things \enquote{binding together},
how can we continue with the process of condensation?

The next thing we want to do with these \enquote{gaseous} profinite sets
(or extremally disconnected sets $\extr$) is to condense the gaseous objects to liquid things such as $\R$ (this to us feels liquid -- there are no gaps and everything flows around)
or, e.g., a connected thing such as $[0,1]$ (this even feels solid as nothing can flow away).
This means we want connected spaces in the form of \enquote{factors} or \enquote{gluings}, and want \enquote{big spaces} in the form of unions.
As an example, the unit interval can be seen as a factor of the profinite set of decimal representations $\prod_{i\ge 0}\{0,\dots, 9\}$
by factoring out the equivalence relation of identifying sequences whose endings are of the form $d9999\dots$ ($d\in\{0,1,\ldots,8\}$)
with the corresponding sequence ending in $(d+1)0000\dots$.
And $\R$ can be built by gluing copies of $[0,1]$ together
(or alternatively can be given by using appropriate decimal representations with Laurent series, which will be important for the liquid analytic ring structure).
But just like we went from $\fin$ to $\profin$ by adjoining limits, disjoint unions and factor maps are just all colimits.
Hence the naive way to define condensed sets would be to take the free cocompletion of $\profin$.
However, again we want to keep some colimits that we are already happy with -- the finite coproducts and the factors by those closed equivalence relations such that the factor is profinite.
The latter precisely coincide with images of morphisms, or coequalizers of kernel pairs.
Note that we do not want to preserve \emph{all} finite colimits.
E.g., the colimit of the irrational torus rotation should be something like $\R/\Q$ (but somehow not with the indiscrete topology) but surely not just a point.
And analogously,
the colimit of equivalence relations that produce $[0,1]$ from $\beta\N$ should not be forced to yield a profinite set,
but rather indeed result in the unit interval.
Hence, the canonical choice is to define the category of condensed sets, $\cond$,
as the universal cocompletion that preserves precisely those colimits, again using the purely categorical construction of this category.
This in fact is a possible definition of $\cond$, as we will see later.

In fact we recover a wonderful category where (essentially by definition) all the limits and colimits behave as well as we want them to.
But nevertheless, we will be able to embed (fully faithfully) the \enquote{correct} topological spaces (we have long forgotten the ones that we do not want anymore),
and will be able to do a lot of classical theory using this category in place of topological spaces. --
It is not a discrete theory!
In particular, we will develop notions of compactness, separation, and later also completeness, allowing many of the (useful, conceptual) \enquote{topological} considerations
(while not needing to think about the T$3\frac{1}{2}$ condition ever again).

This category still contains objects that behave like $\R$ or $\R P^{2}$ but also strange factors like $\R/\R_d$.
So, what have we won compared to the topological situation?
Note that in the last subsection we noted that non-Hausdorff quotients like $\R/\R_d$ should not be done as topological spaces.
Hence it is good news that here we will not just recover the usual topological notion of factor spaces and coproducts;
and in fact we will not even be able to have bad (highly non-separated) topological spaces in our theory.
In other words, in good, separated situations both worlds completely coincide, but the non-Hausdorff situation is handled completely differently.
In fact, condensed sets in general can not be identified with topological spaces (and their morphisms not necessarily, but often, with maps of sets),
instead forming a far more flexible theory,
e.g., by allowing \enquote{complicated structures on a single point}.
Nevertheless, as we will see later, they are still (almost) accessible by classical topological spaces and topological notions.

Of course, we need a more concrete description of condensed sets.
In fact, the intuition described in the following can also alternatively be thought of as the motivation of the definition of \enquote{condensed sets}.
The crucial idea is that of \enquote{testing spaces}.
We have, inside $\cond$, a class of objects that we understand well, the profinite sets.
Now any other object should be defined via its interaction with respect to the class of \enquote{good} objects.
I.e., a condensed set $T$ assigns to every profinite set $S$ the \enquote{set of continuous functions} $T(S)$ from $S$ to $T$.
(Indeed, it is true that $T(S)=\hom_{\cond}(S,T)$ for any condensed $T$ and profinite $S$.)
Hence for the moment think of $T(S)$ as $T(S)=C(S,T)$.
(For condensed sets emerging from topological spaces, this will consequently be simply a true equation.)
For example,
\begin{itemize}
  \item $T(\ast)$ consists of the continuous functions from a point to $T$, hence the \idx{underlying set}.
        (And $T(\ast\sqcup\ast)$ are pairs of points etc.)
  \item $T(\alpha\N)$ are the continuous functions from $\alpha\N$ to $T$, hence the \textbf{convergent sequences in $T$}.
  \item $T(\beta\N)$ are the relatively compact sequences in $T$,
  \item $\dots$
\end{itemize}
Generally, one can imagine the convergent nets/filters of a special form to be continuous functions from appropriate compactifications to $T$,
and hence $T$ being the mapping that assigns to every form the convergent nets of that form in $T$.
This clearly should fix a topology, and fits well into the picture of how one would like to define topologies (or generally types of convergence).

Another crucial point is then, for any $S\to S'$ in $\prof$ a map $T(S')\to T(S)$.
For example, $\ast'\to\ast\sqcup\ast'$ induces the map (from pairs of points to points of $T$) that yields the second coordinate.
The map $\{\infty\}\to\alpha\N$ turns into maps from convergent sequences in $T$ to the underlying set simply returning the limit point of the convergent sequence.
(The points of the sequence themselves can, of course, also be extracted via $\{n\}\to\alpha\N$ and the image of $T(\alpha\N)\to\prod_{n\in\N}T(\{n\})$ is the set of convergent sequences.)

Clearly, not every assignment corresponds to a reasonable sort of \enquote{space} being modelled;
in fact, there are some simple observations on spaces of continuous functions that
yield natural conditions we need to impose on $T$.
(There is a zeroth condition.
Namely that $T$ be a functor, so that $\id_{S}\colon S\to S$ induces the identity on $T(S)$ and similarly without compatibility with composition, $T$ would not make sense as a space.)
The other two properties are explained in more detail in the later sections, but can be described as follows.
\begin{enumerate}[(i)]
    \item The functor $T$ should map finite coproduct to finite products (as $C(S\sqcup S', T)=C(S,T)\times C(S',T)$).
    \item The functor $T$ should behave well under factor maps (for any $S\twoheadrightarrow S'$ one can describe the maps in $C(S',T)$ as those maps in $C(S,T)$ which are constant on fibres).
\end{enumerate}
These two requirements (together with some size constraint that can be understood as \enquote{every topological space should have bounded complexity};
e.g., metrisable spaces have complexity $\N$)
yield the final definition of a condensed set.

Hidden in the above approach, we see beautiful and pure incarnations of the following classical principles from functional analysis.
\begin{enumerate}[(i)]
    \item Defining things in Bourbaki-style \enquote{bottom-up} rather than \enquote{top-down} makes them more elegant.
    \item The more disconnected a space is, the better.
    In particular, disconnected compact spaces are preferable to non-compact connected spaces.
The real numbers $\R$ with the Euclidean topology cannot be used to build $\R_d$ with the discrete topology,
but the Stone-\v{C}ech compactification $\beta \R_d$ of $\R_d$ clearly admits $\R$ as a natural factor.
In fact, $\beta\N$ is a space of the nicest possible form,
and whenever something non-compact occurs, Stone-\v{C}ech compactify can be helpful, see, e.g., \cite{hindman2011algebra,schaefer1974Banachlattices, bihlmaier2023compact}.
The best topological spaces are the extremally disconnected compact Hausdorff spaces
(or Stone-\v{C}ech compactifications of discrete spaces)
and everything should be determined by taking factors of these.
    Using weak*-topologies rather than norm-topologies also fits this philosophy.
  \item We should always be able to define topologies via filters or nets:
        namely, by specifying which filters/nets are convergent (and it should not matter which version we use).
    \item The general idea of \enquote{testing spaces} manifests in different ways, which all seem to be combined here.
    \begin{itemize}
        \item Distributions in the sense of having a test class of good objects and defining all objects by their relation to the test objects.
      \item Koopmanism, or better Co-Koopmanism (which is just Opmanism),
            which is the statement that spaces of functions $C(K,\R)$ and pullback operators $T_\phi$ are better than the underlying topological spaces $X$ and the continuous map $\phi$,
        as they have the algebraic properties of the target but encode the topological information of the underlying space.
          See, e.g., \cite{Eisner2016, engel2000one} for a convincing incorporation of this principle.
        \item Unlike classical Koopmanism, these function spaces can be used as well to encode topological information of the target for more complicated targets than $\R$.
            This idea manifests in relative Koopmanism, relative functional analysis, Kaplansky-Banach-spaces or conditional set theory,
            see, e.g., \cite{N.Edeko2024, HenrikKreidler2021, Edeko2023, Jamneshan2013, Jamneshan2020, jamneshan2022uncountable, Jamneshan_2018, carl2019transfer, Hermle2022Halmos}.
    \item Another related idea is the usefulness of the natural embedding of a space into its bidual.
    \end{itemize}
\end{enumerate}

Equipped with these guiding ideas,
we now embark on the journey of developing this theory.

\section{Profinite sets}

The main building blocks of condensed sets are the profinite sets.
In a particular sense one could see the category of condensed sets as a natural enlargement of the category of profinite sets.
In this chapter we study those building blocks, i.e. the category of profinite sets,
in greater detail.

\subsection{Stone duality}\label{subsec:stone_duality}

First we define what exactly we mean by profinite set and establish a connection to well-known topological spaces, namely the Stone spaces.
Then we make some remark about continuous mappings between those.

\begin{definition}[Profin]\label{def:profin}\chaptwo\uses{def:completion}
  Define $\profin$ to be the free pro-completion of the finite sets $\fin$,
  see~\ref{def:completion}.
The objects of $\profin$ are diagrams in $\fin$ indexed by cofiltered small categories,
and morphisms from $D\colon J\to\fin$ to $E\colon I\to \fin$ are given by
\[
    \hom_{\profin}(D,E)=\varprojlim_I\varinjlim_J \hom_\fin(D_j,E_i).
\]
Explicitly this means that morphism $f\colon\varprojlim_{I}D_{i}\varprojlim_{J}E_{j}$
is given by:
for every $i\in I$ an equivalence class $\overline{f_{ij}}$,
where $f_{ij}\colon D_{i}\to E_{j}$.

Equivalently, $\profin$ is the full subcategory of $[\fin^\op,\Set]$ consisting of cofiltered limits of representable objects.
We also write
\[
    "\varprojlim_{j\in J}" D_j
\]
for \(D\) and similarly for \(E\).
\end{definition}

The given description of maps between profinite sets might appear essentially meaningless.
Luckily,
once we have Stone duality,
we can simply think of morphisms between profinite sets as continuous maps.

\begin{definition}[Extr]\label{def:extr}\uses{def:profin, def:special_objects}
Define the category of extremally disconnected spaces $\extr$ to be the full subcategory of $\profin$ consisting of projective objects.
\end{definition}
\begin{remark}
 Note that, although (under Stone duality) they can be seen as compact Hausdorff spaces,
extremally disconnected spaces do not necessarily form compact objects in the sense of~\ref{def:special_objects}.
\end{remark}

\begin{question}
Having seen the description of $\prof$ as a cocompletion,
one wonders if there is a similar description of $\extr$
which directly yields this category from $\fin$.
\end{question}

\begin{theorem}[Stone duality]\label{thm:stone}\uses{def:profin, lem:completions_via_Yoneda}
The following categories are equivalent.
	\begin{enumerate}[(a)]
		\item $\profin$,
		\item $\stone$, consisting of Stone spaces (totally disconnected compact Hausdorff spaces) together with continuous maps,
		\item $\bool^\op$, the dual of the category consisting of Boolean algebras endowed with Boolean algebra homomorphisms.
	\end{enumerate}
\end{theorem}

\begin{proof}
First, we build three functors which induce the desired equivalences.
\begin{center}\begin{tikzcd}
	& \stone \\
	\profin && {\bool^\op}
	\arrow["G", from=1-2, to=2-3]
	\arrow["F", from=2-1, to=1-2]
	\arrow["H", from=2-3, to=2-1]
\end{tikzcd}\end{center}
\begin{enumerate}
\item[F:] Map any formal cofiltered limit $"\varprojlim_i"S_i$ to the cofiltered limit $\varprojlim S_i$ in $\set$, and endow it with the initial topology of all coordinate projections $\varprojlim S_i\to S_i$.
Further, for any
\[
    f\colon "\varprojlim" S_i\to "\varprojlim"T_j,
\]
i.e.
\[
    f=(\bar f_{ij})_j\in \varprojlim_J\varinjlim_I \hom(S_i,T_j),
\]
        we define $F(f)$ as the map $f=(f_j)_j$ induced by the projections $\varprojlim S_i\to S_i\overset{f_{ij}}\to T_j$.

\item[G:] Map every Stone space $X$ to the Boolean algebra of clopen sets,
\[
    G(X)=\mathrm{clopen}(X)=\hom_{\stone}(X,\mathbb{F}_2),
\]
and every morphism $f\colon X\to Y$ to the induced pullback $f^*\colon G(X)\to G(Y)$ (note that we consider $\bool^\op$, hence the direction of the pullback remains the same.)

\item[H:] For a Boolean algebra $A$ consider the index set $J=\{A_j\le A, \,|A|<\infty\}$ consisting of finite subalgebras, endowed with the cofiltration $A_i\to A_j$ for $A_i\subseteq A_j$. Now define $H(A)$ as $"\varprojlim_j" \hom(A_j,\mathbb{F}_2)$. For each opposite Boolean homomorphism $f\colon A\to B$ consider the map
\[
    \varprojlim_j \hom(A_j,\mathbb{F}_2)\to \varprojlim_k \hom(B_k,\mathbb{F}_2)
\]
via the tuple $(h_k)$ where $h_k=\bar f^*$ for $f^*\colon \hom(f(B_k),\mathbb{F}_2)\to \hom(B_k, \mathbb{F}_2)$.
\end{enumerate}

Next, we prove
\begin{enumerate}
\item $HGF\simeq 1_{\profin}$
\item $FHG\simeq 1_\stone$
\item $GFH\simeq 1_{\bool^\op}$
\end{enumerate}
The natural isomorphy of these functors can be checked quite elementary, as we will do in the following. Together these three identities imply the equivalence of all three categories.

Before embarking on these three tasks, we show that the functor \(F\) is faithful.
For this,
consider the map
\[
    \varinjlim_i C(S_i,T) \to C(\varprojlim_i S_i,T),\quad [f_i]\mapsto f_i\circ\pi_i
\]
where \(S_i\) and \(T\) are finite sets, \(I\) is cofiltered, and \(\pi_i\colon\varprojlim_i S_i\to S_i\) denotes the coordinate projection
and we interpret \(\varprojlim S_i\) as a topological space.
It is obvious that this map is well-defined because of the equivalence relation on \(\varinjlim_i C(S_i,T)\).
Moreover this map is injective.

Denote by \(\tau_k^i\) the transition map \(S_k\to S_i\) if it exists and then we abbreviate \(k\geq i\).
We have \(f_k = f_i\tau_k^i\) for \(k\geq i\).
Now for \([f_i]\neq [g_j]\) we have that \(f_k\neq g_k\) for all \(k\geq i,j\).
Let \(Q_k = \{s\in S_k : f_k(s)\neq g_k(s)\}\) for any \(k\) such that \(f_k\in [f_i]\) and \(g_k\in [g_j]\) and otherwise \(Q_k = S_k\).
For \(k\geq i\) we obtain that \(\tau_k^i(Q_k)\subset Q_i\).
Define
\[
    R_k^i = \{ x\in\prod_{i\in I} S_i : \tau_k^i(x_k) = x_i,\; x_k\in Q_k\}.
\]
Then these sets are closed subsets of the compact Hausdorff space \(\prod_{i\in I} S_i\).
For finitely many \(k_1,\ldots,k_n\) and \(i_1,\ldots,i_n\) (such that \(k_l\geq i_k\))
we have that
\[
    \bigcap_{l=1}^n R_{k_l}^{i_l} \neq \emptyset
\]
because we can choose \(k_0\geq k_1,\ldots,k_n\) (since \(I\) is cofiltered)
and \(x_{k_0}\in Q_{k_0}\subset S_{k_0}\) (because \(f_{k_0}\neq g_{k_0}\) if they are defined).
Then the tupel \(x \in \prod_{i\in I} S_i\) with \(x_i = \tau_{k_0}^{i}(x_{k_0})\)
if \(i\) is \(k_0,k_1,\ldots,k_n,i_1,\ldots,i_n\) and \(x_i\in S_i\) arbitrary otherwise
is an element of the intersection of the \(R_{k_l}^{i_l}\).
Thus there exists
\[
    y\in\bigcap_{k\geq i} R_k^i.
\]
This element is obviously an element of \(\varprojlim_{i\in I} S_i\) and it fulfills that \(f_i(y_i) \neq g_i(y_i)\) whenever this is defined.
Therefore \(f_i\pi_i\neq g_i\pi_i\) what was to be shown.

Now, we are ready.
\begin{enumerate}
\item We build two natural morphisms $\eta \colon 1_{\profin}\to HGF$, and $\gamma\colon HGF\to 1_{\profin}$, which are inverse to each other.
Let $S="\varprojlim_I" S_i$.
The coordinate morphism
\[
    \eta^S \colon "\varprojlim_i" S_i\to "\varprojlim_j" \hom(A_j,\mathbb{F}_2)
\]
with $A_j\le \mathrm{clopen}(\varprojlim_i S_i)$ being finite subalgebras is given as follows.
For each $j$ we need a coordinate $\eta^S_j\in \varinjlim_i \hom(S_i, \hom(A_j,\mathbb{F}_2))$.
To this end, we look for an appropriate $S_i$.
Each element $a\in A_j$ is a clopen set of $\varprojlim S_i$.
Hence it is the union of all cylinders $C\subseteq a$.
By compactness this union is finite, i.e.
\[
    a=\bigcup_{\ell,\,\mathrm{fin}} C_\ell
\]
for
\[
    C_\ell=\bigcap_{h\,\mathrm{fin}}\pi_{i_{h}}^{-1}(C^h_\ell)
\]
with $C^h_\ell\subseteq S_{i_h}$.
As the index set is cofiltered, there exists an $S_\ell$ such that $S_\ell$ maps to all the $S_{i_h}$:
\begin{center}\begin{tikzcd}
	& {S_\ell} \\
	{S_{i_1}} & {S_{i_2}} & {S_{i_h}}
	\arrow["{\tau_{\ell}^{i_1}}"', from=1-2, to=2-1]
	\arrow[from=1-2, to=2-2]
	\arrow["{\tau_{\ell}^{i_h}}", from=1-2, to=2-3]
\end{tikzcd}\end{center}
Now by
\begin{center}\begin{tikzcd}
	{\varprojlim S_i} && {S_\ell} \\
	& {S_{i_1}} & {S_{i_2}} & {S_{i_h}}
	\arrow["{\pi_\ell}", from=1-1, to=1-3]
	\arrow["{\pi_{i_1}}"', from=1-1, to=2-2]
	\arrow["{\tau_{\ell}^{i_1}}"', from=1-3, to=2-2]
	\arrow[from=1-3, to=2-3]
	\arrow["{\tau_{\ell}^{i_h}}", from=1-3, to=2-4]
\end{tikzcd}\end{center}
we have
\[
    C_\ell
    =\bigcap \pi_{i_h}^{-1}(C^h_\ell)
    =\bigcap \pi_\ell^{-1}({\tau_{\ell}^{i_h}}^{-1}(C^h_\ell))=\pi_\ell^{-1} (B_\ell)
\]
for the appropriate
\[
    B_\ell=\bigcap {\tau_{\ell}^{i_h}}^{-1}(C^h_\ell)\subseteq S_\ell.
\]
Hence every cylinder is of the Form $C=\pi_\ell^{-1}(B_\ell)$.
Analogously, there exists $i$ such that for all $\ell$ in the decomposition
$a=\bigcup_{\ell,\,\mathrm{fin}} C_\ell$ we factor over $S_i$,
\begin{center}\begin{tikzcd}
	{\varprojlim S_i} && {S_i} \\
	& {S_{\ell}}
	\arrow["{\pi_i}", from=1-1, to=1-3]
	\arrow["{\pi_{\ell}}"', from=1-1, to=2-2]
	\arrow["{\tau_{i}^{\ell}}"', from=1-3, to=2-2]
\end{tikzcd}.\end{center}
Hence
\[
    a=\bigcup C_\ell= \pi_i^{-1}(\bigcup{(\tau_i^{\ell})}^{-1}(B_\ell))=\pi_i^{-1}(B_a)
\]
for appropriate $B_a\subseteq S_i$. Cofiltering over the finitely many $a\in A_j$ we can assume $S_i$ to be such that all $a\in A_j$ are of the form $\pi_i^{-1}(B_a)$ for some $B_a\subseteq S_i$ for the same $i$.

Let \(f_i\colon S_i\to\mathbb{F}_2\) be the constant function with value \(1\)
and \(g_i\colon S_i\to\mathbb{F}_2\) be the indicator function of the image of \(\pi_i\).
Hence we have \(h_i\pi_i = k_i\pi_i\) and by injectivity of the map
\[
    \varinjlim C(S_i,\mathbb{F}_2)\to C(\varprojlim S_i,\mathbb{F}_2)
\]
we have that \(f_i\) and \(g_i\) are equivalent.
Therefore there is \(\tau_k^i\colon S_k\to S_i\) such that \(f_i\tau_k^i = g_i\tau_k^i = 1\).
This implies that \(\tau_k^i(S_k)\subset\im(\pi_i)\).

For this $S_k$ we define $\eta^S_{kj}\colon S_k\to \hom(A_j,\mathbb{F}_2)$ via
\[
    s\mapsto \left(a\mapsto\begin{cases}1&\tau_k^i(s)\in B_a\\ 0&\tau_k^i(s)\not\in B_a\end{cases}\right).
\]
For any \(s\), this defines a Boolean algebra homomorphism.
First remark that
\[
    B^c_a \cap B_{a^c} \subset S_i\setminus\im(\pi_i)
\]
and
\[
    \pi_i^{-1}(B_{a\wedge b}) = \pi_i^{-1}(B_a\cap B_b).
\]
Now because \(\tau_k^i(s)\in\im(\pi_i)\) it follows that
\(\tau_k^i(s)\in B_{a^c}\) iff \(\tau_k^i(s)\notin B_a\)
and \(\tau_k^i(s)\in B_{a\wedge b}\) iff \(\tau_k^i(s)\in B_a\cap B_b\).
In particular \(\eta_{kj}^S(s)\) is a Boolean algebra homomorphism.

In the following we will omit the \(\tau_k^i\)
and assume that the \(B_a\) already fulfill the relations
\(B_{a\wedge b} = B_a\cap B_b\) and \(B_a^c = B_{a^c}\).

Now let $\eta^S_{j}=\overline{\eta_{ij}^S}$ with the usual equivalence relation on a cofiltered index category.
To see that $\eta^S=(\eta^S_j)_j$ induces a morphism as follows
\begin{center}\begin{tikzcd}
	{"\varprojlim\limits_i" S_i} && {"\varprojlim\limits_j" \hom(A_j,\mathbb{F}_2)} \\
	{S_i} && {\hom(A_j, \mathbb{F}_2)}
	\arrow["\eta", dashed, from=1-1, to=1-3]
	\arrow[from=1-1, to=2-1]
	\arrow["{\eta_j^S}"{description}, from=1-1, to=2-3]
	\arrow[from=1-3, to=2-3]
	\arrow["{\eta^S_{ij}}", from=2-1, to=2-3]
\end{tikzcd},\end{center}
we need commutativity in $j$.
For this consider a morphism $A_j\to A_{j'}$ in \(\bool^\op\), i.e. $A_j\supseteq A_{j'}$.
Then \(\tau_j^{j'}\colon\hom(A_j,\mathbb{F}_2)\to\hom(A_{j'},\mathbb{F}_2)\) is given by restriction.
We have to show that
\[
    \eta_j^S = \tau_j^{j'} \eta_{j'}^S.
\]
Now choose representatives \(S_i\) and \(S_{i'}\) for the maps \(\eta_j^S\) and \(\eta_{j'}^S\).
Refine the corresponding $S_i$ and $S_{i'}$ by some $S_{i''}$
and the \(S_{i''}\) with \(S_k\) such that \(\tau_k^{i''}(S_k)\subset\im(\pi_{i''})\).
In a diagram it looks like this:
\begin{center}\begin{tikzcd}
	{S_k} &&& {"\varprojlim" S_i} && {"\varprojlim" \hom(A_j,\mathbb{F}_2)} \\
	\\
	{S_{i''}} &&& {S_i} && {\hom(A_j, \mathbb{F}_2)} \\
	&& {S_{i'}} && {\hom(A_{j'}, \mathbb{F}_2)}
	\arrow["{\tau_k^{i''}}"', from=1-1, to=3-1]
	\arrow["{\pi_k}"', from=1-4, to=1-1]
	\arrow["\eta", dashed, from=1-4, to=1-6]
	\arrow["{{{\pi_{i''}}}}"', from=1-4, to=3-1]
	\arrow["{{{\pi_i}}}", from=1-4, to=3-4]
	\arrow["{{{\eta_j^S}}}"{description}, from=1-4, to=3-6]
	\arrow["{{{\pi_{i'}}}}"', from=1-4, to=4-3]
	\arrow["{{\eta^S_{j'}}}"{description}, from=1-4, to=4-5]
	\arrow["{{{\pi_j}}}", from=1-6, to=3-6]
	\arrow["{{{\pi_{j'}}}}"'{pos=0.3}, shift left, from=1-6, to=4-5]
	\arrow["{{\tau_{i''}^{i}}}", from=3-1, to=3-4]
	\arrow["{{{\tau_{i''}^{i'}}}}"', from=3-1, to=4-3]
	\arrow["{{{\eta^S_{ij}}}}", from=3-4, to=3-6]
	\arrow["{{{\tau_j^{j'}}}}", from=3-6, to=4-5]
	\arrow["{{{\eta_{i'j'}^S}}}", from=4-3, to=4-5]
\end{tikzcd}\end{center}
So now it suffices to show
\[
    \tau_j^{j'}\eta_{ij}^S\tau_{i''}^{i}=\eta^S_{i'j'}\tau_{i''}^{i'}.
\]
Thus let \(s\in S_k\) and \(a\in A_{j'}\).
Then there exist \(B_a^j\subset S_i\) and \(B_a^{j'}\subset S_{i'}\) such that
\[
     \pi_{i''}^{-1}((\tau_{i''}^{i'})^{-1}(B^{j'}_a))
    =\pi_{i'}^{-1}(B^{j'}_a)
    =a
    =\pi_i^{-1}(B^j_a)
    =\pi_{i''}^{-1}((\tau_{i''}^{i})^{-1}(B^{j}_a))
\]
since \(A_{j'}\subset A_j\) and \(a\) is a clopen subset of \(\varprojlim S_i\).
Hence
\[
        \eta^S_{i'j'}\tau_{k}^{i'}(s)(a)
    =   \eta^S_{i'j'}\tau_{i''}^{i'}\tau_k^{i''}(s)(a)
    =   \begin{cases}
            1,\quad \tau_{i''}^{i'}\tau_k^{i''}(s)\in B_a^{j'} \\
            0,\quad \tau_{i''}^{i'}\tau_k^{i''}(s)\notin B_a^{j'}
    \end{cases}
    =   \begin{cases}
            1,\quad \tau_k^{i''}(s)\in (\tau_{i''}^{i'})^{-1}(B_a^{j'}) \\
            0,\quad \tau_k^{i''}(s)\notin (\tau_{i''}^{i'})^{-1}(B_a^{j'})
    \end{cases}
\]
and similarly
\[
        \tau_j^{j'}\eta_{ij}^S\tau_{i''}^{i}(s)(a)
    =   \begin{cases}
            1,\quad \tau_k^{i''}(s)\in (\tau_{i''}^{i})^{-1}(B_a^{j}) \\
            0,\quad \tau_k^{i''}(s)\notin (\tau_{i''}^{i})^{-1}(B_a^{j})
    \end{cases}
\]
Now from \(\tau_k^{i''}(s)\in\im(\pi_{i''})\) it follows that
\(\tau_k^{i''}(s)\in (\tau_{i''}^{i})^{-1}(B_a^j)\)
iff
\(\tau_k^{i''}(s)\in (\tau_{i''}^{i'})^{-1}(B_a^{j'})\).
and the equality  of the two morphism follows.
Hence the morphism $\eta^S$ is well-defined.

Next we show, that it is natural. Therefore consider a morphism $f=(\overline{f_{ij}})_j=f\in \varprojlim_j\varinjlim_i\hom(S_i,T_j)\colon "\varprojlim"S_i\to "\varprojlim" T_j$. We need commutativity of
\begin{center}\begin{tikzcd}
	{S_{i_1}} && {S_i} \\
	& {S_k} &&& {\hom (A_k,\mathbb{F}_2)} \\
	{S_{i_0}} && {"\varprojlim " S_i} & {"\varprojlim " \hom(A_k,\mathbb{F}_2)} \\
	&& {"\varprojlim " T_j} & {"\varprojlim " \hom(B_\ell,\mathbb{F}_2)} \\
	& {T_j} &&& {\hom(B_\ell,\mathbb{F}_2)}
	\arrow[from=1-1, to=1-3]
	\arrow[from=1-1, to=3-1]
	\arrow["{{\eta^S_{ik}}}"', from=1-3, to=2-5]
	\arrow[from=2-2, to=1-1]
	\arrow["{{\tilde f_{k\ell}}}"', from=2-5, to=5-5]
	\arrow["{{f_{i_0j}}}"', from=3-1, to=5-2]
	\arrow[from=3-3, to=1-3]
	\arrow[from=3-3, to=2-2]
	\arrow[from=3-3, to=3-1]
	\arrow["{{\eta^S}}", from=3-3, to=3-4]
	\arrow["f", from=3-3, to=4-3]
	\arrow[from=3-4, to=2-5]
	\arrow["{{\tilde f}}", from=3-4, to=4-4]
	\arrow["{{\eta^T}}", from=4-3, to=4-4]
	\arrow[from=4-3, to=5-2]
	\arrow[from=4-4, to=5-5]
	\arrow["{{\eta^T_{j\ell}}}"', from=5-2, to=5-5]
\end{tikzcd}.\end{center}
where \(\tau_k^{i_1}(S_k)\subset\im(\pi_{i_1})\).
The commutativity can be checked elementwise in \(\fin\) and uses similar arguments as before.
Note that \(\tilde{f}\) is induced by pullback and taking the inverse image of a clopen set.

Hence in total $(\eta^S)$ is a natural transformation.

It remains to show invertibility of $\eta^S$. For this, define $\gamma^S$ via
\begin{center}\begin{tikzcd}
	{"\varprojlim"\hom(A_j,\mathbb{F}_2)} && {"\varprojlim" S_i} \\
	{\hom(A_j,\mathbb{F}_2)} && {S_i}
	\arrow["\gamma", from=1-1, to=1-3]
	\arrow["{\pi_j}", from=1-1, to=2-1]
	\arrow["{\pi_i}", from=1-3, to=2-3]
	\arrow["{\eps\mapsto \inf\eps^{-1}(1)}", from=2-1, to=2-3]
\end{tikzcd}\end{center}
where $A_j$ is chosen to be the Boolean sub-algebra generated by $\{\pi_i^{-1}(s):s\in S_i\}$,
and $\inf\eps^{-1}(1)$ is that unique $t\in S_i$ such that $\eps(\pi_i^{-1}(t))=1$.
This is well-defined, as for any $S_i\to S_k$ we can cofilter in $j$ via
\begin{center}\begin{tikzcd}
	{S_i} && {S_k} \\
	{\hom(A_j, \mathbb{F}_2)} && {\hom(A_k,\mathbb{F}_2)} \\
	& {\hom(A_{j'},\mathbb{F}_2)} \\
	& {"\varprojlim" \hom(A_j,\mathbb{F}_2)}
	\arrow["f", from=1-1, to=1-3]
	\arrow[from=2-1, to=1-1]
	\arrow[from=2-3, to=1-3]
	\arrow[from=3-2, to=2-1]
	\arrow[from=3-2, to=2-3]
	\arrow[from=4-2, to=2-1]
	\arrow[from=4-2, to=2-3]
	\arrow[from=4-2, to=3-2]
\end{tikzcd}.\end{center}
We want to show that $f(s_i)=s_k$, where $s_i$ and $s_k$ are the induced elements of some object \(\eps\) in $\hom (A_{j'},\mathbb{F}_2)$.
This means that \(\eps(\pi_i^{-1}(s_i)) = \eps(\pi_k^{-1}(s_k)) = 1\).
Therefore \(\eps(\pi_i^{-1}(s_i)\cap \pi_k^{-1}(s_k)) = 1\).
But then \(\pi_i^{-1}(s_i)\cap\pi_k^{-1}(s_k)\subset\pi_k^{-1}(f(s_i))\) implies
that \(\eps(\pi_k^{-1}(f(s_i))) = 1\).

It remains to show $\eta^S\gamma^S=1$ and $\gamma^S\eta^S=1$. For the latter, consider
\begin{center}\begin{tikzcd}
	{"\varprojlim" S_i} & {"\varprojlim" \hom(A_j,\mathbb{F}_2)} & {"\varprojlim" S_i} \\
	{S_i} & {\hom(A_j,\mathbb{F}_2)} & {S_i} \\
	s & {(a\mapsto1_{B_a}(s))} & {\inf(\eta_{ij}(s)^{-1}(1))}
	\arrow["\eta", from=1-1, to=1-2]
	\arrow[from=1-1, to=2-1]
	\arrow["\gamma", from=1-2, to=1-3]
	\arrow[from=1-2, to=2-2]
	\arrow[from=1-3, to=2-3]
	\arrow["{\eta_{ij}}"', from=2-1, to=2-2]
	\arrow["{\gamma_{ji}}"', from=2-2, to=2-3]
	\arrow[maps to, from=3-1, to=3-2]
	\arrow[maps to, from=3-2, to=3-3]
\end{tikzcd}\end{center}
which is the identity because \(\eta_{ij}(s)(\pi_i^{-1}(s)) = 1\).
For the first consider
\begin{center}\begin{tikzcd}
	{"\varprojlim"\hom(A_j,\mathbb{F}_2)} & {"\varprojlim" S_i} & {"\varprojlim"\hom(A_j,\mathbb{F}_2)} \\
	{\hom(A_k,\mathbb{F}_2)} & {S_i} & {\hom(A_j,\mathbb{F}_2)}
	\arrow["\gamma", from=1-1, to=1-2]
	\arrow[from=1-1, to=2-1]
	\arrow["\eta", from=1-2, to=1-3]
	\arrow[from=1-2, to=2-2]
	\arrow[from=1-3, to=2-3]
	\arrow["{\gamma_{ki}}"', from=2-1, to=2-2]
	\arrow["{\eta_{ij}}"', from=2-2, to=2-3]
\end{tikzcd}\end{center}
By construction we have \(A_j\subset A_k\).
Furthermore
\[
    \eta_{ij}(\inf(\eps^{-1}(1)))(a)
    =\eta_{ij}(s_\epsilon)(a)
    =\begin{cases}1&s_\eps\in B_a\\0&s_\eps\in B_{a^c}\end{cases}
\]
as well as $a=\pi_i^{-1}(B_a)$ and $\pi_i^{-1}(s_\eps)=\inf \eps^{-1}(1)\in \eps^{-1}(1)$,
hence $s_\eps\in B_a$ is equivalent to $1=\eps(a)$, thus the function agrees with $\eps(a)$.
This means it factors as the identity
\begin{center}\begin{tikzcd}
	{\varprojlim\hom(A_j,\mathbb{F}_2)} & {\varprojlim S_i} & {\varprojlim\hom(A_j,\mathbb{F}_2)} \\
	{\hom(A_k,\mathbb{F}_2)} & {S_i} & {\hom(A_j,\mathbb{F}_2)} \\
	{\hom(A_j,\mathbb{F}_2)}
	\arrow["\gamma", from=1-1, to=1-2]
	\arrow[from=1-1, to=2-1]
	\arrow["\eta", from=1-2, to=1-3]
	\arrow[from=1-2, to=2-2]
	\arrow[from=1-3, to=2-3]
	\arrow["{\gamma_{ki}}"', from=2-1, to=2-2]
	\arrow["{\eta_{ij}}"', from=2-2, to=2-3]
	\arrow[from=2-1, to=3-1]
	\arrow["1"', from=3-1, to=2-3]
\end{tikzcd}\end{center}
and hence the composition induces the identity morphism. This makes us very happy.

\item Put $\phi^X\colon X\to \varprojlim \hom(A_i,\mathbb{F}_2)$,
where $A_i\le \mathrm{clopen}(X)$ finite via
\[
    x\mapsto \left(a\mapsto\begin{cases}1&x\in a\\ 0&x\not\in  a\end{cases}\right).
\]
This is continuous by definition of the topology, compatible with maps $X\to Y$ and hence defines a natural transformation.
Injectivity of this morphism follows by Hausdorff-property and 0-dimensionality of $X$,
and surjectivity follows, as for any compatible tuple $(\eps_i)$ the preimages $\eps_i^{-1}(1)$ are clopen in $X$ and fulfill the finite intersection property, hence admitting $x\in \bigcap_i \eps_i^{-1}(1)$, being mapped to $(\eps_i)$. Thus $\phi$ is a natural isomorphism.

\item Consider a Boolean algebra $A$. Build a Boolean homomorphism
\[
    \Xi\colon A\to \mathrm{clopen}(\varprojlim_j\hom(A_j,\mathbb{F}_2))
\]
via mapping any $a\in A$ to $\pi_{j}^{-1}(a)$, where the $j$-th coordinate is precisely the one corresponding to the finite sub-algebra generated by $a$. It can be checked that this is a natural isomorphism.
\end{enumerate}
\end{proof}

\begin{example}
  Classically,
  we know some Stone spaces and we want to explain how to explicitly describe them as cofiltered limits of finite (discrete) spaces.
  If we can describe a Stone space $X$ as a cofiltered limit $\varprojlim S_{i}$ of finite discrete spaces in $\Top$
  then we automatically know that $X=\varprojlim S_{i}$ in $\stone$
  as $\stone$ is just a full subcategory of $\Top$,
  see~\ref{lem:lim_in_full_sub}.
  \begin{enumerate}
    \item Consider the one point compactification $\alpha\N=\N\cup\{\infty\}$ of the natural numbers.
          This forms a totally disconnected compact Hausdorff space corresponding to the opposite Boolean algebra of the Boolean algebra of finite or cofinite sets.
          (This description of the clopens is clear from the definition of the Alexandrov compactification.)
          We want to show that in $\Top$,
          $\alpha\N$ is the limit
          \[\alpha\N=\varprojlim_n \{0,1,\dots, n, \infty\}\]
          with inclusions as the only transition maps (yielding a cofiltered diagram).

          To see this, consider the cone
\begin{center}\begin{tikzcd}
	{\{1,\dots, n,\infty\}} & {\{1,\dots,m,\infty\}} \\
	& \alpha\N
	\arrow["{n\ge m}", from=1-1, to=1-2]
	\arrow["{\phi_n}", from=2-2, to=1-1]
	\arrow["{\phi_m}"', from=2-2, to=1-2]
\end{tikzcd}\end{center}
where $\phi_n(a)=\begin{cases}a,& a\le n,\\\infty,&\text{else.}\end{cases}$
          For any $X\in\Top$ with
\begin{center}\begin{tikzcd}
	{\{1,\dots, n,\infty\}} & {\{1,\dots,m,\infty\}} \\
	X & \alpha\N
	\arrow["{n\ge m}", from=1-1, to=1-2]
	\arrow["{f_n}", from=2-1, to=1-1]
	\arrow["{f_m}"{pos=0.2}, from=2-1, to=1-2]
	\arrow["\psi"', from=2-1, to=2-2]
	\arrow["{\phi_n}"'{pos=0.2}, from=2-2, to=1-1]
	\arrow["{\phi_m}"', from=2-2, to=1-2]
\end{tikzcd}\end{center}
clearly $\psi(x)\coloneqq \lim f_n(x)$ (topological limit) does the job, is unique and is continuous as for any end piece $[n,\infty]$, the preimage is given by $f_{n-1}^{-1}(\infty)$.
\item The Cantor set
\[C=\{0,1\}^\N =\varprojlim_{n}\{0,1\}^n,\]
 and more generally for infinite $I$ the cube $\{0,1\}^I$ is profinite. This representation is due to the general decomposition of limits into a finite and a cofiltered limit, or the fact that the system of these finite cubes are final in the clopens.
\item Consider $\beta\N$. Note that for any Boolean algebra $B$ the homomorphisms $\hom(B,\mathbb{F}_2)$ are nothing but ultrafilters in $B$. Now by indexing with all finite subalgebras $A_n$ of $2^\N=\mathrm{clopen}(\beta\N)$, we have that the limit system corresponding to $\beta\N$ is given by
\[\varprojlim_{n} \mathcal{U}(A_n).\]
By noting that ultrafilter properties can be glued together from all finite subalgebras to the whole algebra after passing to the concrete limit in $\set$, as
\[\varprojlim \hom(A_n,\mathbb{F}_2)=\hom(\varinjlim A_n ,\mathbb{F}_2)\]
and as in $\bool$ we have $\varinjlim A_n=A$.
          We can further identify this set with all ultrafilters on $\N$, where the topology from $\beta\N$ is the initial topology of all finite partitions of $\N$.
          More generally this shows that any Stone space $X$ can be identified with the set of ultrafilters of clopen sets of $X$, and the topology is the initial topology of the finite quotients (fixed filters).
          The size of $\beta\N$ is clearly $2^{2^{\aleph_{0}}}$.
\end{enumerate}
\end{example}

\begin{lemma}[Profin complete]\label{lem:prof_complete}\uses{def:profin, lem:pro-completion_is_complete, lem:set_bicomplete}
The category $\profin$ is complete and admits finite coproducts.
\end{lemma}

\begin{proof}
This follows directly from the dual statement of proposition 6.1.18 in \cite{kashiwara2005categories} since \(\fin\) admits finite limits.
The second statement follows from corollary 6.1.17 ibidem.
Alternatively, see lemma~\ref{lem:pro-completion_is_complete}.
\end{proof}

\begin{remark}
We could also argue using Stone duality and compute the limit in \(\CHaus\) via Tychonov's theorem.

Note that, however, we cannot use this argumentation for arbitrary coproducts,
as the Stone-\v{C}ech-functor $\beta$ does not preserve 0-dimensional Hausdorff spaces (see \cite[p.~87]{walker2012stone}).
\end{remark}

Next we show, that profinite sets always admit a slightly more convenient description.

\begin{lemma}[Transition maps epimorphisms]\label{lem:transitionepi}\uses{def:profin, thm:stone, lem:limits_along_final_subdiagrams}
One can always assume all the projections as well as all the transition maps in the diagram
in the profinite set to be surjective;
i.e., every profinite set is isomorphic to a profinite set $\varprojlim S_i$ with all the coordinate projections and transition maps being surjective.
Furthermore one can always restrict the diagram of a profinite set $\varprojlim S_i$ along a cofinal subset.

\end{lemma}

\begin{proof}
  Because Stone duality maps finite sets (in $\prof$) to finite Boolean algebras
  and is an equivalence of categories,
  we simply have to note that for $A\in\bool$,
  \[
    A=\varinjlim_{A_{i}\subseteq A\textrm{ fin.}}A_{i}
  \]
  which is clearly a filtered colimit of finite Boolean algebras along monic transition maps with monic \enquote{embeddings} into the colimit.

  Alternatively, note that restricting all the $S_i$ to the images $\pi_i(S)$ does not change the inverse limit in $\CHaus$,
  since the subsets of the product correspond and both carry the same relative topology.
  Then,
  one can replace $S_{i}$ by the intersection of the images of all arrows in the limit diagram pointing to $S_{i}$.

  The last claim is just~\ref{lem:limits_along_final_subdiagrams}.
\end{proof}

\begin{remark}\label{rem:transitionepi}\uses{lem:transitionepi}
In the following (especially in the second part of the next lemma) we will always assume all profinite sets to be of such form.
\end{remark}

Lastly, we describe morphisms between profinite sets.
The lemma can be proven using category-theoretic statements such as \enquote{limits of monics are monic}
ut we choose to give a classical prove showcasing how to translate smoothly between formal profinite sets and their topological incarnations.

\begin{lemma}[2.5.9/3.2.14 in \cite{engelking1989general}]\label{lem:cofinal_map_between_prof}\uses{thm:stone, lem:transitionepi, lem:limits_along_final_subdiagrams}
  Let $f=(\overline{f^j_{\phi(j)}})\colon\varprojlim_{I}S_i\to \varprojlim_J T_j$ be a mapping between two profinite sets such that  $\phi\colon J\to I$ is monotone with cofinal image
  (i.e. for all $i\in I$ there exists $j\in J$ and $\phi(j)\to i$).
\begin{enumerate}
\item If all the $f^j_{\phi(j)}\colon S_{\phi(j)} \to T_j$ are injective, so is $f$.
\item If all the $f^j_{\phi(j)}\colon S_{\phi(j)}\to T_j$ are surjective, so is $f$.
\end{enumerate}
\end{lemma}
\begin{proof}
\begin{enumerate}
        By Stone duality,
        we work in $\stone$ (or $\Top$) so that monomorphicity is just injectivity and similarly for epics.
\item Take $x\ne y\in \varprojlim S_i$.
Then there exists $i\in I$ with $x_i\ne y_i$.
Now by cofinality of $\phi$ choose any $j\in J$ such that $\phi(j)\to i$.
This leads to $x_{\phi(j)}\ne y_{\phi(j)}$,
since
\[
    \tau^{\phi(j)}_i (x_{\phi(j)}) = x_i\ne y_i = \tau^{\phi(j)}_{i}(y_{\phi(j)}).
\]
Thus by injectivity of $f^j_{\phi(j)}$ we have
\[
    f(x)_j= f^j_{\phi(j)}(x_{\phi(j)})\ne f^j_{\phi(j)}(y_{\phi(j)})=f(y)_j
\]
and consequently $f(x)\ne f(y)$.
\item Let $t=(t_j)\in \varprojlim_J T_j$.
For any fixed $t_j$ define
\[
    Z^j\coloneqq \pi_{\phi(j)}^{-1}((f^j_{\phi(j)})^{-1}(t_j))
    \subseteq \pi_{\phi(j)}^{-1}(S_{\phi(j)})=S.
\]
Every $Z^j$ is nonempty, as $f^j_{\phi(j)}$ is surjective,
and we assumed all the $\pi_{\phi(j)}$ to be surjective.
We want to find some $x\in \bigcap_j Z^j$.
For this we show finite intersection property of the $Z^j$.

Hence, let $j_1,j_2\in J$. Then there exists $j\in J$ with $j\to j_1$ and $j\to j_2$ such that
\begin{center}\begin{tikzcd}
	&&& {S_{\phi(j_1)}} &&& {T_{j_1}} \\
	{\varprojlim S_i} && {S_{\phi(j)}} &&& {T_j} \\
	&&& {S_{\phi(j_1)}} &&& {T_{j_2}}
	\arrow[from=1-4, to=1-7]
	\arrow[from=2-1, to=1-4]
	\arrow[from=2-1, to=2-3]
	\arrow[from=2-1, to=3-4]
	\arrow[from=2-3, to=2-6]
	\arrow[from=2-6, to=1-7]
	\arrow[from=2-6, to=3-7]
	\arrow[from=3-4, to=3-7]
\end{tikzcd}.\end{center}
But this means (the equivalence relation on morphisms $S_i$ is defined that way) that there exists $i$ such that
\begin{center}\begin{tikzcd}
	&&& {S_{\phi(j_1)}} &&& {T_{j_1}} \\
	{S_i} && {S_{\phi(j)}} &&& {T_j} \\
	&&& {S_{\phi(j_1)}} &&& {T_{j_2}}
	\arrow[from=1-4, to=1-7]
	\arrow[from=2-1, to=1-4]
	\arrow[from=2-1, to=2-3]
	\arrow[from=2-1, to=3-4]
	\arrow[from=2-3, to=2-6]
	\arrow[from=2-6, to=1-7]
	\arrow[from=2-6, to=3-7]
	\arrow[from=3-4, to=3-7]
\end{tikzcd}.\end{center}
But as $\phi $ is cofinal, there exists $k\in J$ with
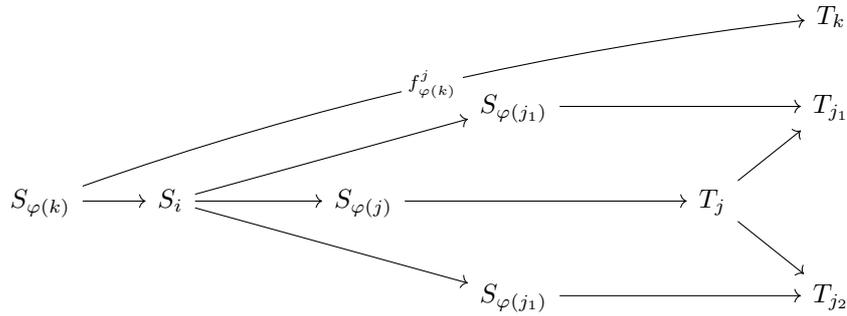
\begin{center}\begin{tikzcd}
	&&&&&&& {T_k} \\
	&&&& {S_{\phi(j_1)}} &&& {T_{j_1}} \\
	{S_{\phi(k)}} & {S_i} && {S_{\phi(j)}} &&& {T_j} \\
	&&&& {S_{\phi(j_1)}} &&& {T_{j_2}}
	\arrow[from=2-5, to=2-8]
	\arrow["{f^j_{\phi(k)}}"{description}, curve={height=-12pt}, from=3-1, to=1-8]
	\arrow[from=3-1, to=3-2]
	\arrow[from=3-2, to=2-5]
	\arrow[from=3-2, to=3-4]
	\arrow[from=3-2, to=4-5]
	\arrow[from=3-4, to=3-7]
	\arrow[from=3-7, to=2-8]
	\arrow[from=3-7, to=4-8]
	\arrow[from=4-5, to=4-8]
\end{tikzcd}\end{center}
and by monotonicity of $\phi$ we deduce
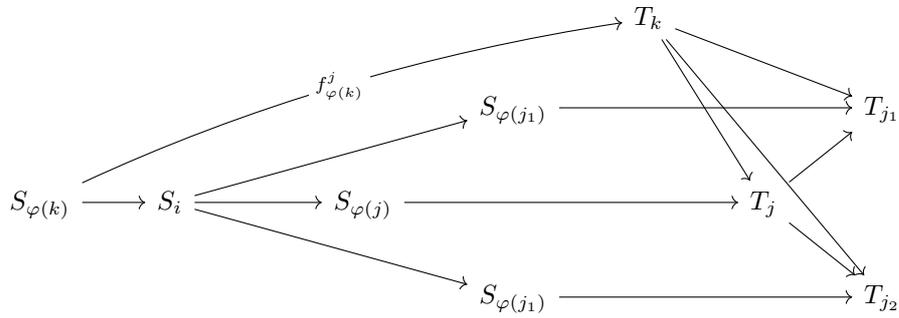
\begin{center}\begin{tikzcd}
	&&&&& {T_k} \\
	&&&& {S_{\phi(j_1)}} &&& {T_{j_1}} \\
	{S_{\phi(k)}} & {S_i} && {S_{\phi(j)}} &&& {T_j} \\
	&&&& {S_{\phi(j_1)}} &&& {T_{j_2}}
	\arrow[from=1-6, to=2-8]
	\arrow[from=1-6, to=3-7]
	\arrow[from=1-6, to=4-8]
	\arrow[from=2-5, to=2-8]
	\arrow["{f^j_{\phi(k)}}"{description}, curve={height=-12pt}, from=3-1, to=1-6]
	\arrow[from=3-1, to=3-2]
	\arrow[from=3-2, to=2-5]
	\arrow[from=3-2, to=3-4]
	\arrow[from=3-2, to=4-5]
	\arrow[from=3-4, to=3-7]
	\arrow[from=3-7, to=2-8]
	\arrow[from=3-7, to=4-8]
	\arrow[from=4-5, to=4-8]
\end{tikzcd}.\end{center}
It remains to remark, that \[Z^k\subseteq Z^{j_1}\cap Z^{j_2}.\] This follows directly, as for any $z\in Z^k$ (i.e. $f^k_{\phi(k)}(\pi_{\phi(k)}(z))=t_k$) fulfills
\[f^{j_1}_{\phi(j_1)}(z)=\tau^{j_1}_k t_k=t_{j_1}\] and
\[f^{j_2}_{\phi(j_2)}(z)=\tau^{j_2}_k t_k=t_{j_2}.\]
Thus there exists $x\in\bigcap Z^j$, which clearly fulfills $f(x)=t$.\qedhere
\end{enumerate}
\end{proof}
\begin{remark}
  Note that without the surjectivity weird things might happen.
  Take the limit $\{1,2\}$ of the left column in the following diagram and let $f$ be given by $1,2\mapsto1$, $3,4\mapsto3$.
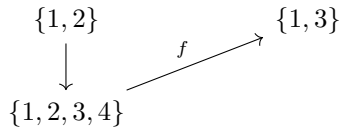
\begin{center}\begin{tikzcd}
	{\{1,2\}} && {\{1,3\}} \\
	{\{1,2,3,4\}}
	\arrow[from=1-1, to=2-1]
	\arrow["f", from=2-1, to=1-3]
\end{tikzcd}\end{center}
By the double limit formula,
this induces the map $\{1,2\}\to\{1,3\}$, $1,2\mapsto1$
which is not surjective.
\end{remark}

\subsection{Weight, size and cutoff cardinals}

In this chapter we relate the cardinality of the diagram \(I\) of a profinite set
\[
    "\varprojlim_{i\in I}" S_i
\]
to the topological weight and the cardinality of the corresponding stone space \(T = \varprojlim_{i\in I} S_i\).
Then we recall important results about topological spaces with certain weights.

\begin{definition}[Strong limit cardinal]\label{def:strong_limit_cardinal}\chaptwo
A \idx{strong limit cardinal} is an infinite cardinal $\kappa$
such that for all $\lambda <\kappa$ we have $2^\lambda<\kappa$.
\end{definition}

The smallest of these cardinals is given by $\aleph_0$,
the next one is the limit (union) of $\aleph_0, 2^{\aleph_0},2^{2^{\aleph_0}},\dots$, and one can proceed similarly from there.
Importantly,
strong limit cardinals are final in all cardinals.

\begin{definition}[Cofinality]\label{def:cofinality}\uses{def:final}\chaptwo
The \idx{cofinality} of a cardinal $\kappa$ is defined as
\[
    \cof(\kappa):=\min\{\lambda: \text{there exists } A\subseteq \kappa \text{ cofinal},\,|A|=\lambda\}.
\]
\end{definition}

\begin{definition}[Weight and size]\label{def:weight}\uses{def:profin,thm:stone}
\begin{enumerate}[(i)]
\item The \idx{size $\kappa$} of a profinite set $T$ is the cardinality of the underlying set, $\kappa=|T|$.
\item The \idx{weight $\lambda$} of a profinite set $T$ is the minimal cardinality $\lambda=|I|$ such that $T\simeq \varprojlim_{i\in I}S_i$.
\item For a topological space $X$, the \idx{topological weight $\omega$}
        is the minimal size of a basis of the topology on $X$.
        \footnote{In this section,
        we never take $\omega$ to meand the least infinite ordinal.}
\end{enumerate}
\end{definition}

\begin{example}\uses{def:profin, thm:stone}\label{ex:profinite}
\begin{enumerate}[(i)]
\item Clearly, $\alpha\N$ has $\lambda=\kappa=|\N|$.
\item The Cantor set has $\lambda=|\N|$ but $\kappa=2^{|\N|}$.
In general $\{0,1\}^I$ has $\kappa=2^{|I|}$ but $\lambda=\omega=|I|$.
\item $\beta\N$ fulfills $\lambda=2^{|\N|}$ but $\kappa=2^{2^{|\N|}}$.
More generally, for any set $I$ the Stone-\v{C}ech compactification $\beta I$
has $\kappa=2^{2^{|I|}}$ and $\lambda = \omega = 2^{|I|}$.
\item Note that for finite discrete sets always $\omega = \kappa$ and $\lambda = 1$ holds.
\end{enumerate}
\end{example}

\begin{remark}[Trivial weight and size estimates]\label{rem:triv_weight_estimates}\uses{def:weight}
Note that for any T1-space $X$ we have $\kappa\le 2^\omega$ and $\omega\le 2^\kappa$.
\end{remark}

\begin{proof}
Consider any base of the topology $\mathcal{V}$.
We obtain an embedding
\[
    X\ni x\mapsto\{U\in \mathcal{V}:x\notin U\}\in\mathcal{P}(\mathcal{V}).
\]
The second inequality is clear.
\end{proof}

\begin{proposition}[Remark of O. Gabber about weights]\label{prop:weight}\uses{lem:weight_clopens, thm:stone, def:weight}
If $\kappa$ is infinite, for any profinite set $T$ the inequality
\[
	\lambda \le \kappa\le 2^\lambda
\]
holds.
\end{proposition}
\begin{proof}
To see the second inequality note that
\[
    T\hookrightarrow \prod_{I}S_i\hookrightarrow \N^I
\]
or alternatively
\[
    \kappa=|\hom(A,\mathbb{F}_2)|\le 2^\lambda
\]
using the lemma below this statement.

For the first inequality consider the set $P=T\times T\setminus \Delta$,
which has cardinality $\kappa^2=\kappa$.
For any distinct pair $p=(t_1, t_2)\in P$ choose a clopen set $A_{p}\subseteq T$ separating $t_1$ and $t_2$ (i.e. $t_1\in A_p\not\ni t_2$).
Now consider
\[
	I=\{F\subseteq P, \quad F \,\text{finite}\},
\]
which is naturally filtered via inclusion and has cardinality $\kappa$.
For each $i\in I$ define the finite subalgebra
\[
	A_i\coloneqq \langle \{A_p : p\in i\}\rangle \le \mathrm{clopen}(T).
\]
Equipped with transition maps $\hom(A_i,\mathbb{F}_2)\to \hom(A_j,\mathbb{F}_2)$ for $j\subset i$ given by restriction,
the sets $\hom(A_i,\mathbb{F}_2)$ form a projective system (i.e. are cofiltered) of finite sets and hence we can form
\[
	S=\varprojlim_I \hom(A_i,\mathbb{F}_2)
\]
and interpret it in $\stone$.

Now for any $i$ the map $T\to \hom(A_i,\mathbb{F}_2)$ (in $\stone$) given by
\[
    t\mapsto \left(a\mapsto \begin{cases}1,& t\in a \\ 0,&\text{else}\end{cases}\right),
\]
is continuous and induces a cone, as it obviously commutes with the transition maps. Hence they induce a map
\[
    f\colon T\to S=\varprojlim_{I}\hom(A_i,\mathbb{F}_2)
\]
in $\stone$.
This map is injective, as for two $t_1\ne t_2$ the value at coordinate $i=\{(t_1,t_2)\}$ evaluated at $A_{(t_1,t_2)}$ is
\[
    f(t_1)_i(A_{(t_1,t_2)})=1\ne 0=f(t_2)_i(A_{(t_1,t_2)}).
\]
Furthermore, any $(a_i)_{i\in I}\in S$ can be identified with $(\min a_i^{-1}(1))_{i\in I}$ where the minimum is taken in each \(A_i\),
which is a family of nonempty clopen sets such that for all finite $J\subset I$
we find an upper bound $j\subset k$ for all $j\in J$, and hence
\[
    \emptyset\ne \min a_k^{-1}(1)\to \bigcap_{i\in J}\min a_i^{-1}(1),
\]
implying that $(\min a_i^{-1}(1))_{i\in I}$ has finite intersection property.
Thus by compactness of $T$ there is an element $\tilde{t}\in \bigcap_{i\in I}\min a_i^{-1}(1)$,
which now obviously fulfills
\[
    f(\tilde{t})=(a_i)_{i\in I}.
\]
Thus $f$ is bijective, which (we are in $\stone$!!) implies
\[
	S\simeq T,
\]
and therefore $\lambda \le \kappa$.
\end{proof}

\begin{remark}
  Note that we just showed that for any point-separating set of finite subalgebras of clopen sets we have $T\simeq \varprojlim \hom(A_i,\f)$.
\end{remark}

\begin{remark}
The estimate $\omega\le |X|$ holds more generally for any locally compact space,
see 3.1.21 and 3.3.6 in \cite{engelking1989general}.
\end{remark}

\begin{lemma}[Weight and size via clopens]\label{lem:weight_clopens}\uses{thm:stone,def:weight}
For any infinite profinite set $T$ we have
\[
    \kappa=|\hom(A,\mathbb{F}_2)|\quad
    \text{ and }\quad\lambda=|A|,
\]
where $A=\mathrm{clopen}(T)$.
Furthermore, $\lambda$ coincides with the topological weight \(\omega\).
\end{lemma}

\begin{proof}
The first statement follows directly when identifying (purely set-theoretically)
\[
    T
    =\varprojlim \hom(A_i,\mathbb{F}_2)
    =\hom(\varinjlim A_i, \mathbb{F}_2)=\hom(A,\mathbb{F}_2).
\]
For the second statement note that clearly $\lambda\le |A|$.

Consider $T=\varprojlim_{i\in I} S_i$, we have to show $|I|\ge |A|$.
But as any clopen set $C\subseteq T$ identifies as a cube $\pi_i^{-1}(N_i)$ with some $N_i\subseteq S_i$,
we have that $A$ embeds into

\[
    \bigsqcup_{i\in I} S_i.
\]
To see that $\lambda$ agrees with the topological weight,
note that in $\stone$ the cardinality of any basis $\mathcal{V}$ of the topology has to be at least as big as $|A|$.
This follows, as every clopen set $C\subseteq T$ can be written as union of elements $V\in \mathcal{V}$.
But as $C$ is compact, we conclude that the union is finite, and hence $A$ embeds into
\[
    \bigsqcup_{i\in\N} \mathcal{V}^i,
\]
thus $|\mathcal{V}|\ge |A|$.
\end{proof}

\begin{definition}[Profinite sets of bounded weight]\label{def:lambda_profinite}\uses{def:profin, def:weight}
For any (in this context always assumed to be infinite) cardinal $\lambda$ define $\prolam$ to be the full subcategory of $\profin$ consisting of profinite sets of weight less than $\lambda$.
\end{definition}
\begin{definition}[Light profinite set]\label{def:light_profinite}\uses{def:weight,def:lambda_profinite}
A profinite set is called \textit{light}, if its weight is at most countable.
\end{definition}

\begin{corollary}[Characterisation of light profinite sets]\label{cor:light_profinite}\uses{def:light_profinite, thm:stone, lem:weight_clopens,lem:limits_along_final_subdiagrams}
The category of light profinite sets is equivalent to
\begin{enumerate}[(a)]
\item sequential limits of finite sets,
\item metrisable totally disconnected compact Hausdorff spaces,
\item the opposite of countable Boolean algebras.
\end{enumerate}
\end{corollary}
\begin{proof}
For the first statement note that any countable codirected set admits a cofinal sequence.
Restricting the diagram along this is a cofinal functor and thereby does not change the colimit (see lemma~\ref{lem:limits_along_final_subdiagrams}).

For the second we only need to know that those (totally disconnected) compact Hausdorff spaces with a countable base of the topology
(i.e. the second-countable compact Hausdorff spaces)
are precisely the metrisable compact Hausdorff spaces.
This follows from Urysohn's metrization theorem.

The weight can be defined as the size of the corresponding Boolean algebra, hence the third assertion is clear.
\end{proof}

\begin{corollary}[Limits and coproducts in $\prolam$]\label{cor:prolam_limits_coproducts}\uses{def:lambda_profinite}
For any infinite cardinal $\lambda$, the category $\prolam$ admits finite limits and finite coproducts.
\end{corollary}

The following property of $\kappa$-profinite sets will be useful later.
One might also proof this purely topological, but this turns out to be quite fuzzy.
\begin{lemma}\label{lem:kappa_prof_is_kappa_cocomp}
	Every $\kappa$-profinite set is a $\kappa$-cocompact object in $\profin$.
\end{lemma}
\begin{proof} By Lemma \ref{lem:kappa_comp_from_comp}, it suffices to show that finite sets are cocompact.
	But this follows by (covariant) Yoneda from $\fin$ to $\profin$.
\end{proof}

\subsection{Topological aspects of profinite sets}

In this section we examine certain topological properties of profinite sets.
This section is not necessary for the understanding of the text
and can thus be skipped.

\begin{lemma}[Universal topological spaces]\label{lem:topologiekrust}\chaptwo
\begin{enumerate}[(i)]
\item Any profinite set of weight $\omega$ is homeomorphic to a subspace of the Cantor cube $\{0,1\}^\omega$.
\item Any compact Hausdorff space of weight $\omega$ is homeomorphic to a subspace of the Tychonov cube $[0,1]^\omega$.
\item Any T0 space of weight $\omega$ is a subspace of the Alexandrov cube $F^\omega$, where $F$ is the Sierpi\`{n}ski space.
\item Any compact Hausdorff space of weight $\omega$ is a continuous image of a subspace of the Cantor cube $\{0,1\}^\omega$. This subspace also has weight $\omega$.
\item Any light profinite set is a retract of $\{0,1\}^{\aleph_0}$.
\end{enumerate}
\end{lemma}

\begin{proof}
\begin{enumerate}[(i)]
\item It suffices to note that the cardinality of all clopen sets is the same as the cardinality of a basis of clopens, as any clopen is finite disjoint union of clopens in the basis. Stone duality yields the desired embedding. See 6.2.16 in \cite{engelking1989general}.
\item Omitted, but works with the embedding
\[x\mapsto (f(x))_{f\colon X\to[0,1]},\]
where it suffices to take $\kappa$-many $f$ (a separating family).
See 2.3.23 and 3.2.5 in \cite{engelking1989general}.
\item See 2.3.26 in \cite{engelking1989general}.
\item We first notice, that any compact Hausdorff space corresponds to a subspace of $[0,1]^\kappa$, hence it suffices to project onto $[0,1]^\kappa$ with $\{0,1\}^\kappa$. But as ${(\{0,1\}^\N)}^\kappa \simeq \{0,1\}^\kappa$ for any $\kappa\ge \aleph_0$, it is legal to just project from $\{0,1\}^\N$ onto $[0,1]$ and glue these together. Such an epimorphism is given e.g. by dyadic representation of real numbers.
See 3.2.2 in \cite{engelking1989general}.
\item As any light profinite set is homeomorphic to a subspace of the Cantor set it suffices to show that subsets of the Cantor set are retracts. Let $X$ be such a subset. Then we can write $X=\varprojlim_{n\in\N}S_n$ where $S_n=\pi_{n}(X)\subseteq \{0,1\}^n$ and $\pi_n$ is the projection onto the first $n$ coordinates, i.e. the $n$-th projection of $\{0,1\}^\N=\varprojlim_n \{0,1\}^n$.
        Note that for any $\{0,1\}^n$ we want to find a surjective map $f_n\colon \{0,1\}^n\twoheadrightarrow S_n$ being the identity on $S_n$ in such a way, that it is compatible with the projections $\{0,1\}^n\to \{0,1\}^m$. For this, let $(x_k)\in \{0,1\}^n$.
        We order $\{0,1\}^n$ lexicographically and for any $x\in \{0,1\}^n$ choose the minimal $s\in S_n\subseteq\{0,1\}^n$ being at least as big as $x$. Trivially, this satisfies the desired properties.\qedhere
\end{enumerate}
\end{proof}

However not every profinite space is the retract of a Cantor cube.
Even a weaker condition is not fulfilled.

\begin{definition}[Dyadic space]\label{def:dyadic}\chaptwo
Spaces that are the continuous image of a Cantor cube are called \emph{dyadic}.
\end{definition}

As we will see, not every profinite set is dyadic,
see also \cite[3.12.12]{engelking1989general} or \cite{Engelking1965CartesianPA}.
In particular, not every profinite set is the retract of a Cantor cube.

\begin{lemma}[Countably many disjoint opens in Cantor cubes]\label{lem:countable_disjoint_opens}\chaptwo
For any Cantor cube $\{0,1\}^\kappa$ there are at most countably many disjoint open subsets.
\end{lemma}

\begin{proof}
Note that $\{0,1\}^\kappa$ is a compact topological group, and hence admits a Haar measure $\mu$.
As Haar measures of open subsets are strictly positive, and the measure is bounded by one and additive on disjoint subsets,
there cannot exist a family of more than $\aleph_0$ many disjoint opens.
\end{proof}

\begin{remark}
One could as well proof this by hand, just using the description of the product topology by cubes.
\end{remark}

\begin{corollary}\label{cor:dyadic_disjoint_opens}\uses{def:dyadic,lem:countable_disjoint_opens}
Every dyadic space has at most $\aleph_0$ disjoint open subsets.
\end{corollary}

\begin{lemma}[Examples for profinite non-dyadic spaces]\label{lem:profinite_non_dyadic}\uses{def:dyadic,lem:countable_disjoint_opens}\chaptwo
$\alpha(\R_d)$ or $\beta\N\setminus\N$ have $2^{\aleph_0}$ many disjoint clopens and hence cannot be dyadic.
\end{lemma}

\begin{proof}
For any real number $r\in [0,1]$ consider the sequence of decimal representations up to some index, $r_n\in \Q$.
Under a bijection $\Q\to\N$ this yields a sequence of natural numbers.
Define $S_r\coloneqq \overline{\{r_n:\, n\in \N\}}\subseteq \beta\N$.
This is a clopen set in $\beta\N$ and the intersection of any two $S_r, S_t$ is finite, as two different real numbers $r,t$ just coincide up to a finite index.
Hence the intersection $S_r\cap S_t$ in $\beta\N\setminus \N$ is empty, as desired.
\end{proof}

We will now investigate whether or not an extremally disconnected space $X$ can be dyadic.
Here the situation is a little bit more complicated.
For \(X\) to be dyadic it is necessary and sufficient to be the retract of some Cantor cube,
since extremally disconnected spaces are projective, see~\ref{lem:CHaus_enough_proj}.
Anyway, this is of no importance, as this situation cannot occur in the infinite case,
see \cite{Engelking1963RemarksOD}.

\begin{theorem}[Extremally disconnected spaces are not dyadic]\label{thm:extremally_disconnected_not_dyadic}\uses{cor:betaN_not_dyadic, lem:real_valued_dyadic, def:dyadic, def:extr}\chaptwo
There exists no infinite dyadic compact extremally disconnected space, hence no extremally disconnected space is a retract of any cantor cube.
\end{theorem}

To prove this, we first show that \(\beta\N\) is not dyadic.

\begin{lemma}[Real valued functions on dyadic spaces]\label{lem:real_valued_dyadic}\uses{def:dyadic}
Let $X$ be any dyadic space, and $f\colon X\to \R$ be continuous.
Then $f$ depends only on countably many coordinates, i.e. there is a countable subset $I$ such that $f$ factors over $\pi_I(X)$ (which in particular is metrisable as a closed subset of the Cantor set).
\end{lemma}

\begin{proof}
It suffices to show this for the cantor cubes $X$.
Here this is a funny application of Stone-Weierstrass.
Denote by $Z$ the set of all functions $Z\to\R$ depending only on countably many coordinates. This set is obviously a subalgebra of $C(X,\R)$ separating points and containing the constant functions.
Furthermore it is closed subspace, as any sequence of functions $f_n\in Z$ with \enquote{supports} $I_n$ converging uniformly has limit only depending on $\bigcup_n I_n$, which is countable.
\end{proof}

\begin{corollary}\label{cor:betaN_not_dyadic}\uses{def:dyadic, lem:real_valued_dyadic}
$\beta\N$ is not dyadic.
\end{corollary}

\begin{proof}
Choose $f\colon \beta\N\to [0,1]$ via $f(n)=1-1/n$.
Then any closed subset $X$ of $\beta\N$ with $f(\beta\N)=f(X)$ has to contain $\N$ and thereby be $\beta\N$, which is not metrisable.
Hence $f$ does not factor over some metrisable subspace of $\beta\N$, which would lead to a contradiction of the last lemma, if $\beta\N$ would be dyadic.
\end{proof}

Now we are able to prove the theorem.

\begin{proof}
Any extremally disconnected infinite compact space admits a surjection on $\beta\N$ by choosing am infinite sequence of disjoint clopens $A_n\subseteq X$ and extending the map $A_n\mapsto \N$ to $\beta\bigcup A_n=\overline{\bigcup A_n}\to \beta\N$. Thus it suffices to know that $\beta\N$ is not dyadic.
\end{proof}

So we have just shown that there exist many profinite spaces which are not dyadic.
If one restricts the spaces to those having countable topological weight,
then the situation becomes much more clear.

\begin{proposition}[Injectives in $\profin$]\label{prop:injective_profinite}\uses{def:profin, lem:topologiekrust}
The injective objects in $\profin$ are precisely the retracts of Cantor cubes.
In particular, all light profinite sets are injective objects in $\profin$.
\end{proposition}

\begin{proof}
It suffices to show injectivity of the set $\{0,1\}$, as $\{0,1\}$ is cogenerating in $\stone$. Here we have to show
\[\hom(S,\{0,1\})\twoheadrightarrow \hom(T,\{0,1\})\]
whenever $S\hookrightarrow T$.
For this we have to note that we can extend any clopen subset $U$ of $S$ to a clopen subset $O$ of $T$ with $O\cap S=U$.
For this note that ${S\setminus U}^c\subseteq T$ is open in $T$ and contains $U$.
As $U$ is compact it can be covered with finitely many clopens $K_i\subseteq T$ with $K_i\subseteq {S\setminus U}^c$.
Taking the union over all $K_i$ yields the desired clopen set.
\end{proof}

In general, it is not true that dyadic spaces are retracts of Cantor cubes,
see the discussion after theorem 16 in \cite{Engelking1965CartesianPA}.

Concerning projective objects, we obtain the following topological lemma.

\begin{lemma}[Projectives in $\CHaus$]\label{lem:CHaus_enough_proj}\uses{lem:adjoints_preserve_projective_objects,def:enough_projectives, lem:projectivity_via_sections, def:weight, lem:topologiekrust, def:extr}
Let $\kappa$ be any cardinal number.
\begin{enumerate}[(i)]
\item Let $X$ be a nonempty compact Hausdorff space of weight $\kappa$.
Then there exists a surjection $\beta\kappa\to X$.
Furthermore, $\beta\kappa$ has weight $2^\kappa$.

\item The category $\CHaus$ has enough projectives and $\profin$ as well.

\item The projective objects in $\CHaus$ are precisely the retracts of Stone Cech compactifications of discrete sets.
They coincide with the extremally disconnected compact Hausdorff spaces.

\item Every nonempty extremally disconnected compact Hausdorff space
of weight $\kappa$ is a retract of $\beta\kappa$.
\end{enumerate}
\end{lemma}

\begin{proof}
\begin{enumerate}[(i)]
\item Since $X$ has weight $\kappa$ there exists a map $\kappa\to X$ with dense image
where we regard $\kappa$ as topological space with the discrete topology.
Then there exists a continuous lift $\beta\kappa\to X$ with dense image,
but the image is closed an so the map is surjective.
The closures of subsets of $\kappa$ form a basis of open sets, hence the weight of $\beta\kappa$ is $2^\kappa$.

\item Follows from (i) and since $\beta\kappa$ is projective in $\CHaus$ by lemma~\ref{lem:projectivity_via_sections}, again using the universal property of $\beta\kappa$.

\item Since $\CHaus$ has enough projectives by (ii), the first claim follows from lemma~\ref{lem:enough-projectives-separating}.
The second claim is a theorem of \cite{Gleason1958}.

\item Since extremally disconnected compact Hausdorff spaces are projective in $\CHaus$ by (iii), the surjection from (i) admits a section.
\end{enumerate}
\end{proof}

\section{Towards condensed sets}

Having studied the building blocks of condensed sets, namely the profinite sets,
we are ready to make a first step towards its definition.

\subsection{Defining k-condensed sets}

In this chapter we install the finitary Grothendieck topology from above on the profinite sets
to be able to talk about sheaves on them.
We present a characterization of these sheaves.

In the following, \(\kappa\) will be an infinite cardinal number.

\begin{definition}[Finitary Grothendieck topology]\label{def:finitary_grothendieck_topology}\uses{def:coverage, def:site, def:profin}
Define a coverage on $\profin$, $\extr$ and $\CHaus$ (with and without cardinal bound) by taking as coverings of objects $T$ finite families
$(S_i\to T)_{i\in F}$ for which
\[
    \coprod_{i\in F} S_i \twoheadrightarrow T
\]
is surjective.
\end{definition}

From now on,
we also install the Grothendieck topology generated by this coverage on \(\profin\), \(\extr\) and \(\CHaus\)
(with and without cardinal bound).

\begin{remark}
Note that a surjective morphism here is exactly an epimorphism,
and even equivalently, a regular epimorphism.
This is the coherent topology from example C2.1.12 (d) in \cite{Johnstone}.
In Chapter D3.3 there are more properties of this topology and its sheaf category.
\end{remark}

\begin{lemma}\label{lem:finitary_forms_coverage}\uses{def:finitary_grothendieck_topology}
This coverage is a coverage in the sense of definition~\ref{def:coverage}.
\end{lemma}

\begin{proof}
Let $(f_i\colon S_i\to T)$ cover $T$ and let $g\colon R\to T$ be a morphism.
Then the pullback family $(f_i^*\colon S_i\times_T R\to R)$ is again a cover in the case we are in $\prof$ or $\CHaus$, since there finite coproducts and regular epimorphisms are stable under pullback.
If we are in $\extr$, take surjections $g_i\colon E_i\to S_i\times_T R$ where $E$ is extremally disconnected.
Then the family $(f_i^*g_i\colon E_i\to R)$ covers $R$ since epimorphisms do not change colimits.
\end{proof}

\begin{remark}
Note that the coverings of the induced Grothendieck topology can be described precisely as the sieves that contain a finite jointly surjective subset,
although they do not necessarily have a finite jointly surjective subset that \emph{generates} the given sieve.

For example consider the Cantor space $C$.
Every map $f\colon C\to \beta\N$ has finite image in $\N$,
as due to cardinality reasons there are no surjections from subsets of $C$ onto $\beta\N$,
and if the image in $\N$ would be infinite,
then the restriction of $f$ to $f^{-1}(\overline{f(C)\cap \N})$ would surject onto $\overline{f(C)\cap \N}\simeq \beta\N$.
Now take any epimorphism $p\colon \beta\N\to C^2$.
The idea is to build any non-finitely generated filter as e.g. all constants $(c_y\colon C\to C, x\mapsto y)_{y\in C}$
and multiply it with anything that ensures infinite image,
e.g. the identity $\mathrm{Id}\colon C\to C$.
Now the sieve generated by $\{c_y\times \mathrm{Id}\mid y\in C\}\cup \{p\}$ contains a surjective morphism but no finite generating family.

Another example is given by taking $i_1,i_2\colon \beta\N\to \beta\N\sqcup\beta\N$
and many non-surjective non-finitely generated maps with contributions in both parts of $\beta\N\sqcup\beta\N$
(e.g. all maps $\{1,2\}\to \beta\N\sqcup\beta\N$ with images in different parts).
This shows that this remark remains valid even in the subcategory of extremally disconnected spaces.
\end{remark}
\begin{lemma}[Initial set of coverings]\label{lem:initial_coverings}
    Consider any profinite set $S$.
    In $\prof$, the set of coverings of $S$ consisting of finite jointly surjective morphisms $(f_i\colon \beta S_d\to S)$ is an initial subset of the category of all coverings.
	If $S$ is $\extr$, then finite clopen decompositions $S=\bigsqcup_{i=1}^n S_i$ form an initial subset of coverings.

    Hence one can restrict the sheafification formula \ref{thm:sheafification} to such covers.
    \end{lemma}
\begin{proof}
    Consider any finite jointly surjective family $(f_i\colon S_i\to S)$.
    Define $g_i$ as the pullback of $f_i$ along $\beta S_d\to S$,

\begin{center}\begin{tikzcd}
	{S_i\times_S\beta S_d} & {\beta S_d} \\
	{S_i} & S
	\arrow["{g_i}", from=1-1, to=1-2]
	\arrow[two heads, from=1-1, to=2-1]
	\arrow[two heads, from=1-2, to=2-2]
	\arrow["{f_i}"', from=2-1, to=2-2]
\end{tikzcd}.\end{center}

    As finite coproducts in $\profin$ are stable under pullback, these glue to a morphism $g\colon \bigsqcup (S_i\times_S\beta S_d)=(\bigsqcup S_i)\times_S\beta S_d\to \beta S_d$

\begin{center}\begin{tikzcd}
	{\bigsqcup S_i\times_S\beta S_d} & {\beta S_d} \\
	{\bigsqcup S_i} & S
	\arrow["{g=\sqcup g_i}", two heads, from=1-1, to=1-2]
	\arrow[two heads, from=1-1, to=2-1]
	\arrow[two heads, from=1-2, to=2-2]
	\arrow["{\sqcup f_i}"', two heads, from=2-1, to=2-2]
\end{tikzcd}\end{center}
this morphism $g$ is epic, since pullbacks of epimorphisms in $\CHaus$ are epimorphisms (this may be checked in $\Set$).
Now, since $\beta S_d$ is projective, there exists a section $s\colon \beta S_d\to \bigsqcup S_i\times_S\beta S_d$ to $g$.

\begin{center}\begin{tikzcd}
	{\bigsqcup S_i\times_S\beta S_d} & {\beta S_d} \\
	{\bigsqcup S_i} & S
	\arrow["g", shift left, two heads, from=1-1, to=1-2]
	\arrow[two heads, from=1-1, to=2-1]
	\arrow["s", shift left, from=1-2, to=1-1]
	\arrow[two heads, from=1-2, to=2-2]
	\arrow["{ f_i}"', from=2-1, to=2-2]
\end{tikzcd}\end{center}
    Defining $A_i=s^{-1}(S_i\times_S \beta S_d)$ yields a clopen decomposition of $\beta S_d$.
    But clopen sets in $\beta S_d$ are precisely of the form $A_i=\overline{B_i}=\beta B_i$ for $B_i\sub S_d$,
	and hence they admit an epimorphism $\beta S_d\twoheadrightarrow A_i$ (since there exists a section $S_d\to B_i$ in $\Set$).

\begin{center}\begin{tikzcd}
	&& {\beta S_d} \\
	{S_i\times_S\beta S_d} & {\beta S_d} & {A_i} \\
	{S_i} & S
	\arrow[two heads, from=1-3, to=2-3]
	\arrow[two heads, from=2-1, to=3-1]
	\arrow["s", from=2-2, to=2-1]
	\arrow[two heads, from=2-2, to=3-2]
	\arrow[hook', from=2-3, to=2-2]
	\arrow["{ f_i}"', from=3-1, to=3-2]
\end{tikzcd}\end{center}
Take $h_i\colon \beta S_d\to S$ to be this composition.
    Then clearly, every $h_i$ factors through $f_i$, and furthermore, they are jointly surjective, as they induce the map $\bigsqcup \beta S_d\twoheadrightarrow \beta S_d\twoheadrightarrow S$.

    If $S$ is extremally disconnected, then we don't need to pullback with $\beta S_d$, but rather have a direct section to $\bigsqcup S_i\twoheadrightarrow S$, inducing a suitable clopen decomposition of $S$.
    \end{proof}
    
    Now we are ready to define $\kappa$-condensed sets.
\begin{definition}[$\kappa$-condensed sets]\label{def:cond_k}\uses{lem:finitary_forms_coverage}
Define the category of $\kappa$-condensed sets, $\condk$, to be $\Sh(\prolak)$,
see definitions \ref{def:sheaf} and \ref{def:D-sheaf},
i.e. contravariant functors $T\colon \prolak\to \set$
such that for every finite covering $(U_i\to U)_{i\in I}$ the diagram 
\begin{center}\begin{tikzcd}
	{T(U)} & {\prod\limits_{i\in I} T(U_i)} & {\prod\limits_{i,j\in I} T(U_i\times_U U_j)}
	\arrow[from=1-1, to=1-2]
	\arrow[shift right, from=1-2, to=1-3]
	\arrow[shift left, from=1-2, to=1-3]
\end{tikzcd}\end{center}
is an equalizer diagram.
\end{definition}

\begin{lemma}[$\kappa$-condensed sets via subbasis]\label{lem:cond_k_on_subbasis}\uses{def:cond_k, lem:stab_sh, def:sheaf}
The $\kappa$-condensed sets are precisely the contravariant functors
$T\colon \prolak\to\set$ such that
\begin{enumerate}[(i)]
\item $T(\emptyset)=\ast$,
\item The canonical map $T(S_1\sqcup S_2)\to T(S_1)\times T(S_2)$ is an isomorphism for all profinite $S_1,\,S_2$, and
\item $T(S)\simeq \{x\in T(S')\mid p_1^*(x)=p_2^*(x)\in T(S'\times_{S}S')\}$
for any surjection $S'\twoheadrightarrow S$,
where $p_1^*=T(p_1)$ and $p_2^*=T(p_2)$ with $p_1,p_2$ being the two projections $S'\times_S S'\to S'$.
\end{enumerate} 
\end{lemma}

\begin{proof}
Let $T$ be a sheaf on $\prolak$.
First of all, the empty covering is a covering of $\emptyset$,
since the empty coproduct is the initial object.
Thus we have the equalizer diagramm 

\begin{center}\begin{tikzcd}
	{T(\emptyset)} & {\prod_\emptyset T(\emptyset)} & {\prod_{\emptyset} T(\emptyset)}
	\arrow[from=1-1, to=1-2]
	\arrow[shift left, from=1-2, to=1-3]
	\arrow[shift right, from=1-2, to=1-3]
\end{tikzcd}.\end{center}
But as $T(\emptyset)\to \prod_\emptyset T(\emptyset)=\ast$ is monic
and every other arrow \(S\to\prod_\emptyset T(\emptyset)\)
equalizes the two arrows between \(\prod_\emptyset T(\emptyset)\),
we know that $T(\emptyset)=\ast$.

As $(S_1\to S_1\sqcup S_2, S_2\to S_1\sqcup S_2)$ is a covering of $S_1\sqcup S_2$,
we have 
\begin{center}\begin{tikzcd}
	{T(S_1\sqcup S_2)} & {T(S_1)\times T(S_2)} & {T(\emptyset)}
	\arrow[from=1-1, to=1-2]
	\arrow[shift left, from=1-2, to=1-3]
	\arrow[shift right, from=1-2, to=1-3]
\end{tikzcd}\end{center}
as $S_1\times_{S_1\sqcup S_2}S_2=\emptyset$.
Therefore $T(S_1\sqcup S_2)\to T(S_1)\times T(S_2)$ is an isomorphism.

For the third property consider an epimorphism $S'\twoheadrightarrow S$, which of course forms a covering of $S$.
Hence the equalizer
\begin{center}\begin{tikzcd}
	{T(S)} & {T(S')} & {T(S'\times_S S')}
	\arrow[from=1-1, to=1-2]
	\arrow[shift right, from=1-2, to=1-3]
	\arrow[shift left, from=1-2, to=1-3]
\end{tikzcd}\end{center}
tells us that $T(S)\to T(S')$ is injective and furthermore that the image precisely consists of the elements $x$ of $T(S')$ such that $p_1^*(x)=p_2^*(x)\in T(S'\times_S S')$.

Now on the other hand consider any contravariant functor $T\colon \prolak\to\set$
fulfilling the three conditions above.
We have to check the sheaf condition on the coverage of finitely jointly surjective families $(f_i\colon U_i\to U)$.
First of all, the canonical map
\[
    T(\bigsqcup U_i)\to \prod T(U_i)
\]
is an isomorphism given by gluing together $e_i^*=T(e_i)$
where $e_i\colon U_i\to \bigsqcup U_i$.
Since $f\colon \bigsqcup U_i\twoheadrightarrow U$ is epic we know that
\begin{center}\begin{tikzcd}
	{T(U)} & {T(\bigsqcup U_i)} & {\prod T(U_i)}
	\arrow[hook, from=1-1, to=1-2]
	\arrow["\sim", from=1-2, to=1-3]
\end{tikzcd}\end{center}
is the equalizer of ${\tilde p_1}^*$ and ${\tilde p_2}^*$ of
\begin{center}\begin{tikzcd}
	{\bigsqcup_i U_i\times_U\bigsqcup_i U_i} & {\bigsqcup_i U_i} \\
	{\bigsqcup_i U_i} & U
	\arrow["{\tilde p_1}", from=1-1, to=1-2]
	\arrow["{\tilde p_2}"', from=1-1, to=2-1]
	\arrow[two heads, from=1-2, to=2-2]
	\arrow[two heads, from=2-1, to=2-2]
\end{tikzcd}.\end{center}

As we are working in $\stone$, finite coproducts are stable under pullback, i.e.,
\[
    \bigsqcup_i U_i\times_U \bigsqcup_i U_i \simeq \bigsqcup_{ij} U_i\times_U U_j,
\]
see example~\ref{ex:CHaus_stable_base_change}.
This yields the resulting equalizer-diagram
\begin{center}\begin{tikzcd}[column sep=small]
	{T(U)} & {\prod T(U_i)} & {T(\bigsqcup U_i)} & {T(\bigsqcup U_i\times_U \bigsqcup U_i)} & {T(\bigsqcup_{i,j} U_i\times_U U_j)} & {\prod T(U_i\times_U U_j)}
	\arrow[hook, from=1-1, to=1-2]
	\arrow["\sim", from=1-2, to=1-3]
	\arrow[shift left, from=1-3, to=1-4]
	\arrow[shift right, from=1-3, to=1-4]
	\arrow["\sim", from=1-4, to=1-5]
	\arrow["\sim", from=1-5, to=1-6]
\end{tikzcd}\end{center}
It remains to check that the maps in this line agree with the maps in the desired equalizer diagram.
This is left to the reader as an exercise.
\end{proof}

\begin{corollary}[Sheaves on $\prof$ via subbasis]\label{cor:sheaves_on_prof_via_subbasis}\uses{lem:cond_k_on_subbasis}
Note that we did not in fact make use of the cardinal bound on \(\prof\),
so the same holds as well for sheaves on \(\prof\), $\extr$ and \(\CHaus\)
with or without cardinal bound.
\end{corollary}

\begin{remark}
We could also reformulate the axioms as follows.
A $\kappa$-condensed set is a contravariant functor $T\colon\prolak\to\set$,
such that:
\begin{enumerate}[(i)]
\item For any finite collection $(S_i)_{i\in I}$ in $\prolak$ the induced map
\begin{center}\begin{tikzcd}
	{T(\bigsqcup_{i\in I}S_i)} & {\prod_{i\in I} T(S_i)}
	\arrow["\sim", from=1-1, to=1-2]
\end{tikzcd}\end{center}
is an isomorphism.
\item For any surjection $f\colon S'\twoheadrightarrow S$ in $\prolak$ the induced map 
\begin{center}\begin{tikzcd}
	{T(S)} & {T(S')}
	\arrow["{f^*}", hook, from=1-1, to=1-2]
\end{tikzcd}\end{center}
is injective with image being precisely the elements $x$ such that $p_1^*(x)=p_2^*(x)$. 
\end{enumerate}
\end{remark}

\subsection{Equivalence to sheaves on Extr or CHaus}

\quot{Suppose you have a topological space that you really like. Then you might want to work in the topos of presheaves on $X$, or the topos of sheaves on $X$.}
{John Baez, in: Topos Theory in a Nutshell\footnote{\url{https://math.ucr.edu/home/baez/topos.html}}}

We have seen that \(\kappa\)-condensed sets are sheaves on \(\prolak\).
In this chapter we will establish an equivalence of this category with the categories of sheaves on \(\CHaus\) and on \(\extr\) with appropriate cardinal bounds.
We will see that on \(\extr\) the sheaf condition becomes particularly simple.

We will apply the comparison lemma~\ref{thm:comparison_lemma} to obtain the desired equivalence.
By lemma~\ref{lem:CHaus_enough_proj}, $\CHaus$ has enough projectives,
so if the cutoff cardinal $\kappa$ has good closure properties,
we obtain that $\extr$ is a dense subsite of $\CHaus$.

\begin{theorem}[Equivalence of sheaves]\label{thm:k_extr_prof_CHaus_equivalence}\uses{def:strong_limit_cardinal,lem:CHaus_enough_proj, thm:comparison_lemma}
The categories of sheaves on 
\begin{enumerate}[(a)]
\item compact Hausdorff spaces $\CHaus_\kappa$ of topological weight less than $\kappa$,
\item profinite sets $\prolak$ of weight less than $\kappa$
\end{enumerate}
are equivalent, and, if $\kappa$ is a strong limit cardinal,
they are equivalent to the category of sheaves on 
\begin{enumerate}[(a)]
\setcounter{enumi}{2}
\item extremally disconnected compact Hausdorff spaces $\extr_\kappa$ of weight less than $\kappa$,
\item Stone-\v{C}ech compactifications $\{\beta I : |I| < \kappa\}$ of discrete sets
    of cardinality less than $\kappa$ with continuous maps.
\end{enumerate}
\end{theorem}

\begin{proof}
The inclusion $\prolak\subset\CHaus_\kappa$ is a dense subsite,
as for any compact Hausdorff space of weight less than $\kappa$,
there is a $\prolak$ set surjecting onto it by lemma~\ref{lem:topologiekrust} (iv)
and example~\ref{ex:profinite} (ii).
Then the comparison lemma~\ref{thm:comparison_lemma} yields the desired equivalence
because $\prolak$ is small.
If $\kappa$ is a strong limit cardinal,
for any profinite set $S$ of weight less than $\kappa$,
there is a surjection $\beta D\twoheadrightarrow S$,
where $D$ has cardinality less than $\kappa$ and hence $\beta D$ has weight less than $\kappa$ by lemma~\ref{lem:CHaus_enough_proj} (i).
Then the comparison lemma yields the equivalences again.
\end{proof}

\begin{remark}
Note that the condition of $\kappa$ being a limit cardinal is necessary,
since the weight of $\beta X_d$ is $2^{|X|}$ in the infinite case.
\end{remark}

\begin{remark}\uses{thm:comparison_lemma, lem:topologiekrust, def:light_profinite}
  As every light profinite set is a retract of the Cantor set,
  a light condensed set is the same thing as a sheaf on the single object site of the Cantor set by the comparison lemma.
  Thus,
  it can be thought of as an underlying set together with an action of $\End(\{0,1\}^{\N})$.
\end{remark}

\begin{remark}
One can even drop the cardinal $\kappa$ and establish the above equivalences of categories of sheaves on $\CHaus$, $\prof$ and $\extr$.
This holds because the needed right Kan extension in the comparison lemma~\ref{thm:comparison_lemma} exists for every sheaf $T\colon\extr^\mathrm{op}\to\Set$ (or on $\prof$)
by lemma~\ref{lem:kan-extension-fully-faithful} since $\extr$ is the colimit of $\extr_\kappa$ along the inclusions over all (strong limit) cardinals (resp. $\prof$).
\end{remark}

Next, we show that on extremally disconnected spaces the sheaf condition is particularly simple.

\begin{proposition}[$\kappa$-condensed sets via extremally disconnected spaces]\label{prop:k_condensed_extremally_disconnected}\uses{def:strong_limit_cardinal, thm:k_extr_prof_CHaus_equivalence, lem:cond_k_on_subbasis, lem:stab_sh}
For a strong limit cardinal $\kappa$ the $\kappa$-condensed sets
are precisely the contravariant functors $T\colon \extr_\kappa\to\set$
that map finite coproducts onto finite products.
\end{proposition}

\begin{proof}
We have to show that such functors indeed form sheaves, as clearly the restriction of sheaves on $\profin_\kappa$ maps finite coproducts to products.
We have to show that for any epimorphism $p\colon S'\twoheadrightarrow S$ the sheaf condition is fulfilled.
Consider any \enquote{compatible} tuple $y\in T(S')$ such that for all $f,g\colon H\to S'$ with $pg=pf$ we have $f^*(y)=g^*(y)$.
We want to show existence of \(z\in T(S)\) such that \(p^*(z) = y\).

Note that by projectivity of $S$ the map $S'\twoheadrightarrow S$ splits,
i.e. has a right inverse $i\colon S\to S'$ with $p\circ i=1_S$.
This implies that $i^*p^*=1$.
Then by choosing \(H = S'\), \(f = ip\) and \(g = 1_{S'}\) we obtain,
since \(pip = p\), that
\[
    p^*(i^*(y)) =y
\]
so \(z = i^*(y)\) does the job.
\end{proof}

\subsection{Embedding topological spaces}

\quot{Many results of pure mathematics, which though likewise apparently fruitless at first, later
become useful in practical science as soon as our mental horizon has been broadened.}{Ludwig Bolzmann (from \cite{Eisner2016})}

Having a \enquote{new concept of topological space},
we wonder how to interpret the classical topological spaces as condensed sets. 

\begin{definition}[Embedding of topological spaces]\label{def:top_to_cond_k}\uses{def:cond_k}
Define a functor $\Top\to\condk$ as follows.
For any topological space $X$ we construct $\underline{X}\in \condk$ via
\[
    S\mapsto C(S,X)
\]
and for $f\colon S\to T$ via $\underline{X}(f)=f^*\colon C(T,X)\to C(S,X)$.
\end{definition}

\begin{lemma}[Topological spaces yield $\kappa$-condensed sets]\label{lem:top_yields_cond_k}\uses{def:top_to_cond_k, lem:cond_k_on_subbasis}
For any topological space $X$, the assignment $\underline{X}$ indeed yields a $\kappa$-condensed set.
\end{lemma}

\begin{proof}
The assignment $S\mapsto C(S,X)$ and $f\mapsto f^*$ clearly induces a contravariant functor,
since composition of continuous functions is continuous.
For any disjoint sum $S=S_1\sqcup S_2$ we have 
\[
    C(S_1\sqcup S_2,X)=C(S_1,X)\times C(S_2,X)
\]
using the $\Hom$-functor in top, and that the coproduct in topological spaces agrees with the disjoint union.
Further, for any epimorphism $\phi\colon S'\twoheadrightarrow S$ we know 
$C(S,X)\hookrightarrow C(S',X)$ via $f\mapsto f\circ \phi$ is injective,
and clearly every $\pi_1^*(f\circ \phi)$ is given by $f\circ\phi\circ\pi_1$ for $\pi_1$ being the coordinate projection of the pullback 
\begin{center}\begin{tikzcd}
	{S'\times_S S'} & {S'} \\
	{S'} & S
	\arrow["{\pi_2}", from=1-1, to=1-2]
	\arrow["{\pi_1}"', from=1-1, to=2-1]
	\arrow["\phi", from=1-2, to=2-2]
	\arrow["\phi"', from=2-1, to=2-2]
\end{tikzcd},\end{center}
thus 
\[
\pi_1^*(\phi^* (f))=\pi_2^*(\phi^*(f)).
\]
On the other hand any $g\colon S'\to X$ that fulfills $g\pi_1=g\pi_2$
is constant on the fibers of $S'\times_S S'$ and therefore by definition of the quotient topology yields a continuous $f$ on $S'/(S'\times_S S')\simeq S$ with
\begin{center}\begin{tikzcd}
	{S'\times_S S'} & {S'} \\
	{S'} & S \\
	&& X
	\arrow["{\pi_2}", from=1-1, to=1-2]
	\arrow["{\pi_1}"', from=1-1, to=2-1]
	\arrow["\phi", from=1-2, to=2-2]
	\arrow["g", from=1-2, to=3-3]
	\arrow["\phi", from=2-1, to=2-2]
	\arrow["g"', from=2-1, to=3-3]
	\arrow["f"{description}, from=2-2, to=3-3]
\end{tikzcd},\end{center}
thus $g=\phi^* f$,
and hence the image of $\phi^*$ is precisely the equalizing set of $\pi_1^*$ and $\pi_2^*$,
as desired.
\end{proof}

\begin{remark}Note that again the cardinal bound was irrelevant for the proof,
and in fact $\underline{X}$ yields a sheaf on the whole category $\prof$.
\end{remark}

\begin{corollary}\label{prof_subcanonical}\uses{ex:subcanonical, lem:top_yields_cond_k}
Every representable presheaf \(\underline{S} = \hom(-,S)\)
for \(S\) profinite is a sheaf.
In particular, the coverage on $\profin$ is \emph{subcanonical}, see example~\ref{ex:subcanonical}.
\end{corollary}

\begin{theorem}[Topologisation of condensed sets]\label{thm:cond_k_to_top}\uses{lem:top_yields_cond_k, def:adjunction}
The functor $X\to\underline{X}$ admits a left adjoint,
which is given by sending a $\kappa$-condensed set $T$ to the \idx{underlying set} $T(\ast)$ endowed with the quotient topology from
\[
    \bigsqcup_{\underline{S}\to T}S\to T(\ast),
\]
where $S$ runs over all $\kappa$-profinite sets.
We denote this topological space by $T(\ast)_{\kappa\mathrm{-top}}$.
\end{theorem}

\begin{proof}
We follow the proof of proposition 1.7 in \cite{scholze2019condensed}.
Denote by $X$ any topological space, and $T$ any $\kappa$-condensed set. 
We want to construct a natural isomorphism
\[
	\hom_\condk(T, \underline{X})\simeq \hom_\Top(T(\ast)_{\kappa\mathrm{-top}}, X).
\]
We prove it in four steps.
\begin{enumerate}
	\item We construct an injective map $\hom(T, \underline{X}) \to \hom_{\mathrm{set}}(T(\ast)_{\kappa-\mathrm{top}}, X)$.
	\item We show its image is a subset of the continuous functions.
	\item We show that the map is surjective.
	\item We show that the map yields a natural isomorphism.
\end{enumerate}
\begin{enumerate}
\item We obviously have a map \(\hom_\condk(T, \underline{X})\to\hom_{\mathrm{set}}(T(\ast), X)\) since \(\underline{X}(\ast) = X\) as sets.
So it suffices to show that all morphisms $T\to \underline{X}$ are determined by its value $T(\ast)\to \underline{X}(\ast)$ on the point \(\ast\).
Note, that for any profinite $S$ and any element $s\in S$ we have an embedding
$i_s\colon\ast{\hookrightarrow} S$, inducing the morphism $\tau$ as follows
\begin{center}\begin{tikzcd}
	{\underline{X}(S)} & {\prod_{s\in S}\underline{X}(\ast)} \\
	{\underline{X}(\ast)}
	\arrow["\tau", from=1-1, to=1-2]
	\arrow["{i_s^*}"', from=1-1, to=2-1]
	\arrow["{\pi_s}", from=1-2, to=2-1]
\end{tikzcd}.\end{center}
Because $i_s^*(f)=f(i_s(\ast))=f(s)$ the morphism $\tau$ maps any $f\in \underline{X}(S)=C(S,X)$ to $(f(s))_{s\in S}$.
Hence \(\tau\) is injective, and nothing else than the inclusion from continuous maps $S\to X$ to all maps $S\to X$.
We simplify the diagram to 
\begin{center}\begin{tikzcd}
	{\underline{X}(S)} & {\prod_{s\in S}X} \\
	X
	\arrow["\tau", from=1-1, to=1-2]
	\arrow["{i_s^*}"', from=1-1, to=2-1]
	\arrow["{\pi_s}", from=1-2, to=2-1]
\end{tikzcd}\end{center}
Now for any morphism $\eta\colon T\to \underline{X}$ the component $\eta(S)$ is fully determined by $\eta(\ast)$ via 
\begin{center}\begin{tikzcd}
	{T(S)} & {\prod_{s\in S} T(\ast)} \\
	{\underline{X}(S)} & {\prod_{s\in S}X}
	\arrow[from=1-1, to=1-2]
	\arrow["{\eta(S)}"', from=1-1, to=2-1]
	\arrow["{\prod \eta(*)}", from=1-2, to=2-2]
	\arrow["\tau", from=2-1, to=2-2]
\end{tikzcd}\end{center}
Hence the assignment $\eta\mapsto \eta(\ast)$ is a good candidate for the natural isomorphism,
i.e. it $\eta\mapsto \eta(\ast)$ is an injective map
$\hom(T,\underline{X})\to \hom_{\Set}(T(\ast),\underline{X}(\ast))$.
\item Now we show that this injection indeed has image inside of continuous maps.
Consider any $f\in T(S)$.
By the Yoneda lemma this is a natural transformation
\begin{center}\begin{tikzcd}
	{\underline{S}(S')} & {\underline S(S'')} \\
	{T(S')} & {T(S'')}
	\arrow[from=1-1, to=1-2]
	\arrow["{f_{S'}}"', from=1-1, to=2-1]
	\arrow["{f_{S''}}", from=1-2, to=2-2]
	\arrow[from=2-1, to=2-2]
\end{tikzcd}.\end{center}
Now, $\eta(S)(f)$ corresponds precisely to the composed natural transformation
\begin{center}\begin{tikzcd}
	{\underline{S}(S')} & {\underline S(S'')} \\
	{T(S')} & {T(S'')} \\
	{\underline{X}(S')} & {\underline X(S'')}
	\arrow[from=1-1, to=1-2]
	\arrow["{f_{S'}}"', from=1-1, to=2-1]
	\arrow["{f_{S''}}", from=1-2, to=2-2]
	\arrow[from=2-1, to=2-2]
	\arrow["{\eta_{S'}}"', from=2-1, to=3-1]
	\arrow["{\eta_{S''}}", from=2-2, to=3-2]
	\arrow[from=3-1, to=3-2]
\end{tikzcd}.\end{center}
But this induces a map 
\begin{center}\begin{tikzcd}
	{\underline{S}(S)} & {\prod_{s\in S} S} \\
	{T(S)} & {\prod_{s\in S} T(\ast)} \\
	{\underline{X}(S)} & {\prod_{s\in S}X}
	\arrow[from=1-1, to=1-2]
	\arrow["{1\mapsto f}"', from=1-1, to=2-1]
	\arrow[from=1-2, to=2-2]
	\arrow[from=2-1, to=2-2]
	\arrow["{\eta(S)}"', from=2-1, to=3-1]
	\arrow["{\prod \eta(*)}", from=2-2, to=3-2]
	\arrow["\tau", from=3-1, to=3-2]
\end{tikzcd}\end{center}
Plugging in $1\in \underline{S}(S)$ tells us precisely
that $\eta(S)(f)\colon S\to X$ being continuous
implies
\[
    S\overset{f}{\to} T(\ast)\overset{\eta(\ast)}{\to} X
\]
has to be continuous. 
Thus $\eta(\ast)$ is continuous with respect to the quotient topology of all $S\to T(\ast)$, which is the topology of $(T(\ast))_{\kappa\mathrm{-top}}$.

Thus $\eta\mapsto \eta(\ast)$ is an injection
\[
    \hom_\condk(T,\underline{X})\to \hom_\Top(T(\ast)_{\kappa\mathrm{-top}},X).
\]

\item To see surjectivity consider any $g(\ast)\in\hom_\Top(T(\ast)_{\kappa-\mathrm{top}},X)$.
We build a natural transformation $g$ out of $g(\ast)$ by defining $g(S)$ as follows: Any $x\in T(S)$ corresponds to a natural transformation (Yoneda):
\begin{center}\begin{tikzcd}
	& x \\
	{\underline{S}(S)} & {T(S)} & {\underline{X}(S)} \\
	{\underline{S}(\ast)} & {T(\ast)} & {\underline{X}(\ast)}
	\arrow["\in", from=1-2, to=2-2]
	\arrow["{x(S)}", from=2-1, to=2-2]
	\arrow[from=2-1, to=3-1]
	\arrow[from=2-2, to=2-3]
	\arrow[from=2-2, to=3-2]
	\arrow[from=2-3, to=3-3]
	\arrow["{x(\ast)}"', from=3-1, to=3-2]
	\arrow["{g(\ast)}"', from=3-2, to=3-3]
\end{tikzcd}.\end{center}
Now define $g(S)(x)$ as the map in the bottom row,
i.e. $g(S)(x) = g(\ast)\circ x(\ast)$. 
This is indeed continuous from $S$ to $X$,
as $g(\ast)$ is continuous for the topology on $T(\ast)_{\kappa\mathrm{-top}}$.
Thus $g(S)$ is well-defined. 

To see that it is a natural transformation,
consider any continuous $\phi\colon S'\to S$. 
This corresponds to a natural transformation $\underline{S'}\to \underline{S}$
which acts by pushforward $\phi_*$.
We are interested in the commutativity of the upper right square:
\begin{center}\begin{tikzcd}
	&& x \\
	{\underline{S}'(S)} & {\underline{S}(S)} & {T(S)} & {\underline{X}(S)} \\
	{\underline{S}'(S')} & {\underline{S}(S')} & {T(S')} & {\underline{X}(S')} \\
	{\underline{S}'(\ast)} & {\underline{S}(\ast)} & {T(\ast)} & {\underline{X}(\ast)}
	\arrow["\in", from=1-3, to=2-3]
	\arrow["{\phi_*}", from=2-1, to=2-2]
	\arrow["{\phi^*}"', from=2-1, to=3-1]
	\arrow["{x(S)}", from=2-2, to=2-3]
	\arrow["{\phi^*}", from=2-2, to=3-2]
	\arrow["{g(S)}", from=2-3, to=2-4]
	\arrow["T\phi", from=2-3, to=3-3]
	\arrow["{\phi^*}", from=2-4, to=3-4]
	\arrow["{\phi_*}", from=3-1, to=3-2]
	\arrow[from=3-1, to=4-1]
	\arrow["{x(S')}", from=3-2, to=3-3]
	\arrow[from=3-2, to=4-2]
	\arrow["{g(S')}", from=3-3, to=3-4]
	\arrow[from=3-3, to=4-3]
	\arrow[from=3-4, to=4-4]
	\arrow["{\phi_*}", from=4-1, to=4-2]
	\arrow["{x(\ast)}", from=4-2, to=4-3]
	\arrow["{g(\ast)}", from=4-3, to=4-4]
\end{tikzcd}\end{center}
We want to compare $\phi^*(g(S)x)$ with $g(S')(T(\phi)x)$.
Clearly, 
\[
	\phi^*g(S)x = \phi^\ast(g(\ast)\circ x(\ast)) = g(\ast)x(\ast)\phi.
\]
But $T(\phi)x\in T(S')$ is precisely the composite natural transformation
$x\circ\phi_\ast$, from the left column into the $T$-column (Yoneda backwards).
Thus $g(S')(T(\phi)x)$ corresponds precisely to the bottom row, i.e. 
\[
	g(S')T(\phi)x = \phi_*(g(\ast)\circ x(\ast)) = g(\ast)x(\ast)\phi.
\]
Thus, $g$ is a natural transformation, and hence the assignment $\eta\mapsto \eta(\ast)$ is bijective.

\item
For the naturality of $\eta\mapsto \eta(\ast)$ in the second argument
consider any continuous map $\phi\colon X\to Y$ and denote by $\underline{\phi}$
the corresponding natural transformation $\underline{X}\to \underline{Y}$.
We want to show that
\begin{center}\begin{tikzcd}
	{\hom(T,\underline{X})} & {\hom(T(\ast)_{\kappa\mathrm{-top}}, X)} \\
	{\hom(T,\underline{Y})} & {\hom(T(\ast)_{\kappa\mathrm{-top}},Y)}
	\arrow[from=1-1, to=1-2]
	\arrow["{\underline{\phi}_{\ast}}", from=1-1, to=2-1]
	\arrow["{\phi_*}", from=1-2, to=2-2]
	\arrow[from=2-1, to=2-2]
\end{tikzcd}\end{center}
is commutative.
For this note that for any $f\colon T\to\underline X$ the assigned map $f(\ast)\colon T(\ast)\to X$ is mapped with $\phi_*$ to $\phi\circ f(\ast)$, and 
$\underline{\phi}\circ f$ has components precisely $\phi \circ f(S)$, hence $(\underline{\phi}\circ f)(\ast)=\phi\circ f(\ast)$.

To see naturality in the first argument consider any $\eta\colon T\to R$.
We want commutativity of
\begin{center}\begin{tikzcd}
	{\hom(T,\underline{X})} & {\hom(T(\ast)_{\kappa\mathrm{-top}}, X)} \\
	{\hom(R,\underline{X})} & {\hom(R(\ast)_{\kappa\mathrm{-top}},X)}
	\arrow[from=1-1, to=1-2]
	\arrow["{\eta^*}", from=1-1, to=2-1]
	\arrow["{\eta(\ast)^\ast}", from=1-2, to=2-2]
	\arrow[from=2-1, to=2-2]
\end{tikzcd}\end{center}
but this is immediate as the map in the right column is precisely the natural transformation in the left column evaluated at $\ast$, as pullback clearly commutes with evaluating a natural transformation, i.e. 
\[
    \eta(\ast)^\ast=\eta^\ast(\ast).
\]
\end{enumerate}
\end{proof}

\begin{remark}
  Again,
  the cardinal bound is irrelevant for the proof.
  Defining $X(*)_{\mathrm{top}}$ in the same way as above (just without cardinal bound on $S$)
  yields an adjunction between $\Top$ and $\Sh(\prof)$.
\end{remark}

\begin{corollary}\label{cor:counit_top_cond_k}\uses{thm:cond_k_to_top}
The counit of the adjunction is given by the continuous bijection
$\underline{X}(\ast)_{\kappa\mathrm{-top}}\to X$, $(*\stackrel{i}{\to}\underline{X})\mapsto i(*)$.
\end{corollary}

\begin{proof}
The counit \((\eps_X)_{X\in\Top}\) arises from taking the identity in
\(\hom_\condk(\underline{X},\underline{X})\) and mapping it via the natural isomorphism of the adjunction into \(\hom_\Top(\underline{X}(\ast)_{\kappa\mathrm{-top}},X)\).
But this is exactly constructed to be the (set theoretic) identity
\(\underline{X}(\ast)_{\kappa\mathrm{-top}}\to X\) which then in turn is continuous.
\end{proof}

\begin{remark}
In particular the \enquote{\(\kappa\)-top topology} on \(X\) is finer than the topology on \(X\).
\end{remark}

\begin{corollary}\label{cor:top_to_cond_k_faithful}\uses{thm:cond_k_to_top, cor:counit_top_cond_k, lem:unit_via_fully_faithfullness}
The functor $X\mapsto\underline{X}$ is faithful from topological spaces to $\kappa$-condensed sets.
\end{corollary}

\begin{proof}
This follows, as the counit is epimorphic, see ~\ref{lem:unit_via_fully_faithfullness}.
Another way to see this is that for $f^*=g^*\colon \underline{X}\to\underline{Y}$
the maps of sets $\underline{X}(\ast)\to \underline{Y}(\ast)$ coincide,
and hence they have to be induced by the same map $X\to Y$.
\end{proof}

\subsection{Finding k-compactly generated spaces}

\quot{All told, this suggests that in $\mathbf{Top}$ we have been studying
the wrong mathematical objects. \\
The right ones are the spaces in $\mathbf{CGHaus}$.}{Saunders Mac Lane, in: \cite[p. 188]{mac2013categories}}

We have seen in the previous chapter that (\(\kappa\)-)compactly generated topologies
play a big role when working with (\(\kappa\)-)condensed sets.
In this chapter we establish some properties of these topological spaces.
All of the results and more on compactly generated spaces can be found in Appendix A of~\cite{Lewis1978}.

\begin{definition}[$\kappa$-compactly generated topological space]\label{def:k_compactly_generated}\uses{def:weight}
Let $\kappa$ be an infinite cardinal number.
Define the category of \idx{$\kappa$-compactly generated topological spaces} $\Top^{\kappa\textrm{-cg}}$
as those topological spaces $X$ that carry the quotient topology of
\[
    \bigsqcup_{S\to X}S\to X,
\]
where $S$ runs over all compact Hausdorff spaces of \emph{weight} less than $\kappa$.
\end{definition}

\begin{remark}
  Because the left adjoint to $\Top\to\condk$, $X\mapsto\underline{X}$ takes values in $\Top^{\kappa\mathrm{-cg}}$,
  the adjunction restricts to an adjunction between $\Top^{\kappa\mathrm{-cg}}$ and $\condk$.
\end{remark}

\begin{example}
  Note that every first countable (e.g. metrisable) space is $\kappa$-compactly generated for every $\kappa$,
  as here continuity may be checked via sequential continuity.
  Further examples include locally compact spaces.
\end{example}

Now, we show, that in the definition of $\kappa$-compactly generated spaces,
instead of testing again \emph{all} compact Hausdorff spaces,
it suffices to test again all profinite sets.

\begin{lemma}\label{lem:comp_gen_quot_top_prof}\uses{ def:weight, def:profin, def:weight, lem:topologiekrust}
Let $\kappa$ be any cardinal.
A topological space $X$ is $\kappa$-compactly generated if and only if it carries the quotient topology of
\[
    \bigsqcup_{S\to X} S\to X,
\]
where $S$ runs over \emph{profinite} sets with weight less than $\kappa$.
Equivalently, the counit of the above adjunction is an isomorphism.
\end{lemma}

\begin{proof}
It suffices to see that for every compact $K$ with weight less than $\kappa$ we find a surjection from some totally disconnected $T$ with weight less than $\kappa$ onto $K$.
This is precisely the statement of the lemma \ref{lem:topologiekrust}.
\end{proof}

\begin{remark}
If $\kappa$ is a strong limit cardinal, in view of lemma~\ref{lem:CHaus_enough_proj},
we could as well restrict to surjections from
Stone-\v{C}ech compactifications of discrete sets with cardinality less than $\kappa$.
In the case of $\kappa$ being the first uncountable cardinal number,
restricting to morphisms from the Cantor set suffices by lemma~\ref{lem:topologiekrust}.
\end{remark}

In particular, lemma~\ref{lem:comp_gen_quot_top_prof} together with the description of the counits of both adjunctions in lemma~\ref{lem:k_compactly_generated_reflective}
and in corollary~\ref{cor:counit_top_cond_k} imply that $\kappa$-compactly generated
spaces form a nice subcategory of $\kappa$-condensed sets.

\begin{corollary}\label{cor:k_compactly_generated_fully_faithful}\uses{lem:comp_gen_quot_top_prof,lem:k_compactly_generated_reflective, cor:counit_top_cond_k, lem:unit_via_fully_faithfullness, lem:comp_of_adjoints}
The functor $X\mapsto\underline{X}$ is fully faithful when restricted to $\kappa$-compactly generated topological spaces.
\end{corollary}

Lastly, we show that the $\kappa$-compactly generated spaces are a coreflective subcategory of $\Top$.

\begin{lemma}\label{lem:k_compactly_generated_reflective}\uses{def:k_compactly_generated, def:reflective_subcategory, thm:FAFT}
The forgetful functor from $\kappa$-compactly generated topological spaces to $\Top$ (and therefore also $\condk$)
admits a right adjoint, sending $X$ to
\[
    X^{\kappa\mathrm{-cg}} = \underline{X}(\ast)_{\kappa\mathrm{-top}}.
\]
The counit is the (not necessarily invertible) identity $\underline{X}(\ast)_{\kappa\mathrm{-top}}\to X$.
Moreover, the counit is an isomorphism precisely if $X$ is $\kappa$-compactly generated.
\end{lemma}

\begin{proof}
We rather give a direct proof of the equality
\[
    C(Y,X)=C(Y,\underline{X}(\ast)_{\kappa\mathrm{-top}})
\]
for all $\kappa$-compactly generated $Y$.
As the topology of $\underline{X}(\ast)_{\kappa\mathrm{-top}}$ is finer than the topology of $X$,
clearly $C(Y,X)\supseteq C(Y,\underline{X}(\ast)_{\kappa\mathrm{-top}})$.
It remains to show that any continuous map $Y\to X$  is continuous with respect to $Y\to \underline{X}(\ast)_{\kappa\mathrm{-top}}$.
As for any $\kappa$-compact set $K$ with map $K\to X$ the map $K\to Y\to X$ is continuous by definition of the topology on $Y$, and thus continuous $K\to Y\to \underline{X}(\ast)_{\kappa\mathrm{-top}}$,
as $\underline{X}(\ast)_{\kappa\mathrm{-top}}$ carries the final topology of all continuous maps $K\to X$.
By definition of the topology on $Y$ this implies $Y\to \underline{X}(\ast)_{\kappa\mathrm{-top}}$ being continuous, and thereby we conclude the equality of both $\hom$ sets.
\end{proof}

\subsection{Eliminating the cardinal bound}

So far, we have only considered condensed sets up to a certain cardinal bound.
In this chapter we lay out how to get rid of this restriction.

\begin{proposition}[$\kappa$-condensed to $\kappa'$-condensed sets]\label{prop:cond_k_to_cond_k_prime}\uses{def:cofinality, lem:adjoints_commute_with_limits, lem:unit_via_fully_faithfullness, thm:computation_of_kan_extensions,prop:filtered_sifted_commutes_set, def:adjunction, def:kan_extensions, thm:sheafification, def:strong_limit_cardinal, prop:k_condensed_extremally_disconnected, thm:comparison_lemma}
Let $\kappa'>\kappa$ be uncountable strong limit cardinals.
Then there is a left adjoint to the restriction $\cond_{\kappa'}\to \cond_\kappa$.
It is given by the functor $\cond_\kappa\to \cond_{\kappa'}$
which sends a condensed set $T$ to the sheafification $T_{\kappa'}$ of its left Kan extension
along the inclusion,
which is given by
\[
    \Tilde S\mapsto \varinjlim_{\Tilde S\to S}T(S).
\]
The functor $T\mapsto T_{\kappa'}$ is fully faithful and commutes with all colimits.
\end{proposition}

\begin{proof}
We follow the proof of 2.9 in \cite{scholze2019condensed}.
Clearly, the restriction of a condensed set on $\prof_{\kappa'}$ to $\prolak$
again yields a condensed set.
Therefore, the restriction functor is given by the upper row in the following diagram.
\begin{center}
\begin{tikzcd}
	\cond_{\kappa'} & {[\prof_{\kappa'},\Set]} & {[\prof_{\kappa},\Set]} & {\condk}
	\arrow[""{name=0, anchor=center, inner sep=0}, "{?}", shift left=2, from=1-1, to=1-2]
	\arrow[""{name=1, anchor=center, inner sep=0}, "\Sh", shift left=2, from=1-2, to=1-1]
	\arrow[""{name=2, anchor=center, inner sep=0}, "{\tau^*}", shift left=2, from=1-2, to=1-3]
	\arrow[""{name=3, anchor=center, inner sep=0}, "{\Lan_\tau}", shift left=2, from=1-3, to=1-2]
	\arrow[""{name=4, anchor=center, inner sep=0}, "\Sh", shift left=2, from=1-3, to=1-4]
	\arrow[""{name=5, anchor=center, inner sep=0}, "{?}", shift left=2, from=1-4, to=1-3]
	\arrow["\dashv"{anchor=center, rotate=90}, draw=none, from=1, to=0]
	\arrow["\dashv"{anchor=center, rotate=90}, draw=none, from=3, to=2]
	\arrow["\dashv"{anchor=center, rotate=-90}, draw=none, from=4, to=5]
\end{tikzcd}
\end{center}
Thus we need to show that the lower row in the diagram,
which is precisely the given map from the proposition,
is indeed its left adjoint.
The three adjunctions indicated in the diagram follow from theorem~\ref{thm:sheafification} and remark~\ref{rem:Kan-adjoint}.
Thus it follows
\[
    \hom((\Sh\circ\Lan_\tau\circ ?) T,S)
    = \hom(\Lan_\tau\circ ? T,?S)
    = \hom(?T,(?S)\circ\tau).
\]
for $T\in\condk$ and $S\in\cond_{\kappa'}$.
Now since $(?S)\circ\tau$ is already a sheaf on $\prolak$,
we further obtain
\[
    \hom(?T,?\circ\Sh((?S)\circ\tau))
    = \hom(T,\Sh((?S)\circ\tau))
    = \hom(T,(\Sh\circ\tau^\ast\circ ?)S)
\]
because the category of sheaves is a full subcategory of presheaves.
The formula for the left Kan extension follows from theorem~\ref{thm:computation_of_kan_extensions}.

The unit of the adjunction is given by mapping any $\kappa$-condensed set to the restriction of the sheafification of its extension.
Sheafification commutes with restriction, since by lemma \ref{lem:initial_coverings} the sheafification formula at some $S\in \prolak$ only uses coverings with elements in $\prolak$, and thereby agrees with the formula of sheafification on $\prolak$.
But since $T$ already was a sheaf on $\prolak$, and the Kan extension does not change the values on $\prolak$ (since $\prolak\hookrightarrow\profin_{\kappa'}$ is fully faithful, using \ref{lem:Kan_extension_extends}),
sheafification does not change the values on $\prolam$.
Hence the unit is an isomorphism and the functor $T\mapsto T_{\kappa'}$ is fully faithful.

Since $T\mapsto T_{\kappa'}$ is a left adjoint, it preserves all colimits by lemma~\ref{lem:adjoints_commute_with_limits}.
\end{proof}

\begin{remark}
  Note that the direction of the arrow $\tilde{S}\to S$ in the index category of the colimit formula is denoted as in $\prolam$ (not the dual)
  although $T$ is, of course, a functor from $\prolak^{\op}$.
\end{remark}

\begin{lemma}\label{lem:cond_k_to_cond_k_prime_alternative}
  For strong limit cardinals $\kappa < \kappa'$,
  identify $\cond_\kappa$ and $\cond_{\kappa'}$ with the corresponding categories of sheaves on $\extr_\kappa$ and $\extr_{\kappa'}$.
\begin{enumerate}[(i)]
\item The embedding $\cond_\kappa\to \cond_\kappa'$ is given by forgetting and Kan extending; no sheafification on $\extr$ is needed.

\item Consider the Yoneda embedding $y^-\colon \extr_\kappa\to\cond_\kappa$
and the inclusion $\extr_\kappa\subset\extr_{\kappa'}$ followed by the Yoneda embedding
$\extr_{\kappa'}\to\cond_{\kappa'}$ which we will call $Y^-$.
Then the left adjoint to the restriction $\cond_{\kappa'}\to\cond_{\kappa}$,
as constructed in the previous proposition, is alternatively given by
left Kan extension of $y^-$ along $Y^-$,
\begin{center}
\begin{tikzcd}
	{\extr_\kappa} && {\cond_\kappa} \\
	\\
	&& {\cond_{\kappa'}}
	\arrow["y", from=1-1, to=1-3]
	\arrow[from=1-1, to=3-3]
	\arrow["{\mathrm{Lan}}", dashed, from=1-3, to=3-3]
\end{tikzcd}
\end{center}

\end{enumerate}
\end{lemma}

\begin{proof}
We prove part (ii).

Theorem~\ref{thm:crit_for_kan} yields existence of the left Kan extension,
since $\extr_\kappa$ is small and $\cond_{\kappa'}$ is cocomplete by lemma~\ref{lem:limits_pointwise}.
Moreover, it is given by
\[
    T\mapsto \varinjlim_{y^E\to T} Y^E.
\]
Since $T$ maps finite coproducts to finite products,
the comma category $y^E\Rightarrow T$ of morphisms from $y^E$ to $T$ (equivalently of elements $x\in T(E)$ by the Yoneda lemma) is sifted as it has finite coproducts:
for $x\in T(E)$ and $y\in T(E')$ the element $(x,y)\in T(E)\times T(E') = T(E\sqcup E')$
is the element above $x$ and $y$.
Henceforth, the pointwise colimit in this formula already is a sheaf.
For $\tilde{S}\in\extr_{\kappa'}$ we obtain
\begin{align*}
         (\varinjlim_{y^E\to T} Y^E)(\tilde{S})
  &=   \varinjlim_{y^E\to T} \hom_{\extr_{\kappa'}}(\tilde{S},E)
    =   \varinjlim_{y^E\to T} \varinjlim_{\tilde{S}\to S}
            \hom_{\extr_{\kappa}}(S,E) \\
    &=   \varinjlim_{\tilde{S}\to S} \varinjlim_{y^E\to T}
      \hom_{\extr_{\kappa}}(S,E)
    =   \varinjlim_{\tilde{S}\to S} \varinjlim_{y^E\to T}
      y^E(S)
      =   \varinjlim_{\tilde{S}\to S} T(S)
\end{align*}
where $\tilde{S}\to S$ runs over all $S\in\extr_{\kappa}$
and the last equality is due to Theorem III.1 in \cite{mac2013categories}.

This also implies (i), since the left side of the equality is a sheaf, proving that the Kan extension already forms a sheaf. 
\end{proof}

\begin{remark}
An alternative (elementary) proof of the fact that on $\extr$ no sheafification is needed will be given later in \ref{lem:sheafification_cond_C}.
\end{remark}

Next, we will show that $\condk\to\cond_{\kappa'}$ also commutes with many limits. 
For this we need the following useful lemma. 
Recall that a diagram is $\lambda$-small, if the set of morphisms has cardinality less than $\lambda$, see definition~\ref{def:diagram}.
A $\lambda$-limit is a limit over a $\lambda$-small diagram.
Recall that a category $\mcI$ is called $\lambda$-filtered,
if the diagonal functor $\mcI\to\mcI^{\mcJ}$ is final for every $\lambda$-small $\mcJ$,
see definition \ref{def:special_limits_and_colimits}~(v).
Equivalently, colimits over $\mcI$ commute with $\lambda$-small limits in $\Set$,
see proposition~\ref{prop:filtered_sifted_commutes_set}.

\begin{lemma}\label{lem:extr-lambda-filtered}
Let $\kappa' > \kappa$ be uncountable strong limit cardinals.
For any extremally disconnected $\tilde{S}$,
the comma category of maps $\tilde{S}\to S$ with $S\in\extr_\kappa$ (i.e., the full subcategory of $(\tilde{S}\downarrow\extr)$) is $\lambda$-cofiltered,
where $\lambda$ is the cofinality of $\kappa$.

The same holds for $\prolak$ instead of $\extr_\kappa$.
\end{lemma}

\begin{proof}
We follow the proof of 2.9 in \cite{scholze2019condensed}.
We have to show that for any index category $(\Tilde S\to S_i)_{i\in I}$ with $S_i\in \extr_\kappa$ of size less than $\lambda$
there is a common upper bound $S\in\extr_{\kappa}$,
i.e., a map $\Tilde S\to S$ such that all $\Tilde S\to S_i$ factor compatibly over $S$,
as shown in the following diagram
\begin{center}
\begin{tikzcd}
	S \\
	{\tilde{S}} & {S_i} \\
	& {S_j}
	\arrow[from=1-1, to=2-2]
	\arrow[from=2-1, to=1-1]
	\arrow[from=2-1, to=2-2]
	\arrow[from=2-1, to=3-2]
	\arrow[from=2-2, to=3-2]
\end{tikzcd}.
\end{center}
A first guess would be $S=\varprojlim S_i$ (taken in $\profin$),
as for this set the universal property precisely yields the map $\Tilde S\to S$ we need.
Now, $\prod_{i\in I}S_i$ is $\kappa$-small,
because all the $S_i$ have weight less than $\kappa$,
and hence their supremum as well since we have less than $\lambda$ many,
and $\lambda$ is the cofinality.
Thus there is a cardinal number $\mu<\kappa$ with $|S_i|<\mu$ for all $i$,
and thereby 
\[
    \left|\prod_i S_i\right|\le \mu^\lambda< \kappa,
\]
as $\kappa$ is a limit cardinal (this estimate can be made sharper).
Hence $\prod_{i\in I}S_i$ is $\kappa$-small and thus as well $S=\varprojlim S_i$ as a closed subset (in case for $\prolak$ we are done here).
Now, take $\beta\kappa$, which is extremally disconnected and projects onto $S$.
As $\kappa$ is a strong limit cardinal, this has weight less than $\kappa$.
By the projectivity of $\Tilde S$ we obtain the desired map:
\begin{center}
\begin{tikzcd}
	{\beta S_d} & S \\
	& {\tilde{S}} & {S_i} \\
	&& {S_j}
	\arrow[two heads, from=1-1, to=1-2]
	\arrow[from=1-2, to=2-3]
	\arrow[dashed, from=2-2, to=1-1]
	\arrow[from=2-2, to=1-2]
	\arrow[from=2-2, to=2-3]
	\arrow[from=2-2, to=3-3]
	\arrow[from=2-3, to=3-3]
\end{tikzcd}.
\end{center}
\end{proof}

\begin{corollary}\label{cor:Lan-lambda-small-limits}
The functor $\cond_\kappa\to\cond_{\kappa'}$ commutes with $\lambda$-small limits, for $\lambda$ the cofinality of $\kappa$.
\end{corollary}
\begin{proof}
Identifying $\cond_\kappa$ and $\cond_{\kappa'}$ with sheaves on $\extr_\kappa$ resp. $\extr_{\kappa'}$, it suffices to show the statement there (where no sheafification is needed).
Now it suffices to show this pointwise for any $\Tilde S\in\prof_{\kappa'}$,
as limits of sheaves are computed pointwise. 
Now, as $\lambda$-filtered colimits in $\Set$ commute with $\lambda$-small limits (see \ref{prop:filtered_sifted_commutes_set}), the Kan extension commutes with these limits, and since no sheafification is needed we are done.
\end{proof}
\begin{remark}
  The same argument hold for Kan extension along $\extr_{\kappa}\hookrightarrow\extr$.
\end{remark}

This implies the following general principle of avoiding the $\kappa$,
\begin{center}
\begin{tikzcd}
	{\Sh(\prof_\kappa)} && {\Sh(\prof)} \\
	\\
	{\Sh(\extr_\kappa)} && {\Sh(\extr),}
	\arrow["{{\Sh\circ \mathrm{Lan}}}", shift left, from=1-1, to=1-3]
	\arrow[shift left, from=1-1, to=3-1]
	\arrow[shift left, from=1-3, to=1-1]
	\arrow[shift right, from=1-3, to=3-3]
	\arrow["{{\mathrm{Ran}}}", shift left, from=3-1, to=1-1]
	\arrow["{{\mathrm{Lan}}}", shift left, from=3-1, to=3-3]
	\arrow["{{\mathrm{Ran}}}"', shift right, from=3-3, to=1-3]
	\arrow[shift left, from=3-3, to=3-1]
\end{tikzcd}
\end{center}
where the left adjoints are labeled,
the right adjoints are restriction and the absence of $\Sh$ in the bottom row is not a typo.

\begin{question}
We are not sure if sheafification is even needed in the \(\prof\) setting.
\end{question}

\section{The category of condensed sets}

We are now ready to define condensed sets.
\begin{definition}[Condensed sets]\label{def:cond}
  Define the category of \idx{condensed sets} as the union (as the described in~\ref{def:cat-unions})
  \[\cond\coloneqq \varinjlim_{\kappa}\cond_\kappa\]
  along all strong limit cardinals $\kappa$, where the transition maps are the fully faithful embeddings described above.
\end{definition}
With this definition at hand,
we will now present many elementary properties and characterisations of this category and of certain subcategories.
In working with condensed sets, it will become clear that they in many ways behave almost like the category of sets,
i.e., have many good topos-theoretic properties.
On the other hand, topological ideas lead to sensible notions in the context of condensed sets.

Many of the observed properties are well-known facts of topoi.
To give a feeling of how to work with condensed sets instead of topological spaces
we prove many results by hand in our case.
For a general overview about topoi and their properties see section~\ref{sec:sheaves-topoi} and especially subsection~\ref{ssec:el-topoi},
as well as the literature there.

\subsection{Equivalent definitions}

In this chapter we establish equivalent definitions of the category of condensed sets.

\begin{theorem}\label{thm:eq_def_cond}\uses{def:cond, lem:completions_via_Yoneda,prop:cond_k_to_cond_k_prime }
The following categories are well defined and equivalent.
Here, $\mcC$ can be taken as $\CHaus$, $\prof$, $\extr$ or $\{\beta I : I\in\Set\}$
with the coverage from definition~\ref{def:finitary_grothendieck_topology}.
\begin{enumerate}[(a)]

\item The colimit $\varinjlim_{\kappa}\cond_\kappa$ of all $\cond_\kappa$ along the ordered class of strong limit cardinal numbers.

\item The full subcategory of $\Sh(\mcC)$ consisting of all objects $T\in\Sh(\mcC)$ such that there exists a strong limit cardinal number $\kappa$ such that
\[
    T\simeq \Sh(\Lan T\lvert_{\mcC_\kappa}),
\]
or, equivalently, $\Sh(\Lan X)$ for any $X\in\Sh(\mcC_\kappa)$.

\item The full subcategory of $\PSh(\extr)$ consisting of all objects $T\in \PSh(\extr)$
mapping finite coproducts to products
such that there exists a strong limit cardinal number $\kappa$ with
\[
    T\simeq \Lan T\lvert_{\extr_\kappa}.
\]

\item The full subcategory of accessible sheaves on $\mcC$.

\item The sifted cocompletion of $\extr$.

\item The cocompletion of $\prof$ which preserves finite coproducts and
coequalizers of kernel pairs.
\end{enumerate}
\end{theorem}

\begin{remark}
Another characterisation we will prove later is that $\cond$ is given by those sheaves $T$ on $\extr$ (or $\profin$) such that there exists a set of $K_i\in \extr$ with an epimorphism
\[\coprod y(K_i)\twoheadrightarrow T\]
in $\Sh(\extr)$.
\end{remark}

To prove this theorem, we first need to establish some results
on the categories of sheaves on $\prof$ and $\extr$.
We need the following basic two facts about $\Sh(\prof)$,
which are not a priori clear, as $\prof$ is not small.

\begin{lemma}[Sheafification on $\profin$]\label{lem:sheafification}\uses{def:profin, def:sheaf, lem:adjoints_commute_with_limits, lem:unit_via_fully_faithfullness, lem:limits_along_final_subdiagrams, def:special_limits_and_colimits}
Let $\mcC$ be one of the categories $\CHaus$, $\prof$, $\extr$ or $\{\beta I : I\in\Set\}$.
Then the inclusion \(\Sh(\mcC)\subset\PSh(\mcC)\) has a left adjoint,
the \emph{sheafification}, which preserves finite limits.

In particular, it preserves colimits, maps epimorphisms to epimorphisms and monomorphisms to monomorphisms and its counit is an isomorphism,
meaning that the sheafification of a sheaf is the sheaf itself.
\end{lemma}

\begin{proof}
Sheafification exists as soon as the site $\mcC$ is small,
since then the colimit formula runs over a small diagram,
see theorem~\ref{thm:sheafification}.
However, $\mcC$ is not small so the coverings of an object in $\mcC$ form a proper class.
But lemma~\ref{lem:initial_coverings} guarantees that it contains an initial subset.
This together with lemma~\ref{lem:limits_along_final_subdiagrams}
implies that the colimit in the formula for sheafification exists.

Since sheafification, as soon as it exists, is a left adjoint (of the inclusion),
it preserves colimits,
and as the formula for calculating the sheafification is over a filtered index set,
sheafification commutes with finite limits.
Because epimorphicity can be phrased as a colimit property by remark~\ref{rem:epi_as_colimit}
and sheafification preserves colimits,
it maps epics to epics.
The same holds for monomorphism,
as these can be expressed as a finite limit.
The last assertion follows because the inclusion
$\Sh(\prof)\subset\PSh(\prof)$
is fully faithful, hence the counit is an isomorphism by lemma~\ref{lem:unit_via_fully_faithfullness}.
\end{proof}

\begin{lemma}[Yoneda for condensed sets]\label{lem:yoneda_in_cond}\uses{def:cond, lem:yoneda, prop:cond_k_to_cond_k_prime}
The category $\mcD$ defined in (b) of the theorem is locally small
and every representable presheaf $\hom_{\prof}(-,S)$ for $S\in\prolak$ is a sheaf that
is the left Kan extension of its restriction to $\prolak$.
In particular, the Yoneda lemma holds.
\end{lemma}

\begin{proof}
Note that the embedding $\cond_\kappa\hookrightarrow \mcD$ given by sheafification of left Kan extension is fully faithful,
as every $\cond_\kappa\hookrightarrow \cond_{\kappa'}$ is fully faithful and natural transformations may be checked on a bounded level.
Alternatively, one could also argue via the unit of the adjunction as in the proof of proposition~\ref{prop:cond_k_to_cond_k_prime}.

Next we show that for any \(S\in\prolak\),
the sheaf \(\hom(-,S)\) actually arises from left Kan extension.
For any $\tilde{S}\in\prof$ we have
\[
    \hom(\tilde{S},S) = \varinjlim_{\tilde{S}\to S'}\hom(S',S)
\]
where $S'$ runs over all $\kappa$-small $S'$.
This is exactly the formula for left Kan extension of the $\hom$ functor $\hom(-,S)$
on $\prolak$.
\end{proof}

\begin{proof}[Proof of the theorem \ref{thm:eq_def_cond}]
The equivalence of (a), (e) and (f):
theorem~\ref{thm:constructing_cocompletions} shows that the sifted cocompletion of $\extr_\kappa$ is exactly $\cond_k$
by lemmas~\ref{prop:k_condensed_extremally_disconnected} and \ref{lem:sind_via_cont},
as is the free cocompletion of $\prolak$ that preserves finite coproducts and coequalizers of kernel pairs, see lemma~\ref{lem:cond_k_on_subbasis}.
By taking the colimit along $\kappa$ of all $\cond_\kappa$ this yields the sifted cocompletion of $\extr$ by lemma~\ref{lem:cocomp_union},
and on the other hand the sifted cocompletion of $\extr$ is constructed this way.
The same argument holds for $\prolak$.
Moreover, lemma~\ref{lem:cocomp_union} also yields the equivalence of (c) to (a).

Equivalence of (a) and (b).
We construct an equivalence of categories as follows:
for any $T_\kappa\in\varinjlim_\kappa\cond_\kappa$ define $T = \Sh(\Lan T_\kappa)$.
For $T$ as in (b), let
\[
    \lambda
= \min\{\kappa : T\simeq \Sh\Lan T\lvert_{\prolak},\,\kappa\text{ strong limit cardinal}\}
\]
and map $T$ to this $T\lvert_{\prolam}$.
Then these functors are essentially inverse to each other:
the one direction is clear, for the other note that similarly to the argument in the proof of proposition~\ref{prop:cond_k_to_cond_k_prime},
the restriction of the sheafification of the left Kan extension of a sheaf is the sheaf itself.
The same argument adapts to $\mcC$ being $\CHaus$, $\extr$ or $\{\beta I : I\in\Set\}$.

Equivalence of (b) and (c).
Let $\mcC = \extr$.
Since the sheaf condition on $\extr$ reduces to $T$ mapping finite coproducts to finite products by lemma~\ref{prop:k_condensed_extremally_disconnected},
we only need to show is that one does not need to sheafify after left Kan extension.
This argument was given in lemma~\ref{lem:cond_k_to_cond_k_prime_alternative}.

Equivalence of (b) and (d).
Every sheaf in the category of (b) is accessible,
because $T\lvert_{\mcC_\kappa}$ is the colimit of representable sheaves (the coverage is subcanonical)
and left Kan extension as well as sheafification are left adjoints,
hence preserve the colimit.
So $T$ is accessible.

If $T$ is accessible, that is, it is the colimit of representable sheaves $S_i\in\prof$,
take the supremum $\kappa$ of all weights of $S_i$, i.e., $S_i\in\prolak$ for every $i$.
By lemma~\ref{lem:yoneda_in_cond} every representable presheaf $S_i$ on $\Sh(\prof)$
is the left Kan extension of its restriction to $\prolak$.
Interchanging the colimit with the left Kan extension yields the result.

\end{proof}

\subsection{(Co)limits of condensed sets}

\begin{proposition}\label{lem:condset_bicomplete}\uses{def:cond, prop:cond_k_to_cond_k_prime}
The category of condensed sets is bicomplete.
\end{proposition}

\begin{proof}
Consider a $\lambda$-small diagram.
For the computation of the (co)limit choose a strong limit cardinal $\kappa$
such that all objects in the diagram are in $\condk$.
Then the colimit, computed in $\cond_\kappa$, is preserved by the left Kan extension
and so is the colimit in $\cond$.
Concerning the limit, further assume that $\kappa$ has cofinality at least $\lambda$.
Then one can compute the limit in $\condk$
and since afterwards all the Kan extensions preserve this limit, corollary~\ref{cor:Lan-lambda-small-limits},
we know that the Kan extension of this $\kappa$-condensed set will be the desired limit.
\end{proof}

\begin{lemma}\label{lem:colimit_computation_cond}\uses{def:cond, thm:eq_def_cond, prop:cond_k_to_cond_k_prime}
In $\Sh(\profin)$, and hence in $\cond$, all limits can be computed pointwise,
as well as all filtered colimits.

In $\Sh(\extr)$, and hence in $\cond$ on extremally disconnected spaces,
all sifted colimits can be computed pointwise.
\end{lemma}

\begin{proof}
The pointwise limit of sheaves is again a sheaf, see corollary~\ref{lem:limits_pointwise}.
In view of lemma~\ref{lem:lim_in_full_sub} it remains to show that this sheaf is actually an object of $\cond$.
Consider a $\lambda$-small diagram $T_i$ in $\cond$.
For any $S\in\prof$,
take a strong limit cardinal $\kappa$ with cofinality at least $\lambda$
such that all $T_i$ are in $\condk$ and $S\in\prolak$.
In $\condk$, the limit of the $T_i$ is computed pointwise,
i.e., $(\varprojlim T_i)(S) = \varprojlim T_i(S)$.
Now the transition $\cond_\kappa\to\cond$ preserves $\lambda$-small limits,
so the limit in $\cond$ is also computed pointwise on $S$.

Since filtered colimits commute with finite limits in \(\Set\),
filtered colimits of sheaves on $\prof$ as presheaves remain sheaves
(as all coverings of the coverage on $\prof$ are finite).
As they furthermore commute with pointwise left Kan extension and sheafification,
filtered colimits may be calculated pointwise.

For the last statement note that on $\Sh(\extr)$ the sheaf condition reduces to commutation with finite products,
and hence this commutes with sifted colimits.
Note that no sheafification of the Kan extension is needed here.
\end{proof}

\begin{remark}
In particular, one can compute the (co)limits in the sheaf category
$\condk = \Sh(\prolak)$ for large enough \(\kappa\).
Furthermore, all (co)limits in \(\cond\)
coincide with the corresponding (co)limit in \(\Sh(\prof)\),
so that they also can be computed there.
\end{remark}

The category $\cond$ has a terminal and initial object.
These are the sheaves
\[
    S\mapsto \ast
\]
and
\[
    S\mapsto\begin{cases}\emptyset,& S\ne \emptyset,\\\ast,&S=\emptyset,\end{cases}
\]
respectively.
(The latter is the sheafification of $S\mapsto\emptyset$.)
Note that any condensed set $T$ other than the initial object has a point,
i.e., nonempty $T(\ast)$.
This is due to the fact that if $T(S)\ne \emptyset$ for some $S\ne \emptyset$
then there exists a map $\ast\to S$,
and thus there must exist a map $\emptyset \ne T(S)\to T(\ast)$ in $\Set$,
which implies $T(\ast)$ to be nonempty.

As we have seen, computing limits and filtered, even sifted, colimits of condensed sets
is particularly easy.
On the point $*$ all colimits can be computed pointwise,
meaning that in any case the underlying set can be immediately understood.

\begin{lemma}[Sheafification on terminal object]\label{lem:Sh_on_terminal}\uses{lem:sheafification}
For any presheaf $F$ on $\profin$ we have
\[\Sh(F)(\ast)=F(\ast),\]
i.e., sheafification does not change the value on the point.
We conclude that the underlying set of any colimit and limit
is always the colimit resp. limit of the underlying sets.
\end{lemma}

\begin{proof}
	Any covering sieve of $\ast$ contains a nontrivial morphism, i.e. a morphism $X\to \ast$ for $X\ne \emptyset$.
	But now there exists an arrow $\ast\to X$, which by the sieve condition implies $\ast\to X\to \ast$ is in the covering.
	But as $\ast$ is terminal, we conclude $1_\ast$ and thus all morphisms to be in the covering.
	Thus there exists only one covering sieve of $\ast$, and hence $\Sh F(\ast)=F(\ast)$
    by the explicit formula for sheafification (see~\ref{thm:sheafification}).
\end{proof}

For general colimits, one has to sheafify the pointwise colimit,
as the following example shows.

\begin{example}\label{ex:colim-not-pointwise}

    Let \(\alpha\in [0,1]\) be an irrational number and denote by \(+\alpha\colon\T\to\T\) the shift/rotation by \(\alpha\) on \(\T = \R/\Z\).
    Consider the coequalizer diagram
    \begin{center}
        \begin{tikzcd}
	       \T && \T && C
	       \arrow["+\alpha", shift left, from=1-1, to=1-3]
	       \arrow["1"', shift right, from=1-1, to=1-3]
	       \arrow["q", two heads, from=1-3, to=1-5]
        \end{tikzcd}
    \end{center}
    taken in \(\cond\).
    By the previous lemma,
    the value of $C$ at the point is given by the corresponding coequalizer diagram
    in $\Set$ and hence given by the projection
    \[
        q_*\colon\T\to\T/\alpha\T,
    \]
    where $\T/\alpha\T$ is the quotient of \(\T\) by the equivalence relation \(s\sim t\)
    if and only if there exists a \(k\in\Z\) such that \(s = k\alpha + t\).
    This is, of course, the smallest equivalence relation such that
    the quotient map \(q_*\colon\T\to\T/\alpha\T\) equalizes \(\alpha\) and \(1\).
    By the general description of coequalizers in \(\set\) it then follows that this quotient is indeed the coequalizer.

    By the sheaf condition, $C(\{0,1\}) = C(\ast)\times C(\ast) = (\T/\alpha\T)^2$.
    But the diagram
\begin{center}
    \begin{tikzcd}
	{\T^2} && {\T^2} && {(\T/\alpha\T)^2}
	\arrow["{(+\alpha,+\alpha)}", shift left, from=1-1, to=1-3]
	\arrow["1"', shift right, from=1-1, to=1-3]
	\arrow["{(q_*,q_*) = q_{\{0,1\}}}", from=1-3, to=1-5]
    \end{tikzcd}
\end{center}
    is not a coequalizer diagram in $\Set$.
    To see this, assume it were and
    let \(S\) be the quotient of \(\T^2\) by the equivalence relation \((s,t)\sim (u,v)\)
    if and only if \((s,t) = (k\alpha + u,k\alpha + v)\) for a \(k\in\Z\).
    Denote by \(r\colon\T^2\to S\) the projection.
    Then \(r\circ(\cdot\alpha,\cdot\alpha) = r\) and so there exists a map
    \(f\colon(\T/\alpha\T)^2\to S\) with \(f\circ(q_*,q_*) = r\).
    But then
    \[
        r(\alpha,0) = f(q_*(\alpha),q_*(0)) = f(q_*(0),q_*(0)) = r(0,0).
    \]
    This is a contradiction.
    So $(\T/\alpha\T)^2$ is not the coequalizer (actually, $S$ is).
    In particular, we cannot compute the coequalizer pointwise in $\Set$.

    In fact,
    there is a much more trivial example given by the binary coproduct of a point with itself,
    $\ast\sqcup\ast$.
    If we didn't need to sheafify,
    we would get
    \[
      (\ast\sqcup\ast)(S)=\ast(S)\sqcup\ast(S)=\{0,1\}.
    \]
    But the sheaf condition implies (in the case $S=\ast\sqcup\ast$)
    \[
      (\ast\sqcup\ast)(\ast\sqcup\ast)=(\ast\sqcup\ast)(\ast)\times(\ast\sqcup\ast)(\ast)
      =\{0,1\}^{2}\neq\{0,1\},
    \]
    so that $\ast\sqcup\ast$ can not be computed pointwise.
\end{example}

As we have seen that one needs to sheafify in general,
we provide the following nice formula for sheafification on $\extr$,
which will be proven later in more generality, see \ref{lem:sheafification_cond_C}.

\begin{lemma}\label{lem:sheafification_cond}
    Consider a presheaf $F\in \PSh(\extr)$.
    Then the sheafification $\Sh(F)$ evaluated on any $S\in \extr$ is given by
    \[\Sh(F)(S)=\varinjlim_{S=\bigsqcup_{i=1}^k T_i}\prod_{i=1}^k F(T_i),\]
    where the colimit runs over all finite clopen partitions of $S$.
\end{lemma}

Next, we exhibit some nice properties of (co)limits in the category of condensed sets
which resemble the category of sets.

\begin{lemma}\label{lem:coprod_disj_stable}\uses{def:cond, lem:colimit_computation_cond}
\begin{enumerate}[(i)]
\item Colimits in $\cond$ are stable under base change.
        In particular, colimits distribute over pullbacks,
        e.g., $(\coprod A_i)\times B = \coprod (A_i\times B)$.
\item Sifted colimits commute with finite products and filtered colimits with finite limits.
\item Coproducts in $\cond$ are disjoint.
\item If \(A_i\hookrightarrow B\) are monic and pairwise disjoint,
then \(\coprod A_i\hookrightarrow B\) is monic as well.
\item If \(A_i\hookrightarrow B_i\) are monic, then \(\coprod A_i\hookrightarrow\coprod B_i\) is monic as well.
In particular, the embeddings into the coproduct of each of its components are monic.
\item If \(B_i\twoheadrightarrow A_i\) are epic, then \(\prod B_i\to\prod A_i\) is epic as well.
In particular, the projections from the product onto its coordinates are epic.
\item Every equivalence relation is effective.
\item If $f\colon A\to B$ is a monomorphism and $g\colon A\to C$, the pushout diagram
\begin{center}
\begin{tikzcd}
	A & B \\
	C & {B\sqcup_A C}
	\arrow["f", hook, from=1-1, to=1-2]
	\arrow["g"', from=1-1, to=2-1]
	\arrow["{g_*}", from=1-2, to=2-2]
	\arrow["{f_*}"', hook, from=2-1, to=2-2]
\end{tikzcd}
\end{center}
        is also a pullback diagram and $f_*$ is a monomorphism;
        and dually for epimorphisms.

\end{enumerate}
\end{lemma}

\begin{proof}
There are at least two ways of obtaining these results.

One can use the fact that $\condk$, as a Grothendieck topos,
is an elementary topos and then use the properties of those, see theorem~\ref{thm:elem_props_topos}.
Since all claims are only made over \emph{small} sets of objects,
choosing $\kappa$ large enough yields the respective assertions in $\cond$
because the left Kan extension preserves the limits and colimits used.

On the other hand, one could argue in the same way as one proves that a Grothendieck topos satisfies Girauds axioms.
All of the properties (i)-(viii) are true in the category of sets,
so they remain true in the category of presheaves $\PSh(\extr_\kappa)$,
since there everything is computed pointwise.
Now sheafification is cocontinuous and preserves finite limits
(and arbitrary limits are computed pointwise anyway),
so the claims still hold in the sheaf category $\condk$.
Then as before one lifts the (co)limits to $\cond$ via Kan extension.
\end{proof}

Note that in $\cond$, general colimits need not preserve monomorphism
and limits need not preserve epimorphisms.
For example, consider the example right after lemma~\ref{lem:colim-preserve-epis}
which can be also interpreted in $\extr$ and then embedded into $\cond$ via Yoneda.
Evaluation at the point again yields the example mentioned
and since on the point, all (co)limits are computed in $\Set$,
we get that the morphism $\emptyset\to \{0,1\}$ is not epimorphic in $\cond$.

The following internal $\hom$ structure will be of huge importance when considering condensed abelian groups.
\begin{lemma}
$\cond$ is a cartesian closed category with internal $\hom$ being given on $S\in\extr$ by
\[\ihom(X,Y)(S)=\hom(X\times S,Y).\]
\end{lemma}
\begin{proof}
  It is well-known that Grothendieck topoi are cartesian closed,
  which yields a right adjoint to $X_{|\kappa}\times-$ from $\condk$ to $\condk$.
  As $X_{|\kappa}\times-$ glue to $X\times-$,
  we get a right adjoint $\ihom(X,-)$ by~\ref{lem:tower-adjs}.
  The stated formula then follows via Yoneda:
  For $Y\in\extr$,
  \[
    \ihom(X,Y)(S)=\hom(S,\ihom(X,Y))=\hom(S\times X,Y).
  \]

  To illustrate how the internal $\hom$ works,
  we also prove the sheaf condition and adjunction by hand.
$\ihom(X,Y)$ is a contravariant functor $\extr\to \Set$, and since coproducts distribute over products, it maps finite coproducts in $\extr$ to finite products,
\[\ihom(X,Y)(S\sqcup T)=\hom(X\times(S\sqcup T),Y)=\hom(X\times S,Y)\times \hom(X\times T,Y).\]
Thus $\ihom(X,Y)$ is a condensed set.
Furthermore,
\[\hom(S\times Y,Z)=\hom(S,\ihom(Y,Z)).\]
Since colimits are stable under base change,
the functor $-\times Y$ is cocontinuous
and thus both sides are cocontinuous, showing that
\[\hom(X\times Y,Z)=\hom(X,\ihom(Y,Z))\]
for all colimits $X$ of objects $S\in \extr$ -- but every condensed set is a colimit of such $S$!
\end{proof}

Note that for compactly generated Hausdorff spaces, the internal $\hom$'s agrees with the classical ones (which carry the compact-open topology).
This essentially follows from the fact that the compact-open topology yields an internal $\hom$ in compactly generated Hausdorff spaces,
and that from there the embedding is fully faithful (and that the right adjoints are unique).
See \cite{Aparicio2021} for more.

\subsection{Morphisms of condensed sets}

Next, we will prove many helpful properties about the morphisms of $\cond$.
For this, we make use of the fact that each \(\cond_\kappa\) is a Grothendieck topos.
These inherit many nice properties of the category of sets.

\begin{proposition}\label{prop:cond-morphisms}
\begin{enumerate}[(i)]
\item All monomorphisms are regular and effective.
\item The category $\cond$ is balanced.
\item The epimorphisms agree with the regular/effective epimorphisms.
    In particular, the pullback of an epimorphism is again an epimorphism.
\item Every morphism $f\colon A\to B$ can be factorized as
\[
        A\twoheadrightarrow M\hookrightarrow B,
\]
    where $M\hookrightarrow B$ is the image and $A\twoheadrightarrow M$ the coimage of $f$.
    If $A\twoheadrightarrow N\hookrightarrow B$ is another such factorisation of $f$
    then there exists a unique isomorphism $M\to N$ respecting the factorization.
    The coimage can be computed as the coequalizer of the kernel pair, and the image is the equalizer of the cokernel pair. 
    Furthermore, the image $M$ is given by the pointwise image
    $S\mapsto f_S(A(S))$ for $S\in\extr$.
\end{enumerate}
\end{proposition}

\begin{proof}
\begin{enumerate}[(i)]
\item Let \(f\colon X\hookrightarrow Y\) be a monomorphism in \(\cond\).
Then choose $\kappa$ so that $X$ and $Y$ arise from their restrictions to $\condk$.
The restriction of $f$ to $\condk$ results in a monomorphism in \(\condk\)
since the embedding $\condk\to\cond$ is faithful.
By remark~\ref{rem:epi_as_colimit}, the kernel pair of $f$ consists of twice the identity on $X$.
In particular, this is also the kernel pair of $f$ in the presheaf category $\PSh(\extr_\kappa)$.
Because there, limits are computed pointwise, $f$ is sectionwise monic,
i.e., every section $f_S\colon X(S)\to Y(S)$ is monic.
But in $\Set$, and hence pointwise in $\PSh(\extr_\kappa)$,
every monomorphism is effective.
As sheafification preserves finite limits,
$f$ is the equalizer of its cokernel pair in $\condk$.
Since left Kan extension preserves finite limits,
the monomorphism \(f\colon X\hookrightarrow Y\)
is the equalizer of two arrows in \(\cond\).

\item Consider a bimorphism $f\colon X\to Y$.
By the previous observation, $f$ is the equalizer of two parallel arrows $g,h$.
But since $f$ is also epimorphic, $g=h$.
This implies that $f$ is an isomorphism.

\item Now to see that in $\cond$ every epimorphism is effective,
consider an epimorphism $f\colon X\twoheadrightarrow Y$.
Let $f(X)$ be the coequalizer of the kernel pair of $f$ in $\cond$.
As it is a reflexive coequalizer, $f(X)$ can be computed pointwise on $\extr$.
Since in $\Set$, every epimorphism is effective
(by example~\ref{example:set-epi-reflective}),
taking the set-theoretic image yields $f(X)(S) = f_S(X(S))$.
The factorisation
\begin{center}
\begin{tikzcd}
	{X\times_YX} & X & {f(X)} \\
	&& Y
	\arrow[shift left, from=1-1, to=1-2]
	\arrow[shift right, from=1-1, to=1-2]
	\arrow["{f^\mid}", two heads, from=1-2, to=1-3]
	\arrow["f"', two heads, from=1-2, to=2-3]
	\arrow["i", hook, from=1-3, to=2-3]
\end{tikzcd}
\end{center}
is given by the pointwise inclusion $f(X)(S) = f_S(X(S)) \subset Y(S)$.
But now $i$ is not only monic, but also epic (after all, $f$ is),
so (ii) implies that $i$ is an isomorphism.
This means that $f$ is the coequalizer of its kernel pair.
The second assertion follows from lemma~\ref{lem:coprod_disj_stable} (i).

\item This follows from the proof of part (xiv) in theorem~\ref{thm:elem_props_topos},
noting that the proof there only uses that every monic and every epic is effective
and that the category is balanced, which we have just established.
The last claim follows from (iii).
\end{enumerate}
\end{proof}

\begin{remark}
As we can see, morphisms behave almost as good as in usual set theory.
One important property of $\Set$ that $\cond$ does not fulfill
is the internal axiom of choice:
In $\Set$,
every epimorphism has a section.
Of course,
it would be catastrophic if this held in $\cond$!
We want to view $\mcT,[0,1]\in\cond$ as \enquote{topological-like} objects and as such,
$[0,1]\to\mcT$ should not have a section.
\end{remark}

Lastly, we give an even easier description of monics and epics.

\begin{proposition}\label{prop:char_epi}\uses{lem:comp-ed-proj}\leanok
  The monics in $\cond$ are precisely the sectionwise\footnote{
    Technically,
    \enquote{sectionwise} and \enquote{pointwise} are of course synonyms
    but for some reason,
    the first is used mainly for morphisms and the latter for objects.
  }
  monomorphisms,
i.e., $f\colon X\to Y$ is monic precisely if for every $S\in \prof$
(or equivalently every $S\in \extr$) the map $f_S\colon X(S)\to Y(S)$ is monic.

A morphism $f\colon X\to Y$ is an epimorphism precisely if it satisfies one of the following equivalent properties.
\begin{enumerate}[(a)]
\item Its sections $f_S$ are epimorphic for all $S\in \extr$.
\item The morphism $f$ is \emph{locally epimorphic}:
for all $S\in \prof$ and for all $h\in Y(S)$
there exists a cover of $S$ such that for all arrows $k\colon K\to S$ in the cover,
there is $g\in X(K)$ with $Y(k)(h)=f_{K}(g)\in Y(K)$.
\item For all $h\colon S\to Y$ with $S$ being compact Hausdorff there exist finitely many $\CHaus$ objects $K_i$ with arrows $k_i\colon K_i\to S$ and $g_i\colon K_i\to X$ such that
\begin{center}
\begin{tikzcd}
	X && Y \\
	\\
	{\coprod K_i} && S
	\arrow["f", from=1-1, to=1-3]
	\arrow["{\amalg g_i}", from=3-1, to=1-1]
	\arrow["{\amalg k_i}"', two heads, from=3-1, to=3-3]
	\arrow["h"', from=3-3, to=1-3]
\end{tikzcd}.
\end{center}
\end{enumerate}
\end{proposition}

\begin{proof}
Clearly, if a natural transformation is sectionwise epic or monic,
it is an epimorphism or a monomorphism, respectively.
A morphism \(f\colon X\to Y\) is monic if and only if its kernel pair
is the identity on \(X\).
Since limits can be calculated pointwise in \(\cond\)
we obtain that pointwise the kernel pair is the identity on \(X(S)\) as well.
This shows that monics are sectionwise monic.

Now let us prove the characterization of epimorphisms.
Let $f$ be an epimorphism.
Because epimorphisms are effective
and reflective coequalizers can be calculated pointwise on $S\in\extr$,
all sections $f_S$ are epimorphic.
On the other hand, if $f$ is sectionwise epic, then $f$ is an epimorphism as well.
The implication from (a) to (b) follows by choosing the cover $\beta S_d\to S$.
The equivalence of (b) and (c) is due to the definition of the finitary Grothendieck topology, using that the Yoneda embedding preserves finite coproducts and epimorphisms
by lemma~\ref{lem:cond_yon_komm}.
Alternatively, use the equivalent description (f) of $\cond$ in theorem~\ref{thm:eq_def_cond}.
Finally, (c) implies (a), as for any extremally disconnected $S$
there is a section $s$ to the epimorphism $k\colon \coprod K_i\twoheadrightarrow S$ and hence for all $h\colon S\to Y$ there exists $g\colon \coprod K_i\to X$
inducing $gs\colon S\to X$ with
\[
    fgs=hks=h
\]
and thus, applying the Yoneda lemma, $f_S$ is surjective.
\end{proof}

\begin{remark}
The characterisation (b) of epimorphisms is a general fact of Grothendieck topoi
and is usually derived using the subobject classifier, see e.g. corollary III.7.5 in \cite{Lane1992SheavesIG}.
\end{remark}

\begin{example}
    The map \(\R_d\to\R\) of condensed sets is not surjective
    because the maps from \(\beta\N\) to \(\R_d\) have finite image
    but there exists a surjection \(\beta\N\to [0,1]\),
    so the map
    \[
        C(\beta\N,\R_d)\to C(\beta\N,\R)
    \]
    is not surjective.
\end{example}

\begin{remark}
Note that although $\cond_\kappa$ admits a subobject classifier, $\cond$ does not, as it could not be left Kan extension from some $\kappa$.
The subobject classifier $\Omega$ would have to satisfy
\[
    \{ A\subset S : A \text{ closed} \} \subseteq \Omega(S)
\]
which cannot be a condensed set, see example~\ref{example:sierpinski}.
\end{remark}

Since we do not have a subobject classifier, this suggests that the category might not be well-powered.
However, it is well-copowered, and later we will see that $\cond$ is indeed also well-powered.

\begin{lemma}
$\cond$ is well-copowered,
and every factor object of a $\kappa$-condensed set is $\kappa$-condensed.

I.e., every object has a small set of (equivalence classes of) factor objects.
\end{lemma}
\begin{proof}
To see this, let $X$ be a condensed set that arises from $\condk$
and let $X\to Y$ be an epimorphism. Then $Y$ is also $\kappa$-condensed (this follows from the characterisation \ref{lem:cond_eq_colim_CHaus}).
Hence all factor objects of $X$ are $\kappa$-condensed objects.
But since $\condk$ is a Grothendieck topos, it is well-copowered.
Alternatively: Because every epimorphism in $\condk$ is sectionwise on $\extr_\kappa$ an epimorphism,
and there is just a set of sections as well as a set of epimorphisms for a fixed section, there is just a set of such epimorphisms.
\end{proof}

In the above proof, we used the following characterisation of condensed sets, which we can finally prove using \ref{prop:char_epi}.

\begin{proposition}\label{lem:cond_eq_colim_CHaus}
  The category of condensed sets is equivalent to the full subcategory of sheaves $T$ on $\prof$,  $\CHaus$, $\extr$ or $\beta\mathbf{Disc}$ such that there exist a set I,
  $K_i\in \prof$ ($i\in I$) and an epimorphism
\[
    \coprod K_i\twoheadrightarrow T.
\]
Furthermore, in this case, $T$ is $\kappa$-condensed for every (uncountable strong limit cardinal) $\kappa>|K_i|$ (for all $i\in I$).
\end{proposition}

\begin{proof}
  By accessibility,
  every condensed set is colimit of (a set of) profinite sets.
Representing this colimit as th coequalizer of two coproducts,
which is always epic, we see that the condensed set is a factor of a coproduct
(in condensed sets!) of profinite sets. On the other hand consider any sheaf $T$ with $\coprod K_i\twoheadrightarrow T$.

As the epimorphism is surjective on extremally disconnected sets,
we know that any map $S\to T$ is induced by a map $S\to \coprod K_i$, and by replacing $\coprod_{i\in I} K_i$ with the filtered limit over finite subdiagrams,
and noting that filtered colimits may be computed pointwise,
\[\left(\coprod_{i\in I} K_i\right)(S)=\varinjlim_{F\sub I\,\mathrm{finite}} \left(\left(\coprod_{i\in F} K_i\right)(S)\right),\]
 we conclude that the map $(S\to \coprod K_i) \in (\coprod K_i)(S)$ factors through a finite sub-coproduct $S\to \coprod_F K_i=K$.
Note that this essentially already contains the proof that $S$ is compact, see \ref{lem:comp-ed-proj}.

Hence the continuous maps $S\to T$ coincide with the colimit $\varinjlim_{S\to K}T(K)$ for $K$ having cardinality at most $\sup |K_i|$. Thus the Kan extension formula is true for $T$, and $T$ is a condensed set. 
\end{proof}

\subsection{Topological spaces as condensed sets}

We have seen that every topological space $X$ induces a $\kappa$-condensed set
via the assignment
\[
    \underline{X}\colon S\mapsto C(S,X).
\]
But when passing to the colimit of all $\kappa$-condensed sets,
i.e., the category of accessible sheaves on $\prof$,
surprisingly,
some topological spaces \emph{that should not exist anyway} indeed do not give condensed sets
under the canonical functor $X\mapsto\underline{X}$.

\begin{example}\label{example:sierpinski}
The Sierpi\`nski space $F=\{a,b\}$ does not yield a condensed set via
\[
    S\mapsto C(S,F)=\mathrm{closed}(S).
\]
This is due to the fact that otherwise there would have to exist a $\kappa$ such that for every $S\in\prof$ the equality
\[
    \mathrm{closed}(S)=\varinjlim_{S\to T}\mathrm{closed}(T)
\]
holds, where $T$ ranges over \(\prof_\kappa\).
But this would mean that a set in $S$ is closed precisely
if it is arises from $\kappa$-profinite sets,
and hence is the intersection of $\kappa$ many clopen subsets.
This follows, as any closed subset $S\to F$ would have to factor over some $S'\in \prof_\kappa$,
and thus correspond to the preimage of a closed subspace in some $\kappa$-profinite set.
But there, the closed subsets are given by intersection of at most $\kappa$ many clopens,
which translates under preimages to the closed subset of $S$.
So any closed subset of any profinite set would be the intersection of \(\kappa\)
many clopen subsets -- for a fixed $\kappa$!
The closed subset \(\{1\}^{2^\kappa} \subset \{0,1\}^{2^\kappa}\) is a counterexample.
\end{example}

\begin{remark}
  Taking the Sierpinksi space in $\condk$, and afterwards Kan extending clearly yields condensed sets.
  Their value at $S$ can be described as all closed sets in $S$ that are intersection of $\kappa$ many clopens.
\end{remark}

However, reasonable topological spaces are still condensed sets.

\begin{proposition}[T1-spaces are condensed sets]\label{prop:T1_is_cond}\uses{def:cond,lem:top_yields_cond_k}\leanok
Let $X$ be a topological space.
Then $\underline{X}\colon S\mapsto C(S,X)$ is a condensed set if and only if $X$ is a T1-space.
\end{proposition}

\begin{proof}
This is  2.15 in \cite{scholze2019condensed}.
If $X$ is not T1, it contains a subset $\{a,b\}$ which is not discrete, and thus does not induce a condensed set by the same argument as for the Sierpinski space above
($\underline{X}(S)$ would have to have a copy of every closed set by defining the two-valued characteristic functions with values in the nondiscrete subspace).

If $X$ is T1, we have to show that there is some $\kappa$ such that for all $\tilde S$ with cardinality bigger than $\kappa$ the formula
\[
    C(\Tilde S,X)=\varinjlim_{\tilde S\to S}C(S,X)
\]
(with $S$ ranging over $\extr_{\kappa}$) holds, i.e.,
that any $f\colon \tilde S\to X$ factors over some $S\to X$ with $S$ being $\kappa$-small.
For this, it suffices to find a surjection $\tilde S\to S$ such that both maps $\tilde S\times_S \tilde S\to \tilde S\to X$ agree, since then
\begin{center}\begin{tikzcd}
	{\tilde S\times_S \tilde S} & {\tilde S} \\
	{\tilde S} & S \\
	&& X
	\arrow[from=1-1, to=1-2]
	\arrow[from=1-1, to=2-1]
	\arrow[from=1-2, to=2-2]
	\arrow[from=1-2, to=3-3]
	\arrow[from=2-1, to=2-2]
	\arrow[from=2-1, to=3-3]
	\arrow[dashed, from=2-2, to=3-3]
\end{tikzcd}\end{center}
which follows from \(S=\tilde{S}/(\tilde{S}\times_S\tilde{S})\).
But this is equivalent to the statement that for any $x\ne y$ in $X$ the preimages
$f^{-1}(x)\times f^{-1}(y)\subseteq \tilde S\times \tilde S$ are disjoint from
$\tilde S\times_S\tilde S$.

We know that $f^{-1}(x)\times f^{-1}(y)$ is a closed and hence compact subset $T$ of $\tilde S\times \tilde{S}$.
For any point $(a,b)\in T$ we know that there is a map $h\colon \tilde{S}\to \{0,1\}$
such that $h(U_a)=0$, $h(U_b)=1$ for neighbourhoods $U_a$ and $U_b$,
and hence $(U_a\times U_b)\cap (\tilde{S} \times_h \tilde{S})=\emptyset$.
This way we can cover $T$, and using compactness we can glue finitely many $h_i$ together
to one $g\colon \tilde S\to \{0,1\}^n$ with $\tilde S \times_g \tilde S\cap T=\emptyset$.
Now we glue these $h$ together for all $(x,y)$
(these are at most $|X|$ many),
to obtain $h\colon \tilde{S}\to \{0,1\}^{|X|}$.
Passing to its image,
we get that \(S = \im h\) fulfills what we look for,
since the weight of \(S\) is bounded by the cardinality of the powerset of \(X\).
\end{proof}

\begin{example}
Consider any space $X$ equipped with the cofinite topology.
Then $X(S)=C(S,X)$ is the set of maps $S\to X$ such that the preimage of every point is closed, i.e.,
can be described as the set $\{(S_{x})_{x\in X}\mid S_{x}\text{closed subset of }S\text{ and }\dot\bigcup_{x\in X} S_x=S\}$
(by using complements of the $S_{y}$,
this implies that every $S_x$ in such a decomposition is the intersection of $|X|$ many opens).

Using the co-$\kappa$-small topology one can build analogous condensed sets.
\end{example}

\begin{remark}
If we are working with $\kappa$-condensed sets,
then this functor has a left adjoint,
namely the functor sending a condensed set $T$ to $T(*)_{\kappa\mathrm{-top}}$,
see theorem~\ref{thm:cond_k_to_top}.
Because we cannot even canonically embed all of $\Top$ into $\cond$,
this adjunction does not hold in full generality anymore.
However,
on T1-spaces,
the functor $X\mapsto\underline{X}$ is still part of an adjunction,
if we restrict to a suitable subcategory of $\cond$, which we will get to know later,
see lemma~\ref{rem:cond-T1-adjunction}.

Moreover, the functor $X\mapsto\underline{X}$ is continuous,
since limits can be computed pointwise and the $\hom$-functor in $\Top$
preserves limits.
\end{remark}

Now we will investigate how the functor $X\mapsto\underline{X}$ behaves
with respect to colimits.

The most important embedding of topological spaces clearly comes from
the Yoneda embedding of $\extr,\prof$ or even $\CHaus$.
We show that the Yoneda embedding yields a good translation between properties of morphisms.

\begin{lemma}[Yoneda in condensed commutes]\label{lem:cond_yon_komm}\uses{lem:yoneda_in_cond,lem:condset_bicomplete, prop:cond-morphisms, prop:char_epi}\leanok
The fully faithful Yoneda embedding
\[
    \prof\hookrightarrow\cond,\quad S\mapsto\underline{S}
\]
preserves arbitrary limits, finite coproducts and coequalizers of kernel pairs of epimorphisms \(S\twoheadrightarrow T\), i.e., quotient maps.
In particular, it preserves images.

Furthermore, it induces bijections between the classes of epimorphisms resp.\ monomorphisms.

If $S\in\prof$ and $X\in\cond$ is a retract of $S$, then $X = \underline{S'}$
for some $S'\in\prof$.
\end{lemma}

\begin{proof}
Clearly, monics and epics between $\underline{X}\to \underline{Y}$ induce monics/epics in $X\to Y$, since the inclusion is fully faithful.
For the converse, consider any arrow $\phi\colon X\to Y$.
If $\phi$ is monic, then $\phi_*(S)\colon \underline{X}(S)\to \underline{Y}(S)$  is injective
for all profinite \(S\) by definition of monomorphism.
Furthermore, for any epimorphism $\phi$ and extremally disconnected $S$ the morphism $\phi_*$ is epimorphic since $S$ is projective in $\prof$.
Hence the morphism $\phi_*$ is sectionwise epimorphic in $\Sh(\extr)$,
and thus clearly epimorphic.
Since the categories of sheaves on $\extr$ and on $\prof$ are equivalent (by restriction),
this means that $\phi_*$ is epimorphic in $\Sh(\prof)$ and thereby in $\cond$.

For the commutation with arbitrary limits note that the contravariant Yoneda embedding always commutes with limits as the Hom-functor does so,
and limits may be computed pointwise.

For the assertion about colimits,
consider first the case of a finite coproduct $\coprod K_i$ in \(\prof\),
which is just the disjoint union with the sum topology.
Then for any condensed set $T$, by the Yoneda lemma and the sheaf condition, it follows that
\[
    \Hom\left(\underline{\coprod K_i},T\right)=T\left(\coprod K_i\right)=\prod T(K_i)=\prod \Hom(\underline{K_i},T),
\]
which precisely is the condition of $\underline{\coprod K_i}$ being the coproduct in $\cond$.

Next, consider any surjection \(f\colon S\twoheadrightarrow T\) in \(\prof\).
Then \(f\) is the coequalizer of its kernel pair
\begin{center}
    \begin{tikzcd}
	{S\times_T S} && S && T
	\arrow["{p_1}", shift left, from=1-1, to=1-3]
	\arrow["{p_2}"', shift right, from=1-1, to=1-3]
	\arrow["f", from=1-3, to=1-5]
\end{tikzcd}
\end{center}
because in \(\prof\) every surjection is a quotient map.
Now we show that
\begin{center}
    \begin{tikzcd}
	{\underline{S\times_T S} \cong \underline{S}\times_{\underline{T}}\underline{S}} && {\underline{S}} && {\underline{T}}
	\arrow["{p_{1,*}}", shift left, from=1-1, to=1-3]
	\arrow["{p_{2,*}}"', shift right, from=1-1, to=1-3]
	\arrow["{f_*}", from=1-3, to=1-5]
\end{tikzcd}
\end{center}
is a coequalizer diagram in \(\cond\),
where the equality on the left hand side arises from the fact that the Yoneda embedding commutes with limits, as we just have proven.

For this consider any morphism \(g\colon\underline{S}\to X\) in \(\cond\)
such that \(g\circ p_{1,*} = g\circ p_{2,*}\).
In view of the third sheaf axiom we obtain that
\begin{center}
    \begin{tikzcd}
	{X(T)} && {X(S)} && {X(S\times_T S)}
	\arrow["{f^*}", from=1-1, to=1-3]
	\arrow["{p_1^*}", shift left, from=1-3, to=1-5]
	\arrow["{p_2^*}"', shift right, from=1-3, to=1-5]
\end{tikzcd}
\end{center}
is an equalizer diagram in \(\Set\).
But we can interpret the morphism \(g\) as an element of \(X(S)\)
and the condition \(g\circ p_{1,*} = g\circ p_{2,*}\) means that \(p_1^*(g) = p_2^*(g)\).
Thus there exists exactly one \(h\in X(T)\) such that \(f^*(h) = g\),
or equivalently \(h\circ f_* = g\).

The last assertion now follows from the fact that the Yoneda embedding preserves images:
if we have $s\colon X\to \underline{S}$ such that $r\colon \underline{S}\to X$ with $rs = 1$,
then $X$ is the image of $sr\colon\underline{S}\twoheadrightarrow X\hookrightarrow\underline{S}$.
Now $sr\colon S\to S$ in $\prof$ has image $S'\in\prof$,
and hence $X = \underline{S}$ because the Yoneda embedding preserves images.
\end{proof}

\begin{remark}\label{rem:cond-yon-chaus}
The Yoneda embedding
\[
    \CHaus\hookrightarrow\cond
\]
has the same properties as the one from $\prof$ in the previous lemma.
So it preserves arbitrary limits, finite coproducts and quotient maps.
\end{remark}

Given a closed equivalence relation $R\subset S\times S$ on a compact Hausdorff space $S$,
one might ask whether
\[
    \underline{S}/\underline{R}=\underline{S/R},
\]
i.e., the Yoneda embedding preserves those quotients.
We can answer this in two ways:
the quotient $S/R$ is again a compact Hausdorff space and its kernel pair is exactly $R$.
Therefore, by the preceding lemma, the above equality holds.
An alternative proof goes as follows.

\begin{lemma}\label{lem:CHaus_to_cond_coe}\uses{prop:T1_is_cond}
Let $S$ be compact Hausdorff and $R\subseteq S\times S$ a closed equivalence relation. Then
\[
    \underline{S}/\underline{R}=\underline{S/R},
\]
where the factor denotes the coequalizer of the two coordinate projections.
\end{lemma}

\begin{proof}
It suffices to show that the sheaf $\underline{S/R}$ is pointwise for every extremally disconnected space $T$ the $\Set$-coequalizer of
\begin{center}\begin{tikzcd}
	\underline{R}(T) & \underline{S}(T)
	\arrow[shift left, from=1-1, to=1-2]
	\arrow[shift right, from=1-1, to=1-2]
\end{tikzcd}.\end{center}
Then the assertion follows as the sheaves are equal
if they agree on extremally disconnected sets,
and there the colimit can be computed in presheaves followed by sheafification,
which in turn is not necessary if one already has the pointwise equality
to the sheaf $\underline{S/R}$.

\begin{center}\begin{tikzcd}
	{C(T,R)} & {C(T,S)} & {C(T,S/R)}
	\arrow[shift left, from=1-1, to=1-2]
	\arrow[shift right, from=1-1, to=1-2]
	\arrow[from=1-2, to=1-3]
\end{tikzcd}\end{center}
is an equalizer.
Thus consider any set \(M\) and any map \(\phi\colon C(T,S)\to M\)
such that

\begin{center}\begin{tikzcd}
	&& M \\
	{C(T,R)} & {C(T,S)} & {C(T,S/R)}
	\arrow[shift left, from=2-1, to=2-2]
	\arrow[shift right, from=2-1, to=2-2]
	\arrow["\phi", from=2-2, to=1-3]
	\arrow[from=2-2, to=2-3]
	\arrow[dashed, from=2-3, to=1-3]
\end{tikzcd},\end{center}
Call the quotient map $q$ or $q_*$.
We are looking for the dashed arrow.
But as $T$ is projective and $S\twoheadrightarrow S/R$, the arrow $C(T,S)\to C(T,S/R)$ is epimorphic. So consider any set-section of this arrow, and compose it with $\phi$ to obtain the desired arrow. It remains to show independence of the section for the uniqueness of the dashed arrow. So consider any $f\in C(T,S/R)$ with two preimages $g,h\in C(T,S)$. Then

\begin{center}\begin{tikzcd}
	& S \\
	T & {S/R}
	\arrow["q", from=1-2, to=2-2]
	\arrow["h"', shift right, from=2-1, to=1-2]
	\arrow["g", shift left, from=2-1, to=1-2]
	\arrow["f"', from=2-1, to=2-2]
\end{tikzcd}\end{center}
and

\begin{center}\begin{tikzcd}
	T \\
	& R & S \\
	& S & {S/R}
	\arrow["k"{description}, dashed, from=1-1, to=2-2]
	\arrow["h", from=1-1, to=2-3]
	\arrow["g"', from=1-1, to=3-2]
	\arrow["{p_2}", from=2-2, to=2-3]
	\arrow["{p_2}"', from=2-2, to=3-2]
	\arrow["q", from=2-3, to=3-3]
	\arrow["q"', from=3-2, to=3-3]
\end{tikzcd}\end{center}
implies that $\phi h=\phi p_2k=\phi p_1 k=\phi g$.
\end{proof}

\begin{remark}
        With this at hand, we have the feeling that this is a good moment to summarise the subtle difference in the commutativity properties of Yoneda from $\profin$ and $\CHaus$.
    \begin{enumerate}[(i)]
       \item The embedding from $\CHaus$ commutes with all factors after equivalence relations, which means the following.
            For a compact Hausdorff space $S$ one can consider any (closed) equivalence relation, that is, a $\CHaus$ subset $R\sub S\times S$ fulfilling the usual conditions (see \ref{def:equivalence_relation}).
            Now we can build the factor (the coequalizer) $S/R$ in $\CHaus$, which may be computed in $\Top$
            on the grounds that the factor in $\Top$ of a compact Hausdorff space by a closed equivalence relation is automatically $\CHaus$ again.
            Note that the kernel pair of this factor map $S\to S/R$ is precisely given by the two canonical maps $R\to S$.
            For short: In $\CHaus$, every equivalence relation is effective.
            Now, the embedding commutes with any such coequalizer, which we summarise as
            \[\underline{S}/\underline{R}=\underline{S/R}.\]

      \item In contrast, from $\profin$ we preserve only those factors of profinite sets by so-called effective equivalence relations,
            i.e., equivalence relations that arise as coequalizers of kernel pairs of an arrow between two profinite sets (see \ref{def:equivalence_relation}).

            This means the following:
            Consider any profinite set $S$ with a profinite equivalence relation $R\sub S\times S$.
            $\profin$ is finitely cocomplete by \ref{lem:pro-completion_is_complete} and so we can consider the factor (coequalizer) in profinite sets,

\begin{center}\begin{tikzcd}
	R & S & T.
	\arrow[shift right, from=1-1, to=1-2]
	\arrow[shift left, from=1-1, to=1-2]
	\arrow[two heads, from=1-2, to=1-3]
\end{tikzcd}\end{center}
            In order to test whether or not this quotient is preserved by the Yoneda embedding,
            we build the $\CHaus$ or topological factor (coequalizer) $S/R$, which a priori is just a compact Hausdorff space.
            Precisely whenever this factor again is totally disconnected, the original equivalence relation in $\prof$ is effective,
            as in this case the equivalence relation constitutes the kernel pair
\begin{center}\begin{tikzcd}
	R & S & {S/R}
	\arrow[shift right, from=1-1, to=1-2]
	\arrow[shift left, from=1-1, to=1-2]
	\arrow[two heads, from=1-2, to=1-3]
  \end{tikzcd}\end{center}
            of $S\twoheadrightarrow S/R$.
            The quotients by this type of equivalence relation are exactly the ones that are preserved.

            Note that the preceding line of reasoning indeed shows the \enquote{exactly} part of the last assertion,
            for as $\CHaus\to\cond$ preserves all coequalizers by (closed) equivalence relations,
            if any such coequalizer in $\prof$ is preserved by $\prof\to\cond$,
            it must also agree with the one in $\CHaus$.

            As we will see in the example after this remark, both situations indeed occur,
            meaning that in $\prof$ not every equivalence relation is effective (arises from an arrow in profinite sets),
            and hence the embedding really does not commute with all quotients by closed equivalence relations.
\end{enumerate}
If one wishes to just remember a single general rule,
the coequalizer by an equivalence relation is preserved \emph{precisely} if the equivalence relation is effective,
regardless of talking about $\CHaus$ or $\prof$.
(But in $\CHaus$, every closed equivalence relation is effective.)

In other words, we see that the quotients by equivalence relations always should be taken in $\CHaus$, and the requirement for an \enquote{effective} equivalence relation in $\profin$
corresponds to preserving the coincidence that sometimes the colimit of $\prof$ spaces in $\CHaus$ is profinite.
\end{remark}

\begin{example}
  For an explicit example illustrating the possible failure of a closed equivalence relation in $\prof$ to be effective,
  consider the profinite set of decimal representations $\prod_{\N}\{0,\dots, 9\}$ and the profinite equivalence relation $R$ induced by identifying decimal representations ending
  in $d999\ldots$ ($d\in\{0,1,\ldots,8\}$) by the corresponding sequence ending in $(d+1)00\dots$.
  (Formally, $R$ can probably most easily be defined as the kernel pair of the canonical map $\prod_\N \{0,\dots, 9\}\to [0,1]$.)
  So we have

  \begin{center}\begin{tikzcd}
	R & {\prod_{\N}\{0,\dots,9\}.}
	\arrow[shift right, from=1-1, to=1-2]
	\arrow[shift left, from=1-1, to=1-2]
\end{tikzcd}\end{center}
The (perfectly reasonable) quotient of $\prod_{\N}\{0,\ldots,9\}$ by this equivalence relation that we obtain in $\cond$ is the quotient in $\CHaus$
(which is precisely what we want to happen), given by the unit interval
        \item \begin{center}\begin{tikzcd}
	R & {\prod_{\N}\{0,\dots,9\}} & {[0,1],}
	\arrow[shift right, from=1-1, to=1-2]
	\arrow[shift left, from=1-1, to=1-2]
	\arrow[two heads, from=1-2, to=1-3]
\end{tikzcd}\end{center}
which, of course, is not totally disconnected ($[0,1]$ is a connected subset with more than one point).

        But since $\profin$ is finitely cocomplete, we could also take the coequalizer in profinite sets:
\begin{center}\begin{tikzcd}
	R & {\prod_{\N}\{0,\dots,9\}} & T.
	\arrow[shift right, from=1-1, to=1-2]
	\arrow[shift left, from=1-1, to=1-2]
	\arrow[two heads, from=1-2, to=1-3]
\end{tikzcd}\end{center}
This coequalizer $T$ in $\prof$, however, is trivial:
Continuous maps from the (connected) unit interval into any profinite set are constant and hence factor through the point

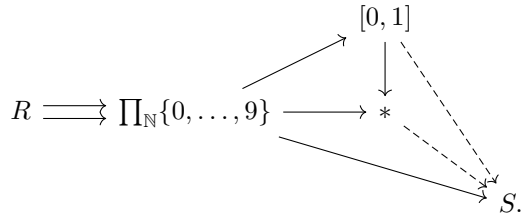
\begin{center}\begin{tikzcd}
	&& {[0,1]} \\
	R & {\prod_{\N}\{0,\dots,9\}} & \ast \\
	&&& S.
	\arrow[from=1-3, to=2-3]
	\arrow[dashed, from=1-3, to=3-4]
	\arrow[shift right, from=2-1, to=2-2]
	\arrow[shift left, from=2-1, to=2-2]
	\arrow[from=2-2, to=1-3]
	\arrow[from=2-2, to=2-3]
	\arrow[from=2-2, to=3-4]
	\arrow[dashed, from=2-3, to=3-4]
\end{tikzcd}.\end{center}
           Therefore, the equivalence relation $R$ in $\prof$ is not effective; the kernel pair of $S\to \ast$ is just given by $S\times S$, not by $R$.
           Hence, we (thankfully!) do not preserve this type of factor of profinite sets taken in $\prof$.
         \end{example}

Having seen the commutativity with coequalizers of kernel pairs, we may wonder about other coequalizers than just quotients by equivalence relations.
Concerning those, we have the following result.

\begin{example}\label{ex:yoneda-colim-not-comm}
To see that the Yoneda embedding $\CHaus\to \cond$ does not preserve all finite colimits

take, e.g., the coequalizer of
\begin{center}
\begin{tikzcd}
	{\mathbb{T}} & {\mathbb{T}}
    \arrow["1"', shift right, from=1-1, to=1-2]
	\arrow["{+\alpha}", shift left, from=1-1, to=1-2]
\end{tikzcd}
\end{center}
for any irrational $\alpha$.
The coequalizer of this diagram in $\CHaus$ is the point $*$
since the orbit of $\alpha$ is dense in $\T$.
However, example~\ref{ex:colim-not-pointwise} shows that the colimit
in $\cond$, evaluated at the point, is equal to $\T/\alpha\T$.

Moreover, the coequalizer of the irrational rotation does not arise from \emph{any} compact Hausdorff space $K$,
since then $\T\to K$ would be surjective in $\CHaus$ and hence a quotient map.
Thus the equivalence relation on $\T$ giving $K$ must contain the equivalence relation generated by $\alpha$ and is closed, hence is the whole product $\T\times\T$.
Therefore $K = *$ and we already have seen that this cannot be the case.

This in particular implies that Yoneda $\CHaus\to\cond$ does not even preserve reflexive coequalizers,
because then by proposition~\ref{prop:criterion_for_continuity} it would preserve all finite colimits.
\end{example}

We now know that the functor $X\mapsto\underline{X}$ preserves certain finite colimits.
The next lemma shows that also some infinite colimits are preserved.

\begin{lemma}\label{inclusions_CHaus_cond}
  Consider a sequence of T1 spaces $X_i$ ($i\in\N$) together with immersions
  \footnote{
    This means that it is a homeomorphism when restricted to its image (equipped with the subspace topology).
  }
  $X_i\hookrightarrow X_{i+1}.$
  Then the colimit $\varinjlim \underline{X_i}$ of the corresponding condensed sets agrees with the condensed set induced by the colimit $X$ of the $X_i$ taken in $\Top$.
\end{lemma}
\begin{proof}
  First of all,
  $X$ is T1 and thus indeed induces a condensed set $\underline{X}$.

  Now it suffices to show that for every $S\in \prof$, the canonical map
  \[
    \varinjlim C(S,X_i)\to C(S,X)
  \]
  is bijective,
  because filtered colimits are computed pointwise.

  For injectivity:
  Any $f\colon S\to X_i$ induces an arrow $S\to X_i\to X$
  and if two such composites are identical, then the morphisms $S\to X_i$ and $S\to X_j$ are equivalent
  (everybody carries subspace topology),
  hence there is only one induced inclusion $S\to X$.

  For surjectivity, we give a proof by contradiction.
  So assume we have a map $f\colon S\to X$ which does not factor through any $X_{i}$.
  Then for every $i\in\N$,
  there is $y_{i}\in S$ with $f(y_{i})\notin X_{i}$ (and therefore $f(y_{i})\notin X_{j}$ for all $j\le i$).
  Without loss of generality,
  assume all $f(y_{i})$ to be distinct.
  Consider $Z\coloneqq\{f(y_{i})\mid i\in\N\}$.
  $Z$ is closed because $Z\cap X_{i}$ is finite for all $i$ (and $X_{i}$ is T1, so finite subsets are closed)
  and $X$ carries the final topology.
  By the same token,
  every subset of $Z$ is closed.
  As $f(S)\sub X$ is compact ($S$ is),
  $Z\sub f(S)$ has an accumulation point $x=f(y_{i^{*}})$.
  Because $\{x\}$ is finite (and $\{x\}$ is not all of $Z$),
  $Z\setminus\{x\}=\{f(y_{i})\mid i^{*}\ne i\in\N\}$ (which, again, is closed)
  still has $x$ as an accumulation point, a contradiction.
\end{proof}

\begin{remark}
  The specialisation of this statement to the $\CHaus$ resp.\ Hausdorff setting can be found in \cite[1.2(4)]{scholze2019Analytic} or \cite[4.3.7]{bhatt2014proetale}
  with a different proof.

  Naturally,
  we wonder whether the conditions of this lemma can be weakened.
  \begin{itemize}
    \item First of all,
          the T1 condition is essential because otherwise $S\mapsto C(S,X_{i})$ does not yield a condensed set.
    \item The injectivity of the transition maps is necessary to ensure that the colimit is T1.
    \item That these transition maps be immersions (which is equivalent to requiring the $X_{i}$ to carry the subspace topology of the colimit)
          was also subtly used in the proof.
          Indeed,
          we used that $f\colon S\to X$ factoring through $X_{i}$ topologically is equivalent $f(S)\sub X_{i}$ (set-theoretically).
    \item That the colimit is sequential (and not just, say, filtered or over any well-ordered set) was used in the proof to show that $Z$ is closed.
          For the filtered case, a specific counterexample will be given below.
  \end{itemize}
\end{remark}

Even if all the $X_{i}$ are $\CHaus$, it is central to take the colimit in $\Top$ rather than $\CHaus$.
\begin{corollary}
If $S_i$ is a sequence of properly increasing $\CHaus$ spaces $S_i\subset S_{i+1}$, then the condensed set
\[S=\varinjlim_i S_i\]
is not $\CHaus$. 
\end{corollary}
\begin{proof}
Assume it were $\CHaus$.
Choose any sequence $x_i\in S_i\setminus S_{i-1}$.
By the universal property of $\beta\N$, there exists an extension $f\colon \beta\N\to S$, which clearly cannot factor through any $S_i$.
This is a contradiction to the last lemma or, alternatively,
the (category-theoretic) compactness of $\beta\N$.
\end{proof}

\begin{remark}
Note that these statements are no longer true for all filtered colimits rather than sequential colimits along inclusions.
For example, consider a first countable compact Hausdorff space $X$
and the filtered colimit of all finite unions
\footnote{Taking finite unions ensures that the diagram is filtered.}
of convergent sequences in $X$,
\[
    \varinjlim_{\substack{(f_i)_i\in X(\alpha\N)^k \\ k\in \N}}
            \bigcup_{i=1}^k f_i(\alpha\N).
\]
If one takes this colimit in $\Top$, one recovers $X$,
since a subset $A$ of $X$ is closed iff it contains all limit points of sequences in $A$.
But the colimit $Y$ in condensed sets is not given by $X$ (if $X$ is uncountable):
As it is filtered,
we can calculate it pointwise,
\[
    Y(S)
    = \varinjlim_{(f_i)_i\in X(\alpha\N)^k,\, k\in \N} \left(\underline{\bigcup_{i=1}^k f_i(\alpha\N)}\right)(S)
    = \varinjlim_{(f_i)_i\in X(\alpha\N)^k}C\left(S,\bigcup_{i=1}^k f_i(\alpha\N)\right).
\]
If $Y = \underline{X}$, then every morphism $S\to Y$ would factor through some $\bigcup_{i=1}^{k}f_{i}(\alpha\N)$, i.e., countable subset of $X$.
But this is not true for the projection $\beta X_d\twoheadrightarrow X$.
\end{remark}

However, arbitrary coproducts are better behaved, and in fact we can always calculate these in $\Top$.

\begin{lemma}\label{lem:cond-top-coprod}
    The embedding $\Top_{\mathrm{T1}}\to \cond, \, X\mapsto \underline{X}$ of T1 spaces into condensed sets preserves coproducts,
    i.e., for a set of T1-spaces $X_{i}$, we obtain
    \[\coprod_{i\in I}\underline{X_i}=\underline{\coprod_{i\in I}X_i}.\]
\end{lemma}

\begin{proof}
  Consider any coproduct of T1-spaces, $\coprod X_i$.
    First, write this as the filtered colimit of its finite subdiagrams,
   \[\coprod_{i\in I}X_i=\varinjlim_{F\sub I\,\mathrm{finite}}\coprod_{i\in F}X_i.\]
    Now, for any $S\in \CHaus$, any map $f\colon S\to\coprod_{i\in I}X_i$ factors through a finite sub-coproduct.
    (This follows by classical topological compactness of $S$, since $f^{-1}(X_i)$ forms a clopen decomposition of $S$.)
   This argument gives the first equality of the following and the second is just that filtered colimits may be computed pointwise:
   \[\left(\underline{\coprod_{i\in I} X_i}\right)(S)=\varinjlim_{F\sub I}\left(\left(\underline{\coprod_{i\in F}X_i}\right)(S)\right)
     =\left(\varinjlim_{F\sub I}\left(\underline{\coprod_{i\in F}X_i}\right)\right)(S).\]
    Hence, it suffices to show that the embedding $X\mapsto \underline{X}$ commutes with finite coproducts.

    For this, we again describe maps $f\colon S\to \coprod_{i=1}^n X_i$ (for $S\in\extr$).
    Any such map induces a finite clopen decomposition of $S$ (by $T_i=f^{-1}(X_i)$).
    Since on any of these $T_i$, the map $f$ maps only into $X_i$, the restriction of $f$ is given by an element of $C(T_i, X_i)\sub \coprod_{j} C(T_i, X_j)$.
    This implies that
    \[C\left(S, \coprod_{i=1}^n X_i\right)=\varinjlim_{S=\coprod_{j=1}^k T_i}\prod_{j=1}^k \coprod_{i=1}^n C(T_j, X_i).\]
    But the right hand side is precisely the formula of sheafification on $\extr$ of the presheaf $\coprod \underline{X_i}$,
    and thus $\underline{\coprod X_i}$ agrees with the colimit $\coprod \underline{X_i}$ in $\cond$.
\end{proof}

Clearly, the discrete world should be contained in the topological world by equipping everything with the discrete topology.
\begin{lemma}
The embedding of sets into condensed sets via $X\mapsto X_d\mapsto \underline{X_d}$ agrees with
    sheafification of the constant presheaf $X\mapsto\Sh(S\mapsto X)$.
    I.e., for any discrete set $X$ the corresponding condensed set on any $S=\varprojlim S_i$ is given by
    \[\underline{X_{d}}(S)=\varinjlim_i X^{S_i}\]
    This embedding $\Set\to\cond$ admits a bicontinuous right adjoint, which is given by evaluating at a point, $T\mapsto T(\ast)$.
    In particular, it is cocontinuous.
\end{lemma}
\begin{proof}
  Clearly, any $f\in C(S,X_d)$ has finite image (as it is compact),
  and hence it is constant on the elements of th induced finite clopen decomposition of $S$,
  i.e.,
  it is an element of
    \[\varinjlim X^{S_i}.\]
    The converse is clear, as any function that is constant on clopens of a clopen decomposition is continuous.
    This shows the fist assertion.
    Looking at the formula of sheafification explains the equality with the sheafification of the constant presheaf with value $X$.

    For the cocontinuity,
    note that mapping $X$ to the constant presheaf with value $X$ induces a bicontinuous functor $\Set\to\PSh(\extr)$.
    But sheafification is also cocontinuous.

    To show that $T\mapsto T(\ast)$ is right adjoint, take any $X\in \Set$.
    We want to show
    \[\hom_{\set}(X,T(\ast))=\hom(\underline{X_d}, T).\]
But rewriting $X=\coprod_{X}\ast$ and using cocontinuity of $X\mapsto \underline{X_d}$, as well as cocontinuity of the $\hom$ in the first coordinate,
    this reduces to the case $X=\ast$.
    But there the result holds by Yoneda,
    \[\hom_{\set}(\ast, T(\ast))=T(\ast)=\hom(\ast,T).\]
\end{proof}

Note that this is similar to the fact that forgetting the topology, $X\mapsto ?X$,
is the right adjoint to the functor equipping a set $Y$ with discrete topology $Y_d$,
\[C(Y_d, X)=\mathrm{Map}(Y,?X).\]

\begin{remark}
    This suggests to define a \idx{discrete condensed set} to a condensed set $X$ as $X_d=(X(\ast))_d$, and to call a condensed set discrete if it is of this form.
\end{remark}
\begin{example}
  We explain in unnecessary detail that $\N=\coprod_{\N}\{n\}$ in $\cond$ in order to showcase the practicality of the \enquote{abstract} (i.e., conceptually clear) arguments.
  The universal property is fulfilled,
  for any family of morphisms $(\eta^{n}\colon\{n\}\to X)_{n\in\N}$,
  sectionwise:
  $\eta^n_S\colon C(S,\{n\})=\{!_{S,\{n\}}\}\mapsto X(S\to \{n\})(\eta^n_{\{n\}}(!_{\{n\},\{n\}}))$.
  But this fully determines a morphism $\N\to X$ given by mapping any locally constant  $f=(n_i)_{S_i}\in C(S,\N)$
  for suitable $\coprod S_i=S$ and $n_i\in \N$ to the unique gluing of $\eta^n_{S_i}(!_{S_{i},\{n\}})$.
  Hence, $\N=\coprod_{N}\{n\}$ in $\cond$.
\end{example}

\begin{remark}
  As the embedding of nice topological spaces admits a left adjoint,
  we expect the limits of discrete sets (in condensed sets) to be computed in topological spaces rather than in $\Set$.
  Indeed, the product of discrete sets is not discrete, as $\profin$ demonstrates;
  i.e., $\Set\to \cond$ is not continuous.
\end{remark}

We close this section by investigating factor objects.
\begin{example}\label{ex:XXd}
Consider any T1 topological space $X$ and its discrete version $X_d\hookrightarrow X$. 
Then for any $S\in \extr$,
\[(X/X_d)(S)=C(S, X)/\sim,\]
where $f\sim g$ precisely if there exist finitely many points $F\sub  X$ such that 
\[A=f^{-1}(X\setminus F)=g^{-1}(X\setminus F)\sub S\]
and 
\[f_{\mid_{A}}=g_{\mid_A}.\]

This can be generalised to discrete subsets $Y_d\sub X$ in the canonical way, identifying $f$ and $g$ if there exists a finite $F\sub A$ with the same property as above.
\end{example}
\begin{proof}
This follows by the formula of sheafification,
\[(X/X_d)(S)=\varinjlim_{S=\bigsqcup_{i=1}^n T_i}\prod_{i=1}^n C(T_i, X)/C(T_i,X_d),\]
i.e., an element is given by continuous functions $f_{i}\colon T_{i}\to X$ (in other words, a continuous function $f\colon S\to X$),
where two functions $f,\,g\colon S\to X$ are identified if there exists a finite clopen decomposition $S=\bigsqcup T_i$
such that for all $i$,
the restrictions either agree or are both contained in $C(T_i,X_d)$ (and thus have finite image).
Taking the union of all images of the restrictions to those $T_i$ where the restrictions of $f$ and $g$ do not agree yields $F$.
\end{proof}

\begin{remark}
  We highly encourage the reader to think about why this is entirely sensible.
\end{remark}


\section{Special subcategories of condensed sets}

Our goal for the next sections is to explain show the following inclusions (arrows with tail)
and equalities (arrows in both directions) of full subcategories of condensed sets.
\begin{center}
\begin{tikzcd}
	&& \mathbf{CGWH} & {\mathrm{T1}} \\
	\extr & \CHaus & {\mathrm{compolog}} & {\mathrm{qs/qs}} \\
	{\mathrm{qcproj}} & {\mathrm{qcqs}} & {\mathrm{qs}} & \cond \\
	{\mathrm{compproj}} & {\mathrm{compqs}} & {\mathrm{Ind}(\CHaus)_{\mathrm{inj}}} & {\mathrm{Ind}(\CHaus)} \\
	{\mathrm{proj}} & {\mathrm{comp}} & {\mathrm{qc}} \\
	{\mathrm{retr}(\coprod\extr)} & {\mathrm{colim}_{\mathrm{fin}}(\extr)} & {\mathrm{epi}(\extr)}
	\arrow[hook, from=1-3, to=1-4]
	\arrow[hook', from=1-3, to=2-3]
	\arrow[hook, from=1-4, to=2-4]
	\arrow[hook, from=2-1, to=2-2]
	\arrow[hook, from=2-2, to=1-3]
	\arrow[tail reversed, from=2-4, to=3-4]
	\arrow[tail reversed, from=3-1, to=2-1]
	\arrow[hook, from=3-1, to=3-2]
	\arrow[tail reversed, from=3-2, to=2-2]
	\arrow[hook, from=3-2, to=3-3]
	\arrow[tail reversed, from=3-3, to=2-3]
	\arrow[hook, from=3-3, to=3-4]
	\arrow[tail reversed, from=4-1, to=3-1]
	\arrow[hook, from=4-1, to=5-1]
	\arrow[tail reversed, from=4-2, to=3-2]
	\arrow[hook, from=4-2, to=5-2]
	\arrow[tail reversed, from=4-3, to=3-3]
	\arrow[hook, from=4-3, to=4-4]
	\arrow[hook, from=4-4, to=3-4]
	\arrow[tail reversed, from=5-1, to=6-1]
	\arrow[hook, from=5-2, to=5-3]
	\arrow[tail reversed, from=5-2, to=6-2]
	\arrow[tail reversed, from=5-3, to=6-3]
	\arrow[hook, from=6-1, to=6-2]
	\arrow[hook, from=6-2, to=6-3]
\end{tikzcd}
\end{center}
\subsection{Projective and compact objects}
First, we will discuss the notions of projective (\ref{def:projectivity}), compact and compact projective (\ref{def:special_objects}) objects.
We start with the following observation. 

\begin{lemma}\label{lem:comp-ed-proj}
    Every compact Hausdorff space is compact and every extremally disconnected space is compact projective. 
\end{lemma}

\begin{proof}
Let \(K\) be a compact Hausdorff space.

By Yoneda,
\[
    \hom(K,\varinjlim T_i)=(\varinjlim T_i)(K).
\]
As every filtered colimit may be computed pointwise on $\CHaus$,
this is equal to
\[
    \varinjlim (T_i(K))=\varinjlim \hom(K,T_i).
\]
Thus \(K\) is a compact object in \(\cond\).

Now let \(S\) be an extremally disconnected space.
We show that \(S\) is compact projective.
This holds, as sifted colimits may be computed pointwise on \(\extr\) using the same Yoneda argument as above.
\end{proof}

This easily implies the following characterisation of projective objects.

\begin{proposition}\label{lem:projective-equiv}\leanok
In condensed sets, the projective objects are precisely those objects $P$ such that for every epimorphism $X\twoheadrightarrow Y$ and morphism $P\to Y$ there exists a lift 
\begin{center}\begin{tikzcd}
	& X \\
	P & Y
	\arrow[two heads, from=2-1, to=1-2]
	\arrow[from=1-2, to=2-2]
	\arrow[from=2-1, to=2-2]
\end{tikzcd}\end{center}
Equivalently, every epimorphism $X\twoheadrightarrow P$ admits a section
\[\begin{tikzcd}
	Y & P
	\arrow[shift left, two heads, from=1-1, to=1-2]
	\arrow[shift left, from=1-2, to=1-1]
\end{tikzcd}.\]

There are enough projectives (the condition of densely generating and separating agree in $\cond$ by \ref{lem:cosep_is_cogen})
and they are given precisely by retracts of coproducts of $\extr$. 
\end{proposition}
\begin{proof}
The first assertion follows, as every epimorphism is regular. 
For the second equivalence see \ref{lem:projectivity_via_sections}.
Since $\extr$ is compact projective and clearly densely generating, the class of projectives in $\cond$ agrees with retracts of coproducts of $\extr$ by \ref{lem:enough-projectives-separating}.
\end{proof}

\begin{question}
Are there more projectives in \(\Sh(\prof)\) than in $\cond$?
\end{question}

\begin{question}
We do not know whether there are injective objects in $\cond$ apart from the terminal object $\ast$, the proof for the nonexistence of injective abelian groups~\ref{lem:no_inj_cab} seems to break down.

By surjecting $\coprod S_i\twoheadrightarrow I$ onto the injective object,
then embedding $\coprod S_i\hookrightarrow\beta(\coprod S_i)$
(everything is computed in $\Top$),
and afterwards extending the surjection,
we see that every injective object admits an epimorphism from a $\CHaus$-space
(in notation that will be introduced later: it is quasicompact).

The canonical candidate for an injective object would be the unit interval.
\end{question}

\begin{proposition}
The compact objects $X\in \cond$ can be described in the following equivalent ways. 
\begin{enumerate}[(a)]
\item Those objects such that $\hom(X,-)$ commutes with filtered colimits.
\item The objects $X$ for which every arrow $X\to \varinjlim_{i\in I}Y_i$ into any limit factors through a finite subdiagram.
\item The finite colimits of $\extr$ or equivalently $\CHaus$.
\end{enumerate}
\end{proposition}

\begin{proof}
  (a) is the definition.
  It is equivalent to (b) by \ref{lem:comp_hom_endl}, for $\cond$ is cocomplete.

Finite colimits of compact objects are compact by \ref{lem:stability-compact-projective}, and thereby every finite colimit of $\CHaus$ is compact.
For the converse, consider any compact $T$. 
Take any colimit representation of $T$ as $\varinjlim K_i$ with $K_i\in \extr$.
Now, the identity $1\in \hom(T,T)=\hom(T,\varinjlim K_i)$ factors through some finite subdiagram $F\subseteq I$,
\[1\colon T\to \varinjlim_{F} K_i\to T,\]
so $T\to\varinjlim_{F}K_{i}$ is split monic.
By showing that it is also epimorphic,
it is an isomorphism (\ref{lem:split_epi_mono_is_iso}).
Indeed, the arrow $T=\varinjlim_I K_i\to \varinjlim_F K_i$ is epimorphic
by the uniqueness assertion of the universal property of $\varinjlim_{F}K_{i}$:
\begin{center}\begin{tikzcd}
	{K_i} & {K_j} & {K_k} \\
	{\varinjlim_I K_i} && {\varinjlim_F K_j} & H.
	\arrow[from=1-1, to=2-1]
	\arrow[from=1-2, to=2-1]
	\arrow[from=1-2, to=2-3]
	\arrow[from=1-3, to=2-1]
	\arrow[from=1-3, to=2-3]
	\arrow[from=2-1, to=2-3]
	\arrow[shift left, from=2-3, to=2-4]
	\arrow[shift right, from=2-3, to=2-4]
\end{tikzcd}\end{center}
\end{proof}

The proof actually shows the following statement.
\begin{corollary}\label{cor:filtered_comp_nicht_comp}
A colimit $\varinjlim K_i$ of compact objects $K_i$ is compact precisely if it is given as the colimit over a finite subdiagram.
\end{corollary}

\begin{example}\label{ex:XXd_not_comp}
For any infinite compact Hausdorff space $X$, the space $X/X_d$ is not compact.
\end{example}
\begin{proof}
Write
\[X_d=\coprod_{x\in X}\{x\}=\varinjlim_{F\sub X\,\mathrm{fin.}} F.\]
Since the subspace $F\sub X$ is discrete, we may identify $X/F$ with a compact Hausdorff space
(Yoneda commutes with factors of closed equivalence relations that induce $\CHaus$ factors).
Because colimits commute with colimits, we obtain
\[X/X_d=\varinjlim_{F\sub X\,\mathrm{fin.}} X/F.\]
But now, assuming $X/X_d$ were compact, this yields a contradiction:
$X/X_d$ is not isomorphic to any of the $X/F$ (it is not compact Hausdorff as it has just a single point).
\end{proof}

Intersecting these classes, we obtain the following description of compact projective objects. 
\begin{proposition}
 The class of compact and projective objects agrees with the compact projective objects and is precisely given by $\extr$. In particular, there are enough compact projectives in $\cond$.

\end{proposition}
\begin{proof}
To see that the compact projective objects are precisely the compact and projective objects, we use \ref{lem:comp-proj}.

Every object in $\extr$ is compact projective by \ref{lem:comp-ed-proj} and 
for the converse, 
every compact and projective object (as a compact object) is a finite colimit of objects in $\extr$,
and hence admits an epimorphism from a finite coproduct of $\extr$'s. 
By projectivity, this admits a section, and by \ref{lem:cond_yon_komm} this implies that the object is a retract of $\extr$ and thus $\extr$ itself. 

\end{proof}

In total, we obtain the following subdiagram. 

\begin{center}\begin{tikzcd}
	{\mathrm{retr}(\coprod\extr)} & \extr & \CHaus & {\mathrm{colim}_{\mathrm{fin}}(\extr)} \\
	{\mathrm{proj}} & {\mathrm{compproj}} && {\mathrm{comp}}
	\arrow[hook', from=1-2, to=1-1]
	\arrow[hook, from=1-2, to=1-3]
	\arrow[hook, from=1-3, to=1-4]
	\arrow[tail reversed, from=2-1, to=1-1]
	\arrow[tail reversed, from=2-2, to=1-2]
	\arrow[hook', from=2-2, to=2-1]
	\arrow[hook, from=2-2, to=2-4]
	\arrow[tail reversed, from=2-4, to=1-4]
\end{tikzcd}\end{center}

Let us now show that all inclusions are strict.
\begin{lemma}
\begin{enumerate}[(i)]
  \item No infinite coproduct of nonempty $\extr$ objects (e.g., infinite discrete sets) is $\extr$,
        showing that most projective objects are not compact projective.
\item Clearly, many $\CHaus$ spaces are not $\extr$.
\item The coequalizer of an irrational torus rotation is an example of a compact object which is not $\CHaus$.
\end{enumerate}
\end{lemma}
\begin{proof}
(i) follows from \ref{lem:cond-top-coprod} and the fact that infinite coproducts of nonempty sets in $\Top$ are never compact.
(iii) is \ref{ex:yoneda-colim-not-comm}.
Concerning (ii), $[0,1]$.
\end{proof}

\begin{remark}
  Among other things,
  we have shown that any compact condensed set $X$ can be written as a finite colimit of compact Hausdorff spaces
  (even extremally disconnected ones).
  Nonetheless,
  it is rarely possible to chose this finite colimit to be a kernel pair
  as this immediately implies $X\in\CHaus$ (which most compact $X$, of course, are not).
  We now prove this last claim.
\end{remark}

\begin{lemma}
Consider $X\twoheadrightarrow Y$ with $X\in \CHaus$. Then 
$Y\in \CHaus$ precisely if $X\times_Y X\in \CHaus$. 
\end{lemma}
\begin{proof}
If $Y\in \CHaus$, then clearly the kernel pair may be computed in $\CHaus$, as the embedding is continuous. 
Conversely, if $X\times_Y X\in \CHaus$, then the quotient by this closed equivalence relation may be computed in $\CHaus$ (Yoneda commutes with this),
and thereby $Y\in \CHaus$.
\end{proof}

\subsection{Quasicompactness}
\subsubsection{Quasicompact objects}
Next, we introduce a notion similar to classical topological compactness.
However, since the name \enquote{compact} is already reserved, one calls these spaces \idx{quasicompact}.

\begin{definition}\label{def:qc}
    Let \(X\) be an object of a category \(\mcC\) where coproducts exist.
    
    \begin{enumerate}[(i)]
\item     A family \((X_i\to X)_{i\in I}\) of morphisms with codomain \(X\) is called a \emph{cover} of \(X\) if the induced map
    \(\coprod X_i\to X\) is a regular epimorphism.

\item An object $X$ in a category $\mcC$ is called \idx{quasicompact},
if every cover $(X_i\to X)_{i\in I}$ has a finite subcover,
i.e., there exists a finite $J\subseteq I$ such that $(X_i\to X)_{i\in J}$ is a cover.
\end{enumerate}
\end{definition}

\begin{proposition}\label{prop:char_qc}
For any object $K\in \cond $ the following conditions are equivalent to $K$ being quasicompact.
\begin{enumerate}[(a)]
\item Every family of morphisms $Y_i\to K$ with $Y_i\in \cond$ which is jointly epimorphic admits a finite jointly epimorphic subfamily.
\item Every family of morphisms $Y_i\to K$ with $Y_i\in \extr$ (or $\prof$ or $\CHaus$) which is jointly epimorphic admits a finite jointly epimorphic subfamily.
\item There is a $\extr$/$\prof$/$\CHaus$ set $S$ with an epimorphism $S\twoheadrightarrow K$.
\item The induced map
    \[
         \varinjlim\hom(K,Y_j)\to \hom(K,\varinjlim Y_j)
    \]
    is injective for every filtered colimit \(\varinjlim Y_j\).
\end{enumerate}
\end{proposition}
\begin{proof}
The first item is clearly equivalent to quasicompactness, since in $\cond$ every epimorphism is regular.

For the next equivalence, assume (b) to follow (a).
Take any jointly surjective family $(Y_i\to K)$. As any $Y_i$ admits an epimorphism $\sqcup_j Y_{ij}\twoheadrightarrow Y_i$ by a set of extremally disconnected sets $Y_{ij}$,
gluing all these together we obtain a jointly surjective family $(Y_{ij}\to K)$.
If finitely many of these are surjective, clearly the corresponding $Y_i\to X$ are epi as well.

The implication from quasicompactness to (c) is clear as finite coproducts of profinites are profinite
and because any condensed set $X$ admits an epimorphism $\coprod S_i\twoheadrightarrow X$, and by quasicompactness finitely many $S_i$ suffice for this.
For the converse, use~\ref{lem:Im_of_qc_is_qc} and quasicompactness of $\extr$'s.

For the last equivalence see \cite[Exp. VI, Theorem 1.23]{Grothendieck1972a}.
\end{proof}

\begin{example}
 E.g., $[0,1]/[0,1]_d$  is quasicompact, as well as
 \[[0,1]\sqcup_{[0,1)} [0,1],\]
 the unit interval with doubled $1$.
\end{example}

\begin{example} $\N$ is not quasicompact.
 For this, let us consider $(\{n\}\hookrightarrow \N)_{n\in\N}$
\begin{center}\begin{tikzcd}
	{\{n\}} & \N \\
	\N
	\arrow[from=1-1, to=1-2]
	\arrow[from=1-1, to=2-1]
	\arrow["{\mathrm{Id}}", from=1-2, to=2-1]
\end{tikzcd}\end{center}
This forms a covering of $\N$, since $\coprod\{i\}=\N$ and the induced map is given by the identity. But clearly there is no finite subcover.
In fact, no infinite coproduct of nonempties is quasicompact.
\end{example}

\begin{corollary}
Every compact object is quasicompact.
\end{corollary}
\begin{proof}
This directly follows from the characterisation.
Another elementary proof of the fact that compact objects are quasicompact goes as follows.
Consider any compact $K$ and any jointly surjective family $(S_i\to K)$. Then, as $\coprod S_i\twoheadrightarrow K$,
by the epimorphism characterisation with $h=1$ there exist finitely many $\CHaus$ spaces $K_i\to K$ being jointly surjective with a morphisms $g\colon\coprod K_i\to \coprod S_i$ such that
\begin{center}\begin{tikzcd}
	{\coprod K_i} & {\coprod S_i} \\
	K
	\arrow[from=1-1, to=1-2]
	\arrow[two heads, from=1-1, to=2-1]
	\arrow[two heads, from=1-2, to=2-1]
\end{tikzcd}\end{center}
By \ref{lem:comp_hom_endl} the map $g$ factors through a finite subdiagram $\coprod_{F} S_i$ inducing the finite subcover $v$
\begin{center}\begin{tikzcd}
	{\coprod K_i} & {\coprod_F S_j} & {\coprod S_i} \\
	& K
	\arrow[from=1-1, to=1-2]
	\arrow[two heads, from=1-1, to=2-2]
	\arrow[from=1-2, to=1-3]
	\arrow["v", from=1-2, to=2-2]
	\arrow[two heads, from=1-3, to=2-2]
\end{tikzcd}\end{center}
But $v$ has to be epimorphic, as $\coprod K_i\twoheadrightarrow K$ is epimorphic.
\end{proof}
Not every quasicompact object is compact. 
\begin{example}\label{ex:indlim_not_comp_aber_qc}
Consider any filtered colimit $\varinjlim K_i$ of compact objects $K_i$. 
If an initial subdiagram of transition maps are epic but none of them isomorphisms, then $\varinjlim K_i$ is quasicompact but not compact.
\end{example}
\begin{proof}
  As all arrows in the diagram are epic and the diagram is filtered,
  every $K_{i}\to\varinjlim K_{i}$ is epic.
  Then $\varinjlim K_{i}$ is quasicompact as a factor of any of the the compact (hence quasicompact) $K_{i}$.
Assume that the colimit were isomorphic to one of the $K_i$.
Then the coordinate embedding of every $K_i\to K_j$ would be split monic, and thus an isomorphism, which is a contradiction to the assumption.
Now, by \ref{cor:filtered_comp_nicht_comp}, the object $\varinjlim K_i$ is not compact.
\end{proof}
Another class of examples are the following. 
\begin{example}\label{lem:XXd_not_compact}
  For any infinite $\CHaus$ space $X$, the condensed set $X/X_d$ is quasicompact (it admits the epimorphism $X\twoheadrightarrow X_d$) but not compact (by \ref{ex:XXd_not_comp}).
\end{example}

\begin{example}
The space $\R/\R_d$ is not quasicompact.
\end{example}
\begin{proof}
We choose the cover $([-n,n]\to \R\to \R/\R_d)$.
This indeed yields a cover. 
To see that there is no finite subcover, it suffices to show that for any $n$ the map
\[[-n,n]\to \R\to \R_d\]
is not epimorphic. 

For this it suffices to see that it is not epimorphic on $\beta\N$. 
Take any injection $\N\to [n+1,n+2]$ (e.g., a strictly monotonic convergent sequence).
This extends to $\beta\N\to [n+1,n+2]$ by local compactness of $\R$.
But clearly, any sequence whose image is contained in $[-n,n]$ has disjoint image, and since the image is not finite we conclude that this map $\beta\N\to \R$ is not equivalent to a map with image in $[-n,n]$ (with respect to the equivalence relation on $(X/X_d)(\beta\N)$ described in \ref{ex:XXd}), showing that $[-n,n](\beta\N)\to \R/\R_d(\beta\N)$ is not epimorphic. 

This implies that $\R/\R_d$ is not quasicompact by definition.
\end{proof}

\begin{question}
  Which other topological spaces yield compact or quasicompact objects in $\cond$?
  (We think it might not be all that hard to prove that classical compactness of a T1 space is equivalent to the corresponding condensed set being (quasi)compact.)
  Evidently, quasicompactness of $\underline{X}$ for T1 $X$ implies classical topological compactness.
\end{question}

\begin{lemma}
The compact projective objects in $\cond$ agree with the quasicompact projective objects.
\end{lemma}
\begin{proof}
Note that in the proof of the equivalence of compact projective objects with $\extr$ we just used the fact that there is an epimorphism from some $\extr$ onto every compact projective object. 
This holds precisely for quasicompact objects.
\end{proof}

\begin{remark}
Note that compact Hausdorff spaces are in general not compact or quasicompact objects in the category $\CHaus$. 

For an explicit example, consider $\beta\N=\coprod_{i\in \N}\{i\}$ (the coproduct is taken in $\CHaus$!), and the family
of embeddings $\{i\}\hookrightarrow \N\hookrightarrow \beta\N$.
This is a cover in $\CHaus$, since it induces the identity $\beta\N\to\beta\N$, but there is no possibility to choose a finite subcover. 
Passing to $\cond$ resolves this issue, as here $\coprod_{i\in \N}\{i\}=\N$,
and thereby the family $(\{i\}\hookrightarrow \beta\N)$ is not a cover, as they glue to the map
\[\N\hookrightarrow \beta\N.\]

More generally, any infinite coproduct $S=\coprod T_i$ in $\CHaus$ of nonempty sets is neither compact nor quasicompact in $\CHaus$.
(Of course, that infinity coproducts of non-initial objects are not (quasi)compact holds in many categories of interest.)
\end{remark}

This leads to the following picture of strict inclusions. 

\begin{center}\begin{tikzcd}
	{\mathrm{retr}(\coprod\extr)} & \extr & \CHaus & {\mathrm{colim}_{\mathrm{fin}}(\CHaus)} & {\mathrm{epi}(\CHaus)} \\
	{\mathrm{proj}} & {\mathrm{qcproj}} && {\mathrm{comp}} & {\mathrm{qc}} \\
	& {\mathrm{compproj}}
	\arrow[hook, from=1-2, to=1-3]
	\arrow[hook, from=1-3, to=1-4]
	\arrow[hook, from=1-4, to=1-5]
	\arrow[tail reversed, from=2-1, to=1-1]
	\arrow[tail reversed, from=2-2, to=1-2]
	\arrow[hook, from=2-2, to=2-1]
	\arrow[hook, from=2-2, to=2-4]
	\arrow[tail reversed, from=2-4, to=1-4]
	\arrow[hook, from=2-4, to=2-5]
	\arrow[tail reversed, from=2-5, to=1-5]
	\arrow[tail reversed, from=3-2, to=2-2]
\end{tikzcd}\end{center}

\begin{lemma}\label{lem:Im_of_qc_is_qc}
Images of quasicompact objects are quasicompact.
\end{lemma}
\begin{proof}
Consider any epimorphism $X\twoheadrightarrow Y$ with $X$ being quasicompact. We want to show that $Y$ is quasicompact. For this consider any jointly surjective family $(K_i\to Y)$ and take the pullbacks and distribute these through the coproduct, and use that pullbacks of epimorphisms are epic.
\begin{center}\begin{tikzcd}
	{\coprod K_i\times_Y X} & {\coprod K_i} \\
	X & Y
	\arrow[two heads, from=1-1, to=1-2]
	\arrow[two heads, from=1-1, to=2-1]
	\arrow[two heads, from=1-2, to=2-2]
	\arrow[two heads, from=2-1, to=2-2]
\end{tikzcd}\end{center}
But now finitely many of these $K_i$ suffice, i.e. the left coproduct can be chosen to be finite.
\begin{center}\begin{tikzcd}
	{\coprod K_i\times_Y X} & {\coprod K_i} \\
	X & Y
	\arrow[two heads, from=1-1, to=1-2]
	\arrow[two heads, from=1-1, to=2-1]
	\arrow[two heads, from=1-2, to=2-2]
	\arrow[two heads, from=2-1, to=2-2]
\end{tikzcd}\end{center}
But now, since the left arrows are epic, the arrow $\sqcup_{\mathrm{fin}}K_i\to Y$ has to be epic as well. Hence $Y$ is quasicompact.
\end{proof}

\begin{remark}
This yields the proof of the fact that every object with epimorphism from $\extr$ is quasicompact since epimorphisms agree with their image and images of quasicompact objects are quasicompact.
\end{remark}

\begin{theorem}[Tychonov]\label{thm:tychonov}
Arbitrary products and finite colimits of quasicompact objects in $\cond$ are quasicompact.
\end{theorem}

\begin{proof}
Consider any product $\prod T_i$ of quasicompact objects \(T_i\) and cover each of them with a single $\prof$ object, $K_i\twoheadrightarrow T_{i}$.
By lemma~\ref{lem:coprod_disj_stable} (vii), we obtain a surjection
\[
    \prod K_i \twoheadrightarrow \prod T_i.
\]
Since Yoneda embedding commutes with limits and \(\prof\) is complete,
we conclude that \(\prod T_i\) is quasicompact.

Next, let \(T_1,\ldots,T_n\) be quasicompact.
Then there exist profinite \(K_1,\ldots,K_n\) such that for every \(i\) we have a surjection \(K_i\twoheadrightarrow T_i\).
Then by lemma~\ref{lem:colim-preserve-epis} this yields a surjection
\[
    K_1\sqcup\ldots\sqcup K_n \twoheadrightarrow T_1\sqcup\ldots\sqcup T_n.
\]
Now the Yoneda embedding also commutes with finite coproducts,
so \(K_1\sqcup\ldots\sqcup K_n\) is profinite and \(T_1\sqcup\ldots\sqcup T_n\) quasicompact.

For an arbitrary finite colimit we represent the colimit as coproducts followed by a coequalizer.
Then we know that the coproducts (since they are finite) are both quasicompact
and hence we have a surjection from one of the quasicompact coproducts onto the finite colimit.
Therefore, the finite colimit must be quasicompact.
\end{proof}

\begin{remark}
Not every limit of quasicompact objects is quasicompact.

For this, consider the inclusion
\[
    (0,1)\hookrightarrow [0,1]
\]
in \(\cond\) (which arises from the usual inclusion and both intervals are equipped with their canonical topologies).
Since monics are effective, we know that
\begin{center}
\begin{tikzcd}
	{(0,1)} && {[0,1]} && {[0,1] \sqcup _{(0,1)}[0,1]}
	\arrow[hook, from=1-1, to=1-3]
	\arrow[shift left, from=1-3, to=1-5]
	\arrow[shift right, from=1-3, to=1-5]
\end{tikzcd}
\end{center}
is an equalizer diagram,
where \([0,1]\sqcup _{(0,1)}[0,1]\) denotes the cokernel pair of the inclusion \((0,1)\hookrightarrow [0,1]\).
Now \([0,1]\) is quasicompact and \([0,1]\sqcup _{(0,1)}[0,1]\) is as well since we have a surjection
\[
    [0,1]\sqcup [0,1]\twoheadrightarrow [0,1]\sqcup _{(0,1)}[0,1]
\]
from the quasicompact object \([0,1]\sqcup [0,1]\).
To see that \((0,1)\) is not quasicompact, we give a proof by contradiction:
if \((0,1)\) were quasicompact then there would exist a profinite set \(S\) with surjection \(S\twoheadrightarrow (0,1)\).
Since \((0,1)\) is compactly generated Hausdorff this surjection arises from a continuous map \(S\to (0,1)\) which has to be surjective.
A contradiction.
\end{remark}

\begin{proposition}\label{prop:qc-subobj-prof}
Quasicompact condensed subobjects of profinite sets are precisely given by closed profinite subsets.
\end{proposition}
\begin{proof}
  Take $S$ to be profinite and $T$ a quasicompact condensed subobject of $S$.
  By the previous lemma,
  there exists a profinite $S'$ with epimorphism $f\colon S'\twoheadrightarrow T\subset S$.
  The composition is an arrow between two profinite sets and hence can be interpreted as a continuous map.
  But this induces a mono-epi factorisation in $\CHaus$ with the usual topological image $f(S')$.
  But since this lifts to a mono-epi factorisation in $\cond$,
  we conclude by the uniqueness of the mono-epi factorisation, that $T\simeq f(S')$ is compact Hausdorff.

\end{proof}

\subsubsection{Quasicompact morphisms}
Using the definition of quasicompactness, we can define a notion of \enquote{proper} map (i.e., preimages of compact sets are compact).

Note that every continuous function $X\to Y$ between topological spaces with $Y$ being compactly generated Hausdorff is closed.

\begin{definition}\label{lem:qc-map-mono}
A morphism $f\colon X\to Y$ in $\cond$ is called \emph{quasicompact} if for any quasicompact $S$
with map $g\colon S\to Y$, the fiber product $S\times_Y X$ is quasicompact.
\end{definition}

\begin{remark}
  An object \(X\) is quasicompact iff \(X\to *\) is quasicompact.
  Furthermore,
  by surjecting on $S$ from some $\extr$/$\prof$/$\CHaus$ space,
  the condition can be weakened to $S$ being one of these.
\end{remark}

\begin{remark}
Another common definition of quasicompact morphism is to define a quasicompact morphism to be a quasicompact object in the commacategory. 
However, these notions do not coincide as the identity $\N\to\N$ shows.
This is a quasicompact map, but not quasicompact in the commacategory $\cond/\N$. 
To see this, first note that $(\N,\mathrm{Id})=\coprod (\{n\}\hookrightarrow \N)$ in the commacategory.
To see this, take any $h\colon X\to \N$ with maps $x_n\colon \{n\}\mapsto X$.
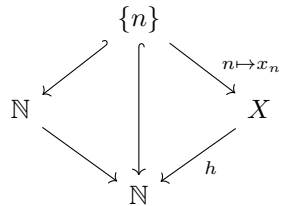
\begin{center}\begin{tikzcd}
	& {\{n\}} \\
	\N && X \\
	& \N
	\arrow[hook', from=1-2, to=2-1]
	\arrow["{n\mapsto x_n}"{pos=0.7}, from=1-2, to=2-3]
	\arrow[hook, from=1-2, to=3-2]
	\arrow[from=2-1, to=3-2]
	\arrow["h", from=2-3, to=3-2]
\end{tikzcd}\end{center}
Consider the natural candidate $g\colon \N\to X$ constructed with the universal property of $\N$ in $\cond$.
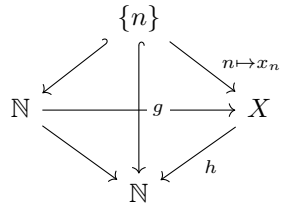
\begin{center}\begin{tikzcd}
	& {\{n\}} \\
	\N && X \\
	& \N
	\arrow[hook', from=1-2, to=2-1]
	\arrow["{n\mapsto x_n}"{pos=0.7}, from=1-2, to=2-3]
	\arrow[hook, from=1-2, to=3-2]
	\arrow["g"{description, pos=0.6}, from=2-1, to=2-3]
	\arrow[from=2-1, to=3-2]
	\arrow["h", from=2-3, to=3-2]
\end{tikzcd}\end{center}
the only thing to show is, that $hg=1_{\N}$. But as $g x_n=1$ we know that $g(h(n))=n$ on any constant function $n\in \N(S)$, and hence by uniqueness of gluings on any locally constant function $f\in \N(S)$. Thus, $hg=1$ and  $\N$ indeed is the coproduct in the comma category.

Hence we note that $\mathrm{Id}\colon (\N,\mathrm{Id})\to (\N,\mathrm{Id})$ is the gluing of all the $n\colon (\{n\},i_n)\to (\N, \mathrm{Id})$, and clearly epimorphic. Furthermore, the identity is a reflective epimorphism, as clearly the fibre product of the identity with the identity yields two times the identity as kernel pair, and the coequalizer of the identity with the identity is the identity. Hence we conclude that $(\{n\},i_n)$ is a covering of $(\N,\mathrm{Id})$.
But for any finite collection $N=\{n_i\}\sub\N$ clearly $(N,i_N)\to (\N,\mathrm{Id})$ is not epimorphic, as one could consider e.g. $(N\sqcup \N\sqcup\N, i_N\sqcup \mathrm{Id}\sqcup\mathrm{Id})$ and $f,g$ both mapping $N$ to the first $N$, but $f$ mapping $\N/\N$ to the first $\N$ and $g$ to the second $\N$.
\begin{center}\begin{tikzcd}
	N & {\N=N\sqcup\N\setminus N} && {N\sqcup\N\sqcup\N} \\
	\\
	&& \N
	\arrow["{i_N}", hook, from=1-1, to=1-2]
	\arrow["{i_N}"', from=1-1, to=3-3]
	\arrow["{f=\mathrm{Id}_N\sqcup i_1}"{pos=0.4}, shift left, from=1-2, to=1-4]
	\arrow["{g=\mathrm{Id}_N\sqcup i_2}"', shift right, from=1-2, to=1-4]
	\arrow["{\mathrm{Id}}", from=1-2, to=3-3]
	\arrow["{i_N\sqcup\mathrm{Id}\sqcup\mathrm{Id}}", from=1-4, to=3-3]
\end{tikzcd}\end{center}

Note that at first it feels weird to have a terminal object ($(\N,\mathrm{Id})$ is terminal in $\cond/\N$) with arrows towards it that are not epimorphic,
but if one considers $\emptyset\to \{\ast\}$ in $\Set$,
this feeling disappears.

So we have seen, that no finite subset of morphisms covers $(\N,\mathrm{Id})$ in the commacategory,
hence the identity $(\N,\mathrm{Id})\in \cond/\N$ does not yield a quasicompact object of $\cond/\N$.
\end{remark}

\begin{example}
    Every isomorphism is quasicompact,
    since the pullback of an isomorphism along any morphism is again an isomorphism.
\end{example}

\begin{corollary}\label{lem:qc-map-prop-1}
    Quasicompact morphisms are stable under base change and the composition of quasicompact morphisms is again quasicompact.
    In particular, if \(f\colon X\to Y\) is quasicompact and \(Y\) is quasicompact, then \(X\) is quasicompact.
\end{corollary}

\begin{proof}
    This follows directly from the previous two lemmas and the pasting lemma for pullback squares.
\end{proof}
\begin{remark}\label{rem:preimage_of_closed_closed}
Note the following familiar formulation of the fact \enquote{quasicompact morphisms are stable under base change}:
Preimages of closed subsets are closed subsets.
\end{remark}

\begin{remark}\label{lem:qc-map-prop-2}
Further properties one expects from quasicompact maps (see, e.g., \cite{Vakil2023}) include that
quasicompact morphisms are preserved under products,
    i.e., if \(X\to X'\) and \(Y\to Y'\) are quasicompact,
    then so is \(X\times Y\to X'\times Y'\).
    Furthermore, if we have morphisms of condensed sets sucht that \(fg = h\),
    and if \(h\) is quasicompact and the diagonal of \(f\) is quasicompact, then \(g\) is quasicompact.
\end{remark}

Having a definition of \enquote{proper} map, we can define closed subobjects (note that we are in a \enquote{compactly generated} setting).

\begin{definition}
    Let \(X\) be a condensed set.
    A subobject \(Y\hookrightarrow X\) is called \idx{closed} if the monic is quasicompact.
\end{definition}

\begin{remark}
The behaviour between closed and quasicompact subobjects is extremely analogous to the comparison between closed subsets and compact subsets of topological spaces:
\begin{enumerate}
\item Not every closed subobject is quasicompact; for example, since the identity is quasicompact, every condensed set is closed in itself,
    but not necessarily quasicompact. 
    However, if \(X\) is quasicompact, then by corollary~\ref{lem:qc-map-prop-1},
    the closed subobjects of $X$ are always quasicompact.
\item Not every quasicompact subobject is closed, 
consider e.g. 
\begin{center}\begin{tikzcd}
	{(0,1)} & {[0,1]} \\
	{[0,1]} & {[0,1]\sqcup_{(0,1)} [0,1]}
	\arrow[hook, from=1-1, to=1-2]
	\arrow[hook', from=1-1, to=2-1]
	\arrow[hook', from=1-2, to=2-2]
	\arrow[hook, from=2-1, to=2-2]
\end{tikzcd}\end{center}
together with the fact that pushouts along monics are pullbacks.

Requiring some \enquote{Hausdorff} condition developed in the next section, will lead to a result of the form \enquote{if $X$ is Hausdorff, then every compact subspace is closed}, see~\ref{lem:qc_closed_in_qs}.
\end{enumerate}
\end{remark}

\begin{definition}
    We say that a condensed set \(X\) has the property T1 or is a \idx{T1 space},
    if all maps \(\ast\to X\) are quasicompact.
    (I.e., all points are closed.)
\end{definition}

\begin{lemma}\label{lem:cond-T1-qc-points}s
If $X$ is a T1 space, then $\underline{X}$ is a T1 condensed set.

If $T$ is a T1 condensed set, then $T(\ast)_{\mathrm{top}}$ (defined as in~\ref{thm:cond_k_to_top}) is a compactly generated T1 space.
\end{lemma}
\begin{proof}

  Let $X\in\Top_{\mathrm{T1}}$
  To see that maps from points are quasicompact,
  we can show that for any profinite $S$ with $S\to X$, the fibre product $S\times_X \{x\}$ is profinite.
  But this is a closed subset of $S$, as the embedding $Y\mapsto \underline{Y}$ commutes with limits, and thus $\underline{S\times_X\{x\}}=\underline{S}\times_{\underline{X}}\underline{S}$.

  For the converse,
  $T(\ast)_\mathrm{top}$ is compactly generated, so it suffices to show the T1 property.
  Take any point $x\in T(\ast)$.
  Using Yoneda,
  we can interpret this as a map from a point to $T$, $\{x\}\to T$.
  Since the final topology installed on $T(\ast)_\mathrm{top}$ tells us that $\{x\}$ is closed precisely if it has a closed preimage under every morphism $S\to T$ from profinite $S$,
  consider any such $S\to T$.
  Then the preimage in $S$ corresponds to $S\times_T \{x\}\subseteq S$. But since the map is quasicompact by assumption, $S\times_T \{x\}$ is a quasicompact subobject of $S$.
  By~\ref{prop:qc-subobj-prof},
  this implies that it is given by a closed subspace of $S$.
\end{proof}

Using this, we can find the correct counterpart to the adjunctions between $\Top$ and $\condk$.

\begin{corollary}\label{rem:cond-T1-adjunction}
  We obtain an adjunction $X\to \underline{X}$ from T1 spaces to T1 condensed sets with left adjoint $T\mapsto T(\ast)_\mathrm{top}$,
  which furthermore is a faithful embedding and fully faithful when restricted to compactly generated T1 spaces.
 
  In particular,
  for any compactly generated T1 space $X$ and any T1 space $Y$ the identity
\[C(X,Y)=C(\underline{X}(\ast)_{\mathrm{top}},Y)=\hom(\underline{X},\underline{Y})\]
holds.
\end{corollary}

\subsection{Quasiseparatedness}
\quot{We must mention one more general nonsensical result about this direct limit topology.}{L. Waelbrock in \cite[p.~45]{Waelbrock1971}}

Having seen T1 as a first separation axiom and compactness, the next step is to define a suitable replacement for Hausdorff spaces.
The idea is to use the characterisation of Hausdorff spaces as those spaces with closed diagonal; or equivalently use the characterisation of weak Hausdorff spaces.
Note that the following are general notions from \cite{Grothendieck1973,Grothendieck1972, Grothendieck1972a}, specialized to our situation.

\begin{lemma}\label{lem:qs-equivalence}
A condensed set $X$ is called \idx{quasiseparated} if it fulfills one of the following equivalent conditions 
\begin{enumerate}[(a)]
\item The diagonal $\Delta\colon X\to X\times X$ is quasicompact.
\item For all quasicompact $S,T$ with maps $T\to X$ and $S\to X$
the fibre product \(S\times_X T\) is quasicompact.
\item For all $S\in \profin$  and all $f,g\colon S\to X$ the set $S\times_X S$ is profinite.
\end{enumerate}
\end{lemma}
We give an extremely detailed proof.
This is not at all necessary but we want to showcase that everything we say can, of course, be proven in such an elementary fashion.
\begin{proof}
\((a)\implies (b)\):
Let \(S,T\) be quasicompact and \(f\colon S\to X\), \(g\colon T\to X\) morphisms.
Then by Tychonov's theorem~\ref{thm:tychonov}, the product
\(P = S\times T\) is still quasicompact and we have a map
\(\phi:=f\times g\colon P \to X\times X\), given by
the universal property of the following diagram:
\begin{center}
    \begin{tikzcd}
	S &&& X \\
	{S\times T} &&& {X\times X} \\
	T &&& X
	\arrow["f", from=1-1, to=1-4]
	\arrow["{p_S}"', from=2-1, to=1-1]
	\arrow["{f\times g}", dashed, from=2-1, to=2-4]
	\arrow["{p_T}", from=2-1, to=3-1]
	\arrow["{p_1}"', from=2-4, to=1-4]
	\arrow["{p_2}", from=2-4, to=3-4]
	\arrow["g"', from=3-1, to=3-4]
\end{tikzcd}.
\end{center}
By assumption and lemma~\ref{lem:qc-map-mono}, the pullback
\begin{center}
\begin{tikzcd}
	{P\times_{X\times X}X} && X \\
	\\
	P && {X\times X}
	\arrow["{p_X}", from=1-1, to=1-3]
	\arrow["{p_P}"', from=1-1, to=3-1]
	\arrow["\Delta", from=1-3, to=3-3]
	\arrow["\phi"', from=3-1, to=3-3]
\end{tikzcd}
\end{center}
is quasicompact.
We show that the square
\begin{center}
    \begin{tikzcd}
	{P\times_{X\times X}X} && T \\
	\\
	S && X
	\arrow["{p_T\circ p_P}", from=1-1, to=1-3]
	\arrow["{p_S\circ p_P}"', from=1-1, to=3-1]
	\arrow["g", from=1-3, to=3-3]
	\arrow["f"', from=3-1, to=3-3]
\end{tikzcd}
\end{center}
is also a pullback diagram,
since then \(P\times_{X\times X}X\cong S\times_X T\) is quasicompact.
First of all, we have the commutativity of the diagram:
\[
    f\circ p_S\circ p_P
    = p_1\circ\phi\circ p_P
    = p_1\circ\Delta\circ p_X
    = p_X
    = p_2\circ\Delta\circ p_X
    = p_2\circ \phi\circ p_P
    = g\circ p_T\circ p_P,
\]
where we have used the universal property of \(\phi\),
\(\phi\circ p_P =\Delta\circ p_X\) and the universal property of the diagonal.
Now let \(Y\) be another condensed set with morphisms
\(h\colon Y\to S\), \(k\colon Y\to T\) such that
\(f\circ h = g\circ k\).
The maps \(h\) and \(k\) give a (unique) map
\(h\times k\colon Y\to S\times T\).
Then we have
\[
    p_1\circ\Delta\circ fh = f\circ h = f\circ p_S\circ (h\times k) = p_1\circ\phi\circ (h\times k)
\]
and
\[
        p_2\circ\Delta\circ fh = fh = g\circ k = g\circ p_T\circ (h\times k) = p_2\circ\phi\circ (h\times k),
\]
which means \(\Delta\circ fh = \phi\circ (h\times k)\).
Therefore, there exists a unique map \(\psi\colon Y\to P\times_{X\times X} X\)
such that the diagram
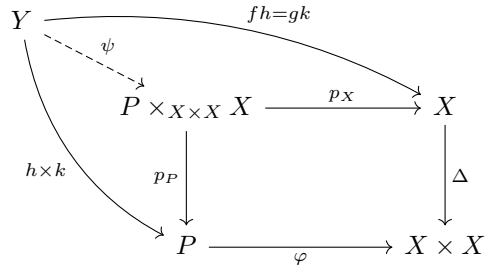
\begin{center}
    \begin{tikzcd}
	Y \\
	& {P\times_{X\times X}X} && X \\
	\\
	& P && {X\times X}
	\arrow["\psi", dashed, from=1-1, to=2-2]
	\arrow["{fh = gk}", curve={height=-18pt}, from=1-1, to=2-4]
	\arrow["{h\times k}"', curve={height=18pt}, from=1-1, to=4-2]
	\arrow["{p_X}", from=2-2, to=2-4]
	\arrow["{p_P}"', from=2-2, to=4-2]
	\arrow["\Delta", from=2-4, to=4-4]
	\arrow["\phi"', from=4-2, to=4-4]
\end{tikzcd}
\end{center}
is commutative.
But the composition \(p_P\circ \psi = h\times k\) exactly means
(by applying the further projections to \(S\) and \(T\))
that the diagram
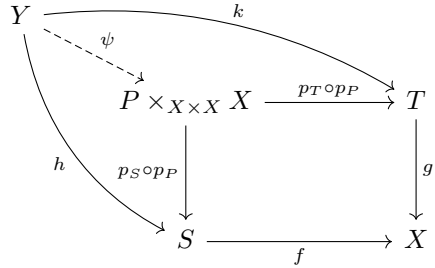
\begin{center}
    \begin{tikzcd}
	Y \\
	& {P\times_{X\times X}X} && T \\
	\\
	& S && X
	\arrow["\psi", dashed, from=1-1, to=2-2]
	\arrow["k", curve={height=-18pt}, from=1-1, to=2-4]
	\arrow["h"', curve={height=18pt}, from=1-1, to=4-2]
	\arrow["{p_T\circ p_P}", from=2-2, to=2-4]
	\arrow["{p_S\circ p_P}"', from=2-2, to=4-2]
	\arrow["g", from=2-4, to=4-4]
	\arrow["f"', from=4-2, to=4-4]
\end{tikzcd}
\end{center}
is commutative.
The uniqueness of \(\psi\) follows follows from the uniqueness of \(\psi\) above:
For any \(\gamma\colon Y\to P\times_{X\times X}X\) which makes the diagram commutative,
we have
\[
    p_X\gamma = fp_Sp_P\gamma = fh
\]
since \(p_X = fp_Sp_P\) (see above).
Furthermore, the equations \(p_Sp_P\gamma = h\) and \(p_Tp_P\gamma = k\)
imply
\[
    p_P\gamma = h\times k.
\]
Therefore \(\gamma\) also makes the pullback of \(\phi\) and \(\Delta\) commutative and by uniqueness of \(\psi\) it follows \(\gamma = \psi\).

\((b)\implies (c)\): This is immediate since every profinite set is quasicompact and hence \(S\times_X S\) is quasicompact.
It is also profinite as quasicompact subobject of the profinite set \(S\times S\) (since Yoneda commutes with products), see proposition~\ref{prop:qc-subobj-prof}.

\((c)\implies (a)\):
We show that the assertion of lemma~\ref{lem:qc-map-mono} holds.
So let \(S\) be any profinite set and \(\phi\colon S\to X\times X\) a morphism.
Define \(f = p_1\circ\phi\) and \(g = p_2\circ\phi\) maps from \(S\) to \(X\).
By assumption the pullback
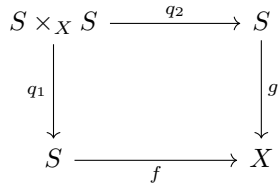
\begin{center}
    \begin{tikzcd}
	{S\times_X S} && S \\
	\\
	S && X
	\arrow["{q_2}", from=1-1, to=1-3]
	\arrow["{q_1}"', from=1-1, to=3-1]
	\arrow["g", from=1-3, to=3-3]
	\arrow["f"', from=3-1, to=3-3]
\end{tikzcd}
\end{center}
is profinite.
Therefore, the pullback
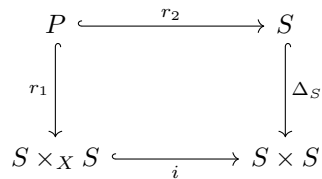
\begin{center}
    \begin{tikzcd}
	{P} && S \\
	\\
	{S\times_X S} && {S\times S}
	\arrow["{r_2}", hook, from=1-1, to=1-3]
	\arrow["{r_1}"', hook, from=1-1, to=3-1]
	\arrow["\Delta_S", hook, from=1-3, to=3-3]
	\arrow["i"', hook, from=3-1, to=3-3]
\end{tikzcd}
\end{center}
can be taken in \(\prof\) or equivalently in \(\cond\),
and hence \(P\) is profinite.
We show that the square
\begin{center}
    \begin{tikzcd}
	P && X \\
	\\
	S && {X\times X}
	\arrow["{fr_2}", from=1-1, to=1-3]
	\arrow["{r_2}"', from=1-1, to=3-1]
	\arrow["\Delta", from=1-3, to=3-3]
	\arrow["\phi"', from=3-1, to=3-3]
\end{tikzcd}
\end{center}
is a pullback diagram since then \(P\cong S\times_{X\times X} X\) is quasicompact
which means that the diagonal \(\Delta\colon X\to X\times X\) is quasicompact.

So let \(Y\) be a condensed set with morphisms \(h\colon Y\to S\) and \(k\colon Y\to X\) such that
\(\phi h = \Delta k\).
But then
\[
    fh = p_1\phi h = p_1\Delta k = k = p_2\Delta k = p_2\phi h = gh
\]
and therefore there is exactly one map
\(l\colon Y\to S\times_X S\) such that the diagram
\begin{center}
    \begin{tikzcd}
	Y \\
	& {S\times_X S} && S \\
	\\
	& S && X
	\arrow["l", dashed, from=1-1, to=2-2]
	\arrow["h", curve={height=-6pt}, from=1-1, to=2-4]
	\arrow["h"', curve={height=6pt}, from=1-1, to=4-2]
	\arrow["{q_2}", from=2-2, to=2-4]
	\arrow["{q_1}"', from=2-2, to=4-2]
	\arrow["g", from=2-4, to=4-4]
	\arrow["f"', from=4-2, to=4-4]
\end{tikzcd}
\end{center}
is commutative.
Now for the coordinate projections \(p'_1\colon S\times S\to S\)
and \(p'_2\colon S\times S\to S\) it follows that
\[
    p'_1 i l = q_1 l = h = p'_1\Delta_S h
\]
and
\[
    p'_2 i l = q_2 l = h = p'_2\Delta_S h
\]
which imply
\[
    il = \Delta_Sh.
\]
In particular there is a unique map \(\psi\colon Y\to P\)
such that the diagram
\begin{center}
    \begin{tikzcd}
	Y \\
	& P && S \\
	\\
	& {S\times_X S} && {S\times S}
	\arrow["\psi", dashed, from=1-1, to=2-2]
	\arrow["h", curve={height=-6pt}, from=1-1, to=2-4]
	\arrow["l"', curve={height=6pt}, from=1-1, to=4-2]
	\arrow["{r_2}", from=2-2, to=2-4]
	\arrow["{r_1}"', from=2-2, to=4-2]
	\arrow["{\Delta_S}", from=2-4, to=4-4]
	\arrow["i"', from=4-2, to=4-4]
\end{tikzcd}
\end{center}
is commutative.
But now
\[
    fr_2\psi = fh = p_1\phi h = p_1\Delta k = k
\]
which means that the diagram
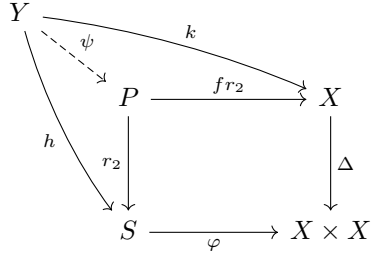
\begin{center}
    \begin{tikzcd}
	Y \\
	& P && X \\
	\\
	& S && {X\times X}
	\arrow["\psi", dashed, from=1-1, to=2-2]
	\arrow["k", curve={height=-6pt}, from=1-1, to=2-4]
	\arrow["h"', curve={height=6pt}, from=1-1, to=4-2]
	\arrow["{fr_2}", from=2-2, to=2-4]
	\arrow["{r_2}"', from=2-2, to=4-2]
	\arrow["\Delta", from=2-4, to=4-4]
	\arrow["\phi"', from=4-2, to=4-4]
\end{tikzcd}
\end{center}
is commutative.

The uniqueness of \(\psi\) follows as in the first part of the proof.
For any other \(\gamma\colon Y\to P\) with
\(r_2\gamma = h\), \(fr_2\gamma = k\)
we have
\[
    q_1r_1\gamma = p'_1ir_1\gamma = p'_1\Delta_Sr_2\gamma = r_2\gamma = h
\]
and by the same argument \(q_2r_1\gamma = h\).
This means that \(r_1\gamma = r_1\psi\) and by uniqueness of \(\psi\) it follows \(\psi = \gamma\).
\end{proof}

\begin{example}
For every weak Hausdorff space $X$, the induced $\underline{X}$ is quasiseparated. 

This is due to the fact that for $f,g \colon S\to \underline{K}$
for profinite $S$ and weak Hausdorff $X$ we can identify
\[
    S\times_{\underline{X}} S= \underline{S\times_X S},
\]
and the second is a compact Hausdorff space by weak Hausdorffness of $X$
(replacing $X$ by the $\CHaus$ image of $S\times_{X}S$ in $X$,
which does not change the pullback).

We will encounter this again in~\ref{lem:cgwh}
\end{example}

\begin{example}
\begin{enumerate}[(i)]
\item Every monic from a non-quasicompact object into quasicompact one has a non-quasiseparated cokernel pair (e.g., $(0,1)\sub [0,1]$)
\item For any subobject $X\to Y$ with $X\in \CHaus$, if the kernel pair is not $\CHaus$, then $Y$ is not quasiseparated.
\end{enumerate}
\end{example}
\begin{proof}
\begin{enumerate}
\item This is immediate, since every monic $X\hookrightarrow Y$ induces a pullback square
\begin{center}\begin{tikzcd}
	X & Y \\
	Y & {Y\sqcup_X Y.}
	\arrow[hook, from=1-1, to=1-2]
	\arrow[from=1-1, to=2-1]
	\arrow[from=1-2, to=2-2]
	\arrow[hook, from=2-1, to=2-2]
\end{tikzcd}\end{center}
\item This is clear by definition of the kernel pair.\qedhere
\end{enumerate}
\end{proof}
\begin{theorem}\label{thm:qs_chars}
For a condensed set $X$ the following assertions are equivalent:
\begin{enumerate}[(a)]
    \item $X$ is quasiseparated.
    \item For any map \(f\colon S\to X\) with \(S\) compact Hausdorff, the image $f(S)$ is again compact Hausdorff.
\item $X$ is filtered colimit of compact Hausdorff spaces with injective transition maps, i.e.,
        an object of $\mathrm{Ind}(\CHaus)$,
        possessing injective transition maps and the double limit formula for any two such objects holds.
 
\item $X$ is the filtered colimit of all its $\qcqs$ subobjects. 
\end{enumerate}
\end{theorem}
\begin{proof}
\((a)\implies (b)\):
First, let \(X\) be quasiseparated.
Then the image is given by the quotient of the closed equivalence relation,
\begin{center}\begin{tikzcd}
	{S\times_XS} & S & f(S)
	\arrow[shift right, from=1-1, to=1-2]
	\arrow[shift left, from=1-1, to=1-2]
	\arrow[two heads, from=1-2, to=1-3]
\end{tikzcd}.\end{center}
But as $X$ is quasiseparated, $S\times_X S$ is quasicompact.
Furthermore, $S\times_X S\hookrightarrow S\times S$,
so it is a quasicompact subobject of the compact Hausdorff space $S\times S$
and thereby compact Hausdorff.
Thus $f(S)$ is the quotient of a compact Hausdorff space by a closed equivalence relation,
and as the embedding preserves these, $f(S)$ itself is compact Hausdorff.

\((b)\implies (c)\):
Next, take any condensed set $T$ with all images from $\CHaus$ into $T$ being $\CHaus$.
Then,
\[
    T=\varinjlim_{S\to T}S,
\]
where the index category ranges over all profinite $S$,
ordered by $f\le g$ if $f=gh$ for some $h$.
But as every morphism has a mono-epi factorisation,
and the image of any $f\colon S\to T$ is compact Hausdorff,
hence the monomorphisms $S\hookrightarrow T$ are final in the index set.
This means that
\[
    T=\varinjlim_{S\hookrightarrow T}S.
\]
But the monics are filtered:
For any two $f\colon S\hookrightarrow T$, $g\colon S'\hookrightarrow T$
consider $g\colon I\hookrightarrow T$ via
\begin{center}\begin{tikzcd}
	S & {S\sqcup S'} & {S'} \\
	& I \\
	& T
	\arrow[from=1-1, to=1-2]
	\arrow[hook', from=1-1, to=3-2]
	\arrow[two heads, from=1-2, to=2-2]
	\arrow[from=1-3, to=1-2]
	\arrow[hook, from=1-3, to=3-2]
	\arrow["g", hook', from=2-2, to=3-2]
\end{tikzcd},\end{center}
and clearly any pair of parallel arrows is given by
\begin{center}\begin{tikzcd}
	S & {S'} \\
	T
	\arrow[shift left, from=1-1, to=1-2]
	\arrow[shift right, from=1-1, to=1-2]
	\arrow[hook, from=1-1, to=2-1]
	\arrow[hook, from=1-2, to=2-1]
\end{tikzcd}\end{center}
which implies, as the right arrow is monic, that the parallel arrows agree. Hence quasiseparated sets are given by filtered colimits with injective transition maps.
To show the double limit formula just note that any compact Hausdorff space is compact, and thus
\begin{align*}
    \hom(\varinjlim_I S_i,\varinjlim_J K_j)=\varprojlim_I \hom(S_i,\varinjlim_J K_j)=\varprojlim_I\varinjlim_J \hom(S_i,K_j).
\end{align*}

\((c)\implies (a)\):
For the other direction first note that $\CHaus$ spaces are quasiseparated,
and then that the quasiseparated condensed sets are closed under filtered colimits with injective transition maps.
For that take any filtered colimit with injective transition maps $T=\varinjlim K_i$ of quasiseparated $K_i$.
Take any maps $f,g\colon S\to T$ for profinite $S$.
As $S$ is compact, these maps factor through some finite index set,
and hence (by cofilteredness) through a single $K_i$,
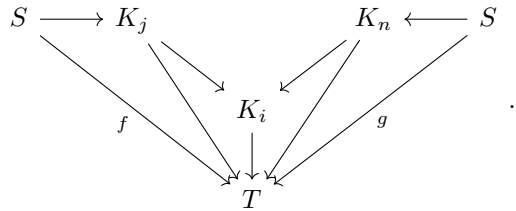
\begin{center}\begin{tikzcd}
	S & {K_j} && {K_n} & S \\
	&& {K_i} \\
	&& T
	\arrow[from=1-1, to=1-2]
	\arrow["f"', from=1-1, to=3-3]
	\arrow[from=1-2, to=2-3]
	\arrow[from=1-2, to=3-3]
	\arrow[from=1-4, to=2-3]
	\arrow[from=1-4, to=3-3]
	\arrow[from=1-5, to=1-4]
	\arrow["g", from=1-5, to=3-3]
	\arrow[from=2-3, to=3-3]
\end{tikzcd}.\end{center}
We want to show that the pullback
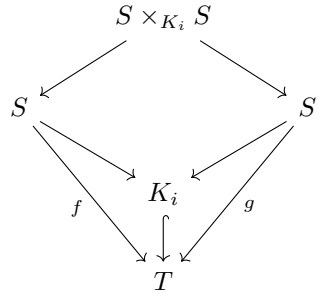
\begin{center}\begin{tikzcd}
	& {S\times_{K_i} S} \\
	S && S \\
	& {K_i} \\
	& T
	\arrow[from=1-2, to=2-1]
	\arrow[from=1-2, to=2-3]
	\arrow[from=2-1, to=3-2]
	\arrow["f"', from=2-1, to=4-2]
	\arrow[from=2-3, to=3-2]
	\arrow["g", from=2-3, to=4-2]
	\arrow[hook, from=3-2, to=4-2]
\end{tikzcd}\end{center}
is the pullback of $f,g$.
Note that $K_i\to T$ is monic, as all the transition maps in the diagram are monic.
But this directly implies the desired statement, since anything equalizing $f,g$
equalizes the two maps $S\to K_i$ composed with the monic $K_i\to K$
and thus as well equalizes the maps $S\to K_i$.
Hence, the pullback $S\times_T S$ can be identified with the pullback $S\times_{K_i} S$ which is quasicompact since \(K_i\) is quasiseparated.

\((c)\implies (d)\):
This is an easy consequence of~\ref{thm:qcqs_is_Chaus}.
\end{proof}

\begin{remark}
In the implication $(a)\implies(b)$, we could have argued by considering the mono-epi factorisation of \(f\),
\[
    f\colon S\twoheadrightarrow f(S)\hookrightarrow X.
\]
But then \(f(S)\) is quasicompact as it admits an epic of a compact Hausdorff space,
and it is quasiseparated as it is a subobject of the quasiseparated condensed set \(X\).
Therefore, \(f(S)\) is qcqs and so it is \(\CHaus\), see below.

But this uses that subobjects of qs are qs, in whose proof we used the lemma itself. However, one could also prove this other result using a different argument (essentially using that pullbacks are preserved under monics), allowing this alternative proof.
\end{remark}

\begin{remark}
  Note that in fact the proof shows that for any two filtered colimits of compact Hausdorff spaces,
  the double limit formula of $\Ind(\CHaus)$ holds.
  However, $\cond$ is not $\Ind(\prof)$, but rather admits the latter as a full subcategory.
  The reason for this is that the category of arrows $S\to X$ is in general not filtered as it fails the second condition of filtered categories.
This becomes true after replacing the homomorphisms with injections.
\end{remark}
\begin{theorem}\label{thm:qcqs_is_Chaus}
The following three full subcategories of \(\cond\) coincide:
\begin{enumerate}[(a)]
    \item \(\CHaus\),
    \item $\qcqs$, the full subcategory of quasicompact quasiseparated condensed sets,
    \item $\mathbf{compqs}$, the full subcategory of compact quasiseparated condensed sets.
\end{enumerate}
\end{theorem}
\begin{proof}
We have seen that compact Hausdorff spaces are qcqs, so let us show the converse.
Take any qcqs set $X$.
Then it is quasicompact and hence admits a surjection $S\twoheadrightarrow X$ from some profinite $S$.
But as any surjection is its image, and the image is compact Hausdorff,
$X$ is compact Hausdorff.

We know that any compact condensed set is quasicompact (corollary after ~\ref{prop:char_qc}), so $\mathbf{compqs}$ is a subcategory of $\qcqs$.
On the other hand, let \(X\) be qcqs. We want to show that it is compact (not using the cheap trick to say that it is compact Hausdorff and these are compact).

Then there exists a surjection \(S\twoheadrightarrow X\) with \(S\) profinite.
Taking the kernel pair \(S\times_X S\), which is compact Hausdorff since \(X\) is quasiseparated, and then since \(X\) is the coequalizer of its kernel pair,
we obtain that \(X\) is a finite colimit of \(\CHaus\) and hence compact,
see proposition ~\ref{prop:char_qc}.
\end{proof}

\begin{lemma}
Subobjects, arbitrary limits and coproducts of quasiseparated condensed sets are quasiseparated. 
Filtered colimits along inclusions of quasiseparated condensed sets are also quasiseparated.
\end{lemma}
\begin{proof}
Consider any quasiseparated set $X$ and subobject $Y\hookrightarrow X$.
Then for any $\CHaus$ set $S$ with morphism $S\to Y$,
the image of $S\to Y$ coincides with the image of $S\to Y\hookrightarrow X$,
and hence is $\CHaus$.

Next, consider any limit $X=\varprojlim X_i$ for quasiseparated $X_i$. Then for any maps $f,g\colon S\to X$, as limits commute with limits, we have that
\[
    S\times_X S=\varprojlim S\times_{X_i}S
\]
is a limit of profinite sets, and hence profinite since \(\prof\) is complete and Yoneda preserves limits.
Hence $X$ is quasiseparated.

Let \(X = \coprod_i X_i\) be a coproduct of quasiseparated \(X_i\).
Now for any $S\in \CHaus$ with map $S\to X$ the map factors through
a finite subdiagram (by compactness of \(S\))
\[
    S\to \coprod_{\mathrm{fin}} X_i\hookrightarrow \coprod X_i,
\]
and hence the image coincides with its image in a compact Hausdorff space
(Yoneda commutes with finite coproducts)
and hence is compact Hausdorff.

For the last assertion, see, e.g., the proof of~\ref{thm:qs_chars}.
\end{proof}
However, not every colimit of quasiseparated sets is quasiseparated,
as $[0,1]\sqcup_{(0,1)}[0,1]$ demonstrates.
We give another class of examples of this common phenomenon.

\begin{example}
Consider any filtered colimit $\varinjlim K_i$ of compact objects $K_i$. 
If an initial subdiagram of transition maps are epic but not isomorphisms, then $\varinjlim K_i$ not quasiseparated (and is quasicompact but not compact). 
\end{example}
\begin{proof}
By \ref{ex:indlim_not_comp_aber_qc} the object $\varinjlim K_i$ is quasicompact but not compact, and since $\qcqs=\mathbf{compqs}$, it cannot be quasiseparated.
\end{proof}

Next, we will show the classical version of Tychonov's theorem.

\begin{proposition}[Classical Tychonov]
The full subcategory \(\qcqs\subset\cond\) is complete.
In particular, \(\CHaus\) is complete.
\end{proposition}
Of course,
this follow from $\CHaus$ being complete and $\CHaus\to\cond$ being continuous.
But in fact,
we demonstrate how this classical result (completeness $\CHaus$) can be proven \enquote{intrinsically} in the condensed framework.
\begin{proof}
For the first claim, it suffices to show that if \(X_i\) are qcqs their limit \(\varprojlim X_i\) in \(\cond\) is qcqs.
By the previous lemma we know that \(\varprojlim X_i\) is quasiseparated.

Since products of $\qc$ objects are $\qc$, it hence just remains to show that equalizer of $\qcqs$ objects are $\qc$. 
For this represent an equalizer 
\begin{center}\begin{tikzcd}
	E & X & Y
	\arrow[hook, from=1-1, to=1-2]
	\arrow[shift left, from=1-2, to=1-3]
	\arrow[shift right, from=1-2, to=1-3]
\end{tikzcd}\end{center}
as pullback 
\begin{center}\begin{tikzcd}
	E & X \\
	{X\times_Y X} & {X\times X}
	\arrow[from=1-1, to=1-2]
	\arrow[from=1-1, to=2-1]
	\arrow[from=1-2, to=2-2]
	\arrow[from=2-1, to=2-2]
\end{tikzcd}.\end{center}
But since $X\times X$ is quasiseparated, $X$ is quasicompact and $X\times_Y X$ is quasicompact by quasiseparatedness of $Y$, the space $E$ is quasicompact. 
This implies the statement.
\end{proof}

\begin{lemma}\label{lem:qs-ev-point}
  If $X$ is quasiseparated, then morphisms $Y\to X$ can be differentiated on a point,
  meaning that the functor $X\mapsto X(\ast)$ from $\qs$ to $\Set$ is faithful.
\end{lemma}
\begin{proof}
See, e.g., \cite[p.~9]{scholze2019Analytic}.

We need to show that the natural transformation
\[
    \ev\colon\hom_{\cond}(-,Y)\Rightarrow\hom_{\Set}(-(*),Y(*))
\]
given by evaluation of a natural transformation $f\colon X\to Y$
at the point $f_*\colon X(*)\to Y(*)$
is injective for every $X\in\cond$.
Since evaluation on a point is bicontinuous by lemma~\ref{lem:Sh_on_terminal},
it suffices to assume that $X\in\extr$ by writing $X = \varinjlim S_i$ as the colimit
of objects $S_i\in\extr$.
Then using that limits preserve monomorphisms, see lemma~\ref{lem:colim-preserve-epis},
the monomorphisms of $\ev$ for $S_i$ glue together to a monomorphism
$\ev\colon\hom_{\cond}(X,Y)\to\hom_{\Set}(X(*),Y(*)$.
Now because $Y\in\qs$ is the filtered colimit of $K_j\in\CHaus$ and $X\in\extr$ is compact,
it suffices to assume that $Y\in\CHaus$ because the filtered colimit preserves monomorphisms in $\Set$ by lemma~\ref{prop:filtered_sifted_commutes_set}.
But for $X,Y\in\CHaus$ the claim is clear since the embedding is fully faithful
by the Yoneda lemma.
\end{proof}

\begin{remark}
Note that in the proof we did not make use of the fact that the colimit for $Y$
is along injective transition map,
so the above result holds also for all of $\Ind(\CHaus)$.
\end{remark}

We furthermore obtain the following surprising result. 
\begin{proposition}
The category $\cond$ is well-powered, and a subobject of any $\kappa$-condensed set is $\kappa$-condensed.
\end{proposition}
\begin{proof}

First assume that the condensed set $X$ is in $\extr_\kappa$.
Then every subobject $T\hookrightarrow X$ is quasiseparated and thus the filtered colimit of its $\qcqs$ subobjects. 
But these $\qcqs$ subobjects are also subobjects of $X$, and in particular $\kappa$-condensed (since $\qcqs\simeq \CHaus$). But this implies that $T$ is $\kappa$-condensed (since colimits of $\kappa$-condensed sets remain $\kappa$-condensed) and since there is just a set of $\qcqs$ subobjects of $X$, there is just a set of combinations of subobjects, and hence a set of such $T$.

Now, consider any $\kappa$-condensed set $X$. 
Then there exist $S_i\in \extr_\kappa$ with an epimorphism $\coprod S_i\twoheadrightarrow X$.

For any subobject $A\hookrightarrow X$, we consider the pullbacks

\begin{center}\begin{tikzcd}
	{S_i\times_X A} & {S_i} \\
	A & X
	\arrow[hook, from=1-1, to=1-2]
	\arrow[from=1-1, to=2-1]
	\arrow[from=1-2, to=2-2]
	\arrow[hook, from=2-1, to=2-2]
\end{tikzcd}\end{center}
Now, since coproducts commute with pullbacks, and pullbacks of epimorphisms are epimorphisms, we conclude 
\begin{center}\begin{tikzcd}
	{\coprod A\times_X S_i} & {\coprod S_i} \\
	A & X
	\arrow[hook, from=1-1, to=1-2]
	\arrow[two heads, from=1-1, to=2-1]
	\arrow[two heads, from=1-2, to=2-2]
	\arrow[hook, from=2-1, to=2-2]
\end{tikzcd}.\end{center}
Remark that any two different $A,A'$ induce different subobjects $\bigsqcup S_i\times A\ne \bigsqcup S_i\times A'$, as otherwise the composition $\bigsqcup S_i\times A\to \coprod S_i\to X$ would have two distinct mono-epi factorisations. 
\begin{center}\begin{tikzcd}
	& {\coprod A'\times_X S_i=\coprod A\times_X S_i} & {\coprod S_i} \\
	& A & X \\
	{A'}
	\arrow[hook, from=1-2, to=1-3]
	\arrow[two heads, from=1-2, to=2-2]
	\arrow[two heads, from=1-2, to=3-1]
	\arrow[two heads, from=1-3, to=2-3]
	\arrow[hook, from=2-2, to=2-3]
	\arrow[hook, from=3-1, to=2-3]
\end{tikzcd}\end{center} 
Hence it suffices to see that there is just a set of (automatically $\kappa$-condensed) subobjects of $S_i$ (or of $\coprod S_i$). But this was the first case.

Thus we conclude, that $\bigsqcup S_i\times A$ is $\kappa$-condensed, and thus also $A$. 
\end{proof}
\begin{remark}\label{rem:saft_in_cond}
  Note that $\cond$ fulfills everything for the special adjoint functor theorem SAFT except of the small cogenerating set.
  $\condk$, however,
  (as a any Grothendieck topos) is perfectly well suited for SAFT
  and one can use this method to procure adjoints from/to $\cond$.
\end{remark}
\begin{question}
Does every continuous functor $\cond\to \mcC$ for locally small $\mcC$ admit a left adjoint?
\end{question}

\begin{lemma}\label{lem:qc_closed_in_qs}
In quasiseparated sets, every quasicompact subobject is closed.
\end{lemma}
\begin{proof}
Consider any quasicompact subobject $K$ of a quasiseparated condensed set $Y$, $K\hookrightarrow Y$. 
We want to show that it is closed, i.e., that for every profinite $S\to Y$ the pullback
\begin{center}\begin{tikzcd}
	{S\times_Y K} & S \\
	K & Y
	\arrow[from=1-1, to=1-2]
	\arrow[from=1-1, to=2-1]
	\arrow[from=1-2, to=2-2]
	\arrow[hook, from=2-1, to=2-2]
\end{tikzcd}\end{center}
is profinite.
But clearly, since $Y$ is quasiseparated, the image $W$ of $S\to Y$ is a qcqs subobject of $Y$,
inducing a factorisation (note that the union of two qc subobjects is qc, since there is a canonical epimorphism from the coproduct into the union, and the coproduct is qc),
\begin{center}\begin{tikzcd}
	{S\times_Y K} & S \\
	& W \\
	K & {K\cup W} & Y
	\arrow[from=1-1, to=1-2]
	\arrow[from=1-1, to=3-1]
	\arrow[two heads, from=1-2, to=2-2]
	\arrow[from=1-2, to=3-3]
	\arrow[hook, from=2-2, to=3-2]
	\arrow[hook, from=2-2, to=3-3]
	\arrow[hook, from=3-1, to=3-2]
	\arrow[hook, from=3-2, to=3-3]
\end{tikzcd}\end{center}
Since monics do not change the pullback, we can identify $S\times_Y K$ with $S\times_{K\cup W} K$. 
But these spaces are all compact Hausdorff spaces, and since the embedding is continuous this implies $S\times_{K\cup W} K$ to be qcqs. 
\end{proof}
In fact, closed subobjects of quasiseparated sets are quite approachable by classical topology.
\begin{proposition}
For any quasiseparated condensed set $X$, the closed subobjects of $X$ are equivalent to closed subspaces $W\sub X(\ast)_{\mathrm{top}}$ via sending
    $W\sub X(\ast)_{\mathrm{top}}$ to the condensed set given as
    \[S\mapsto X(S)\times_{\mathrm{Map}(S,X(\ast))}\mathrm{Map}(S,W)\]
(i.e., $Z(S)$ are those morphisms $S\to X$ whose underlying maps of sets land inside $W$.)
\end{proposition}
\begin{proof}
    This is \cite[4.13]{scholze2019Analytic}.
    Since $X$ is quasiseparated, it is a filtered union of $\CHaus$ along inclusions.
    Now, since the left adjoint commutes with colimits, and colimits are stable under pullback (and filtered colimits may be computed pointwise),
    we may show the statement for a single of these subobjects.
    Hence we may assume that $X$ is compact Hausdorff. 
    Here closed subobjects are quasicompact and in fact again $\CHaus$. But the embedding $\CHaus\to \cond$ commutes with images,
    this reduces us to the classical result that for compact Hausdorff spaces,
    continuous injective maps are already immersions.
\end{proof}

One might wonder how \enquote{far away} from the quasiseparated condensed set any given condensed set can be.
In fact,
every condensed set is quite accessible by quasiseparated ones.
\begin{lemma}
Any condensed set is a quotient of a qs set by a qs equivalence relation.
\end{lemma}
\begin{proof}
First note that every condensed set $X$ admits an epimorphism from some coproduct
$\coprod K_i$ of compact Hausdorff $K_i$.
Rewriting this, we obtain that $X$ is the colimit of
\begin{center}\begin{tikzcd}
	{\coprod K_i\times_X\coprod K_i} & {\coprod K_i.}
	\arrow[shift right, from=1-1, to=1-2]
	\arrow[shift left, from=1-1, to=1-2]
\end{tikzcd}\end{center}
But $\coprod K_i$ is quasiseparated by the lemma before.
Furthermore, the fibre product
\[
    \coprod K_i\times_X\coprod K_i=\coprod_{i,j} K_i\times_X K_j
\]
is quasiseparated, as every $K_i\times_X K_j$ is a subobject of the quasiseparated object $K_i\times K_j$.
\end{proof}

Having a notion of \enquote{Hausdorff spaces}, we can ask for a \enquote{Hausdorffification}.
\begin{lemma}
    There exists a \idx{quasiseparation}, being the left adjoint to $\qs\hookrightarrow\cond$.
    The unit $X\to \qs(X)$ is an epimorphism; and $X \mapsto \qs(X)$ commutes with finite products.
    \end{lemma}
\begin{proof}
    This is \cite[4.14]{scholze2019Analytic}.
    
    One essentially does the canonical thing: Every condensed set $X$ is a quotient of a $\qs$ set $X'$ by a $\qs$ equivalence relation $R\sub X'\times X'$.
Now one takes the closure of the equivalence relation, by using the above equivalence of closed subspaces and quasicompact injections. In fact, for every $X\to Y$ with quasiseparated $Y$, the induced morphism $X'\times_Y X'\sub X'\times X'$ is a quasicompact injection containing $R$. 
This easily implies the adjunction. 
For the commutativity with finite products and more details, see~4.14 in \cite{scholze2019Analytic}.
\end{proof}
\begin{question}
	How does the functor assigning to every condensed set $X$ the filtered colimit of all its $\qcqs$ subobjects fit into this story?
\end{question}

\begin{example}
For any $X\in\cond$, the Hausdorffification of $X/X_d$ ($=X/X(*)_{d}$) is given by $\ast$.
Indeed, for any qs $Y$,
\[\hom(X/X_d,Y)=Y(\ast).\]
\end{example}
\begin{proof}
We have that $X/X_d$ is the coequalizer of the two arrows
\begin{center}\begin{tikzcd}
	{X_d\times X_d} & X & {X/X_d,}
	\arrow["{i\pi_1 }", shift left=2, from=1-1, to=1-2]
	\arrow["{i\pi_2}"', shift right, from=1-1, to=1-2]
	\arrow[from=1-2, to=1-3]
\end{tikzcd}\end{center}
which implies that for any $Y\in \cond$ that
\[\hom(X/X_d,Y)=\mathrm{eq}(i\pi_1^{*},i\pi_2^{\ast}\colon\hom(X,Y)\to \hom(X_d\times X_d,Y)).\]
But as $Y$ is qs,
morphisms into $Y$ are determined on points
so that we may replace $\hom(X,Y)$ and $\hom(X_{d}\times X_{d},Y)$ by the set of all maps $X_{d}(\ast)\to Y(\ast)$ coming from morphisms (and similarly for $X_{d}\times X_{d}$).
Then, clearly, the equalizer is given by $Y(\ast)$.
This implies that $\ast$ is the quasiseparation of $X/X_{d}$, e.g., by Yoneda.
(It fulfills the universal property.)
\end{proof}

We close with some useful classes of non-quasiseparated objects.
\begin{lemma}
\begin{enumerate}[(i)]
\item Every quasicompact object that is not $\CHaus$ is not $\qs$.
\item Every nontrivial object with $T(\ast)=\ast$ is not $\qs$.
\end{enumerate}
\end{lemma}
\begin{proof}
The first assertion is clear, since $\qcqs=\CHaus$. 

For the second item, note that if such a $T$ is quasiseparated,
it admits at most one map from any qs $X$.
But any $S\in\prof$ is qs,
and so $T(S)=\hom(S,T)$ has at most (and hence, precisely, as $T$ is nonempty) one point.
Alternatively,
there is always precisely one morphism $T\to\ast$ and as $T(\ast)=\ast$,
there is (by Yoneda) also precisely one morphism $\ast\to T$.
And as there is only one endomorphism of $\ast$ and only one on $T$ (as $T$ is qs, two different ones would have to be distinguishable on a point),
these to arrows are mutually inverse.
\end{proof}

\begin{corollary}
The irrational torus rotation yields an example of a condensed set which is compact but not quasiseparated.

Every space $X/X_d$ is not quasiseparated if $X$ is nondiscrete.
\end{corollary}

\subsubsection{Compactological spaces and CGWH}
Next, we find the topological counterpart to quasiseparatedness.

\begin{definition}[\cite{Waelbrock1971}]
A bornology on a set $X$ is a collection $\mathcal{B}$ of subsets, called \textit{small} subsets, such that
\begin{itemize}
    \item every finite subset is small,
    \item subsets of small sets are small,
    \item finite unions of small subsets are small.
\end{itemize}
A compactological space is a set $X$ endowed with a compactology, i.e. it is equipped with a (automatically CGWH) topology, and a bornology such that
\begin{itemize}
    \item A set is closed precisely if its intersection with every closed small set is closed
    \item Small sets are contained in small compact Hausdorff subsets.
\end{itemize}

Morphisms of compactological spaces are continuous functions such that the image of every small set is small.
\end{definition}
\begin{example}
The topology of compactological spaces is CGWH , and every CGWH space can canonically be made to a compactological space.
\end{example}
\begin{proof}
This can be found e.g. in \cite{Commelin2023}.
\end{proof}
\begin{proposition}
The full subcategory of quasiseparated sets is equivalent to the category of compactological spaces.
\end{proposition}
\begin{proof}
This is a result in \cite{Waelbrock1971}, where he shows the equivalence to the Ind completion of $\CHaus$ with injective transition maps.
\end{proof}
\begin{lemma}\label{lem:cgwh}
  If $T$ is a quasiseparated condensed set, then $T(\ast)_{\mathrm{top}}$ is compactly generated weakly Hausdorff.
  Conversely, any compactly generated space $X$ with $\underline{X}$ being quasiseparated is CGWH.
\end{lemma}
\begin{proof} This is 2.16 of \cite{scholze2019condensed}
If $T$ is quasiseparated, by a previous lemma the image of any compact Hausdorff space in $T$ is compact Hausdorff, and hence $T(\ast)_{\mathrm{top}}$ is CGWH.
If on the other hand the image $X$ is compactly generated, the embedding is fully faithful, and hence any maps $f,g: S\to X$ can be interpreted in topological spaces. But now as the embedding commutes with limits clearly the fibre product agrees with the topological fibre product $S\times_X S$. This in turn is compact Hausdorff for all $f,g$ precisely if the space $X$ is weak Hausdorff.
\end{proof}
This implies, that some complicated compactological spaces might be identified with easier CGWH spaces when passing from the condensed setting to the topological. 
In particular this explains why the embedding of $T1$-spaces is not full -  there are some non compactly generated spaces, but we check continuity only against the compactly generated topology.

This completes the picture 
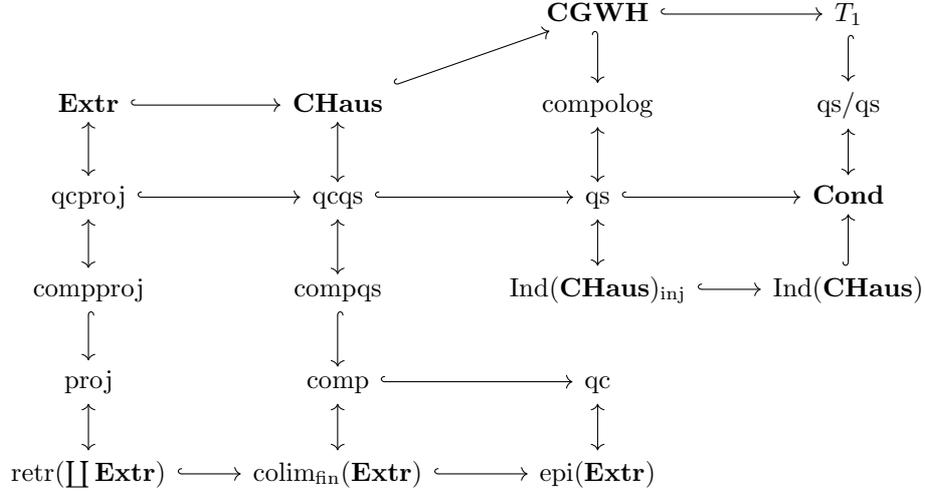
\begin{center}
\begin{tikzcd}
	&& \mathbf{CGWH} & {T_1} \\
	\extr & \CHaus & {\mathrm{compolog}} & {\mathrm{qs/qs}} \\
	{\mathrm{qcproj}} & {\mathrm{qcqs}} & {\mathrm{qs}} & \cond \\
	{\mathrm{compproj}} & {\mathrm{compqs}} & {\mathrm{Ind}(\CHaus)_{\mathrm{inj}}} & {\mathrm{Ind}(\CHaus)} \\
	{\mathrm{proj}} & {\mathrm{comp}} & {\mathrm{qc}} \\
	{\mathrm{retr}(\coprod\extr)} & {\mathrm{colim}_{\mathrm{fin}}(\extr)} & {\mathrm{epi}(\extr)}
	\arrow[hook, from=1-3, to=1-4]
	\arrow[hook', from=1-3, to=2-3]
	\arrow[hook, from=1-4, to=2-4]
	\arrow[hook, from=2-1, to=2-2]
	\arrow[hook, from=2-2, to=1-3]
	\arrow[tail reversed, from=2-4, to=3-4]
	\arrow[tail reversed, from=3-1, to=2-1]
	\arrow[hook, from=3-1, to=3-2]
	\arrow[tail reversed, from=3-2, to=2-2]
	\arrow[hook, from=3-2, to=3-3]
	\arrow[tail reversed, from=3-3, to=2-3]
	\arrow[hook, from=3-3, to=3-4]
	\arrow[tail reversed, from=4-1, to=3-1]
	\arrow[hook, from=4-1, to=5-1]
	\arrow[tail reversed, from=4-2, to=3-2]
	\arrow[hook, from=4-2, to=5-2]
	\arrow[tail reversed, from=4-3, to=3-3]
	\arrow[hook, from=4-3, to=4-4]
	\arrow[hook, from=4-4, to=3-4]
	\arrow[tail reversed, from=5-1, to=6-1]
	\arrow[hook, from=5-2, to=5-3]
	\arrow[tail reversed, from=5-2, to=6-2]
	\arrow[tail reversed, from=5-3, to=6-3]
	\arrow[hook, from=6-1, to=6-2]
	\arrow[hook, from=6-2, to=6-3]
\end{tikzcd}
\end{center}

\subsubsection{Quasiseparated morphisms}

Next, in analogy to quasicompact morphisms, we can define quasiseparated morphisms; inspired by \cite[8.3.1]{Vakil2023}.
However, since we do not yet have any concrete application of these morphisms, we omit the proofs of the following results.

\begin{conjecture}[About quasiseparated morphisms]
\begin{itemize}
\item 
    For a morphism \(X\to Y\) in \(\cond\) the following are equivalent:
    \begin{enumerate}[(a)]
        \item for any map \(S\to Y\) with \(S\) quasiseparated the pullback \(S\times_Y X\) is quasiseparated,
        \item the diagonal \(\Delta_Y\colon X\to X\times_Y X\) is quasicompact.
    \end{enumerate}
    We call these morphisms quasiseparated.
\item 
    Quasiseparated morphisms are stable under base change and composition.
    In particular they are preserved by products.
    Furthermore, if we have morphisms of condensed sets sucht that \(fg = h\),
    and if \(h\) is quasiseparated (and the diagonal of \(f\) is quasiseparated) then \(g\) is quasiseparated.
\item 
    A condensed set \(X\) is quasiseparated iff \(X\to *\) is quasiseparated, where \(*\) is the terminal object in \(\cond\).
\item 
    If \(f\colon X\to Y\) is a quasiseparated morphism with section \(s\), then \(s\) is quasicompact.
\item 
    If \(f\colon X\to Y\) is a morphism with \(X\) quasiseparated, then \(f\) is quasiseparated.
\end{itemize}
\end{conjecture}

\subsubsection{Further questions about condensed topology}

We end this chapter with some more questions that we were not able to answer yet.
\begin{question}
\begin{enumerate}
  \item (How) should one define open subobjects of condensed sets?
  \item Is there any use in open subobjects and do they form a complete distributive lattice (a locale)?
  \item Do the closed subobjects form a complete Heyting algebra?
  \item Using the topos theoretic approach of \cite{leroy2013theorielamesuredans}, can one develop a condensed topological measure theory?
  \item There are canonical translations of other separation properties (at least for T0, sober, T1, R0, T2, etc.).
        Are there interesting results about these, or do the differences become redundant?
  \item Is there a natural version of the condition of a normal space that yields a categorical version of the Tietze/Urysohn extension theorem?
        (We have the feeling that this is a heavily order-theoretical result.)
        What are the injective objects in $\cond$?
  \item The object $\{0,1\}$ seems to be cogenerating in $\qs$ (by the explicit description of morphisms as morphisms on underlying sets).
        Can we thereby infer any classical results?
  \item By the special adjoint functor theorem, there is a Stone-\v{C}ech compactification on $\cond$, meaning a left adjoint to $\CHaus\hookrightarrow \cond$ (explicitly described by composing quasiseparation, topologisation and afterwards classical Stone-\v{C}ech compactification).
        Does there also exist a profinitisation or extremalisation?
        How much of \cite{hindman2011algebra} has a natural translation to the condensed world?
\end{enumerate}
\end{question}

\clearpage{\thispagestyle{empty}\cleardoublepage}

\chapter[Homological Algebra]{Homological Algebra}
Having seen some basic category and some topological theory in form of condensed sets, our goal is to combine algebra and topology.
For us, coming from the analytic side, it hence is first necessary to collect some basic theory of algebraic category theory,
and first understand how a purely algebraic category theory would work, before being able to apply these algebraic ideas to topological situations.

Thus, the aim of this chapter is to describe our current understanding of some basic categorical homological algebra,
which will be needed to understand the combination of condensed structures with algebraic theories.

\section{Elementary homological algebra}
\subsection{Abelian categories}
One of the most important concepts in category theory and especially homological algebra are \idx{abelian categories}.
Roughly, these are categories behaving like categories of $R$-modules for some Ring $R$, e.g., using $R=\Z$, the category of abelian groups.
\footnote{We suggest for now to have the example of abelian groups in mind, and forget about categories of $R$-modules.
Especially the case of $R$ being a field might be sometimes slightly misleading, as here some even stronger properties hold and some important subtleties disappear,
as, e.g., every short exact sequence is split, as every object is injective.}

There is a whole zoo of different weakenings of the notion of abelian categories; \textbf{preadditive, semiadditive, additive, pseudoabelian, semiabelian, quasiabelian, preabelian, abelian,} \dots
We will focus on the notions relevant or helpful (to us).
Most of the results in this section can be found easily online (e.g. on the Stacks project), or e.g. in \cite{Freyd1964, Weibel1994} or any other introductory text on homological algebra.

There are two versions of defining abelian categories, the first one is by requiring stronger conditions on limits and colimits,
and the second one is by equipping the $\hom$ sets with more algebraic structure.

\begin{definition}[Direct sums and the $0$-object]\label{def:direct_sum}\uses{def:category, def:limit, def:limit, def:special_limits_and_colimits}\chapthree
A \idx{$0$-object} in a category $\mcA$ is both a terminal and an initial object.
    Having a $0$-object, one can define for any two objects $A,B\in\mcA$ a unique $0$-morphism as the composition

\begin{center}\begin{tikzcd}
	A & 0 & B.
	\arrow[from=1-1, to=1-2]
	\arrow[from=1-2, to=1-3]
\end{tikzcd}\end{center}

    In any category with a $0$-object, we define a \idx{direct sum} or \idx{biproduct} of two objects $A,B$ to be a simultaneous product and coproduct in the following way.

    We say a direct sum consists of an object $A\oplus B$ with morphisms $i_A\colon A\to A\oplus B$, $i_B\colon B\to A\oplus B$, $\pi_A\colon A\oplus B\to A$, $\pi_B\colon A\oplus B\to B$ so that
    $(A\oplus B, i_A,i_B)$ is a coproduct, $(A\oplus B, \pi_A,\pi_B)$, and furthermore,

\begin{center}\begin{tikzcd}
	& B \\
	A & {A\oplus B} & B \\
	& A
	\arrow["0", from=2-1, to=1-2]
	\arrow["{i_A}"{pos=0.6}, from=2-1, to=2-2]
	\arrow["{1_A}"', from=2-1, to=3-2]
	\arrow["{\pi_B}"'{pos=0.3}, from=2-2, to=1-2]
	\arrow["{\pi_A}"'{pos=0.3}, from=2-2, to=3-2]
	\arrow["{1_B}"', from=2-3, to=1-2]
	\arrow["{i_B}"', from=2-3, to=2-2]
	\arrow["0", from=2-3, to=3-2]
\end{tikzcd}\end{center}
commutes.
   Equivalently, one has a universal arrow $A\sqcup B\to A\times B$ induced by taking identities $A\to A$ and $B\to B$ and zeros $A\to B$, $B\to A$.
Then $A\sqcup B$ can be promoted to a biproduct if and only if this arrow is an isomorphism.
 The category admits all biproducts if and only if this is an isomorphism $(-)\sqcup (-)\to (-)\times(-)$ (automatically natural).

    A category with a zero object and all finite biproducts is called \idx{semiadditive}.
\end{definition}
On the enriched route to abelian categories, the first definition is the following.
\begin{definition}[Preadditive category]\label{def:preadditive}\uses{def:category}\chapthree
  A category is called \idx{preadditive} or \idx{$\Ab$-enriched},
  if it is enriched over abelian groups, meaning that the $\hom$ functor can be equipped with the structure of
  abelian groups $\hom(A,B)$ (i.e., that it factors over $\Ab\to \Set$) and composition is bilinear.
  We will thus freely regard $\hom(A,B)$ as either a set or an abelian group, depending on the context.\footnote{Viewing abelian groups as abelian group objects in $\Set$ helps a lot to be relaxed while doing this.}
  Some authors do also require a $0$-object in a preadditive category.
  It is easy to see that, assuming a $0$-object, the $0$-morphisms in a preadditive category are precisely the neutral elements of $\hom(A,B)$, explaining the name of $0$-object and $0$-morphism.
\end{definition}

Surprisingly, these properties are quite similar.
Given a semiadditive category, one can define addition of two morphisms
$f,g\colon A\to B$ by first gluing $1_A$ and $1_A$ (using that $A\oplus A$ is a product) to a morphism $A\to A\oplus A$ (called the diagonal morphism), and afterwards composing with the coproduct gluing of $f$ and $g$,
\[f+g\colon A\to A\oplus A\to B.\]
A so-called \idx{Eckmann-Hilton} argument shows that this defines a commutative monoid structure on the $\hom$ set (for which composition is bilinear).
Hence the only thing that misses to a preadditive category are additive inverses.

On the other hand, if a finite product or coproduct in any preadditive category exists, then it automatically is a biproduct, and the biproducts are precisely those diagrams

\begin{center}\begin{tikzcd}
	& B \\
	A & {A\oplus B} & B \\
	& A
	\arrow["0", from=2-1, to=1-2]
	\arrow["{i_A}"{pos=0.6}, from=2-1, to=2-2]
	\arrow["{1_A}"', from=2-1, to=3-2]
	\arrow["{\pi_B}"'{pos=0.3}, from=2-2, to=1-2]
	\arrow["{\pi_A}"'{pos=0.3}, from=2-2, to=3-2]
	\arrow["{1_B}"', from=2-3, to=1-2]
	\arrow["{i_B}"', from=2-3, to=2-2]
	\arrow["0", from=2-3, to=3-2]
\end{tikzcd}\end{center}
with $i_A\pi_A+i_B\pi_B=1_{A\oplus B}$.

This explains the following equivalent definitions.
\begin{definition}[Additive category]\label{def:additive}\uses{def:direct_sum, def:preadditive}
An \idx{additive category} is a preadditive category with finite products (equivalently, finite coproducts).
Equivalently, it is a semiadditive category with additive inverses for all morphisms.

A functor $F$ between additive categories is said to be \idx{additive}, if it preserves direct sums, or, equivalently,
    if for any $A,B$, the morphism $\hom(A,B)\to\hom(FA,FB)$ is an abelian group homomorphism.
\end{definition}
The following examples are meant to help distinguish the difference between these concepts.
\begin{example}
    \begin{itemize}
	  \item A \idx{ring} is nothing but a preadditive category with one object.
			Only the trivial Ring is an additive category.
	  \item The category of groups is not preadditive or semiadditive,
			as clearly neither the $\hom$ sets form abelian groups nor do coproducts (free products of groups) and products (product groups) agree
			(compare the cardinality of $F_2$ and $\Z^2$).
	  \item A good example of preadditive but not additive category is the category of odd-dimensional \idx{vector spaces}.
			(In this document, vector spaces are often implicitly meant to be always over either $\R$ or $\C$.)
\end{itemize}
\end{example}

Another important aspect of \enquote{abelian} structures is that the bilinearity of composition with addition simplifies the consideration of kernels and cokernels,
which in the presence of an additive category may be computed easily.

Considering any cone

\begin{center}\begin{tikzcd}
	C & A \\
	A & B
	\arrow["{g_1}", from=1-1, to=1-2]
	\arrow["{g_2}"', from=1-1, to=2-1]
	\arrow["f", from=1-2, to=2-2]
	\arrow["f"', from=2-1, to=2-2]
\end{tikzcd},\end{center}
the equality $fg_1=fg_2$ is equivalent to $f(g_1-g_2)=0=f(g_1-g_2)$, and thus we conclude that to any kernel pair

\begin{center}\begin{tikzcd}
	{A\times_BA} & A \\
	A & B,
	\arrow["{k_1}", from=1-1, to=1-2]
	\arrow["{k_2}"', from=1-1, to=2-1]
	\arrow["f", from=1-2, to=2-2]
	\arrow["f"', from=2-1, to=2-2]
\end{tikzcd}\end{center}
The difference $k=k_1-k_2$ induces an equalizer

\begin{center}
\begin{tikzcd}
	{A\times_BA} & A & B.
	\arrow["k"', from=1-1, to=1-2]
	\arrow["0"', shift right, from=1-2, to=1-3]
	\arrow["f", shift left, from=1-2, to=1-3]
\end{tikzcd}\end{center}
Conversely, any such equalizer $k\colon E\to A$ induces a pullback square via

\begin{center}\begin{tikzcd}
	{E} & A \\
	A & B.
	\arrow["k", from=1-1, to=1-2]
	\arrow["0"', from=1-1, to=2-1]
	\arrow["f", shift left, from=1-2, to=2-2]
	\arrow["f"', shift right, from=2-1, to=2-2]
\end{tikzcd}\end{center}
This explains the following definition.
\begin{definition}[(Co)kernels and (co)images]\label{def:ker_im}\uses{def:additive,def:kernel_pair, def:image-coimage}\chapthree
    Let $\mcA$ be an additive category, and $f\colon A\to B$.
    \begin{itemize}
        \item The \idx{kernel} $\ker(f)$ of $f$ is the equalizer of $f$ with $0$.
        By abuse of notation, we will denote the corresponding object, as well as the morphism, by $\ker(f)$.
        This is equivalent to $(\ker(f),0)$ being a kernel pair of $f$
        \item The \idx{cokernel} $\coker(f)$ of $f$ is the coequalizer of $f$ with $0$, which again is equivalent to $(0,\coker(f))$ being a cokernel pair of $f$.
        \item The \idx{image} of $f$ is the kernel of its cokernel, $\Im(f)=\ker(\coker(f))$.
        We will see that this will agree with the image as defined in \ref{def:image-coimage}.
        \item The \idx{coimage} of $f$ is the cokernel of its kernel, $\Coim(f)=\coker(\ker(f))$.
    \end{itemize}
\end{definition}

\begin{remark}
    Note that a morphism $f\colon A\to B$ in any preadditive category with $0$ is monic precisely if its kernel is $0\to A$,
    and is epic precisely if its cokernel is $B\to 0$.
\end{remark}
\begin{remark}
Analogously to above, in any preadditive category, equalizers

\begin{center}\begin{tikzcd}
	C & A & B
	\arrow["k", from=1-1, to=1-2]
	\arrow["f"', shift right, from=1-2, to=1-3]
	\arrow["g", shift left, from=1-2, to=1-3]
\end{tikzcd}\end{center}

are always kernels (an old name for equalizers, especially in germany, are \enquote{difference kernels}).

\begin{center}\begin{tikzcd}
	C & A & B
	\arrow["k", from=1-1, to=1-2]
	\arrow["0"', shift right, from=1-2, to=1-3]
	\arrow["{g-f}", shift left, from=1-2, to=1-3]
\end{tikzcd}\end{center}
    And dually, coequalizers are always cokernels.
\end{remark}
As, so far, we have only required finite (co)products to exist, there is no reason why kernels and cokernels should exist.
\begin{definition}[Preabelian]
An additive category $\mcA$ is called \idx{preabelian}, if all kernels and cokernels exist, or, equivalently, if it is finitely bicomplete.
\end{definition}
\begin{lemma}[Image factorisation in preabelian categories]\label{lem:preab_monoepi}\uses{def:additive, def:ker_im}
    The image and coimage of any morphism $f\colon A\to B$ in a preabelian category induce a unique factorisation as
    \[A\to \Coim(f)\to \Im(f)\to B,\]
    which we call \idx{image factorisation}.
\end{lemma}

We note that this morphism $\Coim(f)\to\Im(f)$ is not always an isomorphism.
Surprisingly, in general it does not even have to be a bimorphism (such categories are called \idx{semiabelian}).
For an explicit example, see, e.g. \url{https://math.stackexchange.com/questions/3631067/example-of-a-pre-abelian-category-but-not-a-semi-abelian-category}.

In many functional analytic categories this not an isomorphism.
Consider, for example, the dense embedding $i\colon\ell^1\hookrightarrow\ell^2$ in the category of Banach spaces.
This is a bimorphism, and hence has trivial kernel and cokernel, thus $\Im(i)=\ell^2$, $\Coim(i)=\ell^1$.
Now the morphism $\Coim(i)\to\Im(i)$ is simply $i$, which is clearly not an isomorphism.
This defect is basically due to the category not being balanced; in contrast, in many algebraic categories
as in abelian groups and $R$-modules the morphism $\Coim(f)\to\Im(f)$ indeed is an isomorphism.

This shows that preabelian categories do not yet have the strength to provide \enquote{new} algebraic tools for categories like Banach spaces.
For this, we need abelian categories.
\begin{definition}[Abelian category]\label{def:abelian}\uses{def:category, def:direct_sum, def:ker_im}\chapthree
There are many equivalent definitions of abelian categories.

Probably the shortest way is to define an \idx{abelian category} as a preabelian category $\mcA$ such that $\Coim(f)\to\Im(f)$ is always an isomorphism.

Equivalently, an abelian category is a preabelian category such that
\begin{itemize}
    \item every monic is a kernel of any morphism, or equivalently of its cokernel, making every monomorphism its image,
and
    \item every epic is a cokernel of any morphism or equivalently of its kernel, making every epimorphism its own coimage.
\end{itemize}
Spelling this out, we have the following \enquote{preabelian} formulation.
A category $\mcA$ is abelian precisely if
\begin{itemize}
    \item for any two objects $A,B$ the set $\hom(A,B)$ is an abelian group, fulfilling a two-sided distributivity of composition and addition,
    \item finite products in $\mcA$ exist,
    \item $\mcA$ has a $0$-object,
    \item kernels and cokernels exist, and
  	\item the morphism $\Coim(f)\to \Im(f)$ is always an isomorphism.
\end{itemize}

Following the \enquote{semiabelian} approach, one can reformulate the axioms of an abelian category in pure (co)limit form:
    A category $\mcA$ is abelian precisely if
    \begin{itemize}
        \item it admits finite direct sums (including a $0$-object),
        \item it admits kernels and cokernels (thus being finitely bicomplete), and
        \item every monic is a kernel and every epic is a cokernel.
\end{itemize}
    This formulation makes it easy to state and understand many robustness properties of abelian categories, as, e.g., under passage to functor categories etc.
\end{definition}
$\mcA$ will (here and elsewhere) usually denote an abelian category, although we try to always mention the assumption.

See \cite{Junhan2019} for a proof of this, and \cite{Freyd1964, Mitchell1965} or
section 12.3 in \cite[09SE]{SPA2023} for more.

Abelian categories share many of the nice properties of elementary topoi,
although for completely different reasons (or at least we do not know any explanation and the intersection is rather restricted).
\begin{theorem}[First properties of abelian categories]\label{thm:properties_of_abelian_categories}\uses{def:abelian, def:image-coimage}
Any abelian category $\mcA$ fulfills the following.
 \begin{itemize}
          \item Pullbacks and pushouts of monics/epics are monic/epic (as they are regular categories, see regular category in nLab.
          \footnote{\url{https://ncatlab.org/nlab/show/regular+category}}
          )
          Every pullback diagram along an epimorphism is also a pushout diagram, and any pushout diagram along a monomorphism is also a pullback diagram.
           More generally, pullbacks preserve kernels, and pushouts preserve cokernels.
        \item Every monic is regular and every epic is regular.
        \item The category is balanced.
        \item The image is indeed an image in the sense of \ref{def:image-coimage}, inducing an up to ismorphy unique
        mono-epi-factorisation
         \[A\twoheadrightarrow \Im(f)\hookrightarrow B\]
        \item Noetherian isomorphism theorems hold.
    \end{itemize}
\end{theorem}

\begin{lemma}
    In any abelian category $\mcA$, for any object $X$ there is a bijection between factor objects of $X$ and quotient objects.
\end{lemma}
\begin{proof}
  In an abelian category, any monic $A\sub X$ determines an epic by considering some cokernel of the embedding
  (For those already aware of this notion, they sit in a short exact sequence

\begin{center}\begin{tikzcd}
	0 & A & X & {X/A} & 0
	\arrow[from=1-1, to=1-2]
	\arrow[hook, from=1-2, to=1-3]
	\arrow[two heads, from=1-3, to=1-4]
	\arrow[from=1-4, to=1-5]
\end{tikzcd}.)\end{center}
    As cokernels do not change under precomposition with isomorphisms, this induces a well-defined map between subobjects (as isomorphy classes of monics) and factor objects (as isomorphy classes of epics).
    Note that we do not need an axiom of choice here, as after taking quotients the cokernel object is uniquely determined.

    Analogously, taking kernels induces a welldefined map from factor objects to subobjects.
    As any epimorphism is a cokernel of its kernel, and as every monomorphism is a kernel of its cokernel, taking equivalence classes of kernels is a two-sided inverse map.
    Hence this determines a bijection between (isomorphism classes) of subobjects and factor objects.
\end{proof}
\begin{remark}\label{rem:tiny-in-ab}
In an abelian category, the notions of tiny object and compact projective object agree, 
i.e., $\hom(P,-)$ commutes with all sifted colimits precisely if it commutes with all small colimits.
This is due to the fact that every colimit may be decomposed into a sifted colimit and a finite coproduct, and now the finite coproduct is equivalently a finite product, hence always preserved by $\hom(P,-)$.
\end{remark}

\begin{example}
    Note that, e.g., the category of Rings is not abelian, as $\Z\hookrightarrow \Q$ is a bimorphism but not an isomorphism.
    However, forgetting the additional structure and considering $\Z\hookrightarrow \Q$ as morphism of $\Z$-modules (abelian groups), this morphism is no longer epic.

    \end{example}
 In fact, $R$-modules describe the most important examples of abelian categories, and many statements that are true for $R$-modules have a straight forward generalisation to abelian categories.
    This philosophy is supported by the following theorem, known as the \idx{Freyd-Mitchell embedding}.
	In fact,
	it can be made precise in such a way that one can actually prove theorems about arbitrary abelian categories just by providing them for the categories are $R$-modules.

\begin{proposition}\label{lem:freyd_mitchell}\uses{def:abelian, thm:properties_of_abelian_categories}
  Any small abelian category is a full subcategory of the category $\ModR$ of $R$-modules for some ring $R$.
  Furthermore, all finite limits and colimits of the category of this abelian category may be computed in $\ModR$.
\end{proposition}
The structure preserving functors between abelian categories are called \idx{exact functors}.
This terminology will become clear later, when studying exact sequences.

\begin{definition}[Exact functor]\label{def:exact_functor}\uses{def:abelian}
	A functor between two abelian categories is called \idx{additive}, if it preserves direct sums (see the discussion in preadditive category).

	It is called \idx{right exact} if it preserves all finite colimits, and \idx{left exact} if it preserves all finite limits.
	The functor is called \idx{exact}, if it preserves all finite limits and finite colimits.
	Any left/right exact functor is additive.
\end{definition}

Often, abelian categories have further strong properties, which we will introduce now.
The most famous list of properties such a category may posses are the \idx{Grothendieck AB-conditions}, originally formulated in Grothendieck's famous Tohoku paper \cite{Grothendieck1957}.
\begin{definition}[Grothendieck Ab-conditions]\label{def:AB_conditions}\uses{def:abelian, def:special_limits_and_colimits}
An abelian category may or may not fulfill the following so-called Grothendieck axioms.
\begin{enumerate}
    \item[AB3] It admits all coproducts (thus all colimits).
    \item[AB3*] It admits all products (thereby all limits).

    \item[AB4] It is cocomplete and coproducts of monomorphisms are monomorphisms.
    \item[AB4*] It is complete and products of epimorphisms are epimorphisms.
    \item[AB5] It is cocomplete and filtered colimits commute with all finite limits.
    \item[AB5*] It is complete and cofiltered limits commute with all finite colimits.
    \item[AB6] It is cocomplete and filtered colimits distribute over products, meaning that for any set $\{I_j:j\in J\}$ of filtered categories and diagrams $i_j\mapsto A_{i_j}$ of form $I_j$, the natural map
    \[ \varinjlim_{(i_j)_{j\in J}\in \prod_j I_j}\prod_{j\in J} A_{i_j}\to \prod_{j\in J}\varinjlim_{i_j\in I_j} A_{i_j}\]
    is an isomorphism.
    \item[AB6*] It is complete and cofiltered limits distribute over coproducts, meaning that for any set $\{I_j:j\in J\}$ of cofiltered categories and diagrams $i_j\mapsto A_{i_j}$ of form $I_j$, the natural map
    \[\coprod_{j\in J}\varprojlim_{i_j\in I_j} A_{i_j}\to  \varprojlim_{(i_j)_{j\in J}\in \prod_j I_j}\coprod_{j\in J} A_{i_j} \]
    is an isomorphism.
\end{enumerate}
AB0 is the existence of a 0-object.
AB1 is the existence of kernels and cokernels, and AB2 is the condition that $\Coim(f)\to\Im(f)$ is an isomorphism.
\end{definition}
\begin{lemma}[AB-conditions of Abelian Groups]\label{lem:AB_conditions_of_abelian_groups}\uses{def:AB_conditions}
    The category of abelian groups satisfies AB3, AB3*, AB4, AB4*, AB5 and AB6.
\end{lemma}
See \cite[1.6]{Weibel1994} for more.
\begin{lemma}[Sheaves with values in abelian categories]\label{lem:sheaves_with_values_in_abelian_categories}\uses{def:abelian, def:sheaf}
Any category of sheaves with values in abelian groups forms an abelian category
satisfying AB3, AB3*, AB4, AB5.

	Furthermore, any category of sheaves with values in an abelian category is abelian.
\end{lemma}
\subsection{Symmetric monoidal categories}\label{ssec:mon-cats}
\quot{I hail a semigroup when I see one, and I seem to see them everywhere!}
{Einar Hille}

Another essential property we can observe in many categories is a tensor product.
There is an excellent general theory for working tensor products, which leads to monoidal categories, enriched categories and further.
However, we will restrict ourselves to the rudimentary version of just explaining (symmetric) monoidal categories.
For more on \idx{monoidal categories}, \idx{braided monoidal categories} and such, see, e.g., \cite{baez2023, baez1997introduction}.

\begin{definition}[Symmetric strict monoidal category]\label{def:strict_monoidal_cat}\uses{def:category}\chapthree
A \idx{strict monoidal category} is a strict 2-category with one object (recall the similarity to a monoid being a category with one object).

Equivalently,\footnote{after using a special case of a process called \idx{delooping}} it is a category with a functor $\otimes\colon \mcC\times \mcC\to\mcC$, called the \idx{tensor product}, that is associative

\begin{center}\begin{tikzcd}
	{\mcC\times\mcC\times\mcC} & {\mcC\times \mcC} \\
	{\mcC\times \mcC} & \mcC
	\arrow["{\otimes\times 1_\mcC}", from=1-1, to=1-2]
	\arrow["{1\times \otimes}"', from=1-1, to=2-1]
	\arrow["\otimes", from=1-2, to=2-2]
	\arrow["\otimes"', from=2-1, to=2-2]
\end{tikzcd}\end{center}
and admits a neutral object $e\in \mcC$ (called the \idx{tensor unit}), fulfilling $e\otimes -=-\otimes e=1_\mcC$.

A \idx{symmetric strict monoidal category} or \idx{permutative category} is a strict monoidal category such that there exists a natural isomorphism $\gamma\colon X\times Y\to Y\times X$ such that

\begin{center}\begin{tikzcd}
	{X\otimes e} & X & {X\otimes Y} & {Y\otimes X} \\
	& {e\otimes X} && {X\otimes Y} \\
	& {X\otimes Y\otimes Z} && {Z\otimes X\otimes Y} \\
	&& {X\otimes Z\otimes Y}
	\arrow["{1_X}", from=1-1, to=1-2]
	\arrow["\gamma"', from=1-1, to=2-2]
	\arrow["\gamma", from=1-3, to=1-4]
	\arrow["1"', from=1-3, to=2-4]
	\arrow["\gamma", from=1-4, to=2-4]
	\arrow["{1_X}"', from=2-2, to=1-2]
	\arrow["\gamma", from=3-2, to=3-4]
	\arrow["{1\otimes \gamma}"', from=3-2, to=4-3]
	\arrow["{\gamma\otimes 1}"', from=4-3, to=3-4]
\end{tikzcd}\end{center}
    See \cite[3.1]{Elmendorf2006}.
\end{definition}
Note that from a 2-categorical point of view it is a priori not canonical to expect the \enquote{strict} commutativity of any diagram.
Rather, one should impose the existence of natural 2-isomorphisms $(-\otimes-)\otimes-\to -\otimes(-\otimes-)$, fulfilling some identities.
This leads to \idx{unitors}, \idx{associators} and (weak) \idx{monoidal categories}.
\begin{definition}[Symmetric monoidal category]\label{def:symmetric_monoidal}\uses{def:category}\chapthree
    A \idx{symmetric monoidal category} is a category $\mcC$ together with
\begin{itemize}
    \item a \idx{tensor product} $\otimes\colon \mcC\times \mcC\to\mcC$,
    \item a \idx{tensor unit} $e\in\mcC$,
    \item a natural isomorphism
\[a_{X,Y,Z}\colon (X\otimes Y)\otimes Z\to X\otimes(Y\otimes Z),\]
called the \idx{associator},
    \item two natural isomorphisms
    \[\lambda_X\colon e\otimes X\to X,\qquad \rho_X\colon X\otimes e\to X\]
    called the left resp. right \idx{unitors},
\end{itemize}
    for which we further demand the following commutative diagrams:
    \begin{itemize}
        \item The \idx{triangle identity},
      \begin{center}\begin{tikzcd}
	{(A\otimes e)\otimes B} && {A\otimes(e\otimes B)} \\
	& {A\otimes B}
	\arrow[from=1-1, to=2-2]
	\arrow[from=1-1, to=1-3]
	\arrow[from=1-3, to=2-2]
\end{tikzcd}\end{center}
        \item And the \idx{pentagonal identity}
        \begin{center}\begin{tikzcd}
	& {(A\otimes B)\otimes(C\otimes D)} \\
	&& {((A\otimes B)\otimes C)\otimes D} \\
	{A\otimes(B\otimes(C\otimes D))} \\
	&& {(A\otimes(B\otimes C))\otimes D} \\
	& {A\otimes ((B\otimes C)\otimes D)}
	\arrow["{1_A\otimes a_{B,C,D}}", from=5-2, to=3-1]
	\arrow["{\alpha_{A,B,C\otimes D}}"', from=1-2, to=3-1]
	\arrow["{\alpha_{A\otimes B, C,D}}"', from=2-3, to=1-2]
	\arrow["{a_{A,B \otimes C, D}}", from=4-3, to=5-2]
	\arrow["{a_{A,B,C}\otimes 1_D}", from=2-3, to=4-3]
\end{tikzcd}.\end{center}
        \end{itemize}
    A \idx{braided monoidal category} further has a natural isomorphism
    \[B_{X,Y}\colon X\otimes Y\to Y\otimes X,\]
    called the \idx{braiding}, which has to fulfill
    the \idx{hexagonal identities}
\begin{center}\begin{tikzcd}
	& {(X\otimes Y)\otimes Z} & {X\otimes(Y\otimes Z)} \\
	{(Y\otimes X)\otimes Z} &&& {(Y\otimes Z)\otimes X} \\
	& {Y\otimes(X\otimes Z)} & {Y\otimes (Z\otimes X)}
	\arrow[from=1-2, to=2-1]
	\arrow[from=1-2, to=1-3]
	\arrow[from=2-1, to=3-2]
	\arrow[from=3-2, to=3-3]
	\arrow[from=1-3, to=2-4]
	\arrow[from=3-3, to=2-4]
\end{tikzcd}\end{center}
and
\begin{center}\begin{tikzcd}
	& {X\otimes (Y\otimes Z)} & {(X\otimes Y)\otimes Z} \\
	{X\otimes (Z\otimes Y)} &&& {(Y\otimes X)\otimes Z} \\
	& {(X\otimes Z)\otimes Y} & {Y\otimes (X\otimes Z)}
	\arrow[from=1-2, to=2-1]
	\arrow[from=1-2, to=1-3]
	\arrow[from=2-1, to=3-2]
	\arrow[from=3-2, to=3-3]
	\arrow[from=1-3, to=2-4]
	\arrow[from=3-3, to=2-4]
\end{tikzcd}.\end{center}

A \idx{symmetric monoidal category} further fulfills
 \[B_{Y,X}B_{X,Y}=1_{X\times Y}.\]
\end{definition}

\begin{remark} By Mac Lane's coherence theorem every (symmetric) weak monoidal category is equivalent (in the relevant sense) to a (symmetric) strict monoidal category.
Moreover, by Mac Lanes coherence theorem, the pentagonal and triangle identity already imply that every formal diagram constructed of brackets and units commutes.
Thus we often prefer the simpler notion of strict monoidal category, and omit some isomorphisms.

    Further remark that the passage from strict to weak monoidal categories is a special case of the passage from strict 2-categories to bicategories.
\end{remark}

Next, we define (symmetric) monoidal functors as those functors that are preserving this additional structure on a category, again everything up to isomorphisms satisfying appropriate identities.

\begin{definition}[Monoidal functor]\label{def:monoidal_functor}\uses{def:symmetric_monoidal}
A \idx{monoidal functor} $F\colon \mcC\to\mcD$ between monoidal categories is a functor equipped with
\begin{itemize}
    \item a natural isomorphism $\Psi\colon F(X)\otimes F(Y)\to F(X\otimes Y)$
    \item an isomorphism $\phi\colon e\to F(e)$
    \end{itemize}
such that
\begin{center}\begin{tikzcd}[column sep=tiny]
	& {F(X\otimes Y)\otimes F(Z)} & {F((X\otimes Y)\otimes Z)} \\
	{(F(X)\otimes F(Y))\otimes F(Z)} &&& {F(X\otimes(Y\otimes Z))} \\
	& {F(X)\otimes(F(Y)\otimes F(Z))} & {F(X)\otimes F(Y\otimes Z)}
	\arrow[from=2-1, to=1-2]
	\arrow[from=2-1, to=3-2]
	\arrow[from=1-2, to=1-3]
	\arrow[from=3-2, to=3-3]
	\arrow[from=1-3, to=2-4]
	\arrow[from=3-3, to=2-4]
\end{tikzcd},\end{center}
as well as
\begin{center}\begin{tikzcd}
	{e\otimes F(X)} & {F(X)} \\
	{F(e)\otimes F(X)} & {F(e\otimes X)}
	\arrow[from=1-1, to=2-1]
	\arrow[from=2-1, to=2-2]
	\arrow[from=1-1, to=1-2]
	\arrow[from=1-2, to=2-2]
\end{tikzcd}\end{center}
and
\begin{center}\begin{tikzcd}
	{F(X)\otimes e} & {F(X)} \\
	{F(X)\otimes F(e)} & {F(X\otimes e)}
	\arrow[from=1-1, to=2-1]
	\arrow[from=2-1, to=2-2]
	\arrow[from=1-1, to=1-2]
	\arrow[from=1-2, to=2-2]
\end{tikzcd}\end{center}
commute.

We call $F$ \idx{symmetric monoidal}, if $\mcC$ and $\mcD$ are symmetric monoidal and
\begin{center}\begin{tikzcd}
	{F(X)\otimes F(Y)} & {F(Y)\otimes F(X)} \\
	{F(X\otimes Y)} & {F(Y\otimes X)}
	\arrow[from=1-1, to=1-2]
	\arrow[from=1-1, to=2-1]
	\arrow[from=2-1, to=2-2]
	\arrow[from=1-2, to=2-2]
\end{tikzcd}\end{center}
commutes.
\end{definition}

Although this definition (at least to one of the authors) seems very unnatural and complicated,
in practice, monoidal categories are given canonically, and all these isomorphisms and relations are left implicit.
Luckily,
by now we are accustomed to leaving such compatibilities implicit.
Indeed,
we have often written $\prod_{i=1}^{n}X_{i}$ and never told anyone where to put the parentheses.
This does not cause any issues as long as our constructions are natural in any sense of the word.
(Imagine a functor between monoidal categories that satisfies $F(A\otimes B)\cong F(A)\otimes F(B)$ that does not fulfill the hexagonal identity above.
This should cause headaches.)

\begin{example}
  \begin{enumerate}
    \item In any category with finite products, the product yields a symmetric monoidal structure on the category.
		  Analogously with coproducts.
	\item The usual tensor product of abelian groups yields a symmetric monoidal structure on the category of abelian groups.
	\item Analogously, for any commutative (unital) ring $R$, the usual $R$-module tensor product induces a symmetric monoidal structure on $\ModR$.
  \end{enumerate}
\end{example}

A very important natural structure one might have in monoidal categories is an internal $\hom$.
\begin{definition}[Internal $\hom$]
Let $\mcC$ be a (symmetric) monoidal category.
We call a right adjoint to $-\otimes B$, if existent, the \idx{internal hom} $\ihom(B,-)$, yielding the \idx{currying} isomorphism
\[\hom(A\otimes B,C)\simeq \hom(A, \ihom(B,C)).\]
	If this exists (naturally) for all $B$, we call the category \idx{closed (symmetric) monoidal}.
The unit $\ihom(A,B)\otimes A\to B$ of the adjunction is called \idx{evaluation map}.
\end{definition}
\begin{example}
  \begin{enumerate}
	\item Exponential objects are precisely the internal hom with respect to the symmetric monoidal structure induced by the product.
	\item For abelian groups and $R$-modules they are the usual abelian group resp. $R$-module structures on the $\hom$ sets, induced by pointwise algebraic operations.
	\item The category of endofunctors on any category becomes strict symmetric monoidal with respect to the composition.
	\item Equipping the tensor product of abelian groups with compact open topology yields an internal hom for locally compact abelian groups with tensor product being the usual product.
  \end{enumerate}
\end{example}

\begin{remark}
  The term \enquote{internal $\hom$} is always reasonable for functors with this adjunction property.
  Indeed,
  one can always extract the \enquote{underlying} $\hom$ by taking $\hom(e,-)$ (where $e$ denotes the tensor unit),
  \[
	\hom(e,\ihom(A,B))=\hom(e\otimes A,B)=\hom(A,B).
  \]
\end{remark}

\begin{lemma}[$\Sh(\mcC,\mcA)$ monoidal]\label{lem:sh-mon}
  For any (symmetric) monoidal category $\mcA$ and any site $\mcC$ such that sheafification exists,
  the category $\Sh(\mcC,\mcA)$ admits a unique (symmetric) monoidal structure such that sheafification is monoidal
  (where we equip the presheaf category with pointwise tensor product).
  It is given by first defining the tensor product in presheaves as $(F\otimes_{\PSh} G)(C)=F(C)\otimes G(C)$ and afterwards defining $F\otimes G=\Sh(F\otimes_{\PSh} G)$.
\end{lemma}

\subsubsection{Monoids and modules}\label{subsubsec:monoids_and_modules}
\begin{definition}[Monoid]\label{def:monoid} A \idx{monoid} in a monoidal category $(\mcC,\otimes,e)$ consists of an object $M\in \mcC$
	together with a morphism $m\colon M\otimes M\to M$ and a morphism $\eps\colon e\to M$ such that associativity and unitality hold, i.e.,

\begin{center}\begin{tikzcd}
	& {(M\otimes M)\otimes M} && {M\otimes (M\otimes M)} \\
	{M\otimes M} &&&& {M\otimes M} \\
	&& M
	\arrow[from=1-2, to=1-4]
	\arrow[from=1-2, to=2-1]
	\arrow[from=1-4, to=2-5]
	\arrow[from=2-1, to=3-3]
	\arrow[from=2-5, to=3-3]
\end{tikzcd}\end{center}
	and

\begin{center}\begin{tikzcd}
	{e\otimes M} & {M\otimes M} & {M\otimes e} \\
	& M
	\arrow[from=1-1, to=1-2]
	\arrow[from=1-1, to=2-2]
	\arrow[from=1-2, to=2-2]
	\arrow[from=1-3, to=1-2]
	\arrow[from=1-3, to=2-2]
\end{tikzcd}\end{center}
	commute.
Morphisms between monoids $M$ and $N$ are given by morphisms $M\to N$ such that

\begin{center}\begin{tikzcd}
	{M\otimes M} & {N\otimes N} \\
	M & N
	\arrow[from=1-1, to=1-2]
	\arrow[from=1-1, to=2-1]
	\arrow[from=1-2, to=2-2]
	\arrow[from=2-1, to=2-2]
\end{tikzcd}\end{center}
and

\begin{center}
	\begin{tikzcd}
	e & M \\
	& N
	\arrow[from=1-1, to=1-2]
	\arrow[from=1-1, to=2-2]
	\arrow[from=1-2, to=2-2]
\end{tikzcd}
\end{center}
commute.
This defines the \idx{category of monoids} in $\mcC$, $\mathrm{Mon}(\mcC,\otimes)$.\footnote{The unit is usually omitted from notation. We do so also for monoids in monoidal categories.}

	If $\mcC$ is braided, a \idx{commutative monoid} in $\mcC$ is a monoid that furthermore fulfills

\begin{center}\begin{tikzcd}
	{M\otimes M} & {M\otimes M} \\
	M
	\arrow["{B_{M,M}}", from=1-1, to=1-2]
	\arrow[from=1-1, to=2-1]
	\arrow[from=1-2, to=2-1]
\end{tikzcd}\end{center}
\end{definition}
One could also define this in terms of a symmetric monoidal functor $(\Fin^{\mathrm{part}},\times)\to (\mcC,\otimes)$.

This is a very useful and strong notion, as can be seen from the following examples.
\begin{example}
	\begin{enumerate}
		\item Monoids in the category of set equipped with product are just the usual monoids,
		\[\mathrm{Mon}(\Set,\times)=\Mon\]
		\item Monoids in the category of abelian groups with the usual tensor product are precisely the unital Rings
		\[\mathrm{Mon}(\Ab,\otimes)=\Ring.\]
		\item Monoids in $\C$-vector spaces are $\C$-algebras,
		 \[\mathrm{Mon}(\Vect_{\C},\otimes)=\mathbf{Alg}_{\C}.\]
\end{enumerate}
\end{example}
Having a notion of monoids, we can infer a canonical definition of \emph{modules over a given monoid}.
\begin{definition}
Consider a monoid $(R,m)$ in $(\mcC,\otimes)$.
Then a \idx{left $R$-module} (or \idx{left $R$-algebra}) in $\mcC$ is an object $V$ with
\idx{scalar multiplication} or \idx{action} $\lambda\colon R\otimes V\to V$
such that $R$ induces an \idx{action} on $V$, meaning that

\begin{center}\begin{tikzcd}
	{R\otimes R\otimes V} & {R\otimes V} \\
	{R\otimes V} & V
	\arrow["{\id\otimes \lambda}", from=1-1, to=1-2]
	\arrow["{m\otimes\id}"', from=1-1, to=2-1]
	\arrow["\lambda", from=1-2, to=2-2]
	\arrow["\lambda"', from=2-1, to=2-2]
\end{tikzcd}\end{center}

and

\begin{center}\begin{tikzcd}
	{e\otimes V} && {R\otimes V} \\
	& V
	\arrow[from=1-1, to=1-3]
	\arrow[from=1-1, to=2-2]
	\arrow[from=1-3, to=2-2]
\end{tikzcd}\end{center}
commute.

Morphisms between two $R$-modules $V$ and $W$ are given by morphisms $V\to W$ such that

\begin{center}\begin{tikzcd}
	{R\otimes V} & {R\otimes W} \\
	V & W
	\arrow[from=1-1, to=1-2]
	\arrow[from=1-1, to=2-1]
	\arrow[from=1-2, to=2-2]
	\arrow[from=2-1, to=2-2]
\end{tikzcd}\end{center}
commutes.

This defines the \idx{category of $R$-modules} in $\mcC$, $\Mod_R(\mcC)$.

A dual definition yields the notions of \idx{right modules over $R$}.

\end{definition}
\begin{remark}\label{rem:restr_of_scalars}
	\begin{enumerate}[(i)]
		\item
			Every homomorphism between monoids $R\to S$ induces a functor $\Mod_S(\mcC)\to \Mod_R(\mcC)$, the so-called \idx{restriction of scalars},
			by endowing an $S$-modules $V$ with the $R$-module structure from
\begin{center}\begin{tikzcd}
	{R\otimes V} & {S\otimes V} & V.
	\arrow[from=1-1, to=1-2]
	\arrow[from=1-2, to=1-3]
\end{tikzcd}\end{center}
			That this indeed induces an $R$-module structure (and that it yields a functor) is straightforward sind $R\to S$ is a monoid homomorphism.
		\item
	Furthermore, if $R$ is a commutative monoid and the symmetric monoidal category $\mcC$ admits reflexive coequalizers and these are preserved coordinatewise by the tensor product,
	then the category $\Mod_R$ has a canonical induced symmetric monoidal structure by defining the tensor product $V\otimes_R W$ of two $R$-modules as the coequalizer of
\begin{center}\begin{tikzcd}
	{V\otimes R\otimes W} & {V\otimes W,}
	\arrow[shift right, from=1-1, to=1-2]
	\arrow[shift left, from=1-1, to=1-2]
\end{tikzcd}\end{center}
	 making restriction of scalars symmetric monoidal.
	 The restriction functor in fact admits a left adjoint, the so-called \idx{base change} along a monoid morphism $R\to S$.
		 Base change can be described as the functor
	\[-\otimes_R S\]
after canonically exhibiting $S$ as an $R$-module.
	  \item If furthermore $\mcC$ admits an internal $\hom$ structure,
			we can procure the \idx{internal $\hom$ of $R$-modules} $\ihom_R(-,-)$ by defining the internal $\hom$ $\ihom_R(V,W)$ as the equalizer of two arrows
\begin{center}\begin{tikzcd}
	{\ihom(V,W)} & {\ihom(R\otimes V, W)}
	\arrow[shift left, from=1-1, to=1-2]
	\arrow[shift right, from=1-1, to=1-2]
\end{tikzcd}\end{center}
	where the first arrow is given by pullback of the action of $R$ on $V$ and the second one by pushforward of the action on $W$ (see \ref{prop:G-mod-ihom} for more on this).
		The functor $\ihom_R(S,-)$ forms a right adjoint to restriction of scalars along $R\to S$.
\end{enumerate}
\end{remark}
\begin{example}
	\begin{enumerate}
\item 	Every monoid homomorphism $R\to S$ installs an $R$-module structure on $S$.
\item Clearly, the usual modules over a ring are a special case of this definition.
\item Sets equipped with an action by a monoid $H$ form another example.
\item In the next chapter, we will see many more examples of these categories.
	\end{enumerate}
\end{example}

\subsection{Some diagram lemmas}
We collect some elementary lemmas of homological algebra to get a feeling for typical \enquote{diagram chasing} statements.
Note that most of the results in this section can be proven immediately (but much more abstractly) using spectral sequences.

\begin{definition}[Exact sequences]\label{def:chain_complex}\uses{def:category, def:ker_im, def:abelian, thm:properties_of_abelian_categories}
	Let $\mcA$ be an abelian category.
A \idx{sequence} $C^\bullet=((C^i)_{i\in\Z},(\delta_i)_{i\in\Z})$ is a diagram in $\mcA$ whose form is the totally ordered set $\Z$.
Explicitly, $(C^\bullet, \delta^\bullet)$ looks like

\begin{center}\begin{tikzcd}
	\dots & {C^{-2}} & {C^{-1}} & {C^0} & {C^1} & {C^2} & \dots
	\arrow[from=1-1, to=1-2]
	\arrow["{\delta^{-2}}", from=1-2, to=1-3]
	\arrow["{\delta^{-1}}", from=1-3, to=1-4]
	\arrow["{\delta^0}", from=1-4, to=1-5]
	\arrow["{\delta^1}", from=1-5, to=1-6]
	\arrow[from=1-6, to=1-7]
\end{tikzcd}\end{center}

	A sequence is called \textbf{(cochain) complex}\index{complex}\index{cochain complex}, if $\delta^i\circ \delta^{i-1}=0$ for all $i$, or equivalently, if $\Im(\delta^{i-1})\subseteq\ker(\delta^i)$.

For any cochain complex $C^\bullet$ define
\begin{itemize}
	\item the $i$-\idx{cochains} as $C^i$
	\item the $i$-\idx{cocycles} by $Z^i=\ker\delta^i\subseteq C^i$,
	\item the $i$-\idx{coboundaries} by $B^{i}=\Im(\delta^{i-1})$
  \item the \idx{cohomology groups} $H^i(C^\bullet)$ of any complex $C^\bullet$ is the quotient $\ker(\delta^i)/\Im(\delta^{i-1})$ induced as the cokernel of the inclusion $\Im(\delta_{i-1})\hookrightarrow \ker(\delta_i)$.
		\footnote{Note that, of course, $H^{i}(C^{\bullet})$ are objects of $\mcA$ and these often have nothing to do with groups.
		Nonetheless,
		this is totally standard terminology (to which one quickly gets used).}
\end{itemize}
	A complex is called \idx{exact sequence} or \idx{acyclic}, if $H^i(C^\bullet)=0$ for all $i$.

We say that the complex is \idx{bounded to the left}, if there exists $k\in\Z$ such that $C^i=0$ for all $i\le k$, and \idx{positively graded}, if in addition $k$ can be taken to be $\ge0$.
We say that the complex is \idx{bounded to the right}, if there exists $k\in\Z$ such that $C^i=0$ for all $i\ge k$, and \idx{negatively graded}, if $k\le 0$.
	It is \idx{bounded}, if it is both bounded to the left and the right.
In these cases, one often omits $0$'s in the notation and, e.g., simply writes

\begin{center}\begin{tikzcd}
	0 & {C^0} & {C^1} & {C^2} & \dots
	\arrow[from=1-1, to=1-2]
	\arrow[from=1-2, to=1-3]
	\arrow[from=1-3, to=1-4]
	\arrow[from=1-4, to=1-5]
\end{tikzcd}\end{center}
	A complex is said to be \idx{short}, if it is of the form

\begin{center}\begin{tikzcd}
	0 & A & B & C & 0.
	\arrow[from=1-1, to=1-2]
	\arrow[from=1-2, to=1-3]
	\arrow[from=1-3, to=1-4]
	\arrow[from=1-4, to=1-5]
\end{tikzcd}\end{center}
	The important special case of this complex being exact is called a \idx{short exact sequence}.
\end{definition}
\begin{remark}
	We are using cohomological grading as it seems natural to us to use the usual orientation of the integers $\Z$.
	However, reversing all arrows, one obtains definitions for chain complexes, boundaries, cycles and homologies.
	We will often freely switch indexing conventions; as, e.g., when working with projective resolutions, or more generally complexes which are bounded to the left,
	it seems more natural to use positive indices to the left than using many negative numbers.
	We hope that this helps more in having a focus on important ideas than creates confusion.
Usually, homological indexing is indicated by subscripts and cohomological notation by superscripts.
We try our best to do this consequently.
\end{remark}

\begin{definition}[Chain morphisms and category]\label{def:complex}\uses{def:chain_complex}
	A morphism $f\colon A^\bull\to B^\bull$ between two complexes is a natural transformations of the corresponding sequences,

\begin{center}\begin{tikzcd}
	  \dots & {A^{i-1}} & {A^i} & {A^{i+1}} & \dots \\
	\dots & {B^{i-1}} & {B^i} & {B^{i+1}} & \dots
	\arrow[from=1-1, to=1-2]
	\arrow[from=1-2, to=1-3]
	\arrow["{f^{i-1}}", from=1-2, to=2-2]
	\arrow[from=1-3, to=1-4]
	\arrow["{f^i}", from=1-3, to=2-3]
	\arrow[from=1-4, to=1-5]
	\arrow["{f^{i+1}}", from=1-4, to=2-4]
	\arrow[from=2-1, to=2-2]
	\arrow[from=2-2, to=2-3]
	\arrow[from=2-3, to=2-4]
	\arrow[from=2-4, to=2-5]
\end{tikzcd}\end{center}
	This defines the category of chain complexes $\Ch(\mcA)$ and this is an abelian category.
For any full subcategory $\mcB\sub \mcA$ define $\Ch(\mcB)$ to be the full subcategory of complexes in $\mcA$ consisting of complexes whose objects lie in $\mcB$.
The left bounded/right bounded/positively graded/negatively graded/bounded complex category $\Ch^+(\mcA)$/$\Ch^-(\mcA)$/$\Ch^{\ge 0}$/$\Ch^{\le 0}$/$\Ch^b(\mcA)$
are defined as the full subcategories spanned by left bounded/right bounded/etc. complexes.
\end{definition}

\begin{lemma}
$\Ch(\mcA)$ is abelian and all limits and colimits may be computed pointwise (if they exist pointwise in $\mcA$).

If $\mcA$ admits enough (compact) projectives, so does $\Ch(\mcA)$, where a generating class is, e.g., given by bounded complexes of (compact) projectives.
\end{lemma}
\begin{proof}
See, e.g., \cite[chap.~10 and 11]{Murfet2006a} for this and more.
\end{proof}
\begin{remark}
	Any natural transformation between two complexes induces morphisms between the corresponding cohomologies.
	Checking this is an excellent exercise for working with universal properties.
	We give the details for the reader's convenience.
	(Usually, such arguments are done using either Freyd-Mitchell directly or so-called elements.)

	In a first step note that $f_i$ maps the kernel of $A_i\to A_{i+1}$ to the kernel of $B_{i}\to B_{i+1}$, which can be seen by

\begin{center}
	\begin{tikzcd}
	&&& {\ker \delta^i_A} \\
	\dots & \dots && {A^i} && {A^{i+1}} & \dots \\
	&& {\ker(\delta^i_B)} \\
	\dots & \dots && {B^i} && {B^{i+1}} & \dots
	\arrow[from=1-4, to=2-4]
	\arrow["0", from=1-4, to=2-6]
	\arrow[dashed, from=1-4, to=3-3]
	\arrow[from=2-1, to=2-2]
	\arrow[from=2-2, to=2-4]
	\arrow[from=2-2, to=4-2]
	\arrow["{\delta^i_A}", from=2-4, to=2-6]
	\arrow["{f^i}", from=2-4, to=4-4]
	\arrow[from=2-6, to=2-7]
	\arrow["{f^{i+1}}", from=2-6, to=4-6]
	\arrow[from=3-3, to=4-4]
	\arrow[from=4-1, to=4-2]
	\arrow[from=4-2, to=4-4]
	\arrow["{\delta^i_B}", from=4-4, to=4-6]
	\arrow[from=4-6, to=4-7]
\end{tikzcd}\end{center}
	since the composition $\ker \delta^i_A\to A^i\to B^i$ factors over $\ker(\delta ^i_B)$ since composed with $\delta_B^i$ we can \enquote{walk} over $A^{i+1}$ and hence see that the arrow is $0$.

	Hence $f^i$ restricts to a morphism $\ker\delta_A^i\to \ker\delta_B^i$.
	To obtain from this a morphism of cohomologies, we need to see that $f_i$ maps the image of $\delta_A^{i-1}$ into the image of $\delta_B^{i-1}$.
	For this, consider the image factorisation

\begin{center}\begin{tikzcd}
	&& {\Im(\delta_A^{i-1})} \\
	\dots & {A^{i-1}} && {\ker \delta^i_A} \\
	&& {\Im(\delta_B^{i-1})} \\
	\dots & {B^{i-1}} && {\ker\delta^i_B}
	\arrow[from=1-3, to=2-4]
	\arrow[dashed, from=1-3, to=3-3]
	\arrow[from=2-1, to=2-2]
	\arrow[from=2-2, to=1-3]
	\arrow["{\delta_A^{i-1}}", from=2-2, to=2-4]
	\arrow[from=2-2, to=4-2]
	\arrow["{f^i}", from=2-4, to=4-4]
	\arrow[from=3-3, to=4-4]
	\arrow[from=4-1, to=4-2]
	\arrow[from=4-2, to=3-3]
	\arrow["{\delta_B^{i-1}}", from=4-2, to=4-4]
\end{tikzcd}\end{center}
	The desired dashed arrow can be obtained by applying lemma \ref{lem:square_with_diag} to the square induced by the two dashed arrows

\begin{center}\begin{tikzcd}
	&& {\Im(\delta_A^{i-1})} \\
	\dots & {A^{i-1}} && {\ker \delta^i_A} \\
	&& {\Im(\delta_B^{i-1})} \\
	\dots & {B^{i-1}} && {\ker\delta^i_B}
	\arrow[from=1-3, to=2-4]
	\arrow["{h_2}"', dashed, from=1-3, to=4-4]
	\arrow[from=2-1, to=2-2]
	\arrow[two heads, from=2-2, to=1-3]
	\arrow["{\delta_A^{i-1}}", from=2-2, to=2-4]
	\arrow["{h_1}"', dashed, from=2-2, to=3-3]
	\arrow[from=2-2, to=4-2]
	\arrow["{f^i}", from=2-4, to=4-4]
	\arrow[hook, from=3-3, to=4-4]
	\arrow[from=4-1, to=4-2]
	\arrow[from=4-2, to=3-3]
	\arrow["{\delta_B^{i-1}}", from=4-2, to=4-4]
\end{tikzcd}\end{center}

	Having this, we see that it induces a morphism of the factors (which are the cohomologies by definition),

\begin{center}\begin{tikzcd}
	&& {\Im(\delta_A^{i-1})} \\
	\dots & {A^{i-1}} && {\ker \delta^i_A} \\
	&& {\Im(\delta_B^{i-1})} && {\ker\delta_A^i/\im(\delta_A^{i-1})} \\
	\dots & {B^{i-1}} && {\ker\delta^i_B} \\
	&&&& {\ker(\delta_B^{i})/\im(\delta_B^{i-1})}
	\arrow[from=1-3, to=2-4]
	\arrow[from=1-3, to=3-3]
	\arrow[from=2-1, to=2-2]
	\arrow[from=2-2, to=2-4]
	\arrow[from=2-2, to=4-2]
	\arrow[from=2-4, to=3-5]
	\arrow["{f^i}", from=2-4, to=4-4]
	\arrow[from=3-3, to=4-4]
	\arrow[dashed, from=3-5, to=5-5]
	\arrow[from=4-1, to=4-2]
	\arrow[from=4-2, to=4-4]
	\arrow[from=4-4, to=5-5]
\end{tikzcd}\end{center}
	Here we use the universal property of the upper factor (cokernel) $\ker^i_A/\im(\delta^{i-1}_A)$, by seeing that the morphism $\im(\delta_A^{i-1})\to H^i(B^\bullet)$ is zero, using the factorisation.

\end{remark}

Any functor that preserves the $0$ object maps chain complexes to chain complexes,
in particular any, additive functor thus extends to a functor between the complex categories.
However, it does not need to preserve exactness.
\begin{lemma}[Characterisation of left/right exact functors]\label{lem:characterisation_of_exact_functors}\uses{def:exact_functor, def:complex}
A functor $F$ between abelian categories is right exact (in the sense of being finitely cocontinuous, \ref{def:exact_functor}), precisely if it maps
	short exact sequences

\begin{center}\begin{tikzcd}
	0 & A & B & C & 0
	\arrow[from=1-1, to=1-2]
	\arrow[from=1-2, to=1-3]
	\arrow[from=1-3, to=1-4]
	\arrow[from=1-4, to=1-5]
\end{tikzcd}\end{center}
	to exact sequences,

\begin{center}\begin{tikzcd}
	{F(A)} & {F(B)} & {F(C)} & 0
	\arrow[from=1-1, to=1-2]
	\arrow[from=1-2, to=1-3]
	\arrow[from=1-3, to=1-4]
\end{tikzcd}\end{center}
	meaning that the image of the chain complex under $F$,

\begin{center}\begin{tikzcd}
	0 & {F(A)} & {F(B)} & {F(C)} & 0
	\arrow[from=1-1, to=1-2]
	\arrow[from=1-2, to=1-3]
	\arrow[from=1-3, to=1-4]
	\arrow[from=1-4, to=1-5]
\end{tikzcd}\end{center}
is a chain complex which does not have to be exact at the place $F(A)$ (i.e., $F(A)\to F(B)$ is not necessarily monic).

	Similarly, $F$ is left exact precisely if for any short exact sequence, the left part

\begin{center}\begin{tikzcd}
	0 & {F(A)} & {F(B)} & {F(C)}
	\arrow[from=1-1, to=1-2]
	\arrow[from=1-2, to=1-3]
	\arrow[from=1-3, to=1-4]
\end{tikzcd}\end{center}
	of the image chain remains exact.

	Combining these properties, we see that a functor is exact if and only if it maps exact sequences to exact sequences.
	More generally, a functor is exact precisely if it preserves all homology groups.
\end{lemma}
\begin{proof}
	See, e.g., \cite[1.11.2]{Borceux2008a} or lemma 12.7.2 in \cite[010N]{SPA2023}.
\end{proof}

\begin{remark}
  Note the following easy facts for exact sequences:
  \begin{itemize}
		\item \enquote{$0$ at left means monic and $0$ at right means epic}:
		In any exact sequence

\begin{center}\begin{tikzcd}
	\dots & A & B & C & D & \dots
	\arrow[from=1-1, to=1-2]
	\arrow[from=1-2, to=1-3]
	\arrow[from=1-3, to=1-4]
	\arrow[from=1-4, to=1-5]
	\arrow[from=1-5, to=1-6]
\end{tikzcd}\end{center}
		$(A\to B)=0$ if and only if $B\hookrightarrow C$, and $(C\to D)=0$ precisely if $B\twoheadrightarrow C$.
		\item For any exact sequence with

\begin{center}\begin{tikzcd}
	\dots & 0 & A & 0 & \dots
	\arrow[from=1-1, to=1-2]
	\arrow[from=1-2, to=1-3]
	\arrow[from=1-3, to=1-4]
	\arrow[from=1-4, to=1-5]
\end{tikzcd}\end{center}
We conclude $A=0$.
		\item In any exact sequence

\begin{center}\begin{tikzcd}
	\dots & A & B & C & D & \dots,
	\arrow[from=1-1, to=1-2]
	\arrow[from=1-2, to=1-3]
	\arrow[from=1-3, to=1-4]
	\arrow[from=1-4, to=1-5]
	\arrow[from=1-5, to=1-6]
\end{tikzcd}\end{center}
		0ne has that $B\to C$ is an isomorphism if and only if $(A\to B)=(C\to D)=0$.
\end{itemize}
	\end{remark}
\begin{example}
A first prominent example of a short exact sequence is any direct sum

\begin{center}\begin{tikzcd}
	0 & A & {A\oplus B} & B & 0
	\arrow[from=1-1, to=1-2]
	\arrow["{i_A}", from=1-2, to=1-3]
	\arrow["{\pi_B}", from=1-3, to=1-4]
	\arrow[from=1-4, to=1-5]
\end{tikzcd}.\end{center}
	Such an exact sequence (and, of course, any sequence isomorphic to one of these) is called \idx{split exact sequence}.

	Another instructive example of an exact sequence of abelian groups is

\begin{center}\begin{tikzcd}
	0 & {2\Z} & \Z & {\Z/2\Z} & 0
	\arrow[from=1-1, to=1-2]
	\arrow[hook, from=1-2, to=1-3]
	\arrow[two heads, from=1-3, to=1-4]
	\arrow[from=1-4, to=1-5]
\end{tikzcd}\end{center}

\end{example}

\subsection{Elementary homological lemmas}
Next, we will collect some elementary diagram chase lemmas.

First, we will discuss how to recognise those short exact sequences which are split.
One important tool for this is the following 5-lemma.
\begin{lemma}[5-lemma]\label{lem:five_lemma}\uses{def:complex}
	Assume given two exact rows (we omit the irrelevant continuations to left and right)
	with a natural transformation (a morphism between the sequences) between them

\begin{center}\begin{tikzcd}
	A & B & C & D & E \\
	{A'} & {B'} & {C'} & {D'} & {E'}
	\arrow[from=1-1, to=1-2]
	\arrow["f", from=1-1, to=2-1]
	\arrow[from=1-2, to=1-3]
	\arrow["g", from=1-2, to=2-2]
	\arrow[from=1-3, to=1-4]
	\arrow["h", from=1-3, to=2-3]
	\arrow[from=1-4, to=1-5]
	\arrow["k", from=1-4, to=2-4]
	\arrow["\ell", from=1-5, to=2-5]
	\arrow[from=2-1, to=2-2]
	\arrow[from=2-2, to=2-3]
	\arrow[from=2-3, to=2-4]
	\arrow[from=2-4, to=2-5]
\end{tikzcd}.\end{center}
	Then we have the following \idx{5-lemma}:

	If $g$ and $k$ are isomorphisms, $f$ is epic and $k$ is monic, then $h$ is an isomorphism,

\begin{center}\begin{tikzcd}
	A & B & C & D & E \\
	{A'} & {B'} & {C'} & {D'} & {E'}
	\arrow[from=1-1, to=1-2]
	\arrow[two heads, from=1-1, to=2-1]
	\arrow[from=1-2, to=1-3]
	\arrow["\simeq", from=1-2, to=2-2]
	\arrow[from=1-3, to=1-4]
	\arrow["h", from=1-3, to=2-3]
	\arrow[from=1-4, to=1-5]
	\arrow["\simeq", from=1-4, to=2-4]
	\arrow[hook, from=1-5, to=2-5]
	\arrow[from=2-1, to=2-2]
	\arrow[from=2-2, to=2-3]
	\arrow[from=2-3, to=2-4]
	\arrow[from=2-4, to=2-5]
\end{tikzcd}.\end{center}
	This statement decomposes into two \idx{4-lemmas}:

\begin{center}\begin{tikzcd}
	A & B & C & D \\
	{A'} & {B'} & {C'} & {D'}
	\arrow[from=1-1, to=1-2]
	\arrow[two heads, from=1-1, to=2-1]
	\arrow[from=1-2, to=1-3]
	\arrow[hook, from=1-2, to=2-2]
	\arrow[from=1-3, to=1-4]
	\arrow["h", from=1-3, to=2-3]
	\arrow[hook, from=1-4, to=2-4]
	\arrow[from=2-1, to=2-2]
	\arrow[from=2-2, to=2-3]
	\arrow[from=2-3, to=2-4]
\end{tikzcd}\end{center}
implies $h$ monic, and

\begin{center}\begin{tikzcd}
	B & C & D & E \\
	{B'} & {C'} & {D'} & {E'}
	\arrow[from=1-1, to=1-2]
	\arrow[two heads, from=1-1, to=2-1]
	\arrow[from=1-2, to=1-3]
	\arrow["h", from=1-2, to=2-2]
	\arrow[from=1-3, to=1-4]
	\arrow[two heads, from=1-3, to=2-3]
	\arrow[hook, from=1-4, to=2-4]
	\arrow[from=2-1, to=2-2]
	\arrow[from=2-2, to=2-3]
	\arrow[from=2-3, to=2-4]
\end{tikzcd}\end{center}
implies $h$ epic.
\end{lemma}
This easily implies the following \idx{splitting lemma}:
\begin{lemma}[Splitting lemma]\label{lem:splitting_lemma}\uses{lem:five_lemma}
	Consider a short exact sequence

\begin{center}\begin{tikzcd}
	0 & A & B & C & 0
	\arrow[from=1-1, to=1-2]
	\arrow["m", from=1-2, to=1-3]
	\arrow["p", from=1-3, to=1-4]
	\arrow[from=1-4, to=1-5]
\end{tikzcd}.\end{center}
	Then the following assertions are equivalent (explaining the term \enquote{split exact sequence})
	\begin{enumerate}
		\item[a)] The sequence is split exact, i.e. there is an isomorphism $h\colon B\to A\oplus C$ with

\begin{center}\begin{tikzcd}
	&& {A\oplus C} \\
	0 & A & B & C & 0
	\arrow["{\pi_C}", from=1-3, to=2-4]
	\arrow[from=2-1, to=2-2]
	\arrow["{i_A}", from=2-2, to=1-3]
	\arrow[from=2-2, to=2-3]
	\arrow["h"', from=2-3, to=1-3]
	\arrow[from=2-3, to=2-4]
	\arrow[from=2-4, to=2-5]
\end{tikzcd}\end{center}

		\item[b)] The sequence is \idx{right split}, meaning that $p$ is a split epimorphism,

\begin{center}\begin{tikzcd}
	0 & A & B & C & 0
	\arrow[from=1-1, to=1-2]
	\arrow[from=1-2, to=1-3]
	\arrow["p", shift left, from=1-3, to=1-4]
	\arrow["s", shift left, dashed, from=1-4, to=1-3]
	\arrow[from=1-4, to=1-5]
\end{tikzcd}\end{center}
		\item[c)] The sequence is \idx{left split}, meaning that $m$ is a split monomorphism,

\begin{center}\begin{tikzcd}
	0 & A & B & C & 0
	\arrow[from=1-1, to=1-2]
	\arrow["m", shift left, from=1-2, to=1-3]
	\arrow["t", shift left, dashed, from=1-3, to=1-2]
	\arrow[from=1-3, to=1-4]
	\arrow[from=1-4, to=1-5]
\end{tikzcd}\end{center}
	\end{enumerate}
	\end{lemma}

Having understood split exact sequences, next we combine short exact sequences with chain complexes.
\begin{lemma}[\idx{Snake lemma}]\label{lem:snake_lemma}\uses{def:complex}
	For any diagram
\begin{center}
\begin{tikzcd}
	& A & B & C & 0 \\
	0 & {A'} & {B'} & {C'}
	\arrow["f", from=1-2, to=1-3]
	\arrow["a"', from=1-2, to=2-2]
	\arrow["g", from=1-3, to=1-4]
	\arrow["b"', from=1-3, to=2-3]
	\arrow[from=1-4, to=1-5]
	\arrow["c"', from=1-4, to=2-4]
	\arrow[from=2-1, to=2-2]
	\arrow["{f'}"', from=2-2, to=2-3]
	\arrow["{g'}"', from=2-3, to=2-4]
\end{tikzcd}\end{center}
with exact rows, we obtain an exact sequence
\begin{center}
  \begin{tikzcd}
	\ker a\ar[r] & \ker b \ar[r] & \ker c\ar[r,"\delta"] & \coker a\ar[r] & \coker b\ar[r] & \coker c.
  \end{tikzcd}
\end{center}
\end{lemma}

The \idx{connecting homomorphism} $\delta$ is the reason for the snake lemma's name, as it moves like a snake from top right to bottom left in

\begin{center}\begin{tikzcd}
	& {\ker a} & {\ker b} & {\ker c} \\
	& A & B & C & 0 \\
	0 & {A'} & {B'} & {C'} \\
	& {\coker a} & {\coker b} & {\coker c}
	\arrow[from=1-2, to=1-3]
	\arrow[from=1-2, to=2-2]
	\arrow[from=1-3, to=1-4]
	\arrow[from=1-3, to=2-3]
	\arrow[from=1-4, to=2-4]
	\arrow[from=1-4, to=4-2]
	\arrow["f", from=2-2, to=2-3]
	\arrow["a"', from=2-2, to=3-2]
	\arrow["g", from=2-3, to=2-4]
	\arrow["b"', from=2-3, to=3-3]
	\arrow[from=2-4, to=2-5]
	\arrow["c"', from=2-4, to=3-4]
	\arrow[from=3-1, to=3-2]
	\arrow["{f'}"', from=3-2, to=3-3]
	\arrow[from=3-2, to=4-2]
	\arrow["{g'}"', from=3-3, to=3-4]
	\arrow[from=3-3, to=4-3]
	\arrow[from=3-4, to=4-4]
	\arrow[from=4-2, to=4-3]
	\arrow[from=4-3, to=4-4]
\end{tikzcd}.\end{center}
\begin{remark}
	Due to our TikZ-skills we should probably call it \enquote{diagonal lemma}.
\end{remark}

This is a special case of the very important principle \enquote{short exact sequences of complexes give a long exact sequence of homologies}, known as the \idx{zig-zag lemma} or often also called snake lemma.
\begin{lemma}[Zig-Zag lemma]\label{lem:zig_zag_lemma}\uses{def:complex}
	Consider any short exact sequence $0\to A^\bullet\to B^\bullet\to C^\bullet\to 0$ in $\Ch(\mcA)$, i.e.,

\begin{center}\begin{tikzcd}
	& \vdots & \vdots & \vdots \\
	0 & {A^i} & {B^i} & {C^i} & 0 \\
	0 & {A^{i+1}} & {B^{i+1}} & {C^{i+1}} & 0 \\
	0 & {A^{i+2}} & {B^{i+2}} & {C^{i+2}} & 0 \\
	& \vdots & \vdots & \vdots
	\arrow[from=1-2, to=2-2]
	\arrow[from=1-3, to=2-3]
	\arrow[from=1-4, to=2-4]
	\arrow[from=2-1, to=2-2]
	\arrow[from=2-2, to=2-3]
	\arrow[from=2-2, to=3-2]
	\arrow[from=2-3, to=2-4]
	\arrow[from=2-3, to=3-3]
	\arrow[from=2-4, to=2-5]
	\arrow[from=2-4, to=3-4]
	\arrow[from=3-1, to=3-2]
	\arrow[from=3-2, to=3-3]
	\arrow[from=3-2, to=4-2]
	\arrow[from=3-3, to=3-4]
	\arrow[from=3-3, to=4-3]
	\arrow[from=3-4, to=3-5]
	\arrow[from=3-4, to=4-4]
	\arrow[from=4-1, to=4-2]
	\arrow[from=4-2, to=4-3]
	\arrow[from=4-2, to=5-2]
	\arrow[from=4-3, to=4-4]
	\arrow[from=4-3, to=5-3]
	\arrow[from=4-4, to=4-5]
	\arrow[from=4-4, to=5-4]
\end{tikzcd}\end{center}
	with exact rows and complexes as columns.

	Then there exists a long exact sequence of cohomologies

\begin{center}\begin{tikzcd}
	& \dots & \dots \\
	{H^{i-1}(A^\bull)} & {H^{i-1}(B^\bull)} & {H^{i-1}(C^\bull)} \\
	{H^i(A^\bull)} & {H^i(B^\bull)} & {H^i(C^\bull)} \\
	{H^{i+1}(A^\bull)} & {H^{i+1}(B^\bull)} & {H^{i+1}(C^\bull)} \\
	\dots & \dots
	\arrow[from=1-2, to=1-3]
	\arrow[from=1-3, to=2-1]
	\arrow[from=2-1, to=2-2]
	\arrow[from=2-2, to=2-3]
	\arrow[from=2-3, to=3-1]
	\arrow[from=3-1, to=3-2]
	\arrow[from=3-2, to=3-3]
	\arrow[from=3-3, to=4-1]
	\arrow[from=4-1, to=4-2]
	\arrow[from=4-2, to=4-3]
	\arrow[from=4-3, to=5-1]
	\arrow[from=5-1, to=5-2]
\end{tikzcd}\end{center}
where the morphisms $H^{i}(C^\bull)\to H^{i+1}(A^\bull)$ are called \idx{connecting homomorphism}.
	\end{lemma}

	To see that this implies the snake lemma,
	consider the following short exact sequence of chain complexes.

\begin{center}\begin{tikzcd}
	& \vdots & \vdots & \vdots \\
	& 0 & 0 & 0 \\
	0 & {\ker g} & B & C & 0 \\
	0 & {A'} & {B'} & {\coker f'} & 0 \\
	& 0 & 0 & 0 \\
	& \vdots & \vdots & \vdots
	\arrow[from=1-2, to=2-2]
	\arrow[from=1-3, to=2-3]
	\arrow[from=1-4, to=2-4]
	\arrow[from=2-2, to=2-3]
	\arrow[from=2-2, to=3-2]
	\arrow[from=2-3, to=2-4]
	\arrow[from=2-3, to=3-3]
	\arrow[from=2-4, to=3-4]
	\arrow[from=3-1, to=3-2]
	\arrow["f", from=3-2, to=3-3]
	\arrow["a"', from=3-2, to=4-2]
	\arrow["g", from=3-3, to=3-4]
	\arrow["b"', from=3-3, to=4-3]
	\arrow[from=3-4, to=3-5]
	\arrow["c"', from=3-4, to=4-4]
	\arrow[from=4-1, to=4-2]
	\arrow["{f'}"', from=4-2, to=4-3]
	\arrow[from=4-2, to=5-2]
	\arrow["{g'}"', from=4-3, to=4-4]
	\arrow[from=4-3, to=5-3]
	\arrow[from=4-4, to=4-5]
	\arrow[from=4-4, to=5-4]
	\arrow[from=5-2, to=5-3]
	\arrow[from=5-2, to=6-2]
	\arrow[from=5-3, to=5-4]
	\arrow[from=5-3, to=6-3]
	\arrow[from=5-4, to=6-4]
\end{tikzcd}\end{center}

Further homological lemmas include the \idx{salamander lemma} with special case being the \idx{$3\times 3=9$-lemma}.
\subsection{Spectral sequences}
\quot{I am a bit panic-stricken by this flood of cohomology, but have borne up courageously.
Your spectral sequence seems reasonable to me.
	(I thought I had shown that it was wrong in a special case, but I was mistaken, on the contrary it works remarkably well.)}
{J.-P. Serre to A. Grothendieck, \cite[p.~38]{Serre2001}, translation taken from \cite{Vakil2023}.}

Having seen many different homological \enquote{diagram chase} lemmata with bigger and bigger diagrams,
we present a general and powerful machinery behind this kind of result.
See \cite{Vakil2023,Vakila, Vakil, Gelfand2003, Weibel1994} for more.
\begin{definition}[Double complex]\label{def:double_complex}\uses{def:complex}
A \idx{double complex} in an abelian category $\mcA$ is a diagram

\begin{center}\begin{tikzcd}
	& \vdots & \vdots & \vdots \\
	\dots & {A_{j+1,i+1}} & {A_{j,i+1}} & {A_{j-1,i+1}} & \dots \\
	\dots & {A_{j+1,i}} & {A_{j,i}} & {A_{j-1,i}} & \dots \\
	\dots & {A_{j+1,i-1}} & {A_{j,i-1}} & {A_{j-1,i-1}} & \dots \\
	& \vdots & \vdots & \vdots
	\arrow[from=1-2, to=2-2]
	\arrow[from=1-3, to=2-3]
	\arrow[from=1-4, to=2-4]
	\arrow[from=2-1, to=2-2]
	\arrow[from=2-2, to=2-3]
	\arrow[from=2-2, to=3-2]
	\arrow["{r_{j,i+1}}", from=2-3, to=2-4]
	\arrow["{d_{j,i+1}}", from=2-3, to=3-3]
	\arrow[from=2-4, to=2-5]
	\arrow["{d_{j-1,i+1}}", from=2-4, to=3-4]
	\arrow[from=3-1, to=3-2]
	\arrow["{r_{j+1,i}}", from=3-2, to=3-3]
	\arrow[from=3-2, to=4-2]
	\arrow["{r_{j,i}}", from=3-3, to=3-4]
	\arrow["{d_i^{j}}", from=3-3, to=4-3]
	\arrow[from=3-4, to=3-5]
	\arrow[from=3-4, to=4-4]
	\arrow[from=4-1, to=4-2]
	\arrow[from=4-2, to=4-3]
	\arrow[from=4-2, to=5-2]
	\arrow[from=4-3, to=4-4]
	\arrow[from=4-3, to=5-3]
	\arrow[from=4-4, to=4-5]
	\arrow[from=4-4, to=5-4]
\end{tikzcd}\end{center}
	where every small square \idx{anticommutes}, i.e.,
	\[d_{j-1,i+1}r_{j,i+1}=-r_{j,i}d_{j,i+1}.\footnote{The convention to take anticommuting diagrams rather than commuting ones is from \cite{Vakil2023}
		and should help in avoiding many $(-1)^j$ factors and alternating sums.
	We also tried to use mixed indexing conventions $A_i^j$ but this turned out to be way more confusing.}\]

The (direct sum) \idx{total complex} $\tot(A_{\bullet,\bullet})$ (or $\tot^{\oplus}(A_{\bullet,\bullet})$)
of a double complex $A_{\bullet,\bullet}$ (assuming small cocompleteness, i.e., AB3, if $A_{\bullet,\bullet}$ is unbounded in the relevant direction, and writing infinite coproducts as infinite direct sums) is the complex

\begin{center}\begin{tikzcd}
	\dots & {\bigoplus_{j+i=k}A_{j,i}} & {\bigoplus_{j+i=k-1}A_{j,i} }& \dots
	\arrow[from=1-1, to=1-2]
	\arrow["{t_k}", from=1-2, to=1-3]
	\arrow[from=1-3, to=1-4]
\end{tikzcd}\end{center}
where the $t_k$ is the direct sum of all $d_{j,i}+r_{j,i}$ with $j+i=k$.
\end{definition}
\begin{example}[Spectral sequence of a double complex]\label{def:spec_seq_of_double}
  We give a step-by-step definition of the spectral sequence of a double complex.
  This is intended to provide some intuition for the more general (but heavily inspired by this special case)
  definition of a spectral sequence.

Consider a double complex $(A_{\bullet,\bullet})$.
	Define the \idx{spectral sequence with downward orientation} associated to $(A_{\bullet,\bullet})$ to be the iterative sequence $(E_{\bullet,\bullet}^r)_{r\in\N_0}$ obtained in the following way:
	\begin{itemize}
		\item Define the $0$'th page as the sequence of downward chains,

\begin{center}\begin{tikzcd}
	\vdots & \vdots & \vdots & \vdots & \vdots \\
	\dots & {E^0_{j+1,i+1}=A_{j+1,i+1}} & {E^{0}_{j,i+1}=A_{j,i+1}} & {E^0_{j-1,i+1}=A_{j-1, i+1}} & \dots \\
	\dots & {E_{j+1,i}^0=A_{j+1,i}} & {E^0_{j,i}=A_{j,i}} & {E^0_{j-1,i}=A_{j-1,i}} & \dots \\
	\dots & {E^{0}_{j+1,i-1}=A_{j+1,i-1}} & {E^0_{j,i-1}=A_{j,i-1}} & {E^{0}_{j-1,i-1}=A_{j-1,i-1}} & \dots \\
	\vdots & \vdots & \vdots & \vdots & \vdots
	\arrow[from=1-1, to=2-1]
	\arrow[from=1-2, to=2-2]
	\arrow[from=1-3, to=2-3]
	\arrow[from=1-4, to=2-4]
	\arrow[from=1-5, to=2-5]
	\arrow[from=2-1, to=3-1]
	\arrow[from=2-2, to=3-2]
	\arrow[from=2-3, to=3-3]
	\arrow[from=2-4, to=3-4]
	\arrow[from=2-5, to=3-5]
	\arrow[from=3-1, to=4-1]
	\arrow[from=3-2, to=4-2]
	\arrow[from=3-3, to=4-3]
	\arrow[from=3-4, to=4-4]
	\arrow[from=3-5, to=4-5]
	\arrow[from=4-1, to=5-1]
	\arrow[from=4-2, to=5-2]
	\arrow[from=4-3, to=5-3]
	\arrow[from=4-4, to=5-4]
	\arrow[from=4-5, to=5-5]
\end{tikzcd}\end{center}
		by simply forgetting the existence of the rightwards morphisms.
		\item Define the next page, $E^1_{\bullet, \bullet}$ by taking homologies of the chain complexes in $E^0$,
		and then observing that the rightward morphisms induce morphisms of the homologies as follows:

\begin{center}
\begin{tikzcd}
	\dots & {E^1_{j+1,i+1}=H_{j+1,i+1}(E^0_{\bullet,\bullet})} & {E^1_{j,i+1}=H_{j,i+1}(E^0_{\bullet,\bullet})} & {E^1_{j-1,i+1}} & \dots \\
	\dots & {E^1_{j+1,i}=H_{j+1,i}(E^0_{\bullet,\bullet})} & {E^1_{j,i}=H_{j,i}(E^0_{\bullet,\bullet})} & {E^1_{j-1,i}} & \dots \\
	\dots & {E^1_{j+1,i-1}} & {E^1_{j,i-1}} & {E^1_{j-1,i-1}} & \dots
	\arrow[from=1-1, to=1-2]
	\arrow[from=1-2, to=1-3]
	\arrow[from=1-3, to=1-4]
	\arrow[from=1-4, to=1-5]
	\arrow[from=2-1, to=2-2]
	\arrow[from=2-2, to=2-3]
	\arrow[from=2-3, to=2-4]
	\arrow[from=2-4, to=2-5]
	\arrow[from=3-1, to=3-2]
	\arrow[from=3-2, to=3-3]
	\arrow[from=3-3, to=3-4]
	\arrow[from=3-4, to=3-5]
\end{tikzcd}\end{center}
		\item Defining again $E^2_{j,i}=H_{j,i}(E^1_{\bullet,\bullet})$ we obtain (by a snake lemma type argument) a sequence of complexes

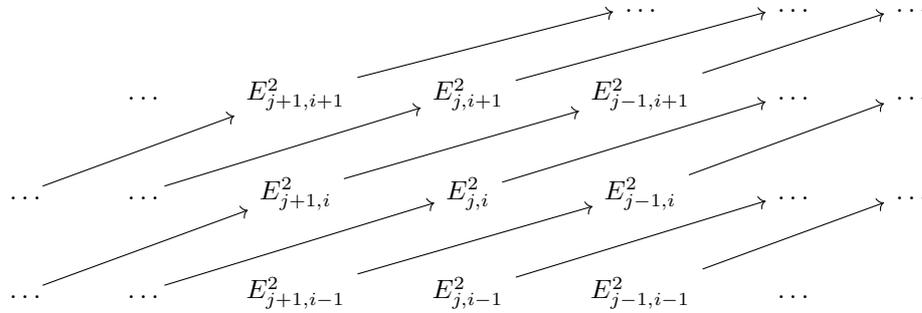
\begin{center}
\begin{tikzcd}
	&&&& \dots & \dots & \dots \\
	& \dots & {E^2_{j+1,i+1}} & {E^2_{j,i+1}} & {E^2_{j-1,i+1}} & \dots & \dots \\
	\dots & \dots & {E^2_{j+1,i}} & {E^2_{j,i}} & {E^2_{j-1,i}} & \dots & \dots \\
	\dots & \dots & {E^2_{j+1,i-1}} & {E^2_{j,i-1}} & {E^2_{j-1,i-1}} & \dots
	\arrow[from=2-3, to=1-5]
	\arrow[from=2-4, to=1-6]
	\arrow[from=2-5, to=1-7]
	\arrow[from=3-1, to=2-3]
	\arrow[from=3-2, to=2-4]
	\arrow[from=3-3, to=2-5]
	\arrow[from=3-4, to=2-6]
	\arrow[from=3-5, to=2-7]
	\arrow[from=4-1, to=3-3]
	\arrow[from=4-2, to=3-4]
	\arrow[from=4-3, to=3-5]
	\arrow[from=4-4, to=3-6]
	\arrow[from=4-5, to=3-7]
\end{tikzcd}
\end{center}
		\item iteratively, define $E^{k+1}_{j,i}=H_{j,i}(E^{k}_{\bullet, \bullet})$ with morphisms (again being induced by a snake lemma type argument) $E^{k+1}_{j,i}\to E^{k+1}_{j-k,i-1+k}$,
		the directions move as follows:

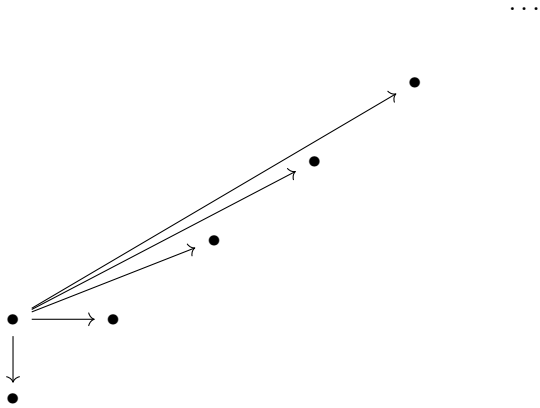
\begin{center}\begin{tikzcd}
	&&&&& \dots \\
	&&&& \bullet \\
	&&& \bullet \\
	&& \bullet \\
	\bullet & \bullet \\
	\bullet
	\arrow[from=5-1, to=2-5]
	\arrow[from=5-1, to=3-4]
	\arrow[from=5-1, to=4-3]
	\arrow[from=5-1, to=5-2]
	\arrow[from=5-1, to=6-1]
\end{tikzcd}.\end{center}
	\end{itemize}

	Dually, define the \idx{spectral sequence with rightward orientation} associated to $A_{\bullet,\bullet}$ via
\begin{itemize}
		\item 0'th page

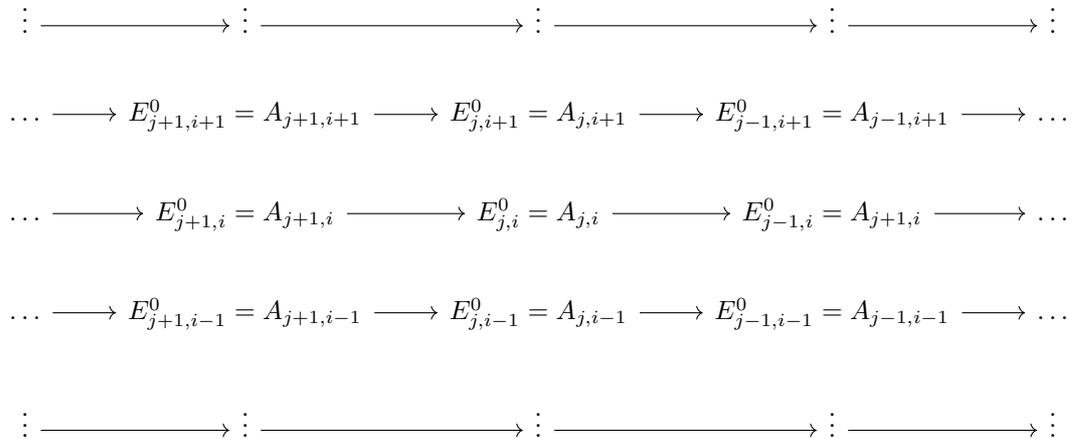
\begin{center}
\begin{tikzcd}
	\vdots & \vdots & \vdots & \vdots & \vdots \\
	\dots & {E^0_{j+1,i+1}=A_{j+1,i+1}} & {E^{0}_{j,i+1}=A_{j,i+1}} & {E^0_{j-1,i+1}=A_{j-1,i+1}} & \dots \\
	\dots & {E_{j+1,i}^0=A_{j+1,i}} & {E_{j,i}^0=A_{j,i}} & {E^0_{j-1,i}=A_{j+1,i}} & \dots \\
	\dots & {E^{0}_{j+1,i-1}=A_{j+1,i-1}} & {E^0_{j,i-1}=A_{j,i-1}} & {E_{j-1,i-1}^0=A_{j-1,i-1}} & \dots \\
	\vdots & \vdots & \vdots & \vdots & \vdots
	\arrow[from=1-1, to=1-2]
	\arrow[from=1-2, to=1-3]
	\arrow[from=1-3, to=1-4]
	\arrow[from=1-4, to=1-5]
	\arrow[from=2-1, to=2-2]
	\arrow[from=2-2, to=2-3]
	\arrow[from=2-3, to=2-4]
	\arrow[from=2-4, to=2-5]
	\arrow[from=3-1, to=3-2]
	\arrow[from=3-2, to=3-3]
	\arrow[from=3-3, to=3-4]
	\arrow[from=3-4, to=3-5]
	\arrow[from=4-1, to=4-2]
	\arrow[from=4-2, to=4-3]
	\arrow[from=4-3, to=4-4]
	\arrow[from=4-4, to=4-5]
	\arrow[from=5-1, to=5-2]
	\arrow[from=5-2, to=5-3]
	\arrow[from=5-3, to=5-4]
	\arrow[from=5-4, to=5-5]
\end{tikzcd}\end{center}
\item pages $E^{k+1}_{j,i}=H_{j,i}(E^{k}_{\bullet,\bullet})$ with induced arrows $E^{k+1}_{j,i}\to E_{k+1}^{j+k-1, i-k}$, first page

\begin{center}
\begin{tikzcd}
	\vdots & \vdots & \vdots \\
	{E^1_{j+1,i+1}} & {E^1_{j,i+1}} & {E^1_{j-1,i+1}} \\
	{E^1_{j+1,i}} & {E^1_{j,i}} & {E^1_{j-1,i}} \\
	{E^1_{j+1,i-1}} & {E^1_{j,i-1}} & {E^1_{j-1,i-1}} \\
	\vdots & \vdots & \vdots
	\arrow[from=1-1, to=2-1]
	\arrow[from=1-2, to=2-2]
	\arrow[from=1-3, to=2-3]
	\arrow[from=2-1, to=3-1]
	\arrow[from=2-2, to=3-2]
	\arrow[from=2-3, to=3-3]
	\arrow[from=3-1, to=4-1]
	\arrow[from=3-2, to=4-2]
	\arrow[from=3-3, to=4-3]
	\arrow[from=4-1, to=5-1]
	\arrow[from=4-2, to=5-2]
	\arrow[from=4-3, to=5-3]
\end{tikzcd},\end{center}
	second page
\begin{center}
\begin{tikzcd}
	& \vdots & \vdots & \vdots \\
	& {E^2_{j+1,i+1}} & {E^2_{j,i+1}} & {E^2_{j-1,i+1}} \\
	\bullet & {E^2_{j+1,i}} & {E^2_{j,i}} & {E^2_{j-1,i}} \\
	\bullet & {E^2_{j+1,i-1}} & {E^2_{j,i-1}} & {E^2_{j-1,i-1}} \\
	\bullet & \vdots & \vdots & \vdots \\
	\bullet & \bullet & \bullet
	\arrow[from=1-2, to=3-1]
	\arrow[from=1-3, to=3-2]
	\arrow[from=1-4, to=3-3]
	\arrow[from=2-2, to=4-1]
	\arrow[from=2-3, to=4-2]
	\arrow[from=2-4, to=4-3]
	\arrow[from=3-2, to=5-1]
	\arrow[from=3-3, to=5-2]
	\arrow[from=3-4, to=5-3]
	\arrow[from=4-2, to=6-1]
	\arrow[from=4-3, to=6-2]
	\arrow[from=4-4, to=6-3]
\end{tikzcd}\end{center}
	and generally directions
\begin{center}\begin{tikzcd}
	&&&& \bullet & \bullet \\
	&&&& \bullet \\
	&&& \bullet \\
	&& \bullet \\
	& \bullet \\
	\dots
	\arrow[from=1-5, to=1-6]
	\arrow[from=1-5, to=2-5]
	\arrow[from=1-5, to=3-4]
	\arrow[from=1-5, to=4-3]
	\arrow[from=1-5, to=5-2]
\end{tikzcd}\end{center}
\end{itemize}
	Informally, in every page one successively takes homology of the last page and obtains induced arrows that walk more and more parallel to the diagonal.
\end{example}
We generalize this situation by allowing isomorphisms in every step, and by forgetting that these pages come from a double complex and that the arrows are induced by some snake lemma argument.
\begin{definition}[Spectral sequence]\label{def:spectral_sequence}
A \idx{downward oriented spectral sequence} (starting with page $a$) is given by
	\begin{itemize}
		\item For every $r\ge a$ a family of objects $(E^r_{j,i})_{j,i\in\Z}$.
		\item Maps $d^r_{j,i}\colon E^r_{j,i}\to E^r_{j-r,i-1+r}$ inducing chain complexes, i.e. $d^r_{j-r,i-1+r}d^r_{j,i}=0$ for all $i,j$.
		\item Isomorphisms
		\[E^{r+1}_{j,i}\simeq H_{j,i}(E^r_{j,i})=\ker(d^r_{j,i})/\Im(d^r_{j+r,i-r+1})\]
\end{itemize}
	A \idx{rightward oriented spectral sequence} (starting with page $a$) is given by
	\begin{itemize}
		\item For every $r\ge a$ a family of objects $(E^r_{j,i})_{j,i\in\Z}$.
		\item Maps $d^r_{j,i}\colon E^r_{j,i}\to E^r_{j-1+r,i-r}$ inducing chain complexes, i.e. $d^r_{j-1+r,i-r}d^r_{j,i}=0$ for all $i,j$.
		\item Isomorphisms
		\[E^{r+1}_{j,i}\simeq H_{j,i}(E^r_{j,i})=\ker(d^r_{j,i})/\Im(d^r_{j-r+1,i+r})\]
\end{itemize}

	We say that a spectral sequence is \idx{bounded}, if for every $k$ there are only finitely many $i,j$ with $i+j=k$ and $E^a_{j,i}\ne 0$.
Clearly, every upper left quadrant spectral sequence fulfills this.
\end{definition}
\begin{definition}[Limit of spectral sequence]\label{def:lim_of_spectral_seq}\uses{def:spectral_sequence}
	Assume that our abelian category is bicomplete.
	Then, given any spectral sequence $E^\bullet_{\bullet,\bullet}$, define $E^\infty_{\bullet,\bullet}$ via
	\[B^\infty_{j,i}=\bigcup_{r\ge a}B_{j,i}^r,\qquad Z^\infty_{j,i}=\bigcap_{r\ge a}Z^r_{j,i}, \quad E^\infty_{j,i}=Z^\infty_{j,i}/B^\infty_{j,i}.\]
Where $Z^r$ are the kernels in the $r$'th page and $B^r$ are the images.

	For any complex $G_\bullet$ we say that $E^\bullet_{\bullet,\bullet}$ \idx{approaches} $G_\bullet$, if
	for every $k$ there exists a \idx{filtration} of $B_k$, being a sequence of inclusions
	\[\dots\sub F_{p-1}G_k\sub F_pG_k\sub F_{p+1} G_k\sub \dots\sub G_k\]
	such that
	\begin{itemize}
		\item there are isomorphisms $E_{j,i}^\infty\simeq F_j G_{j+i} / F_{j-1} G_{j+i}$ for all $p,q$ (this condition is called \idx{weak convergence} of the spectral sequence)
		\item $G_n=\bigcup_p F_p H_n$ and $0=\bigcap_p F_p H_n$.
\end{itemize}
\end{definition}
\begin{definition}[Regular spectral sequence]\label{def:reg_spec_seq}\uses{def:spectral_sequence}
A spectral sequence $E^\bullet_{\bullet,\bullet}$ is called \idx{regular spectral sequence}, if for any $i,j$ there exists a $r\ge a$ such that for all $k\ge r$
	one has $d^k_{j,i}=0$.

	A spectral sequence $E^\bullet_{\bullet,\bullet}$ \idx{converges} to a complex $B_\bullet$, if
	\begin{itemize}
		\item it approaches $B_\bullet$,
		\item it is regular, and
		\item for all $n$n the equality $B_n=\varprojlim_{p} (B_n/F_pB_n)$ holds.
	\end{itemize}
	We denote this by
	\[E^\bullet_{j,i}\implies B_{i+j}.\]
\end{definition}

We suggest to read \cite[chap.~5]{Weibel1994} for many criteria for spectral sequences to converge.

\begin{theorem}[Convergence of spectral sequence of double complex]\label{thm:spectral_sequences}\uses{def:spectral_sequence}
	The associated spectral sequence to an upper left (or dually lower right) quadrant double complex $A_{j,i}$ (meaning a double complex with all $A_{j,i}$ outside of the upper left quadrant) converges,
	as for all places $(i,j)$ after enough pages the spectral sequence always has only 0-morphisms at this place.
	It converges to the homology of the total complex (note that this is independent of the orientation!),
	\[E^\bullet_{j,i}\implies H_{i+j}(\tot(A_{j,i})).\]

	Spelling this out (taking the downwards oriented spectral sequence), it means that for the family $E^\infty_{\bullet,\bullet}$

\begin{center}\begin{tikzcd}
	&&& \dots \\
	&& \dots & {E^\infty_{0,2}} \\
	& \dots & {E^\infty_{1,1}} & {E^\infty_{0,1}} \\
	\dots & {E^\infty_{2,0}} & {E^\infty_{1,0}} & {E^\infty_{0,0}}
\end{tikzcd}\end{center}
we have
	\begin{itemize}
		\item $H^0(\tot(A_{\bullet,\bullet}))=E^\infty_{0,0}$
		\item There is a short exact sequence
\begin{center}\begin{tikzcd}
	0 & {E^\infty_{0,1}} & {H_1(\tot(A_{\bullet,\bullet}))} & {E^\infty_{1,0}} & 0
	\arrow[from=1-1, to=1-2]
	\arrow[from=1-2, to=1-3]
	\arrow[from=1-3, to=1-4]
	\arrow[from=1-4, to=1-5]
\end{tikzcd}\end{center}
		\item There is an object $F_{1}H_2(A_{\bullet,\bullet})$ and a short exact sequence

\begin{center}\begin{tikzcd}
	0 & {E^\infty_{0,2}} & {F_1H_2(\tot(A_{\bullet,\bullet}))} & {E^\infty_{1,1}} & 0
	\arrow[from=1-1, to=1-2]
	\arrow[from=1-2, to=1-3]
	\arrow[from=1-3, to=1-4]
	\arrow[from=1-4, to=1-5]
\end{tikzcd}\end{center}
		with another short exact sequence
\begin{center}\begin{tikzcd}
	0 & {F_1H_2(\tot(A_{\bullet,\bullet}))} & {H_2(\tot(A_{\bullet,\bullet}))} & {E^\infty_{2,0}} & 0
	\arrow[from=1-1, to=1-2]
	\arrow[from=1-2, to=1-3]
	\arrow[from=1-3, to=1-4]
	\arrow[from=1-4, to=1-5]
\end{tikzcd}\end{center}
			\item \dots
\end{itemize}
Using the rightwards oriented associated spectral sequence we obtain a similar result (but with subobjects and quotient objects reversed).
\footnote{A useful tip from \cite{Vakil2023} to remember which direction the filtration works: The differentials on later and later pages point deeper and deeper into the filtration.}

There are similar convergence results for other quadrants, see Chap 5.6 in \cite{Weibel1994}.
\end{theorem}

To explain the general schematic of the \enquote{$F_{i}H$'s},
we explicitly write down all the exact sequences arising from the fifth diagonal, omitting the $\infty$'s on the $E$'s and the \enquote{$(\tot(A_{\bullet,\bullet}))$} after each $H_{i}$:

\begin{center}\begin{tikzcd}
	{(0} & 0 & {E_{0,4}} & {E_{0,4}} & {0)} \\
	0 & {E_{0,4}} & {F_1H_4} & {E_{1,3}} & 0 \\
	0 & {F_1H_4} & {F_2H_4} & {E_{2,2}} & 0 \\
	0 & {F_2H_4} & {F_3 H_4} & {E_{3,1}} & 0 \\
	0 & {F_3 H_4} & {H_4} & {E_{4,0}} & 0 \\
	{(0} & {H_4} & {H_4} & 0 & {0)}
	\arrow[from=1-1, to=1-2]
	\arrow[from=1-2, to=1-3]
	\arrow[from=1-3, to=1-4]
	\arrow[from=1-4, to=1-5]
	\arrow[from=2-1, to=2-2]
	\arrow[from=2-2, to=2-3]
	\arrow[from=2-3, to=2-4]
	\arrow[from=2-4, to=2-5]
	\arrow[from=3-1, to=3-2]
	\arrow[from=3-2, to=3-3]
	\arrow[from=3-3, to=3-4]
	\arrow[from=3-4, to=3-5]
	\arrow[from=4-1, to=4-2]
	\arrow[from=4-2, to=4-3]
	\arrow[from=4-3, to=4-4]
	\arrow[from=4-4, to=4-5]
	\arrow[from=5-1, to=5-2]
	\arrow[from=5-2, to=5-3]
	\arrow[from=5-3, to=5-4]
	\arrow[from=5-4, to=5-5]
	\arrow[from=6-1, to=6-2]
	\arrow[from=6-2, to=6-3]
	\arrow[from=6-3, to=6-4]
	\arrow[from=6-4, to=6-5]
  \end{tikzcd}\end{center}
\begin{remark}
  Note that in the above formulation,
  $F_{i}H_{k}(\tot(A_{\bullet,\bullet}))$ is just a name and can refer to any object of $\mcA$.
  The reason for this notation is that we get a filtration of $H_{k}(\tot(A_{\bullet,\bullet}))$.
\end{remark}

Spectral sequences arising from double complexes are not the only instance of spectral sequences.
In particular, the following spectral sequence associated to a filtration to us seems relevant for applications.
	\begin{example}[Spectral sequence associated to a filtration]\label{ex:spec_seq_of_filtered}
	A \idx{filtration of a complex} $C_\bullet$ is a family of chain subcomplexes
	\[\dots\sub F_{p-1}C_\bullet\sub F_p C_\bullet\sub \dots\sub C_\bullet.\]
	Any filtration determines a spectral sequence starting with $E^0_{j,i}=F_{j}C_{i+j}/F_{j-1}C_{i+1}$.
	For a detailed explanation of this construction we refer to \cite{Gelfand2003} and \cite{Weibel1994}.
	\end{example}
See, e.g., \cite[Theorem 5.5.10]{Weibel1994} for a criterion for spectral sequences associated to filtrations to converge:
\begin{theorem}[Complete convergence theorem]\label{thm:complete_convergence}\uses{def:spectral_sequence}
If the filtration is complete and exhaustive and the corresponding spectral sequence is regular and bounded above
	or if the filtration is bounded, then the associated spectral sequence converges to $H_\bullet(C_\bullet)$.
\end{theorem}
This is a powerful tool for computing homologies through filtrations.
Further easy application of spectral sequences include the following lemmas.
\begin{corollary}
Using spectral sequences of double complexes one can easily prove e.g.
	\begin{itemize}
		\item The snake lemma (or the Zig-zag lemma), using

\begin{center}\begin{tikzcd}
	& 0 & 0 & 0 \\
	0 & A & B & C & 0 \\
	0 & {A'} & {B'} & {C'} & 0 \\
	& 0 & 0 & 0
	\arrow[from=1-2, to=2-2]
	\arrow[from=1-3, to=2-3]
	\arrow[from=1-4, to=2-4]
	\arrow[from=2-1, to=2-2]
	\arrow[from=2-2, to=2-3]
	\arrow[from=2-2, to=3-2]
	\arrow[from=2-3, to=2-4]
	\arrow[from=2-3, to=3-3]
	\arrow[from=2-4, to=2-5]
	\arrow[from=2-4, to=3-4]
	\arrow[from=3-1, to=3-2]
	\arrow[from=3-2, to=3-3]
	\arrow[from=3-2, to=4-2]
	\arrow[from=3-3, to=3-4]
	\arrow[from=3-3, to=4-3]
	\arrow[from=3-4, to=3-5]
	\arrow[from=3-4, to=4-4]
\end{tikzcd}\end{center}
	\item The five lemma, using the double complex

\begin{center}\begin{tikzcd}
	& 0 & 0 & 0 & 0 & 0 \\
	0 & A & B & C & D & E & 0 \\
	0 & {A'} & {B'} & {C'} & {D'} & {E'} & 0 \\
	& 0 & 0 & 0 & 0 & 0
	\arrow[from=1-2, to=2-2]
	\arrow[from=1-3, to=2-3]
	\arrow[from=1-4, to=2-4]
	\arrow[from=1-5, to=2-5]
	\arrow[from=1-6, to=2-6]
	\arrow[from=2-1, to=2-2]
	\arrow[from=2-2, to=2-3]
	\arrow[from=2-2, to=3-2]
	\arrow[from=2-3, to=2-4]
	\arrow[from=2-3, to=3-3]
	\arrow[from=2-4, to=2-5]
	\arrow[from=2-4, to=3-4]
	\arrow[from=2-5, to=2-6]
	\arrow[from=2-5, to=3-5]
	\arrow[from=2-6, to=2-7]
	\arrow[from=2-6, to=3-6]
	\arrow[from=3-1, to=3-2]
	\arrow[from=3-2, to=3-3]
	\arrow[from=3-2, to=4-2]
	\arrow[from=3-3, to=3-4]
	\arrow[from=3-3, to=4-3]
	\arrow[from=3-4, to=3-5]
	\arrow[from=3-4, to=4-4]
	\arrow[from=3-5, to=3-6]
	\arrow[from=3-5, to=4-5]
	\arrow[from=3-6, to=3-7]
	\arrow[from=3-6, to=4-6]
\end{tikzcd}\end{center}
	\end{itemize}
\end{corollary}

\subsection{Projective resolutions}
Next, we will consider projective resolutions, allowing us to describe an object as a chain complex of nice projective objects.

\begin{lemma}[Projectivity as exactness of $\hom$ functor]\label{lem:projectivity_as_exactness_of_Hom}\uses{def:exact_functor, def:special_objects, thm:properties_of_abelian_categories, lem:projectivity_via_sections}
  In any abelian category,
  all of the following conditions on $P$ are equivalent to $P$ being projective.
\begin{enumerate}[(a)]
\item For every epimorphism $A\twoheadrightarrow B$ and morphism $P\to B$ there is a lift 
\begin{center}\begin{tikzcd}
	& A \\
	P & B.
	\arrow[two heads, from=1-2, to=2-2]
	\arrow[dashed, from=2-1, to=1-2]
	\arrow[from=2-1, to=2-2]
\end{tikzcd}\end{center}
\item Every epimorphism $A\twoheadrightarrow P$ splits 
\begin{center}\begin{tikzcd}
	A & P.
	\arrow[shift left, two heads, from=1-1, to=1-2]
	\arrow[shift left, dashed, from=1-2, to=1-1]
\end{tikzcd}\end{center}

\item The functor $\hom(P,-)$ is exact.

\item Every short exact sequence
\begin{center}\begin{tikzcd}
	0 & A & B & P & 0
	\arrow[from=1-1, to=1-2]
	\arrow[from=1-2, to=1-3]
	\arrow[from=1-3, to=1-4]
	\arrow[from=1-4, to=1-5]
\end{tikzcd}\end{center}
	splits.
	\end{enumerate}
\end{lemma}
\begin{proof}
Every epimorphism is regular, thus the first item follows.
For the second characterisation note that epimorphisms are stable under pullback and use \ref{lem:projectivity_via_sections}.
The functor $\hom(P,-)$ commutes with finite coproducts and cokernels, hence the statement follows.

	The last equivalence follows by the splitting lemma, combined with the characterisation of projectivity via splittings.
\end{proof}

\begin{lemma}[Projective resolutions of objects]\label{def:projective_resolution}\uses{def:abelian, def:complex, def:enough_projectives}
	In any abelian category with enough projectives, every object $X$ admits a \idx{projective resolution}, i.e., there exists an exact sequence

\begin{center}\begin{tikzcd}
	\dots & {P_{2}} & {{P_1}} & {P_0} & X & 0
	\arrow[from=1-1, to=1-2]
	\arrow[from=1-2, to=1-3]
	\arrow[from=1-3, to=1-4]
	\arrow[from=1-4, to=1-5]
	\arrow[from=1-5, to=1-6]
\end{tikzcd}\end{center}
	such that all the $P_i$ are projective objects.\footnote{Note that we switched to homological indexing conventions, as this seems more natural in this context}
\end{lemma}
\begin{proof} This may, e.g., iteratively be constructed by projecting onto kernels,

\begin{center}\begin{tikzcd}
	\vdots \\
	{\ker(p_{1})} & {P_{1}} \\
	& {\ker(p_0)} & {P_0} \\
	&& X & 0
	\arrow[from=1-1, to=2-1]
	\arrow[from=2-1, to=2-2]
	\arrow["{p_{1}}", two heads, from=2-2, to=3-2]
	\arrow[hook, from=3-2, to=3-3]
	\arrow["{p_0}", two heads, from=3-3, to=4-3]
	\arrow[from=4-3, to=4-4]
\end{tikzcd}\end{center}
\end{proof}
as there are enough projectives.
\begin{remark}
	Note that this construction can be made functorial (assuming a strong enough axiom of choice), as for any $f\colon X\to Y$ and projective resolutions $P_X^\bullet$, $P_Y^\bullet$,
	by projectivity there is an iteratively constructed extension via

\begin{center}\begin{tikzcd}
	\dots & {P_1^X} & {\ker p_0} & {P_0^X} & X & 0 \\
	\dots & {P_1^Y} && {P_0^Y} & Y
	\arrow[from=1-1, to=1-2]
	\arrow[two heads, from=1-2, to=1-3]
	\arrow[dashed, from=1-2, to=2-2]
	\arrow[from=1-3, to=1-4]
	\arrow[from=1-3, to=2-4]
	\arrow["{p_0}"', two heads, from=1-4, to=1-5]
	\arrow[dashed, from=1-4, to=2-4]
	\arrow[from=1-5, to=1-6]
	\arrow[from=1-5, to=2-5]
	\arrow[from=2-1, to=2-2]
	\arrow[from=2-2, to=2-4]
	\arrow[from=2-4, to=2-5]
  \end{tikzcd}\end{center}

(In many applications,
there is a functorial choice for $P\twoheadrightarrow X$.
Using this,
the construction above becomes functorial, without any axiom of choice.)

	Thus, if an abelian category $\mcA$ has enough projectives we can consider the assignment $\mcA\to \Ch(\mcP)$ (by chopping off the last term), where $\mcP$ denotes the full subcategory of projective objects.

\end{remark}

Next, we see that the choice of the projective resolution, as well as the extension of functions to projective resolutions is almost unique.
\begin{lemma}[Uniqueness of projective resolutions]\label{lem:proj_res_unique}\uses{def:projective_resolution}
Suppose one has a projective resolution (indeed, a chain complex with all $P_i$ projective suffices)

\begin{center}\begin{tikzcd}
	\dots & {P_2^X} & {P_1^X} & {P_0^X} & X &0
	\arrow[from=1-1, to=1-2]
	\arrow[from=1-2, to=1-3]
	\arrow[from=1-3, to=1-4]
	\arrow[from=1-4, to=1-5]
	\arrow[from=1-5, to=1-6]
\end{tikzcd}\end{center}
and another resolution

\begin{center}\begin{tikzcd}
	\dots & {P_2^Y} & {P_1^Y} & {P_0^Y} & Y & 0
	\arrow[from=1-1, to=1-2]
	\arrow[from=1-2, to=1-3]
	\arrow[from=1-3, to=1-4]
	\arrow[from=1-4, to=1-5]
	\arrow[from=1-5, to=1-6]
\end{tikzcd}\end{center}
	as well as a morphism $f\colon X\to Y$.
Then this extends to a chain morphism (as seen in the remark above), and this is unique in the following sense:
For any two extensions $\eta,\nu$

\begin{center}\begin{tikzcd}
	\dots & {P_2^X} & {P_1^X} & {P_0^X} & X & 0 \\
	\dots & {P_2^Y} & {P_1^Y} & {P_0^Y} & Y & 0
	\arrow[from=1-1, to=1-2]
	\arrow[from=1-2, to=1-3]
	\arrow["{\eta_2}", shift left, from=1-2, to=2-2]
	\arrow["{\nu_2}"', shift right, from=1-2, to=2-2]
	\arrow[from=1-3, to=1-4]
	\arrow["{\eta_1}", shift left, from=1-3, to=2-3]
	\arrow["{\nu_1}"', shift right, from=1-3, to=2-3]
	\arrow[from=1-4, to=1-5]
	\arrow[from=1-5, to=1-6]
	\arrow["{\eta_0}", shift left, from=1-4, to=2-4]
	\arrow["{\nu_0}"', shift right, from=1-4, to=2-4]
	\arrow["f", from=1-5, to=2-5]
	\arrow["0"', from=1-6, to=2-6]
	\arrow[from=2-1, to=2-2]
	\arrow[from=2-2, to=2-3]
	\arrow[from=2-3, to=2-4]
	\arrow[from=2-4, to=2-5]
	\arrow[from=2-5, to=2-6]
\end{tikzcd}\end{center}
there exists a family $D_i$ of morphisms $D_i\colon P_i^X\to P_{i+1}^Y$

\begin{center}\begin{tikzcd}
	\dots & {P_2^X} & {P_1^X} & {P_0^X} & X & 0 \\
	\dots & {P_2^Y} & {P_1^Y} & {P_0^Y} & Y & 0
	\arrow[from=1-1, to=1-2]
	\arrow[from=1-2, to=1-3]
	\arrow["{D_2}"'{pos=0.4}, from=1-2, to=2-1]
	\arrow["{\eta_2}", shift left, from=1-2, to=2-2]
	\arrow["{\nu_2}"', shift right, from=1-2, to=2-2]
	\arrow[from=1-3, to=1-4]
	\arrow["{D_1}"'{pos=0.4}, from=1-3, to=2-2]
	\arrow["{\eta_1}", shift left, from=1-3, to=2-3]
	\arrow["{\nu_1}"', shift right, from=1-3, to=2-3]
	\arrow[from=1-4, to=1-5]
	\arrow["{D_0}"'{pos=0.4}, from=1-4, to=2-3]
	\arrow["{\eta_0}", shift left, from=1-4, to=2-4]
	\arrow["{\nu_0}"', shift right, from=1-4, to=2-4]
	\arrow[from=1-5, to=1-6]
	\arrow["f", from=1-5, to=2-5]
	\arrow["0"', from=1-6, to=2-6]
	\arrow[from=2-1, to=2-2]
	\arrow[from=2-2, to=2-3]
	\arrow[from=2-3, to=2-4]
	\arrow[from=2-4, to=2-5]
	\arrow[from=2-5, to=2-6]
\end{tikzcd}\end{center}
such that (putting $D_{-1}=0$)  for all $i\ge 0$
\[\delta_{i+1}^Y D_i+D_{i-1}\delta_i^X=\eta_i-\nu_i.\]
\end{lemma}
\begin{proof}See \cite[2.2.6]{Weibel1994} or \cite[III.1.3]{Gelfand2003}.
	\end{proof}
This lack of uniqueness is one possible motivation to define \idx{chain homotopy}.
\begin{definition}[Chain homotopy]\label{def:chain_homotopy}\uses{def:complex}\chapthree
Consider any two chain complexes $A_\bullet, B_\bullet$ with two chain morphisms $\eta, \nu\colon A_\bullet\to B_\bullet$.
	A \idx{chain homotopy} between $\eta$ and $\nu$ is a family $(D_i)_{i\in\Z}$ of morphisms $D_i\colon A_i\to B_{i+1}$

\begin{center}\begin{tikzcd}
	\dots && {A_{i+1}} && {A_i} && {A_{i-1}} && \dots \\
	\\
	\dots && {B_{i+1}} && {B_i} && {B_{i-1}} && \dots
	\arrow[from=1-1, to=1-3]
	\arrow["{\delta_{i+1}^A}", from=1-3, to=1-5]
	\arrow[from=1-3, to=3-1]
	\arrow["{\eta_{i+1}}", shift left, from=1-3, to=3-3]
	\arrow["{\nu_{i+1}}"', shift right, from=1-3, to=3-3]
	\arrow["{\delta_i^A}", from=1-5, to=1-7]
	\arrow["{D_i}"{description}, from=1-5, to=3-3]
	\arrow["{\eta_i}", shift left, from=1-5, to=3-5]
	\arrow["{\nu_i}"', shift right, from=1-5, to=3-5]
	\arrow["{\delta_{i-1}^A}", from=1-7, to=1-9]
	\arrow["{D_{i-1}}"{description}, from=1-7, to=3-5]
	\arrow["{\eta_{i-1}}", shift left, from=1-7, to=3-7]
	\arrow["{\nu_{i-1}}"', shift right, from=1-7, to=3-7]
	\arrow[from=1-9, to=3-7]
	\arrow[from=3-1, to=3-3]
	\arrow["{\delta^B_{i+1}}", from=3-3, to=3-5]
	\arrow["{\delta^{B}_i}", from=3-5, to=3-7]
	\arrow["{\delta_{i-1}^B}", from=3-7, to=3-9]
\end{tikzcd}\end{center}
such that for all $i$ we have
	\[D_{i-1}\delta_i^A+\delta_{i+1}^BD_i=\eta_i-\nu_i,\]
	meaning that the \enquote{parallelogram} adds up to the difference of the short direct paths:

\begin{center}\begin{tikzcd}
	&& {A_i} && {A_{i-1}} \\
	\\
	{B_{i+1}} && {B_i}
	\arrow["{\delta_i^A}", from=1-3, to=1-5]
	\arrow["{D_i}"{description}, from=1-3, to=3-1]
	\arrow["{\eta_i}", shift left, from=1-3, to=3-3]
	\arrow["{\nu_i}"', shift right, from=1-3, to=3-3]
	\arrow["{D_{i-1}}"{description}, from=1-5, to=3-3]
	\arrow["{\delta^B_{i+1}}", from=3-1, to=3-3]
\end{tikzcd}\end{center}

	Two chain maps are called \idx{chain homotopic}, if there exists a chain homotopy between them.
	A chain map $\eta\colon A_{\bullet}\to B_\bullet$ is called \idx{homotopy equivalence} if there exists $\mu\colon B_\bullet\to A_\bullet$ such that
	$\eta\mu$ is chain homotopic to $1_{B_\bullet}$ and $\mu\eta$ is chain homotopic to $1_{A_\bullet}$.
\end{definition}
\begin{remark}
  Critically,
  $(D_{i})_{i\in\Z}$ does not have to be (and usually, is not) a chain map.

  Note that any additive functor preserves chain homotopy.
\end{remark}
Identifying these morphisms (in the sense of~\ref{def:elementary_constructions_of_categories}), we obtain a fully faithful embedding into the \idx{homotopy category of projectives} $K(\mcP)$.
\begin{definition}[Homotopy category]\label{def:homotopy_cat}\uses{def:chain_homotopy, def:elementary_constructions_of_categories}
For any abelian category $\mcA$ define the \idx{homotopy category of chain complexes} $K(\mcA)$ as the quotient category by chain homotopy,
	i.e., for any chain complexes $A_\bullet, B_\bullet$, define $\hom_{K(\mcA)}(A_\bullet,B_\bullet)=\hom_{\Ch(\mcA)}(A,B)/\sim$
where $\eta\sim\nu$ if and only if there exists a chain homotopy from $\eta$ to $\nu$.

For any full subcategory $\mcB\sub \mcA$ define $K(\mcB)$ to be the full subcategory of $K(\mcA)$ consisting of complexes whose objects lie in $\mcB$.

The left bounded/right bounded/positively graded/negatively graded/bounded homotopy categories $K^+(\mcA)$/$K^-(\mcA)$/$K^{\ge 0}$/$K^{\le 0}$/$K^b(\mcA)$
are defined as the full subcategories spanned by bounded/right bounded/etc. complexes of $K(\mcA)$.
\end{definition}
\begin{remark}
Note that homotopy equivalent complexes form isomorphic objects in $K(\mcA)$, as homotopy equivalences become isomorphisms.
	Furthermore, $\hom_{\Ch(\mcA)}(A_\bullet, B_\bullet)\twoheadrightarrow \hom_{K(\mcA)}(A_\bullet,B_\bullet)$ induces an abelian group structure on $\hom_{K(\mcA)}$, making $K(\mcA)$ an additive category.
\end{remark}

\begin{corollary}[Embedding into homotopy category]\label{cor:emb_into_homotopy}\uses{def:homotopy_cat}
Let $\mcA$ be an abelian category with a generating family $\mcP$ of projective objects.
By taking projective resolutions, one obtains fully faithful embeddings $\mcA\to K(\mcP)^-\to K(\mcP)$.
\end{corollary}
We will discuss this category more in the next section.

First, we will see how to glue projective resolutions together, using the following \idx{horseshoe lemma}.
\begin{lemma}[Horseshoe lemma]\label{lem:horseshoe_lemma}\uses{def:complex}
Let $\mcA$ be an abelian category with enough projectives, and

\begin{center}\begin{tikzcd}
	0 & A & B & C & 0
	\arrow[from=1-1, to=1-2]
	\arrow[from=1-2, to=1-3]
	\arrow[from=1-3, to=1-4]
	\arrow[from=1-4, to=1-5]
\end{tikzcd}\end{center}
a short exact sequence in $\mcA$.
For any two projective resolutions $P_\bullet^A\to A$, $P_\bullet^C\to C$ there exists a projective resolution $P_\bullet^B\to B$ of $B$, such that

\begin{center}\begin{tikzcd}
	& \vdots & \vdots & \vdots \\
	0 & {P_2^A} & {P_2^B} & {P_2^C} & 0 \\
	0 & {P_1^A} & {P_1^B} & {P_1^C} & 0 \\
	0 & {P_0^A} & {P_0^B} & {P_0^C} & 0 \\
	0 & A & B & C & 0\qedhere
	\arrow[from=1-2, to=2-2]
	\arrow[from=1-3, to=2-3]
	\arrow[from=1-4, to=2-4]
	\arrow[from=2-1, to=2-2]
	\arrow[from=2-2, to=2-3]
	\arrow[from=2-2, to=3-2]
	\arrow[from=2-3, to=2-4]
	\arrow[from=2-3, to=3-3]
	\arrow[from=2-4, to=2-5]
	\arrow[from=2-4, to=3-4]
	\arrow[from=3-1, to=3-2]
	\arrow[from=3-2, to=3-3]
	\arrow[from=3-2, to=4-2]
	\arrow[from=3-3, to=3-4]
	\arrow[from=3-3, to=4-3]
	\arrow[from=3-4, to=3-5]
	\arrow[from=3-4, to=4-4]
	\arrow[from=4-1, to=4-2]
	\arrow[from=4-2, to=4-3]
	\arrow[from=4-2, to=5-2]
	\arrow[from=4-3, to=4-4]
	\arrow[from=4-3, to=5-3]
	\arrow[from=4-4, to=4-5]
	\arrow[from=4-4, to=5-4]
	\arrow[from=5-1, to=5-2]
	\arrow[from=5-2, to=5-3]
	\arrow[from=5-3, to=5-4]
	\arrow[from=5-4, to=5-5]
\end{tikzcd}\end{center}
commutes, and has exact rows.
\end{lemma}

\begin{proof}
	Construct this iteratively as follows:

\begin{center}\begin{tikzcd}
	& \vdots & \vdots & \vdots \\
	0 & {P_1^A} & {P_1^A\oplus P_1^C} & {P_1^C} & 0 \\
	0 & {\ker p_0^A} & {\ker p_0^B} & {\ker p_0^C} & 0 \\
	0 & {P_0^A} & {P_0^A\oplus P_0^{B}} & {P_0^B} & 0 \\
	0 & A & B & C & 0
	\arrow[from=1-2, to=2-2]
	\arrow[from=1-3, to=2-3]
	\arrow[from=1-4, to=2-4]
	\arrow[from=2-1, to=2-2]
	\arrow[from=2-2, to=2-3]
	\arrow[two heads, from=2-2, to=3-2]
	\arrow[from=2-2, to=3-3]
	\arrow[from=2-3, to=2-4]
	\arrow[dashed, two heads, from=2-3, to=3-3]
	\arrow[from=2-4, to=2-5]
	\arrow[dashed, from=2-4, to=3-3]
	\arrow[two heads, from=2-4, to=3-4]
	\arrow[from=3-1, to=3-2]
	\arrow[from=3-2, to=3-3]
	\arrow[hook, from=3-2, to=4-2]
	\arrow[from=3-3, to=3-4]
	\arrow[hook, from=3-3, to=4-3]
	\arrow[from=3-4, to=3-5]
	\arrow[hook, from=3-4, to=4-4]
	\arrow[from=4-1, to=4-2]
	\arrow[from=4-2, to=4-3]
	\arrow["{p_0^A}", two heads, from=4-2, to=5-2]
	\arrow[from=4-2, to=5-3]
	\arrow[from=4-3, to=4-4]
	\arrow["{p_0^B}", dashed, from=4-3, to=5-3]
	\arrow[from=4-4, to=4-5]
	\arrow[dashed, from=4-4, to=5-3]
	\arrow["{p_0^C}", from=4-4, to=5-4]
	\arrow[from=5-1, to=5-2]
	\arrow[from=5-2, to=5-3]
	\arrow[from=5-3, to=5-4]
	\arrow[from=5-4, to=5-5]
\end{tikzcd}\end{center}
	\end{proof}

\begin{lemma}[Cartan Eilenberg resolutions of complexes]\label{def:cartan_eilenberg}\uses{def:abelian, def:complex, def:enough_projectives}
Any complex $C^\bullet\in\Ch^{-}(\mcA)$\footnote{Here, we use homological and cohomological notation to distiguish rows and colums.} in an abelian category $\mcA$ with enough projectives admits a double complex

\begin{center}\begin{tikzcd}
	& \vdots & \vdots & \vdots \\
	\dots & {P_2^{i-1}} & {P_2^{i}} & {P_2^{i+2}} & \dots \\
	\dots & {P_1^{i-1}} & {P_1^{i}} & {P_1^{i+1}} & \dots \\
	\dots & {P_0^{i-1}} & {P_0^{i}} & {P_0^{i+1}} & \dots \\
	\dots & {C^{i-1}} & {C^i} & {C^{i+1}} & \dots \\
	\dots & 0 & 0 & 0 & \dots
	\arrow[from=1-2, to=2-2]
	\arrow[from=1-3, to=2-3]
	\arrow[from=1-4, to=2-4]
	\arrow[from=2-2, to=2-3]
	\arrow[from=2-2, to=3-2]
	\arrow[from=2-3, to=2-4]
	\arrow[from=2-3, to=3-3]
	\arrow[from=2-4, to=3-4]
	\arrow[from=3-1, to=3-2]
	\arrow[from=3-2, to=3-3]
	\arrow[from=3-2, to=4-2]
	\arrow[from=3-3, to=3-4]
	\arrow[from=3-3, to=4-3]
	\arrow[from=3-4, to=3-5]
	\arrow[from=3-4, to=4-4]
	\arrow[from=4-1, to=4-2]
	\arrow[from=4-2, to=4-3]
	\arrow[from=4-2, to=5-2]
	\arrow[from=4-3, to=4-4]
	\arrow[from=4-3, to=5-3]
	\arrow[from=4-4, to=4-5]
	\arrow[from=4-4, to=5-4]
	\arrow[from=5-1, to=5-2]
	\arrow[from=5-2, to=5-3]
	\arrow[from=5-2, to=6-2]
	\arrow[from=5-3, to=5-4]
	\arrow[from=5-3, to=6-3]
	\arrow[from=5-4, to=5-5]
	\arrow[from=5-4, to=6-4]
\end{tikzcd}\end{center}
such that
\begin{itemize}
	\item all rows are complexes,
	\item all columns are projective resolutions,
	\item if $C^i=0$, then $P_j^i=0$ for all $j$,
\item for any $i$, the induced sequences

\begin{center}
\begin{tikzcd}
	\vdots \\
	{Z^i(P_2^\bullet)} \\
	{Z^i(P_1^\bullet)} \\
	{Z^i(P^\bullet_0)} \\
	{Z^i(C^\bullet)}
	\arrow[from=1-1, to=2-1]
	\arrow[from=2-1, to=3-1]
	\arrow[from=3-1, to=4-1]
	\arrow[from=4-1, to=5-1]
\end{tikzcd} and
\begin{tikzcd}
	\vdots \\
	{B^i(P_2^\bullet)} \\
	{B^i(P_1^\bullet)} \\
	{B^i(P^\bullet_0)} \\
	{B^i(C^\bullet)}
	\arrow[from=1-1, to=2-1]
	\arrow[from=2-1, to=3-1]
	\arrow[from=3-1, to=4-1]
	\arrow[from=4-1, to=5-1]
\end{tikzcd}
	and \begin{tikzcd}
	\vdots \\
	{H^i(P_2^\bullet)} \\
	{H^i(P_1^\bullet)} \\
	{H^i(P^\bullet_0)} \\
	{H^i(C^\bullet)}
	\arrow[from=1-1, to=2-1]
	\arrow[from=2-1, to=3-1]
	\arrow[from=3-1, to=4-1]
	\arrow[from=4-1, to=5-1]
\end{tikzcd}\end{center}
	are projective resolutions (recall that $Z$ denotes the cycles (kernels), and $B$ are the boundaries (images)).
\end{itemize}
	Such a resolution is called \idx{Cartan-Eilenberg resolution}.
\end{lemma}
\begin{proof}
	See \cite[proof of III.7.7]{Gelfand2003}.
\end{proof}
Often (but not always), this choice is functorial, see, e.g., the discussion in \cite[2.3.2]{Riehl2014}.
However, it is \enquote{unique egough} to induce a functor $\Ch^{\le 0}(\mcA)\to K^{\le 0}(\mcP)$, although this does not need to lift to $\Ch(\mcP)$.

\begin{definition}[Cartan-Eilenberg deformation]\label{def:cartan_eilenber_deform}\uses{def:cartan_eilenberg}
We say, that an abelian category $\mcA$ with enough projectives admits a \idx{Cartan-Eilenberg deformation},
if Cartan-Eilenberg resolution composed with taking total complexes induces a functor $P\colon \Ch^{\le0}(\mcA)\to \Ch^{\le0}(\mcA)$
together with a natural transformation $(\eta_{A_\bullet})_{A_\bullet\in \Ch^{\le0}(\mcA)}\colon P\Rightarrow 1$ with all $\eta_{\mcA}\colon P(A_\bullet)\to A_\bullet $ inducing isomorphisms on homologies.
\end{definition}
\begin{example}
The category of $R$-modules admits a Cartan-Eilenberg deformation, see \cite[chap. 2.3]{Riehl2014}.
\end{example}
\begin{lemma}[Cartan Eilenberg quasiisomorphic]\label{lem:cart_eil_isomorphic}\uses{def:cartan_eilenberg}
	The induced morphism $\tot(P_\bullet^\bullet)\to A_\bullet$ of a Cartan Eilenberg resolution of $A_\bullet$ induces isomorphisms on homologies.
\end{lemma}
\begin{proof} Clear by a spectral sequence argument: Take the the downward oriented spectral sequence associated to $P_\bullet^\bullet$; then $E^1$ has just $A_\bullet$ as lowest row and only 0's apart from this.
\end{proof}

\begin{remark} In fact, one can even find projective resolutions of unbounded complexes, see, e.g., \cite{Spaltenstein1988ResolutionsOU} or \cite[Prop.~71]{Murfet2006b}.
\end{remark}


\section{Derived categories}
\quot{I find it very agreeable to stick all sorts of things, which are not much fun when taken
individually, together under the heading of derived functors.}
{A.\ Grothendieck, letter to J.-P.\ Serre, February 18, 1955}

For more on homological algebra see \cite{Riehl2014, Gelfand2003, Weibel1994} or \cite{hatcher2002algebraic, Dowker1950, forster1999lectures, Vakil2023, Niemiro, thomas2000derived, Grothendieck1957, Hartshorne1997,Godemont1958, Vakil, Urbanik2019}.

Having seen how chain complexes work, we wonder which chain complexes carry \enquote{essentially the same information} and hence should be identified.
The goal of this is to collapse enough redundant information to transform the category of chain complexes into something like a natural extension of the abelian category itself.

As an example, consider any object $X\in\mcA$.
This clearly contains the same information as the complex with just $X$ at the $0$'th place

\begin{center}\begin{tikzcd}
	\dots & 0 & X & 0 & \dots
	\arrow[from=1-1, to=1-2]
	\arrow[from=1-2, to=1-3]
	\arrow[from=1-3, to=1-4]
	\arrow[from=1-4, to=1-5]
\end{tikzcd}\end{center}
which we also denote by $X[0]$.

Using this, we can identify $\mcA$  with a full subcategory of the category $\Ch(\mcA)$ of complexes in $\mcA$.
But clearly, the category $\Ch(\mcA)$ has way more flexibility, whenever one does not just want to consider a single object at once.

For the moment, assume the following ($\Ch^{-}$-minded) intuition behind homologies:
While the images \enquote{create} information, the kernels \enquote{kill} information.
The chain condition now tells us that when taking a stroll through our complex, going from left to right, we kill more information than we have preserved from the last step.
The homology of a chain complex at a point, $H^i(C^\bullet)$, now is the \enquote{really killed information} at time step $i$.

There is also a dual ($\Ch^{+}$-minded) point of view.
In it,
we start with a huge baggage of information and then let it drip into the $H^{i}$ at time $i$.
In taking the step from $A^{i-1}$ to $A^{i}$,
we gain a basic supply of (new) possibilities to encode our information.
The chain condition says that we cannot get back wasted opportunities to encode our information.

A complex hence can also be read as: \enquote{at time $i$ you have hidden the information of $H^i$}.

An exact sequence thus essentially carries no information -- the chain maps kill precisely the amount of information (the kernel) that was created in the last step.

This leads to an interesting concept of which chain complexes should \enquote{carry essentially the same information}.
For example, suppose one wants to understand a factor space $X/Y$.
Then clearly, it should be possible to obtain every fact about $X/Y$ by just looking at the embedding $Y\hookrightarrow X$.
This means, the two chain complexes (in the rows)

\begin{center}\begin{tikzcd}
	0 & 0 & {X/Y} & 0 \\
	0 & Y & X & 0
	\arrow[from=1-1, to=1-2]
	\arrow[from=1-2, to=1-3]
	\arrow[from=1-3, to=1-4]
	\arrow[from=2-1, to=1-1]
	\arrow[from=2-1, to=2-2]
	\arrow[from=2-2, to=1-2]
	\arrow[hook, from=2-2, to=2-3]
	\arrow[from=2-3, to=1-3]
	\arrow[from=2-3, to=2-4]
	\arrow[from=2-4, to=1-4]
\end{tikzcd}\end{center}
both should represent \enquote{the same object}.
This can be made precise by looking at the homologies of both complexes -- they agree!
\begin{itemize}
  \item In the upper line all homologies are $0$, except for the one at place of $X/Y$ (let us put this at place $0$):
		As $X/Y\to 0$ has kernel $X/Y$ and  $0\to X/Y$ has image $0$, and thus $H^0=\ker(X/Y\to 0)/\Im(0\to X/Y)=X/Y$.
\item In the second row, at $Y$, we have $\ker(Y\to X)=0$ and $\im(0\to Y)=0$, hence homology $0$.
	At $X$, we have $\ker(X\to 0)=X$ and $\Im(Y\to X)=Y$, hence homology $X/Y$.
\end{itemize}
As another example, assume that one wants to understand an object $Y$, and one knows how to embed this object into a simpler object $X$ that one understands.
It would only be fair if one could recover all information about $Y$ by understanding the larger object $X$ and the quotient $X/Y$.
\begin{center}\begin{tikzcd}
	0 & Y & 0 & 0 \\
	0 & X & {X/Y} & 0
	\arrow[from=1-1, to=1-2]
	\arrow[from=1-2, to=1-3]
	\arrow[from=1-2, to=2-2]
	\arrow[from=1-3, to=1-4]
	\arrow[from=1-3, to=2-3]
	\arrow[from=2-1, to=2-2]
	\arrow[from=2-2, to=2-3]
	\arrow[from=2-3, to=2-4]
\end{tikzcd}\end{center}
Morally, again both rows should carry the same information -- the space $Y$ has precisely the information that is contained in the quotient $X\to X/Y$.
Again, one checks that the homologies agree, since in the second line,
the homology at $X$ is $\ker(X\to X/Y)/\Im(0\to X)\simeq Y$,
and at $X/Y$ we have $\ker(X/Y\to 0)=X/Y$ and $\Im(X\to X/Y)=X/Y$, leading to $\ker/\Im=0$.

Combining these, one sees that if one has an arrow $\phi\colon X\to Y$, the information of the two complexes

\begin{center}
\begin{tikzcd}
	0 & 0 & X & 0 & 0 \\
	0 & A & Y & B & 0
	\arrow[from=1-1, to=1-2]
	\arrow[from=1-2, to=1-3]
	\arrow[from=1-3, to=1-4]
	\arrow["\phi", from=1-3, to=2-3]
	\arrow[from=1-4, to=1-5]
	\arrow[from=2-1, to=2-2]
	\arrow[hook, from=2-2, to=2-3]
	\arrow[from=2-3, to=2-4]
	\arrow[two heads, from=2-4, to=2-5]
\end{tikzcd}\end{center}

agree, precisely if $\phi$ induces an isomorphism
\[X\simeq \ker(Y\to B)/\Im(A\to Y)=\ker(Y\to B)/A.\]

As a further example, suppose one loves to work with projective objects (for example, if the last section has had an enchanting effect)
and knows how to resolve an object $X$,

\begin{center}\begin{tikzcd}
	\dots & {P_2} & {P_1} & {P_0} & X & 0,
	\arrow[from=1-1, to=1-2]
	\arrow[from=1-2, to=1-3]
	\arrow[from=1-3, to=1-4]
	\arrow[from=1-4, to=1-5]
	\arrow[from=1-5, to=1-6]
\end{tikzcd}\end{center}
then the projective resolution should carry all the information about the object $X$, i.e., the two complexes

\begin{center}\begin{tikzcd}
	\dots & {P_2} & {P_1} & {P_0} & 0 \\
	\dots & 0 & 0 & X & 0
	\arrow[from=1-1, to=1-2]
	\arrow[from=1-2, to=1-3]
	\arrow[from=1-2, to=2-2]
	\arrow[from=1-3, to=1-4]
	\arrow[from=1-3, to=2-3]
	\arrow[from=1-4, to=1-5]
	\arrow[from=1-4, to=2-4]
	\arrow[from=1-5, to=2-5]
	\arrow[from=2-1, to=2-2]
	\arrow[from=2-2, to=2-3]
	\arrow[from=2-3, to=2-4]
	\arrow[from=2-4, to=2-5]
\end{tikzcd}\end{center}
should be \enquote{the same}.
Again, one checks that by definition of a resolution, all the \enquote{higher homologies} $H^i$ are $0$, and $H^0=\ker(P_0\to 0)/\Im(P_1\to P_0)=X$.

These examples suggest the following idea:
Identify two complexes if all their homology groups agree.

This is not entirely what we want, as this point of view would completely ignore the transition maps, and hence identify every complex with a trivial complex of cohomologies,

\begin{center}\begin{tikzcd}
	\dots & {C^{i-1}} & {C^i} & {C^{i+1}} & \dots \\
	\dots & {H^{i-1}(C^\bullet)} & {H^{i}(C^\bullet)} & {H^{i+1}(C^\bullet)} & \dots
	\arrow[from=1-1, to=1-2]
	\arrow[from=1-2, to=1-3]
	\arrow[from=1-3, to=1-4]
	\arrow[from=1-4, to=1-5]
	\arrow["0", from=2-1, to=2-2]
	\arrow["0", from=2-2, to=2-3]
	\arrow["0", from=2-3, to=2-4]
	\arrow["0", from=2-4, to=2-5]
\end{tikzcd}\end{center}

As an example, these two chain complexes (of $\C[x,y]$-modules) should not be identified, as they are completely different ways of \enquote{attaching the information},
see~\cite{thomas2000derived}.

\begin{center}\begin{tikzcd}
	0 & {\C[x,y]\oplus \C[x,y]} & {\C[x,y]} & 0 \\
	0 & {\C[x,y]} & \C & 0
	\arrow[from=1-1, to=1-2]
	\arrow["{(x,y)}", from=1-2, to=1-3]
	\arrow[from=1-3, to=1-4]
	\arrow[from=2-1, to=2-2]
	\arrow["0", from=2-2, to=2-3]
	\arrow[from=2-3, to=2-4]
\end{tikzcd}\end{center}

What we should observe is that in the examples above, the isomorphism of the $H^i$ is induced by a morphism of chain complexes.
Two chain complexes where the same amount of information floats around, but is stacked together completely differently, should not be identified.
This leads to the following idea:

\textbf{Morphisms between chain complexes that induce isomorphisms on all homology groups should be isomorphisms.}

This will yield the \idx{derived category}.
First, let us introduce a name for such morphisms.
\begin{definition}[Quasiisomorphism]\label{def:quasiiso}\uses{def:complex}\chapthree
	A morphism between two complexes, $\eta\colon A^\bullet\to  B^\bullet$ is called a \idx{quasiisomorphism}, if all the induced arrows between homologies, $H^i(A^\bullet)\to H^i(B^\bullet)$, are isomorphisms.
\end{definition}

\begin{lemma}[Chain homotopies are quasiisomorphisms]\label{lem:chain_homotopies_are_quasiisomorphisms}\uses{def:chain_homotopy, def:quasiiso}
If two morphisms $\eta,\nu$ between chain complexes $A_\bullet, B_\bullet$ are chain homotopic, then they induce the same morphisms on homologies.
	In particular, any homotopy equivalence $\eta\colon A_\bullet\to B_\bullet$ is a quasiisomorphism.
\end{lemma}
\begin{proof}
See \cite[III.1.2]{Gelfand2003}.
\end{proof}

\begin{definition}[Derived category]\label{def:derived}\uses{def:complex, def:quasiiso}
Consider an abelian category $\mcA$.
Then the \idx{derived category} $D(\mcA)$ of $\mcA$ is the \idx{localisation} of $\Ch(\mcA)$ (or equivalently $K(\mcA)$) at quasiisomorphisms.
This means it is a category $D(\mcA)$ together with a functor $\Ch(\mcA)\to D(\mcA)$ mapping quasiisomorphisms to isomorphisms, such that for any other category $\mcD$ and functor $\Ch(\mcA)\to\mcD$ that maps quasiisomorphisms to isomorphism,
there exists a unique $D(\mcA)\to\mcD$ with

\begin{center}\begin{tikzcd}
	{\Ch(\mcA)} & {D(\mcA)} \\
	& \mcD
	\arrow[from=1-1, to=1-2]
	\arrow[from=1-1, to=2-2]
	\arrow[dashed, from=1-2, to=2-2]
\end{tikzcd}\end{center}
This category exists and can be constructed explicitly by iteratively adding formal inverses to those morphisms which need to have inverses, see \cite[III.2]{Gelfand2003}.

The left bounded/right bounded/ bounded derived categories $D^{+/-/b}(\mcA)$ are defined by localising $\Ch^{+/-/b}(\mcA)$ at quasiisomorphisms.
The positively graded/negatively graded/concentrated in degree 0 derived category, $D^{\ge 0}, D^{\le 0}, D^0$ are defined analogously.
The category $D^-(\mcA)$ is equivalent to the full subcategories of $D(\mcA)$ spanned by complexes with $H^i(C^\bullet)=0$ for all $i\ge i_0(C^\bullet)$
(or the corresponding conditions for the other cases) (\cite[III.2.5]{Gelfand2003}).
\end{definition}
\begin{lemma}[Description of the derived category]\label{lem:descr_derived}\uses{def:derived}
  One can describe the derived category $D(\mcA)$ as follows: The objects are complexes $C_\bullet$, and morphisms $C_\bullet\to D_\bullet$ are represented by \idx{roofs},
  which are objects $A_\bullet$ together with a quasiisomorphism $q\colon A_\bullet\to C_\bullet$ and a chain morphism $A_\bullet\to C_\bullet$,

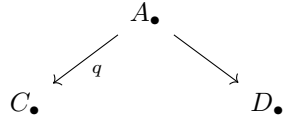
\begin{center}\begin{tikzcd}
	& {A_\bullet} \\
	{C_\bullet} && {D_\bullet}
	\arrow["q", from=1-2, to=2-1]
	\arrow[from=1-2, to=2-3]
\end{tikzcd}\end{center}

Composition of two roofs is defined as follows: for any two roofs $A\to B$ and $B\to C$ there exists a roof connecting the \enquote{tops},
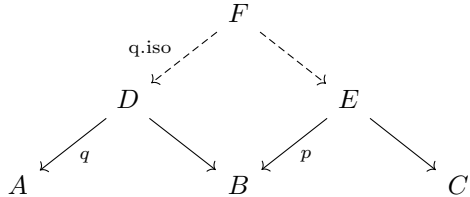
\begin{center}\begin{tikzcd}
	&& F \\
	& D && E \\
	A && B && C
	\arrow["{\mathrm{q.iso}}"'{pos=0.6}, dashed, from=1-3, to=2-2]
	\arrow[dashed, from=1-3, to=2-4]
	\arrow["q", from=2-2, to=3-1]
	\arrow[from=2-2, to=3-3]
	\arrow["p", from=2-4, to=3-3]
	\arrow[from=2-4, to=3-5]
\end{tikzcd}\end{center}
where the parallelogram in the middle commutes up to chain homotopy.

Note that there is an equivalence relation on the roofs, making the representatives non-unique.
\end{lemma}
\begin{proof}
See \cite[2.1.6]{Riehl2014}.
\end{proof}
\begin{theorem}[Elementary properties of derived categories]\label{thm:elementary_properties_of_derived_categories}\uses{def:derived}
The derived category is usually not abelian, however, it still satisfies some useful properties.
\begin{itemize}
	\item The derived category $D(\mcA)$ is additive.
	\item If $\mcA$ is (co)complete and (co)products are exact, then $\mcD(\mcA)$ and $K(\mcA)$ admit (co)products and they may be computed pointwise in $\Ch(\mcA)$.  
	\item The embedding $\mcA\to D(\mcA)$, $A\mapsto A[0]$ exhibits an equivalence between $\mcA$ and the full subcategory of $D(\mcA)$ consisting of complexes with $H^i=0$ for all $i\ne 0$ (we call these complexes \idx{homologically concentrated in degree 0}).
\end{itemize}
\end{theorem}
\begin{proof}
A good reference is \cite[prop.~41--46]{Murfet2006b}.
\end{proof}
The next theorem is of great interest in concretely describing the abstractly defined derived category.
\begin{proposition}[Derived category via projectives]\label{prop:derived_via_proj}\uses{def:derived, def:homotopy_cat, lem:proj_res_unique}
	Let $\mcA$ be an abelian category with a generating class of projectives $\mcP$.
	Then
	\[D^-(\mcA)\simeq K^-(\mcP)\]
\end{proposition}
\begin{proof}
  This can be found, e.g., in \cite{Gelfand2003}, using a Cartan-Eilenberg resolution.
  The main point is to see that quasisomorphisms between complexes of projectives already form chain homotopies.
\end{proof}

\subsubsection{Derived functors}
The next important concept are \idx{derived functors}, which describe the \enquote{best} way to extend a given functor $\mcA\to\mcB$ between abelian categories to the derived categories $D(\mcA)\to D(\mcB)$
and special cases of which give intensely studied classical notions.
\begin{definition}[Derived functor]\label{def:derived_functor}\uses{def:derived, def:functor, def:kan_extensions}
	Consider any functor $F\colon\Ch(\mcA)\to\Ch(\mcB)$ between abelian categories, and the localisations $q_{\mcA}\colon \Ch(\mcA)\to D(\mcA)$, $q_{\mcB}\colon \Ch(\mcB)\to D(\mcB)$.
	Then the \idx{left derived functor} of $F$ is the \emph{right} Kan extension of $q_{\mcB} F$ along $q_{\mcA}$:
	\[LF=\Ran_{q_\mcA}(q_{\mcB}F)\]
Explicitly, $LF$ is equipped with a natural transformation $\eps\colon LF q_{\mcA}\Rightarrow q_{\mcB} F$

\begin{center}
\begin{tikzcd}
	\Ch(\mcA) & \Ch(\mcB) \\
	{D(\mcA)} & {D(\mcB)}
	\arrow["F", from=1-1, to=1-2]
	\arrow["{q_{\mcA}}"', from=1-1, to=2-1]
	\arrow["{q_{\mcB}}", from=1-2, to=2-2]
	\arrow["\eps"'{pos=0.6}, Rightarrow, from=2-1, to=1-2]
	\arrow["LF"', from=2-1, to=2-2]
\end{tikzcd}
\end{center}
such that for all other functors $G\colon D(\mcA)\to D(\mcB)$ with natural transformations $\mu\colon G q_{\mcA}\Rightarrow q_{\mcB}F$, there exists a unique $\nu\colon G\Rightarrow LF$ with

\begin{center}
\begin{tikzcd}
	{G q_\mcA} & {q_{\mcB}F} \\
	& {(LF)q_{\mcA} }
	\arrow["\mu",Rightarrow, from=1-1, to=1-2]
	\arrow["{\nu.q_{\mcA}}"'{pos=0.3},Rightarrow,  from=1-1, to=2-2]
	\arrow["\eps"',Rightarrow, from=2-2, to=1-2]
\end{tikzcd}\end{center}

	The \idx{right derived functor} of $F$ is the \emph{left} Kan extension of $q_{\mcB} F$ along $q_{\mcA}$:
		\[RF=\Lan_{q_\mcA}(q_{\mcB}F)\]
		Again, part of the data is a natural transformation $\eta\colon q_{\mcB} F\Rightarrow RFq_{\mcA}$ such that for any $G\colon D(\mcA)\to D(\mcB)$ with natural transformation $\mu\colon q_{\mcB} F\Rightarrow G q_{\mcA}$,
		there exists a unique $\nu\colon  RF\Rightarrow G$ with

\begin{center}\begin{tikzcd}
	{q_{\mcB} F} & {(RF)q_{\mcA}} \\
	& {Gq_{\mcA}}
	\arrow["\eta",Rightarrow, from=1-1, to=1-2]
	\arrow["\mu"',Rightarrow, from=1-1, to=2-2]
	\arrow["\nu",Rightarrow, dashed, from=1-2, to=2-2]
\end{tikzcd}
\end{center}

Note that, equivalently, we can identify such a functor with a functor $\Ch(\mcA)\to D(\mcB)$ sending quasiisomorphisms to isomorphisms, which we will sometimes do implicitly.

	For any functor $F\colon \mcA\to \mcB$ define the right/left derived functor of $F$ as the derived functors of the canonical pointwise extension $F_\bullet\colon \Ch(\mcA)\to \Ch(\mcB)$.
\end{definition}
In some cases, the derived functors lift to $\Ch(\mcA)$.
\begin{definition}[Point-set derived functors]\label{def:point_set_derived}\uses{def:derived_functor}
	A \idx{left point-set derived functor} of some $F\colon\Ch(\mcA)\to\Ch(\mcB)$ is a functor $L_p F\colon \Ch(\mcA)\to \Ch(\mcB)$ mapping quasiisomorphisms to isomorphism (or equivalently, a functor $D(\mcA)\to \Ch(\mcB)$),
	with natural transformation $\lambda\colon L_p F\Rightarrow F$ such that $q_{\mcB}.\lambda\colon q_{\mcB} L_p F\Rightarrow q_{\mcB} F $ is a left derived functor of $F$.
	Dually we define \idx{right point-set derived functors}.
\end{definition}

Note that we do not claim that the derived functor of a given functor exists.
In fact, often they do not exist.
In the following, we will find criteria for them to exist and methods of computing them.

\begin{remark}
We define completely analogously the derived functors on the bounded versions of the derived categories.

Note that the derived functor on the bounded derived category agrees with the restriction of the unbounded derived functor, see, e.g., \cite[lemma 14]{Murfet2006c}
\end{remark}

The first trivial case is given if the functor preserves quasiisomorphisms.
In this case, the functor \enquote{survives} the localisation and indeed can be extended naively.
\begin{lemma}[Derived functors of exact functors]\label{lem:derived_functors_of_exact_functors}\uses{def:derived_functor, def:derived, def:exact_functor, lem:characterisation_of_exact_functors}
For any exact functor $F\colon \mcA\to\mcB$ the left and right derived functors agree and are induced by the pointwise derived functors,
which can be computed by pointwise applying $F$.
\end{lemma}

\begin{proof}
This is due to the fact that an exact functor induces a unique functor making the following commute

\begin{center}\begin{tikzcd}
	{\Ch(\mcA)} & {\Ch(\mcB)} \\
	{D(\mcA)} & {D(\mcB)}
	\arrow["{F_\bullet}", from=1-1, to=1-2]
	\arrow["{q_{\mcA}}"', from=1-1, to=2-1]
	\arrow["{q_{\mcB}}", from=1-2, to=2-2]
	\arrow[dashed, from=2-1, to=2-2]
\end{tikzcd}\end{center}
\end{proof}
For the following, see, e.g., \cite[2.4]{Weibel1994}.
\begin{proposition}\label{thm:elementary_properties_of_derived_functors}\uses{def:derived_functor, def:pointwise_kan_extensions}
The derived functor of an additive functor is additive.
\end{proposition}

Next, we will approach the computation of derived functors via easier objects, and in particular in the presence of enough projectives (see~\ref{def:cartan_eilenberg}, \ref{def:cartan_eilenber_deform}).
\begin{definition}
	A  \idx{left deformation} of $\Ch(\mcA)$ is an endofunctor $Q\colon \Ch(\mcA)\to \Ch(\mcA)$ with natural transformation $Q\Rightarrow 1$, whose components are quasiisomorphisms.
	A  \idx{right deformation} of $\Ch(\mcA)$ is an endofunctor $Q\colon \Ch(\mcA)\to \Ch(\mcA)$ with natural transformation $ 1\Rightarrow Q$, whose components are quasiisomorphisms.

	For any deformation we have an equivalence $D(\Im(Q))\simeq D(\mcA)$ for $\Im(Q)$ being the essential image of $Q$.

A left/right deformation $Q\colon\Ch(\mcA)\to\Ch(\mcA)$ is a \idx{left/right deformation for a functor} $F\colon \Ch(\mcA)\to \Ch(\mcB)$ if the functor restricted to the essential image preserves quasiisomorphisms.
	We say that $F\colon\Ch(\mcA)\to \Ch(\mcB)$ is (left/right) \idx{deformable}, if a (left/right) deformation exists for $F$.
	\end{definition}
\begin{theorem}
	If $F\colon\Ch(\mcA)\to\Ch(\mcB)$ is left deformable by some $Q$, then $FQ$ is a point-set left derived functor of $F$.
	\end{theorem}
\begin{proof} This is \cite[2.2.8]{Riehl2014}.
	\end{proof}
\begin{proposition}
	A derived functor of a deformable functor is an \emph{absolute} Kan extension.
	In particular, this absolute Kan extension preserves adjunctions, explaining that adjunctions may be derived, inducing derived adjunctions.
\end{proposition}
\begin{proof} \cite[2.2.13]{Riehl2014}.
\end{proof}
\begin{lemma}
  Any additive functor maps quasiisomorphism between complexes with projective objects to quasiisomorphisms.
\end{lemma}
\begin{proof} As any quasiisomorphism between complexes of projective objects is a chain homotopy, this is clear.
\end{proof}
\begin{example}[Computing derived functors]
  If a Cartan-Eilenberg deformation $P$ for an abelian category $\mcA$ with enough projectives $\mcP$ exists,
  any additive functor $F\colon \mcA\to\mcB$  has induced functor $F^\bullet\colon \Ch^{\ge 0}(\mcA)\to \Ch^{\ge 0}$ which is deformed by Cartan-Eilenberg deformation,
  and thus admits a point-set left-derived functor, computed by

\begin{center}\begin{tikzcd}
	{\Ch^{\ge 0}(\mcA)} & {\Ch^{\ge 0}(\mcP)} & {\Ch^{\ge 0}(\mcB).}
	\arrow["P", from=1-1, to=1-2]
	\arrow["{F^\bullet}", from=1-2, to=1-3]
\end{tikzcd}\end{center}
In general, if the abelian category has enough projectives but not necessarily a functorial choice of resolutions,
one can still define a derived functor of any additive functor $F\colon\mcA\to \mcB$ by this procedure,
which then does not have to yield a point-set derived functor:
	\begin{itemize}
		\item Take the complex $A_{\bullet}$
			and build a Cartan-Eilenberg resolution\footnote{For once, we depart from our indexing convention to avoid negative numbers as well as indices separated by comma.
			Lord, have mercy upon us.}
\begin{center}\begin{tikzcd}
	& \vdots & \vdots & \vdots \\
	\dots & {P2^2} & {P_1^2} & {P_0^2} & 0 \\
	\dots & {P_2^1} & {P_1^1} & {P_0^1} & 0 \\
	\dots & {P_2^0} & {P_1^0} & {P_0^0} & 0 \\
	\dots & {A_{2}} & {A_1} & {A_0} & 0 \\
	& 0 & 0 & 0.
	\arrow[from=1-2, to=2-2]
	\arrow[from=1-3, to=2-3]
	\arrow[from=1-4, to=2-4]
	\arrow[from=2-1, to=2-2]
	\arrow[from=2-2, to=2-3]
	\arrow[from=2-2, to=3-2]
	\arrow[from=2-3, to=2-4]
	\arrow[from=2-3, to=3-3]
	\arrow[from=2-4, to=2-5]
	\arrow[from=2-4, to=3-4]
	\arrow[from=3-1, to=3-2]
	\arrow[from=3-2, to=3-3]
	\arrow[from=3-2, to=4-2]
	\arrow[from=3-3, to=3-4]
	\arrow[from=3-3, to=4-3]
	\arrow[from=3-4, to=3-5]
	\arrow[from=3-4, to=4-4]
	\arrow[from=4-1, to=4-2]
	\arrow[from=4-2, to=4-3]
	\arrow[from=4-2, to=5-2]
	\arrow[from=4-3, to=4-4]
	\arrow[from=4-3, to=5-3]
	\arrow[from=4-4, to=4-5]
	\arrow[from=4-4, to=5-4]
	\arrow[from=5-1, to=5-2]
	\arrow[from=5-2, to=5-3]
	\arrow[from=5-2, to=6-2]
	\arrow[from=5-3, to=5-4]
	\arrow[from=5-3, to=6-3]
	\arrow[from=5-4, to=5-5]
	\arrow[from=5-4, to=6-4]
\end{tikzcd}\end{center}
		Delete the lowest row to obtain a double complex of projectives which has quasiisomorphic total complex to $A_\bullet$.
		This does not change the value of the derived functor, $LF(A_\bullet)=LF(\mathrm{tot}(P_\bullet^\bullet))$.
		\item Apply the functor $F$ pointwise to the double complex of projectives (as $F$ is additive this commutes with taking the total complex).

\begin{center}\begin{tikzcd}
	& \vdots & \vdots & \vdots \\
	\dots & {FP_2^2} & {FP_1^2} & {FP_0^2} & 0 \\
	\dots & {FP_2^1} & {FP_1^1} & {FP_0^1} & 0 \\
	\dots & {FP_2^0} & {FP_1^0} & {FP_0^0} & 0 \\
	\dots & 0 & 0 & 0 &
	\arrow[from=1-2, to=2-2]
	\arrow[from=1-3, to=2-3]
	\arrow[from=1-4, to=2-4]
	\arrow[from=2-1, to=2-2]
	\arrow[from=2-2, to=2-3]
	\arrow[from=2-2, to=3-2]
	\arrow[from=2-3, to=2-4]
	\arrow[from=2-3, to=3-3]
	\arrow[from=2-4, to=2-5]
	\arrow[from=2-4, to=3-4]
	\arrow[from=3-1, to=3-2]
	\arrow[from=3-2, to=3-3]
	\arrow[from=3-2, to=4-2]
	\arrow[from=3-3, to=3-4]
	\arrow[from=3-3, to=4-3]
	\arrow[from=3-4, to=3-5]
	\arrow[from=3-4, to=4-4]
	\arrow[from=4-1, to=4-2]
	\arrow[from=4-2, to=4-3]
	\arrow[from=4-2, to=5-2]
	\arrow[from=4-3, to=4-4]
	\arrow[from=4-3, to=5-3]
	\arrow[from=4-4, to=4-5]
	\arrow[from=4-4, to=5-4]
\end{tikzcd}\end{center}
		This yields $LF(\mathrm{tot}(P_\bullet^\bullet))=\mathrm{tot}(F(P_\bullet^\bullet))=\mathrm{tot}(F(P_i^j))_{i\in \Z}^{j\in\Z}$.
		\item Take the total complex of the resulting double complex,

\begin{center}\begin{tikzcd}
	\dots & {FP_2^0\oplus FP_1^1\oplus FP_0^2} & {FP_1^0\oplus FP_0^1 } & {FP_0^0} & 0
	\arrow[from=1-1, to=1-2]
	\arrow[from=1-2, to=1-3]
	\arrow[from=1-3, to=1-4]
	\arrow[from=1-4, to=1-5]
\end{tikzcd}\end{center}
\end{itemize}
Clearly, instead of the first two steps one could also find any other complex which is quasiisomorphic to $A_\bullet$ and consists only of projectives, and act pointwise on this complex.

Instead of taking the total complex in the last step,
it is also possible to directly calculate the homologies of the resulting complex as the limit of the spectral sequence of the double complex,
as this converges (it is a left upper quadrant double complex).
\end{example}

\subsubsection{Triangulated Categories and Truncations}
    \quot{We worked on this book [Methods of Homological Algebra \cite{Gelfand2003}] with the disquieting feeling that the development of homological algebra is currently in a state of flux,
    and that the basic definitions and constructions of the theory of triangulated categories, despite their widespread use, are of only preliminary nature (...).
        There is no doubt that similar thoughts have occured to the founders of the theory, and to everyone who has seriously worked with it; the absence of a monographic exposition is one of the symptoms.}{
Introduction to \cite{Gelfand2003}.}

Another important, although initially scary-looking, property of derived categories is that they are \idx{triangulated categories},
a notion originally introduced by Verdier, see \cite[1.1.2.5]{Lurie2017} for the following definition.

\begin{definition}
A \idx{triangulated category} consists of
\begin{enumerate}
	\item an additive category $\mcD$,
	\item a so-called \idx{shift functor} $\Sigma\colon \mcD\to\mcD,\, X\mapsto X[1]$, inducing an equivalence of $\mcD$ with $\mcD$,
	\item a class of so-called \idx{distinguished triangles}, being some diagrams of the form
\begin{center}\begin{tikzcd}
	X & Y & Z & {X[1]}
	\arrow[from=1-1, to=1-2]
	\arrow[from=1-2, to=1-3]
	\arrow[from=1-3, to=1-4]
\end{tikzcd}\end{center}
	for some $X,Y,Z\in \mcD$.
\end{enumerate}
	We further require the following axioms.
	\begin{enumerate}
		\item \begin{itemize}
				  \item For every morphism $f\colon X\to Y$, there exists a distinguished triangle using $f$,

\begin{center}\begin{tikzcd}
	X & Y & Z & {X[1]}
	\arrow["f", from=1-1, to=1-2]
	\arrow[from=1-2, to=1-3]
	\arrow[from=1-3, to=1-4]
\end{tikzcd}\end{center}
			  \item 	Distinguished triangles are stable under isomorphisms.
This means that if
\begin{center}\begin{tikzcd}
	X & Y & Z & {X[1]}
	\arrow[from=1-1, to=1-2]
	\arrow[from=1-2, to=1-3]
	\arrow[from=1-3, to=1-4]
\end{tikzcd}\end{center}
				  is a distinguished triangle and
\begin{center}\begin{tikzcd}
	{X'} & {Y'} & {Z'} & {X'[1]}
	\arrow[from=1-1, to=1-2]
	\arrow[from=1-2, to=1-3]
	\arrow[from=1-3, to=1-4]
\end{tikzcd}\end{center}
				  an isomorphic triangle, meaning that there exist vertical isomorphisms

\begin{center}\begin{tikzcd}
	X & Y & Z & {X[1]} \\
	{X'} & {Y'} & {Z'} & {X'[1]}
	\arrow[from=1-1, to=1-2]
	\arrow["\simeq"', from=1-1, to=2-1]
	\arrow[from=1-2, to=1-3]
	\arrow["\simeq"', from=1-2, to=2-2]
	\arrow[from=1-3, to=1-4]
	\arrow["\simeq"', from=1-3, to=2-3]
	\arrow["\simeq"', from=1-4, to=2-4]
	\arrow[from=2-1, to=2-2]
	\arrow[from=2-2, to=2-3]
	\arrow[from=2-3, to=2-4]
\end{tikzcd},\end{center}
				  then the second triangle is also distinguished.
				  \item for every $X\in \mcD$, the triangle
\begin{center}\begin{tikzcd}
	X & X & 0 & {X[1]}
	\arrow["{\mathrm{id}_X}", from=1-1, to=1-2]
	\arrow[from=1-2, to=1-3]
	\arrow[from=1-3, to=1-4]
\end{tikzcd}\end{center}
				  is distinguished.
\end{itemize}
		\item A triangle
\begin{center}\begin{tikzcd}
	X & Y & Z & {X[1]}
	\arrow["f", from=1-1, to=1-2]
	\arrow["g", from=1-2, to=1-3]
	\arrow["h", from=1-3, to=1-4]
\end{tikzcd}\end{center}
		is distinguished precisely if the triangle
\begin{center}\begin{tikzcd}
	Y & Z & {X[1]} & {Y[1]}
	\arrow["g", from=1-1, to=1-2]
	\arrow["h", from=1-2, to=1-3]
	\arrow["{f[1]}", from=1-3, to=1-4]
\end{tikzcd}\end{center}
		is distinguished.
		\item Morphisms between distinguished triangles may always be completed as follows
\begin{center}\begin{tikzcd}
	X & Y & Z & {X[1]} \\
	{X'} & {Y'} & {Z'} & {X'[1]}
	\arrow[from=1-1, to=1-2]
	\arrow[from=1-1, to=2-1]
	\arrow[from=1-2, to=1-3]
	\arrow[from=1-2, to=2-2]
	\arrow[from=1-3, to=1-4]
	\arrow[dashed, from=1-3, to=2-3]
	\arrow[from=1-4, to=2-4]
	\arrow[from=2-1, to=2-2]
	\arrow[from=2-2, to=2-3]
	\arrow[from=2-3, to=2-4]
\end{tikzcd}\end{center}
		\item For three distinguished triangles of the form
\begin{center}\begin{tikzcd}
	X & Y & {Y/X} & {X[1]} \\
	Y & Z & {Z/Y} & {Y[1]} \\
	X & Z & {Z/X} & {X[1]}
	\arrow["f", from=1-1, to=1-2]
	\arrow["u", from=1-2, to=1-3]
	\arrow["d", from=1-3, to=1-4]
	\arrow["g", from=2-1, to=2-2]
	\arrow["v", from=2-2, to=2-3]
	\arrow["{d'}", from=2-3, to=2-4]
	\arrow["{g\circ f}", from=3-1, to=3-2]
	\arrow["w", from=3-2, to=3-3]
	\arrow["{d''}", from=3-3, to=3-4]
\end{tikzcd}\end{center}
		there always exists a fourth distinguished triangle
\begin{center}\begin{tikzcd}
	{Y/X} & {Z/X} & {Z/Y} & {(Y/X)[1]}
	\arrow["\phi", from=1-1, to=1-2]
	\arrow["\psi", from=1-2, to=1-3]
	\arrow["\theta", from=1-3, to=1-4]
\end{tikzcd}\end{center}
		with
\begin{center}\begin{tikzcd}
	X && Z && {Z/Y} && {(Y/X)[1]} \\
	& Y && {Z/X} && {Y[1]} \\
	&& {Y/X} && {X[1]}
	\arrow["{g\circ f}", from=1-1, to=1-3]
	\arrow["f", from=1-1, to=2-2]
	\arrow["v", from=1-3, to=1-5]
	\arrow["w", from=1-3, to=2-4]
	\arrow["\theta", from=1-5, to=1-7]
	\arrow["{d'}", from=1-5, to=2-6]
	\arrow["g", from=2-2, to=1-3]
	\arrow["u", from=2-2, to=3-3]
	\arrow["\psi", from=2-4, to=1-5]
	\arrow["{d''}", from=2-4, to=3-5]
	\arrow["{u[1]}"', from=2-6, to=1-7]
	\arrow["\phi", from=3-3, to=2-4]
	\arrow["d", from=3-3, to=3-5]
	\arrow["{f[1]}"', from=3-5, to=2-6]
\end{tikzcd}\end{center}
\end{enumerate}
	An additive functor between triangulated categories is called triangulated, if it sends distinguished triangles to distinguished triangles.
\end{definition}
There is a sizeable theory of triangulated categories, which we will not present here, see, e.g., \cite{Gelfand2003, Weibel1994, kashiwara2005categories} for more.
However, derived categories naturally fulfill this condition.
\begin{lemma}
 The derived category $D(\mcA)$ (as well as the bounded versions to which the shift operation restricts) of any abelian category $\mcA$ is triangulated,
 where the shift functor is given by shifting,
 and distinguished triangles are defined by taking all triangles that are isomorphic to the induced triangle
 (by the snake lemma) of a short exact sequences of complexes.
	The homotopy category $K(\mcA)$ is also triangulated, but here the distinguished triangles are more complicated.
\end{lemma}
\begin{proof} See \cite[III.3.4, cor.~7 and prop.~8 of~IV.2]{Gelfand2003}.
Another good reference is \cite[prop.~42, 44]{Murfet2006b}
\end{proof}

\begin{lemma}[Short exact sequence produces long exact sequence]\label{lem:short_to_long}\uses{def:derived}
If they exist, derived functors of additive functors are triangulated and taking their $i$'th homology is a cohomological functor, which means the following.

Given any additive functor $F\colon\mcA\to\mcB$ in an abelian category with enough projectives (and Cartan-Eilenberg deformation),
a short exact sequence $0\to A_\bullet\to B_\bullet\to C_\bullet\to 0$ of complexes in nonnegative homological degree induces a long exact sequence
\begin{center}\begin{tikzcd}
	& \dots & {H_2(LFC)} \\
	{H_1(LFA)} & {H_1(LFB)} & {H_1(LFC)} \\
	{H_0(LFA)} & {H_0(LFB)} & {H_0(LFC)} \\
	0
	\arrow[from=1-2, to=1-3]
	\arrow[from=1-3, to=2-1]
	\arrow[from=2-1, to=2-2]
	\arrow[from=2-2, to=2-3]
	\arrow[from=2-3, to=3-1]
	\arrow[from=3-1, to=3-2]
	\arrow[from=3-2, to=3-3]
	\arrow[from=3-3, to=4-1]
\end{tikzcd}\end{center}
If the functor is left exact, then $H_0(LF X)=F(X)$ for all objects $X\in \mcA$.
\end{lemma}
\begin{proof}This is \cite[2.3.1]{Riehl2014} and also stated in \cite[13.3.1]{kashiwara2005categories}.
	\end{proof}

\begin{remark}
  We  believe that an analogous statement holds for the unbounded derived category,
  at least in the case that $F$ is right exact,
  see~\ref{thm:calculating_derived_functors_via_truncations}.
  However,
  we could not find a reference nor a simple proof.
\end{remark}

Next, we will approach the relation of the unbounded derived category with its bounded versions. 
	For the following definition see, e.g., \cite[12.3.1]{kashiwara2005categories}.
\begin{definition}[Truncations]\label{def:truncations}\uses{def:complex}
For any complex $A^\bullet$ and $k\in \Z$ define the
\begin{itemize}
		\item \idx{stupid truncations} $\sigma^{\ge n}A^\bullet$ as the subcomplex

\begin{center}\begin{tikzcd}
	{\sigma^{\ge n}A^\bullet} & \dots & 0 & {A^n} & {A^{n+1}} \\
	{A^\bullet} & \dots & {A^{n-1}} & {A^n} & {A^{n+1}}
	\arrow[hook, from=1-1, to=2-1]
	\arrow[from=1-2, to=1-3]
	\arrow[from=1-3, to=1-4]
	\arrow[from=1-3, to=2-3]
	\arrow[from=1-4, to=1-5]
	\arrow[from=1-4, to=2-4]
	\arrow[from=1-5, to=2-5]
	\arrow[from=2-2, to=2-3]
	\arrow[from=2-3, to=2-4]
	\arrow[from=2-4, to=2-5]
\end{tikzcd}\end{center}
		and $\sigma^{\le n}A^\bullet$ via
\begin{center}
\begin{tikzcd}
	{\sigma^{\le n}A^\bullet} & \dots & {A^{n-1}} & {A^n} & 0 & \dots \\
	{A^\bullet} & \dots & {A^{n-1}} & {A^n} & {A^{n+1}} & \dots
	\arrow[from=1-2, to=1-3]
	\arrow[from=1-3, to=1-4]
	\arrow[from=1-4, to=1-5]
	\arrow[from=1-5, to=1-6]
	\arrow[from=2-1, to=1-1]
	\arrow[from=2-2, to=2-3]
	\arrow[from=2-3, to=1-3]
	\arrow[from=2-3, to=2-4]
	\arrow[from=2-4, to=1-4]
	\arrow[from=2-4, to=2-5]
	\arrow[from=2-5, to=1-5]
	\arrow[from=2-5, to=2-6]
\end{tikzcd}\end{center}
		\item define the \idx{Postnikov truncations} (or \idx{canonical truncations}) $\tau^{\ge n}$

\begin{center}
\begin{tikzcd}
	{\tau^{\ge n}A^\bullet} & \dots & 0 & {\coker(d)} & {A^{n+1}} & \dots \\
	{A^\bullet} & \dots & {A^{n-1}} & {A^n} & {A^{n+1}} & \dots
	\arrow[from=1-2, to=1-3]
	\arrow[from=1-3, to=1-4]
	\arrow[from=1-4, to=1-5]
	\arrow[from=1-5, to=1-6]
	\arrow[from=2-1, to=1-1]
	\arrow[from=2-2, to=2-3]
	\arrow[from=2-3, to=1-3]
	\arrow["d"', from=2-3, to=2-4]
	\arrow[from=2-4, to=1-4]
	\arrow[from=2-4, to=2-5]
	\arrow[from=2-5, to=1-5]
	\arrow[from=2-5, to=2-6]
\end{tikzcd}\end{center}
and $\tau^{\le n}A^\bullet$ via
\begin{center}
\begin{tikzcd}
	{\tau^{\le n}A^\bullet} & \dots & {A^{n-1}} & {\ker(d)} & 0 & \dots \\
	{A^\bullet} & \dots & {A^{n-1}} & {A^n} & {A^{n+1}} & \dots
	\arrow[from=1-1, to=2-1]
	\arrow[from=1-2, to=1-3]
	\arrow[from=1-3, to=1-4]
	\arrow[from=1-3, to=2-3]
	\arrow[from=1-4, to=1-5]
	\arrow[from=1-5, to=1-6]
	\arrow[from=2-2, to=2-3]
	\arrow[from=2-3, to=2-4]
	\arrow[from=2-4, to=1-4]
	\arrow["d", from=2-4, to=2-5]
	\arrow[from=2-5, to=1-5]
	\arrow[from=2-5, to=2-6]
\end{tikzcd}\end{center}
\end{itemize}
\end{definition}
\begin{lemma}[Complex as limit of truncations]\label{lem:complex_via_truncations}\uses{def:truncations}
The Postnikov truncation $\tau^{\le n}$ is the left adjoint to the inclusion $D^{\le n}(\mcA)\to D(\mcA)$.
	The Postnikov truncation $\tau^{\ge n}$ is the right adjoint to the inclusion $D^{\ge n}(\mcA)\to D(\mcA)$.
	One furthermore has
\[A^\bullet=\varprojlim_{n\to \infty}\sigma^{\le n}A^\bullet\]
\[A^\bullet=\varinjlim_{n\to \infty}\tau^{\le n}A^\bullet\]
\[A^\bullet=\varinjlim_{n\to -\infty}\sigma^{\ge n}A^\bullet\]
\[A^\bullet=\varprojlim_{n\to -\infty}\tau^{\ge n}A^\bullet\]
\end{lemma}
\begin{proof}
	The statements about the adjunctions of Postnikov truncations are, e.g., \cite[1.2.3.4]{Lurie2017}.
	
	The limits of the stupid and the Postnikov truncations can be found in \cite{Murfet2006b}, and the adjunctions with the Postnikov as well.
	There,
	prop.~33 is also interesting.
\end{proof}

\begin{question}
  Do the stupid truncations give the resp.\ right/left adjoints?
\end{question}
See, e.g., in \cite[sec.~4.5]{Murfet2006b} for the following reduction of \enquote{bounded} complexes to complexes concentrated in one degree.
\begin{lemma}
For any abelian category $\mcA$ consider a class of objects $\mcP\sub \mcA$ containing $0$. 
Now, every triangulated subcategory $\mcS\sub D(\mcA)$ containing $P[i]$ for all $P\in \mcP$ and $i\in \Z$ contains all bounded complexes of objects in $\mcP$.
\end{lemma}
\begin{proof}
This is \cite[lemma~79]{Murfet2006b}.
 The main trick of the decomposition is given by  looking at the mapping cone of the morphism
\begin{center}\begin{tikzcd}
	0 & {P^0} & 0 & 0 & 0 \\
	0 & {P^1} & {P^2} & {P^3} & {P^4} & 0.
	\arrow[from=1-1, to=1-2]
	\arrow[from=1-2, to=1-3]
	\arrow[from=1-2, to=2-2]
	\arrow[from=1-3, to=1-4]
	\arrow[from=1-3, to=2-3]
	\arrow[from=1-4, to=1-5]
	\arrow[from=1-4, to=2-4]
	\arrow[from=1-5, to=2-5]
	\arrow[from=2-1, to=2-2]
	\arrow[from=2-2, to=2-3]
	\arrow[from=2-3, to=2-4]
	\arrow[from=2-4, to=2-5]
	\arrow[from=2-5, to=2-6]
\end{tikzcd}\end{center}

\end{proof}
\begin{corollary}
If $\mcA$ admits enough projectives, then a (under the induced triangulated category) generating class of projectives in $D(\mcA)$ is given by $P[i]$ for projectives $P\in\mcA$ and $i\in \Z$.
\end{corollary}

\begin{conjecture}[Calculating derived functors via truncations]\label{thm:calculating_derived_functors_via_truncations}\uses{def:truncations, def:derived_functor, def:derived}
Consider any cocontinuous functor $F\colon \mcA\to\mcB$ between abelian categories with enough projectives and its point-set left derived functor $LF\colon D^-(\mcA)\to D^-(\mcB)$.
	Then
\[LF(A^\bullet)=\varinjlim LF(\tau^{\le n}A^\bullet)\]
defines the point-set left derived functor.
\end{conjecture}
\begin{proof}[Justification]
We did not find an explicit reference for this. 
However, it should follow by the construction of the projective resolution of unbounded complexes, which is given by the filtered colimit of resolutions of the truncations, 
as well as the fact that the functor commutes with filtered colimits on level of chain complexes.
\end{proof}
\begin{remark}
This seems to be related to \cite[lemma 77]{Murfet2006b}.
\end{remark}

Another frequently useful way to compute derived functors is via spectral sequences.
\begin{lemma}[Grothendieck spectral sequence]\label{lem:groth_spec_seq}\uses{def:spectral_sequence}
	Let $\mcA$, $\mcB, \mcC$ be abelian categories with enough projectives, $G\colon\mcA\to\mcB$ and $F\colon\mcB\to\mcC$ right exact functors
	, and assume that
	all objects $GA$ for projective $A\in \mcA$ are $F$-\idx{acyclic}, meaning that $LF(GA[0])$ has homology concentrated in degree $0$.
	Then there is a (upper left quadrant) spectral sequence with second page
	\[E^2_{j,i}=H^j(LF(H^iLG(A)))\implies H_{j+i}L(FG).\]
	\end{lemma}
	\begin{proof}
	  This is a nice proof, using a Cartan-Eilenberg resolution of $G$ of a projective resolution of $A$.
	  Then use both directions of the corresponding spectral sequences.
	\end{proof}
	
\begin{definition}[Some examples of derived functors]\label{def:examples_of_derived_functors}\uses{def:derived_functor, def:derived}
Consider any abelian category $\mcA$ with enough (compact) projectives.
\begin{itemize}
	\item The most important derived functor is the \idx{derived hom}, $\RHom\colon D(\mcA)^{\op}\times D(\mcA)\to D(\Ab)$.
This is pointwise the right derived functor of the left exact functor $\hom(-,X)\colon \mcA^{\op}\to \Ab$. 

	For objects $A,B\in \mcA$ define the \idx{homology groups of $A$ with coefficients in $B$} as
	\[\Ext^i(A,B):=H^i(A,B):=H^i(\RHom(A[0],B[0]))\]
	These groups are called the $i$-th extensions of $B$ by $A$.

	And more generally for complexes $A,B\in D(\mcA)$ define the extensions
	\[\Ext^i(A,B)=H^i(\RHom(A,B))=\hom_{D(\mcA)}(A,B[i])\]
	\item If the abelian category is symmetric monoidal, the \idx{derived tensor product} $\otimes^L\colon D(\mcA)\times D(\mcA)\to D(\mcA)$ on $D(\mcA)$ is the left derived
	\[\otimes^L=L\otimes\]
	Note that this installs again a symmetric monoidal structure on $D(\mcA)$ (see~\cite[13.4]{kashiwara2005categories})
  \item If the abelian category is closed symmetric monoidal, there is an \idx{internal derived hom} $R\ihom\colon D(\mcA)^{op}\times D(\mcA)\to D(\mcA)$,
		defined as the
	right derived of $\ihom\colon\mcA^\op\times\mcA\to \mcA$.
	 Note that as the derived functors are given by absolute Kan extensions, they preserve adjunctions, and hence $\otimes^L\dashv R\ihom$.
\end{itemize}
\end{definition}
The following general description of Ext groups explains the name extensions, see, e.g., \cite[3.4.6]{Weibel1994} or, for the case $i=1$, \cite[III.5.2]{Gelfand2003}.
\begin{lemma}\label{lem:ext_as_extensions}
The abelian groups $\Ext^i(A,B)=\RHom(A[0],B[i])$ may be described as equivalence classes of exact sequences starting with $B$, having $i$ many terms in between, and ending with $A$.
The equivalence relation for $i=1$ is given by isomorphisms of short exact sequences.
\end{lemma}
\begin{remark}
	The functor $\hom$ usually is regarded as a bifunctor, and there is a whole theory of derived bifunctors.
	We refer to, e.g., \cite[13.4]{kashiwara2005categories} for more.
	Also note that in the $\infty$-categorical setting (which we introduce next), this becomes easier, since there, $\mcD(A\times B)=\mcD(A)\times \mcD(B)$.
\end{remark}

\section{Infinity categories}
\quot{The way I think of it these days, the mysterious way quantum theory slammed into physics in the early 20th century was just nature's way of telling us we'd better learn n-category theory.}
{John Baez, The Tale of n-Categories, Week77, 1996}

In this chapter, we introduce some basic notions of infinity category theory, which allows us to define a more natural notion of derived category and much more.
It seems that many experts prefer to work in the $\infty$-categorical setting, and that sooner or later one should switch to $\infty$-categories,
as there the theorems and definitions seem to be \enquote{the correct ones}.

See also \cite{Bergner2009, Hoermann2019} for a gentle introduction to $(\infty,1)$-categories, \cite{Lurie2024, SAG, Lurie2017, Lurie2009, mair2021animated, riehl2023type, Riehl_Verity_2022, Cisinski2025} or Kerodon
\footnote{\url{https://kerodon.net/}}.

The main idea of $\infty$-categories is to replace ordinary categories with categories that allow morphisms between morphisms, afterwards morphisms between those etc.
There are several ways to introduce $\infty$-categories (by which we always mean $(\infty,1)$-categories), which are all equivalent in the relevant sense.
We stick to the intuitive and, crucially, furthest developed version via quasicategories and simplicial sets.

\begin{warning}
  Note that here at latest, it seems reasonable to work with a sufficient amount of universes (ZFCU) to get rid of set-theoretic problems.
We hence will mostly ignore the set theoretic issues in this section.
See \cite[section 1.2.15]{Lurie2009} for many possible solutions.
\end{warning}
\begin{definition}[Simiplicial sets]\label{def:simplicial_sets}\uses{def:category}\chapthree
Define $\Delta$ to be the category whose objects are finite nonempty ordinals, $[n]=\{0, 1,\dots, n\}$, and whose objects are non-decreasing maps.
    The category $\sSet$ of \idx{simplicial sets} is defined as $\PSh(\Delta)$ (or, equivalently, the free cocompletion of $\Delta$).
More generally, for any category $\mcC$, define the \idx{simplicial objects} in $\mcC$ as $\Fun(\Delta^\op,\mcC)$.

The \idx{standard $n$-simplex} $\Delta_n$ is the presheaf represented by $[n]$, i.e., the simplicial set
    \[[m]\mapsto \hom_{\Delta}([m],[n]).\]
The \idx{standard $k$'th horn} $\Lambda^k_n$ is the simplicial subset of $\Delta_n$ defined by \enquote{removing the face opposite to the $k$'th vertex} (and, consequently, the interior).
    Explicitly, for $0\le i\le n$, define the \idx{elementary face operators} as the unique injective map avoiding $i$,
\[\delta^i_n\colon [n-1]\to [n], \quad j\mapsto \begin{cases}j, &j<i,\\ j+1,&j\ge i.\end{cases} \]
Via Yoneda, these give maps $\delta^i\colon \Delta_{n-1}\to \Delta_n$.
    Now, define the $k$'th horn do be the pushout of all faces except the $k$'th
    \[\Lambda^k_n=\bigcup_{i\ne k}\delta^{i}_{n}(\Delta_{n-1})=\Delta_{n-1}\sqcup_{\Delta_{n-2}}\cdots\sqcup_{\Delta_{n-2}}\Delta_{n-1}.\]
    An \idx{inner horn} is a horn $\Lambda^k_n$ with $0<k<n$.
\end{definition}
\begin{example}
    Take any simplicial set $S$.
    \begin{itemize}
        \item One can view the set $S([0])$ as the \idx{0-morphisms} or objects of the simplicial set, which we draw as points.

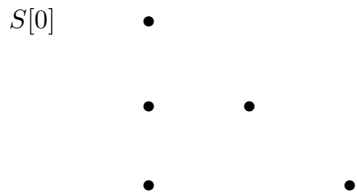
\begin{center}\begin{tikzcd}
	{S[0]} & \bullet \\
	& \bullet & \bullet \\
	& \bullet && \bullet
\end{tikzcd}\end{center}
        \item The set $S([1])$ may be interpreted as the \idx{1-morphisms} or simply morphisms of the simplicial set. 
        The induced arrow of $[0]\to [1], \, 0\mapsto 0$ assigns to every arrow a domain $S([1])\to S([0])$
        and the induced arrow of $[0]\to [1],\, 0\mapsto 1$ assigns to every arrow a codomain.
        Hence we can draw them as arrows, and obtain a picture

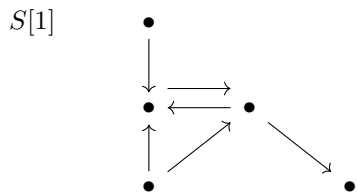
\begin{center}\begin{tikzcd}
	{S[1]} & \bullet \\
	& \bullet & \bullet \\
	& \bullet && \bullet
	\arrow[from=1-2, to=2-2]
	\arrow[shift left=3, from=2-2, to=2-3]
	\arrow[from=2-3, to=2-2]
	\arrow[from=2-3, to=3-4]
	\arrow[from=3-2, to=2-2]
	\arrow[from=3-2, to=2-3]
\end{tikzcd}\end{center}
        \item The value $S[2]$ can be described as \idx{commutative triangles}: There are three morphisms $[0]\to [2]$.
         But to every such triangle there are also three corresponding arrows: $\delta^i_2\colon [1]\to [2]$.
        Leaving the details aside, this boils down to $S[2]$ being a set of triangles between \enquote{composable morphisms}.
         We draw this as natural isomorphism; drawing a surface would be better.

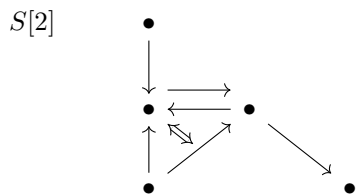
\begin{center}\begin{tikzcd}
	{S[2]} & \bullet \\
	& \bullet & \bullet \\
	& \bullet && \bullet
	\arrow[from=1-2, to=2-2]
	\arrow[shift left=3, from=2-2, to=2-3]
	\arrow[from=2-3, to=2-2]
	\arrow[from=2-3, to=3-4]
	\arrow[from=3-2, to=2-2]
	\arrow[""{name=0, anchor=center, inner sep=0}, from=3-2, to=2-3]
	\arrow[shorten <=2pt, Rightarrow, 2tail reversed, from=0, to=2-2]
\end{tikzcd}\end{center}
        Informally, elements of $S[2]$ are statements of the form \enquote{$h$ is a possible composite of $f$ and $g$}.
      \item The elements of $S[3]$, or 3-isomorphisms, are pyramids (3-simplices), filling out a space between four \enquote{composable} triangles,
            and telling us which triangles (2-morphisms) compose to essentially the same 2-isomorphisms.
        \item $\dots$
    \end{itemize}
    We see, that we obtain a \idx{space of morphisms} (in the sense of an actual CW-complex from topology, the so-called geometric realization), describing the nature of our simplicial set.
    Compare this view with the classical view on 1-categories as graphs and 2-categories as graphs with 2-morphisms.
\end{example}
The main idea now is the following: While composition in a 1-category is strict in the sense that there exists only one possible composite,
in an $\infty$-category we require composition not to be unique, but rather unique up to $2$-isomorphisms (these are the elements $S[2]$), and these $2$-isomorphisms should be unique up to $3$-isomorphisms etc.

But what we still want to require is the possibility to compose morphisms with appropriate (co)domains (although this may not be unique).
\begin{definition}[Infinity categories]\label{def:inf_cat}\uses{def:category, def:simplicial_sets}\chapthree
An \idx{$\infty$-category} (other names  for this include \idx{quasi-category} or \idx{weak Kan complex}) is a simplicial set $S$ such that for every inner horn
    $\Lambda^k_n$ in $S$ there exists a filler, meaning
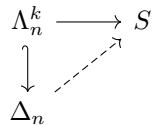
\begin{center}\begin{tikzcd}
	{\Lambda^k_n} & S \\
	{\Delta_n}
	\arrow[from=1-1, to=1-2]
	\arrow[hook', from=1-1, to=2-1]
	\arrow[dashed, from=2-1, to=1-2]
\end{tikzcd}\end{center}

\end{definition}
    This precisely encapsulates the idea of being able to compose all morphisms that should be composable.
 For example, consider a horn $\Lambda_2^1\to S$.
This, by Yoneda, consists of two morphisms $a, b\in S[1]$ with the codomain of the first being the domain of the second morphism,

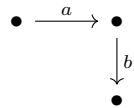
\begin{center}\begin{tikzcd}
	\bullet & \bullet \\
	& \bullet
	\arrow["a", from=1-1, to=1-2]
	\arrow["b", from=1-2, to=2-2]
\end{tikzcd}\end{center}
This admitting a filler means that there is a morphism $c\in S[1]$ with the same domain as $a$ and the same codomain as $b$ and a triangle $t\in S[2]$ witnessing the \enquote{commutativity}
of the following diagram.

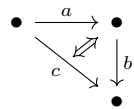
\begin{center}\begin{tikzcd}
	\bullet & \bullet \\
	& \bullet
	\arrow["a", from=1-1, to=1-2]
	\arrow[""{name=0, anchor=center, inner sep=0}, "c"', from=1-1, to=2-2]
	\arrow["b", from=1-2, to=2-2]
	\arrow[shorten <=2pt, Rightarrow, 2tail reversed, from=0, to=1-2]
\end{tikzcd}\end{center}
Playing around with the same notion for 3-horns gives a quite good feeling for the horn filling condition:
Given three maps $\bullet\stackrel{a}{\to}\bullet\stackrel{b}{\to}\bullet\stackrel{c}{\to}\bullet$,
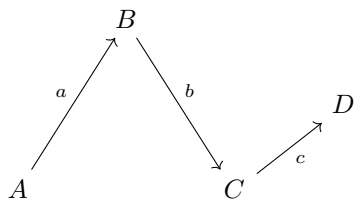
\begin{center}
\begin{tikzcd}
	& B \\
	&&& D \\
	A && C
	\arrow["b", from=1-2, to=3-3]
	\arrow["a", from=3-1, to=1-2]
	\arrow["c"', from=3-3, to=2-4]
\end{tikzcd}
\end{center}

we can first infer the existence of composites $f$ of $a$ with $b$ and $g$ of $b$ with $c$ by the horn filling property just discussed
together with 2-isomorphisms witnessing these commutativities.
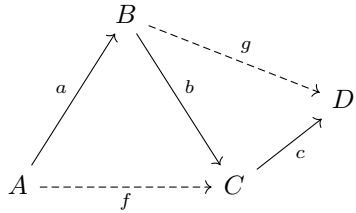
\begin{center}
\begin{tikzcd}
	& B \\
	&&& D \\
	A && C
	\arrow["g", dashed, from=1-2, to=2-4]
	\arrow["b", from=1-2, to=3-3]
	\arrow["a", from=3-1, to=1-2]
	\arrow["f"', dashed, from=3-1, to=3-3]
	\arrow["c"', from=3-3, to=2-4]
\end{tikzcd}
\end{center}

This is a horn of the form $\Lambda_{3}^{2}$.
A filling for such a horn consists of a simultaneous composite $h$ of $a$ and $g$ as well as $c$ and $f$ (together with appropriate 2-isos)
and an interior for the emerging pyramid (tetrahedron).

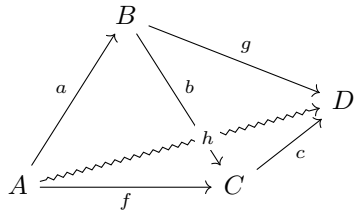
\begin{center}
\begin{tikzcd}
	& B \\
	&&& D \\
	A && C
	\arrow["g", from=1-2, to=2-4]
	\arrow["b", from=1-2, to=3-3]
	\arrow["a", from=3-1, to=1-2]
	\arrow["h"{description, pos=0.6}, squiggly, from=3-1, to=2-4]
	\arrow["f"', from=3-1, to=3-3]
	\arrow["c"', from=3-3, to=2-4]
\end{tikzcd}
\end{center}

This filled interior witnesses the associativity, i.e., the higher coherence between the two ways of composing the three arrows $a,b,c$.

Similarly, the inner horn filling property for $\Lambda_{n}^{i}$ gives coherent ways to compose $n$ many arrows.

We note that it is not reasonable to require horn filling for \textbf{all} (not just inner) horns, as e.g., the horn filling condition for $\Lambda_2^0$ would look like

\begin{center}\begin{tikzcd}
	\bullet & \bullet \\
	& \bullet
	\arrow["a", from=1-1, to=1-2]
	\arrow[""{name=0, anchor=center, inner sep=0}, "c"', from=1-1, to=2-2]
	\arrow["b", dashed, from=1-2, to=2-2]
	\arrow[shorten <=2pt, Rightarrow, dashed, 2tail reversed, from=0, to=1-2]
\end{tikzcd}\end{center}
\begin{remark}
  Note that the idea is to have all higher morphisms invertible (in order to reasonably draw them as surfaces without \enquote{direction}),
  unlike the usual 2-categorical thought of 2-morphisms as directed arrows between arrows.
Thus, the name $(\infty,1)$-category is more precise (the $1$ indicates the order such that all morphisms of higher order are required to be invertible).
\end{remark}

Having a new concept of categories, let us embed the classical concept of 1-categories into this.
\begin{definition}[(Simplicial) nerve of a category]\label{def:nerve}\uses{def:inf_cat}
  For any classical (small) category $\mcC$ define the \idx{nerve} $\mcN(\mcC)$ of $\mcC$ by the yoneda embedding $\Cat\to [\Cat^{\op},\Set]$ followed by restriction
  $[\Cat^{\op},\Set]\to[\Delta^{\op},\Set]$ to the family of categories $\Delta$.
Explicitly, this is given as the $\infty$-category induced by the
 mapping \[[n]\mapsto \{ (f_1,\dots, f_n)\,\text{composable}\}=\mathrm{Fun}([n]\to \mcC).\]
    This yields a fully faithful embedding $\Cat\to \Cat_\infty$ (to be defined later).
\end{definition}
\begin{remark}
    Nerves of ordinary categories are characterised as those $\infty$-categories that admit a \textbf{unique} inner horn filling, see \cite[1.1.2.2.]{Lurie2009}.
\end{remark}

\begin{definition}[Homotopy category of infinity category]\label{def:hom_cat_of_inf_cat}
To any $\infty$-category $\mcC$ (in fact, any simplicial set) there is an associated $1$-category $\mathrm{h}\mcC$, called the \idx{homotopy category}
such that the nerve functor $\Cat\to \sSet$ is right adjoint to the homotopy category functor $\sSet\to \Cat$.
\end{definition}
Many explicit constructions of the homotopy category can be found in \cite{Lurie2009}, after 1.1.4.3.
For the adjunction,
see \cite[1.2.3.1]{Lurie2009}.
\begin{lemma}
  Infinity categories $\mcC$ for which $\mathrm{h}\mcC$ is a groupoid are precisely those for which the Kan lifting condition is fulfilled for all horns (not just inner horns).
    These categories are called \idx{anima} (or \idx{$\infty$-grupoids} or \idx{spaces} or \idx{Kan complexes}).

    These ae the objects the category $\Ani$ as full subcategory of $\Cat_\infty$ (to be defined later).
Moreover, the inclusion admits a right adjoint, and, in a higher categorical sense, a left adjoint.
\end{lemma}
\begin{proof}
      See \cite[1.2.5.2--6]{Lurie2009}.
\end{proof}

\begin{warning}[Basic definitions in $\infty$-category theory]\label{def:basic_infty}\uses{def:inf_cat}
      Any concept we have developed in the first and third chapter for ordinary categories has an appropriate $\infty$-categorical version.
    We will thus just use them like in ordinary category theory without giving formal definitions in this text.
\end{warning}
\subsubsection{Some advanced constructions on $\infty$-categories}
\quot{This general nonsense, you see, to make it a mathematical statement would be kind of awkward. [\dots]
  So I don't think I will try to make it really mathematics.
  What is important is to get the general idea and have the correct -- how shall I say -- thought patterns for these situations so as to not use time and energy doing involved things.}
{A.Grothendieck, 1973}

In this subsection, we will try to collect some advanced constructions from \cite{Lurie2009, Lurie2017, SAG} which seem to be related to our theory.
Note that we are just at the beginning of understanding this theory, and hence this chapter should be read with caution. 
We recommend to read the corresponding sections of the references.

\begin{example}
  The functors from one $\infty$-category to another form an $\infty$-category $\Fun(\mcC,\mcD)$ (1.2.7.3 in \cite{Lurie2009}).

  The $\infty$-category of all $\infty$-categories is the nerve of the ordinary category obtained by taking the hom-sets to be the largest Kan complex contained in $\Fun(\mcC,\mcD)$
  (the \enquote{largest Kan komplex} is the right adjoint to $\Ani\to \Cat_\infty$).
  This is up to set theory (the range is over all $\infty$-categories in a smaller universe), and usually one would, just as for $\Cat$, regard this as an $\infty$-bicategory, see \cite[3.0.0.4]{Lurie2009}.
  Objects of $\Cat_{\infty}$ are $\infty$-categories, morphisms are functors, 2-morphisms are equivalences between functors, etc.

  The $\infty$-category $\Cat_\infty$ is bicomplete in the appropriate sense, \cite[3.3.3--4]{Lurie2009}.
\end{example}

Next, we explain the appearance of chain complexes in homological algebra.
\begin{theorem}[Dold-Kan correspondence]\label{thm:dold_kan_correspondence}\uses{def:simplicial_sets, def:chain_complex}
Let $\mcA$ be an abelian category.
Then there is an equivalence of categories of negatively graded chain complexes in $\mcA$ and simplicial objects in $\mcA$,
    \[\Ch^{\le 0}(\mcA)\simeq \Fun(\Delta^\op,\mcA).\]
\end{theorem}
\begin{proof} See, e.g., \cite[1.2.3]{Lurie2017}.
    \end{proof}
This shows that the appearance of chain complexes is totally natural from an $\infty$-categorical perspective.
We will explore this further by first defining the appropriate $\infty$-categorical version categories coming up in abelian/derived/triangulated category theory from \cite[1.1.1.9]{Lurie2017}.
\begin{definition}[Stable $\infty$-cat]\label{def:stable_infty}\uses{def:inf_cat}
    An $\infty$-category is called \idx{stable $\infty$-category}, if
\begin{itemize}
    \item it admits a $0$-object,
    \item for every morphism $g\colon X\to Y$ there exists a \idx{fiber}, which is a pullback square (in the $\infty$-categorical sense)

\begin{center}\begin{tikzcd}
	Z & X \\
	0 & Y
	\arrow[from=1-1, to=1-2]
	\arrow[from=1-1, to=2-1]
	\arrow["g", from=1-2, to=2-2]
	\arrow[from=2-1, to=2-2]
\end{tikzcd}\end{center}
    and a \idx{cofiber}, which is a pushout square

\begin{center}\begin{tikzcd}
	X & Y \\
	0 & Z
	\arrow["g", from=1-1, to=1-2]
	\arrow[from=1-1, to=2-1]
	\arrow[from=1-2, to=2-2]
	\arrow[from=2-1, to=2-2]
\end{tikzcd}\end{center}
    \item a \enquote{triangle}
\begin{center}\begin{tikzcd}
	X & Y \\
	0 & Z
	\arrow[from=1-1, to=1-2]
	\arrow[from=1-1, to=2-1]
	\arrow[from=1-2, to=2-2]
	\arrow[from=2-1, to=2-2]
\end{tikzcd}\end{center}
is a pushout diagram precisely if it is a pullback diagram.
\end{itemize}

    The second requirement is the analogon of images being isomorphic to coimages.

    Stable $\infty$-categories are finitely bicomplete (\cite[1.1.3.4]{Lurie2017}).
    An \idx{exact} functor between stable $\infty$-categories is a finitely limit preserving functor
    (or equivalently(!) finitely colimit preserving or finitely bicontinuous, see \cite[1.1.4.1]{Lurie2017})
    \end{definition}

A general passage from $1$-categories to $\infty$-categories can be described via \idx{animation}, being a similar construction to the passage to derived categories.
\begin{definition}[Animation]
Recall the process of 1-animation \ref{def:one_ani}.

    Instead of taking the 1-cocompletion, we can now rather take the $\infty$-categorical cocompletion, leading to an $\infty$-category.
    Define $\Ani(\mcC)$ as the $\infty$-category freely generated by $\mcC^{\mathrm{cp}}$ under sifted colimits (i.e., $\sInd(\mcC^{\mathrm{cp}})$ in the $\infty$-categorical sense).
    Note that sifted colimits in the $\infty$-setting are given by filtered colimits and geometric realisations, the latter meaning colimits over $\Delta^{\op}$.

Again, this agrees with the full $\infty$-subcategory of the $\infty$-functor category $[{\mcC^{\mathrm{cp}}}^{\op}, \Ani(\Set)]$ generated under sifted colimits by representables.
    If $\mcC^{\mathrm{cp}}$ is small, this agrees with all functors mapping finite coproducts in $\mcC^{\mathrm{cp}}$ to products in $\Ani(\Set)$.

    See \cite[11.1]{scholze2019Analytic} and the references therein for more.
  
\end{definition}
\begin{remark}
  There exists an analogous theory of (co)completions of categories in the $\infty$-categorical setting, giving similar theorems to the ones in \ref{ssec:completions}.
  See, e.g., \cite[5.5.8]{Lurie2009} for more.
\end{remark}
\begin{example}
    \begin{itemize}
        \item $\Ani(\Set)=\Ani$ is the $\infty$-category of $\infty$-groupoids.
        \item $\Ani(\Ab)$ is the $\infty$-derived category $\mcD_{\ge 0}(\Ab)$.
\end{itemize}
\end{example}

Here, the $\infty$-derived category is defined as follows.
We recommend reading \cite[1.3.2]{Lurie2017} for this.

\begin{definition}[$\infty$-derived cat]\label{def:inf_der}\uses{def:inf_cat}
    Let $\mcA$ be an abelian category with enough projectives $\mcA_{\mathrm{proj}}$.
    Then the category $\Ch^{-}(\mcA_\mathrm{proj})$ can be regarded as a differential graded category, see \cite[1.3.2.1]{Lurie2017}.
    A \idx{differental graded category} is -- roughly -- a category enriched over chain complexes, for a precise definition see \cite[1.3.1.1]{Lurie2017}.
    For any such category, one can define an $\infty$-category called the \idx{differential graded nerve} $\mcN_{dg}(\mcC)$, see \cite[1.3.1.6]{Lurie2017}.
    The objects of $\mcN_{dg}(\mcC)$ are the objects of $\mcC$, and the 1-morphisms are morphisms in $\mcC$ (with $df=0$), see \cite[1.3.1.8]{Lurie2017}.
    Now, define the \idx{$\infty$-derived ($\infty$-)category} of $\mcA$ as $\mcN_{dg}(\Ch^-(\mcA_{\mathrm{proj}}))$.

    The $\infty$-derived category is stable, see \cite[I.13.10]{Lurie2024} or \cite[1.3.2.18]{Lurie2017}.
    Furthermore, one can characterise this category by a universal property similar to the usual definition of derived categories, see \cite[1.3.3.2]{Lurie2017}.

    In fact, by \cite[1.3.4.4]{Lurie2017}, we can define the $\infty$-derived category $D^-(\mcA)$ as the localisation (in $\infty$-categories) of $\Ch^-(\mcA)$ by quasiisomorphisms.

   By \cite[11.3]{Hoermann2019}, this agrees with the $\{\Delta^{\op}\}$-continuous (a more involved version of reflexive coequalizers) cocompletion (in $\infty$-categories) of the projective objects in $\mcA$.
\end{definition}
The $\infty$-derived category roughly consists of the following data (see \cite[1.3.2]{Lurie2017}):
\begin{itemize}\item Objects are right-bounded complexes of projective objects.
    \item Morphisms are maps of chain complexes.
    \item 2-morphisms are chain homotopies.
    \item $\dots$ 
    \end{itemize}

\begin{theorem}[Comparison with classical derived category]\label{thm:inf_der_to_der}\uses{def:inf_der}
   The homotopy category of the $\infty$-derived category agrees with the usual derived category, see \cite[1.3.2.999999999]{Lurie2017}.

   \[\mathrm{h}\mcD^-(\mcA)\simeq D^-(\mcA)\]

\end{theorem}
For more on this see \cite[1.3]{Lurie2017}.
\begin{lemma}[$\infty$-derived commutes with diagrams]\label{lem:inf_der_commutes}\uses{def:inf_der}
    For any category $\mcC$ and small category $\mcI$, one has
    \[\mcD(\mcC^{\mcI})\simeq \mcD(\mcC)^{\mcI}.\]
    \end{lemma}
    \begin{proof}
      Without proof, this is \cite[9.2]{Hoermann2019}.
    \end{proof}

Next, we want to see how to recover the whole derived category from the right-bounded derived category. This is done via spectrum objects, see, e.g., \cite[1.4]{Lurie2017}.

\begin{definition}[Spectrum objects]\label{def:spec_obj}\uses{def:inf_cat}
   In any finitely bicomplete $\infty$-category $\mcC$ with $0$-object, one can define the \idx{suspension} of any object $X$ to be the pushout

\begin{center}\begin{tikzcd}
	X & 0 \\
	0 & {\Sigma X}
	\arrow[from=1-1, to=1-2]
	\arrow[from=1-1, to=2-1]
	\arrow[from=1-2, to=2-2]
	\arrow[from=2-1, to=2-2]
\end{tikzcd}\end{center}
and the \idx{loop space} as the pullback

\begin{center}\begin{tikzcd}
	{\Omega X} & 0 \\
	0 & X
	\arrow[from=1-1, to=1-2]
	\arrow[from=1-1, to=2-1]
	\arrow[from=1-2, to=2-2]
	\arrow[from=2-1, to=2-2]
\end{tikzcd}.\end{center}
    Clearly, this is functorial, yielding the functors $\Sigma$ and $\Omega$.
    Furthermore these two functors are adjoint, $\Sigma\dashv \Omega$.
    Now, the $\infty$-category is stable precisely if $\Sigma$ and $\Omega$ form an equivalence (see \cite[1.4.2.27]{Lurie2017}).

    Define the \idx{stabilisation} $S(\mcC)$ of $\mcC$ as the limit (in $\Cat_\infty$) of

\begin{center}
    \begin{tikzcd}
	\dots & \mcC & \mcC & \mcC
	\arrow["\Omega", from=1-1, to=1-2]
	\arrow["\Omega", from=1-2, to=1-3]
	\arrow["\Omega", from=1-3, to=1-4]
\end{tikzcd}.\end{center}
    The objects in $S(\mcC)$ are called \idx{spectrum objects}.
\end{definition}
\begin{lemma}
    For the category $\mcD_{\ge 0}(\mcA)$ for an abelian category $\mcA$, the functor $\Omega$ is precisely given by right shifting, and $\Sigma$ by shift to the left.
    The stabilisation $S(\mcD)$ agrees with the unbounded $\infty$-derived category $\mcD(\mcA)$ (see \cite[C.1.2.9]{SAG}).
    
    In particular, any complex $X\in \mcD(\mcA)$ is a canonical limit of its Postnikov truncations, and these are just the coordinate projections.
 \end{lemma}
See, e.g., \cite[C.3.6.3]{SAG} or \cite[5.5.6.23]{Lurie2009} for more on this and other constructions.

We have the following very elegant $\infty$-categorical analog of symmetric monoidal categories from \cite[3.3/4]{Scholze2023}.
This is just a very special case of the vast theory of $\infty$-operads, for which refer to \cite[2--4]{Lurie2017}
\begin{definition}[Symmetric monoidal $\infty$-cat]\label{def:symm_mon_infty}\uses{def:inf_cat}
  Define the category $\Fin^{\mathrm{part}}$ as the category whose objects are finite sets and whose morphisms are partially defined maps,
  i.e., maps defined on a subset of the codomain.

A \idx{commutative monoid} in an $\infty$-category $\mcC$ is a functor
$F\colon \Fin^{\mathrm{part}}\to \mcC$ such that for any finite $I$, the morphism
\[F(I)\to \prod_{i\in I}F(\{i\})=X(\ast)^I\]
being induced by all $I\to \{i\}$ mapping $i\to i$ and being defined only there, is an isomorphism.

A \idx{symmetric monoidal $\infty$-category} is a commutative monoid in $\Cat_\infty$.
\end{definition}

\section{Universal algebra}\label{sec:univ-alg}
\quot{
  What the functor is an algebra?
  Is a ring an algebra?
  What if it isn't commutative?
  Does an algebra have to be associative?
  What about Lie algebras?
  Is a group an algebra?
  But then,
  why do we also have algebraic groups?
  There seems to be a variety of algebras,
  but an algebraic variety is not an algebra;
  it's no even a class of algebras ---
  well, unless you're talking about a variety of algebras.
}{Sheafification of G, in Algebra -- It's not what you think it is!%
\footnote{\url{youtu.be/bD8kjpynF6A}}}
So far, we have seen the categorical development of a theory that is suitable to replace the theory of abelian groups or $R$-modules.
However, there certainly are more algebraic structures than those two. For example, Groups, Rings, algebras and many more do not fit in this framework.
So what \emph{is} algebra?

\quot{
  Algebra is the offer made by the devil to the mathematician.
  The devil says:
  \enquote{I will give you this powerful machine, it will answer any question you like.
  All you need to do is give me your soul:
  give up geometry and you will have this marvelous machine.}
}{Sir M. Atiyah, in \cite{Atiyah2001}}

\subsubsection*{Algebraic or Lawvere theories}

For an appropriate answer to this question we take a step back and think about how to put algebraic structure on the objects of a given category.
For example, we will show how to form $\Ab(\mcC)$,
the category of abelian group objects (and homomorphic $\mcC$ maps),
but also $\Grp(\mcC)$, the group objects in $\mcC$, and many more will fit into this picture.
This process of \enquote{algebraisation} will turn out to be surprisingly similar to that of condensation (which will be explained in~\ref{sec:cond_cats} and can be read first)
and the two commute.

So, what is an algebraic structure?
When recapitulating the common undergraduate definitions of (abelian) groups, (commutative) rings and modules,
a pattern emerges.
There is always an \enquote{underlying set} on which a number of different (mostly binary) operations are defined
and these operations are required to satisfy some universally quantified equations,
for example,
$\forall x,y,z\in R\colon x\cdot(y+z)=x\cdot y+x\cdot z$.
A first objection to this perspective might be that many algebraic structures have specified \textit{constants}.
Luckily,
these fit neatly into this perspective:
a constant is simply an operation that takes no inputs (and therefore necessarily always gives the same output),
a so-called \textit{nullary operation} (compare unary, binary, trinary, etc.).
A second objection might be that often one encounters existence quantifiers,
for example, part of the usual definition of an abelian group is $\forall x\in A\colon\exists y\in A\colon x+y=0$.
This can also easily be remedied by requiring corresponding unary operations.
In this case,
we need $A\to A,\,a\mapsto -a$ and the earlier requirement turns into $\forall x\in A\colon x+(-x)=0$,
which is universally quantified.
A third objection could be that whenever there is scalar multiplication involved,
there are actually two underlying sets.
In the case of modules,
there is a ring (which has an underlying set $R$)
and an abelian group (with underlying set $M$).
This is also easily fixed,
but for clearness of exposition we postpone this generalisation.

Guided by these ideas,
the following definition seems plausible.

\begin{definition}[preliminary]
  A type of algebraic structure consists of a set $\Omega$ together with an arity map $\Omega\to\N_{0},\,F\mapsto [F]$
  and a set $\mcE$ of equations of terms formed from the operations in $\Omega$.
\end{definition}

\begin{example}
  In our toy example of abelian groups,
  we have $\Omega=\{+,-,0\}$, $[+]=2$, $[-]=1$, $[0]=0$
  and
  \[
    \mcE=\{X+(Y+Z)=(X+Y)+Z,\,X+Y=Y+X,\,X+(-X)=0\}.
  \]
\end{example}

In order to make this notion precise,
one needs to inductively define what a term is.
We won't bother with this hassle here as we will never formally use this notion.

With this idea in mind,
an actual algebraic structure $A$ (also commonly referred to as a model, an algebra or an instantiation of a type of algebraic structure)
is classically a set $U(A)$ together with maps $F_{A}\colon U(A)^{[F]}\to U(A)$ for each $F\in\Omega$ such that the equations of $\mcE$ are fulfilled for all inputs.
Homomorphisms $A\to B$ of algebraic structures with underlying sets $U(A)$ resp.\ $U(B)$ are maps $f\colon U(A)\to U(B)$ of the underlying sets preserving all structure;
i.e., for all $F\in\Omega$ and $x=(x_{1},\ldots,x_{[F]})\in U(A)^{[F]}$,
$f(F_{A}(x))=F_{B}(f(x_{1}),\ldots,f(x_{[F]}))$,
or more concisely: for all $F\in\Omega$,$f\circ F_{A}=F_{B}\circ f^{[F]}$.
Clearly,
in this way,
for every $(\Omega,\mcE)$,
we obtain a category with objects algebraic structure of type $(\Omega,\mcE)$ and morphisms homomorphisms as just described.

This is how universal algebra has been done for a long time!
Using these definitions,
one can prove things such as the existence of free algebras (left adjoint to the forgetful underlying set functor)
or that isomorphisms of algebras are bijective homomorphisms.

Before we do things like that,
we want to free ourselves of the constraint of only working \enquote{within $\Set$}.
For example,
we want to know what a topological group is.
To achieve this,
we only have to replace the terms \enquote{underlying set} and \enquote{map} in the earlier definition.

\begin{definition}[preliminary]
  Given a type of algebraic structure $(\Omega,\mcE)$ and a category $\mcC$,
  an \textit{$(\Omega,\mcE)$-object $A$ in $\mcC$} consists of an object $U(A)\in\mcC$%
  \footnote{More precisely, an object $U(A)\in\mcC$ for which the appropriate finite powers exist.}
  and for each $F\in\Omega$ a morphism $F_{A}\colon U(A)^{[F]}\to U(A)$
  such that for every equation in $\mcE$, the corresponding diagram commutes.
\end{definition}

\begin{remark}
  The last part of definition is intentionally vague.
  For $F=+$ in the case of abelian groups,
  $X+(Y+Z)=(X+Y)+Z$ turns into the diagram
  \begin{center}
    \begin{tikzcd}
      U(A)^3\ar[rr,"\id_{U(A)}\times F_A"]\ar[d,"F\times\id_{U(A)}"] & & U(A)^2\ar[d,"F_A"]\\
      U(A)^2\ar[rr,"F_A"] & & U(A)
    \end{tikzcd}
  \end{center}
  and $X+0=X$ (and $G=0$, $[G]=0$) becomes
  \begin{center}
    \begin{tikzcd}
      {U(A)^1\times U(A)^0}\ar[rr,"\id_{U(A)}\times G_A"] & & {U(A)^2}\ar[d,"F_A"]\\
      {U(A)^1}\ar[u,"{(\id_{U(A)},test)}"]\ar[rr,"{\id_{U(A)}}"] & & U(A).
    \end{tikzcd}
  \end{center}
  Again,
  we won't formalise this simply because the definition is not in final form.
\end{remark}

\begin{example}
  \begin{itemize}
    \item A topological abelian group is an abelian group object in $\Top$.
    \item A Lie group is nothing but a group object in the category of differentiable manifolds and smooth maps.
    \item An algebraic group is a group object in the category of algebraic varieties (by which we mean special schemes).
    \item A group group is an abelian group and therefore the fundamental group of a topological group is abelian,
          having the funny (and otherwise nontrivial) consequence that the figure eight ($S^{1}\vee S^{1}$) can not be endowed with a (continuous) group structure
          (as such a structure would be transported to $\pi_{1}(S^{1}\vee S^{1})$ as $\pi_{1}$ is finitely product preserving).
  \end{itemize}
  Elliptic curves can also be seen as group objects.
\end{example}

\begin{remark}
  Of course,
  there are still classical types of algebraic structure which are not captured by the sketched formalism of operations and universal equations.
  Most notably,
  there is the case of fields,
  where inversion would have to be a partially defined function.
  Surely,
  there are no complaints about this formalism not capturing structures such as finitely presented modules.
  There is simply no reason to expect this specific phenomenon to fit into such a general framework and a finitely presented module
  is most reasonably viewed as a specific flavour of module,
  which in turn is a sort of structure that we can form.
  If  If on squints a little,
  it is certainly possible to similarly simply see fields as a flavour of rings (such as noetherian rings or integral domains).
\end{remark}

\begin{remark}
  There is one notable conceptual infelicity in our current framework.
  Consider the following type of algebraic structure:
  We take $\Omega=\{*,+,-,0\}$ with $[*]=3$ and the other arities as before
  and
  \[
    \mcE=\{X*Y*Z=(X+Y)+Z,\,X+(Y+Z)=(X+Y)+Z,\,X+Y=Y+X,\,X+(-X)=0\}.
  \]
  Obviously,
  the resulting algebras are going to be abelian groups.
  The operation $*$ is entirely superfluous.
  Therefore,
  we would like to somehow identify these to types with each other.
  However,
  it is not at all clear how one would go about doing this.
  After all,
  it can be arbitrarily difficult to figure out whether or not two types are equivalent in this sense
  (just as figuring out whether or not two groups given by generators and relations are isomorphic is non-computable).
\end{remark}

\begin{definition}[preliminary]
  A morphism of $(\Omega,\mcE)$-objects $A$, $B$ in $\mcC$ is a ($\mcC$-)morphism
  $f\colon U(A)\to U(B)$ such that,
  for each $F\in\Omega$,
  the following diagram commutes:
  \begin{center}
    \begin{tikzcd}
      U(A)^{[F]}\ar[rr,"f^{[F]}"]\ar[d,"F_A"] & & U(B)^{[F]}\ar[d,"F_B"]\\
      U(A)\ar[rr,"f"] && U(B).
    \end{tikzcd}
  \end{center}
\end{definition}

This is extremely reminiscent of the definition of a natural transformation.
And indeed,
algebras for a type of algebraic structures behave a lot like functors.
Clearly,
their target category needs to be $\mcC$.
Figuring out an apt definition of the domain category
(but more importantly, having this realisation at all)
is work originally due to Lawvere,
whence these categories are called \textbf{Lawvere theories}.
In a Lawvere theory, we need a \enquote{generic object} $x$;
we will interpret the evaluation at $x$ as the underlying object $U(A)$.
Then we need finite products of $x$ (including $x^{0}$).
We will interpret the morphisms $x^{n}\to x$ as $n$-ary operations.
(And accordingly,
the morphisms $x^{n}\to x^{m}$ are to be interpreted as $m$-tuples of $n$-ary operations).
In this sense,
we do not differentiate the \enquote{basic} operations (those of $\Omega$ in our current language)
and the terms formed from these.
Though at first unusual,
this is a big advantage as it completely solves the problem of the last remark.

\begin{definition}\label{def:lawvere_theory_ss}
  A \idx{single-sorted Lawvere theory} is a category $\mcT$ with finite products and
  an essentially surjective product-preserving functor $M\colon\Fin^{\op}\to\mcT$.
  We call $x:=M([1])$ the generic object of $\mcT$.
  A $\mcT$-object (or -algebra\footnote{Unfortunately, this terminology is absolutely standard in the literature and we will not attempt to avoid it.})
  in a category $\mcC$ is a (finite) product-preserving functor $A\colon \mcT\to\mcC$.
  We call $U(A)=A(x)$ the underlying object of $A$ (and in practice also denote it by $A$ or $?A$).

  A (homo)morphism of $\mcT$-algebras $A$, $B$ is a natural transformation $f\colon A\to B$
  (which is uniquely determined by the morphism $U(f):=f_{x}\colon U(A)\to U(B)$).

  We obtain a \textbf{category of $\mcT$-algebras in $\mcC$}
  which might reasonably be denoted by $\Fun_{\times}(\mcT,\mcC)$ or $[\mcT,\mcC]_{\times}$.
  This category is also often called $\mcT(\mcC)$.
\end{definition}

Now that we are beginning to give the definitions we will employ to prove theorems,
it is high time we repayed our debt regarding multiple underlying sets of an algebraic structure (such as in the case of modules).
For this,
we now give the first official definition of this section.

\begin{definition}\label{def:lawvere_theory_ms}
  A (multi-sorted) \idx{Lawvere theory} consists of
  \begin{itemize}
	\item a finite set of sorts $S$,
    \item a category $\mcT$ with finite products, and
    \item an essentially surjective product-preserving functor
          \[M\colon (\Fin/S)^{\op}=[S,\Fin^{\op}]=[S,\Fin]^{\op}\to\mcT.\footnote{Note that $[S,\Fin]^{\op}=[S^{\op},\Fin^{\op}]=[S,\Fin^{\op}]$ since $S$ is a discrete category.}\]
  \end{itemize}

  A \idx{$\mcT$-object} (or \idx{$\mcT$-algebra}) in a category $\mcC$ is a (finite) product-preserving functor $A\colon \mcT\to\mcC$.

  For $s\in S$,
  we identify with $s$ the object of $[S,\Fin^{\op}]$ given by $\{*\}$ on $s$ and $\emptyset$ everywhere else.
  We call $\{M(s)\mid s\in S\}$ the generic objects of $\mcT$ (and again denote $M(s)$ by $s$).

  For a $\mcT$-object $A$,
  we call $U_{s}(A)=A(s)$ the underlying objects of $A$ (where $s\in S$).
  A (homo)morphism of $\mcT$-algebras $A$, $B$ is a natural transformation $f\colon A\to B$.
  Note that this natural transformation is uniquely determined by the morphisms $U_{s}(f):=f_{s}\colon U_{s}(A)\to U_{s}(B)$.

  We obtain a \idx{category of $\mcT$-algebras} in $\mcC$, which we usually denote by $\mcT(\mcC)$; other reasonable notations are $\Fun_{\times}(\mcT,\mcC)$ or $[\mcT,\mcC]_{\times}$.
\end{definition}

\begin{remark}
  \begin{enumerate}[(i)]
    \item
  The essential surjectivity of $M$ says that the generic objects of $\mcT$ generate $\mcT$ under finite products.

    \item
  We often assemble the forgetful functors $U_{s}$ for each sort into a single one by currying
  $U\colon\mcT(\mcC)\to[S,\mcC]$.

    \item
  When $\mcT$ is clear from to context,
  write $F\colon x_{1}\times\cdots\times x_{n}\to y$
  to mean $F\in\Hom_{\mcT}(M(\{1,\ldots,n\},x),y)$ (for $x\in S^{n}$, $y\in S$).

  \end{enumerate}
\end{remark}

Breaking this definition down,
a Lawvere theory is a category $\mcT$ generated under finite products by a finite set of objects of $\mcT$,
i.e., its objects can be indexed by $\N_{0}^{k}$.
A $\mcT$-algebra in $\mcC$ is just a finitely product preserving functor $\mcT\to\mcC$.
This simultaneously encodes the underlying objects (by evaluating at the generic objects of $\mcT$) and the algebraic structure (given by the values of the morphisms).

\begin{remark}
  Note that there is a strong similarity to the construction of $\cond(\mcC)$ of~\ref{chap:cond_alg}.
  There,
  we use product preserving functors $\extr^{\op}\to\mcC$.
  In particular,
  a condensed set is (essentially) a product preserving functor $\extr^{\op}\to\Set$.
  Of course,
  $\extr$ is not at all generated under finite products by a finite set of objects and therefore not a Lawvere theory.
\end{remark}

\begin{proposition}
  For any category $\mcC$ and any Lawvere theory $\mcT$,
  the forgetful functor $U\colon\mcT(\mcC)\to[S,\mcC]$ is \idx{conservative},
  i.e.,
  a homomorphism of $\mcT$ algebras is an isomorphism if and only if the underlying $\mcC$ maps are isomorphisms.
\end{proposition}
\begin{proof}
  In the single-sorted case,
  this is an easy exercise.
  But then in generalises swiftly.
\end{proof}

\begin{remark}
  This is summarised in the slogan \enquote{\emph{In algebra}, any bijective homomorphism is an isomorphism}.
  This proposition is also the reason for the common atrocity in courses on linear algebra to \emph{define} the term \enquote{isomorphism} this way.
\end{remark}

\begin{proposition}
  For any Lawvere theory $\mcT$,
  the forgetful functor $\mcT(\Set)\to\Set$ has a left adjoint.
\end{proposition}
\begin{proof}
  For a proof,
  see~\cite[4.4]{Brandenburg2017}.
\end{proof}
\begin{remark}
  This also holds in much greater generality,
  in particular in any elementary topos,
  and most notably in $\cond$.
\end{remark}

\begin{remark}\label{rem:expl-law-th}
  The attentive reader may have noticed that we have not actually defined,
  for example,
  the Lawvere theory of abelian groups.
  It is not hard to do so.
  The $\hom$ set $\hom(x^{n},x)$ for the generic object $x$ is just $\Z[\{1,\ldots,n\}]$.
  Similarly for rings (with polynomial algebras), groups (free groups), modules (giving a 2-sorted algebraic theory), $R$-modules (for any ring $R$ giving a 1-sorted theory) etc.
\end{remark}

\begin{proposition}
  The forgetful functor $U\colon \mcT(\mcC)\to[S,\mcC]$ creates limits.
\end{proposition}
\begin{proof}
  Limits commute with limits and products are limits.
\end{proof}

\begin{remark}
  Of course,
  the forgetful functor usually does not create colimits.
  Consider the case of groups objects in $\Set$.
  The coproduct of two groups $G,H$ is the so-called \idx{free product} $G\ast H$.
  Its underlying set can be described by a sort of \enquote{word construction} involving the elements of $G$ and $H$
  and has nothing to do with the disjoint union of the underlying sets.
  (In fact, the coproduct of two abelian groups in $\Grp$ is the same as in $\Ab$ and so has underlying set the \emph{product} of the two underlying sets.)

  So in general,
  there is no reason to expect $\mcT(\mcC)$ to have many colimits.
  In practice, however, it does.
  This is mainly due to the fact that $\mcT(\mcC)$ is mostly of interest for exceptionally well-behaved $\mcC$ such as elementary topoi.
  In this case,
  $\mcT(\mcC)$ is bicomplete as soon as $\mcC$ (so for example, in the cases of $\Set$ and $\cond$).
  There are obviously many intermediate results of the form \enquote{If $\mcC$ has such and such (co)limits, then $\mcT(\mcC)$ has this and that (co)limit}
  but in the literature,
  these are only spelled out for monads,
  not Lawvere theories.
\end{remark}

This is on possible motivation to pass to monads but we also give a \enquote{pedagogical} introduction.
Let us mention now that there are theorems giving exact correspondences between Lawvere theories and monads.
For more on this, see~\cite{Berger_2012}.

\subsubsection{Monads}\label{subsubsec:monads}

\quot{I hail a [monad] when I see one, and I seem to see them everywhere!\\Friends have observed, however, that there are mathematical objects which are not [monads].}{Einar Hille, 1948 (referring to semigroups, freely adapted by the authors)}

Now,
we introduce a different point of view on all things algebraic.
For our purposes (and as is usually done when introducing monads) we restrict our aim to what we have begun to call \emph{single-sorted} types of algebraic objects.
Surprisingly,
despite this restriction,
the theory of monads will be of greater use than that of Lawvere theories when discussing condensed modules.

Let's stick to our toy example of abelian group objects.
For concreteness,\footnote{Pun intended.}
let us consider abelian group objects in $\Set$, i.e., abelian groups.
Given any set $M$,
we can form the set $\Z[M]$ of formal expressions involving the elements of $M$.
For example,
for $m_{1},m_{2},m_{3}\in M$,
we have $m_{1},3m_{1}-7m_{2},-m_{1}+m_{2}-m_{3}\in\Z[M]$.%
\footnote{Later, it will be very important to us how such formal expressions become equal.
  To avoid confusion:
  Yes, technically, there is already an identification happening in this step:
  We regard $m+m$ and $2m$ as the same elements of $\Z[M]$.
  However,
this a purely notational matter and has nothing to do with the actual non-trivial identifications that will come up soon.}
These formal expressions are meaningless by design.
For $M=\{\clubsuit,\spadesuit,\heartsuit,\diamondsuit\}$,
we have $\clubsuit-7\diamondsuit,-3\heartsuit\in\Z[M]$.
By giving these formal expressions meaning (for example, by saying that $\clubsuit-7\diamondsuit$ evaluates to $\heartsuit$),
we are endowing $M$ with the structure of an abelian group.
Evidently,
we need to ask for some coherence in this procedure.
For example,
the formal expression $\clubsuit\in\Z[M]$ should evaluate to $\clubsuit\in M$.
On the other hand,
for any formal expression of formal expressions,
such as $7(\clubsuit-\spadesuit)-3(2\heartsuit+\diamondsuit)\in\Z[\Z[M]]$,
the two ways of evaluating it should coincide.
I.e., it should not matter whether or not we first evaluate $\clubsuit-\spadesuit$ (giving, say, $\spadesuit$) and $2\heartsuit+\diamondsuit$ (say, to $\diamondsuit$)
and then evaluate $7\spadesuit-3\diamondsuit$ or just evaluate $7\clubsuit-7\spadesuit-6\heartsuit+3\diamondsuit$.

In analyzing (and attempting to generalise) this,
we see that we have a functor $\Z\colon\Set\to\Set$ (in general, an endofunctor $T\colon \mcC\to\mcC$)
and an algebraic structure for this functor consists of a set $M$ and a map $\Z[M]\to M$ (an object $A$ and a morphism $m\colon TA\to A$).
Now we want to require the two properties of this evaluation map that we just discussed.

The first is that $M\to\Z[M]\to M$ is the identity on $M$.
To generalize this,
we actually need more data.
Specifically,
we also take a natural transformation $\nu\colon 1_{\mcC}\to T$ as part of our monad and then require
\[
  A\stackrel{\nu_{A}}{\longrightarrow}TA\longrightarrow A
\]
to be identity.

For the other,
i.e., that the two maps $\Z[\Z[M]]\to\Z[M]\to M$ are equal,
we also subtly used structure associated to the functor $\Z[-]$ that we did not make explicit.
Namely,
there is a \enquote{formal} way of reducing formal expressions of formal expressions to (simple) formal expressions.
In other words,
there is a natural transformation $\mu\colon T\circ T\to T$ and we require the following square to commute.
\begin{center}
  \begin{tikzcd}
    TTA\ar[rr,"T(m)"]\ar[d,"\mu_A"] && TA\ar[d,"m"]\\
    TA\ar[rr,"m"] && A.
  \end{tikzcd}
\end{center}

Although we have not yet defined what a monad is exactly,
we already know that it consists of an endofunctor $T\colon \mcC\to\mcC$ and two natural transformations $\nu\colon 1_{\mcC}\to T$ and $\mu\colon TT\to T$.
This is already enough to define the notion of an algebra for a monad.

\begin{definition}\label{def:alg-over-monad}
  Given a monad $(T,\nu,\mu)$ on $\mcC$ (see below for a precise definition),
  an \textbf{algebra for the monad $(T,\nu,\mu)$}\index{algebra for a monad} or $T$-algebra
  consists of an object $A\in\mcC$ and a map $m\colon TA\to A$ such that the following diagrams commute
  \begin{center}
    \begin{tikzcd}
      TTA\ar[r,"T(m)"]\ar[d,"\mu_A"] & TA\ar[d,"m"]&& A\ar[dr,"1_A"]\ar[r,"\nu_A"] & TA\ar[d,"m"] \\
      TA\ar[r,"m"] & A &&&A.
    \end{tikzcd}
  \end{center}
  A morphism of $T$-algebras $A$, $m\colon TA\to A$, $B$, $n\colon TB\to B$ is a map $f\colon A\to B$ (in $\mcC$) such that the following diagram commutes.
  \begin{center}
    \begin{tikzcd}
      TA\ar[r,"Tf"]\ar[d,"m"] & TB \ar[d,"n"]\\
      A\ar[r,"f"] & B.
    \end{tikzcd}
  \end{center}
  We call the resulting category the \idx{Eilenberg-Moore-category} of $T$ and write it as $\EM(T)$.
  It comes equipped with a forgetful functor to $\mcC$.
\end{definition}

What we have so far pushed under the rug is that $\nu$ and $\mu$ of course need to satisfy coherences of their own.
For example, $T\stackrel{T.\nu}{\longrightarrow}TT\stackrel{\mu}{\longrightarrow}T$ should be the identity.
These can be swiftly summed up using the language of monoidal categories.

\begin{definition}
  A \idx{monad} over $\mcC$ is a monoid\footnote{in the sense of monoidal categories, see \ref{def:monoid}} in $(\End(\mcC),\circ,1_{\mcC})$, the monoidal category of endofunctors $\End(\mcC)=[\mcC,\mcC]$ of $\mcC$
  with composition as \enquote{tensor product}.
  As such,
  it consists of a functor $T\colon \mcC\to\mcC$,
  a natural transformation $\mu\colon T\circ T\to T$ called the \idx{monad multiplication}
  and a natural transformation $\eta\colon1_{\mcC}\to T$ called the \idx{monad unit}.

  We call $T$ a \idx{finitary monad}, if it commutes with filtered colimits.
\end{definition}

\begin{remark}
  For an excellent intuitive explanation for why \emph{finitary} monad is the closest to our intuition for free algebraic structures,
  see \emph{Algebra -- It's not what you think it is!}%
  \footnote{\url{youtu.be/bD8kjpynF6A}}
   from \texttt{17:54} on.
\end{remark}

\begin{proposition}\label{prop:alg-th-is-monad}
  Given an $S$-sorted Lawvere theory $\mcT$ and a category,
  there is a finitary monad $T_{\mcT}$ on $\Set^{S}$ such that $\EM(T)\simeq\mcT(\Set)$
  and the forgetful functors agree.
\end{proposition}

\begin{remark}\label{rem:ab-group-monad}
  In fact, this proposition holds in much greater generality.
  For precise statements,
  see~\cite{Berger_2012}.

  Explicitly, for the Lawvere theory $\mcT$ of abelian groups,
  $T_{\mcT}=\Z[-]$ and similarly for any ring $R$ in place of $\Z$.
\end{remark}

\begin{proposition}
  For any monad $T$ on any category $\mcC$,
  the forgetful functor $\EM(T)\to\mcC$ has a left adjoint.
\end{proposition}
\begin{proof}
  The left adjoint is given by
  \[
    C\mapsto (TC,m:=\mu_{C}\colon TTC\to TC),\qquad f\mapsto Tf.
  \]
\end{proof}

\begin{example}\label{ex:rmod_as_monad}
  One of the most central examples of monads is in the context of monoidal categories.
  Let $(\mcC,\otimes)$ be a monoidal category and $R$ a monoid in $\mcC$.
  Then,
  $R\otimes-\colon \mcC\to\mcC$ is a monad with monad multiplication and unit canonically induced by multiplication and unit of $R$.
  Algebras for this monad are precisely (left) $R$-modules.
  The reader is highly encouraged to take a moment to consider all the theorems in this section for this monad,
  in particular in the case that $R\otimes-$ is cocontinuous as this is often the case.
  (Most importantly, we will have a tensor product on condensed abelian groups that is componentwise cocontinuous.)
\end{example}

\begin{proposition}
  For any category $\mcC$ and monad $T$ on $\mcC$,
  the forgetful functor $?\colon \EM(T)\to\mcC$ is conservative.
\end{proposition}
\begin{proof}
  Clear, as $f^{-1}n=mTf^{-1}$
  so that $f^{-1}$ is again a $T$-homomorphism.
\end{proof}

\begin{proposition}\label{prop:monad-create-limits}
  For any category $\mcC$ and monad $T$ on $\mcC$,
  the forgetful functor $?\colon \EM(T)\to\mcC$ creates limits.
\end{proposition}
\begin{proof}
  This is just a much simpler version of the proof of the following proposition.
\end{proof}

\begin{proposition}\label{prop:monad-create-colimits}
  For a category $\mcC$ and a monad $T$ on it that preserves all colimits of form $\mcI$,
  the forgetful functor $?\colon\EM(T)\to\mcC$ creates colimits of form $\mcI$.
  In particular,
  for finitary monads,
  the forgetful functor preserves filtered colimits.
\end{proposition}
\begin{proof}
  Take a diagram $A\colon\mcI\to \EM(T)$ for which $?\circ A$ has a colimit in $\mcC$.
  By assumption,
  in the following diagram,
  the canonical map induced by identities is an isomorphism:
  \begin{center}
    \begin{tikzcd}
      T\varinjlim_{\mcI}?A_i \ar[rr,<->,"\sim"] && \varinjlim_{\mcI}T?A_i\\
      T?A_i \ar[u,"T(\iota_{?A_i})"] \ar[rr,"1_{T?A_i}"] && T?A_i. \ar[u,"\iota_{t?A_i}"]
    \end{tikzcd}
  \end{center}
  We may thus define a $T$-algebra structure on $\varinjlim_{\mcI}A_{i}$
  by the map $T\varinjlim_{\mcI}A_{i}\to \varinjlim_{\mcI}A_{i}$ induced by this isomorphism and the maps $TA_{i}\to A_{i}$
  (which define the $T$-algebra structure on the $A_{i}$) in the canonical way:
  \begin{center}
    \begin{tikzcd}
      && T?A_i \ar[rr]\ar[d] && ?A_i \ar[d]\\
      T\varinjlim_{\mcI} ?A_i \ar[rr,"\sim"] && \varinjlim_{\mcI}T?A_i \ar[rr] && \varinjlim_{\mcI} ?A_i.
    \end{tikzcd}
  \end{center}
  To check that this indeed defines a $T$-algebra structure,
  i.e.,
  that both ways of resolving $TT\varinjlim_{\mcI}?A_{i}$ coincide.
  By testing this on each $?A_{i}\varinjlim_{\mcI}?A_{i}$,
  this reduced the corresponding property of maps $TT?A_{i}\to ?A_{i}$.
  But of course,
  the maps $T?A_{i}\to ?A_{i}$ define $T$-algebra structures on $?A_{i}$.

  Henceforth,
  we are only left with proving that the resulting $T$-algebra indeed satisfies the universal property.
  Uniqueness is clear.
  For existence,
  take the map in $\mcC$ induced by the universal property of $\varinjlim_{\mcI}?A_{i}=?\varinjlim_{\mcI}A_{i}$.
  This map is a $T$-algebra homomorphism essentially by construction.
\end{proof}

\begin{remark}
  Evidently,
  the following stronger statement is also true:
  Any diagram $A\colon \mcI\to \EM(T)$
  for which $?A$ and $T?A$ have colimits in $\mcC$ and $T$ maps the former to the latter
  has a colimit in $\EM(T)$ with underlying object $\varinjlim_{\mcI}?A_{i}$.

  As with algebraic theories,
  we see that limits are simpler than colimits.
  However,
  in contrast to the case of algebraic theories,
  there are much more refined statements giving colimits in $\EM(T)$.
  For a discussion with pointers to the literature,
  see colimits in categories of algebras on nLab%
  \footnote{\url{https://ncatlab.org/nlab/show/colimits+in+categories+of+algebras}}.
  We give an example from~\cite{Linton_1969} next.
  For more, see~\cite{Barr1984}.
\end{remark}

\begin{theorem}[Linton]
  Let $T$ be a monad on a cocomplete category $\mcC$.
  Then if $\EM(T)$ has reflexive coequalizers (e.g., when $T$ preserves reflexive coequalizers),
  $\EM(T)$ is cocomplete.
\end{theorem}

\begin{corollary}
  For a finitary monad $T$ on a category $\mcC$,
  the left adjoint $L$ to the forgetful functor $U\colon\EM(T)\to\mcC$
  sends compact objects to compact objects.
\end{corollary}
\begin{proof}
  Let $X\in\mcC$ be a compact object and consider a filtered colimit $\varinjlim Y_{i}\in\EM(T)$.
  Since $T$ commutes with filtered colimits,
  $U$ does, too, and we get
  \begin{align*}
    \Hom_{\EM(T)}(LX,\varinjlim Y_{i})
    &=\Hom_{\mcC}(X,U(\varinjlim Y_{i}))
      =\Hom_{\mcC}(X,\varinjlim UY_{i})\\
    &=\varinjlim\Hom_{\mcC}(X,UY_{i})=\varinjlim\Hom_{\EM(T)}(LX,Y_{i})
  \end{align*}
  so that $LX$ is indeed compact.
\end{proof}

\begin{remark}
  Analogously,
  if $T$ preserves sifted colimits,
  the left adjoint sends compact projectives to compact projectives.
\end{remark}

The following proposition will be of much importance later.
\begin{proposition}\label{prop:free-r-sift}
  The functor $R[-]\colon\Set\to\Set$ preserves sifted colimits for all rings $R$.
\end{proposition}
\begin{proof}
  As $R[-]$ is finitary (\ref{rem:ab-group-monad})
  so that we only have to show that it preserves reflexive coequalizers by~\ref{lem:comm_with_sifted}.
  But in fact,
  any finitary monad over $\Set$ preserves reflexive coequalizers,
  as is easy to see.
  For a spelling out of a (weird) proof,
  see theorem~2.5 in reflexive coequalizers on nLab%
  \footnote{\url{https://ncatlab.org/nlab/show/reflexive+coequalizer}}.
\end{proof}

\clearpage{\thispagestyle{empty}\cleardoublepage}


\thispagestyle{empty}

\chapter[Condensed Algebra]{Condensed Algebra}\label{chap:cond_alg}

In this chapter, we discuss some aspects of the condensed theory mixed with algebra.
We start by defining general condensed categories, and develop special properties those categories have, depending on the starting category.

In particular, if one takes the condensation of an algebraic category, many algebraic properties remain, and lead to a \enquote{topologised} version of the algebraic theory.
This is where the condensed theory really starts to shine.
Most mathematicians (except maybe algebraic geometers) are totally used to topological versions of algebraic structures being highly pathological
and would never expect the categorical properties to be almost as good \emph{or even better} than the respective \enquote{discrete} versions.
So much so that we are certain everyone can immediately think of at least three such pathologies.
As functional analysts spending all day handling objects that carry both topological and algebraic structure,
the fact that in the condensed setup,
those two finally come together in a harmonious way\footnote{Tempting one to call this real harmonic analysis.}
is thrilling (especially in combination with a nice notion of completeness given by analytic rings).

We will in particular focus on condensed $R$-modules (first for classical, then for condensed rings $R$) and condensed abelian groups.
After having seen some general properties of these categories, 
we will discuss completeness, in form of analytic rings.

Using these results, in the last section we will start developing some ergodic theory.

\section{Condensed categories}\label{sec:cond_cats}
In this section, we will discuss general condensed categories and how the process of \enquote{condensing} a category preserves many important characteristics.
\begin{definition}[Condensed 1-categories]\label{def:cond_C}\uses{def:cond}\chapfour
    Let $\mcC$ be any category.
We define the \idx{condensation of $\mcC$}, $\cond(\mcC)$ as the category of contravariant functors $X\colon\extr\to \mcC$ such that
\begin{itemize}
    \item $X$ maps finite coproducts in $\extr$ to finite products in $\mcC$
    \item $X$ is the pointwise left Kan extension of its restriction onto some $\extr_\kappa$.
\end{itemize}
    Furthermore, for any cardinal $\kappa$ define $\cond_\kappa(\mcC)$ as the category of contravariant functors $\extr_\kappa\to \mcC$ sending finite coproducts to products.
\end{definition}
However, if $\mcC$ does not admit filtered colimits, then not necessarily every Kan extension is pointwise,
and if filtered colimits do not commute with finite products, then the Kan extension does not need to preserve the property of being product preserving.

\begin{lemma} If $\mcC$ admits filtered colimits, then every Kan extension of any functor $T\colon \extr_\kappa^{\op}\to\mcC$ along inclusions $\extr_\kappa^{\op}\to \extr_{\kappa'}^{\op}$ or the inclusion $\extr_\kappa^{\op}\to \extr^{\op}$ exists and is pointwise, and extends the functor $T$.
\end{lemma}
\begin{proof}
    Let $T\colon \extr_\kappa^{\op}\to \mcC$ and consider the embedding $\tau_\kappa^{\kappa'}\colon \extr_\kappa^{\op}\to \extr_{\kappa'}^{\op}$.
    Then the Kan limit formula reads as
    \[(\Lan_{\tau_\kappa^{\kappa'}}T)(S)=\varinjlim_{S\to K\in \extr_\kappa}T(K).\]
    This colimit is over a small $\mathrm{cof}(\kappa)$-filtered diagram (see \ref{lem:extr-lambda-filtered}), and hence exists.
    Thus the pointwise left Kan extension exists (\ref{thm:computation_of_kan_extensions}).
    As the embedding $\tau_\kappa^{\kappa'}$ is fully faithful,
    the restriction of the Kan extension gives back the original functor (\ref{lem:Kan_extension_extends}).

    The result for $\extr_\kappa\to \extr$ follows by passage to the union, see \ref{lem:kan-extension-fully-faithful}.
\end{proof}
\begin{lemma}\label{lem:kan_on_extr_sh}
    If filtered colimits commute with finite products in $\mcC$, then for any finite product preserving functor $T\colon \extr_\kappa^{\op}\to\mcC$
with pointwise left Kan extension along $\extr_\kappa^{\op}\to \extr_{\kappa'}^{\op}$ or $\extr_\kappa^{\op}\to \extr^{\op}$,
    this Kan extension still yields a finite product preserving functor.
\end{lemma}

\begin{proof}
Note that $\kappa$ is an infinite cardinal (in particular a limit ordinal), so it has infinite cofinality.
Hence the diagram of the left Kan extension along the inclusion
$\tau\colon\extr_\kappa^{\op}\to\extr^{\op}_{\kappa'}$ is always filtered by lemma~\ref{lem:extr-lambda-filtered} (note that here we are in the opposite category),
and thereby colimits over this diagram commute with finite products in $\mcC$.
Fix $A,B\in \extr_{\kappa'}^{\op}$.
We want to show
\[\Lan_{\tau}T(A\sqcup B)=\Lan_{\tau} T(A)\times \Lan_{\tau} T(B).\]
Let us first look at the diagram on the left hand side,
\[\Lan_{\tau} T(A\sqcup B)=\varinjlim_{A\sqcup B\to S} T(S) \]
Here the diagram, denoted as in $\extr_\kappa$ (so now it is cofiltered),
is given as follows, where $S$ is taken over all $S\in \extr_{\kappa}$.
\begin{center}\begin{tikzcd}
	& {A\sqcup B} \\
	S && {S'} \\
	{T(S)_{h,g}} && {T(S')_{h',g'}}
	\arrow["{{(h,g)}}"', from=1-2, to=2-1]
	\arrow["{{(h',g')}}", from=1-2, to=2-3]
	\arrow["f"', from=2-1, to=2-3]
	\arrow["{{T(f)}}", from=3-3, to=3-1]
\end{tikzcd}.\end{center}

    We want to identify the diagram on the right hand side with a initial subdiagram of this.
    For this we use that filtered colimits commute with finite products and that $T$ commutes with finite products.
    Now the right hand side of the equations is given by
\[
    (\Lan_{\tau} T (A))\times (\Lan_{\tau} T(B))
    =\varinjlim_{A\to S}T(S)\times \varinjlim_{B\to S'} T(S')
\]
Now we can pull out the colimit one for one to obtain
\[
    \varinjlim_{(A\to S, B\to S')}T(S\sqcup S'),
\]
where the diagram is over the product form of the index categories $A\to S$ and $B\to S'$,
i.e., its objects consist of tuples $(g\colon A\to S, h\colon B\to S')$
and the morphisms are given by $(f\colon S\to S_0,f'\colon S'\to S_0')$ such that the diagrams
\begin{center}\begin{tikzcd}
	& A &&& B \\
	S && {S_0} & {S'} && {S'_0} \\
	& {T(S\sqcup S')} &&& {T(S_0\sqcup S_0')}
	\arrow["g"', from=1-2, to=2-1]
	\arrow["{{g_0}}", from=1-2, to=2-3]
	\arrow["h"', from=1-5, to=2-4]
	\arrow["{{h_0}}", from=1-5, to=2-6]
	\arrow["f"', from=2-1, to=2-3]
	\arrow["{{f'}}"', from=2-4, to=2-6]
	\arrow["{{T(f\sqcup f')}}"', from=3-2, to=3-5]
\end{tikzcd}\end{center}
are commutative.
Gluing such a tuple $(g\colon A\to S, h\colon B\to S')$ we assign to it the morphism
$g\sqcup h\colon A\sqcup B\to S\sqcup S'$.

This assignment yields a functor between the index categories from
the colimit on the right to the colimit on the left.
So in view of lemma~\ref{lem:limits_along_final_subdiagrams} it remains to show
that this functor is initial (we are working in $\extr_\kappa$ and $T$ reverses the
direction of the morphisms).

Note that both index categories are cofiltered.
Hence the zig-zag condition (ii) in the definition~\ref{def:final} of initial subcategory is fulfilled.

Thus take any morphism $f\colon A\sqcup B\to S$.
Then two copies of $f$ yield a morphism $f\sqcup f\colon A\sqcup B\to S\sqcup S$
that glues to $f$ under the surjection $S\sqcup S\twoheadrightarrow S$,
as depicted in the following diagram
\begin{center}
\begin{tikzcd}
	& {A\sqcup B} \\
	S && {S\sqcup S}
	\arrow["f"', from=1-2, to=2-1]
	\arrow["{f\sqcup f}", from=1-2, to=2-3]
	\arrow[two heads, from=2-3, to=2-1]
\end{tikzcd}.
\end{center}
This is precisely condition (i) in the definition of initial.

Lastly, we need to check that $\Lan_{\tau} T(\emptyset)=\ast$.
This follows elementary by noting that the identity on $\emptyset$ is initial in
the diagram of all morphism $\emptyset\to S$,
and hence the colimit is computed by evaluating $T$ there, so
\[
    \Lan_{\tau }T (\emptyset)=\varinjlim_{\emptyset \to S} T(S)=T(\emptyset)=\ast.
\]
\end{proof}

The property that filtered colimits commute with finite products is not fulfilled in every category, as the example~\ref{ex:filtered-colimit-limit-not-commute} shows.

\begin{lemma}\label{lem:cond_as_colimit}
Let $\mcC$ have filtered colimits that commute with finite products.
Left Kan extension $T\in\cond_\kappa(\mcC)$ along $\tau\colon\extr_\kappa\hookrightarrow\extr_{\kappa'}$ and $\tau\colon\extr_\kappa\to\extr$
yields an object of $\cond_{\kappa'}(\mcC)$ or of $\cond(\mcC)$, respectively.
The induced functor $\Lan_{\tau}(-)\colon\cond_\kappa(\mcC)\to\cond_{(\kappa')}(\mcC)$
is the fully faithful left adjoint to the restriction $T\mapsto T\circ \tau$.
The unit of this adjunction is given by the identity morphism
$T\to \res\Lan_\tau T$.

In particular, 
\[
    \cond(\mcC)=\varinjlim_{\kappa}\cond_\kappa(\mcC).
\]
\end{lemma}
\begin{proof}
	This is basically just the statement that if all Kan extensions along some $\tau$ exist, they form a left adjoint to the pullback functor $\tau^*$, see \ref{rem:Kan-adjoint}.

    Consider the identity $T\to (\Lan_\tau T)\tau$.
    For any $F\in \cond_{\kappa'}$ and $T\to F\tau$ we need to show the existence of a unique $g\colon \Lan_\tau T\to F$ with

	\begin{center}\begin{tikzcd}
	T & {(\Lan_\tau T)\tau} \\
	& {F\tau }
	\arrow["\id", from=1-1, to=1-2]
	\arrow["f"', from=1-1, to=2-2]
	\arrow["{g.\tau}", dashed, from=1-2, to=2-2]
\end{tikzcd}\end{center}

    But this is precisely the universal property of the left Kan extension:
    $F$ is a functor from $\extr_{\kappa'}$ to $\mcC$, $f$ a natural transformation $T\to F\tau$, and by the universal property of $\Lan$, there exists a unique $g\colon \Lan_\tau T\to F$,
    with $g.\tau$ making

\begin{center}\begin{tikzcd}
	{\extr_\kappa} && \mcC \\
	\\
	{\extr_{\kappa'}}
	\arrow[""{name=0, anchor=center, inner sep=0}, "T", from=1-1, to=1-3]
	\arrow["\tau"', hook, from=1-1, to=3-1]
	\arrow[""{name=1, anchor=center, inner sep=0}, "{\Lan_\tau T}"{description, pos=0.3}, from=3-1, to=1-3]
	\arrow[""{name=2, anchor=center, inner sep=0}, "F"{description, pos=0.3}, shift right=5, from=3-1, to=1-3]
	\arrow["f"{description, pos=0.4}, shorten <=5pt, shorten >=5pt, Rightarrow, from=0, to=2]
	\arrow["1"{description, pos=0.4}, shorten >=9pt, Rightarrow, from=0, to=3-1]
	\arrow[shorten <=2pt, shorten >=2pt, Rightarrow, from=1, to=2]
\end{tikzcd}\end{center}
commute.
To see $\cond(\mcC)=\varinjlim_{\kappa}\cond_\kappa(\mcC)$,
we use the obvious identification of objects:
Every $\kappa$-condensed object can be Kan-extended to a condensed object,
and every condensed object can canonically be restricted to being $\kappa$-condensed.
The $\hom$ isomorphism follows from
\[
    \hom(\Lan T,\Lan S)=\hom(T, (\Lan S)_{|}=\hom(T, S).
\]
\end{proof}

\begin{lemma}[$\condk(\mcC)$ as sheaves]\label{lem:cond_C_as_sheaves}\uses{def:cond_C}
Let $\kappa$ be a strong limit cardinal.
We equip the categories $\extr_\kappa$, $\prolak$ and $\CHaus_\kappa$ with the
Grothendieck topology from definition~\ref{def:finitary_grothendieck_topology}.
Then the following categories are equivalent:
\begin{enumerate}[(a)]
    \item the $\kappa$-condensation of $\mcC$, $\cond_\kappa(\mcC)$,
    \item the category $\Sh(\extr_\kappa,\mcC)$ of $\mcC$-valued sheaves on $\extr_\kappa$
    \item the full subcategory of functors $T\colon\prolak^\mathrm{op}\to\mcC$
    mapping finite coproducts and coequalizers of kernel pairs in $\prolak$
    to finite products and to equalizers of kernel pairs, resp., in $\mcC$,
    \item the category $\Sh(\prolak)$ of $\mcC$-valued sheaves on $\prolak$,
    \item the full subcategory of functors $T\colon\CHaus_\kappa^\mathrm{op}\to\mcC$
    mapping finite coproducts and coequalizers of kernel pairs  in $\CHaus_\kappa$
    to finite products and to equalizers of kernel pairs, resp., in $\mcC$,
    \item the category $\Sh(\CHaus_\kappa)$ of $\mcC$-valued sheaves on $\CHaus_\kappa$.
\end{enumerate}
\end{lemma}

\begin{proof}
The equivalences of (b), (d) and (f) follow from the enriched comparison lemma~\ref{thm:comparison_lemma}.
The corresponding equivalences (a)$\iff$(b), (c)$\iff$(d) and (e)$\iff$(f)
are all proven similarly, so we restrict the proof to the equivalence of (a) and (b).

Let $T\colon\extr_\kappa^\mathrm{op}\to\mcC$ be a sheaf.
Then the empty sieve is a cover of the empty set.
Therefore, $T(\emptyset) = \varprojlim_{S\in\emptyset} T(S) = \ast$
where $\ast$ is the terminal object of $\mcC$.

Now let $S,S'\in\extr_\ast$.
We will show that the induced map
\begin{center}
\begin{tikzcd}
	& {T(S)\times T(S')} \\
	{T(S)} && {T(S')} \\
	& {T(S\sqcup S')}
	\arrow[from=1-2, to=2-1]
	\arrow[from=1-2, to=2-3]
	\arrow[dashed, from=3-2, to=1-2]
	\arrow[from=3-2, to=2-1]
	\arrow[from=3-2, to=2-3]
\end{tikzcd}
\end{center}
is an isomorphism.
For if $C\in\mcC$ has morphisms $C\to T(S)$ and $C\to T(S')$,
then for $S_0$ such that
\begin{center}
\begin{tikzcd}
	{S_0} & S \\
	{S'} & {S\sqcup S'}
	\arrow[from=1-1, to=1-2]
	\arrow[from=1-1, to=2-1]
	\arrow[from=1-2, to=2-2]
	\arrow[from=2-1, to=2-2]
\end{tikzcd}
\end{center}
is commutative,
then $S_0 = \emptyset$, $T(S_0) = \ast$ and hence the diagram
\begin{center}
\begin{tikzcd}
	\ast & {T(S)} \\
	{T(S')} & {T(S\sqcup S')}
	\arrow[from=1-2, to=1-1]
	\arrow[from=2-1, to=1-1]
	\arrow[from=2-2, to=1-2]
	\arrow[from=2-2, to=2-1]
\end{tikzcd}
\end{center}
is commutative.
Similarly, if $f,g\colon S_0\to S$ are parallel morphisms
with $e_S\circ f = e_S\circ g$ then $f = g$
since the embedding $e_S\colon S\to S\sqcup S'$ is monic.
Therefore $f^\ast = g^\ast$ .
Since $T(S\sqcup S') = \varprojlim_{E\in J} T(E)$ where $J$ is the cover of $S\sqcup S'$
generated by the embeddings $S\to S\sqcup S'$, $S'\to S\sqcup S'$,
there exists a unique morphism $C\to T(S\sqcup S')$ making the diagram
\begin{center}
\begin{tikzcd}
	& C \\
	{T(S)} && {T(S')} \\
	& {T(S\sqcup S')}
	\arrow[from=1-2, to=2-1]
	\arrow[from=1-2, to=2-3]
	\arrow[dashed, from=1-2, to=3-2]
	\arrow[from=3-2, to=2-1]
	\arrow[from=3-2, to=2-3]
\end{tikzcd}
\end{center}
commutative.

For the other implication, let $T\colon\extr_\kappa^\mathrm{op}\to\mcC$
be a functor preserving finite products.
Let $(f_i\colon S_i\to S)_{i=1}^n$ be a cover of $S$,
i.e. the induced morphism $f\colon\bigsqcup_{i=1}^n S_i \to S$ is epic.
Since $S$ is extremally disconnected,
there exists a retraction $s$ with $f\circ s = 1_S$.
We need to show that $T(S)$ is the limit of the diagram generated by the given cover of $S$.
Let $C\in\mcC$ with projections $p_E$ for every morphism $E\to S$
that factors over two (not necessarily distinct) $f_i$'s.
Since $T(\bigsqcup_{i=1}^n S_i)$ is the product of the $T(S_i)$,
the projections $p_{S_i}$ glue to a morphism
\[
    \varphi\colon C\to T\left(\bigsqcup_{i=1}^n S_i\right)
\]
such that $e_i^\ast\circ\varphi = p_{S_i}$
where $e_i\colon S_i\to\bigsqcup_{i=1}^n S_i$ is the embedding.
We will show that $s^\ast\circ\varphi\colon X\to T(S)$ is the desired morphism.

Let $A_i = s^{-1}(S_i)$ and denote by $s_i\colon A_i\to S_i$ the restriction of $s$ to $A_i$.
Then we have $f_is_if_i = f_i$ since $f_is_i = 1_{A_i}$.
In particular, $f_i^\ast = (s_i\circ f_i)^\ast\circ f_i^\ast$,
which implies $p_{S_i} = (s_i\circ f_i)^\ast\circ p_{S_i}$
by definition of the diagram.
Therefore,
\[
    f_i^\ast\circ s^\ast\circ\varphi
    = (s\circ f_i)^\ast\circ\varphi
    = (e_i\circ s_i\circ f_i)^\ast\circ\varphi
    = (s_i\circ f_i)^\ast\circ e_i^\ast\circ\varphi
    = (s_i\circ f_i)^\ast\circ p_{S_i}
    = p_{S_i}.
\]
Now denote by $j_i\colon A_i\to S$ the inclusions.
For uniqueness of $s^\ast\circ\varphi$ let $\psi\colon C\to T(S)$ any other morphism
making the diagram commutative
We have that
$f_i^\ast\circ s^\ast\circ\varphi = f_i^\ast\circ\psi$
and by composing with $s_i^\ast$ we obtain that
$j_i^\ast\circ s^\ast\circ\varphi = s_i^\ast\circ f_i^\ast\circ s^\ast\circ\varphi
= s_i^\ast\circ f_i^\ast\circ \psi = j_i^\ast\circ\psi$
Since $A = \bigsqcup_{i=1}^n A_i$, we obtain that $T(S)$ is the product
$\prod_{i=1}^n T(A_i)$.
This implies the desired equality $s^\ast\circ\varphi = \psi$.
\end{proof}

Passing to the colimit along all $\kappa$ of the categories from the previous lemma,
we obtain the following characterization of $\cond(\mcC)$.

\begin{corollary}
The following categories are equivalent.
\begin{enumerate}[(a)]
    \item The condensation $\cond(\mcC)$ of $\mcC$,
    \item the colimit $\varinjlim_\kappa\condk(\mcC)$,
    \item the full subcategory of sheaves on $\extr$ with values in $\mcC$
    that are left Kan extension of its restriction to some $\extr_\kappa$,
\end{enumerate}
\end{corollary}

\begin{proof}
The equivalence of (a) and (b) is lemma~\ref{lem:cond_as_colimit}.
The equivalence of (a) and (c) follows from the equivalence of (a) and (b) in lemma~\ref{lem:cond_C_as_sheaves}.
\end{proof}

\begin{question}
  One might hope for the result that the following are also equivalent characterisations:
  \begin{itemize}
    \item the full subcategory of functors $T\colon\prof^\mathrm{op}\to\mcC$
    mapping finite coproducts and coequalizers of kernel pairs in $\prof$
    to finite products and to equalizers of kernel pairs, resp., in $\mcC$,
    that are left Kan extension of their restriction to some $\prolak$,
    \item the full subcategory of sheaves on $\prof$ with values in $\mcC$
    that are left Kan extension of their restriction to some $\prolak$,
    \item the full subcategory of functors $T\colon\CHaus^\mathrm{op}\to\mcC$
    mapping finite coproducts and coequalizers of kernel pairs in $\CHaus$
    to finite products and to equalizers of kernel pairs, resp., in $\mcC$,
    that are left Kan extension of their restriction to some $\CHaus_\kappa$,
    \item the full subcategory of sheaves on $\CHaus$ with values in $\mcC$
    that are left Kan extension of their restriction to some $\CHaus_\kappa$.
  \end{itemize}
  We do not know whether the Kan extension along $\prof_\kappa\to \prof$ (or even $\prof_{\kappa}\to\prof_{\kappa'}$)
  takes sheaves to sheaves.
  Whatever the case may be,
  we believe the characterisations can be made, sheafifying (see below) if necessary.
\end{question}

\begin{lemma}[Sheafification]\label{lem:sheafification_cond_C}\uses{def:cond_C}
    Assume that $\mcC$ admits filtered colimits, finite products and that filtered colimits commute with finite products.

Then sheafification of presheaves $F\in [\extr^{\op}, \mcC]$ exists and can be computed by the usual $F^+$ formula, which in our case reduces to
\[\Sh(F)(S)=\varinjlim_{S=\coprod_{i\in I\,\mathrm{fin.}}S_i}\prod F(S_i).\]
with $S_i\in \extr$.
Analogous statements hold for $\cond_\kappa$.

Furthermore, sheafification commutes with Kan extensions, i.e., for any $F\in [\extr_\kappa,\mcC]$, the following diagram commutes

\begin{center}\begin{tikzcd}
	{[\extr_\kappa^{\op},\mcC]} & {\cond_\kappa} \\
	{[\extr_{\kappa'}^{\op},\mcC]} & {\cond_{\kappa'}} \\
	\vdots & \vdots \\
	{\varinjlim[\extr^{\op}_\kappa,\mcC]} & \cond
	\arrow["\Sh", from=1-1, to=1-2]
	\arrow["\Lan"', from=1-1, to=2-1]
	\arrow["\Lan", from=1-2, to=2-2]
	\arrow["\Sh"', from=2-1, to=2-2]
	\arrow["\Lan"', from=2-1, to=3-1]
	\arrow["\Lan", from=2-2, to=3-2]
	\arrow["\Sh"', from=4-1, to=4-2]
\end{tikzcd}\end{center}

\end{lemma}
\begin{proof}

  We show that $\Sh(F)$ as defined in the lemma fulfills the sheaf conditions.
	First we explain how to sheafify a morphism.

	 For any $\phi\colon S\to S'$ the map $\Sh(F)(\phi)\colon \varinjlim_{S'=\coprod S_i'} \prod F(S_i')\to \varinjlim_{S=\coprod S_i}\prod F(S_i)$
    is induced by taking for any $S'=\coprod S_i'$ the pullback decomposition $S_i\coloneqq \phi^{-1}(S_i')$, and
    restricting $\phi_i\colon S_i\to S_i'$, applying $F$ to obtain morphisms $F(\phi_i)\colon F(S_i')\to F(S_i)$, gluing them to a morphism
    $\prod F(S_i')\to \prod F(S_i)$ and afterwards embedding into $\varinjlim_{S=\coprod S_i} F(S_i)$.

    For this, first note that the empty cover of the empty set is initial in the diagram of all covers of $\emptyset$.
    Hence
    \[\Sh(F)(\emptyset)=\varinjlim_{\emptyset=\sqcup_{i\in I} S_i}\prod_{I}F(S_i)=\prod_{\emptyset} F(\emptyset)=\ast.\]
    Now, take $X\sqcup Y$ in $\extr_\kappa$.
    We can refine any decomposition $X\sqcup Y=\coprod_{i=1}^n A_i$ by defining $B_i=A_i\cap X$ and $C_i=A_i\cap Y$ to obtain a finer covering
    \[X\sqcup Y=\coprod_{i=1}^n B_i\coprod_{j=1}^n C_j.\]
    This forms an initial subdiagram, and hence
    \begin{align*}
        \Sh(F)(X\sqcup Y)&=\varinjlim_{X\sqcup Y=\coprod_{i=1}^n A_i}\prod F(A_i)\\
        &=\varinjlim_{X=\coprod B_i,\, Y=\coprod C_j} \prod F(B_i)\times \prod F(C_j)\\
        &=\varinjlim_{X=\coprod B_i}\varinjlim_{Y=\coprod C_j} \prod F(B_i)\times \prod F(C_j)\\
        &=\left(\varinjlim_{X=\coprod B_i} \prod F(B_i)\right)\times \left(\varinjlim_{Y=\coprod C_j}\prod F(C_j)\right)\\
        &=\Sh(F)(X)\times \Sh(F)(Y).
        \end{align*}

    For any $F$, define the unit $\eta_F\colon F\to \Sh(F)$ as $\eta_F(S)$ to be the injection of the coordinate $S=\coprod_{i=1}^1 S$ into the colimit,
    \[\eta_F(S)\colon F(S)\to \varinjlim_{S=\coprod S_i} \prod F(S_i).\]

   This is natural, as for $\phi\colon S\to S'$

\begin{center}\begin{tikzcd}
	S & {F(S)} & {\Sh(F)(S)} \\
	{S'} & {F(S')} & {\Sh(F)(S')}
	\arrow["\phi"', from=1-1, to=2-1]
	\arrow["{\eta_S}", from=1-2, to=1-3]
	\arrow["{F(\phi)}", from=2-2, to=1-2]
	\arrow[from=2-2, to=2-3]
	\arrow["{\Sh(F)(\phi)}"', from=2-3, to=1-3]
\end{tikzcd}\end{center}
commutes.

    For the universal property, consider any sheaf $G$ and natural transformation $\eps\colon F\to G$, we want a unique natural transformation $\phi\colon \Sh(F)\to G$ with

\begin{center}
\begin{tikzcd}
	& G \\
	F & {\Sh(F)}
	\arrow["\eps", from=2-1, to=1-2]
	\arrow["\eta"', from=2-1, to=2-2]
	\arrow["\phi"', dashed, from=2-2, to=1-2]
\end{tikzcd}
\end{center}
Consider any $S\in \extr$.
    Then this induces an arrow $\phi_S$ via

\begin{center}\begin{tikzcd}
	&& {G(S)} & {\prod G(S_i)} \\
	\\
	{F(S)} && {\varinjlim\prod F(S_i)} & {\prod F(S_i)}
	\arrow["\simeq"', from=1-4, to=1-3]
	\arrow["{\eps_S}", from=3-1, to=1-3]
	\arrow["{\eta_S}"', from=3-1, to=3-3]
	\arrow["{\phi_S}"', dashed, from=3-3, to=1-3]
	\arrow["{\prod\eps_{S_i}}"', from=3-4, to=1-4]
	\arrow[from=3-4, to=3-3]
\end{tikzcd}\end{center}
which clearly commutes.
    Again, the naturality can be seen easily.
    For the uniqueness take any $\phi'_S$ and note that it suffices to show commutativity of

\begin{center}\begin{tikzcd}
	{G(S)} & {\prod G(S_i)} \\
	\\
	{\varinjlim\prod F(S_i)} & {\prod F(S_i)}
	\arrow["\simeq"', from=1-2, to=1-1]
	\arrow["{\phi_S'}"', dashed, from=3-1, to=1-1]
	\arrow["{\prod\eps_{S_i}}"', from=3-2, to=1-2]
	\arrow[from=3-2, to=3-1]
\end{tikzcd}\end{center}

    which follows by reversing the isomorphism $G(S)=\prod G(S_i)$ and noting that the outer morphisms do not depend on $\phi'$ and commute, as they commute for $\phi$.

\begin{center}\begin{tikzcd}
	&& {G(S)} & {\prod G(S_i)} \\
	\\
	{F(S)} && {\varinjlim\prod F(S_i)} & {\prod F(S_i)}
	\arrow["\simeq", from=1-3, to=1-4]
	\arrow["{\eps_S}"{description}, from=3-1, to=1-3]
	\arrow["{\eta_S}"{description}, from=3-1, to=3-3]
	\arrow["{\phi_S'}"', dashed, from=3-3, to=1-3]
	\arrow["{\prod\eps_{S_i}=\prod \phi_{S_i}'\eta_{S_i}}"', from=3-4, to=1-4]
	\arrow[from=3-4, to=3-3]
\end{tikzcd}.\end{center}

    To see commutativity with left Kan extensions, note both compositions of $\Lan$ and $\Sh$ form left adjoints to $\cond_{\kappa'}\to [\extr_\kappa^{\op},\mcC]$.

\begin{center}\begin{tikzcd}
	{[\extr_\kappa^{\op},\mcC]} & {\cond_\kappa} \\
	{[\extr_{\kappa'}^{\op},\mcC]} & {\cond_{\kappa'}} \\
	\vdots & \vdots \\
	{\varinjlim[\extr^{\op}_\kappa,\mcC]} & \cond
	\arrow[""{name=0, anchor=center, inner sep=0}, "\Sh", curve={height=-12pt}, from=1-1, to=1-2]
	\arrow[""{name=1, anchor=center, inner sep=0}, "\Lan"', curve={height=12pt}, from=1-1, to=2-1]
	\arrow[""{name=2, anchor=center, inner sep=0}, "{?}"{description}, curve={height=-12pt}, from=1-2, to=1-1]
	\arrow[""{name=3, anchor=center, inner sep=0}, "\Lan", curve={height=-12pt}, from=1-2, to=2-2]
	\arrow[""{name=4, anchor=center, inner sep=0}, "{?}"', curve={height=12pt}, from=2-1, to=1-1]
	\arrow[""{name=5, anchor=center, inner sep=0}, "\Sh"', curve={height=12pt}, from=2-1, to=2-2]
	\arrow["\Lan"', curve={height=12pt}, from=2-1, to=3-1]
	\arrow[""{name=6, anchor=center, inner sep=0}, "{?}", curve={height=-12pt}, from=2-2, to=1-2]
	\arrow[""{name=7, anchor=center, inner sep=0}, "{?}"{description}, curve={height=12pt}, from=2-2, to=2-1]
	\arrow["\Lan", curve={height=-12pt}, from=2-2, to=3-2]
	\arrow[curve={height=12pt}, from=3-1, to=2-1]
	\arrow[curve={height=-12pt}, from=3-2, to=2-2]
	\arrow["\Sh"', curve={height=12pt}, from=4-1, to=4-2]
	\arrow["{?}"', curve={height=12pt}, from=4-2, to=4-1]
	\arrow["\dashv"{anchor=center}, draw=none, from=1, to=4]
	\arrow["\dashv"{anchor=center, rotate=-90}, draw=none, from=0, to=2]
	\arrow["\dashv"{anchor=center, rotate=-180}, draw=none, from=3, to=6]
	\arrow["\dashv"{anchor=center, rotate=90}, draw=none, from=5, to=7]
\end{tikzcd}\end{center}
commutativity follows by uniqueness of adjoints.
    \end{proof}

\begin{definition}
	Having a sheafification, we can find a canonical functor $\mcC\to \cond(\mcC)$,
	mapping a $C\in\mcC$ to the \idx{discrete condensed object} $C_d$ given by sheafification of the constant (this clearly fulfills the Kan condition) presheaf.
\end{definition}
\begin{question}
	We have a canonical candidate for a left adjoint to the discrete objects, given by evaluation at a point, $T\mapsto T(\ast)$.
	However, we are unsure whether this in general yields a left adjoint.
\end{question}
\begin{lemma}[Exactness of sheafification]\label{lem:sh_exact}
	If filtered colimits in $\mcC$ exist and are exact, then filtered colimits may be computed pointwise on $\extr$.
	Furthermore, sheafification is an exact functor.
\end{lemma}
\begin{proof}
	Consider any filtered diagram $(X_i)$ in $\cond(\mcC)$.
	The colimit is computed by sheafification of the pointwise colimit, thus it suffices to show that the pointwise filtered colimit again yields a sheaf.
	But this follows with
	\[\varinjlim X_i(S\sqcup T)=\varinjlim X_i(S)\times X_i(T)=(\varinjlim X_i(S))\times (\varinjlim X_i(T)).\]
	Sheafification clearly is right exact as a left adjoint an thus even cocontinuous.
	Hence, take any finite limit of some $X_i$.
	Inspecting the formula for sheafification yields
	\[\Sh(\varprojlim X_i)(S)=\varinjlim_{S=\sqcup T_j} \prod_j \varprojlim_i X_i(T_j)=\varprojlim_{i}\varinjlim_{S=\sqcup T_j}\prod_j X_i(T_j)=\varprojlim \Sh(X_i)(S).\]
	\end{proof}

\begin{definition}[Condensing as 2-functor]\label{def:cond_C_functor}\uses{def:cond_C}
       Clearly, any functor $F\colon\mcC\to \mcD$ induces a functor $\cond_\kappa(F)\colon \cond_\kappa(\mcC)\to\cond_\kappa(\mcD)$ by forgetting, pushforward, followed by sheafification;
\[\cond_\kappa(F)(T)=\Sh(F\circ T).\]
To any natural transformation $\eta\colon F\to G$ between functors $F,G\colon \mcC\to\mcD$, we can associate a natural transformation $\cond_\kappa(\eta)\colon \cond_\kappa(F)\to\cond_\kappa(G)$ via
\[\cond_\kappa(\eta)_{T}=\Sh(\eta.T)\colon \Sh(FT)\to \Sh(GT).\]

By passing to the colimit we obtain a 2-functor $\cond\colon \Cat\to\Cat$.
This yields a 2-adjunction

\begin{center}\begin{tikzcd}
	\Cat && \Cat
	\arrow[""{name=0, anchor=center, inner sep=0}, "\cond"', curve={height=12pt}, from=1-1, to=1-3]
	\arrow[""{name=1, anchor=center, inner sep=0}, "{\varinjlim[\extr_\kappa^{\op}, -]}", curve={height=-12pt}, from=1-1, to=1-3]
	\arrow[""{name=2, anchor=center, inner sep=0}, "{?}", shift left=4, shorten <=3pt, shorten >=3pt, Rightarrow, from=0, to=1]
	\arrow[""{name=3, anchor=center, inner sep=0}, "\Sh", shift left=4, shorten <=3pt, shorten >=3pt, Rightarrow, from=1, to=0]
	\arrow["\dashv"{anchor=center, rotate=180}, draw=none, from=3, to=2]
\end{tikzcd}\end{center}
\end{definition}
\begin{remark}\label{rem:cond-monad}
  Just like any other 2-functor $\Cat\to\Cat$,
  $\cond$ thus preserves monads.
  The 2-functoriality might also directly imply the following statement,
  but we give an elementary proof.
\end{remark}
\begin{lemma}[Condensing preserves adjunctions]\label{lem:cond_pres_adj}\uses{def:cond_C, lem:sheafification_cond_C}
Again, let $\mcC$ and $\mcD$ be categories with filtered colimits that commute with finite products.

For any adjunction $L\dashv R$, there is an induced adjunction $\cond_\kappa(L)\dashv\cond_\kappa(R)$.
$\cond(R)$ is given by forgetting and pullback, no sheafification is necessary.
Taking colimits induces an adjunction $\cond(L)\dashv \cond(R)$.
\end{lemma}
 \begin{proof}
   Consider any adjunction $L\dashv R\colon \mcC\to\mcD$, and denote forgetful functors by $?$.
   Then this induces by pushforward an adjunction between $L_*\colon \mcC^{\extr_\kappa^\op}\to\mcD^{\extr_\kappa^{\op}}$ and $R_*$ (\ref{lem:adj_in_functor_cats}).
   By composing with the sheafification adjunction,
   we obtain an adjunction $\Sh\circ L_*\colon \mcC^{\extr_\kappa^{\op}}\to \Sh(\extr)$ to $R_* \circ ?$ (precomposed with the forgetful functor).

    As the right adjoint commutes with limits, it sends sheaves on $\extr$ to sheaves, and hence $R_* T\in \cond_\kappa(\mcD)$ for all $T\in \cond_\kappa(\mcC)$.
     Thus, $R_*?=?\Sh\circ R_*?\colon \cond_\kappa(\mcD)\to[\extr_\kappa^{\op},\mcC]$ and hence the adjunction restricts to an adjunction
    \[\cond_\kappa(L)=\Sh\circ L_* ?\vdash \cond_\kappa(R)= \Sh R_* ?.\]
Explicitly, for $F\in \cond_\kappa(\mcC)$ and $G\in \cond_\kappa(\mcD)$ one obtains
	 \begin{align*}
		 \hom(\Sh L_* ? F, G)\simeq \hom(?F, ?\Sh R_* ?G)\simeq \hom(F,\Sh R_*?G).
		 \end{align*}

	 Hence we obtain a tower of adjunctions

\begin{center}\begin{tikzcd}
	{\cond_\kappa(\mcC)} && {\cond_\kappa(\mcD)} \\
	{\cond_{\kappa'}(\mcC)} && {\cond_{\kappa'}(\mcD)} \\
	\vdots && \vdots \\
	{\cond(\mcC)} && {\cond(\mcD).}
	\arrow[""{name=0, anchor=center, inner sep=0}, "{\Sh \,L_*}"{description}, curve={height=-12pt}, from=1-1, to=1-3]
	\arrow["\Lan"', curve={height=6pt}, from=1-1, to=2-1]
	\arrow[""{name=1, anchor=center, inner sep=0}, "{R_*}"{description}, curve={height=-12pt}, from=1-3, to=1-1]
	\arrow["\Lan", curve={height=-6pt}, from=1-3, to=2-3]
	\arrow["{\mathrm{restr}}"', curve={height=6pt}, from=2-1, to=1-1]
	\arrow["{\Sh\, L_*}"', curve={height=-12pt}, from=2-1, to=2-3]
	\arrow[curve={height=6pt}, from=2-1, to=3-1]
	\arrow["{\mathrm{restr}}", curve={height=-6pt}, from=2-3, to=1-3]
	\arrow["{R_*}"{description}, curve={height=-12pt}, from=2-3, to=2-1]
	\arrow[curve={height=6pt}, from=2-3, to=3-3]
	\arrow[curve={height=6pt}, from=3-1, to=2-1]
	\arrow[curve={height=6pt}, from=3-3, to=2-3]
	\arrow[curve={height=-18pt}, from=4-1, to=4-3]
	\arrow[curve={height=-18pt}, from=4-3, to=4-1]
	\arrow["\dashv"{anchor=center, rotate=-90}, draw=none, from=0, to=1]
\end{tikzcd}\end{center}
The commutativity follows from the commutativity of the $R_*$ and forgetting, and uniqueness of left adjoints.
	 To see that the colimits still form an adjunction, see~\ref{lem:tower-adjs}.
 \end{proof}
\begin{warning}
	Note that we do not claim that the global right adjoint is given by $R_*$, as this a priori has no reason to commute with left Kan extension.
	The left adjoint in turn can always be described as $\Sh L_* ?=\cond(L)$, as left adjoints preserve left Kan extensions.
	If however $R_*?$ preserves the Kan extensions $\Lan_{\tau} F$, then the right adjoint on $\cond$ may be computed as $R_*?$.
\end{warning}
\begin{lemma}
	If $\mcC$ is cocomplete and filtered colimits commute with finite products, then $\cond(\mcC)$ is cocomplete and colimits can be computed by taking pointwise colimits followed by sheafification.

	If for all limit cardinals $\kappa$, $\kappa$-filtered colimits in $\mcC$ exist and commute with $\kappa$-small limits,
	then limits in $\cond(\mcC)$ exist as soon as the pointwise limits exist and can be computed pointwise.
\end{lemma}
\begin{proof}
	As sheafification is cocontinuous, it suffices to show that we can compute the colimits pointwise in $\varinjlim [\extr_\kappa^{\op},\mcC]$.
	Take any diagram $D\colon \mcI\to\Cond(\mcC)$ in $\Cond(\mcC)$ and choose $\kappa\ge \kappa_i$ for $\kappa_i=\min\{\lambda: D_i\in \cond_\lambda\}$ and all $i$ (using that the diagram is small).
	By forgetting, we can hence interpret the diagram in $[\extr^{\op}_\kappa, \mcC]$; we can compute the colimit pointwise there.
	As the left Kan extension commutes with colimits, and furthermore preserves pointwise colimits (as it is a pointwise colimit),
	this implies that in $\varinjlim  [\extr_\kappa^{\op},\mcC]$, colimits may be computed pointwise.

	For the calculation of limits take any diagram $D$, and take $\kappa$ to be bigger than all $\kappa_i$ to interpret the diagram in $\cond_\kappa$ and
	furthermore take $\kappa$ such that is has larger cofinality than $|I|$ (use, e.g., $\beth_{\max\{\kappa,|I|\}}$).
	Note that again, it suffices to check that for this $\kappa$, the $\Lan_{\tau}$ (for $\tau\colon \extr_\kappa\to \extr$) commutes with
	pointwise limits of size $|I|$.

	But this follows by looking at the Kan extension formula:
	\[\varprojlim_{I}(\Lan_{\tau}D_i(S))=\varprojlim_I\varinjlim_{S\to K}D_i(K),\]
	noting that the index category is $|I|$-filtered, and that such colimits commute with $|I|$-small limits,
	we obtain equality to
	\[\varinjlim_{S\to K}\varprojlim_I D_i(K)=(\Lan_{\tau}\varprojlim D_i)(S).\]
	\end{proof}
\begin{remark}
  Note that we do not know whether the condition \enquote{$\kappa$-filtered colimits commute with $\kappa$-small limits for all $\kappa$}
  in the last lemma is already equivalent to just assuming that filtered colimits are exact.
\end{remark}
\begin{corollary}
	Assume that $\kappa$-filtered colimits in $\mcC$ exist and commute with $\kappa$-small limits.
	If filtered colimits in $\mcC$ are $\mcM$-exact (commute with limits of forms in $\mcM$), then sheafification is $\mcM$-continuous.
\end{corollary}

Having \enquote{condensation} as some kind of process of adding topology to an arbitrary category, we note that \idx{algebraic theories}, as introduced in~\ref{sec:univ-alg},
form a quite similar construction.
These form a method of talking about \idx{algebric objects} in an arbitrary category.

When given a reasonably well behaved category $\mcC$,
we now know how to obtain a category ($\cond(\mcC)$) whose objects we can imagine as objects of $\mcC$ with a topological structure (and continuous $\mcC$-maps).
The next step will be to additionally endow the objects of $\cond(\mcC)$ with algebraic structure,
e.g., with the structure of an abelian group,
to obtain condensed abelian groups in $\mcC$.
\begin{definition}[Condensed algebraic theory]\label{def:condensed_algebraic_theory}\uses{def:lawvere_theory_ms, def:cond_C}\chapfour
For any category $\mcC$ where filtered colimits exist and commute with finite products,
and any algebraic theory $\mcT$ with filtered colimits commuting with finite products in $\mcT$,
the category of $\mcT$-objects in $\cond(\mcC)$ are precisely given by the condensation of the $\mcT$-objects in $\mcC$.
    \[\mcT(\cond(\mcC))\simeq \cond(\mcT(\mcC))\]
\end{definition}
\begin{proof}
	First we show that the classes of product preserving functors agree.

	Currying in the large closed cartesian category $\Cat$ yields
	\[\Fun(\extr_\kappa^{\op},\Fun(\mcT,\mcC))\simeq \Fun(\extr_\kappa^{\op}\times \mcT, \mcC)\simeq \Fun(\mcT, \Fun(\extr_\kappa^{\op},\mcC)).\]

	As finite products may be taken pointwise in both, $\Fun(\extr_\kappa^{\op}, \mcC)$ and $\Fun(\mcT, \mcC)$,
	The classes of finitely product preserving functors agree, i.e. $T\in \Fun(\extr_\kappa\times\mcT,\mcC)$ is finitely product preserving if and only if
	for any fixed $S\in \extr_\kappa$ and $n\in \mcT$, the induced functors
	\[T^n\colon \extr_\kappa^{\op}\to \mcC\]
	and
	\[T(S)\colon \mcT\to \mcC\]
	are finitely product preserving and furthermore the induced functors
	\[\Fun(\extr^{\op},\mcT(\mcC)), \, S\mapsto T(S)\]
	and
	\[\Fun(\mcT,\cond_\kappa(\mcC)),\, n\mapsto T^n\]
	are finitely product preserving.
	This can be seen by considering
	\[T^{\prod_i n_i^{k_i}}(\prod S_j)=\prod_{i,j} T^{n_i}(S_j).\]

	Hence
	\[\cond_\kappa(\mcT(\mcC))\simeq \mcT(\cond_\kappa(\mcC)).\]

	Doing the same with $\extr$ instead of $\extr_\kappa$, we see that the classes of product preserving functors agree, and it just remains to check the Kan condition.

	Now we show that the Kan-condition works well.

	For this, take any Kan extension of some $T\in \cond(\mcT(\mcC))$ along $\tau\colon \extr_\kappa\to \extr$.
	As filtered colimits in $\mcT(\mcC)$ are computed pointwise,
	\[(\Lan_{\tau} T)^n(S)=(\varinjlim_{S\to S'}T(S'))^n= \varinjlim_{S\to S'} T^n(S')\]
	implying that if $T\in \cond(\mcT(\mcC))$, then $T\in \mcT(\cond(\mcC))$.
	Conversely, use a large enough $\kappa$ such that all $T^n$ are left Kan extension from the same $\extr_\kappa$, and read the above equation in the converse direction.
	Note that this $\kappa$ can be found, as $\mcT$ is essentially small.
\end{proof}

\begin{corollary}
    The underlying condensed set functor $\cT\cond\to\cond$
    has a left adjoint.
    In particular,
    sheafification of $\cond\cT(\Set)$ may be computed in $\cond$.
\end{corollary}

\begin{remark}
  From this,
  one can go searching for properties of $\mcC$ and $\mcT$ that imply that $\cT\cond(\mcC)$ has enough (compact) projectives
  and that in the case $\mcC=\Set$,
  such a class of generators are given by the free structures on $\extr$.
\end{remark}

Next, we will slowly drift towards condensing abelian categories.
Note that in semiadditive categories, finite products commute with filtered colimits, as they can be identified with finite coproducts.

\begin{proposition}[Condensed semiadditive pointwise]\label{prop:cond_semiadd_pointw}\uses{def:cond_C}
    For any semiadditive category $\mcA$ where $\kappa$-filtered colimits exist and commute with $\kappa$-small limits,
    all colimits in $\cond(\mcA)$ may be computed pointwise (limits may be computed pointwise anyways).

    In particular it remains semiadditive and epi-/monomorphisms are precisely the pointwise epics/monics.
\end{proposition}
    \begin{proof}
  The sheaves are precisely the functors from $\extr$ commuting with finite products.
  But as finite products agree with finite coproducts (here we use the semiadditivity),
  this condition clearly commutes with all limits and all colimits.
  This means that no sheafification is needed, and thus the colimits may be computed pointwise.
\end{proof}
\begin{corollary}[Condensed abelian categories]\label{prop:cond_ab_cats}\uses{def:cond_C}
    If $\mcA$ is an abelian category, then $\cond(\mcA)$ is abelian.
Furthermore, $\cond(\mcA)$ fulfills the same Grothendieck axioms as $\mcA$.
\end{corollary}

\begin{proposition}[Condensed monoidal categories]\label{prop:cond_monoidal}\uses{def:cond_C}
    If $\mcA$ is (symmetric) monoidal, then $\cond(\mcA)$ is (symmetric) monoidal and sheafification is (symmetric) monoidal.
\end{proposition}
\begin{proof}
  Every $\condk(\mcA)$ is (symmetric) monoidal by~\ref{lem:sh-mon}.
  These monoidal structures survive the passage to the colimit $\condk$.
\end{proof}
\begin{question}
    If furthermore the category $\mcA$ admits an internal $\hom$, does $\cond(\mcA)$ admit an internal $\hom$ as well?
\end{question}
\begin{proposition}
  If $\mcA$ is generated by compact projectives,
  then $\cond(\mcA)$ is also generated under compact projectives (under finite colimits).
\end{proposition}
\begin{proof}
  This is \cite[11.8]{scholze2019Analytic}.
\end{proof}

This implies that for abelian $\mcA$ with enough compact projectives,
$\cond(\mcA)$ also has such and thus $D(\cond(\mcA))$ can be constructed neatly.
In general,
for any abelian $\mcA$,
$\cond(\mcA)$ is abelian as well,
there are two canonical ways to construct a \enquote{condensed derived category},
namely, one can derive $\cond(\mcA)$ or condense $D(\mcA)$.
However, it is not true that $D(\cond(\mcA))\simeq \cond(D(\mcA))$.
This defect is one of the possible reasons to pass to the $\infty$-categorical setting (see \cite[chap.~11]{scholze2019Analytic}).
\begin{definition}[Condensing $\infty$-categories]\label{def:cond_infty}
    Let $\mcC$ be an $\infty$-category with all small colimits.
For $\kappa$ being an strong limit cardinal, define $\cond_\kappa(\mcC)$ as the $\infty$-category of contravariant functors $\extr_\kappa\to \mcC$.
Define the \idx{$\infty$-condensation} $\cond(\mcC)$ by
    \[\cond(\mcC)=\varinjlim_{\kappa}\cond_\kappa(\mcC).\]
\end{definition}

\begin{lemma}\label{lem:cond-enouh-proj}
  Condensing commutes with animation.
This means for any category $\mcC$ which is generated by $\mcC^{\mathrm{cp}}$, $\cond(\mcC)$ is also generated under compact projectives and
\[\Ani(\cond(\mcC))\simeq \cond(\Ani(\mcC)).\]
    \end{lemma}
\begin{proof} This is \cite[11.8]{scholze2019Analytic}, and essentially is the $\infty$-categorical version of condensation commuting with (general) algebraic theories.
    \end{proof}
\begin{example}
    \[\mcD_{\ge 0}(\cond(\Ab))=\Ani(\cond(\Ab))=\cond(\Ani(\Ab)).\]
    This induces a $\infty$-category $\mcD_{\ge 0}(\mcA)$ of $\mcA$-modules in $\Ani(\cond(\Ab))$ for any Ring in $\Ani$.
    This can be extended to a stable $\infty$-category $\mcD(\mcA)$, see 12.0/12.5 in \cite{scholze2019Analytic} or Appendix C in \cite{SAG}.
\end{example}


\section{Condensed modules}

When given a reasonably well behaved category $\mcC$,
we now know how to obtain a category ($\cond(\mcC)$) whose objects we can imagine as objects of $\mcC$ with a topological structure (and continuous $\mcC$-maps).
Now, we will further investigate special cases of additionally endowing the objects of $\cond(\mcC)$ with algebraic structure,
e.g., with the structure of an abelian group,
to obtain condensed abelian groups in $\mcC$.
In particular,
$\cond(\Ab(\Set))=\Ab(\cond(\Set))$,
i.e., condensed abelian groups are the same thing as abelian group objects in $\Cond$
(and analogously for other algebraic structures such as rings, modules, etc.).

In this section,
we lay out the delightful theory of condensed modules.
Naturally,
we can form condensed abelian groups.
There are at least two equivalent ways of viewing condensed abelian groups:
as abelian group objects\footnote{See section \ref{sec:univ-alg} for an in depth explanation of this notion.
  For large parts of the present section,
  one can simply think of abelian group objects as objects $A$ together with three maps,
  the unit $0\colon *\to A$,
  negation $-\colon A\to A$ and multiplication $m\colon A\times A\to A$
  fulfilling translations of the usual axioms of ablien groups.
}
in $\cond$ (i.e., $\Ab(\cond)$) or as, literally, condensed abelian groups,
i.e., accessible sheaves on $\extr$ with values in $\Ab$ (i.e., $\cond(\Ab)$),
and these two notions coincide (\ref{def:condensed_algebraic_theory}):
\[
  \cond(\Ab)\simeq\Ab(\cond).
\]
We call this category $\cab$.
In fact,
we will immediately consider condensed $R$-modules
where $R\in\Ring$ is a classical ring.
In the first subsection,
we investigate the categories $\cond(\ModR)$ for rings $R$
and discover that they are extremely well-behaved.
One can work with them almost as comfortably as with classical modules.
This, of course, forms a stark contrast to the setting of topological modules which have all sorts of problems.
In the commutative case,
we also get a nice tensor product and an internal $\hom$.
The tensor product is given as sheafification of $S\mapsto M(S)\otimes N(S)$.
(An even more explicit description will be given shortly.)

Clearly,
as functional analysts,
we do not want our rings to be discrete.
(After all,
$\R$ and $\C$ are of utmost importance for us and we certainly do not want to endow them with a discrete structure.)
Building especially on the case $R=\Z$,
we thus discuss modules over condensed rings.
What is a condensed ring?
Well, what is a ring?
If one thinks of a ring as a set with some operations satisfying some equations,
then one can view condensed rings as condensed sets together with some operations satisfying some equations.
In that case,
we are effectively talking about ring objects in $\cond$ (i.e., $\Ring(\cond)$).
Again,
this is the same as $\cond(\Ring)$.
If, on the other hand,
one thinks of a ring as an abelian group $R$ together with homomorphisms $R\otimes R\to R$ etc.
satisfying the axioms of a monoid,
then one can think so of condensed rings:
monoids in $\cab$ with respect to $\otimes$.
These two (or three) notions coincide:
\[
  \cond(\Ring)\simeq\Ring(\cond)\simeq\Mon(\cab,\otimes).
\]
We call this category $\CondRing$.
In this category we find $\R$, $\C$, $\C[T]$ and alike
but it also includes $\Z$ and $\R_{\textrm{disc}}$.
The viewpoint $\CondRing=\Mon(\cab,\otimes)$ is best when considering modules.
After all:
What is a condensed module over a condensed ring $R\in\CondRing$?
It's a condensed abelian group $M$ together with a homomorphism $R\otimes M\to M$ satisfying certain equations (that translate to commutativity of some diagrams).
This notion was studied in $\ref{ssec:mon-cats}$
We will call the resulting category $\CondModR$.
Again,
there will be a tensor product in the commutative case given as sheafification of $S\mapsto M(S)\otimes_{R(S)} N(S)$
with respect to which $\CondModR$ is closed symmetric monoidal.
As in classical algebra,
a homomorphism of condensed rings $R'\to R$ endows any $R$-module (in particular $R$ itself)
with the structure of an $R'$-module,
inducing a functor $\CondModR\to\CondMod_{R'}$.
For $R$, $R'$ commutative,
this functor has a left adjoint given by $M\mapsto M\otimes_{R'}R$ and a right adjoint $M\mapsto\ihom_{R'}(R,-)$.
The case $R_{d}\to R$ will be of importance soon.

In the pursuit of building up the theory,
we will want to consider different points of view when it comes to modules.
For example,
seeing them as algebras over the monad $R\otimes-$ will be useful.
For now,
however,
we want to focus on an important point about condensed modules over discrete rings
and will try to keep this discussion as elementary as possible.\footnote{
  To this end,
  we will also ignore \enquote{$\kappa$ issues}.
  In other words,
  we will phrase everything as if condensed sets were just sheaves of sets on $\extr$.
}

Unfolding the definition somewhat,
a condensed module $M\in\cond(\ModR)$ for a classical ring $R$ consists of a sheaf $M$ of abelian groups on $\extr$
together with,
for each $S\in \extr$,
a map $R\otimes M(S)\to M(S)$
defining an $R$-module structure.
(And these need to be compatible along maps $S'\to S$ in $\extr$ but this is besides the point.)

Given a condensed ring $R'\in\CondRing$, a module $M\in\CondModR$ is a sheaf $M$ of abelian groups on $\extr$
together with maps $R(S)\otimes M(S)\to M(S)$ ($S\in\extr$).
This is all very reasonable.
There was however one point that bugged us for some time.
Naturally,
one would expect that for $R\in\Ring$ and $R_{d}=R_{\textrm{disc}}\in\CondRing$
the corresponding discrete condensed ring,
the following categories are equivalent:
$\CondMod_{R_{d}}\simeq\cond(\ModR)$.
And in fact,
this turned out to be the point.
What confused us,
was that
by definition,
for any $M\in\CondMod_{R_{d}}$ and $S\in\extr$,
$M(S)$ has a given $R_{d}(S)$-module structure,
but an $M\in\cond(\ModR)$ a priori only has a $R=R_{d}(*)$-module structure.
In fact,
it will turn out that the $R_{d}(S)$-module structure comes from the sheaf condition (see~\ref{rem:rs-mod-magic}).

\subsection{Condensing classical modules}

We collect the information of the last section specialised to $\cond(\ModR)$
and pave the way for the study of condensed $R$-modules.

We can combine the theorems of the last section to categories of the nicest possible sort -- algebraic theories over set that are abelian categories.
However, these already form a quite restrictive class of examples.

\begin{proposition}
	Any category of single-sorted algebraic theory objects over set $\mcT(\Set)$ that also forms an abelian category is automatically of the form $\ModR$ for a (not necessarily commutative) ring $R$.
\end{proposition}
\begin{proof}
  This is~\cite[B.8]{Adamek2010a}.
\end{proof}

Because of this,
we will now focus on condensed $R$-modules.
As explained in the introduction to this section,
there are two notions of \enquote{condensed $R$-modules}.
In this subsection,
we will take a look at the case of classical rings.
I.e.,
we will investigate the categories $\cond(\ModR)$ for classical (purely algebraic) rings $R$,
leaving the case of modules over condensed rings for later.

As an important special case,
we have $R=\Z$.
Since $\Z$-modules are better known as abelian groups,
we will usually write $\cab$ instead of $\cond(\Ab)=\cond(\Mod_{\Z})$.

\begin{corollary}[of \ref{lem:cond_as_colimit}, \ref{def:condensed_algebraic_theory} and \ref{prop:alg-th-is-monad}]
  The category $\cond(\ModR)$ of \idx{condensed $R$-modules} is equivalently given as any of the following.
  \begin{enumerate}[(a)]
    \item The category of contravariant functors
    \[
	  M\colon \extr\to \ModR
	\]
    mapping finite coproducts to finite products
	such that there exists a strong limit cardinal $\kappa$ with
    \[
	  M(X)=\varinjlim_{X\to S,\, S\in \extr_\kappa} M(S)
	\]
    for all $X\in \extr$ (where the colimit can be taken equivalently in $\ModR$ or $\Set$).
    \item The category
          \[\varinjlim \cond_\kappa(\ModR).\]
    \item The category of $R$-module objects in $\cond$ in the sense of the (single-sorted) Lawvere theory of $R$-modules,
          see~\ref{rem:expl-law-th}.
    \item The (Eilenberg-Moore) category of algebras for the monad $R[-]$,
          where $R[X]$ ($X\in\cond$) is
          given by sheafification of $S\mapsto R[X(S)]$,
          the free $R$-module on $X(S)$.
  \end{enumerate}
\end{corollary}

\begin{corollary}[of \ref{lem:sheafification_cond_C}, \ref{prop:monad-create-limits}, \ref{prop:monad-create-colimits} and \ref{prop:free-r-sift}]
  Sheafification in $\cond(\ModR)$ exists and agrees with the sheafification in $\cond$,
  i.e.,
  for any presheaf $M$ in $[\extr^{\op}, \ModR]$ that is the Kan extension of its restriction to $\extr_{\kappa}$ for some strong limit cardinal $\kappa$,
  one can compute the left adjoint to $\cond(\ModR)\to \varinjlim[\extr^{\op}_\kappa, \ModR]$ by forgetting the $R$-module structure on each $M(S)$ and
  using sheafification in $\cond$,

\begin{center}\begin{tikzcd}
	{\varinjlim [\extr_\kappa^\op, \ModR]} & {\varinjlim [\extr_\kappa^\op, \Set]} \\
	\cond(\ModR) & \cond
	\arrow["{?}", from=1-1, to=1-2]
	\arrow[""{name=0, anchor=center, inner sep=0}, "\Sh"', curve={height=12pt}, from=1-1, to=2-1]
	\arrow[""{name=1, anchor=center, inner sep=0}, "\Sh"', curve={height=12pt}, from=1-2, to=2-2]
	\arrow[""{name=2, anchor=center, inner sep=0}, "{?}"', curve={height=12pt}, from=2-1, to=1-1]
	\arrow["{?}"', from=2-1, to=2-2]
	\arrow[""{name=3, anchor=center, inner sep=0}, "{?}"', curve={height=12pt}, from=2-2, to=1-2]
	\arrow["\dashv"{anchor=center}, draw=none, from=0, to=2]
	\arrow["\dashv"{anchor=center}, draw=none, from=1, to=3]
\end{tikzcd}\end{center}

	On $\extr$, we can compute the sheafification explicitly by
	\[\Sh(M)(S)=\varinjlim_{S=\sqcup_{j=1}^k T_j}\prod_{j=1}^k M(T_j).\]
\end{corollary}

\begin{corollary}[of \ref{lem:sh_exact}]
    The sheafification functor
    \[\varinjlim [\extr^{\op}_\kappa, \ModR]\to \cond(\ModR)\]
    is exact.
\end{corollary}

\begin{corollary}[of \ref{prop:cond_semiadd_pointw}]\label{lem:lim_incondAb}
    All limits and colimits in $\cond(\ModR)$ exist and can be taken pointwise, i.e., in $[\extr^\op, \ModR]$.
    In particular,
	epi-/monomorphisms are precisely the pointwise epi-/monomorphisms.
\end{corollary}

\begin{corollary}[of \ref{prop:cond_ab_cats}]
  $\cond(\ModR)$ is a bicomplete abelian category (AB3, AB3*),
  coproducts are exact (AB4) and filtered colimits are exact (AB5).
\end{corollary}

\begin{corollary}[of \ref{prop:cond_ab_cats}]
    $\condAb$ is a bicomplete (AB3, AB3*) abelian category, all coproducts and products are exact (AB4, AB4*), filtered colimits are exact (AB5) and filtered colimits distribute over products (AB6).
\end{corollary}

In particular this implies that
\begin{itemize}
    \item $\cond(\ModR)$ is balanced.
    \item There is an image factorisation.
    \item All epics and monics are regular, and thus the regular projective objects and projective objects agree.
    \item Pullbacks of epics are epics and pushouts. Pushouts of monics are monics and pullbacks.
    \item The notions of densely generating and separating agree.
\end{itemize}

\begin{proposition}
  The forgetful functor $\cond(\ModR)\to\cond$ is right adjoint.
  Its left adjoint $R[-]$ is given on $X\in\cond$ as the sheafification of
  \[
	S\mapsto R[\hom(S,X)]=R[X(S)]=\bigoplus_{X(S)}R,
  \]
  for short: $R[X]=\Sh(R[X(-)])$.
\end{proposition}
\begin{proof}
  This is trivially trivial,
  for example as we have already seen that $\cond(\ModR)$ is the Eilenberg-Moore category for exactly this monad.

  We also give an ad hoc argument.
  Since $\ModR\to\Set$ creates limits (it is preserves them as a right adjoint,
  but both $\ModR$ and $\Set$ are complete),
  the sheaf condition of any given object of $\varinjlim_{\kappa}[\extr^{\op}_{\kappa},\ModR]$
  can be checked in $\varinjlim_{\kappa}[\extr^{\op}_{\kappa},\Set]$.
  In other words,
  the composite
  \[
	\cond(\ModR)\to \varinjlim_{\kappa}[\extr^{\op}_{\kappa},\ModR]\to \varinjlim_{\kappa}[\extr^{\op}_{\kappa},\Set]
  \]
  lands in $\cond$.
  But then the claim follows
  because both $\cond(\ModR)\to\varinjlim_{\kappa}[\extr^{\op}_{\kappa},\ModR]$ and $[\extr^{\op}_{\kappa},\ModR]\to [\extr^{\op}_{\kappa},\Set]$
  have left adjoints.
  As these are given by sheafification resp.\ $F\mapsto(T\mapsto R[F(T)])$,
  we get the given explicit description.
\end{proof}

\begin{proposition}\label{prop:r-con-enogh-cproj}
  For $S\in\extr$,
  $\hom_{\cond(\ModR)}(R[S],-)$ commutes with all limits and colimits.
  In particular,
  $R[S]$ is compact projective.
  Indeed, this yields a class of compact projective generators.
\end{proposition}
\begin{proof}
  As all limits and colimits in $\cond(\ModR)$ can be taken pointwise,
  the first assertion can be deduced by the usual Yoneda argument:
  Let $\lim M_{i}\in\cond(\ModR)$ be any limit or colimit.
  Then
  \begin{align*}
    \hom_{\cond(\ModR)}(R[S],\lim M_{i})
    &=\hom_{\cond}(S,?\lim M_{i})
      =(?\lim M_{i})(S)\\
    &=(\lim M_{i})(S)=\lim M_{i}(S)=\lim (?M_{i})(S)\\
    &=\lim\hom_{\cond}(S,?M_{i})=\lim\hom_{\cond(\ModR)}(R[S],M_{i}).
  \end{align*}
  One can even avoid Yoneda and just use compact projectivity of $S$ (i.e., that $\hom(S,-)$ commutes with sifted colimits):
  Given any sifted colimit $\varinjlim M_{i}\in\cond(\ModR)$,
  we can compute
  \begin{align*}
	\hom_{\cond(\ModR)}(R[S],\varinjlim M_{i})
	&=\hom_{\cond}(S,?\varinjlim M_{i})
	  =\hom_{\cond}(S,\varinjlim ?M_{i})\\
	&=\varinjlim\hom_{\cond}(S,?M_{i})
	  =\varinjlim\hom_{\cond(\ModR)}(R[S],M_{i}).
  \end{align*}
  But then,
  take a look at~\ref{rem:tiny-in-ab}.

  Another argument could be constructed along the lines~\ref{lem:adjoints_preserve_projective_objects}.

  So let us turn to the last claim.
  we give to arguments.
  One can be found in~\ref{lem:cond-enouh-proj}.
  Alternatively,
  we show that it is a separating class.
  For this,
  let $f\neq g\colon M\to N$ in $\cond(\ModR)$.
  Because $\extr$ is separating in $\cond$ (and because of the faithfulness of $\cond(\ModR)\to \cond$),
  there is map $h\colon S\to M$ for which $fh\neq gh$.
  But then the induced map $R[S]\to M$ also distinguishes $f$ and $g$ by the universal property of $R[S]$
  (i.e., by the adjunction with the forgetful functor).
\end{proof}

\begin{proposition}\label{prop:trace-pie}
  For $M\in\cond(\ModR)$,
  the evaluation (or trace) map $R[M]\to M$ is a split epimorphism.
\end{proposition}
\begin{proof}
  This is an immediate consequence of the definition of an algebra over a monad.
  (The triangle in~\ref{def:alg-over-monad} exhibits its split epimorphicity.)

  We also give an elementary proof here.
  Obviously,
  $R[M(S)]\to M(S)$ is epic for all $S$.
  Therefore $(S\mapsto R[X(S)])=:R_{\textrm{PSh}}[M]\to M$ is epic in $\varinjlim[\extr_{\kappa}^{\op},\Set]$.
  But then $R[M]\to M$ is also epic as a map of presheaves
  because it factors $R_{\textrm{PSh}}[M]\to M$:
  \begin{center}
    \begin{tikzcd}
      R_{\textrm{PSh}}[M] \ar[r,two heads] \ar[d] &  M\\
      R[M]\ar[ur] &
    \end{tikzcd}
  \end{center}
  (after all, $M$ is a sheaf).
  As such,
  it is of course also epic as a map in $\cond(\ModR)$.
  Along these lines,
  one can also find the section.
\end{proof}

\begin{lemma}\label{lem:T1Ab_to_CondAb} For any T1 topological $R$-module $A$ (and hence Hausdorff), the corresponding condensed set is pointwise an $R$-module,
  and thus
$\underline{M}$ is a condensed $R$-module.

The embedding $M\mapsto \underline{M}$ commutes with all limits and finite coproducts.
\end{lemma}
\begin{proof}
The first assertion is clear, as we can describe condensed $R$-modules as sheaves with values in $\ModR$.
The continuity claims follow by comparing limits/colimits to the ones in condensed $R$-modules and T1 topological $R$-modules.
\end{proof}

\begin{corollary}\chapfour
  Any condensed $R$-module is cokernel of some map $\bigoplus R[T_i]\to \bigoplus R[S_j]$ ($T_{i},S_{j}\in\extr$).
\end{corollary}
\begin{proof}
  As shown, for any $A\in \cond(\ModR)$,
  there exists an epimorphism $\bigoplus R[S_j]\twoheadrightarrow A$,
  which is effective and hence the cokernel of its kernel $B\hookrightarrow \bigoplus R[S_j]\twoheadrightarrow A$.
  Choosing an epimorphism $\bigoplus R[T_i]\twoheadrightarrow B$ yields the correct morphism $\bigoplus R[T_i]\twoheadrightarrow B\hookrightarrow \bigoplus R[S_j]$.

  For another simple argument,
  write $A=\varinjlim R[S_{i}]$ and write this colimit as a coequalizer of coproducts as in~\ref{lem:lim_via_special_lim}.
\end{proof}

When $R$ is commutative,
it is to be expected that the closed symmetric monoidal structure on $\ModR$ can survive the condensation process.
Of course,
we already know that symmetric monoidal structures survive condensation (\ref{prop:cond_monoidal}),
but here,
it does so particularly smoothly.

\begin{theorem}\label{lem:computation_tensor_hom}\chapfour
  Let $R$ be a commutative ring.
  Then the category $\cond(\ModR)$ has a symmetric monoidal structure $\otimes_{R}$
  with respect to which $R[X]$ is flat for all $X\in\cond$,
  which is cocontinuous in both arguments (separately)
  and for which $R[-]\colon (\cond,\times)\to(\cond(\ModR),\otimes_{R})$
  is symmetric monoidal.
  That is to say: $R[X\times Y]=R[X]\otimes_{R} R[Y]$.
  It is uniquely determined by the latter two properties.
  Furthermore,
  this tensor product represents bilinear maps and its tensor unit is $R[\ast]=R_{d}$.

  In fact,
  this symmetric monoidal structure is closed, i.e.,
  there is an internal $\hom$ functor
  \[
	\ihom\colon\quad\cond(\ModR)^{\op}\times\cond(\ModR)\to\cond(\ModR)
  \]
  which is partially right adjoint to $\otimes_{R}$,
  meaning that for all $A,B,C\in\cond(\ModR)$
  the following holds naturally in each coordinate:
  \[
	\hom_{\cond(\ModR)}(A\otimes_{R} B,C)=\hom_{\cond(\ModR)}(A,\ihom(B,C)).
  \]

  In general,
  the tensor product $A\otimes_{R} B$ can be calculated as sheafification of the pointwise tensor product or more explicitly by
  \[
	(A\otimes_{R} B)(S)=A(S)\otimes_{R_{d}(S)}B(S)
  \]
  with $R_{d}(S)=C(S,R_{\mathrm{disc}})$ ($S\in\extr$).
  The internal hom $\ihom(A,B)$ can be computed by
  \[
	\ihom_{R}(A,B)(S)=\hom_{\cond(\ModR)}(R[S]\otimes_{R} A,B)
  \]
  for any $S\in \extr$.

  Lastly, sheafification is symmetric monoidal, meaning for any presheaves $A,B\in\varinjlim[\extr_{\kappa}^{\op},\ModR]$,
	\[\Sh(A)\otimes_{\cond(\ModR)}\Sh(B)=\Sh(A\otimes_{\PSh}B).\]
\end{theorem}

The rest of this subsection will be devoted to proving this theorem.

\begin{remark}
  We will often use the term \enquote{accessible} for sheaves and presheaves $A\in\Sh(\extr,\ModR)$ or $B\in[\extr^{\op},\ModR]$.
  By this,
  we mean that they can be written as (small) colimits of $R[T]$ ($T\in\extr$) resp.\ $S\mapsto R[T(S)]$ ($T\in\extr$)
  and this is equivalent to $A$ resp.\ $B$ lying in $\varinjlim[\extr_{\kappa}^{\op},\ModR]$.
\end{remark}

\begin{definition}\label{def:sym-mon-on-psh}
  Let $R$ be a commutative ring.
  We define two tensor products on $[\extr^{\op},\ModR]$,
  \[
    A\otimes_{R,\PSh}B\coloneqq (S\mapsto A(S)\otimes_{R}^{\mathrm{alg}}B(S)),\qquad
    A\otimes_{R}B\coloneqq \Sh(A\otimes_{R,\PSh}B).
  \]
  Clearly, both define symmetric monoidal structures on $[\extr^{\op},\ModR]$.

  We also the define two functors
  \[
    \ihom_{R,\PSh},\ihom_{R}\colon[\extr^{\op},\ModR]^{\op}\times [\extr^{\op},\ModR]\to [\extr^{\op}]
  \]
  by
  \begin{align*}
    \ihom_{R,\PSh}(A,B)(S)&\coloneqq \hom_{[\extr^{\op},\ModR]}(A\otimes_{R,\PSh}R[S],B),\\
    \ihom_{R}(A,B)(S)&\coloneqq \hom_{[\extr^{\op},\ModR]}(A\otimes_{R}R[S],B).
  \end{align*}
\end{definition}

\begin{remark}\label{rem:tp-ihom-cocont}
  Because the algebraic tensor product $\otimes_{R}^{\mathrm{alg}}$ is componentwise cocontinuous,
  $\otimes_{R,\PSh}$ is.
  But sheafification is also cocontinuous,
  so the same holds for $\otimes_{R}$.

  As all limits and colimits can be computed pointwise,
  both $\ihom_{R,\PSh}$ and $\ihom_{R}$ are (co)continuous in the usual way.
\end{remark}

\begin{lemma}\label{lem:acc-ot-acc}
  For accessible presheaves $A,B\in\varinjlim[\extr_{\kappa}^{\op},\ModR]$,
  the tensor product $A\otimes_{R}B$ is also accessible.
\end{lemma}
\begin{proof}
  This is a direct consequence of the last remark.
  Just write both $A$ and $B$ as colimits of $R[S_{i}]$'s resp.\ $R[T_{j}]$'s
  and then pull out both colimits in $A\otimes_{R}B$.
  And for $S_{i},T_{j}\in\extr$,
  $R[S_{i}\times T_{j}]=R[S_{i}]\otimes_{R}R[T_{j}]$,
  can be verified explicitly.
  Indeed,
  for any $S,\,T\in\cond$ and $k\in\extr$,
  we have by the sheafification formula (used in the beginning and at the end)
  \begin{align*}
    &(R[S]\otimes_{R} R[T])(K)\\
    &=\varinjlim_{\sqcup_{\textrm{fin.}}K_{i}=K}\bigoplus_{i}R[S](K_{i})\otimes_{R}^{\textrm{alg}}R[T](K_{i})\\
    &=\varinjlim_{\sqcup_{\textrm{fin.}}K_{i}=K}\bigoplus_{i}
      \left(\left(\varinjlim_{\sqcup_{\textrm{fin.}}K_{ij}=K_{i}}\bigoplus_{j}R[S(K_{ij})]\right)\otimes_{R}^{\textrm{alg}}
      \left(\varinjlim_{\sqcup_{\textrm{fin.}}K'_{ik}=K_{i}}\bigoplus_{k}R[T(K'_{ik})]\right)\right)\\
    &=\varinjlim_{\sqcup_{\textrm{fin.}}K_{i}=K}\bigoplus_{i}\varinjlim_{\sqcup_{\textrm{fin.}}K_{ij}=K_{i}}\,\varinjlim_{\sqcup_{\textrm{fin.}}K'_{ik}=K_{i}}
      \left(\bigoplus_{j}R[S(K_{ij})]\otimes_{R}^{\textrm{alg}}\bigoplus_{k}R[T(K'_{ik})]\right)\\
    &=\varinjlim_{\sqcup_{\textrm{fin.}}K_{i}=K}\bigoplus_{i}\varinjlim_{\sqcup_{\textrm{fin.}}K_{ij}=K_{i}}\,\varinjlim_{\sqcup_{\textrm{fin.}}K'_{ik}=K_{i}}
      \bigoplus_{j,k}R[S(K_{ij})]\otimes_{R}^{\textrm{alg}}R[T(K'_{ik})]\\
    &=\varinjlim_{\sqcup_{\textrm{fin.}}K_{i}=K}\bigoplus_{i}\varinjlim_{\sqcup_{\textrm{fin.}}K_{ij}=K_{i}}\,\varinjlim_{\sqcup_{\textrm{fin.}}K'_{ik}=K_{i}}
      \bigoplus_{j,k}R[S(K_{ij})\times T(K'_{ik})]\\
    &\stackrel{(\ast)}{=}\varinjlim_{\sqcup_{\textrm{fin.}}K_{i}=K}\bigoplus_{i}\varinjlim_{\sqcup_{\textrm{fin.}}K_{ij}=K_{i}}
      \bigoplus_{j}R[S(K_{ij})\times T(K_{ij})]\\
    &=\varinjlim_{\sqcup_{\textrm{fin.}}K_{i}=K}\bigoplus_{i}\varinjlim_{\sqcup_{\textrm{fin.}}K_{ij}=K_{i}}
      \bigoplus_{j}R[S\times T(K_{ij})]\\
    &\stackrel{(\cdot)}{=}\varinjlim_{\sqcup_{\textrm{fin.}}K_{i}=K}\varinjlim_{\substack{(K_{ij})_{i},\\\forall i\colon\sqcup_{\textrm{fin.}}K_{ij}=K_{i}}}
      \bigoplus_{i}\bigoplus_{j}R[S\times T(K_{ij})]\\
    &=\varinjlim_{\sqcup_{\textrm{fin.}}K_{h}=K}\bigoplus_{h}R[S\times T(K_{h})]\\
    &=R[S\times T](K),\\
  \end{align*}
  where at $(\ast)$ we have noticed that the two colimits we combine are over the same index category (call it $I$)
  and the diagonal functor $I\to I\times I$ is final, and at $(\cdot)$ we have used the property AB6 of the classical $\ModR$.

  For $R[S_{i}\times T_{j}]=R[S_{i}]\otimes_{R}R[T_{j}]$,
  one could alternatively argue that both represent bilinear maps in a suitable sense.
\end{proof}
\begin{remark}\label{rem:acc-to-acc-k}
  As is clear from the proof, if both $A$ and $B$ can be written as colimits of $R[S_{i}]$'s resp.\ $R[T_{j}]$'s
  with $R[S_{i}],\,R[T_{j}]\in\extr_{\kappa}$,
  so can $A\otimes_{R} B$.
\end{remark}

\begin{lemma}\label{lem:ihom-ot-adj-on-acc}
  For accessible $A\in\varinjlim[\extr_{\kappa}^{\op},\ModR]$
  and $B,C\in[\extr^{\op},\ModR]$,
  the equation of the hom-tensor-adjunction for $B$ holds, i.e.,
  \[
    \hom_{[\extr^{\op},\ModR]}(A\otimes_{R}B,C)=\hom_{[\extr^{\op},\ModR]}(A,\ihom_{R}(B,C)).
  \]
\end{lemma}
\begin{proof}
  We first prove the statement for the case $A=R[S]$ ($S\in\extr$).
  But then, by Yoneda,
  \begin{align*}
    &\quad \hom_{[\extr^{\op},\ModR]}(R[S]\otimes_{R}B,C)\\
    &=\ihom_{R}(B,C)(S)=\hom_{[\extr^{\op},\Set]}(S,\ihom_{R}(B,C))\\
    &=\hom_{[\extr^{\op},\ModR]}(R[S],\ihom_{R}(B,C)).
  \end{align*}
  Now,
  by~\ref{rem:tp-ihom-cocont},
  both sides of the equation are cocontinuous in $A$, so the claim follows for accessibles.
\end{proof}

\begin{lemma}\label{lem:ihom-yields-acc}
  For $A,B\in\varinjlim[\extr_{\kappa}^{\op},\ModR]$,
  $\ihom_{R}(A,B)$ is accessible
  and $\ihom_{R}(A,-)$ is right adjoint to $-\otimes_{R}A$ as functors from and to $\varinjlim[\extr_{\kappa}^{\op},\ModR]$.
\end{lemma}
\begin{proof}
  Let $A$ be the left Kan extension from its restriction to $\extr_{\kappa}$.
  Then for $\kappa'>\kappa$,
  $-\otimes_{R}A$ is a functor $[\extr_{\kappa}^{\op},\ModR]\to [\extr_{\kappa}^{\op},\ModR]$ by~\ref{rem:acc-to-acc-k} and
  similarly for $\ihom_{R}$.
  But by~\ref{lem:ihom-ot-adj-on-acc},
  these functors are adjoint.
  Thus,
  \[
    -\otimes_{R}A\colon\quad\varinjlim[\extr_{\kappa}^{\op},\ModR]\to\varinjlim[\extr_{\kappa}^{\op},\ModR]
  \]
  has a right adjoint.
  But $\ihom_{R}(A,-)$ fulfills the adjunction formula!
  Therefore,
  by Yoneda,
  $\ihom_{R}(A,-)$ actually has to send accessibles to accessibles.
\end{proof}

\begin{lemma}\label{lem:ihom-sheaf}
  For any $A,B\in[\extr^{\op},\ModR]$,
  $\ihom_{R}(A,B)$ is a sheaf.
\end{lemma}
\begin{proof}
  We simply use that $R[-]$,
  as any left adjoint,
  preserves colimits and~\ref{rem:tp-ihom-cocont}:
  \begin{align*}
    &\quad\ihom_{R}(A,B)(S\sqcup T)
      =\hom_{[\extr^{\op},\ModR]}(A\otimes_{R}R[S\sqcup T],B)\\
    &=\hom_{[\extr^{\op},\ModR]}\big((A\otimes_{R}R[S])\oplus(A\otimes_{R}R[T]),B\big)\\
    &=\hom_{[\extr^{\op},\ModR]}(A\otimes_{R}R[S],B)\times\hom_{[\extr^{\op},\ModR]}(A\otimes_{R}R[S],B)\\
    &=\ihom_{R}(A,B)(S)\times\ihom_{R}(A,B)(T)
  \end{align*}
  for any $S,T\in\extr$.
\end{proof}

\begin{lemma}\label{lem:ihoms-agree}
  For any sheaf $B\in\Sh(\extr,ModR)$
  and any $A\in[\extr^{\op},\ModR]$,
  \[
    \ihom_{R}(A,B)=\ihom_{R,\PSh}(A,B).
  \]
\end{lemma}
\begin{proof}
  As $B$ is a sheaf and $\Sh(\extr,\ModR)$ is a full subcategory of $[\extr^{\op},\ModR]$,
  \begin{align*}
    \ihom_{R}(A,B)(S)
    &=\hom_{[\extr^{\op},\ModR]}(A\otimes_{R}R[S],B)\\
    &=\hom_{\Sh(\extr,\ModR)}(\Sh(A\otimes_{R,\PSh}R[S]),B)\\
    &=\hom_{[\extr^{\op},\ModR]}(A\otimes_{R,\PSh}R[S],B)
      =\ihom_{R,\PSh}(S)
  \end{align*}
  for all $S\in\extr$.
\end{proof}

\begin{lemma}\label{lem:sh-ot-wherever-you-want}
  For any $A,B\in\varinjlim[\extr_{\kappa}^{\op},\ModR]$,
  \[
    A\otimes_{R}B=\Sh(A)\otimes_{R}B
  \]
  and as direct consequences,
  also $A\otimes_{R}B=A\otimes_{R}\Sh(B)=\Sh(A)\otimes_{R}\Sh(B)$.
\end{lemma}
\begin{proof}
  We simply use the adjunctions,~\ref{lem:ihoms-agree} and calculate
  \begin{align*}
      &\quad\hom_{\cond(\ModR)}(A\otimes_{R}B,C)=\hom_{\cond(\ModR)}(\Sh(A\otimes_{R,\PSh} B),C)\\
      &=\hom_{[\extr^{\op},\ModR]}(A\otimes_{R,\PSh} B, C)= \hom_{[\extr^{\op},\ModR]}(A, \ihom_{R,\PSh}(B,C))\\
      &=\hom_{[\extr^{\op},\ModR]}(A,\ihom_{R}(B,C))=\hom_{\cond(\ModR)}(\Sh(A), \ihom_{R}(B,C))\\
      &=\hom_{\cond(\ModR)}(\Sh(A)\otimes B, C),
  \end{align*}
  so that the claim follow by (the dual version of) Yoneda.
\end{proof}

\begin{remark}
  The attentive reader has probably noticed that
  we have also used a $(-\otimes_{R,\PSh}B)$-$(\ihom_{R,\PSh}(B,-))$-adjunction.
  In fact,
  the details for this go through analogously to the case or $\otimes_R$ and $\ihom_{R}$.
\end{remark}

\begin{corollary}\label{cor:sh-symm-mon}
  Sheafification
  \[
    \Sh\colon\quad\big(\varinjlim[\extr_{\kappa}^{\op},\ModR],\otimes_{R,\PSh}\big)\to\big(\cond(\ModR),\otimes_{R}\big)
  \]
  is symmetric monoidal.
\end{corollary}

\begin{lemma}\label{lem:flat}
  For any $X\in\cond$,
  $R[X]$ is flat,
  i.e.,
  $R[X]\otimes_{R}-$ is exact.
\end{lemma}
\begin{proof}
  $R[X]\otimes_{R}-=\Sh(S\mapsto R[X(S)])\otimes_{R}-=\Sh((S\mapsto R[X(S)])\otimes_{R,\PSh}-)$
  is a composite of exact functors.
  Here,
  $R[X(-)]\otimes_{R,\PSh}$ is exact as all $R[X(S)]$ are flat (even free).
\end{proof}

\begin{lemma}\label{lem:r-symm-mon}
  The functor $R[-]\colon(\cond,\times)\to(\cond(\ModR),\otimes_{R})$ is symmetric monoidal.
\end{lemma}
\begin{proof}
  We can simply compute that for any $X,Y\in\cond$,
  \begin{align*}
    R[X]\otimes_{R}R[Y]
    &=\Sh(S\mapsto R[X(S)]\otimes_{R}^{\mathrm{alg}}R[Y(S)])\\
    &=\Sh(S\mapsto R[X(S)\times Y(S)])\\
    &=\Sh(S\mapsto R[(X\times Y)(S)])=R[X\times Y].
  \end{align*}
  as, classically, $R[-]\colon(\Set,\times)\to(\ModR,\otimes_{R}^{\mathrm{alg}})$ is symmetric monoidal
  (and limits in $\cond$ can be computed pointwise).
\end{proof}

\begin{remark}\label{rem:rs-mod-magic}
  In $\varinjlim[\extr^{\op}_{\kappa},\ModR]$ with $\otimes_{R,\PSh}$,
  we have $R_{d}$ as a monoid by pointwise multiplication
  \[
    R_{d}(S)\otimes_{R}^{\mathrm{alg}}R_{d}(S)\to R_{d}(S)
  \]
  given by actual $R_{\mathrm{disc}}$-pointwise multiplication
  \[
    C(S,R_{\mathrm{disc}})\otimes_{R}^{\mathrm{alg}}C(S,R_{\mathrm{disc}})\to C(S,R_{\mathrm{disc}}).
  \]
  As $R_{d}$ is a sheaf and sheafification is symmetric monoidal,
  this induces a $\otimes_{R}$-monoid structure on $R_{d}$.
  Along the same lines,
  we see that any $A\in\cond(\ModR)$ is an $R_{d}$-module.
  In particular,
  it is pointwise an $R_{d}(S)$-module,
  so that the formula in the next lemma makes sense.
\end{remark}

\begin{lemma}\label{lem:formula-tp}
  The tensor product $\otimes_{R}$ on $\cond(\ModR)$ can be calculated by
  \[
    (A\otimes_{R}B)(S)=A(S)\otimes_{R_{d}(S)}^{\mathrm{alg}}B(S)
  \]
  for any $S\in\extr$.
\end{lemma}
\begin{proof}
  As can be checked elementarily,
  $S\mapsto A(S)\otimes_{R_{d}(S)}^{\mathrm{alg}}B(S)$ is already a sheaf.
\end{proof}

\begin{proof}[Proof of~\ref{lem:computation_tensor_hom}]
  Taken together, we have proven all the claims except uniqueness (which is clear) and the tensor unit.
  That $R_{d}$ is the tensor unit is clear by the explicit formula and $R[\ast]=\Sh(S\mapsto R[\ast(S)]=R)=R_{d}$
\end{proof}

\subsection{Modules over a condensed ring}

Of course,
by general nonsense above we have a category $\cab=\cond(\Ab)$ of condensed abelian groups
with a symmetric monoidal tensor product $\otimes=\otimes_{\Z}$.
A condensed ring is nothing but a monoid in this monoidal category (see \ref{def:monoid}),
i.e.,
a condensed abelian group $R$ together with maps $R\otimes R\to R$ and $\Z\to R$ satisfying associativity and unitality.
Again,
by abstract nonsense,
this is the same as ring objects in $\cond$ which in turn is the same as condensed rings.
That is,
\[
  \Mon(\condAb,\otimes)=\Ring(\cond)=\cond(\Ring)=:\CondRing.
\]
So far, so clear.

Now,
given a condensed ring $R$,
there are several obvious notions of condensed $R$-module
that all coincide.
For starters,
using algebraic theories or operads we have a functor
\[
  \CondMod=\cond(\Mod)\to\cond(\Ring)=\CondRing
\]
given by condensing the functor
\begin{align*}
  \Mod&\to\Ring,\\
  (R,M)&\mapsto R,\\
  (\varphi,f)&\mapsto\varphi.
\end{align*}
Recall that a morphism $(\varphi,f)\colon(R,M)\to(S,N)$ in $\Mod$ is given by a ring homomorphism
$\varphi\colon R\to S$ and an abelian group homomorphism $f\colon M\to N$ with $f(rm)=\varphi(r)f(m)$ for all $r\in R$, $m\in M$.
Equivalently,
one could ask for $f$ to be an $S$-module homomorphism $M\otimes_{R}S\to N$.
Following these lines of thought,
we can define $\CondModR$ as the subcategory of $\CondMod$ with objects $(R,-)$ and morphisms $(\id_{R},-)$ --
the fiber over $R$ and $\id_{R}$, so to say.

Thinking more about the monoidal structure on $\cab$,
we could also define an $R$-module to be a condensed abelian group together with an action by $R$.
In fact,
let us do so now to fix terminology and notation, being specialised from \ref{subsubsec:monoids_and_modules}.

\begin{definition}\label{def:condmodr}
  A \idx{condensed module} over a condensed ring $R$ consists of a condensed abelian group $M$ together with an action by $R$.
  Concretely,
  we have a morphism $\alpha\colon R\otimes M\to M$ such that
  \begin{center}
	\begin{tikzcd}
	  R\otimes R\otimes M \ar[r, "R\otimes\alpha"]\ar[d,"\mu_R\otimes M"] & R\otimes M\ar[d,"\alpha"]\\
	  R\otimes M\ar[r,"\alpha"] & M
	\end{tikzcd}
  \end{center}
  commutes
  and the composite $M\stackrel{\nu_{R}\otimes M}{\longrightarrow}R\otimes M\stackrel{\alpha}{\longrightarrow}M$
  is the identity on $M$;
  here,
  we let $\nu_{R}\colon\Z=1_{\otimes}\to R$ denote $R$'s unit.

  A morphism of condensed $R$-modules is an equivariant morphism of condensed abelian groups.
  That is,
  a morphism $\varphi\colon M\to N$ such that
  \begin{center}
	\begin{tikzcd}
	  R\otimes M\ar[r,"R\otimes \varphi"]\ar[d,"\alpha_M"] & R\otimes N\ar[d,"\alpha_N"]\\
	  M \ar[r,"\varphi"] & N
	\end{tikzcd}
  \end{center}
  commutes.

  We obtain a category $\CondModR$ together with a forgetful functor $?\colon\CondModR\to\cab$.
\end{definition}

Clearly, the following holds, see e.g. \ref{ex:rmod_as_monad}.
\begin{proposition}
    The category of condensed $R$-modules is given by the Eilenberg Moore category
  $\CondModR=\EM(T_{R})$ for the monad
    \[T_{R}=R\otimes-\colon\cab\to\cab,\, T_{R}T_{R}\to T_{R}\]
    given by the multiplication $R\otimes R\to R$ and $\id_{\cab}\to T_{R}$ being the unit of $R$.
\end{proposition}

Unravelling what a morphism $R\otimes M\to M$ is,
we get the following characterisation.
\begin{proposition}
  Let $R\in\CondRing$.
  A condensed $R$-module is a condensed abelian group $M$ together with an $R(S)$-module structure on $M(S)$ for each $S\in\extr$
  (which is just a special homomorphism of abelian groups $R(S)\otimes M(S)\to M(S)$)
  that is compatible in the sense that for every $f\colon S'\to S$,
  \begin{center}
	\begin{tikzcd}
	  R(S)\otimes M(S) \ar[rr]\ar[d,"R(f)\otimes M(f)"] && M(S)\ar[d,"M(f)"]\\
	  R(S')\otimes M(S') \ar[rr] && M(S')
	\end{tikzcd}
  \end{center}
  commutes.
\end{proposition}
\begin{proof}
  Since $R\otimes M=\Sh(?R\otimes_{\PSh}?M)$,
  a map $R\otimes M\to M$ corresponds to a map $?R\otimes_{\PSh}?M\to ?M$.
  Such a map is given by compatible maps $R(S)\otimes M(S)\to M(S)$.
  Using the same translation,
  one sees that a map $R\otimes M\to M$ defining an $R$-module structure on $M$ corresponds to such compatible maps defining $R(S)$-module structures.
\end{proof}

\begin{remark}
    Clearly, by installing the pointwise $\Ab$-tensor product on $\PSh(\extr, \Ab)$, we see that a monoid in
    $\cond$ can be considered as a monoid in $\PSh(\extr,\Ab)$.
    Conversely, using sheafification, we see that every monoid in $\PSh(\extr,\Ab)$ that has underlying object already in $\cab$
    installs a monoid in $\cab$.
    This leads to an appropriate version of \enquote{$R$-valued presheaves} as $R$-modules in $\PSh$, having pointwise values in $\Mod_{R(S)}$ as well as a sheafification.
\end{remark}

Because $\CondModR=\EM(T_{R})$,
$\CondModR\to\cab$ has a left adjoint, the base change $R\otimes-\colon\cab\to\CondModR,$,
we can unpack this somewhat to get an explicit \enquote{formula} and more, see e.g. \ref{rem:restr_of_scalars}.

\begin{corollary}
    \begin{enumerate}[(i)]
        \item  The left adjoint to $\CondModR\to \cab$ can be computed as
  \[
	A\mapsto A\otimes R =\Sh\left(S\mapsto R(S)\otimes A(S)\right)=(S\mapsto R(S)\otimes_{\Z(S)} A(S))
  \]
        \item There exists a monoidal structure on $\CondModR$, which is given by defining $A\otimes_R B$
        as the coequalizer of
\begin{center}\begin{tikzcd}
	{A\otimes R\otimes B} & {A\otimes B}
	\arrow[shift right, from=1-1, to=1-2]
	\arrow[shift left, from=1-1, to=1-2]
\end{tikzcd}\end{center}
        But since the coequalizer may be computed pointwise, and we have an explicit formula for sheafification of modules over discrete rings as $\Z$, this reduces to $(A\otimes_R B)(S)$ being the coequalizer

\begin{center}\begin{tikzcd}
	{A(S)\otimes_{\Z(S)} R(S)\otimes_{\Z(S)} B(S)} & {A(S)\otimes_{\Z(S)} B(S)}
	\arrow[shift right, from=1-1, to=1-2]
	\arrow[shift left, from=1-1, to=1-2]
\end{tikzcd}\end{center}
in algebraic $\Z(S)$-modules.
        But this is precisely the algebraically induced tensor product (using that base change composes, we can interpret this over $\Z$ as well),
        \[(A\otimes_R B)(S)=A(S)\otimes_{R(S)}B(S).\]
        \item There exists an internal $\hom$ in $\CondModR$, being partially right adjoint to $\otimes_{R}$, and it can be computed as the equalizer of

\begin{center}\begin{tikzcd}
	{\ihom(A,B)} & {\ihom(A\otimes R,B )}
	\arrow[shift right, from=1-1, to=1-2]
	\arrow[shift left, from=1-1, to=1-2]
\end{tikzcd}\end{center}
        Again by the explicit formulae, the $S$-valued points are just the algebraically $R(S)$-linear maps from $A(S)$ to $B(S)$.
    \end{enumerate}
    In particular, this implies that base change $A\mapsto R[A]=A\otimes R$ is bicontinuous from $\cab\to\CondModR $.

\end{corollary}
We close with two more wonderful properties of $\CondModR$, showing that in fact the transition from $\Ab$ to $\CondModR$ can be done without any loss.
\begin{corollary}[Of \ref{subsubsec:monads} or the corollary before]
    \begin{enumerate}[(i)]
        \item Limits and colimits may be computed in $\cab$, in particular $\CondModR$ is bicomplete and abelian, and fulfills the same Grothendieck axioms as $\Ab$.
        \item The family of $R[S]=\Z[S]\otimes R$ for extremally disconnected $S$ forms a family of compact projectives, and they generate since the trace map $R[A]\to A$ is epimorphic (for $A\in \cab$).
    \end{enumerate}
\end{corollary}

\subsection{Computations in condensed abelian groups}

For condensed abelian groups, we can carry out some further computations and explicit arguments, which are supposed to help in gaining an intuition for condensed abelian groups.

Firstly, we can explicitly calculate the free objects $\Z[S]$ for profinite sets $S$ as follows.

\begin{lemma}\label{lem:expl-zs}\chapfour
For an $S=\varprojlim S_i$ in $\prof$ with $S_i$ finite we have
\[
    \Z[S] = \bigcup_n\varprojlim_i\Z[S_i]_{\ell^1\le n} \sub \varprojlim\Z[S_i]
\]
as condensed abelian groups.
Here, $\Z[S_i]_{\ell^1\le n}$ are those $\Z$ combinations
indexed by the elements in $S_i$ such that their $\ell^1$-norm is $\le n$.
They in particular are finite, hence $\varprojlim_i\Z[S_i]_{\ell^1\le n}$
is profinite and thus $\Z[S]$ is induced by a
totally disconnected Hausdorff abelian group.
\end{lemma}

\begin{proof}
This is \cite[2.1]{scholze2019Analytic}.

Let $S\in\prof$ be represented as $S=\varprojlim S_i$.
The projections $S\to S_i$ induce maps $\Z[S]\to \Z[S_i]$,
which glue to a map
\[
    \phi\colon \Z[S]\to \varprojlim_i \Z[S_i].
\]
We show that this map $\phi$ is monic.
Consider $K\in\extr$ and $f\in \ker\phi_K$.
By the construction of adjunctions, $\Z[S]$ is the sheafification of the presheaf
$K\mapsto \Z_{\mathrm{alg}}[C(K,S)]$.
An element $f\in \Z[S](K)$ can be written (by the formula for sheafification) as
\[
    f=(f^j)_{j=1,\dots, n}\in \prod \Z_{\mathrm{alg}}[C(T_j,S)]
\]
where $\coprod T_j=K$ is a clopen partition of $K$.
Hence we can represent $f^j$ as formal sum
\[
    f^j=\sum_{i=1}^{k_j} a_{i,j} [f_{i}^j].
\]
with $f_{i}^j\in C(T_j,S)$.

For a partition $M=(N_{k})_k$ of $\{1,\ldots,k_j\}$ define the closed subset
\[
    E_M^j
    =\bigcap_{N\in M}\{t\in T_j\colon f_n^j(t)=f_{n'}^j(t)\text{ for all } n,n'\in N\}
\]
of $T_j$ and, if it exists, choose an element
\[
    t_M^j\in
    \bigcap_{N\in M}\bigcap_{k\in N^c}
    \{t\in T_j\colon f_n^j(t) = f_{n'}^j(t) \ne f_k^j(t)\text{ for all } n,n'\in N\}.
\]
Now take only those $E^j_M$ for which a $t^j_M$ exists.
Together, they are a covering of $T_j$ by closed sets.
Running over all $j$, we obtain a finite cover of $K$
\[
    \coprod_{M\in C} E_M \twoheadrightarrow K
\]
in $\prof$, indexed by a set $C$.
This in turn yields a monomorphism
\[
    \Z[S] \hookrightarrow \prod_{M\in C}\Z[S](E_M).
\]
Hence it suffices to show that the image of $f$ is zero in the codomain of this monic.
This image is given by
\[
    f = \prod_{M\in C}\sum_{i=1}^{k_j}a_{i,j}[f_{i}^j\rvert_{E_M^j}]
\]
But for all $n, n'\in N$ and $N\in M$, $M\in C$,
the functions $f^j_n$ and $f^j_{n'}$ agree.
So we define $f_N = f_n\rvert_{E_M^j}$ for $n\in N$ (which is independent of $n$),
and obtain that
\[
    f = \prod_{M^j\in C}\sum_{N\in M^j}(\sum_{n\in N}a_{i,j})[f_N].
\]
Therefore it is enough to show that for any fixed $M^j\in C$ and $N\in M^j$
all $\sum_{n\in N}a_{i,j}$ are zero.
To achieve this we compute the value at the chosen element $t_M^j$.
Plugging in yields
\[
    \sum_{i}a_{i,j}[f_{i}^j(t^j_M)]=0
\]
because the diagram
\begin{center}
\begin{tikzcd}
	{\Z[S](K)} & {\varprojlim\limits_i\Z[S_i](K)} \\
	{\Z[S](\ast)=\Z_{\mathrm{alg}}[S]} & {\varprojlim\limits_i \Z[S_i](*) = \varprojlim\limits_i \Z_{\mathrm{alg}}[S_i]}
	\arrow["{\phi_K:f\mapsto 0}", from=1-1, to=1-2]
	\arrow["{f\mapsto f(t^j_M)}"', from=1-1, to=2-1]
	\arrow[from=1-2, to=2-2]
	\arrow[hook, from=2-1, to=2-2]
\end{tikzcd}
\end{center}
is commutative.
But the values $f_i^j(t^j_M)=f_{i'}^j(t^j_M)$ agree precisely if $i,{i'}\in N$
for some $N\in M$, so we have
\[
    \sum_{N\in M}(\sum_{i\in N}a_{i,j})[f_N(t^j_M)]=0
\]
By definition of $M$, those values are distinct in $S$,
so this implies
\[
    \sum_{i\in N}a_{i,j}=0
\]
for all $N\in M$ and hence finally that $f=0$ in $\Z[S]$.
We have shown that
\[
    \Z[S] \hookrightarrow \varprojlim_i \Z[S_i]
\]
is a monomorphism.

The addition on $\Z[S_i]_{\le n}$ maps to $\Z[S_i]_{\le 2n}$
and thus
$\varprojlim \Z[S_i]_{\le n}\times\varprojlim \Z[S_i]_{\le n} \to \varprojlim \Z[S_i]_{\le 2n}$.
So the condensed set
\[
    \bigcup_n\varprojlim_i \Z[S_i]_{\le n}
\]
is a condensed abelian group,
The map $S\to\bigcup_n \varprojlim_i \Z[S_i]_{\le n}$
lifts by freeness of $\Z[S]$ to a map of condensed abelian groups
\[
    \Z[S]\to\bigcup_n \varprojlim_i \Z[S_i]_{\le n}.
\]
Restricted to $S$, the embedding
\[
    S\to\bigcup_n\varprojlim_i \Z[S_i]_{\le n}\to \varprojlim_i \Z[S_i]
\]
agrees with directly embedding
\[
    S\to \varprojlim_i \Z[S_i].
\]
By uniqueness of the induced arrow that lifts the arrow from $S$ to $\Z[S]$,
the diagram
\begin{center}
\begin{tikzcd}
	S && {\varprojlim\limits_i \Z[S_i]} \\
	& { \bigcup\limits_n \varprojlim_i \Z[S_i]_{\le n}}
	\arrow[hook, from=1-1, to=1-3]
	\arrow[from=1-1, to=2-2]
	\arrow[from=2-2, to=1-3]
\end{tikzcd}
\end{center}
is commutative.
In particular, there is a monomorphism
$\Z[S]\hookrightarrow  \bigcup_n \varprojlim_i \Z[S_i]_{\le n}$.
Hence it suffices to show that this map is sectionwise epimorphic on $\extr$.
For $n\in \N$ and a partition $S_i$ of $S$, consider the map
\[
    (S\times \{-1,0,1\})^n\to \Z[S_i]_{\le n},\quad
    (x_i,a_i)_{1\le i\le n}\mapsto \sum a_i [x_i].
\]
It is surjective and clearly commutes with the transition maps $\Z[S_i]\to \Z[S_j]$.
Therefore, it induces a morphism
\[
    S^{n}\times \{-1,0,1\}^{n}\to \varprojlim_i \Z[S_i]_{\le n}
\]
which is surjective as cofiltered limit of surjections in $\prof$.
Moreover, $S^{n}\times \{-1,0,1\}^{n}$ can also be embedded into $\Z[S]$.
All of those morphism fit into the commutative diagram
\begin{center}\begin{tikzcd}
	{S^n \times \{-1,0,1\}^n} & {\varprojlim\limits_i \Z[S_i]_{\le n}} \\
	{\Z[S]} & {\bigcup\limits_m\varprojlim\limits_i\Z[S_i]_{\le m} }
	\arrow[two heads, from=1-1, to=1-2]
	\arrow[from=1-1, to=2-1]
	\arrow[hook, from=1-2, to=2-2]
	\arrow[hook, from=2-1, to=2-2]
\end{tikzcd}.\end{center}
We conclude that the image contains $\varprojlim\Z[S_i]_{\le_n}$ for every $n\in\N$,
and thus all of the union
\[
    \bigcup_m\varprojlim\Z[S_i]_{\le m}.
\]
\end{proof}

\begin{lemma}\chapfour\label{comp_proj_gen}
The category $\cab$ is generated by compact projective objects, and a generating family of such is $\Z[T]$ for extremally disconnected $T$.
In fact, $\hom_{\cab}(\Z[T],-)$ commutes with all limits and colimits.
\end{lemma}
\begin{proof}
  We already know this, see~\ref{prop:r-con-enogh-cproj}.

  We give a quick explicit argument.
  For $M\in\cab$,
  take a surjection $\bigsqcup S_{i}\twoheadrightarrow M$,
  apply the left adjoint $\Z$ to obtain $\oplus_{i}\Z[S_{i}]=\Z[\sqcup S_{i}]\twoheadrightarrow\Z[M]$
  and compose with the trace map $\Z[M]\twoheadrightarrow M$ (\ref{prop:trace-pie}).
  It is, of course, also possible to simply exploit the universal property of $\Z[\sqcup S_{i}]$.

  Alternatively,
  we give a direct argument from~\cite[2.2]{scholze2019condensed} showing that the $\Z[T]$ are generating.
	For this consider any $M\in\cond$ and $\kappa$ such that $M=\mathrm{Lan}_\kappa M|_\kappa\in \cond_\kappa$.
	Now consider the set of  $\cond_\kappa$ subobjects $N\hookrightarrow M$ such that there exists a projection $\bigoplus \Z[T_i]\twoheadrightarrow N$.
	They clearly form a partially ordered set, and for any totally ordered chain $M_i$ of subobjects with epimorphism $\bigoplus_{j\in J_i} \Z[T_j]\twoheadrightarrow M_i\hookrightarrow M$ their union $\bigoplus_{j\in \bigcup_i J_i}\Z[T_j]\to M$ defines a morphism with image $K\hookrightarrow M$ containing every $M_i$, as by the universal property of the image there exists a unique monic
\begin{center}\begin{tikzcd}
	{\bigoplus_{j\in J_i}\Z[T_j]} & {\bigoplus_{j\in \bigcup J_i}\Z[T_j]} \\
	{M_i} & K \\
	& M
	\arrow[hook, from=1-1, to=1-2]
	\arrow[two heads, from=1-1, to=2-1]
	\arrow[two heads, from=1-2, to=2-2]
	\arrow["{\exists!}", dashed, hook, from=2-1, to=2-2]
	\arrow[hook, from=2-1, to=3-2]
	\arrow[hook, from=2-2, to=3-2]
\end{tikzcd}\end{center}
Hence, every chain of such subobjects admits an upper bound.
Thus, by Zorn's lemma there exists a maximal subobject $N$ with epimorphism $\bigoplus\Z[T_i]\twoheadrightarrow N$.
	But if $N\ne M$, then $M/N\ne 0$,  hence $M/N(S)\ne 0$ for some $S\in \extr$.
	But this means there is a nontrivial morphism $S\to M/N$, i.e. a nontrivial homomorphism of condensed abelian groups $\Z[S]\to M/N$.
	Lifting this by projectivity of $\Z[S]$ to $M$ there is a nonzero homomorphism
\begin{center}\begin{tikzcd}
	& M \\
	{\Z[S]} & {M/N}
	\arrow[two heads, from=1-2, to=2-2]
	\arrow[from=2-1, to=1-2]
	\arrow[from=2-1, to=2-2]
\end{tikzcd}.\end{center}
We claim that the induced morphism
\[\Z[S]\oplus \bigoplus \Z[T_i]\to M\]
has image strictly larger than $N$.
	Clearly it contains $N$, and if the image would be precisely $N$, then the composition with $M\to M/N$ would be $0$.
	Hence $N$ is not maximal, and therefore the maximal subobject with surjection is $M|_\kappa$ itself.
	As the left adjoint (given by left Kan extension) preserves colimits, this induces an epimorphism to $M$.
\end{proof}
\begin{lemma}
    The compact projective objects in $\cab$ are precisely the retracts of $\Z[\beta I]$ for some set $I$.
	For such a retract, the remarkable property $P\oplus \Z[\beta I]\simeq \Z[\beta I]$ holds.
\end{lemma}
\begin{proof}
    The first statement follows as there are enough compact projectives, and $\Z[\beta I]$ is a generating class of such.
    See 3.5 in \cite{scholze2022complex} for the second statement.
    \end{proof}

It is an open question whether all compact projective objects are of the form $\Z[S]$, see \cite[3.6]{scholze2022complex}.
\begin{lemma}\label{lem:no_inj_cab}
  There are no nonzero injective condensed abelian groups.
\end{lemma}
\begin{proof}
  We give a proof by Scholze given in a mathoverflow answer.%
  \footnote{\url{https://mathoverflow.net/questions/352448/are-there-enough-injectives-in-condensed-abelian-groups}}

	Take any injective $I$.
    We first show that there is a compact Hausdorff $K$ with $K\twoheadrightarrow I$.
	For this take some surjection $\sqcup_{i\in J}S_j\twoheadrightarrow I $,
    lift this by the universal property of $\Z[\sqcup S_{i}]=\oplus_{i}\Z[S_{i}]$ to a map
	\[\bigoplus_{j\in J}\Z[S_j]\twoheadrightarrow I,\]
	where $\bigoplus_{j\in J}\Z[S_j]\simeq \Z[\sqcup_{j\in J}S_j]$ follows as $\Z[-]$ is a left adjoint.
	By using the embedding $\Z\hookrightarrow b\Z$ of $\Z$ into the \idx{Bohr-compactification} of $\Z$
    (or alternatively any other compact Hausdorff topological abelian group such as the $p$-adic integers),
	we obtain a monomorphism
	\[\bigoplus_{i\in J}\Z[S_j]\hookrightarrow \prod_{j\in J} b\Z[S_j] \eqqcolon K.\]
    (Monomorphicity can be tested in $\Set$ or $\Ab$.)
	As $I$ is injective, this lifts to a morphism
\begin{center}\begin{tikzcd}
	& K \\
	I & {\bigoplus \Z[S_j]}
	\arrow[from=1-2, to=2-1]
	\arrow[hook', from=2-2, to=1-2]
	\arrow[two heads, from=2-2, to=2-1]
\end{tikzcd}\end{center}
which automatically has to be epimorphic.
Hence $I$ is quasicompact with epimorphism $K\twoheadrightarrow I$.
Now consider
	\[\bigoplus_{2^J}I\hookrightarrow \prod_{2^J}I,\]
	and extend the sum map $\bigoplus_{2^J} I\to I$
\begin{center}\begin{tikzcd}
	& {\prod_{2^J} I} \\
	I & {\bigoplus_{2^J} I}
	\arrow[from=1-2, to=2-1]
	\arrow[hook', from=2-2, to=1-2]
	\arrow[from=2-2, to=2-1]
\end{tikzcd}\end{center}
Next, we want to see that this factors through a subdiagram of cardinality $\le J$.
\begin{center}\begin{tikzcd}
	{\prod_A I} & {\prod_{2^J} I} \\
	I & {\bigoplus_{2^J} I}
	\arrow[from=1-1, to=2-1]
	\arrow[from=1-2, to=1-1]
	\arrow[hook', from=2-2, to=1-2]
	\arrow[from=2-2, to=2-1]
\end{tikzcd}\end{center}
For this,
it suffices to show that $\prod_{2^{J}}K\to\prod_{2^{J}}I\to I$
factors through a sufficiently small subdiagram of $\prod_{2^{J}}K$,
because then by definition of the arrow $\prod_{2^{J}}K\to I$,
this induces a factoring $\prod_{2^{J}}K\to\prod_{A}I\to I$,
but $\prod_{2^{J}}K\to\prod_{2^{J}}I$ is epic.

To see that $\prod K\to I$ factors through a sufficiently small subdiagram,
consider the Kan extension formula,
giving us $|J|$-profinite $S$ together with a factorisation $\prod K\to S\to I$.
But the arrow $\prod K\to S$ already factors through a sufficiently small subdiagram.
(There are many topological ways to see this,
e.g., by noting that the cardinal bound on the weight on $S$
even gives a good cardinal bound on the number of points of $S$.
However, the probably simplest way is to use $\kappa$-cocompactness
of $S$ by lemma~\ref{lem:kappa_prof_is_kappa_cocomp} together with
remark~\ref{rem:kappa_commutes_into_subdiag}.)

Now, take any $i\in 2^J\setminus A$.
Then the addition map $I\to I$ is the identity, but the composition through the product is $0$,
\begin{center}\begin{tikzcd}
	{\prod_A I} & {\prod_{2^J} I} \\
	I & I
	\arrow[from=1-1, to=2-1]
	\arrow["{\pi_A}"', from=1-2, to=1-1]
	\arrow["0"', from=2-2, to=1-1]
	\arrow["{\tau_i}"', hook', from=2-2, to=1-2]
	\arrow["\id", from=2-2, to=2-1]
\end{tikzcd}.\end{center}
Thus $I=0$.
\end{proof}

\begin{lemma}[Breen-Deligne resolution]\chapfour
There exists a functorial resolution of every condensed abelian group $A$ of the form
\[\bigoplus_{\mathrm{fin}}\Z[A^{a_i}]\to \dots \to \Z[A^2]\to \Z[A]\to A\to 0.\]
	The last map is induced by the identity on $A$ and applying the adjunction.
	The map $\Z[A\times A]\to \Z[A]$ is induced by the addition
	\[[(a,b)]\mapsto [a+b]-[a]-[b],\]
	which is defined by using the pullback of the projections $A^2\to A$ and the addition map $A\times A\to A$.
\end{lemma}

\begin{proof}
This is \cite[4.10]{scholze2019condensed}.
\end{proof}

\begin{lemma}[Quasiseparated condensed abelian groups]\label{lem:qs_tensor_hom}\chapfour
For any quasiseparated condensed set $X$, the free abelian group $\Z[X]$ is quasiseparated (as a condensed set).

For any condensed abelian group $A$, the quasiseparation of the underlying condensed set $A$ yields a quasiseparated condensed abelian group, inducing a left adjoint to $\Ab(\qs)\sub \cab$.

For any $A\in \cab$ and $B\in \Ab(\qs)$, we have that $\ihom(A,B)$ is quasiseparated.
In particular, the category $\Ab(\qs)$ of quasiseparated condensed abelian groups closed under internal homs,
and taking tensor products followed by quasiseparation induces a unique closed symmetric monoidal
structure on $\Ab(\qs)$ such that quasiseparation is symmetric monoidal.
\end{lemma}
\begin{proof}
	First, we remark that for any compact Hausdorff space $K$, the free abelian group $\Z[K]$ can be written as the countable union $\bigcup \Z[K]_{\le n}$,
	where $\Z[K]_{\le n}$ are compact Hausdorff (see also Ex.\ 2.3 in \cite{scholze2019Analytic}).
	In particular, it is quasiseparated.

	Now take any quasiseparated condensed set $X$ and write it as the filtered colimit along inclusions of its qcqs subobjects $K_i$.
	Since $\Z[-]$ commutes with colimits, clearly $\Z[X]=\varinjlim \Z[K_i]$, and since monics $K_i\hookrightarrow K_j$ are sectionwise monic,
	they induce monics $\Z_{\mathrm{alg}}[K_i(S)]\hookrightarrow \Z_{\mathrm{alg}}[K_j(S)]$.
	Since sheafification is exact, this already implies that all the $\Z[K_i]\to \Z[K_j]$ are monic.
	Hence $\Z[X]$ is the filtered colimit (these may be computed in $\cond$) along inclusions of the quasiseparated sets $\Z[K_i]$, and hence $\Z[X]$ is quasiseparated.

	Since quasiseparation commutes with finite products, the quasiseparation of any abelian group object yields an condensed abelian group
    (clearly, this works with condensed $\mcT$-objects for general Lawvere theories),
	being the left adjoint of $\Ab(\qs)\sub \cab$.

	As $\ihom(A,B)$ is contravariant and cocontinuous in the first coordinate and since limits of quasiseparated sets are quasiseparated,
	it suffices to show the assertion for $A=\Z[S]$ with $S\in \extr$.
	But here
	\[\ihom(\Z[S], B)(T)=\hom_{\cab}(\Z[S\times T],B)=\hom_{\cond}(S\times T, B).\]
	Since $S\times T$ is compact, for any quasiseparated $B=\varinjlim K_i$ we obtain
	\[\hom_{\cond}(S\times T, B)=\varinjlim \hom_{\cond}(S\times T, K_i),\]
	and as $K_i\to K_j$ is monic, clearly $\hom_{\cond}(S\times T,K_i)\hookrightarrow \hom_{\cond}(S\times T, K_j)$.
	Hence $\ihom(\Z[S],B)$ is the filtered colimit along inclusions of $\hom_{\cond}(S\times -,K_i)=\ihom_{\cond}(S, K_i)$.
	But the latter is equipped with compact open topology and thereby quasiseparated.
\end{proof}
\begin{remark}
		Note that a similar statement can be made about the internal hom in condensed sets by the exact same reasoning.

		Furthermore remark that we do not know whether a $\cab$ tensor product of quasiseparated condensed sets is automatically quasiseparated.
        We wonder how smoothly this generalises, for example,
        to condensed modules over arbitrary condensed rings.
	\end{remark}
\begin{remark}
	There are canonical alternative definitions of (quasi)compact abelian groups and quasiseparated abelian groups, by noting that the definition of quasiseparated and quasicompact objects work in any category.
	However, as we want condensed sets to be the replacement of topological spaces and want to allow an interplay between \enquote{compact subsets} of \enquote{topological abelian groups},
	we stick to the term quasicompact and quasiseparated to mean the corresponding notions in condensed sets.
	\end{remark}

Next, we see that the internal hom of good topological abelian groups corresponds with the internal hom in $\cab$.
\begin{lemma}\chapfour
Let $A, B$ be Hausdorff topological abelian groups and $A$ be compactly generated.
The internal hom $\ihom(\underline{A},\underline{B})$ corresponds to the topological one with the compact-open topology, i.e.,
\[\ihom(\underline{A},\underline{B})\simeq \underline{\hom(A,B)},\]
where $\hom(A,B)$ is equipped with compact open topology.
\end{lemma}
\begin{proof}
We want to show that \[\ihom(\underline{A},\underline{B})\simeq \underline{\hom(A,B)}.\]
By Yoneda, for any $S\in \extr$, the elements of $\ihom(\underline{A},\underline{B})(S)$ are given by
maps of condensed sets $S\to\ihom(A,B)$.
 \[\ihom(A,B)(S)=\hom_{\cond}(S,\ihom(A,B))=\hom_{\cab}(\Z[S],\ihom(A,B))\simeq \hom_{\cab}(A\otimes \Z[S],B)\]

 Next, use the epimorphism $\Z[A]\twoheadrightarrow A$, and tensor with $\Z[S]$ to obtain an epimorphism
 (as tensoring is a left adjoint, it preserves epimorphisms) $\Z[A\times S]\twoheadrightarrow A\otimes \Z[S]$.
 Pullback by this morphism gives a map
\[\ihom(A,B)(S)=\hom(A\otimes\Z[S], B)\to \hom_{\cab}(\Z[A\times S],B)\simeq \hom_{\cond}(A\times S, B).\]
This shows that any map in $\hom(A\otimes \Z[S],B)$ corresponds to a map $A\times S\to B$, which can be seen as a map of topological spaces,
as $A\times S$ is compactly generated and thus $A\times S=(\underline{A\times S})(\ast)_{\mathrm{top}}$.

Clearly, as pullback by an epimorphism is a monic procedure (by definition of epimorphicity), this map is a monomorphism.
To see that it is an isomorphism, consider any $f\colon \Z[A]\otimes\Z[S]\to B$.
This is equivalently a map $A\times S\to B$.
Using pullback by $A\times \{s\}\to A\times S$, we see that for every $s\in S$, the map
$f_s\colon A\to B$ is a group homomorphism.

We need to show that $f$ factors through $A\otimes \Z[S]$, i.e., we need to show that $f$ vanishes on the kernel of $p\colon \Z[A]\otimes \Z[S]\to A\otimes \Z[S]$.
But by considering the Breen-Deligne resolution of $A$ and tensoring with the projective object $\Z[S]$ we see that the kernel can be described as

\begin{center}\begin{tikzcd}
	\dots & {\Z[A\times A]\otimes \Z[S]} && {\Z[A]\otimes \Z[S]} & {A\otimes\Z[S]} & 0
	\arrow[from=1-1, to=1-2]
	\arrow["{([a+b]-[a]-[b])\otimes 1}", from=1-2, to=1-4]
	\arrow["p", from=1-4, to=1-5]
	\arrow[from=1-5, to=1-6]
\end{tikzcd}\end{center}
Hence we need to show that the composition

\begin{center}\begin{tikzcd}
	{\Z[A\times A \times S]} &&& {\Z[A\times S]} & {A\otimes\Z[S]}
	\arrow["{[(a+b,s)]-[(a,s)]-[(b,s)]}", from=1-1, to=1-4]
	\arrow["f", from=1-4, to=1-5]
\end{tikzcd}\end{center}
vanishes.
This is equivalently a map $A\times A\times S\to B$, described by
\[(a,b,s)\mapsto f_s(a+b)-f_s(a)-f_s(b)=0.\]
Hence $f$ factors through $A\otimes \Z[S]$, and we obtain isomorphisms (natural in $S$)
\[\ihom(A,B)(S)\simeq \hom_{\cond}(A\times S,B).\]

Lastly, note that continuous maps from $S$ to $\underline{\hom(A,B)}$ with compact open topology are equivalently
continuous maps $A\times S\to B$, by the cartesian closedness of locally compact abelian groups.
\[\hom_{\cond}(A\times S,B)\simeq \hom_{\cond}(S,\underline{\hom(A,B)})=\underline{\hom(A,B)}(S).\]

This in total shows that there is a natural isomorphism
\[\ihom(A,B)\to\underline{\hom(A,B)}.\]

\end{proof}
\begin{lemma}
    A condensed abelian group $A$ is flat precisely if $A(S)$ is torsion-free for all $S\in \extr$ (equivalently $S\in\prof$).
    \end{lemma}
\begin{proof} This is \cite[3.4]{scholze2022complex}.

The flat abelian groups are precisely the torsion free abelian groups.
Hence, clearly, if $A(S)$ is torsion-free for all $S\in \extr$, then tensoring on level of presheaves with $A$ is exact.
As sheafification is exact, this implies that tensoring on the level of sheaves remains exact.

Conversely, consider any $S\in \extr$ and any flat condensed abelian group $A$.
It suffices to see that for $i\colon n\Z\hookrightarrow \Z$ the arrow
\[A(S)\otimes n\Z\to A(S)\otimes \Z\]
is monic.

We know that by flatness of $A$ (and computing everything pointwise) and the embedding $n\Z\hookrightarrow \Z$ as discrete abelian groups being monic, the arrow
\[A(S)\otimes C(S,n\Z)\to A(S)\otimes C(S,\Z), \quad a\otimes f\mapsto a\otimes i f\]
is monic as a map of abelian groups.
But using constant functions $k$, this implies that the constant function
\[a\otimes k\in A(S)\otimes C(S,\Z)\]
is nonzero, which implies that the value
\[a\otimes k\in A(S)\otimes \Z\]
is nonzero.

This is essentially the factorisation
\begin{center}\begin{tikzcd}
	{A(S)\otimes n\Z} & {A(S)\otimes C(S,n\Z)} \\
	{A(S)\otimes \Z} & {A(S)\otimes C(S,\Z)}
	\arrow[hook, from=1-1, to=1-2]
	\arrow[hook', from=1-1, to=2-1]
	\arrow[hook, from=1-2, to=2-2]
	\arrow[hook, from=2-1, to=2-2]
\end{tikzcd}\end{center}
and hence $A(S)$ is indeed torsion-free.

Alternatively, we could have used that $C(S,\Z)$ and thus also $C(S,n\Z)$ are free abelian groups (with obvious identification of summands), and decomposed

\begin{center}\begin{tikzcd}
	{A(S)\otimes C(S,n\Z)} & {A(S)\otimes C(S,\Z)} \\
	{\bigoplus (A(S)\otimes n\Z)} & {\bigoplus (A(S)\otimes \Z)}
	\arrow[hook, from=1-1, to=1-2]
	\arrow["\simeq"', from=1-1, to=2-1]
	\arrow["\simeq"', from=1-2, to=2-2]
	\arrow["{\bigoplus 1\otimes i}"', from=2-1, to=2-2]
\end{tikzcd}\end{center}
Hence $A(S)$ is torsion-free for all $S\in \extr$.

Clearly, as subgroups of torsion-free abelian groups are torsion free,
for any profinite set $S$ the map $A(S)\hookrightarrow A(\beta S_d)$ shows that the values on all profinite sets are torsion-free.
\end{proof}

\begin{remark}Note that if every sections functor $X\mapsto X(S)$ were essentially surjective, the above argument could have been simplified.
	But unfortunately, e.g., $X(2)=X(\ast)\times X(\ast)$ must be decomposable and hence cannot be, e.g., $\Z$.
	\end{remark}

We close with a quite canonical question we are not quite sure about.
\begin{question}[Adjoint functor theorem for $\cab$]\chapfour
    Is it true that a functor from $\cab$ to some category has a left adjoint if and only if it is continuous?
\end{question}

\subsubsection{Derived condensed abelian groups}

As we have seen,
$\cab$ is an abelian category of the best sort (closed symmetric monoidal structure, generated by compact projectives, etc.).
Therefore (especially because of it being generated by compact projectives),
we can employ the whole powerful machinery of homological algebra.

\begin{lemma}\chapfour
Consider the derived category of condensed abelian groups, $D(\cab)$.
It is a triangulated additive category with symmetric monoidal tensor product, internal $\hom$ and is generated (as triangulated category)
by compact projectives $P[0]$ for compact projectives $P\in \cab$.
\end{lemma}
\begin{proof}
The derived category of any abelian category with the corresponding properties fulfills these properties.
\end{proof}

We can explicitly calculate internal $\RHom$'s via the following spectral sequence,
which can be found in \cite[4.8]{scholze2019condensed} or \cite[4.1]{Aparicio2021}.
\begin{proposition}\label{prop:spec_seq_condab}
    Let $A$, $B$ be condensed abelian groups.
    Then there is a spectral sequence with
    \[E_1^{i,j}=\prod_{k=1}^{n_i}H^j(A^{r_{i,j}}\times S,M)\]
    converging to $\underline{\mathrm{Ext}}^{i+j}(A,M)(S)$,
    where $H^{j}(-,-)=H^{j}(\RHom(\Z[-],-))$.
    This is functorial in $A,M $ and $S$.
\end{proposition}
\begin{proof}
    As $\ihom(A,B)(S)=\hom(A\otimes \Z[S],B)$, we can derive both sides and obtain
    \[\iRHom(-,-)(S)=R(\Hom(-\otimes \Z[S],-)).\]
    Using that derived $\hom$ as a derived bifunctor may be computed by fixing one variable,
    or by an argument such as using exactness of $-\otimes \Z[S]$ and the Grothendieck spectral sequence, we can see that this in fact coincides with
    \[\iRHom(-,-)(S)=\RHom(-\otimes \Z[S],-).\]
    This reduces the problem to calculating the $\RHom$.

    Consider a Breen-Deligne resolution of $A$,

\begin{center}\begin{tikzcd}
	\dots & {\bigoplus_{k=1}^{n_2}\Z[A^{r_{2,k}}]} & {\bigoplus_{k=1}^{n_1}\Z[A^{r_{1,k}}]} & {\bigoplus_{k=1}^{n_0}\Z[A^{r_{0,k}}]} & 0 \\
	& \dots & 0 & A & 0
	\arrow[from=1-1, to=1-2]
	\arrow[from=1-2, to=1-3]
	\arrow[from=1-3, to=1-4]
	\arrow[from=1-3, to=2-3]
	\arrow[from=1-4, to=1-5]
	\arrow[from=1-4, to=2-4]
	\arrow[from=1-5, to=2-5]
	\arrow[from=2-2, to=2-3]
	\arrow[from=2-3, to=2-4]
	\arrow[from=2-4, to=2-5]
\end{tikzcd}\end{center}
    Now, as $\Z[S]$ is flat, i.e., tensoring with $\Z[S]$ is exact, we obtain the quasiisomorphism

\begin{center}\begin{tikzcd}
	\dots & {\bigoplus_{k=1}^{n_2}\Z[A^{r_{2,k}}]\otimes \Z[S]} & {\bigoplus_{k=1}^{n_1}\Z[A^{r_{1,k}}]\otimes \Z[S]} & {\bigoplus_{k=1}^{n_0}\Z[A^{r_{0,k}}]\otimes \Z[S]} & 0 \\
	& \dots & 0 & {A\otimes \Z[S]} & 0
	\arrow[from=1-1, to=1-2]
	\arrow[from=1-2, to=1-3]
	\arrow[from=1-3, to=1-4]
	\arrow[from=1-3, to=2-3]
	\arrow[from=1-4, to=1-5]
	\arrow[from=1-4, to=2-4]
	\arrow[from=1-5, to=2-5]
	\arrow[from=2-2, to=2-3]
	\arrow[from=2-3, to=2-4]
	\arrow[from=2-4, to=2-5]
\end{tikzcd}\end{center}
    Now, as $\Z[-]$ is symmetric monoidal, the above row is isomorphic to

\begin{center}\begin{tikzcd}
	\dots & {\bigoplus_{k=1}^{n_2}\Z[A^{r_{2,k}}\times S]} & {\bigoplus_{k=1}^{n_1}\Z[A^{r_{1,k}}\times S]} & {\bigoplus_{k=1}^{n_0}\Z[A^{r_{0,k}}\times S]} & 0
	\arrow[from=1-1, to=1-2]
	\arrow[from=1-2, to=1-3]
	\arrow[from=1-3, to=1-4]
	\arrow[from=1-4, to=1-5]
\end{tikzcd}\end{center}

    Hence, $\Rhom(A,B)$ can be computed by applying the $\RHom(-,B)$ to this chain complex.
   If $A$ is extremally disconnected, this may be done pointwise, as in this case $\Z[A]$ is projective.
    To compute this in general, consider any Cartan Eilenberg resolution $P_{\bullet,\bullet}$

\begin{center}\begin{tikzcd}
	&& \vdots & \vdots \\
	& \dots & {P_{1,1}} & {P_{0,1}} & 0 \\
	\dots & {P_{2,0}} & {P_{1,0}} & {P_{0,0}} & 0 \\
	\dots & {\bigoplus_{k=1}^{n_2}\Z[A^{r_{2,k}}\times S]} & {\bigoplus_{k=1}^{n_1}\Z[A^{r_{1,k}}\times S]} & {\bigoplus_{k=1}^{n_0}\Z[A^{r_{0,k}}\times S]} & 0
	\arrow[from=1-3, to=2-3]
	\arrow[from=1-4, to=2-4]
	\arrow[from=2-2, to=2-3]
	\arrow[from=2-2, to=3-2]
	\arrow[from=2-3, to=2-4]
	\arrow[from=2-3, to=3-3]
	\arrow[from=2-4, to=2-5]
	\arrow[from=2-4, to=3-4]
	\arrow[from=2-5, to=3-5]
	\arrow[from=3-1, to=3-2]
	\arrow[from=3-2, to=3-3]
	\arrow[from=3-2, to=4-2]
	\arrow[from=3-3, to=3-4]
	\arrow[from=3-3, to=4-3]
	\arrow[from=3-4, to=3-5]
	\arrow[from=3-4, to=4-4]
	\arrow[from=3-5, to=4-5]
	\arrow[from=4-1, to=4-2]
	\arrow[from=4-2, to=4-3]
	\arrow[from=4-3, to=4-4]
	\arrow[from=4-4, to=4-5]
\end{tikzcd}\end{center}

    We apply $\hom(-,B)$ pointwise to obtain a double complex whose spectral sequences converge to $\RHom(A\otimes\Z[S], B)$.
    The resulting double complex is

\begin{adjustbox}{max width = \textwidth}
  \begin{tikzcd}
	&& \vdots & \vdots \\
	& \dots & {\hom(P_{1,1},B)} & {\hom(P_{0,1},B)} & 0 \\
	\dots & {\hom(P_{2,0},B)} & {\hom(P_{1,0},B)} & {\hom(P_{0,0},B)} & 0 \\
	\dots & 0 & 0 & 0 & 0
	\arrow[from=1-3, to=2-3]
	\arrow[from=1-4, to=2-4]
	\arrow[from=2-2, to=2-3]
	\arrow[from=2-2, to=3-2]
	\arrow[from=2-3, to=2-4]
	\arrow[from=2-3, to=3-3]
	\arrow[from=2-4, to=2-5]
	\arrow[from=2-4, to=3-4]
	\arrow[from=2-5, to=3-5]
	\arrow[from=3-1, to=3-2]
	\arrow[from=3-2, to=3-3]
	\arrow[from=3-2, to=4-2]
	\arrow[from=3-3, to=3-4]
	\arrow[from=3-3, to=4-3]
	\arrow[from=3-4, to=3-5]
	\arrow[from=3-4, to=4-4]
	\arrow[from=3-5, to=4-5]
	\arrow[from=4-2, to=4-3]
	\arrow[from=4-3, to=4-4]
	\arrow[from=4-4, to=4-5]
\end{tikzcd}\end{adjustbox}
This allows us to compute the internal $\RHom$ in many cases, see also \cite{Aparicio2021} or \cite[chap.~4]{scholze2019condensed} for more details.
As every column yields a projective resolution $P_{i,\bullet}$ of ${\bigoplus_{k=1}^{n_i}\Z[A^{r_{i,k}}\times S]}$,
the homology of the columns precisely compute the homology
\[H_j\left(\RHom({\bigoplus_{k=1}^{n_i}\Z[A^{r_{i,k}}\times S]},B)\right),\]
which is the definition of
\[H_j\left({\bigoplus_{k=1}^{n_i}\Z[A^{r_{i,k}}\times S]}, B\right).\]
Thus, consider the downwards orientated spectral sequence.
$E^1_{\bullet,\bullet}$ looks like

\begin{adjustbox}{max width=\textwidth}
  \begin{tikzcd}
	&& \dots & \dots \\
	& \dots & {H_1({\bigoplus_{k=1}^{n_1}\Z[A^{r_{1,k}}\times S]},B)} & {H_1({\bigoplus_{k=1}^{n_0}\Z[A^{r_{0,k}}\times S]},B)} & 0 \\
	\dots & {H_0({\bigoplus_{k=1}^{n_2}\Z[A^{r_{2,k}}\times S]},B)} & {H_0({\bigoplus_{k=1}^{n_1}\Z[A^{r_{1,k}}\times S]},B)} & {H_0({\bigoplus_{k=1}^{n_0}\Z[A^{r_{0,k}}\times S]},B)} & 0 \\
	\dots & 0 & 0 & 0 & 0
	\arrow[from=2-2, to=2-3]
	\arrow[from=2-3, to=2-4]
	\arrow[from=2-4, to=2-5]
	\arrow[from=3-1, to=3-2]
	\arrow[from=3-2, to=3-3]
	\arrow[from=3-3, to=3-4]
	\arrow[from=3-4, to=3-5]
	\arrow[from=4-2, to=4-3]
	\arrow[from=4-3, to=4-4]
	\arrow[from=4-4, to=4-5]
\end{tikzcd}\end{adjustbox}
Now, as $\RHom$ (and taking limits and colimits) is additive, $H_i(-,B)$ commutes with finite direct sums, and
we obtain
\[E_{i,j}^1=H_i({\bigoplus_{k=1}^{n_j}\Z[A^{r_{j,k}}\times S]},B)=\bigoplus_{k=1}^{n_j}H_i(\Z[A^{r_{j,k}}\times S],B)\]
    \end{proof}

For locally compact abelian groups, this is particularly easy, as we have the following structural result for LCA groups.
\begin{lemma}[Decomposition of LCA-groups]
    Every LCA group $A$ can be decomposed as
    \[A=\R^{n_A}\oplus T_A,\]
    where the torsion part $T_A$ admits a compact open subgroup.
    This can be rephrased as $T_A$ being an extension of a discrete group $D_A$ by a compact abelian group $K_A$,
     i.e., there is a short exact sequence (in $\cab$!)

\begin{center}\begin{tikzcd}
	0 & K_A & T_A & D_A & 0
	\arrow[from=1-1, to=1-2]
	\arrow[from=1-2, to=1-3]
	\arrow[from=1-3, to=1-4]
	\arrow[from=1-4, to=1-5]
\end{tikzcd}\end{center}

	The discrete abelian group $D_A$ fits into an exact sequence

\begin{center}\begin{tikzcd}
	0 & {\Z[Q_A]=\bigoplus_{Q_A}\Z} & {\Z[R_A]=\bigoplus_{R_A} \Z} & D_A & 0
	\arrow[from=1-1, to=1-2]
	\arrow[from=1-2, to=1-3]
	\arrow["{\psi_A}", from=1-3, to=1-4]
	\arrow[from=1-4, to=1-5]
\end{tikzcd},\end{center}

	Pontryagin dually, the compact abelian group $K_A$ fits into an exact sequence

\begin{center}\begin{tikzcd}
	0 & K_A & {\prod_I \R/\Z} & {\prod_J\R/\Z} & 0
	\arrow[from=1-1, to=1-2]
	\arrow[from=1-2, to=1-3]
	\arrow["\phi_A", from=1-3, to=1-4]
	\arrow[from=1-4, to=1-5]
\end{tikzcd}\end{center}

\end{lemma}
\begin{proof}
        The fact that locally compact abelian groups can be decomposed in this way in the category of locally compact abelian groups is, e.g., \cite[2.2]{Hoffmann2007}.
    We need to see that the exact sequence in locally compact abelian groups
    \begin{center}\begin{tikzcd}
	0 & K_A & T_A & D_A & 0
	\arrow[from=1-1, to=1-2]
	\arrow[from=1-2, to=1-3]
	\arrow[from=1-3, to=1-4]
	\arrow[from=1-4, to=1-5]
\end{tikzcd}\end{center}
    remains exact in condensed abelian groups.
    Since the embedding $\CHaus\to \cond$ is left exact (and thus, so is $\Ab(\CHaus)\to\cab$), it remains to show surjectivity of $T_A\to D_A$.
    But this follows since $(T_A)_d \to (D_A)_d$ is epimorphic ($\Set\to\cond$ admits a right adjoint)
    and the corresponding map $(T_{A})_{d}\to D_{A}$ may be factorized as $ (T_A)_d\to T_A\to D_A$.
    The other cases can be handled similarly or, e.g., using the criterion 2.18 in \cite{Aparicio2021}.
\end{proof}

Now, we calculate explicitly some $\RHom$'s; essentially this is \cite[4.3]{scholze2019condensed}.
\begin{lemma}
	\begin{enumerate}[(i)]
		\item $\iRHom(\Z,B)=B[0]$ for all $B$,
		\item $\iRHom(\R, D)=0$ for $D$ discrete,
		\item $\iRHom(\prod_I\T,D)=\bigoplus_I D[-1]$ for $D$ discrete,
		\item $\iRHom(\prod_I \T,\R)=0$,
		\item $\iRHom(\prod_I \T, \T)=\bigoplus_I \Z[0]$,
		\item $\iRHom(\R,\R)=\R[0]$, and
		\item $\iRHom(\R, \T)=\R[0]$.
\end{enumerate}	\end{lemma}
\begin{proof}
\begin{enumerate}[(i)]

\item	To see the first statement, note that $\ihom(\Z,B)=B$ since $\Z[\ast]=\Z$.
	This implies
	\begin{align*}\ihom(\Z,B)(S)&=\hom_{\cond}(S,\ihom(\Z,B))=\hom_{\cab}(\Z[S],\ihom(\Z,B))
	\\&=\hom(\Z[S]\otimes \Z[\ast],B)=\hom_{\cab}(S\times \ast,B)=B(S).\end{align*}

\item We only sketch the arguments for (ii), see \cite[4.3]{scholze2019condensed} for more details.
Using the spectral sequence from \ref{prop:spec_seq_condab},
	\[E_1^{i,j}=\prod_{k=1}^{n_i}H^j(\R^{r_{i,j}}\times S,D)\Rightarrow \iRHom(\R,D)\]

	We show that $0\to \R$ induces an isomorphism $H^j(\R^r\times S,D)\to H^j(\ast\times S, D)$; this implies
	\[\prod_{k=1}^{n_i}H^j(\R^{r_{i,j}}\times S,D)=\prod_{k=1}^{n_i}H^j(\ast\times S,D)\Rightarrow \iRHom(0,D)\]
	but as the limit is unique this means
	\[0=\iRHom(0,D)=\iRHom(\R, D).\]

	Thus it suffices to show $H^{j}()\R^r\times S,D)\simeq H^j(\ast\times S,D)$.
	But this follows by knowing that the condensed cohomology on these locally compact Hausdorff spaces agrees with the classical {\v C}ech cohomology, and using known results there.
        See 3.2 and 3.4 \cite{scholze2019condensed} for a detailed proof of this, our rough intuition is that first we move the colimit $\R^n=\varinjlim [-n,n]^r$
        out of the colimit in the first coordinate of the $\RHom$,
        and afterwards the result for $[-n,n]^r$ is classical as it is homotopy-equivalent to a point $\ast$.

\item 	Now the result (iii) follows for finite products by resolving

\begin{center}\begin{tikzcd}
	0 & \Z & \R & \T & 0
	\arrow[from=1-1, to=1-2]
	\arrow[from=1-2, to=1-3]
	\arrow[from=1-3, to=1-4]
	\arrow[from=1-4, to=1-5]
\end{tikzcd}\end{center}
and using the induced long exact sequence
\begin{center}\begin{tikzcd}
	\dots & {H^{i+1}(\T,D_B)} & {H^i(\Z,D_B)} & {H^i(\R,D_B)} & {H^i(\T,D_B)} & \dots
	\arrow[from=1-2, to=1-1]
	\arrow[from=1-3, to=1-2]
	\arrow[from=1-4, to=1-3]
	\arrow[from=1-5, to=1-4]
	\arrow[from=1-6, to=1-5]
\end{tikzcd}\end{center}
as well as (i) and (ii).

For arbitrary products one again uses comparison to {\v Cech}-cohomology, see 3.4 in \cite{scholze2019condensed}.

\item To see that $\iRHom(\prod T, \R)=0$, again we only sketch parts of the argument from 3.4 in \cite{scholze2019condensed}.
One can show (3.3. in \cite{scholze2019condensed}), that $\iRHom(A, \R)(S)$ is computed by

\begin{center}\begin{tikzcd}
	0 & {\bigoplus_{j=1}^{n_0} C(A^{r_{0,j}}\times S, \R)} & {\bigoplus_{j=1}^{n_1}C(A^{r_{1,j}}\times S, \R)} & \dots
	\arrow[from=1-1, to=1-2]
	\arrow[from=1-2, to=1-3]
	\arrow[from=1-3, to=1-4]
\end{tikzcd}\end{center}
Now, we use that multiplication with $2$ on $A$ is bounded, but on $\R$ unbounded; see for the rest of this nice argument.\cite{scholze2019condensed}.

\item Resolving $\T$, $\iRHom(\prod \T, \T)$ fits in

\begin{center}\begin{tikzcd}
	{H^i(\prod\T, \R)} & {H^i(\prod \T,\T)} & {H^{i+1}(\prod \T, \Z)} & {H^{i+1}(\prod\T, \R)}
	\arrow[from=1-1, to=1-2]
	\arrow[from=1-2, to=1-3]
	\arrow[from=1-3, to=1-4]
\end{tikzcd}\end{center}
		and thus, using (iii) and (iv) we conclude $\iRHom(\prod_I \T, \T)=\bigoplus_I \Z[0]$.

\item Now $\iRHom(\R,\R)=\R[0]$ follows by resolving the left coordinate, obtaining

\begin{center}\begin{tikzcd}
	{H^{i-1}(\T, \R)} & {H^i(\Z, \R)} & {H^i(\R, \R)} & {H^i(\T, \R)} & {H^{i+1}(\Z, \R)}
	\arrow[from=1-1, to=1-2]
	\arrow[from=1-2, to=1-3]
	\arrow[from=1-3, to=1-4]
	\arrow[from=1-4, to=1-5]
\end{tikzcd}\end{center}
\item We lastly solve $\iRHom(\R, \T)=\R[0]$ by resolving $\T$, inducing

\begin{adjustbox}{max width = \textwidth}
\begin{tikzcd}
	\dots & {H^1(\R,\R)} & {H^1(\R, \T)} & {H^0(\R, \Z)} & {H^0(\R, \R)} & {H^0(\R, \T)} & 0 \\
	& 0 & {H^1(\R,\T)} & 0 & \R & {H^0(\R,\T)} & 0
	\arrow[from=1-1, to=1-2]
	\arrow[from=1-2, to=1-3]
	\arrow[from=1-3, to=1-4]
	\arrow[from=1-4, to=1-5]
	\arrow[from=1-5, to=1-6]
	\arrow[from=1-6, to=1-7]
	\arrow[from=2-2, to=2-3]
	\arrow[from=2-3, to=2-4]
	\arrow[from=2-4, to=2-5]
	\arrow[from=2-5, to=2-6]
	\arrow[from=2-6, to=2-7]
\end{tikzcd}\end{adjustbox}
	\end{enumerate}
	\end{proof}

Using this, we can compute all internal $H^i(A,B)$ for locally compact abelian groups $A,B$, as explained in \cite{Aparicio2021}.

For any two locally compact abelian groups $A, B$ we can compute $\underline{\RHom}(A,B)$ as follows.
The functor $\iRHom(-,-)$ is additive in both coordinates, hence
	\[\iRHom(A, B)=\iRHom(\R^{n_A}, \R^{n_B})\oplus \iRHom(\R^{n_A}, T_B)\oplus \iRHom(T_A,\R^{n_B})\oplus \iRHom(T_A,T_B).\]
Now we treat those cases individually.

\begin{lemma}
	\[\iRHom(\R^{n_A},\R^{n_B})=\R^{n_{A}n_B}[0]\]
\end{lemma}
\begin{proof}
	As $\iRHom$ is additive in both coordinates,
	\[\iRHom(\R^{n_A},\R^{n_B})=\iRHom(\R,\R)^{n_A n_B}.\]
	But $\iRHom(\R,\R)$ can be computed using the resolution
\begin{center}\begin{tikzcd}
	0 & \Z & \R & {\R/\Z} & 0
	\arrow[from=1-1, to=1-2]
	\arrow[from=1-2, to=1-3]
	\arrow[from=1-3, to=1-4]
	\arrow[from=1-4, to=1-5]
\end{tikzcd},\end{center}
	which induces the long exact sequence

\begin{center}\begin{tikzcd}
	\dots & {H^{i+1}(\T,\R)} & {H^i(\Z,\R)} & {H^i(\R,\R)} & {H^i(\T,\R)} & \dots
	\arrow[from=1-2, to=1-1]
	\arrow[from=1-3, to=1-2]
	\arrow[from=1-4, to=1-3]
	\arrow[from=1-5, to=1-4]
	\arrow[from=1-6, to=1-5]
\end{tikzcd}.\end{center}
	As $\RHom(\T,\R)=0$, we conclude
	\[\iRHom(\R, \R)\simeq \iRHom(\Z,\R)\simeq \R[0]\].
\end{proof}

Next, we use the resolution of $T_B$ to solve the next term.
\begin{lemma}
	\[\iRHom(\R^{n_A},T_B)=(\ker \phi_*)^{n_A}[0],\]
	where $\phi_*$ is the induced morphism
	\[\prod_{I_B}\R\to \prod_{J_B}\R\]
	induced by the decomposition of $K_B$ into tori.
\end{lemma}
\begin{proof}
	First, note that $\iRHom(\R^{n_A},T_B)=(\iRHom(\R, T_B))^{n_A}$

The short exact sequence

\begin{center}\begin{tikzcd}
	0 & {K_B} & {T_B} & {D_B} & 0
	\arrow[from=1-1, to=1-2]
	\arrow[from=1-2, to=1-3]
	\arrow[from=1-3, to=1-4]
	\arrow[from=1-4, to=1-5]
\end{tikzcd}\end{center}
induces the long exact sequence (a distinguished triangle)

\begin{adjustbox}{max width = \textwidth}
\begin{tikzcd}
	\dots & {H^{i-1}(\R,D_B)} & {H^i(\R,K_B)} & {H^i(\R,T_B)} & {H^i(\R,D_B)} & {H^{i+1}(\R,K_A)} & \dots
	\arrow[from=1-1, to=1-2]
	\arrow[from=1-2, to=1-3]
	\arrow[from=1-3, to=1-4]
	\arrow[from=1-4, to=1-5]
	\arrow[from=1-5, to=1-6]
	\arrow[from=1-6, to=1-7]
\end{tikzcd}\end{adjustbox}
As $\iRHom(\R, D_B)=0$, we conclude
	\[\iRHom(\R, T_B)=\iRHom(\R, K_B)\]

	For $\iRHom(\R, K_B)$, we use the resolution of the compact abelian group $K_B$ into tori
	\begin{center}\begin{tikzcd}
	0 & K_B & {\prod_{I_B} \T} & {\prod_{J_B}\T} & 0
	\arrow[from=1-1, to=1-2]
	\arrow[from=1-2, to=1-3]
	\arrow[from=1-3, to=1-4]
	\arrow[from=1-4, to=1-5]
\end{tikzcd}\end{center}

		This yields $\iRHom(\R, K_B)=0$ for $i>0$ and

\begin{center}\begin{tikzcd}
	0 & {H^0(\R, K_B)} & {\prod_{I_B}\R} & {\prod_{J_B} \R} & 0
	\arrow[from=1-1, to=1-2]
	\arrow[from=1-2, to=1-3]
	\arrow[from=1-3, to=1-4]
	\arrow[from=1-4, to=1-5]
\end{tikzcd}\end{center}
\end{proof}

\begin{lemma}
	\[\iRHom(T_A, \R^{n_B})=(\ker {\psi_A}^*)^{n_B}[0]+ (\coker \psi_A^*)^{n_B}[-1]\]
	where $\psi_A^*$ is the induced arrow

\begin{center}\begin{tikzcd}
	{\bigoplus_{Q_A}\R} & {\bigoplus_{R_A}\R}
	\arrow["{\psi_A}^*"', from=1-2, to=1-1]
\end{tikzcd}\end{center}
	\end{lemma}
\begin{proof}
Again $\iRHom(T_A, \R^{n_B})=\iRHom(T_A,\R)^{n_B}$.
	By resolving $T_A$ into $D_A$ and $K_A$,
	we reduce to $\iRHom(K_A, \R)$ and $\iRHom(D_A,\R)$.
But since $K_A$ can be resolved by products of tori, and $\iRHom(\prod T, \R)=0$,
	we conclude $\iRHom(K_A,\R)=0$ and hence $\iRHom(T_A,\R)\simeq \iRHom(D_A, \R)$.
	Using the resolution of $D_A$ into direct sums of $\Z$
	as
\begin{center}\begin{tikzcd}
	0 & {\bigoplus_{Q_A}\Z} & {\bigoplus_{R_A}\Z} & {D_A} & 0
	\arrow[from=1-1, to=1-2]
	\arrow["{\psi_A}", from=1-2, to=1-3]
	\arrow[from=1-3, to=1-4]
	\arrow[from=1-4, to=1-5]
\end{tikzcd}\end{center}
	This leads to $H^i(T_A,\R)=0$ for $i>1$ and

\begin{center}\begin{tikzcd}
	0 & {H^1(T_A,\R)} & {\bigoplus_{Q_A}\R} & {\bigoplus_{R_A}\R} & {H^0(T_A,\R)} & 0
	\arrow[from=1-2, to=1-1]
	\arrow[from=1-3, to=1-2]
	\arrow["{\psi_A}^*"', from=1-4, to=1-3]
	\arrow[from=1-5, to=1-4]
	\arrow[from=1-6, to=1-5]
\end{tikzcd}\end{center}

\end{proof}

\begin{lemma}
	The last term,
	\[\iRHom(T_A,T_B)\]
	can be computed as below, and in particular fulfills $H^i=0$ for $i>1$.

\end{lemma}
\begin{proof}
	We resolve $T_A$ into $D_A$ and $K_A$, and compute $\iRHom(D_A,T_B)$ and $\iRHom(K_A,T_B)$.
	\begin{itemize}
		\item $\iRHom(D_A,T_B)$ can be computed by resolving $D_A$ and hence is given by $H^i(D_A,T_B)=0$ for $i>1$ and

\begin{center}\begin{tikzcd}
	0 & {H^1(D_A,T_B)} & {\bigoplus_{Q_A}T_B} & {\bigoplus_{R_A}T_B} & {H^0(D_A,T_B)} & 0
	\arrow[from=1-2, to=1-1]
	\arrow[from=1-3, to=1-2]
	\arrow["{\psi_A}^*"', from=1-4, to=1-3]
	\arrow[from=1-5, to=1-4]
	\arrow[from=1-6, to=1-5]
\end{tikzcd}.\end{center}
		Hence $\iRHom(D_A,T_B)=(\ker {\psi_A}^*)^{n_B}[0]+ (\coker \psi_A^*)^{n_B}[-1]$.
		\item $\iRHom(K_A, T_B)$ can be computed by first resolving $K_A$, leading to

\begin{center}
\begin{tikzcd}
	{H^i(\prod_{I_A}\T, T_B)} & {H^i(\prod_{J_A}\T, T_B)} & {H^i(K_A,T_B)} & {H^{i-1}(\prod_{I_A}\T, T_B)}
	\arrow["{\phi_A^*}"', from=1-2, to=1-1]
	\arrow[from=1-3, to=1-2]
	\arrow[from=1-4, to=1-3]
\end{tikzcd}\end{center}
		Thus it suffices to compute $\iRHom(\prod\T, T_B).$
		But for this we again can resolve $T_B$, and knowing that the discrete part is $\iRHom(\prod \T, D_B)=\bigoplus D_B[-1]$,
		and knowing that $\iRHom(\prod_I\T,\T)=\bigoplus_I \Z[0]$, to compute $H^i(\prod_{J_A}\T, T_B)$ we consider the exact diagram

\adjustbox{max width = \textwidth}{\begin{tikzcd}
	&& {H^{i}(\prod_{J_A}\T, D_B)} \\
	& {H^i(\prod_{I_A}\T, T_B)} & {H^i(\prod_{J_A}\T, T_B)} & {H^i(K_A,T_B)} & {H^{i-1}(\prod_{I_A}\T, T_B)} \\
	{H^{i}(\prod_{J_A}\T, \prod_{I_B}\T)} & {H^i(\prod_{J_A}\T, \prod_{J_B}\T)} & {H^i(\prod_{J_A}\T, K_B)} & {H^{i-1}(\prod_{J_A}\T,\prod_{I_B}\T)} \\
	&& {H^{i}(\prod_{J_A}\T, K_B)}
	\arrow[from=1-3, to=2-3]
	\arrow["{\phi_A^*}"', from=2-3, to=2-2]
	\arrow[from=2-3, to=3-3]
	\arrow[from=2-4, to=2-3]
	\arrow[from=2-5, to=2-4]
	\arrow["{\phi_B^*}", from=3-1, to=3-2]
	\arrow[from=3-2, to=3-3]
	\arrow[from=3-3, to=3-4]
	\arrow[from=3-3, to=4-3]
\end{tikzcd}}
		\begin{itemize}
		\item For $i=0$ this is

\begin{center}\begin{tikzcd}
	&& 0 \\
	& {H^0(\prod_{I_A}\T, T_B)} & {H^0(\prod_{J_A}\T, T_B)} & {H^i(K_A,T_B)} & 0 \\
	{\prod_{I_B}\bigoplus_{J_A}\T} & {\prod_{J_B}\bigoplus_{J_A}\T} & {H^0(\prod_{J_A}\T, K_B)} & 0 \\
	&& 0
	\arrow[from=1-3, to=2-3]
	\arrow["{\phi_A^*}"', from=2-3, to=2-2]
	\arrow[from=2-3, to=3-3]
	\arrow[from=2-4, to=2-3]
	\arrow[from=2-5, to=2-4]
	\arrow["{\phi_B^*}", from=3-1, to=3-2]
	\arrow[from=3-2, to=3-3]
	\arrow[from=3-3, to=3-4]
	\arrow[from=3-3, to=4-3]
\end{tikzcd}\end{center}

and we conclude $H^0(\prod_{J_A}\T, T_B)=\coker {\phi_B}_*$.
		Analogously, we conclude $H^0(\prod_{I_A}\T, T_B)=\coker{\phi_B}_{*}$, in total yielding

\begin{center}\begin{tikzcd}
	{\prod_{I_B}\bigoplus_{I_A}\T} & {\prod_{I_B}\bigoplus_{J_A}\T} \\
	{\prod_{J_B}\bigoplus_{I_A}\T} & {\prod_{J_B}\bigoplus_{J_A}\T} \\
	{\coker{\phi_B}_*} & {\coker{\phi_B}_*} & {H^0(K_A,T_B)=\ker \phi_A^*} & 0 \\
	0 & 0
	\arrow["{{\phi_B}_*}", from=1-1, to=2-1]
	\arrow["{{\phi_B}_*}", from=1-2, to=2-2]
	\arrow[from=2-1, to=3-1]
	\arrow[from=2-2, to=3-2]
	\arrow[from=3-1, to=4-1]
	\arrow["{\phi_A^*}", from=3-2, to=3-1]
	\arrow[from=3-2, to=4-2]
	\arrow[from=3-3, to=3-2]
	\arrow[from=3-4, to=3-3]
\end{tikzcd}\end{center}

		\item For $i=1$ the above diagram reduces to calculating $H^1(K_A,T_B)$ via

\begin{adjustbox}{max width = \textwidth}
\begin{tikzcd}
	0 & 0 & 0 & 0 & 0 \\
	{\bigoplus_{I_A}D_B} & {\bigoplus_{J_A} D_B} & {H^1(K_A,T_B)} & {\coker{\phi_B}_{*}} & {\coker {\phi_B}_*'} \\
	{H^i(\prod_{I_A}\T, T_B)} & {H^1(\prod_{J_A}\T, T_B)} & {H^1(K_A,T_B)} & {H^{0}(\prod_{I_A}\T, T_B)} & {H^0(\prod_{J_A}\T,T_B)} \\
	0 & 0 & 0 & 0 & 0
	\arrow[from=1-1, to=2-1]
	\arrow[from=1-2, to=2-2]
	\arrow[from=1-3, to=2-3]
	\arrow[from=1-4, to=2-4]
	\arrow[from=1-5, to=2-5]
	\arrow[from=2-1, to=3-1]
	\arrow[from=2-2, to=2-1]
	\arrow[from=2-2, to=3-2]
	\arrow[from=2-3, to=2-2]
	\arrow[from=2-3, to=3-3]
	\arrow[from=2-4, to=2-3]
	\arrow[from=2-4, to=3-4]
	\arrow[from=2-5, to=2-4]
	\arrow[from=2-5, to=3-5]
	\arrow[from=3-1, to=4-1]
	\arrow["{\phi_A^*}"', from=3-2, to=3-1]
	\arrow[from=3-2, to=4-2]
	\arrow[from=3-3, to=3-2]
	\arrow[from=3-3, to=4-3]
	\arrow[from=3-4, to=3-3]
	\arrow[from=3-4, to=4-4]
	\arrow[from=3-5, to=3-4]
	\arrow[from=3-5, to=4-5]
\end{tikzcd}\end{adjustbox}
			\item For $i\ge 2$ we obtain $H^i(K_A, T_B)=0$.
		\end{itemize}

\end{itemize}
Now, knowing $\iRHom(D_A,T_B)$ and $\iRHom(K_A,T_B)$ we obtain $\iRHom(T_A,T_B)$ via

\begin{center}\begin{tikzcd}
	{H^{i+1}(K_A,T_B)} & {H^i(D_A,T_B)} & {H^{i}(T_A,T_B)} & {H^i(K_A,T_B)} & {H^{i-1}(D_A,T_B)}
	\arrow[from=1-2, to=1-1]
	\arrow[from=1-3, to=1-2]
	\arrow[from=1-4, to=1-3]
	\arrow[from=1-5, to=1-4]
\end{tikzcd}\end{center}
\end{proof}

This shows that the category of locally compact abelian groups can be seen as part of the condensed theory, and that there, things really can be computed elementarily.

\begin{remark}
  Surely, many aspects of the preceding discussion can be translated swiftly to $\CondModR$, at least for nice enough $R$ (for example, discrete or $\R$).
\end{remark}


\section{Analytic rings}
So far, we have seen the topological concepts of compactness being generalised to quasicompactness,
as well as the replacement for Hausdorffness being quasiseparatedness.
In this section we study another important concept in topology -- completeness.

The first major problem one has to deal with when finding a proper replacement for completeness, is that classically,
completeness of topological spaces is not an intrinsic property of the space.
For talking about completeness, we need additional information
in form of some uniformity on the space.
In particular, this means that a given space is supposed to admit multiple versions of completeness.

Some sort of uniformity comes up whenever one has additional algebraic structure because relative closeness of points can be compared independently of their absolute position.
This means that whenever we have an addition map that gives a way to \enquote{move} the space around, the concept of $x$ and $y$ being close together can be compared to
$x-a$ and $y-a$ being close together.
This explains why we stick to finding a notion of completeness for abelian groups or more generally $R$-modules.

\subsection{From classical completeness to preanalytic rings}
The notion of completeness is heavily set-theoretic and hence we have to reformulate it in a way that is more accessible to categorification.
This is done best by not working with a uniformity as a system of sets, but rather with pseudometrics.
A particularly important class of complete objects is that of Banach spaces, so we will first concentrate on these.

\begin{warning}[Skip if you're a functional analyst]
 As usual in functional analysis, the term \enquote{vector space} implicitly means that the base field is $\R$ or $\C$.
 Usually, one does not care which one of these two is used, but we will mostly stick to $\R$-vector spaces.
\end{warning}

The first step in a reformulation of completeness of Banach spaces is to get rid of Cauchy sequences.
This can be done as follows.
\begin{lemma} A (semi)normed vector space $(E,\|\cdot\|)$ is complete precisely if for any sequence $(x_n)\sub E$
    that is absolutely summable,
    \[\sum\|x_n\|<\infty,\]
    there exists a sum
    \[\sum x_n=\lim_{n\to \infty} \sum_{k\le n}x_k.\]
    \end{lemma}
\begin{proof}
    This can be found, e.g., in \cite[I.1.8]{werner2006funktionalanalysis}.
  \end{proof}
This reformulation has the advantage that it avoids the notion of Cauchy sequence.
However, we would rather like $(x_n)$ to be a convergent sequence than some absolutely summable sequence,
since a convergent sequence is something we can
describe neatly as a continuous function $\alpha\N\to E$.

Hence we reformulate this a bit further.
\begin{lemma}
    A normed vector space $E$ is complete precisely if for every convergent sequence $(y_n)\sub E$,
    and every absolutely summable sequence $a_i\in \ell^1(\R)$, the sum
    \[\sum a_i y_i\]
    exists.
\end{lemma}
\begin{proof}
    First assume that $E$ is complete.
    Then clearly for any $a_i\in \ell^1$ and convergent $(y_n)\sub E$, the convergence of $y_n$ implies in particular boundedness of $\|y_n\|$ by some constant $C$,
    and hence the sequence $(x_n)=(a_ny_n)$ is absolutely summable by
    \[\sum \|x_n\|=\sum |a_n|\|y_n\|\le C\sum |a_n|<\infty.\]
    Hence the sum $\sum x_n=\sum a_n y_n$ exists.

    For the converse implication, take any absolutely summable sequence $x_n$ in $E$.
    We need to show existence of $\sum x_n$.
The idea is to use the fact that there is no \enquote{minimal growth rate} in $\ell^1$, and use this to split $x_n$ into a absolutely summable sequence of scalars
(approximately $a_n=\|x_n\|$) and taking this \enquote{out} of the sequence $x_n$, $y_n=x_n/a_n$.

To make this precise first we show that for any sequence $b_n\in \ell^1$ there exists a decomposition
    \[b_n=a_n c_n\]
    with $(a_n)\in \ell^1$ and $(d_n)\in c_0$.
    This can be achieved by a standard \enquote{analysis 1} tactic as follows.
    For every $k\in \N$, there is a minimal index $N_k$ such that
    \[\sum_{n\ge N_k} |b_n|\le 2^{-2k}.\]
    Define
    \[d_n= 2^{-k} \qquad \text{ with } N_k\le n< N_{k+1}\]
    Then clearly $c_n\to 0$, and for
    \[a_n=b_n/d_n=b_n 2^k\]
    we obtain
    \[\sum_{N_{k}\le n< N_{k+1}} |a_n|= 2^{k}\sum_{N_k\le n< N_{k+1}} |b_n|\le 2^k \sum_{n\ge N_k}|b_n|\le 2^{-k}\]
    And thereby
    \[\sum |a_n|\le \sum_{k\in \N} 2^{-k}=1<\infty.\]

    Using this, consider any sequence $x_n$ in $E$ which is absolutely summable (hence in particular bounded), and define
    $a_n$, $d_n$ to be such a decomposition of $b_n=\|x_n\|$.

    Now, define $y_n=(d_n x_n)/\|x_n\|$ (and $y_n=0$ if $x_n=0$).
    Then $\|y_n\|=d_n$ (and $0$ if $\|x_n\|=0$), hence $y_n\in c_0$, and furthermore
    \[a_{n}y_n=a_n d_n x_n/\|x_n\|=x_n.\]

Thus this is of the form of a convergent sequence weighted by a $\ell^1$-sequence.
    This implies the existence of
    \[\sum x_n=\sum a_n y_n\]
and thereby the completeness.
\end{proof}

This characterisation has the advantage that we have a good use of convergent sequences as just maps $\{\alpha\N\to E\}=E(\alpha\N)$.
But clearly, as soon as one leaves the normed setting and moves to non-metrisable spaces, the ability to sum just sequences will not suffice.
Hence we need another shift in our point of view that is generalisable to spaces other than $\alpha\N$.

Now, the key observation is to understand the appearance of $\ell^1$-sequences.
The central idea is the following:
\begin{center}
    \textbf{One can view an element of $\ell^1$ as a (finite signed Radon) measure on $\N$.}
\end{center}
This sheds a completely different light on completeness -- the sums clearly are just integrals with respect to these measures,
and \textbf{completeness can be described as the ability to integrate.}

\begin{lemma}
A normed vector space $E$ is complete precisely if for any continuous map $f\colon \alpha\N\to E$,
and any measure $\mu\in \M_{\mathrm{Radon}}(\alpha\N)$, the Bochner integral
\[\int f\mathrm{d}\mu\]
exists.
\end{lemma}
\begin{proof}
    For any such measure $\mu$ define $(a_n)=(\mu(\{n\}))\in \ell^1$, and combine this with $a_\infty=\mu(\{\infty\})$.
    Then for any $f\colon \alpha\N\to E$, the integral clearly is given by
    \[\int f \mathrm{d}\mu =a_\infty\lim_{n\to \infty}f(n)+\sum a_n f(n)\]
    and hence existent as soon as $E$ is complete.

    Conversely, it suffices to use measures with $\mu(\{\infty\})=0$.
    \end{proof}

Clearly, the ability to integrate with respect to some measures is one of the key aspects of analysis and functional analysis, and this viewpoint will be the guiding idea of analytic rings.
However, the general theory of vector-valued integrals is surprisingly subtle,
as, e.g., especially in the locally convex setting we need to distinguish weak integrals, Pettis integrals, Bochner integrals, etc.\
(see, e.g., \cite{diestel1978vector} and the references therein for vector integration).
This fact boils down to the map $\mu\mapsto \int f \mathrm{d}\mu$ being continuous with respect to different topologies on the space of measures.

Hence, especially if we want a non-metrisable theory and consider measures on more complicated spaces, we need to include the \textbf{space of measures} into our considerations.
There are a priori many different ways to topologise this space, as it is described as the dual space of the space of continuous functions (by Riesz),
\[\mcM_{\mathrm{Radon}}(K)=C(K)'=\hom(K, \R).\]

Now, a functional analytically minded person would expect the space of Radon measures to be equipped with the weak*-topology (or Mackey topology or even weak or norm topology),
but having seen the internal $\hom$ of condensed abelian groups, it is clear that the correct topology to install is the compact open topology (which also fits well with Arzela-Ascoli):
\[
    \mcM_{\mathrm{Radon}}(K)=\ihom(K,\R).
\]

Let us compare this topology to the classical functional analytic ones.

 As the weak*-convergence is given by uniform convergence on finite sets and the compact open topology by uniform convergence on compact sets, clearly
    the compact open topology is finer.
Because the internal dual in condensed $R$-modules fulfills $\ihom(E',\R)\simeq E$ (as we will see later), this is also true in locally convex spaces, and thereby
the Mackey topology is finer than the compact open topology.

Since the Mackey topology as well as the weak*-topology have dual space $E$, they have compact unit ball.
This implies that the compact open topology makes the unit ball compact.
These topologies characterise dual spaces, see \cite{kaijser1977note}.
As they are comparable, we note that on norm-bounded subsets all three topologies induce comparable compact Hausdorff topologies, and hence they agree.
Thus the compact open topology is really similar to the weak*-topology.

Furthermore, the compact open topology can be described as the final topology with respect to the embeddings
\[E'_{\|\cdot\|\le n}\hookrightarrow E',\]
and hence
\[E'_{KO}=\bigcup E'_{\le n}.\]

Using this, we have the following reformulation of completeness of Banach spaces.
\begin{lemma}
A normed vector space $E$ is complete precisely if for all $f\colon \alpha\N \to E$ there is a unique continuous linear extension

\begin{center}\begin{tikzcd}
	{\mcM_{\mathrm{Radon}}(\alpha\N)} & E \\
	\alpha\N
	\arrow[dashed, from=1-1, to=1-2]
	\arrow[hook, from=2-1, to=1-1]
	\arrow["f"', from=2-1, to=1-2]
\end{tikzcd}\end{center}
    where $\mcM_{\mathrm{Radon}}(\alpha\N)$ is equipped with compact open topology.
\end{lemma}

Finally, we take the step out of the metrisable theory by allowing varying $S$ instead of $\alpha\N$.

\begin{lemma}\label{lem:locally_convex_radon_comp}
Consider any Hausdorff complete locally convex topological vector space $E$.
Let $S\in \prof$, and install the compactly generated topology on $\mcM_{\mathrm{Radon}}(S)$ induced by the compact unit ball (the compact open topology).
Then for any map $f\colon S\to E$
there exists a unique continuous linear map $\mu\mapsto \int f d\mu$ from $\mcM_{\mathrm{Radon}}(S)\to E$,
extending the continuous Dirac embedding $S\to \mcM_{\mathrm{Radon}}(S)$,

\begin{center}\begin{tikzcd}
	{\mcM_{\mathrm{Radon}}(S)} & E \\
	S
	\arrow[dashed, from=1-1, to=1-2]
	\arrow[hook, from=2-1, to=1-1]
	\arrow["f"', from=2-1, to=1-2]
\end{tikzcd}.\end{center}

\end{lemma}
\begin{proof}
    This is 3.4 in \cite{scholze2019Analytic}.

    The uniqueness is clear by density of Dirac measures in the unit ball with respect to weak*-topology, which implies density of Dirac measures with
    compact open topology.
    By definition of final topologies, this implies that finite linear combinations of Dirac measures are dense in $\mcM_{\mathrm{Radon}}(S)$.
    Hence the extension is unique, if it exists (as all the spaces are Hausdorff).

    Consider $f\colon S\to E$ for $S\in \prof$ and complete locally convex $E$.
    We want to show that there exists a unique extension to $\mcM_{\mathrm{Radon}}(S)$.
    By definition of the final topology of the embeddings, it suffices to find a continuous linear extension
    to those radon measures with norm $\le 1$.

    We now define the integral in the classical Bochner/Riemann way by approximating with step functions.

    For this we need to partition $S$ into finer and finer \enquote{intervals}.
    Luckily, $S$ is profinite, hence $S=\varprojlim S_i$ (with surjective transition maps) does precisely this.
    So choose any such clopen partition $S_i$, meaning that $T_k=\pi_i^{-1}(k)\ne \emptyset$ for $k\in S_i$ and $\pi_i\colon S\to S_i$ partition $S$,
    \[S=\bigsqcup_{k\in S_i}T_k\]
    Now, for any such partition $S_i$ we pick \enquote{base points} $t_i(k)$ of the \enquote{intervals} $T_k$.

    The correct way to state this is to take $t_i$ to be a section of $S\to S_i$,

\begin{center}\begin{tikzcd}
	S \\
	{S_i}
	\arrow[two heads, from=1-1, to=2-1]
	\arrow["{t_i}", shift left=3, from=2-1, to=1-1]
\end{tikzcd}\end{center}
    meaning that for any $k\in S_i$ the value $t_i(k)$ is contained in the clopen subset $T_k=\pi_i^{-1}(k)$ of $S$ that belongs to the index $k$ of the clopen partition.

    Similarly to the case of defining Riemann integrals by evaluating the function in the middle of each intervall or the right or left of the intervall,
    this choice does not really matter in the end.

    Now we define the $i$-th approximation of the integral of $f$ along the measure $\mu$ (again in complete analogy to the Riemann/Bochner integral)
    \[v_i(\mu)=\sum_{k\in S_i} f(t_i(k))\mu(T_k)\in E\]
    We want to define
    \[\int f \, d\mu=\lim_{i\in I} v_i\in E\]
    Now, the index category is filtered, in fact a directed set, and hence a net.
    Hence it suffices to show that this net is Cauchy.

    For any  $0$-neighbourhood (w.l.o.g.\ being absolutely convex absorbant or simply a barrel), we need to find an index $i_0$
    such that for all $i, j\ge i_0$ one has $v_i-v_j\in U$.
    Extend $f-f\colon S\times S\to E$ via $(f-f)(t,r)=f(t)-f(r)$.

    Now, $(f-f)^{-1}(U)$ is a neighbourhood of $\Delta\sub S\times S$, and as $\{T_k\times T_{k'} \, k\in S_i, \, k'\in S_{i'}\}$ forms a basis of the topology of
    $S\times S$, there is a open cover of $\Delta$ contained in $(f-f)^{-1}(U)$ with sets of the form $\{T_k\times T_{k'}, \, k\in S_i, \, k'\in S_{i'}\}$.
    As $\Delta$ is compact, this admits a finite subcover.
    As the diagram $S_i$ is filtered (alternative language: the clopen partitions allow finite refining),
    there exists $i_0\in I$ such that all the sets $T_{k}\times T_{k'}$ in the covering can be written with $k,k'\in S_{i_0}$.
    Using that the sets $T_k$  for $k\in S_{i_0}$ are pairwise disjoint, we know that for every $k\in S_{i_0}$, the set $T_k\times T_k$ is in the cover (as those sets cover the diagonal)
    This implies that for any $t_{i_0}(k), t_{i_0}(k)'\in T_k$ the difference $f(t_{i_0})(k))-f(t_{i_0}(k)')\in U$.
   Now, take any $i,j\ge i_0$.
    We want to show that
        \[v_{i}(\mu)-v_{j}(\mu)=\sum_{k\in S_i} f(t_i(k))\mu(T_k^i)- \sum_{k'\in S_j} f(t_j(k'))\mu(T_{k'}^j)\in U\]
    for any $\mu$ of norm $\le 1$.
    Refine $i$ and $j$ to $h\ge i,j$, and define for any $\ell\in S_h$ the value $i(\ell)$ as the image of $S_h\to S_i$ and $j(\ell)$ as the image under $S_h\to S_j$.
    By finite additivity of $\mu$ we can use the partition $T_k^{i}=\bigsqcup_{i(\ell)=k} T_\ell^{h}$ to see $\mu(T_k^i)=\sum_{i(\ell)=k}\mu(T_\ell^{h})$, and thus
    \[\sum_{k\in S_i}f(t_i(k))\mu(T_k^i)=\sum_{\ell \in S_h}f(t_{i}(i(\ell)))\mu(T_\ell^h)\]
Analogously, we obtain
\[\sum_{k'\in S_j}f(t_j(k'))\mu(T_{k'}^j)=\sum_{\ell\in S_h}f(t_j(j(\ell)))\mu(T_\ell^h).\]
    In total this yields
    \[v_i(\mu)-v_j(\mu)=\sum_{\ell\in S_h} (f(t_i(i(\ell)))-f(t_j(j(\ell))))\mu(T_\ell^h).\]
    Now we note that for every $\ell\in S_h$,  $(i(\ell))$ and $(j(\ell))$ are mapped to the same $p\in S_{i_0}$ (as the diagram is thin).
    Hence $t_i(i(\ell))\in T_{i(\ell)}^i\sub T_p^{i_0}$ and $t_j(j(\ell))\in T_{j(\ell)}^j\sub T_p^{i_0}$,
    i.e., they are contained in the same clopen set of the partition $S_{i_0}$.
    This implies that for all $\ell\in S_h$, the difference
    \[d_\ell=f(t_i(i(\ell)))-f(t_j(j(\ell)))\in U\]
    is contained in $U$.
Now, since $U$ is absolutely convex and $\sum_{\ell \in S_h}|\mu(T_\ell^{h})|\le 1$, this finally implies that
    \[v_i(\mu)-v_j(\mu)=\sum_{\ell\in S_h} d_\ell \mu(T_\ell^h)\in U.\]

    So we have shown that $v_i$ forms a Cauchy net in $E$.
    By completeness of $E$ thereby the integral
    \[\int fd\mu=\lim_{i\in I}v_i\in E\]
    exists.
    Linearity of this map is clear as every approximation step is linear.

    It remains to show continuity of $\mu\mapsto \int f d\mu$.
    This follows from the independence of $i_0$ from $\mu$.
    In more detail, the topology on the measures is made in such a way that for any convergent net $(\mu_{i})$
and any fixed clopen partition $S_p$ of $S$, the values $\mu_i(T_k)$ converge for all $k\in S_p$.
    Hence consider any net $\mu_i\to 0$ (w.l.o.g.\ of bounded norm), and any absolutely convex $0$-neighbourhood $U$.
    Choose any closed absolutely convex $0$-neighbourhood $W$ with $W-W\sub U$.
    Choose $i_0$ large enough so that for all $j\ge i_0$ and $\mu$ with norm $\le 1$ the value
    \[v_{i_0}(\mu)-v_j(\mu)\in W\]
    (and hence also $\int fd\mu-v_i\in W$ by closedness of $W$) furthermore,
    as $W$ is absorbant, for every $k\in S_{i_0}$ there exists $\eps_k$ with $\eps_k f(t_{i_0}(k))\in W$.

    Now choose $\ell_0$ such that $|S_{i_0}|\mu_\ell(T_k^{i_0})\le \eps_k$ for all $k\in S_{i_0}$ and all $\ell\ge \ell_0$.
    We conclude for every $\ell\ge \ell_0$
    \[ v_{i_0}(\mu_\ell)=\sum_{k\in S_{i_0}} f(t_{i_0}(k))\mu_\ell(T_k^{i_0}) \in W\]
    and thereby
    \[\int f d\mu_\ell=\underbrace{\int f d\mu_\ell -v_{i_0}(\mu_\ell)}_{\in W}+\underbrace{v_{i_0}(\mu_\ell)}_{\in W}\in U\]
    and thus $\int f d\mu_\ell\to 0$ and the map is continuous.
\end{proof}

Now, from a condensed point of view, one cannot help but notice that the space $\alpha\N$ is profinite,
and that the statement \enquote{for all profinite spaces}
should be reduced to \enquote{for all $S\in \extr$}.

This is precisely the situation we want to generalise.
Hence let us see how concepts of completeness other than Banach spaces fit into this setting.
\begin{remark}
A $\Q_p$-vector space $E$ with topology being induced by a valuation $\|\cdot\|$ 
    is complete precisely if for any $\alpha\N\to E$ there is a unique continuous linear extension

\begin{center}\begin{tikzcd}
	{c_0\oplus \R} & E \\
	\alpha\N
	\arrow[dashed, from=1-1, to=1-2]
	\arrow[hook, from=2-1, to=1-1]
	\arrow["f"', from=2-1, to=1-2]
\end{tikzcd}\end{center}
    where the spaces are equipped with appropriate topologies.
\end{remark}
In this sense, $p$-adic completeness can be captured in the same way, simply changing the spaces of measures from Radon measures to a different one.

Analogously, one could install a trivial notion of completeness as follows.
\begin{remark}
    Consider any condensed ring $R$, then
    by definition of the free objects $R[-]$,
    for any condensed $R$-module and any map $S\to R$ for extremally disconnected $R$,
    there is always a unique continuous extension to

\begin{center}\begin{tikzcd}
	{R[S]} & E \\
	S
	\arrow[dashed, from=1-1, to=1-2]
	\arrow[hook, from=2-1, to=1-1]
	\arrow["f"', from=2-1, to=1-2]
\end{tikzcd}\end{center}
In this sense, any $R$-module is \enquote{trivially complete with respect to the trivial measures}.
\end{remark}
This choice leads to another way of looking at the space of measures -- it is in some sense the free complete object over $S$,
seen as the \enquote{completion} of the Dirac measures (given by $S$ itself).
This yields to the general idea of \enquote{wishing for a certain space to be the completion of spaces of dirac measures}
to obtain corresponding free objects and an induced notion of \enquote{completeness} given by integrability against these measures.

In total, this suggests, that we should stay flexible in the choice of measures if we want to include these different concepts of completeness into our theory.

\subsection{1-categorical analytic rings}
Recall that the main idea of the last section was to define completeness in the following way:
assign a space of measures to every extremally disconnected set, against which we want to be able to integrate.

Note that in the following $R$ is always a condensed ring and $\RMod$ is the category of condensed $R$-modules.
Although this is not strictly necessary except for tensor products, we implicitly assume $R$ to be commutative unital.
We will sometimes leave the term \enquote{condensed} implicit -- everything is condensed.

We are lead to the definition of a pre-analytic ring.
\begin{definition}\label{def:preanalytic_ring}
  A \idx{pre-analytic ring} consists of a condensed ring $R$ together with a functor
  \[\mcM\colon\extr\to \RMod, \, S\mapsto \mcM(S)\]
  and a natural transformation
    \[\delta\colon \id\to ?\mcM, \, \delta_S\colon S\to \mcM(S),\]
    such that $\mcM$ maps finite coproducts in $\extr$ to finite direct sums in $\RMod$.
\end{definition}

\begin{remark}
Unlike in classical topology, we want the information of completeness of modules to be a property of the base ring; this explains the terminology.
\end{remark}
The structure of a pre-analytic ring induces a natural definition of what \enquote{completeness} of modules should mean.
\begin{definition}[Completeness]\label{def:completeness}
Consider any pre-analytic ring $(R,\mcM,\delta)$ and $E\in \RMod$.
    We say that $E$ is $\mcM$-complete (or just \idx{complete}), if for every $f\colon S\to E$ there exists a unique condensed $R$-module homomorphism $h\colon \mcM\to E$ with

\begin{center}\begin{tikzcd}
	{\mcM(S)} & E \\
	S
	\arrow["h", dashed, from=1-1, to=1-2]
	\arrow["{\delta_S}", from=2-1, to=1-1]
	\arrow["f"', from=2-1, to=1-2]
\end{tikzcd}\end{center}
    The full subcategory of complete $R$-modules is denoted by $\cMod\sub \RMod$.
\end{definition}
\begin{remark}
       Equivalently to the diagram,
\begin{center}\begin{tikzcd}
	{\mcM(S)} & E \\
	S
	\arrow["h", dashed, from=1-1, to=1-2]
	\arrow["{\delta_S}", from=2-1, to=1-1]
	\arrow["f"', from=2-1, to=1-2]
\end{tikzcd}\end{center}
      in condensed sets, one could use the diagram

\begin{center}\begin{tikzcd}
	{\mcM(S)} & E \\
	{R[S]}
	\arrow["h", dashed, from=1-1, to=1-2]
	\arrow["{R[\delta_S]}", from=2-1, to=1-1]
	\arrow["{R[f]}"', from=2-1, to=1-2]
\end{tikzcd}\end{center}
       in condensed $R$-modules,
    using the transformation $R[\delta]\colon R\to \mcM, \, R[S]\to \mcM(S)$.

    Alternatively, this can be characterised by $\delta$ inducing the second isomorphy:
    \[\hom_{\cond}(S,E)\simeq \hom_{\RMod}(R[S], E)\simeq \hom(\mcM(S),E).\]
\end{remark}
We can define the analogous notion in the derived category $D(\RMod).$
\begin{definition}
    An object $E\in D(\RMod)$ is $\mcM$-complete, if $\delta$ induces a natural isomorphism
    \[\iRHom(R[S][0],E)\simeq \iRHom(\mcM(S)[0], E).\]
    The full subcategory of complete derived $R$-modules is denoted by $\cDMod\sub \DRMod$.
\end{definition}
\begin{remark}
It may a priori seem quite strong to require the isomorphism of the internal $\RHom$'s, and not just the regular $\RHom$.
However, once we have established the definition of analytic rings, the equality of the internal $\hom$ will be automatic from the $\RHom$ equality.
\end{remark}

This simple definition of a pre-analytic ring is too general to imply good properties for the categories of complete modules.
Hence we need to add another important ingredient: The spaces of measures themselves should be complete.

\begin{definition}

    An \idx{analytic ring} is a pre-analytic ring $(R,\mcM,\delta)$, such that for all index sets $I$, $J$, extremally disconnected sets $S_i,\, i\in I$, $T_j,\, j\in J$ and maps
    \[f\colon \bigoplus_i \mcM(S_i)\to \bigoplus_j \mcM(T_j),\]
    with kernel $K\in \RMod$, the object $K[0]\in D(\RMod)$ is $\mcM$-complete.

\end{definition}
Note that it is essential to use the derived category here.
This lack of \enquote{higher morphisms} in classical functional analysis in the opinion of the authors can be viewed as one of the essential concepts functional analysis lacks,
somewhat holding the theory back from access to more powerful results.

The definition of an analytic ring directly implies wonderful properties of the category of complete modules, see, e.g,. \cite[5.9]{scholze2019condensed} for the central argument,
or \cite[7.5]{scholze2019condensed} or \cite[6.5, 6.14]{scholze2019Analytic} for versions of this theorem.
\begin{theorem}
    Let $(R,\mcM, \delta)$ be an analytic ring.
     The full subcategory $\cMod\sub \RMod$ of complete condensed $R$-modules has the following properties.
    \begin{enumerate}[(i)]
        \item It is an abelian category stable under all limits, colimits and extensions.
        \item $\mcM(S)$ for extremally disconnected $S$ form a family of compact projective generators.
      \item The inclusion admits a left adjoint, the \idx{completion}, which is denoted by $X\mapsto X\otimes_{(R,\mcM)}(R,\mcM)$ or $\mcM(X)$.%
            \footnote{This is just notation.
            In analogy with the base change, $X\otimes_{R}(R,\mcM)$ is probably a more reasonable notation.
            $\mcM(X)$ is the most canonical notation.
            }
        Completion is the unique colimit preserving extension of $R[\delta]\colon R[S]\to \mcM(S)$.
        \item For commutative $R$, there is a unique symmetric monoidal tensor product $\otimes_{(R,\mcM)}$ on
        $\cMod$ making the completion symmetric monoidal.
        We call this the \idx{completed tensor product}.
        \item The category $\cMod$ is closed under internal $\hom$ and internal extensions, and in particular the tensor-hom adjunction holds.
    \end{enumerate}
    Furthermore, the complete $R$-modules $E$ can be equivalently described in the following ways.
    \begin{enumerate}[(a)]
        \item  Those $E\in \RMod$, such that $\delta$ induces an isomorphism
        \[\hom(R[S],E)\simeq \hom(\mcM(S),E)\]
        \item Those $E\in \RMod$, such that $\delta$ induces an isomorphism
        \[\ihom(R[S],E)\simeq \ihom(\mcM(S),E)\]
        \item There exists a map \[f\colon \bigoplus_{i\in I}\mcM(S_i)\to \bigoplus_{j\in J}\mcM(T_j)\]
         such that $E$ is the cokernel of $f$.
    \end{enumerate}
\end{theorem}
\begin{proof}
    First, we note that the subcategory of complete objects is stable under limits.
    For this it suffices to see that for any diagram $E_i\sub \cMod$ and $S\in \extr$, one obtains
    \[\hom(S,\varprojlim E_i)=\varprojlim\hom(S,E_i)\simeq \varprojlim \hom(\mcM(S),E_i)=\hom(\mcM(S),\varprojlim E_i).\]
    The stability under colimits will follow by the description via cokernels, as colimits commute with colimits.
    First, note that using the map $\bigoplus \mcM(T_i)\to 0$,
    we can describe every direct sum of objects of the form $\mcM(T)$ as kernels of maps between such direct sums,
    and hence the condition of analyticity of the ring implies that $\bigoplus\mcM(S_i)$ are complete with respect to $\iRHom$
    (and hence especially $\RHom(-,-)=\iRHom(-,-)(\ast)$ and $\hom=H^{0}\RHom$).

    Next, we will show that $\cMod$ is closed under cokernels.

    For this consider any cokernel $Q$ of a map $X\to Y$ with complete $X$ and $Y$.
    As $\RMod$ has enough compact projectives, one has a projection $\bigoplus R[T_i]\twoheadrightarrow Y$ for extremally disconnected $T_i$.
    Every coordinate $R[T_i]\to Y$ by completeness of $Y$ lifts to $\mcM(T_i)\to E$

\begin{center}\begin{tikzcd}
	{\mcM(T_i)} & Y \\
	{R[T_i]}
	\arrow[from=1-1, to=1-2]
	\arrow[from=2-1, to=1-1]
	\arrow[from=2-1, to=1-2]
\end{tikzcd}.\end{center}
    Gluing these together we obtain

\begin{center}\begin{tikzcd}
	{\bigoplus\mcM(T_i)} & Y \\
	{\bigoplus R[T_i]}
	\arrow[from=1-1, to=1-2]
	\arrow[from=2-1, to=1-1]
	\arrow[from=2-1, to=1-2]
\end{tikzcd}, \end{center}
    and noting that $\bigoplus R[T_i]\to Y$ is epic, $\bigoplus \mcM(T_{i})\to Y$ is epic as well.

    Pullback with this epimorphism does not change the cokernel, as the category is abelian, see \ref{thm:properties_of_abelian_categories}.
\begin{center}\begin{tikzcd}
	{X\times_Y\bigoplus\mcM(T_i)} & {\bigoplus\mcM(T_i)} & Q & 0 \\
	X & Y & Q & 0
	\arrow[from=1-1, to=1-2]
	\arrow[from=1-1, to=2-1]
	\arrow[two heads, from=1-2, to=1-3]
	\arrow[from=1-2, to=2-2]
	\arrow[from=1-3, to=1-4]
	\arrow[from=2-1, to=2-2]
	\arrow[two heads, from=2-2, to=2-3]
	\arrow[from=2-3, to=2-4]
\end{tikzcd}\end{center}
    Renaming $X\times_Y\bigoplus \mcM(T_i)$ as $X$, and noting that this is still contained in $\cMod$, as by assumption all direct sums of $\mcM(T_i)$ are complete, and the limit remains complete,
    we can without loss of generality assume that $Y=\bigoplus \M(T_i)$.
    Next, we surject $\bigoplus\mcM(S_j)\twoheadrightarrow X$.
    As precomposition with epimorphisms does not change cokernels, again we can replace $f$ by $\bigoplus \mcM(S_j)\twoheadrightarrow X\to Y$
    and thus further assume that $X$ is a direct sum of $\mcM[S_j]$.

\begin{center}\begin{tikzcd}
	{\bigoplus\mcM(S_j)} & X & {\bigoplus\mcM(T_i)} & Q & 0
	\arrow[two heads, from=1-1, to=1-2]
	\arrow[from=1-2, to=1-3]
	\arrow[two heads, from=1-3, to=1-4]
	\arrow[from=1-4, to=1-5]
\end{tikzcd}\end{center}
Defining $K$ as the kernel of this composition, we obtain an exact sequence

\begin{center}\begin{tikzcd}
	0 & K & {\bigoplus\mcM(S_j)} & {\bigoplus\mcM(T_i)} & Q & 0
	\arrow[from=1-1, to=1-2]
	\arrow["\phi", from=1-2, to=1-3]
	\arrow[from=1-3, to=1-4]
	\arrow["\psi", from=1-4, to=1-5]
	\arrow[from=1-5, to=1-6]
\end{tikzcd}\end{center}
    But now, as $K$ is a kernel of a map between direct sums of objects $\mcM(S)\to\mcM(T)$,
    by definition of an analytic ring, it is $\iRHom$-complete (and hence in particular $\RHom$-complete), i.e., we obtain
    \[\RHom(R[S],K[0])\simeq \RHom(\mcM(S),K[0]).\]

    \noindent Now we decompose the exact sequence into short exact sequences, by introducing $W=\coker(\phi)=\ker(\psi)$.

\begin{center}\begin{tikzcd}
	0 & K & {\bigoplus\mcM(S_j)} & W & 0 \\
	0 & W & {\bigoplus\mcM(T_i)} & Q & 0
	\arrow[from=1-1, to=1-2]
	\arrow[from=1-2, to=1-3]
	\arrow[from=1-3, to=1-4]
	\arrow[from=1-4, to=1-5]
	\arrow[from=2-1, to=2-2]
	\arrow[from=2-2, to=2-3]
	\arrow[from=2-3, to=2-4]
	\arrow[from=2-4, to=2-5]
\end{tikzcd}\end{center}
The upper row induces for any $S\in \extr$ an exact triangle

\begin{adjustbox}{max width=\textwidth}
\begin{tikzcd}
	{\RHom(R[S],K)} & {\RHom(R[S], \bigoplus \mcM(S_j))} & {\RHom(R[S],W)} & {\RHom(R[S],K)[1]}
	\arrow[from=1-1, to=1-2]
	\arrow[from=1-2, to=1-3]
	\arrow[from=1-3, to=1-4]
\end{tikzcd}
\end{adjustbox}

   \noindent Using the isomorphisms $\RHom(R[S],K)\simeq \RHom(\mcM[S],K)$ as well as 
    \[\RHom(R[S],\bigoplus \mcM(S_j))\simeq \RHom(\mcM[S],\bigoplus \mcM(S_j)),\]
    we see that $\RHom(R[S],W)=\RHom(\mcM(S),W)$; one elementary way to do so would be to spell out the long exact sequence of homologies,

\begin{adjustbox}{max width=\textwidth}\begin{tikzcd}
	{H^i(R[S],K)} & {H^i(R[S], \bigoplus \mcM(S_j))} & {H^i(R[S],W)} & {H^{i-1}(R[S],K)} & {H^{i-1}(R[S],\bigoplus\mcM(S_j))} \\
	{H^i(\mcM(S),K)} & {H^i(\mcM(S), \bigoplus \mcM(S_j))} & {H^i(\mcM(S),W)} & {H^{i-1}(\mcM(S),K)} & {H^{i-1}(\mcM(S), \bigoplus \mcM(S_j))}
\arrow[from=1-1, to=1-2]
	\arrow[from=1-2, to=1-3]
	\arrow[from=1-3, to=1-4]
	\arrow[from=1-4, to=1-5]
	\arrow["\simeq", from=2-1, to=1-1]
	\arrow[from=2-1, to=2-2]
	\arrow["\simeq", from=2-2, to=1-2]
	\arrow[from=2-2, to=2-3]
	\arrow[from=2-3, to=1-3]
	\arrow[from=2-3, to=2-4]
	\arrow["\simeq"', from=2-4, to=1-4]
	\arrow[from=2-4, to=2-5]
	\arrow["\simeq"', from=2-5, to=1-5]
\end{tikzcd}\end{adjustbox}
    and now apply the $5$-lemma.

    Now we look at the short exact sequence

\begin{center}\begin{tikzcd}
	0 & W & {\bigoplus\mcM(T_i)} & Q & 0
	\arrow[from=1-1, to=1-2]
	\arrow[from=1-2, to=1-3]
	\arrow[from=1-3, to=1-4]
	\arrow[from=1-4, to=1-5]
\end{tikzcd}\end{center}
and the corresponding long exact sequence of $\RHom$'s,

\begin{adjustbox}{max width=\textwidth}\begin{tikzcd}
	{H^i(R[S],W)} & {H^i(R[S], \bigoplus \mcM(T_i))} & {H^i(R[S],Q)} & {H^{i-1}(R[S],W)} & {H^{i-1}(R[S],\bigoplus\mcM(T_i))} \\
	{H^i(\mcM(S),W)} & {H^i(\mcM(S), \bigoplus \mcM(S_j))} & {H^i(\mcM(S),Q)} & {H^{i-1}(\mcM(S),W)} & {H^{i-1}(\mcM(S), \bigoplus \mcM(T_i))}
\arrow[from=1-1, to=1-2]
	\arrow[from=1-2, to=1-3]
	\arrow[from=1-3, to=1-4]
	\arrow[from=1-4, to=1-5]
	\arrow["\simeq", from=2-1, to=1-1]
	\arrow[from=2-1, to=2-2]
	\arrow["\simeq", from=2-2, to=1-2]
	\arrow[from=2-2, to=2-3]
	\arrow[from=2-3, to=1-3]
	\arrow[from=2-3, to=2-4]
	\arrow["\simeq"', from=2-4, to=1-4]
	\arrow[from=2-4, to=2-5]
	\arrow["\simeq"', from=2-5, to=1-5]
\end{tikzcd}\end{adjustbox}
    which using the $\RHom$-completeness of $Q$ and of $\bigoplus \mcM(T_i)$, again by the $5$-lemma, implies $\RHom(R[S], Q)\simeq \RHom(\mcM(S), Q)$.
    Hence the category $\cMod$ is stable under cokernels.
Note that the above argument can also be made with $\iRHom$, yielding that for such a cokernel automatically $\iRHom(R[S],Q)\simeq \iRHom(\mcM(S),Q)$ holds.

    Now, as the (full sub-) category is stable under cokernels, we may compute the cokernels in the category $\RMod$ and they agree with the cokernels in $\cMod$.
    This in particular implies that for any complete $Q$, by surjecting with $\bigoplus R[T_i]\twoheadrightarrow Q$,
    and using that every epimorphism is its image, i.e., the cokernel of its kernel (we are in an abelian category), and afterwards surjecting onto the kernel, $Q$ is the cokernel of an arrow

\begin{center}\begin{tikzcd}
	{\bigoplus \mcM(S_j)} & {\bigoplus \mcM(T_i)} & Q & 0
	\arrow[from=1-1, to=1-2]
	\arrow[from=1-2, to=1-3]
	\arrow[from=1-3, to=1-4]
\end{tikzcd}\end{center}
This shows the equivalence of the first and third description of $\cMod$.
As we have seen that all such cokernels are in fact $\iRHom$-complete, this implies the equivalence of all three descriptions of $\cMod$.
Now, the stability under all colimits follows, as for any direct sum $\bigoplus_\ell Q^{\ell}$ with $Q^\ell\in \cMod$ with corresponding representations as cokernel via

\begin{center}\begin{tikzcd}
	{\bigoplus_{j\in J^{\ell}} \mcM(S_j^{\ell})} & {\bigoplus_{i\in I^{\ell}} \mcM(T_i^\ell)} & {Q^\ell} & 0
	\arrow[from=1-1, to=1-2]
	\arrow[from=1-2, to=1-3]
	\arrow[from=1-3, to=1-4]
\end{tikzcd}\end{center}
    We can just glue all these arrows together to obtain $\bigoplus Q^{\ell}$ as a cokernel,

\begin{center}\begin{tikzcd}
	{\bigoplus_\ell\bigoplus_{j\in J^{\ell}} \mcM(S_j^{\ell})} & {\bigoplus_\ell\bigoplus_{i\in I^{\ell}} \mcM(T_i^\ell)} & {\bigoplus_\ell Q^\ell} & 0.
	\arrow[from=1-1, to=1-2]
	\arrow[from=1-2, to=1-3]
	\arrow[from=1-3, to=1-4]
\end{tikzcd}\end{center}
    As colimits commute with colimits (in the language of abelian categories: coproducts are right exact),
    this implies the stability under all coproducts and thus together with cokernels under all colimits.

    This in turn answers the question of compact projectivity of $\mcM(S)$, as $\hom(R[S],-)$ commutes with all colimits in $\RMod$ (which now agree with the ones in $\cMod$),
    and $\delta$ induces a natural isomorphism
\[\hom(\mcM(S),\varinjlim W_i)\simeq \hom(R[S],\varinjlim W_i)\simeq \varinjlim \hom(R[S], W_i)\simeq \varinjlim\hom(\mcM(S),W_i).\]
    By the description of $\cMod$, the class $\mcM(S),\, S\in \extr$, is generating, hence in conclusion we have enough compact projectives.

The stability under extensions follows again by the five lemma, as every object in $\cMod$ by the description above is $\RHom$-complete.
    Explicitly, let $A$, $C$ be complete objects, hence by the description above $\RHom(R[S],A)\simeq \Rhom(\mcM(S),A)$ and $\RHom(R[S],C)\simeq \RHom(\mcM(S),C)$ we obtain for any short exact sequence

\begin{center}\begin{tikzcd}
	0 & A & B & C & 0
	\arrow[from=1-1, to=1-2]
	\arrow[from=1-2, to=1-3]
	\arrow[from=1-3, to=1-4]
	\arrow[from=1-4, to=1-5]
\end{tikzcd}\end{center}
    that $B$ is complete, by

\begin{adjustbox}{max width = \textwidth}\begin{tikzcd}
	{H^{i+1}(R[S],C)} & {H^i(R[S],A)} & {H^i(R[S],B)} & {H^i(R[S], C)} & {H^{i-1}(R[S],A)} \\
	{H^{i+1}(\mcM(S),C)} & {H^i(\mcM(S),A)} & {H^i(\mcM(S),B)} & {H^{i}(\mcM(S),C)} & {H^{i-1}(\mcM(S),A).}
	\arrow[from=1-1, to=1-2]
	\arrow[from=1-2, to=1-3]
	\arrow[from=1-3, to=1-4]
	\arrow[from=1-4, to=1-5]
	\arrow["\simeq", from=2-1, to=1-1]
	\arrow[from=2-1, to=2-2]
	\arrow["\simeq", from=2-2, to=1-2]
	\arrow[from=2-2, to=2-3]
	\arrow[from=2-3, to=1-3]
	\arrow[from=2-3, to=2-4]
	\arrow["\simeq"', from=2-4, to=1-4]
	\arrow[from=2-4, to=2-5]
	\arrow["\simeq"', from=2-5, to=1-5]
\end{tikzcd}\end{adjustbox}

    The analogous argument for $\iRHom$ instead of $\Rhom$ shows stability under internal $\Ext$'s
    (so in particular, stable under $\ihom=H^{0}\Ext$).

    It only remains to show the existence of the completion and the existence of the completed tensor product.%
    \footnote{Probably there is an adjoint functor theorem specialised to situations like the present, but we stick to the elementary setting.}
    For the completion, use that any $A\in \RMod$ can be written as the cokernel of a map

\begin{center}\begin{tikzcd}
	{\bigoplus R[S_j]} & {\bigoplus R[T_i]} & A & 0
	\arrow["\phi",from=1-1, to=1-2]
	\arrow[from=1-2, to=1-3]
	\arrow[from=1-3, to=1-4]
\end{tikzcd}.\end{center}
    Now define the completion of $A$, $\mcM(A)$, to be the cokernel

\begin{center}\begin{tikzcd}
	{\bigoplus \mcM(S_j)} & {\bigoplus \mcM(T_i)} & {\mcM(A)} & 0
	\arrow["\psi", from=1-1, to=1-2]
	\arrow[from=1-2, to=1-3]
	\arrow[from=1-3, to=1-4]
\end{tikzcd}\end{center}

Now
    \[\hom(\mcM(A),?B)=\hom(\coker(\phi), ?B)=\ker(\phi^*\colon \hom(\bigoplus R[T_i],B)\to \hom(\bigoplus R[S_j],B))\]
    which, by completeness of $B$ (clearly holding also for $\bigoplus \mcM(T_i)$ by continuity of the $\hom$ in the first coordinate)
    agrees with
    \[\ker(\psi^*\colon \hom(\bigoplus \mcM(T_i),B)\to \hom(\bigoplus\mcM(S_j),B))=\hom(\coker(\psi), B)=\hom(\mcM(A),B).\]
    To see that this left adjoint sends $R[S]$ to $\mcM(S)$ simply use

\begin{center}\begin{tikzcd}
	0 & {R[S]} & {R[S]} & 0
	\arrow["\phi", from=1-1, to=1-2]
	\arrow[from=1-2, to=1-3]
	\arrow[from=1-3, to=1-4]
\end{tikzcd}\end{center}

    The tensor product now clearly is given by using the tensor product in $\RMod$ and completing,
    \[A\otimes_{(R,\mcM)}B=\mcM(A\otimes B).\]

    It is clearly is colimit preserving, and hence unique as every element is a colimit of $\mcM(T)$'s
    (and as left adjoints are fully determined by their values on objects).
    It fulfills the tensor-$\hom$ adjunction (with respect to the same internal $\hom$) as for any complete $A,B,C\in \cMod$,
    \[\hom(A\otimes_{(R,\mcM)} B,C)=\hom(A\otimes B, ?C)=\hom(A,\ihom(B,C)).\]

    For any $A\in \RMod$, $B\in \cMod$ and any $C\in \cMod$, we compute (as the internal $\hom$'s agree by stability under internal $\hom$'s)
	\begin{align*}
      &\quad\hom_{\cMod}(\mcM(A\otimes B),C)\\
      &=\hom_{\RMod}(A\otimes B, C)\simeq \hom_{\RMod}(A, \ihom_{\RMod}(B,C))\\
      &\simeq \hom_{\RMod}(A,?\ihom_{\cMod}(B,C)) \simeq\hom_{\cMod}(\mcM(A), \ihom(B,C))\\
      &= \hom_{\cMod}(\mcM(A)\otimes_{(R,\mcM)} B, C).
	\end{align*}
    By Yoneda, this implies $\mcM(A\otimes_{\RMod} B)=\mcM(A)\otimes_{(R,\mcM)} B$,
    and by symmetry it follows that
        \[\mcM(A)\otimes_{(R,\mcM)} \mcM(B)=\mcM(\mcM(A)\otimes \mcM(B))=\mcM(A\otimes B),\]
     This is the last missing statement that completion is symmetric monoidal.
\end{proof}

\begin{remark}
  As the three forgetful functor $\cMod\to\RMod$, $\RMod\to\cab$ and $\cab\to\cond$ have left adjoints,
  we have, on any condensed set $X$, a free complete $R$-module, $\mcM(R[X])$.
\end{remark}

This shows that the category of complete $R$-modules is essentially as good as $R$-modules itself, and one obtains an extremely well-behaved theory of completions.
In particular, the unique completed tensor product and that the category is abelian (even with particular good properties as such)
clearly must be surprising to a functional analyst, and we suggest to be impressed.

This statement also survives the passage to derived categories, see again \cite[7.5]{scholze2019condensed} or \cite[6.14]{scholze2019Analytic}.
\begin{theorem}\label{thm:der-an-rings-gut}
  For any analytic ring $(R,\mcC)$,
    the full subcategory of complete derived $R$-modules, $\cDRMod$, fulfills the following conditions.
    \begin{enumerate}
        \item It is stable under all limits and colimits, and $\mcM(S)[i]$ for $i\in \Z$ and $S\in \extr$ forms a generating class of compact projective objects.
      \item There exists a completion $X\mapsto X\otimes_{(R,\mcM)}^L (R,\mcM)$, being the left adjoint to the inclusion,
            which coincides with the left derived functor of the completion of $R$-modules.
        \item If $R$ is commutative, there exists a unique symmetric monoidal tensor product $\otimes^L_{(R,\mcM)}$ such that completion is symmetric monoidal.
        \item It is closed under $\iRHom$'s.
\end{enumerate}
    Furthermore, $\cDRMod$ can be described in the following equivalent ways.
    \begin{enumerate}
        \item The category of objects $E\in \DRMod$ such that $\delta$ induces an isomorphism
        \[\RHom(R[S],E)\simeq \RHom(\mcM(S),E).\]
        \item The category of objects $E\in \DRMod$ such that $\delta$ induces an isomorphism
        \[\iRHom(R[S],E)\simeq \iRHom(\mcM(S),E).\]
        \item Those objects $E\in \DRMod$ such that all $H^i(E)\in \cMod$.
        \item The derived category $D(\cMod)$.
\end{enumerate}
\end{theorem}
\begin{proof}
    This is mainly \cite[5.9, 7.5]{scholze2019condensed}.
    We only sketch the proof that $D(\cMod)\sub D(\RMod)$ is fully faithful.
    Note that $\cMod$ admits enough compact projectives, as well as $\RMod$, and $\RHom$ is computed via projective resolutions.
    Thus, it suffices to show that for any $C\in D(\cMod)$
    \[\RHom_{D(\cMod)}(\mcM(X), C)\simeq \RHom_{D(\RMod)}(\mcM(X),C)\simeq \RHom(R[X],C).\]
    By an approximation argument with Postnikov truncations, see \cite[5.9]{scholze2019condensed},
    it suffices to show that for all $Y\in \cMod$ and $i\ge 0$,
    \[\Ext^i_{\cMod}(\mcM(X),Y)\simeq \Ext^i_{\RMod}(R[X],Y).\]
    As $R[X]$ and $\mcM(X)$ are projective in $\RMod$ resp. $\cMod$, both sides are $0$ for $i>0$, and as $Y\in \cMod$, the equality for $i=0$ follows (note that this argument is essentially also contained in the proof of the analogous 1-categorical statement above.)

    For the rest of the equivalent characterisations and for the derived completions, see \cite[5.9]{scholze2019condensed}.
    These characterisations also imply the stability under all limits and colimits.
    For the derived tensor product and the internal $\RHom$ equality, we refer to \cite[7.5]{scholze2019condensed}.
\end{proof}
Note that the completed derived tensor product does not necessarily have to be the left derived functor of the completed tensor product on $\cMod$,
although in practice, it always is.
We will investigate a criterion for this, which is \cite[7.6]{scholze2019condensed}.
\begin{lemma}\label{lem:cond-lder-tp}
    The derived completed tensor product is given by the left derived functor of the completed tensor product precisely if the following condition holds.%
    \footnote{Note that we are not entirely sure if we understood this properly, see warning 7.6 in \cite{scholze2019condensed} for the original statement.}
    For all $S,T\in \extr$, compute the complex $\mcM^L(S\times T)$ as follows.
    We compute $\mcM^L(S\times T)$ by taking a projective (simplicial) resolution $P_i$ of $S\times T$,
    inducing via Dold-Kan a projective resolution of $\Z[S\times T]$

\begin{center}\begin{tikzcd}
	\dots & {\Z[P_2]} & {\Z[P_1]} & {\Z[P_0]} & {\Z[S\times T]} & 0
	\arrow[from=1-1, to=1-2]
	\arrow[from=1-2, to=1-3]
	\arrow[from=1-3, to=1-4]
	\arrow[from=1-4, to=1-5]
	\arrow[from=1-5, to=1-6]
\end{tikzcd}\end{center}
    and applying $\mcM$ to each degree, i.e., $\mcM^L(S\times T)$ is given by

\begin{center}\begin{tikzcd}
	\dots & {\mcM(P_2)} & {\mcM(P_1)} & {\mcM(P_0)} & 0
	\arrow[from=1-1, to=1-2]
	\arrow[from=1-2, to=1-3]
	\arrow[from=1-3, to=1-4]
	\arrow[from=1-4, to=1-5]
\end{tikzcd}\end{center}

    Now the derived tensor product agrees with the left derived functor of the tensor product precisely if all these complexes are concentrated in degree $0$.
\end{lemma}
\begin{proof}
See Warning 7.6 in \cite{scholze2019condensed}.

    Our understanding is that we just need to show that the left derived functor of the tensor product makes the completion symmetric monoidal.
    This reduces to testing against the generators $R[S]$, i.e., showing that the derived completion of $R[S\times T]$ is given by $\mcM(S)\otimes^L \mcM(T)$.
    But since $\mcM(S)$ and $\mcM(T)$ are compact projective, $\mcM(S)\otimes^L\mcM(T)=(\mcM(S)\otimes \mcM(T))[0]$.
    Computing the derived completion of $R[S\times T]$ can be done by choosing a projective resolution as in the statement of the lemma,
    and showing that it is concentrated in degree $0$ precisely would mean that $\mcM(S\times T)\simeq H^0(\mcM(P_\bullet))[0]=\mcM(S\times T)[0]$.
    Now, since the tensor product makes completion symmetric monoidal, these two agree.
\end{proof}

As a last topic, we will discuss a helpful lemma, which essentially states that considering complete condensed abelian groups is equivalent to considering the resp.\ complete $R$-modules.
One should try to get a feel for why this is totally reasonable.
\begin{lemma}
    For any analytic ring $(R,\mcM, \delta)$, the object $\mcM(\ast)$ is a condensed ring and the $R$-modules $\mcM(S)$ are $\mcM(\ast)$-modules, turning
    $(\mcM(\ast),\mcM,\delta)$ into an analytic ring.
    The forgetful functor induces equivalences $\cMod\simeq \mcM(\ast)\mathrm{-Mod}_{\mcM}$ and
    $\cDRMod\simeq D(\mcM(\ast)\mathrm{-Mod})_{\mcM}.$
    \end{lemma}
\begin{proof}
This is \cite[6.15]{scholze2019Analytic}.
\end{proof}
\begin{corollary}
    The categories of complete $R$-modules can equivalently be described as the $\mcM$-complete abelian groups.
I.e., condensed abelian groups $E$ such that
\[\hom_{\cab}(\Z[S],E)\simeq \hom_{\cab}(\mcM(S),E)\]
are automatically complete condensed $\mcM(\ast)$-modules.
Analogously, a derived complete $\mcM(\ast)$-module is equivalently a complete derived condensed abelian group.
    \end{corollary}
This in particular shows, that the most natural choice of the base ring is given by $\mcM(\ast)$,
and that furthermore the definition of completeness can fully be reduced to working in abelian groups.
\begin{definition}
    A \idx{normalized analytic ring} is an analytic ring $(R,\mcM)$ for which $\delta_\ast\colon R=R[\ast]\to \mcM(\ast)$ is an isomorphism.
\end{definition}

\subsection{Animated analytic rings}
Although analytic rings posses excellent properties, the definition is at heart 1-categorical and in particular the base ring is a non-derived object.
This defect induces problems with general base change, as well as the many subtle switches between the category and the derived category.

As we want to develop a theory that, e.g., contains \idx{relative functional analysis} in the form of base changing the base ring,
we need to allow everything to be equally derived.

To develop the full strength of analytic rings, we have to switch from using the \enquote{wrong} derived category
to the $\infty$-categorical versions.
This makes not only the definitions easier and proofs more structured, but also allows for stronger results.
Recall that for any animated condensed ring $R$, we can define the $\infty$-derived category $\mcD_{\ge0}(R)$ of $R$-modules in condensed animated abelian groups.

\begin{definition}
An \idx{animated analytic ring} is an animated ring $R$ in $\Ani(\cond(\Ab))=\mcD_{\ge 0}(\cab)$ together with a functor
    $\mcM\colon \extr\to \mcD_{\ge 0}(R)$ taking finite coproducts to finite direct sums and a natural transformation $\delta_S\colon S\to \mcM(S)$.
    We further demand that for any sifted colimit $X$\footnote{Note that sifted colimits in $\infty$-categories are generated by geometric realisations and filtered colimits.}
    in $\mcD_{\ge 0}(R)$ of objects of the form $M(S)$,
    $\delta$ induces an isomorphism
\[\ihom_{\mcD_{\ge 0}(R)}(R[S],X)\simeq \ihom_{\mcD_{\ge 0}(R)}(\mcM(S),X)\]
\end{definition}
\begin{remark}
    Note that this version seems to be generalisable to more general categories than just animated rings, e.g., to animations of condensed algebraic categories.
\end{remark}
\begin{remark}
  In \cite[12.5]{scholze2019Analytic} it is remarked, that by passing to spectrum objects one recovers similar results for the unbounded derived $\infty$-category.
In particular, one could replace the $\ihom_{\mcD_{\ge0}(R)}$-isomorphism with an unbounded $\iRHom$-isomorphism.
We guess that this follows by writing the unbounded complex as colimit of its left postnikov truncations,
and commuting the colimit through the $\iRHom$ (as this is a right adjoint, and starts in the opposite category), and further shifting any individual truncated complex far enough to the right.

However, we do not yet understand enough of the theory to explain how this works precisely.
\end{remark}

Now, we have an analog of complete objects.
\begin{definition}
    An object $E\in \mcD_{\ge 0}(R)$ is said to be complete with respect to an animated analytic ring $(R,\mcM, \delta)$, if $\delta$ induces an isomorphism
    \[\ihom_{\mcD_{\ge 0}(R)}(R[S],E)\simeq \ihom_{\mcD_{\ge 0}(R)}(\mcM(S),E).\]
The full subcategory of the complete objects is denoted by $\mcD_{\ge 0}(R,\mcM)\sub \mcD_{\ge 0}(R)$.
\end{definition}
This subcategory again possesses excellent properties.

\begin{proposition}
For an animated analytic ring $(R,\mcM,\delta)$, the full subcategory $\mcD_{\ge 0}(R,\mcM)$ has the following properties.
\begin{enumerate}[(i)]
        \item It is stable under all limits and colimits, and furthermore it is generated under sifted colimits by the compact projective objects $M(S)$ for $S\in \extr$.
    \item There exists a completion, forming the left adjoint to the inclusion, and extending the map $R[S]\mapsto \mcM(S)$.
       \item If $R$ is a condensed animated commutative ring, then $\mcD_{\ge_0}(R)$ is symmetric monoidal, and there exists a unique symmetric monoidal structure on $\mcD_{\ge 0}(R,\mcM)$ such that completion is symmetric monoidal.
\end{enumerate}
    Furthermore we can characterize the complete objects as follows.
    The category $\mcD_{\ge 0}(R,\mcM)$ is prestable, and its heart is given by colimits of those $\pi_0 R$-modules that are given by $\pi_0 \mcM(S)$ for $S\in \extr$.
    The objects $X$ in $\mcD_{\ge 0}(R)$ are complete precisely if all $H_i(X)$ are contained in the heart.

\end{proposition}
\begin{proof}
    This is \cite[12.4]{scholze2019Analytic}.
\end{proof}
One advantage of the generalisation to the $\infty$-categorical setting is the surprising fact that these properties of the subcategory essentially classify analytic rings.
This in particular shows that the definition of analytic rings, which a priori might have seemed quite arbitrary,
in fact could not have been different if one wants to recover these good properties.
\begin{theorem}
    Consider a condensed animated ring $R$ and a full subcategory $\mcD\sub \mcD_{\ge 0}(R)$.
    Then there exists a (necessarily unique) animated analytic ring structure $(R,\mcM)$ on $R$ with $\mcD_{\ge 0}(R,\mcM)=\mcD$ precisely if all of the following hold.
    \begin{enumerate}[(i)]
        \item The subcategory $\mcD$ is stable under all limits and colimits.
        \item The subcategory $\mcD$ is stable under $\ihom_{\mcD_{\ge 0}(\cab)}(\Z[S],-)$ for $S\in \extr$.
        \item The inclusion $\mcD\sub \mcD_{\ge 0}(R)$ admits a left adjoint
    \end{enumerate}
    \end{theorem}
\begin{proof}
  This is \cite[12.20]{scholze2019Analytic}.
\end{proof}

The next powerful feature of animated analytic rings, which in general fails for 1-categorical condensed analytic rings, is base change.
\begin{proposition}
    Let $(R,\mcM)$ be an animated analytic ring, and $C$ a condensed animated ring, together with a map $R\to C$ of condensed animated rings.
Then we can define an induced analytic ring structure $\mcN$ on $C$ by defining
    \[S\mapsto \mcN(S)=C[S]\otimes_{(R,M)}(R,\mcM),\]
    i.e., the space of $C$-measures on $S$ is the completion of the free $C$-module with respect to the completion induced by $(R,\mcM)$.
    This forms an analytic ring.
\end{proposition}
\begin{proof} This is \cite[12.8]{scholze2019Analytic}.
\end{proof}

A sample consequence of this result is the possibility to complexify, i.e., having a notion of completeness of $\R$-modules,
we have an induced notion of completeness of complex vector spaces, and a complexification.
We are convinced that this result will be of huge relevance to relative functional analysis in form of, e.g., Kaplanski-Banach modules,
or $L^0$-modules, see, e.g., \cite{N.Edeko2024, HenrikKreidler2021, Edeko2023, Jamneshan2020, jamneshan2022uncountable, Jamneshan_2018, Hermle2022Halmos}
for the type of theory we believe can be carried out fruitfully in this setting.

Lastly, this also implies that, similar to the 1-categorical setting, changing the base ring to its correct choice does not affect the complete objects, and one can work with normalised analytic rings.
See \cite[12.9/26]{scholze2019Analytic} for more.
\begin{definition}
    There is an equivalence of categories of complete objects with respect to $(R,\mcM)$ and with respect to the analytic ring $(\mcM(\ast),\mcM)$.
    An animated analytic ring is normalized if $R\mapsto \mcM(\ast)$ is an isomorphism.
\end{definition}


\section{Embedding functional analysis}

    \quot{[Algebraic geometry] seems to have acquired
the reputation of being esoteric, exclusive, and very abstract, with adherents who are secretly
plotting to take over all the rest of mathematics. In one respect this last point is accurate.}
{Mumford in \enquote{Curves and their Jacobians}}

    \quot{For some reason, this secret plot has so far stopped short of taking over analysis. The goal
of this course is to launch a new attack, turning functional analysis into a branch of commutative
algebra [...]}
{Scholze in \cite{scholze2019Analytic}, referring to Mumford}

Having seen the general theory of analytic rings, in this section we will find the correct replacement of locally convex complete $\R$-vector spaces.
This has some surprising twists to it, i.e., the naive version, which would essentially lead just to the old locally convex theory,
will be replaced by a beautiful different version, leading to a theory stronger than the classical one.

For classical general functional analysis we recommend \cite{charalambos2013infinite}, regarding ergodic theory \cite{Eisner2016}, for locally convex spaces \cite{schaefer1971Topological} and for Banach lattice theory \cite{schaefer1974Banachlattices}.
Some combination of functional analysis and category theory can be found, e.g., in \cite{sanchez2023homological, castillo2021hitchhiker, Castillo2010, Aviles2016, hassoun2020examplesnonexamplesintegralcategories, Semadeni, Jamneshan2020, Frederick2010, Hoffmann1977}
\subsection{Towards liquid vector spaces}

First, we will categorify the notion of Radon measures.
For this, it is helpful to look at two \enquote{extreme cases} of analytic rings, induced by the two extreme cases of completeness:
trivial and $p$-adic completeness.

\begin{example}
    \begin{itemize}
        \item The trivial structure $S\mapsto \Z[S]=\bigcup_n \varprojlim_i \Z[S_i]_{\le n}$ induces an analytic ring.
        \item The assignment \[S\mapsto \varprojlim \Z[S_i]\] for $S=\varprojlim S_i$ induces the analytic ring of \idx{solid abelian groups}.
    \end{itemize}
\end{example}

Going back to our functional analytic setup, we need to understand Radon measures.
Since on totally disconnected spaces, definitions always become easier, the definition of Radon measures becomes categorifiable.

\begin{definition}
For any profinite set $S=\varprojlim S_i$ the space $\mathcal{M}_1(S)$ of (signed) Radon measures on $S$ equipped with compact open topology agrees with
\begin{align*}
\mathcal{M}_1(S)&=\bigcup_{c>0}\varprojlim_i\mathcal{M}(S)_{\le c}\\
&=\ihom(C(S,\R),\R)\\
\end{align*}
where $\mathcal{M}(S)_{\le c}\coloneqq \varprojlim_i \R[S_i]_{\ell^1\le c}$ denotes the limit over all subsets $\R[S_i]_{\ell^1\le c}$ of $\ell^1$ norm at most $c$.

Furthermore, the classical definition of Radon measures simplifies to being maps $\mu\colon \mathrm{clopen}(S)\to \R$ such that the following hold:
\begin{enumerate}[(i)]
    \item $\mu(\emptyset)=0$,
    \item finite additivity, i.e., $\mu(U_1\sqcup U_2)=\mu(U_1)+\mu(U_2)$ for disjoint clopens $U_1$ and $U_2$,
    \item bounded 1-variation, i.e., there exists $c<\infty$ such that for any finite family of disjoint clopens $(U_j)$,
     \[\sum_j |\mu(U_j)|\le c.\]
\end{enumerate}
More generally, for any compact Hausdorff space $X$ given as $S/S\times_X S$ for a profinite set $S$, the set of signed Radon measures $\mathcal{M}(X)$ on $X$ agrees with
the coequalizer of the arrows $\mathcal{M}(S\times_X S)\to \mathcal{M}(S)$.
\end{definition}

\begin{proof}
This is \cite[3.3]{scholze2019Analytic} or \cite[3.9]{scholze2022complex}.
\end{proof}
\begin{remark}
    One would a priori expect the Radon measures to be the objects $\R[S]$, by analogy to~\ref{lem:expl-zs}.
    This, however, is not equivalent; one has
    \[\R[S]=\bigcup_{n\ge 0}\varprojlim_i \R[S_i]_{\ell^0\le n},\]
    where the $\ell^0$ simply counts the number of nonzero coefficients.
    See \cite[3.10]{scholze2022complex}.
    \end{remark}

The naive guess of defining an analytic ring structure for functional analysis now clearly would be the choice of radon measures.
\begin{definition}
    Define the pre-analytic ring of Radon measures $(\R,\mcM_1, \delta)$ by
    \[\mcM_1(S)=\bigcup_{c>0}\varprojlim R[S_i]_{\ell^1\le c}\]
    and $\delta\colon S\to \mcM_1(S),\, x\mapsto 1[x]$.
\end{definition}
Now clearly, any complete locally convex vector space is $\mcM_{1}$-complete by \ref{lem:locally_convex_radon_comp}.
Interestingly, this does not define an analytic ring, basically due to the following problem with entropy, called \idx{Ribe's extension},
see also \cite{Commelin}.

Consider the map $\ell^1\to \ell^2$, induced by the continuous map $S((x_n))=(x_n\log|x_n|)$, where $0\log(0)\coloneqq 0$.
    This map has the feature that it is \enquote{almost linear} in the sense that the nonlinearity
    \[S(x+y)-S(x)-S(y)\]
    is contained in $\ell^{1}$ for all $x,y\in \ell^1$.
    Hence, it induces a \emph{linear} continuous map between condensed $\R$-modules
    \[S\colon \ell^1\to \ell^2/\ell^1.\]
   Extend $S$ from $\ell^1$ to $\mcM_1(\alpha\N)=\ell^1\oplus \R$ arbitrarily.

    As $\ell^1$ and $\ell^2$ are both Banach spaces, they clearly are $\mcM_1$-complete.
    If the Radon measures formed an analytic ring, then the factor $\ell^2/\ell^1$ still would be complete.
    But now note that $S$ restricted to Radon measures is $0$, i.e., it is an extension of $0$,
    but clearly not identical to $0$ as $S$ does not always land inside $\ell^1$ (e.g., use $x_{n}=1/(n(\log n))^2$).

    This means that there are two continuous $\R$-linear extensions

\begin{center}\begin{tikzcd}
	{\mcM_1(\alpha\N)} \\
	\\
	{\alpha\N} && {\ell^2/\ell^1}
	\arrow["0"', shift right=2, from=1-1, to=3-3]
	\arrow["S", shift left=2, from=1-1, to=3-3]
	\arrow["\delta", from=3-1, to=1-1]
	\arrow["0"', from=3-1, to=3-3]
\end{tikzcd}\end{center}
    meaning that integration would not be uniquely determined in $\ell^2/\ell^1$, and thereby it is not $\mcM_1$-complete.

See, e.g., \cite{sanchez2023homological} for more on these types of extension problems; they are not avoidable in classical theory.
Extensions may always destroy \enquote{infinitesimal amounts} of convexity of the unit balls.

This crucial point in the theory can be explained roughly as follows.
(At least this seems like a good intuition to the authors.)
Although $\ell^{1}$ is perfectly reasonable with respect to Hausdorff behaviour,
it is exposed as \enquote{slightly too big} as soon as one enters the non-Hausdorff regime (think: $\ell^{2}/\ell^{1}$).
Generally, one has many possibilities of installing growth bounds on function spaces, usually in form of Orlicz spaces.
However, any such growth condition is \enquote{closed}, meaning that there is a rather sharp upper bound to the growth rate.
The sharper this upper bound is, the easier the condition of remaining in the space can be destroyed by simple operations.
In fact, here some \enquote{almost trivial} deformation in the form of adding logarithmic growth already suffices to fall out of $\ell^1$.
But as such a deformation can not be distinguished by the linear theory, we have to find something that is stable under such deformations.
This leads to the idea of using an \enquote{open} growth condition in the following form of an inductive limit.
Take a filtered colimit (a union) of spaces, each yielding an individually closed growth condition,
but where an infinitesimal deformation of something satisfying one growth condition is still captured by some slightly bigger, i.e., weaker growth condition.
And most importantly -- don't take the closure afterwards, meaning that one should not try to mix the growth condition (Orlicz functions) to one new growth condition.
Then one will obtain a space that is \enquote{stable} under infinitesimal deformations.

Let us give concrete examples of this idea, which clearly fits the topological idea of
open sets as sets where \enquote{there is always an $\epsilon$ of room} very well, explaining our intuitive use of the terms \enquote{closed} and \enquote{open}.

\begin{itemize}
  \item The first quite prominent and self-explaining example of this superiority of \enquote{open} phenomena are ordinals, large cardinals and universes.
  \item Another instance is given by the notion condensed sets themselves.
        Each category $\cond_\kappa$ is not particularly stable -- there is an absolute growth condition in form of size boundaries on $\extr_\kappa$ and thus the complexity of the topology,
        that is destroyed, e.g., when taking large products.
        But we want a category that is stable under such operations.
        Hence we find a filtered family $\cond_{\kappa}\hookrightarrow\cond_{\kappa'}\hookrightarrow \dots$, and take the colimit of those.
        Now we obtain a wonderful category that is stable under whatever we wanted;
        as for each argument one can pass to a \enquote{bounded} level, do whatever one wishes there and by these operations only lands in another (though quite possibly bigger) $\cond_{\kappa'}$.
        However, it is important not to take the \enquote{closure} of this category in the form of taking all sheaves on $\extr$,
        as otherwise one would again obtain a \enquote{closed} condition.
  \item The next analogy is with the situation of power series in one complex variable.
        If one considers the disc algebra (those power series that converge on the closure of the unit ball),
        then this is not a particularly stable growth condition, as, e.g., the derivative of such a function may well explode at the boundary of the ball.
        But at least, the derivative of such a function can only explode \enquote{linearly} on the boundary,
        so the next guess would be to take those functions that explode at most linearly on the boundary.
        But again, derivatives of those functions would tend quadratically to infinity, leading to the same problem.
        The solution is to pass to the way stabler set of those power series that have a radius of convergence at least one
        (i.e., that are holomorphic on $B_{r}(0)$ for all $r<1$) -- this condition essentially cannot be destroyed anymore.
  \item Another situation quite similar in spirit is the subexponential growth of powers of operators with spectral radius $1$,
        compared to bounding the polynomial growth of the resolvents.
\end{itemize}
Applying this to our situation, the idea is to replace $\ell^1$ (which is slightly too big) as the space of measures by the union of $\ell^p$ for all $p<1$
(we need $p<1$ in order to approximate from the inside),
\[ \varinjlim_{p<1} \ell^p.\]

This might seem strange, as it is a space sitting (very) densely in $\ell^1$.
However, it is critical to not succumb to the habit of the classical functional analytic setting and take the closure.
Taking the closure would lead to $\ell^1$ and one loses.
Instead, we remind the reader that it is often more natural for spaces of measures to be equipped with weak*- or similar topologies and thus not normable.
So using a locally convex inductive limit here does not feel too bad after all.

Note that the appearance of $\ell^p$-norms instead of other sequence spaces might seem arbitrary here at first.
We believe that this can be explained by some intuition coming from valuations, as the $p$-norms are special in the sense that they behave
way better with algebraic structures and hence are precisely the ones that capture some algebraic behaviour.
But as we do not fully understand these more subtle aspects yet,
we suggest the reader to have a look at \cite[chap.~6--10]{scholze2019Analytic}, especially the use of $\Z((T))_{>r}$ seems to give reason as to why $\ell^p$ is the correct choice, see in particular chap.~7.

In some sense, to us this fits the intuition of exhausting any unit ball from the inside with concave $L^p$ balls ($p<1$),
rather than taking the closure (which would be the $L^1$-ball).
Thus we would claim that this can essentially be summed up as the brilliant and quite funny idea to replace local convexity by \enquote{local non-concavity},
leading to a robust growth condition.

\subsection{The liquid ring}
\quot{One Ring to rule them all,\\
  One Ring to find them,\\
  One Ring to bring them all\\
  and in the darkness bind them.}{J.R.R. Tolkien, The Fellowship of the Ring}

The preceding discussion hopefully helps to motivate the following definition of the space of Radon measures with bounded $q$-variation for some $q<1$:

\begin{definition}
    For any $p\in (0,1]$ and profinite $S$, define the \idx{measures of bounded $p$-variation on $S$} to be
    \[\mcM_p(S)=\bigcup_{c>0}\varprojlim \R[S_i]_{\ell^p<c},\]
    where $\R[S_i]_{\ell^p<c}$ are those $\R$-linear combinations $\sum a_i [x_i]$ of the elements in $S_i$ with $\sum |a_i|^p \le c$.
    There are canonical inclusions $\mcM_q(S)\sub \mcM_p(S)$ for $q\le p$.
Using these, for any $p\le 1$, define
    \[\mcM_{<p}(S)=\varinjlim_{q<p} \mcM_{q}(S),\]
    which together with the Dirac embedding yields the \idx{$p$-liquid pre-analytic ring} $(\R,\mcM_{<p})$.
\end{definition}
\begin{remark} The measures in $\mcM_p$ can equivalently be characterised as the usual Radon measures of bounded $p$-variation, and described as those $\mu\colon \mathrm{clopen}(S)\to \R$ with
    \begin{itemize}
        \item $\mu(\emptyset)=0$,
        \item $\mu(A\sqcup B)=\mu(A)+\mu(B)$,
        \item and the property there exists a $c>0$ such that for any finitely many disjoint clopens $A_i$ one has
        \[\sum |\mu(A_i)|^p\le c.\]
\end{itemize}
    \end{remark}

    Of course,
    this definition immediately begs the question of whether or not this pre-analytic ring is an analytic ring.
    Indeed, it is.

\begin{theorem}\label{thm:liquid}
    The $p$-liquid structure indeed forms an analytic ring.
\end{theorem}
Instead of a proof, we give the following remark.
\begin{remark}
    This theorem is by far the deepest theorem in the entire text.
    The proof of this theorem fills the main part of \cite{scholze2019Analytic},
    and in fact is incredibly subtle.

    After being published in the lecture notes \cite{scholze2019Analytic},
Scholze in 2020 writes the following remarkable blog post,%
\footnote{Liquid tensor experiment on \url{https://xenaproject.wordpress.com/2020/12/05/liquid-tensor-experiment/}}
proposing the challenge of computer-verifying the proof of this theorem.

    The following quotes of the post should suffice to get the functional analytic reader to become a little excited about
    the preceding theorem.
\begin{itemize}
\item \emph{I want to make the strong claim that in the foundations of mathematics, one should replace topological spaces with condensed sets (except when they are meant to be topoi - topoi form a separate variant of topological spaces that is useful, and somewhat incomparable to condensed sets).
This claim is only tenable if condensed sets can also serve their purpose within real functional analysis.}

\item  \emph{With this theorem, the hope that the condensed formalism can be fruitfully applied to real functional analysis stands or falls. I think the theorem is of utmost foundational importance, so being 99.9\% sure is not enough.}

\item \emph{If it stands, the theorem gives a powerful framework for real functional analysis, making it into an essentially algebraic theory. [\dots] Generally, whenever one is trying to mix real functional analysis with the formalism of derived categories, this would be a powerful black box. As it will be used as a black box, a mistake in this proof could remain uncaught.}

\item \emph{I spent much of 2019 obsessed with the proof of this theorem, almost getting crazy over it. In the end, we were able to get an argument pinned down on paper, but I think nobody else has dared to look at the details of this, and so I still have some small lingering doubts.}

\item \emph{I think this may be my most important theorem to date. (It does not really have any applications so far, but I’m sure this will change.)}

\item \emph{It may well be the most logically involved statement I have ever proved. (On the other hand, if I want to claim that the theorem on liquid vector spaces makes it possible to black box functional analysis, hard estimates have to be \textit{somewhere}.
 Better be sure the estimates actually work...!)}
\end{itemize}

In fact, as can be seen at \emph{Half a year of the Liquid Tensor Experiment: Amazing developments}%
\footnote{\url{https://xenaproject.wordpress.com/2021/06/05/half-a-year-of-the-liquid-tensor-experiment-amazing-developments/}}
or at \emph{Completion of the Liquid Tensor Experiment} on the Lean community blog%
\footnote{\url{https://leanprover-community.github.io/blog/posts/lte-final/} with bueprint \url{https://leanprover-community.github.io/liquid/}},
by now, the proof has been fully verified in \texttt{Lean} after one and a half years,
which has lead to a lot of attention inside and outside of mathematics, see, e.g. the articles in Nature%
\footnote{Mathematicians welcome computer-assisted proof in ‘grand unification’ theory, \url{https://www.nature.com/articles/d41586-021-01627-2}}
and Quanta magazine%
\footnote{\enquote{Lean} computer program confirms Peter Scholze proof, \url{https://www.quantamagazine.org/lean-computer-program-confirms-peter-scholze-proof-20210728/}}.
\end{remark}
Clearly, we cannot give an appropriate overview of the beautiful proof of this theorem.
However, there are many functional analytic ideas involved, including some very cleverly placed switches between Banach spaces and dual spaces in order to be able to play with norms and compactness.
We suggest the reader to spend some time enjoying the proof in \cite{scholze2019Analytic}, lectures 6--9.

We recap what this implies for the category of liquid complete $\R$-modules (also called \idx{liquid vector spaces}), which is a consequence of the general good properties of analytic rings discussed in the last section.
For the following statement see \cite[3.11]{scholze2022complex} or \cite[6.5]{scholze2019Analytic}.

\begin{corollary}
Consider a condensed abelian group $E$.
Then the following conditions are equivalent formulations of $p$-liquidity of $E$.
\begin{enumerate}[(i)]
    \item For all $S\in \extr$ and $f\colon S\to E$ and $q<p$ there is a unique extension to a condensed $\R$-module homomorphism

\begin{center}\begin{tikzcd}
	{\mcM_q(S)} & E \\
	S
	\arrow[dashed, from=1-1, to=1-2]
	\arrow["\delta", from=2-1, to=1-1]
	\arrow["f"', from=2-1, to=1-2]
\end{tikzcd}\end{center}
\item For all $S\in \extr$ and $f\colon S\to E$ there is a unique extension to a condensed $\R$-module homomorphism

\begin{center}\begin{tikzcd}
	{\mcM_{<p}(S)} & E \\
	S
	\arrow[dashed, from=1-1, to=1-2]
	\arrow["\delta", from=2-1, to=1-1]
	\arrow["f"', from=2-1, to=1-2]
\end{tikzcd}\end{center}
    \item $E$ is cokernel of a map
    \[\bigoplus_j \mcM_{<p}(S_j)\to \bigoplus_i\mcM_{<p}(T_i)\]
\end{enumerate}
Furthermore, the full subcategory of liquid abelian groups $\liq_p\sub \cab$ fulfills the following.
\begin{enumerate}[(i)]
    \item It is an abelian subcategory of $\cab$ stable under all limits, colimits, extensions, internal $\Ext$'s and (therefore) internal $\hom$'s.
    \item There exists liquidification, being the left adjoint to $\liq_p\sub \cab$, and sending $\Z[S]$ to $\mcM_{<p}(S)$.
    In particular, the liquidification of $\Z$ is $\R$.
    \item There is a unique symmetric monoidal completed tensor product on $\liq_p$ making liquidification symmetric monoidal, and still fulfilling the tensor-hom adjunction.
    This in particular implies that all objects in $\liq$ automatically have unique functorial $\R$-module structure.
\end{enumerate}

    The analogous statement for the derived category holds, see~\ref{thm:der-an-rings-gut}.
    In addition, the derived completed tensor product is given by the left derived functor of the tensor product in $\liq$, as the condition in~\ref{lem:cond-lder-tp} is fulfilled.
\end{corollary}
\begin{corollary}[Complexification]
    There is a canonical choice of analytic ring structure on $\C$-vector spaces, induced by base change.
    For this, note that we can interpret everything in the $\infty$-derived category and use base change there.
    This means that
    \[\mcM^{\C}_{<p}(S)=(\C[S])^{\liq}\]
    induces an analytic ring structure on $\C$-vector spaces, and
    the left adjoint mapping $V\mapsto V\otimes_{\R}\C$ yields a complexification and maps onto an equally well behaved category as liquid $\R$ -vector spaces.
\end{corollary}

As every locally convex space is $\mcM_1$-complete, it is in particular liquid.
In fact, we have the following characterisation in the quasiseparated setting from 2.14 in \cite{scholze2022complex}, being quite similar to a $p$-convex version of Krein's theorem (II.4.3 in \cite{schaefer1971Topological})

\begin{definition}
    Let $E$ be a quasiseparated $\R$-module, $K\sub E$ be a quasicompact sub-condensed set, and $q\ge 0$.
    We say that $K$ is \idx{$q$-convex}, if for all $\lambda_1,\dots, \lambda_n$ with
    \[\sum_{i=1}^n\|\lambda_i\|^q\le 1,\]
    and $x_i\in K$, the linear combination stays in $K$,
    \[\sum \lambda_i x_i\in K.\]

    Note that since we assumed everything to be quasiseparated resp. qcqs,
    this is really a classical map between locally compact Hausdorff spaces, and the property reduces to

\begin{center}\begin{tikzcd}
	& E \\
	{(\R^{n})_{\ell^q\le 1}\times K^n} & K
	\arrow[from=2-1, to=1-2]
	\arrow[dashed, from=2-1, to=2-2]
	\arrow[hook, from=2-2, to=1-2]
\end{tikzcd}\end{center}
\end{definition}
\begin{theorem}
Consider any quasiseparated $\R$-module $V$.
Then $V$ is $p$-liquid precisely if any quasicompact subobject of $V$ is contained in a quasicompact $q$-convex subobject (possibly depending on $q$) for every $q<p$.

Moreover, if it is $\mcM_1$-complete then the factorisation property holds for $1$-convex subobjects.
\end{theorem}
\begin{proof}
  This is \cite[2.14]{scholze2022complex}.

First, take any $\mcM_p$-complete quasiseparated vector space $E$, and a quasicompact $K\sub E$.
Since $E$ is quasiseparated, $K$ is compact Hausdorff,
and hence admits a surjection $S\twoheadrightarrow K$ for some $S\in \extr$.
Consider
\[f\colon S\twoheadrightarrow K\hookrightarrow E.\]
By definition of $\mcM_p$-completeness, this map factors (uniquely) $\R$-linearly over

\begin{center}
\begin{tikzcd}
	{\mcM(S)_p} & E \\
	S
	\arrow["T", from=1-1, to=1-2]
	\arrow["\delta", from=2-1, to=1-1]
	\arrow["f"', from=2-1, to=1-2]
\end{tikzcd}
\end{center}
By compactness of $S$, this map factors through $\mcM_p(S)_{\le c}$ for some $c$.
But hence the image $V\sub E$ is contained in the image of the restriction of $T$ in

\begin{center}\begin{tikzcd}
	{(\mcM(S)_p)_{\le c}} & E \\
	S
	\arrow["T", from=1-1, to=1-2]
	\arrow["\delta", from=2-1, to=1-1]
	\arrow["f"', from=2-1, to=1-2]
\end{tikzcd}.
\end{center}
The set $(\mcM(S)_p)_{\le c}$ is a $p$-convex compact set by writing
\[\R^n_{\ell^p\le 1}\times\varprojlim\R[S_i]_{\ell^p\le c}=\varprojlim \R^n_{\ell^p}\times \R[S_i]_{\ell^p \le c}, \]
and now noting that
$\R[S_i]_{\ell^p\le c}$ is $p$-convex by Cauchy-Schwarz,
\[\sum |\lambda_i x_i|^p=\sum |x_i^p||\lambda_i^p|\le (\sum |x_i|^p)(\sum |\lambda_i|^p)\le c,\]
this property is preserved by continuous linear maps, hence $T(\mcM(S)_p)_{\le c}$ is a compact $p$-convex subset of $E$ containing $K$.

The converse needs some more work, we give a fairly detailed sketch of the argument, for a precise proof see \cite{scholze2022complex}.

Let $V$ be a quasiseparated condensed vector space with the factorisation property.
Take $f\colon S\to V$ for some $S\in \extr$, and some $q<p$.
We want to find an $\R$-linear extension $\mcM_{q}(S)\to V$.
Uniqueness is clear by denseness of linear combinations of dirac measures in $\mcM_1$, it remains to show existence.

The idea to do this clearly is to approximate any measure with Dirac measures -- as any measure is a limit of finite linear combinations of Dirac measures.
However, we want a \textbf{continuous} extension to the whole space $\mcM_q(S)$ (or, for the sake of simplicity, to the compact unit ball $\mcM_q(S)_{\le 1}$),
and not just a pointwise extension.
So we need a \textbf{uniform} way to approximate the measure by finite convex combinations of Dirac measures.

This is done in a very clever way in \cite{scholze2019Analytic} and \cite{scholze2022complex} by passing from $\mcM_p$ to a slightly more general space of measures as follows:

Define the analytic ring of \emph{overconvergent Laurent series} (which is a priori very creepy to a functional analyst) in complete analogy to $\mcM_{<p}$ as follows.
For any finite $S$ define a condensed ring with
\[
(\Z((T))_{r})_{\le c}[S](S')
=\{\sum_{n\in \Z, s\in S} a_{n,s} T^n [s]\, \mid\, a_n\in C(S', \Z),\, \sum_{n\in \Z, s\in S}|a_{n,s}|r^n\le c\},
\]

\begin{align*}
\Z((T))_r[S](S')
&=\bigcup_{c>0} (\Z((T))_r[S])_{\le c}(S') \\
&=\{\sum_{n\in \Z,\, s\in S} a_{n,s} T^n[s]\,\mid\, a_{n,s}\in C(S',\Z),\, \sum_{n\in \Z, s\in S}|a_{n,s}|r^n <\infty.\},
\end{align*}

\[\Z((T))_{>r}[S]=\bigcup_{p>r} \Z((T))_{p}[S].\]

In particular, define $\Z((T))_{r}=\Z((T))_r[\ast]$.
Now, for any profinite set $S=\varprojlim S_i$, define
\[\mcM(S,\Z((T))_r)_{\le c}=\varprojlim (\Z((T))_r)_{\le c}\]
and
\[\mcM(S,\Z((T))_r)=\bigcup_{c>0}\mcM(S,\Z((T))_r)_{\le c}.\]

These so-called overconvergent Laurent series yield a generalisation of $\mcM_{p}(S)$ as follows (see Theorem 6.9 in \cite{scholze2019Analytic}).

Plugging in $T=\lambda$ for some fixed $\lambda\in \R$ yields a natural isomorphism
\[\mcM(S,\Z((T))_r)\otimes_{\Z((T))_r}\R=\mcM_q(S),\]
where $\lambda^q=r$.

This is roughly the same as obtaining the real numbers from laurent series by plugging in some $\lambda$, which determines the \enquote{base} of the representation.
E.g., for $\lambda=10$ we obtain $\R$ from laurent series in $\Z$ via classical decimal representation.

The above equality by currying implies that $\R$-linear maps from $\mcM_q(S)$ are uniquely given by
$\Z((T))_r$-linear maps from $\mcM(S,\Z((T))_r)$.

Thus the idea is to extend the function $f$ from above to the space $\mcM(S,\Z((T))_r)$ rather than $\mcM_q(S)$.
The good thing about $\mcM(S,\Z((T))_r)$ now is, that there is a canonical uniform way of approximating measures via Dirac measures, by truncating the allowed exponents in the power series.

Explicitly, for $n\in \N$, consider
\[\mcM(S,\Z((T))_r)_{\le C}\sub \prod_{k\in \Z}\Z[S] T^k\twoheadrightarrow \prod_{k\le n}\Z[S] T^k,\]
and note that this induces an epimorphism from the union
\[\mcM(S,\Z((T))_r)\twoheadrightarrow \prod_{k\le n}\Z[S] T^k .\]

Now, on $\prod_{k\le n}\Z[S] T^k$ the morphism $f$ clearly is canonically expandable,
as by adjunction $f\colon S\to E$ induces a $\Z[T,T^{-1}]$-linear extension
\[\Z[T,T^{-1}][S]\to V,\]
and now we note $\bigoplus_{k\le n}\Z[S] T^k\sub \Z[T,T^{-1}][S]$.
Hence, define $f_n$ via
\[\mcM(S,\Z((T))_r)\to \bigoplus_{k\le n}\Z[S]T^k\to \Z[T,T^{-1}][S]\to E.\]
Now, it suffices to show that $f_n$ converges pointwise to a continuous function, as this limit would correspond to the extension on the whole $\Z[T,T^{-1}][S]$.

For this, by the colimit topology it suffices to show that for every $c>0$, the restriction of $f_n$ onto $\mcM(S,\Z((T))_r)_{\le c}$ converges pointwise.

For this, we use the Cauchy-completeness type lemma directly after this proof.
    We first see that there exists a compact $K\sub V$ with all $f_i$ mapping into $K$;
    for this use that the image of $f$ is compact in $V$, and hence by assumption contained in a $q'$-convex subset $K$ for some $q<q'<p$.
    Now, the image of $\mcM(S,\Z((T))_r)\le c$ under evaluation at $\lambda$ is $\mcM_q(S)$,
    hence the image is contained in the $q$-convex hull of $K$, which is again $K$, as $K$ is $q'$ convex.
Hence all $f_n$ map into $K$.
    For the uniform bound $\lambda_n$ now we use the fact, that $\ell^{q'}$ and $\ell^q$ have different growth rates,
and hence we can find a nullsequence $\lambda_n$ such that any element $\ell^{q}$ multiplied with $1/\lambda_n$ still lands inside $\ell^{q'}$.
E.g., this can be achieved by
\[(f_n-f_m)(\sum c_i T^i)=\sum_{n<i\le m}f(c_i)\lambda^i=a^{-n}\sum_{n<i\le m}f(c_i)a^{n-i}(a\lambda)^i\]
and noting that by choosing $a=1+\eps$ for very small $\eps$, the last part
\[\sum_{n<i\le m}f(c_i)a^{n-i}(a\lambda)^i\]
can be chosen to have small enough growth to still land inside $K$.

Note that this last part is the first argument that heavily uses the $\mcM_{<p}$-structure and can not be made in $\mcM_1$ - for any nullsequence there exists a $\ell^1$-sequence such that $1/\lambda_n$ times the sequence is no longer in $\ell^1$.
In other words, the speed of approximation of measures by dirac measures can get arbitrarily bad, destroying the uniform argument.
\end{proof}
In the proof, we used the following lemma,
which seems to be some kind of Cauchy completeness of spaces of continuous functions with appropriate convergence.
This is \cite[2.17]{scholze2022complex}.
\begin{lemma}
Consider a compact Hausdorff space $T$, a quasiseparated $\R$-module $E$ and any quasicompact $p$-convex $K\sub E$.
Then the following statements hold.
    \begin{enumerate}[(i)]
        \item $\bigcap_{c>0}cK=\{0\}$.

        \item Consider any sequence of continuous maps $f_i\colon T\to K$.
        Assume furthermore that $(f_i)$ is Cauchy in the following sense:
There exists a nullsequence $0\le \lambda_i\sub \R$ such that for all $n$ and all $m\ge n$
         \[f_n(t)-f_m(t)\in \lambda_n K\qquad\text{for all}\quad t\in T.\]
Then there exists a continuous function $f_\infty\colon T\to K$ such that pointwise $\lim_{n\to \infty}f_n(t)=f_\infty(t)$.
    \end{enumerate}
\end{lemma}
\begin{proof}
    Assume that $x\in \bigcap_{c>0}cK$.
    Then $cx\in K$ for all $c>0$, and by $p$-convexity $\R x\sub K$.
    Thereby $\overline{\R x}\sub K$, being a compact Hausdorff $\R$-vector space and thus $0=\R x$, implying $x=0$.

    The other statement can also be verified elementarily, see \cite[2.17]{scholze2022complex}.
\end{proof}

The last estimate in the proof of the theorem above used the $q<q'<r$-interplay, which can not be exploited in $\mcM_1$.
However, as any locally convex complete vector space is $\mcM_1$-complete by a similar approximation argument,
we wonder whether this description can also be made for $\mcM_1$.
\begin{question}
    Does a similar description work for $\mcM_1$-complete quasiseparated $\R$-modules?
    I.e., is a qs $\R$-module $E$ such that any quasicompact $K\sub E$ is contained in an absolutely convex compact $E$,
    automatically $\mcM_1$-complete?
\end{question}

   Every complete locally convex vector space fulfills this property, essentially due to Krein's theorem (see II.4.2 in \cite{schaefer1971Topological}) or by noting that they are $\mcM_1$-complete.
\begin{proposition}
Let $V$ be any complete locally convex vector space and $S$ any compact Hausdorff space with map $f\colon S\to V$.
Then $f$ factors over a compact absolutely convex subset $K\subseteq V$.
\end{proposition}
\begin{proof}
Let $V$ be a locally convex vector space.
If $V$ is complete, then any continuous map $f\colon S\to V$ for any compact Hausdorff space $S$ has a unique extension to a continuous linear map $\mathcal{M}_1(S)\to V$.
This yields a factorisation of $f$ over some absolutely convex compact subset,
as the image of $f$ is norm-bounded and hence contained in some absolutely convex compact ball in $\mcM_1$
\end{proof}
    
\begin{conjecture}
The converse of this proposition holds in an appropriate sense.
\end{conjecture}

\subsubsection{Basic functional analysis}

We first recall some classical theory of complete locally convex vector spaces.

\begin{definition}[Classical Smith spaces]
Define a \textbf{classical Smith space}\index{Smith space (classical)} as a complete locally convex vector space $V$ (always assumed to be Hausdorff) that carries the compactly generated topology of one compact Hausdorff absolutely convex subset $K$,
\[
V=\bigcup_{c>0}cK.
\]
\end{definition}

Using this, we obtain the \textbf{classical Smith duality} which is also some type of \textbf{Pontryagin result}.
\begin{lemma}[Smith duality]
The dual space of a Banach space $E$ equipped with compact open topology is a classical Smith space.

Conversely, the dual space of a classical Smith space $V$ equipped with compact open topology is a Banach space.
 
The canonical map $E\to E''$ is an isomorphism.
\end{lemma}

\begin{proof}
See \cite{Smith1952} or 3.8 in \cite{scholze2019Analytic}.
\end{proof}

\begin{corollary}[\idx{Banach-Alaoglu}]
The unit ball of the dual space of a Banach space is weak*-compact.
\end{corollary}

For general complete locally convex vector spaces we have a similar Pontryagin type result using pairings.

\begin{lemma}
For a complete locally convex vector space $E$ consider the canonical pairing
$E\times E'\to \R$. 
Equipping $E'$ with the induced $\sigma(E,E')$-topology yields
\[
    E\simeq E''.
\]
\end{lemma}

\begin{proof}
This is IV.1.2 in \cite{schaefer1971Topological}.
\end{proof}

One of the central results in the classical locally convex theory is Hahn-Banach, which can be extended to the following Hahn-Banach-Schaefer injectivity result.

\begin{lemma}
For $S\in \extr$, a Banach space $E$ and a normed subspace $F\sub E$ (with subspace topology),
a continuous linear function $f\colon F\to C(S,\R)$ has an extension $E\to C(S,\R)$.
\end{lemma}

\begin{proof}
See II.7.10 in \cite{schaefer1974Banachlattices}.
\end{proof}

Note that those lifting properties often are formulated in terms of $L^\infty$-spaces,
or in terms of the projectivity of the corresponding $L^1$-spaces.
That these theorems agree is essentially due to a special case of \idx{Gelfand duality}.

\begin{lemma}[Gelfand type result]
The space $L^\infty(\Omega,\R)$ for $\Omega$ a (separable) measure space is isomorphic to $C(S,\R)$ for a $S\in \extr$. 
\end{lemma}

We furthermore have the following commutativity results of dual spaces and (co)limits. This is some kind of cocompactness result of $\R$.

\begin{lemma}
Consider a locally convex space $E$ which is a projective limit along epimorphisms of locally convex spaces $E_i$.
Then $E'$ with Mackey topology agrees with the inductive limit of the dual spaces $E_i'$ (also equipped with Mackey topology).

Dually, the weak dual (weak*-topology) of an inductive limit $E=\varinjlim E_i$ agrees with the projective limit of the weak duals $(E_i',\sigma(E_i, E_i'))$.

\end{lemma}
\begin{proof}
This is IV.4.4 and IV.4.5 of \cite{schaefer1971Topological}.
\end{proof}

Finally, we obtain the following strong representation result. 

\begin{lemma}
Every complete locally convex vector space $E$ is isomorphic to a cofiltered limit along epimorphisms of Banach spaces.
\end{lemma}

\begin{proof}
See II.5.4 in \cite{schaefer1971Topological}.
\end{proof}

Note that the theory for filtered colimits along inclusions is more subtle.
We regard this as a possible hidden appearance of the fact that (uncountable)
filtered colimits along inclusions should not be computed
in $\Top$ but rather in $\cond$.

Another important aspect of functional analysis are tensor products.
There are many (fourteen, to be precise) different tensor products on locally convex spaces, but the projective tensor product is universal with respect to bilinear maps.

\begin{lemma}
There is a tensor product $A\otimes_\pi B$ of locally convex vector spaces $A$ and $B$ such that bilinear maps $A\times B\to C$ agree with linear maps $A\otimes_\pi B\to C$.
\end{lemma}

\begin{proof}
This is III.6.2 in \cite{schaefer1971Topological}.
\end{proof}

Now we see how to translate parts of this into the condensed setting.
We will develop some basic functional analysis in $\mcM_1$-complete condensed $\R$-vector spaces,
which is not really restricted to the liquid setting as it remains purely in the quasiseparated setting.
However it gives a feeling for how classical arguments translate.

We define condensed locally convex complete $\R$-vector spaces as follows.

\begin{definition}
Define the category of \textbf{condensed complete locally convex vector spaces} as the full subcategory of quasiseparated $\mcM_1$-complete objects in $\RRMod$.
\end{definition}

\begin{lemma}
The category of condensed locally convex vector spaces is a subcategory of $\RRMod$ stable under all limits, filtered colimits along inclusions and images.
\end{lemma}
\begin{proof}
This is part of 4.6 in \cite{scholze2019Analytic}.
\end{proof}

\begin{lemma}
Clearly, every classical locally convex vector space induces a condensed complete locally condensed vector space. 
\end{lemma}

Analogous to the classical setting, we define Smith spaces and Banach spaces.

\begin{definition}
    A \textbf{(condensed) Smith space}\index{Smith space (condensed)} is an $\mcM_1$-complete condensed $\R$-module $E$ such that there exists an absolutely convex compact Hausdorff $K\sub E$ such that
    \[E=\bigcup_{c>0}cK.\]
\end{definition}

\begin{remark}
A Smith space,  as the union of quasiseparated spaces, is automatically quasiseparated.
\end{remark}

\begin{question}
Is the $\mcM_1$-completeness in the above definition automatic?
Note that the liquidity is automatic.
\end{question}

\begin{lemma}
   The category of Smith spaces agrees with the category of classical Smith spaces.
\end{lemma}
\begin{proof}
    See \cite[4.2]{scholze2019Analytic}.
\end{proof}
\begin{lemma}
Clearly, for any profinite (or $\CHaus$) $S$, $\mcM_1(S)$ is a Smith space.
\end{lemma}
Duality should give us the correct notion of Banach spaces, essentially due to the gauge-norm of $K$.
\begin{definition}
    The \idx{dual vector space} of any condensed $\R$-vector space $E$ is defined as
    \[E'=\ihom(E,\R)\]
\end{definition}
The following example is Ex.~2 resp.~3 on p.~46 in \cite{scholze2022complex}.
\begin{example}
\begin{itemize}
\item For $p\le 1$ we have $\mathcal{M}_p(I)'=c_0(I)$ for countable $I$, and $\ell^p(I)'=\mathcal{M}_\infty(I)$.
  \item For $p> 1$ we have $\mathcal{M}_p(I)'=\ell_q(I)$ for countable $I$, and $\ell^p(I)'=\mathcal{M}_q(I)$ where $1/p+1/q=1$.%
\footnote{Here we employ the funny convention of \cite{scholze2022complex} to use subscripts for equipping the product topology and superscripts for the usual norm topology.}
\end{itemize}
\end{example}

\begin{definition} A \idx{Banach space} is an $\R$-vector space $E$ of the form
    \[F'=\ihom(F,\R)\]
    for some Smith space $F$.
\end{definition}
\begin{lemma}
    This agrees with the classical definition of a Banach space.
\end{lemma}
\begin{proof}
This is part of 4.7 in \cite{scholze2019Analytic}.
\end{proof}

Again, we want some Pontryagin type result (Smith duality).

\begin{lemma}[Pontryagin type result]
    The internal $\hom$ installs a duality between the category of Banach spaces and the category of Smith spaces,
    i.e., for any Banach space $E$ and Smith space $F$ the space $E'$ is Smith and $F'$ is Banach, and furthermore
    \[E''=E\qquad \text{and}\qquad F''=F.\]
\end{lemma}

\begin{proof}
See 4.7 in \cite{scholze2019Analytic} for a nice proof.
\end{proof}

\begin{corollary}[Condensed Banach Alaoglu]
The unit ball of the dual space of a Banach space is $\qcqs$.
\end{corollary}

This duality furthermore implies some Hahn-Banach type result.
\begin{theorem}[Condensed Hahn-Banach-Schaefer]
Every space of the form $C(S,\R)$ for $S\in\extr$ is an injective object with respect to reflective monomorphisms of
Banach spaces, and $\mcM_1$ is projective in Smith spaces.
\end{theorem}
\begin{proof}
See the discussion after 4.7 in \cite{scholze2019Analytic}.

The projectivity of $\mcM_1(S)$ follows by projectivity of $\R[S]$ in $\RRMod$ and the observation
that images in condensed complete locally convex vector spaces may be computed inside $\RRMod$.
Hence every effective epimorphism in complete locally convex vector spaces is an epimorphism in $\RRMod$.
For such colimits we obtain
\[\hom(\mcM_1(S),\varinjlim X_i)=\hom_{\RRMod}(\R[S],\varinjlim X_i).\]
The injectivity now follows by using condensed Smith duality.
\end{proof}

\begin{remark}
However, we note that one has to be careful about the notion of injectivity here.
It is certainly not true that $\R$ is injective with respect to all monics.
Indeed, using the monic of some dual Banach space $E_{\|\cdot\|}'$ equipped with norm topology into
$E'$ with compact open topology, the lifting property

\begin{center}\begin{tikzcd}
	{E'} \\
	{E'_{\|\cdot\|}} & \R
	\arrow[dashed, from=1-1, to=2-2]
	\arrow[hook, from=2-1, to=1-1]
	\arrow[from=2-1, to=2-2]
\end{tikzcd}\end{center}
would imply $(E'_{\|\cdot\|})'\sub E''$ which is certainly not true for non-reflexive $E$.
A simpler example is given by the embedding $\ell^1\hookrightarrow c_0$.
Here, the lifting property would yield $(\ell^1)'=\ell^\infty\hookrightarrow c_0'=\ell^1$.
\end{remark}

Lastly, we obtain a dual representation result to the classical one (which can be found in \cite[3.7, 4.6]{scholze2019Analytic}).

\begin{lemma}
Every $\mcM_1$-complete vector space is the filtered colimit of its Smith subspaces. 
\end{lemma}

Again, like in the classical theory, there exists a universal tensor product.

\begin{lemma}
There is a symmetric monoidal tensor product $\otimes_\pi$ of $\mcM_1$-complete vector spaces representing bilinear maps.
Furthermore, the tensor product is cocontinuous in each coordinate and makes $\mcM_1$ symmetric monoidal.
\end{lemma}

\begin{proof}
This is 4.9 in \cite{scholze2019Analytic}.
\end{proof}

We give an elementary example of an interesting factor space.

\begin{example}
Consider the space $E=\bigcup_{r>0}[-r,r]^\N$ (similar results hold for $\R^\N$ or $\ell^\infty$) and the space $c_f$ of finite sequences equipped with colimit topology,
\[
    c_f=\bigcup_{n\in \N}\R^n.
\]
Both are liquid vector spaces (since they are complete locally convex)
and the injections  $\R^n\hookrightarrow E$ induce an monomorphism
\begin{center}\begin{tikzcd}
	{\phi} & {E}
	\arrow[hook, from=1-1, to=1-2]
\end{tikzcd}\end{center}
We want to study the cokernel $E/c_f$ (in liquid vector spaces, or, equivalently, in condensed abelian groups).

The underlying set of $E/c_f$, $E/c_f(\ast)$,
is given by equivalence classes of bounded sequences by the equivalence relation
which identifies $(x_n)$ with $(y_n)$ if there exists an $N\in\N$
such that $x_k=y_k$ for all $k\ge N$.
We call these equivalence classes \emph{tails}.

The topology on these \emph{tails} is given as follows:
The relatively compact sequences (the value on $\beta\N$) are given by the factor
\[
E/c_f(\beta\N)=E(\beta\N)/c_f(\beta\N)=\varinjlim_n C(\beta\N,E)/C(\beta\N, \R^n)=\varinjlim_n\bigcup_{r>0} C(\beta\N,[-r,r]\N)/C(\beta\N,\R^n).
\]
This can be simplified further to
\[
\bigcup_{r>0}\varinjlim_{n\in \N} ([-r,r]^\N)^\N/([-r,r]^n)^\N.
\]
This means, that the maps $\beta\N$ to $E/c_f$ are given by sequences
of uniformly bounded sequences modulo the equivalence relation of identifying
two sequences of sequences if they differ only in finitely many coordinates.

E.g. the relatively compact sequence
\[
    n\mapsto\begin{cases} e_1 & n\,\text{odd}\\e_2 & n\,\text{even}\end{cases},
\]
where $e_1 = (1,0,\ldots)\in E$ and $e_2 = (0,1,0,\ldots)\in E$,
is equivalent to the constant sequence $0$.

On the other hand consider for example the sequence $(e_n/n)_{n\in\N}$.
It converges to $0$ (as above, $e_n$ is the $n$-th unit vector).
The corresponding underlying sequence of \emph{points} is given by the constant
$0$-sequence $(0)_{n\in\N}$,
but as a \emph{sequence} it is not identified with the constant sequence $(0)_{n\in\N}$.
The same statement holds for the sequence $(e_n)_{n\in\N}$.
In contrast, the sequence $(n\cdot e_n)_n$ is not relatively compact
and neither is $(n\cdot e_1)_n$,
although both are the underlying sequences of points $0$.

Since the value on $\alpha\N$ is given by those elements on $\beta\N$
that are constant on $\beta\N\setminus \N$,
we can produce the convergent sequences of $E/c_f$.
In fact, all the above sequences were convergent.
Corresponding examples of uniformly bounded non-convergent sequences of sequences
can be build accordingly.

This also shows that the quotient is highly not quasiseparated,
and therefore it is necessary to leave the quasiseparated $\mcM_1$-complete setup
and work in liquid vector spaces.
\end{example}

\subsubsection{Some questions}

Next, we sketch a rough conjecture of how to extend these results to the full locally convex strength.

\begin{conjecture}\label{conj:M1_vs_lctvs}
The following categories are equivalent.
\begin{enumerate}[(a)]
\item Condensed complete locally convex vector spaces, i.e., quasiseparated $\mcM_1$-complete $\R$-vector spaces, written as $\RRMod_{\mcM_1}$.
\item Complete locally convex vector spaces (always assumed to be Hausdorff).
\item Dual spaces of all $\mcM_1$-complete $\R$-vector spaces.
\end{enumerate}
Furthermore, the following results hold. 
\begin{enumerate}[(i)]
\item  $\ihom(-,\R)$ installs an autoduality of this category (full Pontryagin).
\item The spaces of the form $C(S,\R)$ for $S\in \extr$ are a densely cogenerating class (representation result) of cocompact injectives (full Hahn-Banach),
\item The spaces $\mcM_1(S)$ give a densely generating class under filtered colimits along inclusions (in $\RRMod$) (dual representation result) of compact projectives (full lifting result).
\item The described category is a full subcategory of $\RRMod$ stable under all limits, images and filtered colimits along inclusions.
These may be taken inside complete locally convex vector spaces (but filtered colimits along inclusions may not necessarily be computed inside topological spaces).
\end{enumerate}
\end{conjecture}

\begin{proof}[Possible approaches]
We rather consider the above statement as a conjecture, since the following arguments have some subtle details which have not been checked yet, and to us it would feel more natural if the category $(a)$ were slightly larger than $(b)$ and $(c)$.
Nevertheless, we sketch a possible approach of proving these kinds of results.

In particular, we think that one can read the statement about inductive limits of dual spaces as internal cocompactness of $\R$ with respect to cofiltered limits along projections.

Define $\RRMod_{\mcM_1}'$ as the full subcategory of dual spaces of $\mcM_1$-complete vector spaces. 
For such a space $F=E'$, represent $E$ as a filtered colimit $E=\varinjlim W_i$ with Smith spaces $W_i$.
Since the internal $\hom$ is continuous in the first coordinate, we obtain
\[
    F=\ihom(\varinjlim W_i, \R)=\varprojlim \ihom(W_i,\R).
\]
Thus $F$ is a cofiltered limit of Banach spaces.
But these can be computed in locally convex space,
showing that $F$ is a cofiltered limit of Banach spaces
and thereby induced by a classical locally convex vector space. 
Conversely, every complete locally convex vector space is a cofiltered limit of Banach spaces along projections,
leading to equivalence of (b) and (c)
(using the injectivity of $\R$ in complete locally convex vector spaces).

Next, we use the internal cocompactness of $\R$ with respect to projective limits along projections in locally convex vector spaces.
This implies that for any inductive limit $Y=\varinjlim \ihom(X_i,\R)$ with Banach spaces $X_i$, we have $Y=\ihom(\varprojlim X_i, \R)$ and thus (a) and (c) are equivalent.

This also implies that duality induces an autoequivalence via
\begin{align*}
    X''
    &=\ihom(\ihom(\varinjlim Y_i,\R),\R) \\
    &=\ihom(\varprojlim \ihom(Y_i,\R),\R) \\
    &=\varinjlim \ihom(\ihom(Y_i,\R),\R) \\
    &=\varinjlim Y_i \\
    &=X
\end{align*}
or 
\begin{align*}
    X''
    &=\ihom(\ihom(\varprojlim X_{i},\R),\R) \\
    &=\ihom(\varinjlim\ihom(X_i,\R),\R) \\
    &=\varprojlim \ihom(\ihom(X_i,\R),\R) \\
    &=\varprojlim X_i \\
    &=X.
\end{align*} 
This shows that (a), (b) and (c) are equivalent, as well as the Pontryagin result (i).

Since images and filtered colimits of $\mcM_1$-complete vector spaces can be computed in condensed $\R$-modules and by the condition of $\mcM_1$-completeness,
the space $\mcM_1(S)$ for $S\in\extr$ is compact projective
(with respect to filtered colimits along inclusions and images)
in $\mcM_1$-complete condensed $\R$-vector spaces
(note that this implies projectivity with respect to reflective epimorphisms,
not all epimorphisms!).
Thus, $C(S,\R)$ is cocompact and injective with respect to the appropriate limits in the category $\RRMod_{\mcM_1}'$. 

Since furthermore every Smith space (in fact every $\mcM$-complete condensed $\R$-vector space) is colimit of objects of the form $\mcM_1(S)$,
every Banach space is limit of spaces of the form $C(S,\R)$,
and this now passes to all quasiseparated $\mcM_1$-complete $\R$-modules
(since limits of limits are limits).
But since all limits of Banach spaces are locally convex,
this yields the last description.
\end{proof}

This shows how some basic functional analysis can be done in $\R$-modules.
We remark that this result also would imply some of the questions and conjectures above.
It also suggests that the subtle differences that not all filtered colimits along inclusions can be computed in $\Top$ (not even of $\CHaus$ spaces!) seems to disappear whenever considering complete locally convex vector spaces.

However, we claim that working in $\liq$ has huge advantages and will lead to a significantly stronger theory than classical functional analysis.
(See, for example, the next section which requires an analytic ring.)

The point where the liquid theory really starts to pay off is whenever one tries to use tensor products, or other homological arguments.
We give the following toy example of how one could define \emph{very abstract Cauchy problems}, due to a joint work with Robert Boehringer.

\begin{definition}
A \idx{closed operator} $A$ on $E\in \liq$ is given by two arrows 
\begin{center}\begin{tikzcd}
	{D(T)} & E
	\arrow["i", shift left, hook, from=1-1, to=1-2]
	\arrow["T"', shift right, from=1-1, to=1-2]
\end{tikzcd}\end{center}
with $D(T)\in \liq$.

Consider a closed arrow $(D(A),A,E)$, modelling a \emph{space equation}
(think, e.g., of the Laplace operator on $E=L^2(\R^3)$),
and a closed arrow $(D(T),T,F)$ modelling a \emph{time equation}
(think of $T$ being the first derivative on $F=C(\R)$).
Then the combined autonomous equation
\enquote{time derivative equals space operator} is given by considering the arrows
\begin{center}\begin{tikzcd}
	{D(A)\otimes D(T)} & {E\otimes F}
	\arrow["{A\otimes i_T}", shift left, from=1-1, to=1-2]
	\arrow["{i_A\otimes T}"', shift right, from=1-1, to=1-2]
\end{tikzcd}.\end{center}
The \emph{space of global solutions} is given by its equalizer.
If now, e.g., $F\hookrightarrow \ihom(G, H)$ we can define \emph{initial conditions} by taking the pullback with a subobject of $\ihom(\ast, H)$ for some point $\ast\to G$.
\end{definition}

\begin{question}
If now $G$ is, e.g., a semigroup object, can one develop a semigroup theory like the one in \cite{engel2000one}?
E.g., is it possible to define abstract Sobolev towers (iterated pushouts and pullbacks of closed operators).
How do, e.g., interpolation spaces fit in this picture?

We believe that this is closely related to the representation theory developed in the next section.
\end{question}

As another example, spectral theory (see, e.g., \cite{grobler1995spectral}) essentially becomes part of homological algebra.
E.g., the following classical result now has a surprisingly different proof, see 5.12 in \cite{scholze2022complex}.

\begin{proposition}
Every endomorphism on a nontrivial Banach space has nonempty compact spectrum.
\end{proposition}

We remark a completely trivial application of the liquid setting to compute the spectrum of an operator.

\begin{lemma}
Consider an operator $T\colon E\to E$ for $E$ a complex $p$-liquid vector space
(e.g., a Banach space).
Furthermore consider a $T$-invariant subspace $F\hookrightarrow E$.
Then $T$ induces operators $T_\mid$ and $T_\slash$ in the spaces $F$ and $E/F$, i.e.,
one obtains a morphism of short exact sequences.
\begin{center}\begin{tikzcd}
	0 & F & E & {E/F} & 0 \\
	0 & F & E & {E/F} & 0
	\arrow[from=1-1, to=1-2]
	\arrow[from=1-2, to=1-3]
	\arrow["{T_\mid}"', from=1-2, to=2-2]
	\arrow[from=1-3, to=1-4]
	\arrow["T", from=1-3, to=2-3]
	\arrow[from=1-4, to=1-5]
	\arrow["{T_{\slash}}", from=1-4, to=2-4]
	\arrow[from=2-1, to=2-2]
	\arrow[from=2-2, to=2-3]
	\arrow[from=2-3, to=2-4]
	\arrow[from=2-4, to=2-5]
\end{tikzcd}\end{center}
Then $\lambda$ is contained in the resolvent $\rho(T)=\C\setminus \sigma(T)$ of $T$
precisely if one of the following conditions hold
\begin{enumerate}[(i)]
    \item $T_\mid-\lambda$ is monic, i.e. $\lambda$ is not in the point spectrum of the restriction.
    \item $T_\slash-\lambda$ is epic (note that we are talking about epics in condensed $\R$-modules).
    \item The connecting homomorphism $\ker(T_\slash-\lambda)\to \coker(T_\mid-\lambda)$ is an isomorphism.
\end{enumerate}
\end{lemma}

\begin{proof}
This is just the snake lemma.
\end{proof}

Using spectral sequences of filtrations one clearly can extend this result,
and generally we believe that classical spectral theoretic questions should become solvable after this transition.
One possible question is the following, see, e.g., \cite{gluck2016growth}, \cite{bihlmaier2023compact} or \cite{grobler1995spectral} or the first authors bachelor's thesis for more.

\begin{question}[Cyclicity problem]
Consider a Banach lattice $E$ and a positive operator $T$ on $E$.
We want to know whether the peripheral spectrum (the part of the spectrum with maximal modulus) is cyclic
(i.e., for $r(T) e^{i\phi}\in \sigma(T)$ it follows that $r(T)e^{ik\phi}\in \sigma(T)$ for all $k\in\Z$).
For this we make the following reductions
\begin{enumerate}[(i)]
\item Rescale to $r(T)=1$.
\item Pass from $E$ to a bidual of an ultraproduct of $E$, hence obtain an eigenvector
   $x\in E$ with $Tx=\lambda x$, $\lambda\in \sigma_{\mathrm{per}}(T)$, as well as order completeness of the Banach lattice
        and a weak*-topology.
\item Define $u=|x|$, and for all $k\in \N$ consider the principal ideal $E_{T^kx-x}$. 
\item Now define the invariant ideal $I=\bigcup_{k\in \N}E_{T^kx-x}$ and take the quotient $E/I$ (of course in liquid spaces). 
\begin{itemize}
\item One sees that $x\notin I$
(using the subexponential growth bound on $\|T^k\|$ induced by the spectral radius of $T$),
and clearly $Tu=u$ in $E/I$, i.e., one has a dominating fixed vector.
Can one show in this highly non-separated setting the classical result that this implies $Tx^k=\lambda^kx^k$ in $E/I$?
(See \cite{lotz1968spektrum} or \cite{gluck2016peripheral}.)
\item Can one prove (e.g., using a spectral sequence of the filtration on $I$) that the point spectrum of the factor operator can be lifted to the spectrum of the whole operator on $E$?
(See \cite{grobler1995spectral} for classical versions of these spectral theoretic considerations.)
\end{itemize}
\end{enumerate}
\end{question}

We end the introduction to condensed functional analysis here and leave it for future work to tackle classical questions using this theory.
See \cite{scholze2019Analytic} or \cite{scholze2022complex} for many tensor product calculations, and much more functional analysis such as DNF spaces, etc.
See also \cite{Schefers2023} for an entire functional calculus in liquid vector spaces, and note that this uses heavily the whole machinery of $\infty$-derived categories.
A stronger, holomorphic functional calculus can be extracted from lecture~6 in \cite{scholze2022complex}.
We believe that understanding those topics will ultimately lead to strong applications wherever functional calculus is used,
hence in particular to harmonic analysis and PDE theory, see, e.g., the book series \cite{hytonen2016analysis}.


\section{Towards ergodic theory}

\quot{Es könnte sein, dass ihr da eine gute Idee habt.}{R. Nagel, 21.01.2025}

In the study of group theory one often studies the properties of a group $G$
by investigating the actions of the group on some set $X$.
Here, an action is a map $G\times X\to X$ mapping $(g,x)$ to $g.x$ such that
$e.x = x$, where $e$ is the identity element of the group,
and $g.(h.x) = (gh).x$ for every $x\in X$.
Phrased differently, one has a group homomorphism $G\to\Aut(X)$, where $\Aut(X)$
is the group of all invertible maps $X\to X$.
Given a group action we obtain such a group homomorphism by mapping $g\in G$ to the map
$x\mapsto g.x$.
On the other hand, given a group homomorphism $f\colon G\to\Aut(X)$ one obtains an action
of $G$ on $X$ via $(g,x)\mapsto f(g)(x)$.
Now one is interested in obtaining properties of $G$ through properties of the action of $G$.
One of the first but nevertheless important insights is that a group acts on itself
by left multiplication.
This means that every group is a subgroup of the automorphism group of a set $X$
(the group itself). We have just established Cayley's theorem.
Another insight is that instead of looking at actions on a set,
one should much rather look at actions on an abelian group,
or more generally, on an $R$-module $M$ for a commutative unital ring $R$.
Then an action of $G$ on $M$ is the same as that of the set-theoretic description
described above, but now one insists that all maps $m\mapsto g.m$ for $g\in G$ are $R$-linear.
Instead of looking at actions $G\times M\to M$ one often considers representations
$G\to\Aut(M) = \GL(M)$. But we have already seen that there is no real difference in those
two points of view, so we will reluctantly switch between them.
An important technique is that of \emph{group cohomology} where one closely
investigates the fixed and orbit spaces of an action on an $R$-module $M$.
This helps in classifying group extensions, amenable groups, etc.
Another result that should be pointed out here is Tannaka duality.
It shows that one can always recover the group from its category of representations.
Tannaka duality could also be called the \emph{Yoneda lemma of representation theory}
and actually is an application thereof.

As in most cases in mathematics, the structure one investigates also has a topology,
as is with many groups $G$.
Then one is interested in studying \emph{continuous} representations of $G$,
i.e., topological $R$-modules $M$ such that the action maps $g\mapsto g.m$ are continuous.
One hopes to recover (topological) properties of the group by looking
at its continuous representations.
Some results can be carried over for compact groups,
e.g. Peter-Weyl and Tannaka-Krein duality.
But in other aspects difficulties arise: e.g.,
short exact sequences of modules with $G$-action do not give rise to long exact
sequences in the cohomology anymore.

In this section we explain how one could develop a theory of continuous group
cohomology in the condensed setting.
We give some applications to ergodic theory.
For the classical theory of group cohomology we refer to chapter 6 in \cite{Weibel1994}
and the lecture notes \cite{Loeh2019}.
A treatment of continuous group cohomology can be found in \cite{Casselman1974} and \cite{Fust2023}.
Some development of condensed group cohomology has already been made in 
\cite{bhatt2014proetale} and in \cite{Brink2023}.

\subsubsection{Group actions}

Let $G$ be a condensed group with multiplication $m\colon G\times G\to G$.
We denote with $0\colon G\to *$ and $0\colon *\to G$ the obvious zero morphisms.
To begin with we give a definition of an action of $G$ on a condensed set $X$
which is just the straightforward application of the classical definition
to the condensed setting.

\begin{definition}\label{def:cond-action}\chapfour
A \idx{condensed $G$-action} on a condensed set $X$ is a morphism
\[
    a\colon G\times X\to X
\]
such that the diagrams
\begin{center}
\begin{tikzcd}
	{G\times G\times X} & {G\times X} & X \\
	{G\times X} & X & {G\times X} & X
	\arrow["{1_G\times a}", from=1-1, to=1-2]
	\arrow["{m\times1_X}"', from=1-1, to=2-1]
	\arrow["a", from=1-2, to=2-2]
	\arrow["{0\times 1_X}"', from=1-3, to=2-3]
	\arrow["{1_X}", from=1-3, to=2-4]
	\arrow["a"', from=2-1, to=2-2]
	\arrow["a"', from=2-3, to=2-4]
\end{tikzcd}
\end{center}
are commutative.
A morphism $f$ between two such sets $X$ and $Y$ is a morphism between condensed sets which preserves the action,
i.e., the diagram
\begin{center}
\begin{tikzcd}
	{G\times X} & {G\times Y} \\
	X & Y
	\arrow["{1_G\times f}", from=1-1, to=1-2]
	\arrow["a"', from=1-1, to=2-1]
	\arrow["b", from=1-2, to=2-2]
	\arrow["f"', from=2-1, to=2-2]
\end{tikzcd}
\end{center}
is commutative.
Denote by $\cond\Set_G$ the category of condensed sets with $G$-action.
\end{definition}

The commutativity of the top square on the left asserts that $g.(h.x) = (gh).x$
and the commutative triangle says that $e.x = x$ for all $g,h\in G$, $x\in X$.
The diagram for a morphism $f\colon X\to Y$ asserts $g.f(x) = f(g.x)$.

\begin{remark}\label{rem:cond-action-monad}
One could alternatively define the endofunctor
\[
    \cond\Set\to\cond\Set,\quad X\mapsto G\times X
\]
which becomes a monad by the canonical morphisms $0\times 1_X\colon X\to G\times X$
and $m\times 1_X\colon G\times G\times X\to G\times X$.
Then $\cond\Set_G$ is just the category of algebras for this monad.
\end{remark}

The interpretation of $\cond\set_G$ as the category of algebras on a monad
directly yields its bicompleteness.

\begin{proposition}\label{prop:cond-action-bicomp}
The category $\cond\Set_G$ is bicomplete and (co)limits can be computed in $\cond\Set$.
\end{proposition}

\begin{proof}
The forgetful functor $?\colon\cond\Set_G\to\cond\Set$ creates limits,
see \cite[Proposition 4.3.1]{Borceux2008a}.
Since $\cond\Set$ is complete, the same holds for $\cond\Set_G$.

On the other hand, the functor $G\times -$ preserves colimits,
since colimits in $\cond\Set$ are stable under base change.
Thus, by \cite[Proposition 4.3.2]{Borceux2008a}, the forgetful functor $?$ also creates
colimits, hence the cocompleteness of $\cond\Set_G$.
\end{proof}

\begin{example}
Let $\alpha\in [0,1]$ be irrational and consider the $\Z$-action on the torus $\T$
given by $(n,t)\mapsto n\alpha + t$.
Then the coequalizer of this action with the trivial action on $\T$ is the coequalizer
described in examples~\ref{ex:colim-not-pointwise} and \ref{ex:yoneda-colim-not-comm}.
Different than in the category of actions of $\Z$ on Hausdorff spaces,
here we actually have a nontrivial \emph{orbit space} even though each orbit
is dense in $\T$.
\end{example}

\begin{remark}
In \cite{bhatt2014proetale} the category $\Sh(\mathbf{B}G_\textrm{pro\'{e}t})$
is introduced, where $G$ is a profinite group.
It is used to study a variant of profinite group cohomology.
We are interested in more general groups, e.g. infinite discrete or compact connected groups,
so our approach differs by looking at actions of the group on the condensed sets itself
rather than on the test objects (the profinite sets).
\end{remark}

\subsubsection{Group representations on linear spaces}

Besides $G$-actions on condensed sets, i.e., representations of $G$ on $\Aut(X)$,
one may also ask for representations on linear spaces.
Classically, a representation of $G$ on an $R$-module $M$ is a group homomorphism
$G\to\GL(M)$ where $\GL(M)$ is the general linear group on $M$,
i.e., all invertible $R$-linear maps $M\to M$.
Equivalently, it is an action $G\times M\to M$ of $G$ on $M$ such that $g.(m+m') = g.m + g.m'$.
Having a group homomorphism $G\to\GL(M)$, there is a unique $R$-linear extension
$R[G]\to\GL(M)$ to the free $R$-module $R[G]$ over $G$.
Now $R[G]$ carries the structure of a (non-commutative) ring:
its multiplication is the unique $R$-linear extension of the multiplication on $G$.
Thus $M$ is a left $R[G]$-module via the assignment
\[
    \left(\sum_{g\in G} r_g g\right)\cdot m = \sum_{g\in G} r_g (g.m).
\]
On the other hand, each left $R[G]$-module $M$ is also an $R$-module by restriction
and the induced map $R[G]\to\GL(M)$ restricts to a representation $G\to\GL(M)$ of $G$
on $M$.
So instead of investigating $G$-representations over $R$,
we look for left $R[G]$-modules.
A left $R[G]$-module is an $R$-module $M$ with a map $R[G]\times M\to M$
(this is just the $R$-linear extension of the action of $G$ on $M$)
that is $R$-bilinear.
Therefore it extends to a unique $R$-linear map $R[G]\otimes_R M\to M$.
This means that a left $R[G]$-module is just a left $R$-module with map
$R[G]\otimes_R M\to M$.
This will be our next definition.

Let $(R,\mcM,\delta)$ be a commutative normalized analytic ring.
To simplify notation, we denote by $-\otimes_R-$ the completed tensor product.
We also denote by $\mcM(-)$ the functor from $\cond\set$ to $\cMod$
which is the left adjoint to the forgetful functor
$\cMod\to\RMod\to\cond\set$ and is given by $\mcM(X) =  R[X]\otimes_{(R,\mcM)}(R,\mcM)$.
It is the unique colimit preserving extension of $S\mapsto\mcM(S)$
and is symmetric monoidal w.r.t. the product on $\cond\set$.
Now $\mcM(G)\coloneq R[G]\otimes_{(R,\mcM)}(R,\mcM)$,
the free complete $R$-module over $G$, is a (not necessarily commutative) condensed ring,
where multiplication is given by the map
\[
    \mcM(G)\otimes_R \mcM(G) = \mcM(G\times G)\to \mcM(G)
\]
in $\cMod$ induced by multiplication $m$ of $G$.
Tensoring with $\mcM(G)$ in $\cMod$ yields a monad on $\cMod$.

\begin{definition}\label{def:G-mod}
We denote by $\cGMod$ the subcategory of $\cMod$ consisting of $\mcM(G)$-modules
and $G$-equivariant morphisms.
\end{definition}

\begin{remark}
Explicitly, the objects of $\cGMod$ are given by complete $R$-modules $M$
and an action (or scalar multiplication) of $G$
\[
    \mcM(G)\otimes_R M\to M,
\]
i.e., an $R$-bilinear map $\mcM(G)\times M\to M$.
The morphisms are those $R$-module homomorphisms $M\to N$ that respect the $G$-action.
So $\cGMod$ is the category of $R$-complete left $G$-modules,
or equivalently the category of left $G$-representations on complete $R$-modules.
\end{remark}

\begin{example}
As explained above, the forgetful functor $\cMod\to\cond\Set$ admits a left adjoint
$\mcM(-)\colon\cond\Set\to\cMod$ which maps a condensed set $X$ to the completion
$\mcM(X)$ of $R[X]$.
Therefore, every condensed set $X$ with $G$-action gives rise to a complete $G$-module $\mcM(X)$
since $\mcM(-)$ is symmetric monoidal w.r.t. the product in $\cond\Set$.
\end{example}

\begin{remark}
If $\mcM$ is the trivial analytic ring structure on $R$
such that $\cMod = \RMod$, then we could have also described $\cGMod = G\textrm{-Mod}$
as the category of $R[G] = \mcM(G)$-modules.
\end{remark}

\begin{warning}
Note that we do not require that $S\mapsto R[G][S]\otimes_{(R,\mcM)}(R,\mcM)$ 
defines an analytic ring structure.
We are working in the 1-categorical setting where this in general need not be the case.
Besides this is not what we want to consider.
For example, if we start with the liquid ring structure,
one has to think about $\cGMod$ as the representations of $G$
on real locally convex complete vector spaces (= liquid vector spaces).
If one works over the discrete analytic ring over $\R$,
then one obtains all representations of $G$ over (condensed) $\R$ vector spaces, etc.
\end{warning}

Nevertheless, we still have plenty of good categorical properties.
Those guarantee that we can compute limits and colimits of $G$-modules
in the category of $R$-modules and then they automatically carry a $G$-module structure.

\begin{proposition}\label{prop:free-G-mod}
The embedding $?\colon\cGMod\to\cMod$ has a left adjoint
$\mcM(G)\otimes_1 -$.
\end{proposition}

This left adjoint is called the \emph{free $G$-module over $M$}.
Explicitly, the left adjoint $\cMod\to\cGMod$ is constructed as follows:
for any complete $R$-module $M$, one can define a $G$-action on $\mcM(G)\otimes_R M$
by the multiplication on $G$, i.e.,
\[
    m\otimes 1_M\colon\mcM(G)\otimes_R\mcM(G)\otimes_R M\to \mcM(G)\otimes_R M.
\]
A morphism $f\colon M\to N$ is mapped to $1_G\otimes f\colon \mcM(G)\otimes M\to\mcM(G)\otimes N$.

\begin{proof}
This is Proposition 4.1.4 in \cite{Borceux2008a} and Theorem VI.3.1 in \cite{mac2013categories}.

The unit of the adjunction is given by the morphism $M\to\mcM(G)\otimes M$
induced by $0\colon *\to G$ and the counit is given by the action of $G$, i.e.,
$\mcM(G)\otimes_R N\to N$.
Then one verifies that they indeed satisfy the triangle equalities
(which follow from the axioms of a monad and of a module homomorphism).
\end{proof}

We denote this functor by $\mcM(G)\otimes_1-$ because objectwise it is given by tensoring with $\mcM(G)$,
but the action only happens on the first component.
This construction is investigated more closely later.

\begin{proposition}\label{prop:G-mod-prop}
The category $\cGMod\subset\cMod$ is a reflective abelian subcategory
closed under all limits and colimits.
\end{proposition}

\begin{proof}
The closedness under limits follows from \cite[Proposition 4.3.1]{Borceux2008a}
since $\cMod$ is complete and the embedding $?\colon\cGMod\to\cMod$ creates limits.

Again, since the functor $\mcM(G)\otimes_R -$
preserves colimits in $\cMod$ and $\cMod$ is cocomplete,
Proposition 4.3.2 in \cite{Borceux2008a} implies that
$\cGMod$ is also closed under colimits and the forgetful functor creates colimits.

The fact that the category is abelian then follows since being an abelian
category can be described purely in terms of (finite) limits and colimits.
\end{proof}

\subsubsection{Constructions of $G$-modules}

Given an $R$-module $A$, one way to obtain a $G$-module
is to construct the free $G$-module over $A$, see proposition~\ref{prop:free-G-mod}.
A much simpler way is obtained by the \emph{trivial action functor} $\triv$.
Given $A$, one obtains a $G$-module by simply defining $g.a \coloneq a$
for every $a\in A$ and $g\in G$.
Explicitly, since $\mcM(*)$ is the tensor unit in $\cMod$, the null morphism
$G\to\ast$ induces a morphism
\[
    \mcM(G)\otimes_R A\to\mcM(\ast)\otimes_R A = A.
\]
This is the trivial action on $A$.
Clearly, any morphism $f\colon A\to B$ between $R$-modules respects the trivial action,
so that $\triv$ becomes a functor from $\cMod$ to $\cGMod$.

Besides the $\Ab$-enrichment, the category $\cGMod$ has also a natural enrichment over $\cMod$.
In the classical setting this is evident.
We need some more observations for this definition in our situation.
Given two $G$-modules $M$ and $N$,
one would like to have the $G$-equivariant morphisms from $M$ to $N$
to also carry an $R$-module and a topological structure,
i.e., $\cGMod$ is enriched over $\cMod$.
To define this enrichment, we first need another construction.
Given a $G$-module $M$ and an $R$-module $A$,
one can construct, similar to proposition~\ref{prop:free-G-mod},
the free $G$-module over $A$ with values in $M$.
The underlying complete $R$-module is given by $M\otimes_R A$ in $\cMod$
and the action of $G$ is given by the action of $G$ on $M$.
We thus have defined a functor
\[
    -\otimes_1 -\colon\cGMod\times\cMod\to\cGMod.
\]

\begin{lemma}\label{lem:free-G-mod-cont}
The functor $-\otimes_1-$ is separately cocontinuous in both variables.
\end{lemma}

\begin{proof}
The composition $?\circ(-\otimes_1-)\colon\cGMod\times\cMod\to\cMod$
is equal to the tensor product functor $(?-)\otimes_R-$.
Since $?$ creates colimits, it suffices to show that the functor $(?-)\otimes_R-$
is separately cocontinuous in both variables.
But this follows from the fact that the tensor product is symmetric
and left adjoint.
\end{proof}

Now we define the enriched Hom in $\cGMod$ as
\[
    \ihom_G(M,N)(S) = \hom_G(M\otimes_1\mcM(S),N)
\]
for $S\in\extr$.

\begin{proposition}\label{prop:G-mod-ihom}
The enriched Hom $\ihom_G$ is a well-defined functor
\[
    \ihom_G\colon\cGMod^\mathrm{op}\times\cGMod\to\cMod,\quad (M,N)\mapsto\ihom_G(M,N)
\]
where $\ihom_G(M,N)(S) = \hom_G(M\otimes_1\mcM(S),N)$ for $S\in\extr$.
Moreover, the following diagram
\begin{center}
\begin{tikzcd}
	{\ihom_G(M,N)} & {\ihom_R(M,N)} & {\ihom_R(\mcM(G)\otimes_R M,N)}
	\arrow[from=1-1, to=1-2]
	\arrow[shift left, from=1-2, to=1-3]
	\arrow[shift right, from=1-2, to=1-3]
\end{tikzcd}
\end{center}
is an equalizer diagram in $\cMod$.
Here, the parallel arrows are given by the pullback of the action of $G$ on $M$
and the other by the pushforward of the action of $G$ on $N$.
\end{proposition}

\begin{proof}
The assignment $S\mapsto \hom_G(M\otimes_1\mcM(S),N)$ maps finite coproducts to finite products since $\mcM(-)$ and $M\otimes_1-$ are cocontinuous.
So it suffices to show the equalizer condition.
Since limits in $\cMod$ are computed pointwise, it suffices to show that the diagram
\begin{center}
\begin{tikzcd}
	{\hom_G(M\otimes_1\mcM(S),N)} & {\hom_R(M\otimes_R\mcM(S),N)} & {\hom_R(\mcM(G)\otimes_R M\otimes_R\mcM(S),N)}
	\arrow[from=1-1, to=1-2]
	\arrow[shift left, from=1-2, to=1-3]
	\arrow[shift right, from=1-2, to=1-3]
\end{tikzcd}
\end{center}
is an equalizer diagram in $\Ab$, and hence in $\Set$.
By the descriptions of equalizers this is equivalent to showing that
\[
    \hom_G(M\otimes_1\mcM(S),N) = \{ f\in\hom_R(M\otimes_R\mcM(S),N) :
            f\circ w_M = w_N\circ(1_G\otimes f) \}
\]
where $w_M$ and $w_N$ are the structure maps of $M$ and $N$, respectively.
But this is exactly the definition of morphism in $\cGMod$.
\end{proof}

By means of this enrichment, the forgetful functor $?\colon\cGMod\to\cMod$
is representable.
Namely, one has an isomorphism
\[
    ? \simeq \ihom_G(\mcM(G),-).
\]
This follows from proposition~\ref{prop:free-G-mod}.
For $M\in\cGMod$ we obtain, when evaluated at $S\in\extr$, that
\[
    ?M(S)=\hom_{R}(\mcM(S),?M)=\hom_G(\mcM(G)\otimes_1\mcM(S),M)=\ihom_G(\mcM(G),M)(S).
\]

Besides looking at $G$-equivariant maps from $M$ to $N$ one might also be interested
in investigating all $R$-linear maps from $M$ to $N$
and then installing a $G$-module structure on them.
This is achieved by defining a map $G\times\hom(M,N)\to\hom(M,N)$ given by
\[
    (g,f)\mapsto \left(m\mapsto g.f(g^{-1}.m)\right).
\]
We make this notion precise in our setting.
Let $M$ and $N$ be $G$-modules.
Now we want to define a morphism
\[
    w\colon\mcM(G)\otimes_R\ihom_R(M,N)\to\ihom_R(M,N)
\]
which gives the $R$-module $\ihom_R(M,N)$
(which is complete since $\cMod$ is closed under internal hom's)
of (internal) $R$-linear maps $M\to N$ the structure of a $G$-module.
Instead of looking for the map $(g,f)\mapsto g.f$ we rather define the map
$(g,f,m)\mapsto (g.f)(m) = g.f(g^{-1}.m)$ for $m\in m$.
Under the Hom-tensor adjunction this map is given by a morphism
\[
    w\colon\mcM(G)\otimes_R\ihom_R(M,N)\otimes_R M\to N.
\]
Using the antidiagonal $\nabla\colon G\to G\times G, g\mapsto (g,g^{-1})$
we obtain $w$ as the composition
\begin{center}
\begin{tikzcd}
	{\mcM(G)\otimes_R\ihom_R(M,N)\otimes_RM} & {\mcM(G)\otimes_R\mcM(G)\otimes_R\ihom_R(M,N)\otimes_RM} \\
	{\mcM(G)\otimes_RN} & {\mcM(G)\otimes_R\ihom_R(M,N)\otimes_RM} \\
	N
	\arrow[from=1-1, to=1-2]
	\arrow[from=1-2, to=2-2]
	\arrow[from=2-1, to=3-1]
	\arrow[from=2-2, to=2-1]
\end{tikzcd}.
\end{center}
Here the first morphism is just the antidiagonal $\nabla$,
followed by the action on $M$ by the inverted coordinate of $G$,
then evaluation at $M$ (the counit of the hom-tensor adjunction),
and lastly the action on $N$.
Using the axioms of the actions of $G$ on $M$ and $N$ one is able to deduce the same
for the morphism $w$, hence it actually is a $G$-action on $\ihom_R(M,N)$.
Thus we have constructed a functor
\[
    \hom^G\colon\cGMod^\mathrm{op}\times\cGMod\to\cGMod,\quad
    (M,N)\mapsto \ihom_R(M,N)
\]
with the action on $\ihom_R(M,N)$ just described.
The value on morphisms is clear.
We will later establish the relation between the two functors $\ihom_G$ and $\hom^G$.

Having seen how the morphisms in $\cGMod$ work,
we now investigate the role of tensor products.
The tensor product $M\otimes_G N$ of two $G$-modules $M$, $N$
should by given by the $R$-linear tensor product $M\otimes_R N$
but now we enforce $(g^{-1}.m)\otimes n = m\otimes (g.n)$ on elementary tensors.
Given two complete $G$-modules $M$, $N$, their $G$-tensor product is constructed by the coequalizer diagram
\begin{center}
\begin{tikzcd}
	{M\otimes_R \mcM(G)\otimes_R N} & {M\otimes_RN} & {M\otimes_GN}
	\arrow[shift left, from=1-1, to=1-2]
	\arrow[shift right, from=1-1, to=1-2]
	\arrow[from=1-2, to=1-3]
\end{tikzcd}
\end{center}
in $\cMod$, where the two parallel arrows are given by
the inverse $g\mapsto g^{-1}$ on $G$ followed by the action of $G$ on $M$ and the action of $G$ on $N$ (without inverting).
The functoriality of colimits thus yields a functor
\[
    -\otimes_G-\colon \cGMod\times\cGMod\to\cMod,\quad (M,N)\mapsto M\otimes_G N.
\]
The unit of this tensor product is $\mcM(G)$:
the zero morphism $G\to *$ induces a morphism $M\otimes_R\mcM(G)\to M$
(since $\mcM(*) = R$ is the tensor unit in $\cMod$)
which has a right inverse $M\to M\otimes_R\mcM(G)$ induced by the other zero morphism $*\to G$.
This implies the arrow $M\otimes_R\mcM(G)\to M$ to coequalize the parallel arrows
and being universal among these.
Moreover, the $G$-linear tensor product preserves colimits in both variables:
clearly, any colimit in $\cGMod$ can be computed in $\cMod$ after applying the forgetful functor.
Then using that the $R$-tensor product preserves colimits (since it has a right adjoint,
the internal $\hom$-functor) and that colimits commute with colimits,
we get the desired cocontinuity.
Let us summarize our findings.

\begin{lemma}\label{lem:G-mod-tensor}
There exists a functor
\[
    -\otimes_G -\colon\cGMod\times\cGMod\to\cMod
\]
which preserves colimits in both variables and $\mcM(G)\otimes_G M \simeq M \simeq M\otimes_G \mcM(G)$ for every $M\in\cGMod$.
\end{lemma}

Fixing $N$, we obtain a functor $-\otimes_G N\colon\cGMod\to\cMod$.
This yields a left adjoint to the previously constructed functor $\hom^G(N,-)$
of $G$-equivariant morphisms.

\begin{proposition}\label{prop:G-mod-tensor-hom}
There is an adjunction
\[
    \ihom_R(M\otimes_G N, A) = \ihom_G(M,\hom^G(N,A_\triv))
\]
natural in $A\in\cMod$ and $M,N\in\cGMod$.
\end{proposition}

\begin{proof}
We have isomorphisms between diagrams
\begin{center}
\begin{tikzcd}
	{\ihom_R(M\otimes_R N, A)} & {\ihom_R(M\otimes_R\mcM(G)\otimes_R N, A)} \\
	{\ihom_R(M,\ihom_R(N, A))} & {\ihom_R(M\otimes_R\mcM(G),\ihom_R(N,A))} \\
	{\ihom_R(M,\ihom_R(N, A))} & {\ihom_R(\mcM(G)\otimes_RM,\ihom_R(N,A))}
	\arrow["{\overline{w}_M^*}", shift left, from=1-1, to=1-2]
	\arrow["{w_N^*}"', shift right, from=1-1, to=1-2]
	\arrow["\simeq"', from=1-1, to=2-1]
	\arrow["\simeq", from=1-2, to=2-2]
	\arrow["{\overline{w}_M^*}", shift left, from=2-1, to=2-2]
	\arrow["{w_N^*}"', shift right, from=2-1, to=2-2]
	\arrow["\simeq"', from=2-1, to=3-1]
	\arrow["\simeq", from=2-2, to=3-2]
	\arrow["{{w}_M^*}", shift left, from=3-1, to=3-2]
	\arrow["{(w_B)_*}"', shift right, from=3-1, to=3-2]
\end{tikzcd}
\end{center}
where $w_M$ and $w_N$ are the action of $G$ on $M$ and $N$, resp.,
and $\overline{w}_M$ denotes the inversion on $G$ followed by the action.
Moreover, $w_B$ is the action on $\ihom_R(N,A) = \hom^G(N,A)$.
Since $\ihom_R(M\otimes_G N, A)$ is the equalizer of the upper row
and $\ihom_G(M,\hom^G(N,A))$ the equalizer of the lower row,
the desired equality follows.
\end{proof}

Summarizing our findings, the relation between the functors $\mcM(G)\otimes_1-$ and
$\hom^G(\mcM(G),-_\triv)$ is given by an adjunction with the forgetful functor.

\begin{proposition}\label{prop:G-mod-adjunction-?}
There are adjunctions
\begin{center}
\begin{tikzcd}
	\cGMod && \cMod
	\arrow[""{name=0, anchor=center, inner sep=0}, "{{?}}"{description}, from=1-1, to=1-3]
	\arrow[""{name=0p, anchor=center, inner sep=0}, phantom, from=1-1, to=1-3, start anchor=center, end anchor=center]
	\arrow[""{name=0p, anchor=center, inner sep=0}, phantom, from=1-1, to=1-3, start anchor=center, end anchor=center]
	\arrow[draw=none, from=1-1, to=1-3]
	\arrow[""{name=1, anchor=center, inner sep=0}, "{{\hom^G(\mcM(G),-_\triv)}}", curve={height=-18pt}, from=1-3, to=1-1]
	\arrow[""{name=1p, anchor=center, inner sep=0}, phantom, from=1-3, to=1-1, start anchor=center, end anchor=center, curve={height=-18pt}]
	\arrow[""{name=2, anchor=center, inner sep=0}, "{{\mcM(G)\otimes_1-}}"', curve={height=18pt}, from=1-3, to=1-1]
	\arrow[""{name=2p, anchor=center, inner sep=0}, phantom, from=1-3, to=1-1, start anchor=center, end anchor=center, curve={height=18pt}]
	\arrow["\dashv"{anchor=center, rotate=-95}, draw=none, from=0p, to=1p]
	\arrow["\dashv"{anchor=center, rotate=-85}, draw=none, from=2p, to=0p]
\end{tikzcd}.
\end{center}
That is, for any $M\in\cGMod$ and $A\in\cMod$ we have that
\[
    \ihom_G(\mcM(G)\otimes_1 A,M) \simeq \ihom_R(A,?M)
\]
and
\[
    \ihom_R(?M,A) \simeq \ihom_G(M,\hom^G(\mcM(G),A_\triv))
\]
natural in $A\in\cMod$ and $M\in\cGMod$.
\end{proposition}

\begin{proof}
The first adjunction is the one from proposition~\ref{prop:free-G-mod}
under the enrichment, so it follows pointwise on $S\in\extr$.
The second adjunction follows from proposition~\ref{prop:G-mod-tensor-hom}
and the fact that $M\otimes_G\mcM(G) = M$.
\end{proof}

We end this part with a Tannaka type result, see also 9.0.0.1 in \cite{SAG}.
Given the forgetful functor $?\colon\cGMod\to\cMod$
one can ask if it is possible to reconstruct the group $G$ from the endomorphisms on this forgetful functor.
We have already seen that the forgetful functor $?$ is representable via $\ihom_G(\mcM(G),-)$.
So using an enriched version of the Yoneda lemma we obtain
\[
    \hom(?,?) = \hom(\ihom_G(\mcM(G),-),\ihom_G(\mcM(G),-)) = \ihom_G(\mcM(G),\mcM(G)) = \mcM(G).
\]

\subsubsection{Group cohomology}

Having the most important tools leads to studying invariants of the group $G$.
Given a representation $M$, one may ask which elements of $M$ are fixed by the action of $G$.
So one asks for the subspace
\[
    \fix M = \{ m\in M : g.m = m \text{ for all } g\in G \}.
\]
This is accomplished by the fixed space functor $\fix$ which is given by
the equalizer of the two parallel morphisms $M\to\ihom_R(\mcM(G),M)$
which are given by mapping $m\in M$ to $g\mapsto g.m$ and $g\mapsto m$, respectively.
In our case, this is accomplished by the of the action of $G$ on $M$,
now instead of a morphism $a\colon\mcM(G)\otimes_R M\to M$ is represented by a morphism $M\to\ihom_R(\mcM(G),M)$,
and the trivial action $t\colon\mcM(G)\otimes M\to M$ on $M$,
which equivalently is a morphism $M\to\ihom(\mcM(G),M)$.
Then we impose on $\fix$ that the diagram
\begin{center}
\begin{tikzcd}
	{\fix M} & M & {\ihom_R(\mcM(G),M)}
	\arrow[from=1-1, to=1-2]
	\arrow["a", shift left, from=1-2, to=1-3]
	\arrow["t"', shift right, from=1-2, to=1-3]
\end{tikzcd}
\end{center}
is an equalizer diagram in $\cMod$.
Since this equalizer diagram is functorial in $M$, $\fix$ is a functor from $\cGMod$ to $\cMod$.

Dually, one can define the orbit functor $\orb$.
Classically, this is obtained by taking the quotient of $M$ by the submodule $DM$
generated by all elements of the form $g.m-m$ for $g\in G$ and $m\in M$,
\[
    \orb M = M/\langle g.m - m : g\in G,\, m\in M \rangle
\]
In our case, the orbit functor assigns to every $G$-module $M$ the coequalizer
\begin{center}
\begin{tikzcd}
	{\mcM(G)\otimes_RM} & M & {\orb M}
	\arrow["a", shift left, from=1-1, to=1-2]
	\arrow["t"', shift right, from=1-1, to=1-2]
	\arrow[from=1-2, to=1-3]
\end{tikzcd}
\end{center}
where again $a$ is the action on $M$ and $t$ is the trivial action.
This construction also yields a functor from $\cGMod$ to $\cMod$.

The fixed space and the orbit functor have alternative descriptions.
First, for any module $M$ one has $\hom(R,M) = M$
since a $R$-linear map $R\to M$ is determined on the unit $1\in R$.
Now if one installs the trivial action on $R$,
then a $G$-equivariant map from $R$ to $M$ (which is just an element of $M$, namely $f(1)$)
needs to also satisfy $g.f(1) = f(g.1) = f(1)$.
Hence $f(1)$ is a fixed point of $M$.
So $\fix M$ is given by all $G$-equivariant maps from $R$, with the trivial action, to $M$.

Dually, an equivalence class $[m]\in M/DM$ gives rise to the element $1\otimes m\in R\otimes_G M$ where $R$ has the trivial action.
For another $n\in [m]$, there exists $g\in G$ such that $g.n = m$ (or a linear combination of such elements).
Thus $1\otimes m = 1\otimes g.n = g^{-1}.1\otimes n = 1\otimes n$
and the map is well-defined.
On the other hand, every such element in the tensor product produces a generator of the
orbit of $G$ on $M$.

\begin{lemma}\label{lem:fix-orb-descr}
We have isomorphisms $\fix\simeq\ihom_G(R_\triv,-)$
and $\orb\simeq R_\triv\otimes_G -$.
\end{lemma}

\begin{proof}
The first natural isomorphism follows from the description of $\ihom_G(R_\triv,M)$
as the equalizer of the two morphisms
\[
    {\ihom_R(R,M)} \to {\ihom_R(\mcM(G)\otimes_R R,M)}
\]
where one is the pullback of the (trivial) action on $R$
and the other the pushforward of the action on $M$.
On the other hand, the fixed space of $M$ is defined to be the equalizer of
the two maps $M\to\ihom_R(\mcM(G),M)$,
one being the action on $M$ and the other the trivial action on $M$.
The identities $\ihom_R(R,M) = R$ and $\mcM(G)\otimes_R R = \mcM(G)$
then imply the claim.
The second isomorphism follows in the same manner from the general description
of $R_\triv\otimes_G M$ as a coequalizer.
\end{proof}

The fixed space functor now relates all $R$-linear maps $M\to N$
to the $G$-equivariant maps $M\to N$:
if $f\colon M\to N$ is a fixed space in $\hom(M,N)$,
then by definition of the action thereof, we have $g.f(g^{-1}.m) = f(m)$,
or equivalently $f(g^{-1}.m) = g^{-1}.f(m)$.
This means that $f$ respects the action of $G$.

\begin{corollary}\label{cor:hom-fix}
There is an isomorphism $\fix\hom^G(M,N) \simeq \ihom_G(M,N)$ natural in $M,N\in\cGMod$.
\end{corollary}

\begin{proof}
This follows from the description of $\ihom_G(M,N)$ in proposition~\ref{prop:G-mod-ihom}
as an equalizer.
\end{proof}

The following proposition summarizes the connection between the fixed space and the orbit functor.

\begin{proposition}\label{prop:G-mod-adjunction-triv}
The following adjunctions hold.
\begin{center}
\begin{tikzcd}
	\cGMod && \cMod
	\arrow[""{name=0, anchor=center, inner sep=0}, "{{\orb}}", curve={height=-18pt}, from=1-1, to=1-3]
	\arrow[""{name=0p, anchor=center, inner sep=0}, phantom, from=1-1, to=1-3, start anchor=center, end anchor=center, curve={height=-18pt}]
	\arrow[""{name=1, anchor=center, inner sep=0}, "{{\fix}}"', curve={height=18pt}, from=1-1, to=1-3]
	\arrow[""{name=1p, anchor=center, inner sep=0}, phantom, from=1-1, to=1-3, start anchor=center, end anchor=center, curve={height=18pt}]
	\arrow[""{name=2, anchor=center, inner sep=0}, "\triv"{description}, from=1-3, to=1-1]
	\arrow[""{name=2p, anchor=center, inner sep=0}, phantom, from=1-3, to=1-1, start anchor=center, end anchor=center]
	\arrow[""{name=2p, anchor=center, inner sep=0}, phantom, from=1-3, to=1-1, start anchor=center, end anchor=center]
	\arrow[draw=none, from=1-3, to=1-1]
	\arrow["\dashv"{anchor=center, rotate=-85}, draw=none, from=0p, to=2p]
	\arrow["\dashv"{anchor=center, rotate=-95}, draw=none, from=2p, to=1p]
\end{tikzcd},
\end{center}
i.e., there are isomorphisms
\[
    \ihom_R(\orb M,A) \simeq \ihom_G(M,A_\triv)
\]
and
\[
    \ihom_G(A_\triv,M) \simeq \ihom_R(A,\fix M)
\]
natural in $A\in\cMod$ and $M\in\cGMod$.
\end{proposition}

\begin{proof}
The first isomorphism follows from the definition of $\orb M$,
since then $\ihom_R(\orb M,A)$ is given by the equalizer of the two arrows
$\ihom_R(M,A)\to\ihom_R(\mcM(G)\otimes_R M,A)$.
But this is exactly $\ihom_G(M,A_\triv)$ by definition.

The second isomorphism follows from the description of $\ihom_G(A_\triv,M)$
as the equalizer of the two parallel arrows
$\ihom_R(A, M)\to\ihom_R(\mcM(G)\otimes_R A,M)$
where we use the trivial action on $A$.
But this is then, by the description of $\fix M$ as an equalizer,
precisely $\ihom_R(A,\fix M)$.
\end{proof}

The proposition also implies that the fixed point functor is left exact
and the orbit functor is right exact.
Therefore it is a natural question to ask for the existence of the right resp. left 
derivation of those functors.
So we need certain conditions on the categories involved to guarantee the existence of those.
First of all, the category $\cMod$ has enough compact projective objects.
There, the objects $\mcM(S)$ for $S\in\extr$ are compact projective generators.
Now we have a free construction on $G$, i.e.,
there exists a left adjoint to the forgetful functor $?\colon\cGMod\to\cMod$.
Namely, it is given by $M\mapsto \mcM(G)\otimes_1 M$.
This left adjoint now maps compact projective objectives in $\cMod$
to compact projectives in $\cGMod$.
Moreover, by a formal argument, the objects $\mcM(G)\otimes_1\mcM(S)$ are
also generating in $\cGMod$, since any epimorphism $\bigoplus_i\mcM(S_i)\to ?M$ in $\cMod$
for $M\in\cGMod$
translates to an epimorphism $\bigoplus_i\mcM(G)\otimes_1\mcM(S_i)\to M$.
In particular, $\cGMod$ has enough compact projectives.
Moreover, direct sums in $\cGMod$ exist and are exact.
So are filtered colimits.
This is due to the fact that all (co)limits can be computed in $\cMod$.
This allows us to derive the (internal) $\hom$-functor and the tensor product in $\cGMod$
on the unbounded derived category and justifies the following definition.

\begin{definition}\label{def:group-cohomology}
The \idx{condensed group cohomology} of $G$ is the right derived functor of $\fix$,
\[
    R\fix\colon D(\cGMod)\to D(\cMod) = D(\RMod)_{\mcM}.
\]
We denote by $H^i(G,M)$ its $i$-th cohomology group
(or rather $\mcM$-complete condensed $R$-module) of $R\fix(M)$.
Dually, the \idx{condensed group homology} of $G$ is the left derived functor of $\orb$,
\[
    L\orb\colon D(\cGMod)\to D(\cMod).
\]
We denote by $H_i(G,M)$ its $i$-th cohomology group.
\end{definition}

For the rest of this chapter, we will focus on group cohomology.
The results for its dual concept then can be carried over verbatim.

\subsubsection{Computing group cohomology}

We give a method of calculating the right derived functor $R\fix$.
Since $\fix \simeq \ihom_G(R_\triv,-)$ it suffices to derive the enriched $\hom$-functor
$\ihom_G(-,-)$ in the contravariant argument at $R_\triv$.
Since $R_\triv$ is in general not projective in $\cGMod$,
we first need to resolve $R_\triv$ by a compact projective resolution.
In a first step, we construct two important resolutions.
The first of them is the so-called \emph{simplicial resolution}.

\begin{proposition}[Simplicial resolution]\label{prop:simp-res}
There exists a quasiisomorphism $S^\bullet\to R_\triv$ where $S^\bullet$ is given by
\begin{center}
\begin{tikzcd}
	\cdots & {\mcM(G^n)} & \cdots & {\mcM(G^2)} & {\mcM(G)} & 0
	\arrow[from=1-1, to=1-2]
	\arrow[from=1-2, to=1-3]
	\arrow[from=1-3, to=1-4]
	\arrow[from=1-4, to=1-5]
	\arrow[from=1-5, to=1-6]
\end{tikzcd}
\end{center}
in $\cGMod$.
\end{proposition}

\begin{proof}
The action of $G$ on each entry $\mcM(G^n)$ is given by simultaneous left multiplication
of $G$ on $G^n$.
The augmentation map $\mcM(G)\to R_\triv$ is the image of the zero map $0\colon G\to\ast$
under the functor $\mcM(-)$.
For $n\geq 2$, the morphism $\mcM(G^n)\to\mcM(G^{n-1})$ is the $R$-linear extension
of the morphism
\[
    G^n\to\mcM(G^{n-1}), \quad (g_1,\ldots,g_n)\mapsto \sum_{i=0}^n (-1)^i (g_1,\ldots,\hat{g}_i,\ldots,g_n)
\]
where $(g_1,\ldots,\hat{g}_i,\ldots,g_n) = (g_1,\ldots,g_{i-1},g_{i+1},\ldots,g_n)\in G^{n-1}$.
(Abstractly, it is $\sum_{i=0}^n (-1)^i \hat{p}_i$ where $\hat{p}_i\colon G^n\to G^{n-1}$ is the projection leaving out the $i$-th component).
Then one easily checks that this is an exact sequence in $\cGMod$.
\end{proof}

The construction reveals why it is called simplicial:
forgetting the linear structure, it is given by the (condensed) simplicial set $(G^n)$
which is regarded as an infinite dimensional simplex with $n$-simplex $G^n$.
The boundary maps are then the projections on each of its subsimplices $G^{n+1}\subset G^n$.
The second resolution is the so-called \emph{Bar construction} or \emph{Bar resolution}.

\begin{proposition}[Bar resolution]\label{prop:bar-res}
There exists a quasiisomorphism $B^\bullet\to R_\triv$ where $B^\bullet$ is given by
\begin{center}
\begin{tikzcd}
	\cdots & {\mcM(G)\otimes_1\mcM(G^{n-1})} & \cdots & {\mcM(G)\otimes_1\mcM(G)} & {\mcM(G)} & 0
	\arrow[from=1-1, to=1-2]
	\arrow[from=1-2, to=1-3]
	\arrow[from=1-3, to=1-4]
	\arrow[from=1-4, to=1-5]
	\arrow[from=1-5, to=1-6]
\end{tikzcd}
\end{center}
in $\cGMod$.
\end{proposition}

\begin{proof}
The augmentation map $\mcM(G)\to R_\triv$ is given by the zero morphism
$0\colon G\to *$. This clearly respects the $G$-action.
The morphism $\delta_1\colon\mcM(G)\otimes_1\mcM(G)\to\mcM(G)$
is induced via the adjunction of proposition~\ref{prop:G-mod-adjunction-?}
by the map $\mcM(G)\to\mcM(G)$ given by $1_G - R[0]$.
Then we have that $\delta_0\circ\delta_1 = 0$ since $R[0]\circ(1_G-R[0]) = R[0]-R[0] = 0$.
Now for $n\geq 2$,
the morphism $\delta_n\colon\mcM(G)\otimes_1\mcM(G^n)\to\mcM(G)\otimes_1\mcM(G^{n-1})$
is given under the adjunction by the map
$\mcM(G^n)\to\mcM(G)\otimes_1\mcM(G^{n-1})$ defined by
\[
    1_{G^n} + \sum_{i=1}^{n-1}(-1)^i(R[0]\otimes\partial_i) + (-1)^n(R[0]\otimes R[p_n]),
\]
where $\partial_i\colon\mcM(G^n)\to\mcM(G^{n-1})$ is the induced map of
$G^n\to G^{n-1}$ given by multiplication of the $i$-th with the $i+1$-th factor
and letting the other coordinates invariant,
and $p_n\colon G^n\to G^{n-1}$ given by projecting onto the first $n-1$ coordinates.
It is easily checked that the constructed chain complex is exact.
\end{proof}

\begin{remark}
Note that in the case where $\mcM$ is the trivial analytic ring
we could have also proved the preceding two proposition using the fact
that both the Bar and the simplicial resolution exist for (abstract) group modules $RG$
and are functorial in $G$ and $R$.
Hence they exist in any presheaf topos and since sheafification is exact,
they also exist in $\GMod$.
However if $\mcM$ is not trivial, completion is not exact anymore and we have to construct the resolutions by hand.
\end{remark}

The two resolutions are closely related.

\begin{proposition}\label{prop:bar-simp}
There exists a chain isomorphism $S^\bullet\to B^\bullet$.
\end{proposition}

\begin{proof}
For $n\in\N_0$, the map $S^n\to B^n$ is the morphism $\mcM(G^n)\to\mcM(G)\otimes_1\mcM(G^{n-1})$ which is induced by the morphism
\[
    G^n\to G^n,\quad (g_1,\ldots,g_n)\mapsto (g_1,g_1^{-1}g_2,\ldots,g_1^{-1}g_n)
\]
under the identification $\mcM(G)\otimes_R\mcM(G^{n-1}) = \mcM(G^n)$.
Its inverse $B^n\to S^n$ is the morphism $\mcM(G)\otimes_1\mcM(G^{n-1})\to\mcM(G^n)$
induced by the morphism
\[
    G^{n-1}\to G^n, \quad (g_1,\ldots,g_{n-1})\mapsto (e_G,g_1,g_1g_2,\ldots,g_1\cdots g_{n-1}).
\]
Then one checks that both morphisms are actually $G$-equivariant chain morphisms
and inverse to each other.
\end{proof}

In the abstract world, those two resolutions are the key to group cohomology
since they are projective.
But we are dealing with topological structures,
so the Bar and the simplicial resolution need not be projective resolutions.
Thus we cannot compute group cohomology by just applying $\ihom_G(-,M)$
pointwise on either of the chain complexes.
Instead we choose a Cartan-Eilenberg resolution
\begin{center}
\begin{tikzcd}
	& \vdots && \vdots & \vdots \\
	\cdots & {P_2^{-n}} & \cdots & {P_2^{-1}} & {P_2^0} & 0 \\
	\cdots & {P_1^{-n}} & \cdots & {P_1^{-1}} & {P_1^0} & 0 \\
	\cdots & {P_0^{-n}} & \cdots & {P_0^{-1}} & {P_0^0} & 0 \\
	\cdots & {\mcM(G)\otimes_1\mcM(G^n)} & \cdots & {\mcM(G)\otimes_1\mcM(G)} & {\mcM(G)} & 0
	\arrow[from=1-2, to=2-2]
	\arrow[from=1-4, to=2-4]
	\arrow[from=1-5, to=2-5]
	\arrow[from=2-1, to=2-2]
	\arrow[from=2-2, to=2-3]
	\arrow[from=2-2, to=3-2]
	\arrow[from=2-3, to=2-4]
	\arrow[from=2-4, to=2-5]
	\arrow[from=2-4, to=3-4]
	\arrow[from=2-5, to=2-6]
	\arrow[from=2-5, to=3-5]
	\arrow[from=3-1, to=3-2]
	\arrow[from=3-2, to=3-3]
	\arrow[from=3-2, to=4-2]
	\arrow[from=3-3, to=3-4]
	\arrow[from=3-4, to=3-5]
	\arrow[from=3-4, to=4-4]
	\arrow[from=3-5, to=3-6]
	\arrow[from=3-5, to=4-5]
	\arrow[from=4-1, to=4-2]
	\arrow[from=4-2, to=4-3]
	\arrow[from=4-2, to=5-2]
	\arrow[from=4-3, to=4-4]
	\arrow[from=4-4, to=4-5]
	\arrow[from=4-4, to=5-4]
	\arrow[from=4-5, to=4-6]
	\arrow[from=4-5, to=5-5]
	\arrow[from=5-1, to=5-2]
	\arrow[from=5-2, to=5-3]
	\arrow[from=5-3, to=5-4]
	\arrow[from=5-4, to=5-5]
	\arrow[from=5-5, to=5-6]
\end{tikzcd}
\end{center}
of the Bar sequence.
This means a double complex $P_\bullet^\bullet$ in $\cGMod$
such that all the $P_i^j$ are projective in $\cGMod$, $P_i^j = 0$ for $i < 0$ or $j > 0$,
the rows are chain complexes and each column $P_\bullet^{-n}$ is a projective resolution
of $\mcM(G)\otimes_1\mcM(G^n)$.
In particular, the induced morphism $\tot(P_\bullet^\bullet)\to B^\bullet$
from the additive total complex to the Bar resolution is a quasiisomorphism,
i.e., it induces isomorphisms $H^i(\tot P_\bullet^\bullet) \simeq H^i(B^\bullet)$ for all $i\in\N_0$.
Since $B^\bullet\to R_\triv$ is a quasiisomorphism,
the composition $\tot P^\bullet_\bullet\to R_\triv$ also is a quasiisomorphism.
By definition of the right derived functor, the value of $R\fix M = R\ihom_G(R_\triv,M)$
is given by the pointwise application of $\ihom_G(-,M)$ to $\tot P^\bullet_\bullet$.

On the other hand, applying $\ihom_G(-,M)$ pointwise to the double complex yields the double complex
$Q_\bullet^\bullet = \ihom_G(P_\bullet^\bullet,M)$.
Since $\ihom_G(-,M)$ preserves finite sums, the total complex of $Q_\bullet^\bullet$
is equal to the pointwise application of $\ihom_G(-,M)$ to the total complex $\tot P_\bullet^\bullet$ of $P_\bullet^\bullet$.
Now the (upwards and leftwards oriented) spectral sequence induced from the complex

\begin{adjustbox}{max width =\textwidth}
\begin{tikzcd}
	& \vdots && \vdots & \vdots \\
	\cdots & {\ihom_G(P_2^{-n},M)} & \cdots & {\ihom_G(P_2^{-1},M)} & {\ihom_G(P_2^0,M)} & 0 \\
	\cdots & {\ihom_G(P_1^{-n},M)} & \cdots & {\ihom_G(P_1^{-1},M)} & {\ihom_G(P_1^0,M)} & 0 \\
	\cdots & {\ihom_G(P_0^{-n},M)} & \cdots & {\ihom_G(P_0^{-1},M)} & {\ihom_G(P_0^0,M)} & 0 \\
	& 0 && 0 & 0
	\arrow[from=2-2, to=1-2]
	\arrow[from=2-2, to=2-1]
	\arrow[from=2-3, to=2-2]
	\arrow[from=2-4, to=1-4]
	\arrow[from=2-4, to=2-3]
	\arrow[from=2-5, to=1-5]
	\arrow[from=2-5, to=2-4]
	\arrow[from=2-6, to=2-5]
	\arrow[from=3-2, to=2-2]
	\arrow[from=3-2, to=3-1]
	\arrow[from=3-3, to=3-2]
	\arrow[from=3-4, to=2-4]
	\arrow[from=3-4, to=3-3]
	\arrow[from=3-5, to=2-5]
	\arrow[from=3-5, to=3-4]
	\arrow[from=3-6, to=3-5]
	\arrow[from=4-2, to=3-2]
	\arrow[from=4-2, to=4-1]
	\arrow[from=4-3, to=4-2]
	\arrow[from=4-4, to=3-4]
	\arrow[from=4-4, to=4-3]
	\arrow[from=4-5, to=3-5]
	\arrow[from=4-5, to=4-4]
	\arrow[from=4-6, to=4-5]
	\arrow[from=5-2, to=4-2]
	\arrow[from=5-4, to=4-4]
	\arrow[from=5-5, to=4-5]
\end{tikzcd}
\end{adjustbox}

converges to the cohomology of the total complex $\tot Q_\bullet^\bullet$.
But this cohomology is just $H^*(G,M)$
since $R\fix M = \ihom_G(\tot P^\bullet_\bullet,M)$.
How does the spectral sequence look like?
Taking cohomology w.r.t. the columns, one obtains the first page $E^{\bullet,\bullet}_1$
of the spectral sequence.
Each column $P_\bullet^{-n}$ is a projective resolution of $\mcM(G)\otimes_1\mcM(G^n)$
and hence the cohomology of the induced sequence $\ihom_G(P_\bullet^{-n},M)$
is just $R^*\ihom_G(\mcM(G)\otimes_1\mcM(G^n),M) = \iExt_G^*(\mcM(G)\otimes_1\mcM(G^n),M)$.
This implies that the first page $E^{\bullet,\bullet}_1$ is given by
\begin{center}
\begin{tikzcd}
	& \vdots & \vdots \\
	\cdots & {\iExt_G^2(\mcM(G)\otimes_1\mcM(G),M)} & {\iExt_G^2(\mcM(G),M)} & 0 \\
	\cdots & {\iExt_G^1(\mcM(G)\otimes_1\mcM(G),M)} & {\iExt_G^1(\mcM(G),M)} & 0 \\
	\cdots & {\iExt_G^0(\mcM(G)\otimes_1\mcM(G),M)} & {\iExt_G^0(\mcM(G),M)} & 0
	\arrow[from=2-2, to=2-1]
	\arrow[from=2-3, to=2-2]
	\arrow[from=2-4, to=2-3]
	\arrow[from=3-2, to=3-1]
	\arrow[from=3-3, to=3-2]
	\arrow[from=3-4, to=3-3]
	\arrow[from=4-2, to=4-1]
	\arrow[from=4-3, to=4-2]
	\arrow[from=4-4, to=4-3]
\end{tikzcd}
\end{center}
The adjunction $(\mcM(G)\otimes_1-)\dashv ?$ over the enrichment of $\cGMod$ translates
to an adjunction $R\ihom_G(\mcM(G)\otimes_1 A,M) = R\ihom_R(A,?M)$ for $A\in\cMod$, $M\in\cGMod$,
and hence the spectral sequence reduces to

\begin{adjustbox}{max width=\textwidth}
\begin{tikzcd}
	& \vdots && \vdots & \vdots \\
	\cdots & {\iExt_R^2(\mcM(G^n),M)} & \cdots & {\iExt_R^2(\mcM(G),M)} & {\iExt_R^2(R,M)} & 0 \\
	\cdots & {\iExt_R^1(\mcM(G^n),M)} & \cdots & {\iExt_R^1(\mcM(G),M)} & {\iExt_R^1(R,M)} & 0 \\
	\cdots & {\iExt_R^0(\mcM(G^n),M)} & \cdots & {\iExt_R^0(\mcM(G),M)} & {\iExt_R^0(R,M)} & 0
	\arrow[from=2-2, to=2-1]
	\arrow[from=2-3, to=2-2]
	\arrow[from=2-4, to=2-3]
	\arrow[from=2-5, to=2-4]
	\arrow[from=2-6, to=2-5]
	\arrow[from=3-2, to=3-1]
	\arrow[from=3-3, to=3-2]
	\arrow[from=3-4, to=3-3]
	\arrow[from=3-5, to=3-4]
	\arrow[from=3-6, to=3-5]
	\arrow[from=4-2, to=4-1]
	\arrow[from=4-3, to=4-2]
	\arrow[from=4-4, to=4-3]
	\arrow[from=4-5, to=4-4]
	\arrow[from=4-6, to=4-5]
\end{tikzcd}
\end{adjustbox}

with $\iExt_R^i(R,M) = 0$ if $i>0$ and $\iExt_R^0(R,M) = M$.
Clearly, if $M\to N$ is a morphism of $G$-modules,
the constructed spectral sequences translate into each other since the construction
was independent of $M$.
To summarize our findings, we have just proven the following.

\begin{corollary}\label{cor:group-cohomology-spec-seq}
There exists a spectral sequence
\[
    E_1^{p,q} = \iExt_R^p(\mcM(G^q),M)
        \Rightarrow H^{p+q}(G,M)
\]
functorial in $M$.
\end{corollary}

All of the cohomology groups $H^i(G,M)$ are actually condensed $R$-modules,
so one is interested in its sections $H^i(G,M)(S)$ for $S\in\extr$.
First of all, we have
\[
    (\fix M)(S) = \hom_R(\mcM(S),\fix M) = \hom_G(\mcM(S)_\triv,M).
\]
In particular, $(R\fix M)(S) = R\hom_G(\mcM(S)_\triv,M)$.
Applying $-\otimes_1\mcM(S)$ to the Bar resolution yields a resolution of
$\mcM(S)_\triv = R_\triv\otimes_1\mcM(S)$ functorial in $S$,
since the functor $-\otimes_1\mcM(S)$ is exact
(we can check exactness in $\cMod$ and there it is just tensoring with $\mcM(S)$
which is projective for $S\in\extr$).
Now the same reasoning as above gives rise to the corresponding spectral sequence
\[
    E_1^{p,q} = \Ext_R^p(\mcM(G^q\times S),M)\Rightarrow H^{p+q}(G,M)(S)
\]
which is functorial in $M$ and $S$.
Furthermore using the adjunction $\hom_R(\mcM(X),A) = \hom(X,?A)$ gives the spectral sequence
\[
    E_1^{p,q} = H_\mathrm{cond}^p(G^q\times S,M)\Rightarrow H^{p+q}(G,M)(S).
\]
Note that here the $\mcM$-completeness on the left hand side is encoded in the condensed cohomology groups $H_\mathrm{cond}^p(G^q\times S,M)$
since $M$ is $\mcM$-complete.
These cohomology groups are fairly reasonable if $G$ is a compact group and $M$ is discrete.
Then they coincide with sheaf cohomology and can be computed,
and in turn the spectral sequence yields the values of the condensed group cohomology groups.
But this will be a task for future investigations.
For now let's turn back to ordinary group cohomology.

\subsubsection{Comparison with continuous group cohomology}

To be able to compare condensed group cohomology with continuous group cohomology,
we give a short overview over the definition.

All topological spaces are assumed to have closed points, i.e., they are $T_1$-spaces
(for linear spaces and groups this is equivalent to being Hausdorff).
This is due to the need to embed the spaces into the condensed setting.
We are working in the category of topological $R$-modules for a topological ring $R$.
Let $G$ be a topological group and $M$ a topological $G$-module.

For $n\in\N$ denote by $C^n(G,M)$ the topological $G$-module given by
all continuous maps $G^{n+1}\to M$ (note the shift in dimension) equipped with
the compact open topology and the induced $G$-action of $G^n$ and $M$.
These form a chain complex via the boundary maps
\[
    \partial^n\colon C^n(G,M)\to C^{n-1}(G,M),\quad \partial^nf(g_0,\ldots,g_n) = \sum_{i=0}^n (-1)^i(g_0,\ldots,\hat{g}_i,\ldots,g_n).
\]
Then the embedding $M\to C^0(G,M)$ mapping $M$ to the constant map gives an exact sequence
\begin{center}
\begin{tikzcd}
	0 & M & {C^0(G,M)} & {C^1(G,M)} & {C^2(G,M)} & \cdots
	\arrow[from=1-1, to=1-2]
	\arrow[from=1-2, to=1-3]
	\arrow[from=1-3, to=1-4]
	\arrow[from=1-4, to=1-5]
	\arrow[from=1-5, to=1-6]
\end{tikzcd}
\end{center}
of topologically injective $G$-modules (whatever that means - it is not relevant for us).
Then taking pointwise the fixed space yields a chain complex $C^\bullet(M)$
\begin{center}
\begin{tikzcd}
	0 & {\fix C^0(G,M)} & {\fix C^1(G,M)} & {\fix C^2(G,M)} & \cdots
	\arrow[from=1-1, to=1-2]
	\arrow[from=1-2, to=1-3]
	\arrow[from=1-3, to=1-4]
	\arrow[from=1-4, to=1-5]
\end{tikzcd}.
\end{center}
This assignment $M\mapsto C^\bullet(M)$ is functorial in $M$.
Now the $i$-th \idx{continuous group cohomology} group $H^i_c(G,M)$ is defined
to be the $i$-th cohomology group of the chain complex $C^\bullet(M)$
where each $H^i_c(G,M)$ carries the quotient topology
(assuming that it is again Hausdorff).

First of all, this description is very explicit.
This is because working with topological $G$-modules does not bring with it the best
categorical properties (they do not form an abelian category)
and so one has to mimic the definition from the discrete case by hand.
Also note the similarity to the simplicial resolution of $R_\triv$.
Actually, we are resolving the point $*$ by powers $G^n$ and then applying the functor
$C(-,M)$.
This gives a resolution of $M$ by $G$-modules.

Nevertheless we can make use of it to compare continuous group cohomology
with condensed group cohomology.

\begin{proposition}\label{prop:cond-cont-group-cohom}
Let $R = \mcM(*)$ be a topological ring and $M$ a topological $R$-module
such that $\underline{M}$ is $\mcM$-complete.
Let $G$ be a topological group such that every power $G^n$ for $n\in\N$ is compactly generated,
and let $G$ act continuously on $M$
Then there exists a morphism
\[
    \underline{C}^\bullet(M) \to R\fix(\underline{M}).
\]
natural in $M$.
\end{proposition}

\begin{proof}
The condition on $G$ guarantees that
\[
    \underline{C^n(G,M)} = \underline{C(G^{n+1},M)} = \ihom(\underline{G}^{n+1},\underline{M})
\]
as condensed sets and furthermore
\[
    \underline{C^n(G,M)} = \hom^G(\mcM(\underline{G}^{n+1}),\underline{M})
\]
as $\mcM$-complete $\underline{G}$-modules.
Since limits can be computed pointwise and are preserved by condensation,
we obtain
\[
    \underline{C}^n(M) = \underline{\fix C^n(G,M)} = \fix\underline{C^n(G,M)} = \ihom_G(\mcM(\underline{G}^{n+1}),\underline{M}).
\]
The boundary maps of $C^\bullet(-)$ are precisely the ones from the simplicial resolution,
so $\underline{C}^\bullet(M) \simeq \ihom_G(S^\bullet,\underline{M})$.
Choosing a projective resolution $P^\bullet\to S^\bullet\to \underline{R}_\triv$ then yields the desired morphism
upon applying $\ihom_G(-,\underline{M})$.
\end{proof}

The question now is: when is this map an isomorphism?
Let us first consider the case when $G$ is discrete.
Then for every $n\in\N$ the free $R$-modules $\mcM(G^n)$ are projective in $\cMod$.
Hence the Bar resolution and so the simplicial, are projective resolutions of $R_\triv$.
In this case, $R\fix(\underline{M})$ is obtained by taking the cohomology of the complex $\ihom_G(S^\bullet,\underline{M})$
which we have just showed to be $\underline{C}^\bullet(M)$.
Since evaluation at the point is exact, for any topological $G$-module $M$ we have
\[
    H^n_c(G,M) \simeq H^n(\underline{G},\underline{M})(*)
\]
for all $n\in\N$.

\begin{question}
If the quotient topology of $H^n_c(G,M)$ happens to be Hausdorff, does then $\underline{H^n_c(G,M)} \simeq H^n(\underline{G},\underline{M})$ hold as well?
\end{question}

Another interesting case is that when $G$ is a compact group.
Then in $\CAb$ it holds that
\[
    H^n_c(G,M) \simeq H^n(\underline{G},\underline{M})(*)
\]
for all topological $\R$-vector spaces $M$, see \cite{Brink2023}.
This should remain true in the category of liquid vector spaces for $M$ a complete locally convex complete vector space
and thus allows for investigations in representations of $G$ on locally convex spaces.

For more comparisons with continuous group cohomology see \cite{Brink2023}.

\subsubsection{Cohomology for dynamical systems}

\begin{center}
\textit{\enquote{Every mathematician is a category theorist, some just don't know yet.}}
\end{center}

\noindent To end this section we show how one can use condensed group cohomology in ergodic theory.
The usage of cocycles and extensions in ergodic theory goes back to \cite{Helson1960} where the equation for a cocycle pops out of a functional equation.
Over the years, cocycles have been proven to be helpful in the construction of ergodic systems out of easier ones,
see \cite{Helson1965, Ellis1965, Parry1969, Jones1972, Parry1974} and \cite{Parry1975}.
Over the years it has been recognized that not only do cocycles give rise to ergodic systems, but they are in a way the only
extensions of a given system by a certain group, see \cite{Zimmer1976} and \cite{Ellis1987}.
This is the so-called \emph{Mackey-Zimmer theorem} we introduce in this chapter.
Closely related is the so called Moore-Schmidt theorem which classifies coboundaries with values in the circle group $\T$.
Both of them are still relevant today, see e.g. \cite{jamneshan2022uncountablemackeyzimmertheorem, Jamneshan2020, jamneshan2022uncountable}.

Despite using terminology from homological algebra, there seems to be no systematic usage of the methods thereof.
For discrete dynamical systems, T. Tao explains in \cite{Tao2008} the origin of the term \emph{cocycle} by means of homological algebra.
Our goal is to expand this example to non-discrete topological dynamical systems (and hence to ergodic theory) with the aid of condensed mathematics,
where we can now use homological algebra on objects that have a topological structure.
This last part is a little exhibition of what might be possible and the first step into that direction.

The following is the definition of cocycles in topological dynamics.

\begin{definition}\label{def:top-dyn-cocycle}
Let $G$ be a (discrete) group that acts on a topological space $X$, and let $U$ be a topological group.
A \idx{continuous cocycle} is a continuous function $\rho\colon G\times X\to U$ such that  
\[
    \rho(hg,x)=\rho(g,x)+\rho(h,g.x)
\]
for all $x\in X$ and $g,h\in G$.
A continuous cocycle $\rho$ is said to be a \idx{continuous coboundary} if there is a continuous function $f\colon X\to U$
such that
\[
    \rho(g,x)=f(g.x)-f(x).
\]
\end{definition}

The set of all cocycles is clearly a group under pointwise addition and
the first {cohomology group} is defined to be the factor of this group by the subgroup of all coboundaries.
Now given a topological dynamical system $G\curvearrowright X$ and an $U$-valued cocycle $\rho$,
we can build a new topological dynamical system.
Its underlying set is the product $U\times X$ and now the action of $G$ on $U\times X$ is given by
\[
    G\times U\times X\to U\times X,\quad (g,u,x)\mapsto (u+\rho(g,x),g.x).
\]
Now the cocycle equation ensures that this actually is a dynamical system
(and this is how often the cocycle equation is justified - not really helpful).
We denote this dynamical system by $U\rtimes_\rho X$.
The group $U$ also acts by left translation on $U\rtimes_\rho X$.
This translation behaves well with the action of $G$:
we have $u.(g.x) = g.(u.x)$ since the action of $G$ is independent of the value at $U$.
One can quotient out the orbits of this action, so that $\orb_U(U\rtimes_\rho X) = X$
and the quotient map $U\rtimes_\rho X\to X$ is $G$-equivariant.
This leads to the definition of an \emph{extension}.

\begin{definition}\label{def:top-dyn-ext}
Let $X$ and $Y$ be topological systems over $G$ and $Y\to X$ a $G$-equivariant map.
Then $Y$ is an \idx{extension} of $X$ by $U$ if
\begin{enumerate}[(i)]
    \item the group $U$ acts faithful on $Y$, i.e., $u.y = y$ for all $y\in Y$ implies $u = e_U$,
    \item for all $y\in Y$, we have $g.(u.y) = u.(g.y)$ for every $g\in G$, $u\in U$,
    \item each fibre of the given map $Y\to X$ is precisely an orbit of $U$ in $Y$.
\end{enumerate}
\end{definition}

The last condition asserts that the orbit space of $Y$ by $U$ is isomorphic to $X$.
Now the \idx{Mackey-Zimmer theorem} relates those two concepts.

\begin{theorem}[Mackey-Zimmer]\label{thm:top-dyn-mackey-zimmer}
Let $X$ be extremally disconnected compact Hausdorff and $U$ a compact group.
If $Y\to X$ is an extension of $X$ by $U$ over $G$ such that $U$ acts \emph{freely} on the compact Hausdorff space $Y$,
then there exists a continuous cocycle $\rho$ such that the diagram
\begin{center}
\begin{tikzcd}
	{U\rtimes_\rho X} && Y \\
	& X
	\arrow["\simeq", from=1-1, to=1-3]
	\arrow[from=1-1, to=2-2]
	\arrow[from=1-3, to=2-2]
\end{tikzcd}
\end{center}
is commutative where the horizontal arrow is an isomorphism.

Two cocycle extensions $U\rtimes_\rho X$ and $U\rtimes_\eta X$ are isomorphic over $X$
precisely if $\rho$ and $\eta$ are equivalent in the first cohomology group.
\end{theorem}

See Theorem 4.9 in \cite{Ellis1987} for a proof.
See also \cite{jamneshan2022uncountablemackeyzimmertheorem} for a slightly more general and ergodic theoretic formulation.

Another important theorem that investigates the first cohomology group is the \idx{Moore-Schmidt theorem}.
It goes as follows.

\begin{theorem}[Moore Schmidt]\label{thm:top-dyn-moore-schmidt}
Let $X$ be an extremally disconnected compact Hausdorff topological system over $G$
and $U$ be a compact abelian group.
Then for an $U$-valued cocycle $\rho$ the following are equivalent.
\begin{enumerate}[(a)]
    \item The cocycle $\rho$ is a coboundary.
    \item For every character $\chi\in\hat{U}$, the cocycle $\chi\circ\rho$ is a coboundary.
\end{enumerate}
\end{theorem}

See Theorem 1.6 in \cite{jamneshan2022uncountable}.

Having reviewed the most important notions and results,
we show how they adapt to our setting.
We start with the notion of cocycle.

Let $G$ be a discrete group.
We have already seen that in this case, condensed group cohomology agrees with continuous group cohomology.
Moreover the Bar resolution is a projective resolution of $R_\triv$ so we can compute cohomology pointwise.
To obtain the notion of cocycle we also need a space $X$ on which $G$ acts.
For this, let $X$ be a compactly generated $T_1$-space and let $U$ be an $R$-module such that $\underline{U}$ is $\mcM$-complete
(e.g. a complete locally convex vector space in the case of the liquid ring,
or just a topological abelian group in the case of the trivial ring structure).
Then $\mcM(\underline{X})$ is in $\cGMod$ and
\[
    M = \hom^G(\mcM(\underline{X}),\underline{U}_\triv)
\]
has the action of $G$ via pullback on $X$ and therefore also is an object of $\cGMod$.
We want to compute the cohomology groups of $G$ with values in this module $M$.
To better keep track of all objects involved, we write $H^i(G,X,U)$ instead of $H^i(G,M)$.
As already explained, the condensed group cohomology with values in $M$ is given by
the chain complex
\begin{center}
\begin{tikzcd}
	0 & {\ihom_G(\mcM(G),M)} & {\ihom_G(\mcM(G)\otimes_1\mcM(G),M)} & \cdots
	\arrow[from=1-1, to=1-2]
	\arrow["{\delta_1^*}", from=1-2, to=1-3]
	\arrow["{\delta_2^*}", from=1-3, to=1-4]
\end{tikzcd}.
\end{center}
Applying the adjunction from proposition~\ref{prop:G-mod-adjunction-?},
this chain complex is isomorphic to the chain complex
\begin{center}
\begin{tikzcd}
	0 & M & {\ihom_R(\mcM(G),M)} & \cdots & {\ihom_R(\mcM(G^n),M)} & \cdots
	\arrow[from=1-1, to=1-2]
	\arrow["{d_1^*}", from=1-2, to=1-3]
	\arrow["{d_2^*}", from=1-3, to=1-4]
	\arrow["{d_{n-1}^*}", from=1-4, to=1-5]
	\arrow["{d_n^*}", from=1-5, to=1-6]
\end{tikzcd}
\end{center}
where now the boundary maps $d^*$ are given by the pullbacks of the maps constructed
in the proof of the Bar resolution.
We will come back to the boundary morphisms later.
Plugging in $M = \hom^G(\mcM(\underline{X}),\underline{U}_\triv)$ and using the hom-tensor adjunction
yields the complex
\begin{center}
\begin{tikzcd}
	0 & {\ihom_R(\mcM(\underline{X}),\underline{U})} & {\ihom_R(\mcM(\underline{G\times X}),\underline{U})} & \cdots
	\arrow[from=1-1, to=1-2]
	\arrow["{{d_1^*}}", from=1-2, to=1-3]
	\arrow["{{d_2^*}}", from=1-3, to=1-4]
\end{tikzcd}.
\end{center}
Using that $\mcM(-)$ is left adjoint then yields
\begin{center}
\begin{tikzcd}
	0 & {\ihom(\underline{X},\underline{U})} & {\ihom(\underline{G\times X},\underline{U})} & \cdots & {\ihom(\underline{G^n\times X},\underline{U})} & \cdots
	\arrow[from=1-1, to=1-2]
	\arrow["{{d_1^*}}", from=1-2, to=1-3]
	\arrow["{{d_2^*}}", from=1-3, to=1-4]
	\arrow["{{d_{n-1}^*}}", from=1-4, to=1-5]
	\arrow["{{d_n^*}}", from=1-5, to=1-6]
\end{tikzcd}.
\end{center}
Now since $G$ is discrete and $X$ is compactly generated, so is $G^n\times X$ for every $n\in\N$
and hence comparison of the internal hom with the classical objects with compact open topology yields the complex
\begin{center}
\begin{tikzcd}
	0 & {\underline{C(X,U)}} & {\underline{C(G\times X,U)}} & \cdots & {\underline{C(G^n\times X,U)}} & \cdots
	\arrow[from=1-1, to=1-2]
	\arrow["{{d_1^*}}", from=1-2, to=1-3]
	\arrow["{{d_2^*}}", from=1-3, to=1-4]
	\arrow["{{d_{n-1}^*}}", from=1-4, to=1-5]
	\arrow["{{d_n^*}}", from=1-5, to=1-6]
\end{tikzcd}
\end{center}
Now inspecting the boundary maps constructed for the Bar resolution more closely,
one can see that the first boundary map $d_1^*$ is given by mapping a function $f\colon X\to U$
to the function $(g,x)\mapsto f(g.x)-f(x)$.
The second boundary map $d_2^*$ is given by mapping a function $f\in\C(G\times X,U)$
to the function $(g,h,x)\mapsto f(h,g.x) - f(hg,x) + f(g,x)$,
and by proceeding further one can construct every boundary map in such a manner.

Now what is $H^i(G,X,U)$ evaluated at the point $*$? 
Since evaluation at the point is exact,
we can compute $H^i(G,X,U)(*)$ by first evaluating the chain complex at the point and then computing its cohomology groups.

For $i=0$ we obtain the submodule of $C(X,U)$ consisting of all function $f\colon X\to U$ such that $f(g.x) = f(x)$,
or equivalently, all function from the orbit of $X$ to $U$.

For $i=1$ we obtain the module that is the quotient of all functions $\rho\colon G\times X\to U$ such that $\rho(h,g.x) - \rho(hg,x) + \rho(g,x) = 0$,
or equivalently,
\[
    \rho(hg,x) = \rho(h,g.x) + \rho(g,x)
\]
(This is exactly the cocycle equation from definition~\ref{def:top-dyn-cocycle}!),
by the submodule of all functions of the form $\rho(x) = f(g.x) - f(x)$ for $f\colon X\to U$.

Thus we have reconstructed the first cohomology group of definition~\ref{def:top-dyn-cocycle} as an actual part of the first condensed cohomology
group of group cohomology with values in $M = \hom^G(\mcM(\underline{X}),\underline{U}_\triv)$.

Now what about the connection made in the Mackey-Zimmer theorem between cocycles and extensions?
In homological algebra, an \emph{extension of $C$ by $A$} is a short exact sequence of the form
\begin{center}
\begin{tikzcd}
	0 & A & B & C & 0
	\arrow[from=1-1, to=1-2]
	\arrow[from=1-2, to=1-3]
	\arrow[from=1-3, to=1-4]
	\arrow[from=1-4, to=1-5]
\end{tikzcd}.
\end{center}
So we want to bring these extensions in relation with the first cohomology group we have just presented.
Since Ext-groups classify extensions, see lemma~\ref{lem:ext_as_extensions}, we have the following trivial statement.

\begin{lemma}[Condensed Mackey-Zimmer]
Let $G$ be a discrete group, $X$ a compactly generated $T_1$-space with continuous $G$-action
and $U$ a $T_1$-topological $R$-module such that $\underline{U}$ is $\mcM$-complete.
Then the first cohomology group $H^1(G,X,U)$ classifies extensions
\begin{center}
\begin{tikzcd}
	0 & {U_\triv} & M & {\mcM(X)} & 0
	\arrow[from=1-1, to=1-2]
	\arrow[from=1-2, to=1-3]
	\arrow[from=1-3, to=1-4]
	\arrow[from=1-4, to=1-5]
\end{tikzcd}
\end{center}
of $\mcM(X)$ by $U$ in the sense that to every cocycle $[\rho]\in H^1(G,X,U)$
there is an extension $U_\triv\oplus_\rho\mcM(X)$ fitting in the middle of the exact sequence such that
two such extensions are isomorphic if and only if $[\rho] = [\eta]$.
On the other hand, every such extension arises in this way.
\end{lemma}

\begin{proof}
We already know that $\Ext^1$ classifies extensions, see lemma~\ref{lem:ext_as_extensions},
so we need only show that the short exact sequences have the form described in the lemma.
This is due to the fact that by corollary~\ref{cor:hom-fix} we have
\[
    R^*\fix\hom^G(\mcM(X),U_\triv) = \iExt_G^*(\mcM(X),U_\triv)
\]
which follows from the adjunction in corollary~\ref{cor:hom-fix}.
Then evaluating at the point yields $H^1(G,X,U) = \Ext_G^1(\mcM(X),U_\triv)$.
\end{proof}

\begin{remark}
This is a generalisation of the classical Mackey-Zimmer theorem in the following two (and many more) directions.
Note that we only impose on $X$ that it is compactly generated.
It need not be extremally disconnected compact Hausdorff
as is the case in the classical theorem.
Moreover such a theorem also generalizes to non-discrete groups $G$,
spaces that are not Hausdorff and much more, due to the general flexibility in the setup.
The continuity of the action is automatically encoded even for complicated $G$, where the classical cocycle equation does make no sense at all.

Second, we can vary the base ring of the modules and
may install completeness properties on them,
allowing to vary from abelian groups up to locally convex $\R$-vector spaces,
all in the same theorem.

One drawback is that $U$ needs to be an abelian group.
This is due to the fact that there is almost no account for group cohomology
with non-abelian coefficients.
This theorem (among others) can be a good starting point to develop such a theory.
\end{remark}

To close this section let us remark that we think that there is a translation of the Moore-Schmidt theorem
into the condensed setting whose proof also is just an abstract manipulation of symbols.
We are certain that by translating much of ergodic theory to the condensed setting one is finally able of unleashing the power
of homological algebra onto this area of mathematics.
Doing so one might be able to develop a strong structure theory within ergodic theory that
does not rely on tedious computations.
See for example \cite{Host2018, host2009, Glasner2003,candela2017,manners2018,manners2018ii, manners2020} for more ergodic theory we think is worth condensing.
Also, the articles \cite{Delvenne2019CategoryTF} and \cite{bader2023higherkazhdanpropertyunitary} contain further possible directions of research.
We are excited what the future will bring.

\cleardoublepage

\thispagestyle{empty}

\quot{
    I can illustrate the [\dots] approach with the [\dots] image of a nut to be opened.\\[\baselineskip]
The first analogy that came to my mind is of immersing the nut in some softening liquid, and why not simply water?\\[\baselineskip]
From time to time you rub so the liquid penetrates better, and otherwise you let time pass. The shell becomes more flexible through weeks and months -- when the time is ripe, hand pressure is enough, the shell opens like a perfectly ripened avocado![\dots]\\[\baselineskip]
A different image came to me a few weeks ago. The unknown thing to be known appeared to me as some stretch of earth or hard marl, resisting penetration [\dots] the sea advances insensibly in silence, nothing seems to happen, nothing moves, the water is so far off you hardly hear it [\dots] yet finally it surrounds the resistant susbstance.
}{A. Grothendieck p.552-3 in \cite{Grothendieck1985}, translation taken from \cite{Vakil2023}}

\clearpage{\thispagestyle{empty}\cleardoublepage}

\bibliographystyle{alpha_mod}
\bibliography{bibliography}



\section*{Overview over categories used}
\begin{center}
\rotatebox{90}{
\begin{tabular}{c | c | c}
    name & objects & morphisms \\
    \hline \hline
    $\Set$ & sets & functions \\
    $\Fin$ & finite sets & functions \\
    $\Delta$ & finite ordinals & monotone functions \\
    $\sSet$ & presheaves on $\Delta$ & natural transformations \\
    \textbf{Cat} & small categories & functors \\
    $\Top$ & topological spaces & continuous maps \\
    \textbf{CGWH} & compactly generated weak Hausdorff spaces & continuous maps \\
    $\CHaus$ & compact Hausdorff spaces & continuous maps \\
    $\stone$ & totally disconnected compact Hausdorff spaces & continuous maps \\
    $\prof$ & profinite sets & continuous maps \\
    $\bool$ & Boolean algebras & Boolean algebra homomorphisms \\
    $\extr$ & extremally disconnected compact Hausdorff spaces & continuous maps \\
    $\bdis$ & Stone \v{C}ech compactifications of discrete sets & continuous maps \\
    $\Grp$ & groups & group homomorphisms \\
    \textbf{Ab} & abelian groups & group homomorphisms \\
    \textbf{Ring} & associative rings with unit & ring homomorphisms \\
    $\ModR$ & left modules over $R$ & module homomorphisms \\
    $\LCA$ & locally compact Hausdorff abelian groups & continuous group homomorphisms \\
    $\EM(T)$ & Eilenberg-Moore category of monad $T$ & $T$-algebra homomorphisms \\
    $\cond$ & condensed sets & natural transformations \\
    $\qc$ & quasicompact condensed sets & natural transformations \\
    $\qs$ & quasiseparated condensed sets & natural transformations \\
    $\qcqs = \CHaus$ & quasicompact quasiseparated condensed sets & natural transformations \\
    $\CondRing$ & condensed rings & natural transformations \\
    $\CondMod$ & condensed modules over condensed rings & natural transformations \\
    $\RMod = \CondModR$ & condensed $R$-modules ($R\in\CondRing$) & natural transformations \\
    $\cond(\ModR)$ & condensation of $\ModR$ & natural transformations \\
    $\cab = \cond(\Mod_\Z)$ & condensed abelian groups & natural transformations \\
    $\cMod$ & complete condensed $R$-modules w.r.t. an analytic ring $\mcM$ & condensed module homomorphisms \\
    $\liq_p$ & $p$-liquid $\R$-vector spaces & condensed linear maps \\
    $\GMod$ & condensed $G$-representations & $G$-equivariant morphisms \\
    $\cGMod$ & condensed $G$-representations, complete w.r.t. an analytic ring $\mcM$ & $G$-equivariant morphisms \\
\end{tabular}}
\end{center}

\printindex

\end{document}